%% file: texte.tex
\documentclass[a4paper,11pt]{smfbook}
\usepackage{amssymb,amsmath,mathrsfs,graphics,graphicx,mathptm}
\usepackage[T1]{fontenc}
\usepackage[english, francais]{babel} 
\usepackage{float,lscape}
\usepackage[all]{xy}

\theoremstyle{plain}
\newtheorem{theoreme}{Th\'eor\`eme}
\newtheorem{corollaire}[theoreme]{Corollaire}

\theoremstyle{definition}
\newtheorem*{exemples}{Exemples}
\newtheorem{remarque}[theoreme]{Remarque}

\theoremstyle{plain}
\newtheorem{thm}{Th\'eor\`eme}[chapter]
\newtheorem{pro}[thm]{Proposition}
\newtheorem{lem}[thm]{Lemme}
\newtheorem{cor}[thm]{Corollaire}

\theoremstyle{definition}
\newtheorem*{defi}{D\'efinition}

\newtheorem*{eg}{Exemple}
\newtheorem*{egs}{Exemples}
\newtheorem{rem}[thm]{Remarque}
\newtheorem{rems}[thm]{Remarques}

\newcommand{\clearemptydoublepage}{%
    \newpage{\pagestyle{empty}\cleardoublepage}}

\def\noteA#1#2#3{{\begin{footnotesize}$#1$\end{footnotesize}} & {\begin{footnotesize}#2\end{footnotesize}} \hfill{\begin{footnotesize}\pageref{#3}\end{footnotesize}} \\}
\def\noteAB#1#2#3{{\begin{footnotesize}$#1$\end{footnotesize}} & {\begin{footnotesize}#2\end{footnotesize}} \hfill{} \\}
\def\noteAC#1#2#3{{} & {\begin{footnotesize}#2\end{footnotesize}} \hfill{\begin{footnotesize}\pageref{#3}\end{footnotesize}} \\}

\def\noteB#1#2{{\begin{small}#1\end{small}} \hfill {\begin{small}\pageref{#2}\end{small}} \\}

\def\transp #1{\vphantom{#1}^{\mathrm t}\! {#1}}

\def\og{\leavevmode\raise.3ex\hbox{$\scriptscriptstyle\langle\!\langle$~}}
\def\fg{\leavevmode\raise.3ex\hbox{~$\!\scriptscriptstyle\,\rangle\!\rangle$}}

\addtocounter{section}{0}             
\numberwithin{equation}{section}       

\begin{document}

\frontmatter

\title{Transformations birationnelles de petit degr\'e}

\author{Dominique \textsc{Cerveau}}

\address{Membre de l'Institut Universitaire de France.
IRMAR, UMR 6625 du CNRS, Universit\'e de Rennes $1,$ $35042$ Rennes, France  \\
Membre de l'ANR BLAN$06$-$3$\_$137237$}
\email{dominique.cerveau@univ-rennes1.fr}

\author{Julie \textsc{D\'eserti}}

\address{Institut de Math\'ematiques de Jussieu, Universit\'e Paris $7,$ Projet G\'eom\'etrie et Dynamique, Site Chevaleret, Case $7012,$ $75205$ Paris Cedex 13, France\\
Membre de l'ANR BLAN$06$-$3$\_$137237$}
\email{deserti@math.jussieu.fr}

\begin{abstract}
\selectlanguage{french}
Depuis la fin du XIX\`eme si\`ecle on sait que toute
transformation birationnelle du plan projectif complexe dans
lui-m\^eme, encore appel\'ee transformation de \textsc{Cremona},
s'\'ecrit comme la compos\'ee de transformations birationnelles
quadratiques; ceci a motiv\'e notre travail qui porte
essentiellement sur ces transformations. Nous \'etablissons des
propri\'et\'es de type alg\'ebriques comme la classification des
groupes \`a un param\`etre
de transformations de \textsc{Cremona} quadratiques ou encore la
lissit\'e de l'espace des transformations birationnelles de degr\'e
$2$ dans l'espace des transformations rationnelles. On peut aussi
voir qu'un nombre fini de transformations de \textsc{Cremona}
quadratiques choisies g\'en\'eriquement engendre un groupe libre.
Par ailleurs nous montrons que si $f$ est une transformation
birationnelle de degr\'e $2$ ou un automorphisme non trivial du
plan projectif complexe, le sous-groupe normal engendr\'e par $f$
est le groupe des transformations de \textsc{Cremona} tout entier;
nous en d\'eduisons que ce groupe est parfait. Nous d\'emontrons
aussi des propri\'et\'es de nature dynamique: en suivant une
id\'ee de \textsc{Guillot} nous implantons  aux transformations
birationnelles de degr\'e $2$ des invariants propres aux feuilletages
ce qui nous permet par exemple d'obtenir l'\'enonc\'e suivant: si
deux transformations de \textsc{Cremona} quadratiques
g\'en\'eriques sont birationnellement conjugu\'ees, elle le sont
lin\'eairement; nous nous int\'eressons aussi \`a la pr\'esence ou
non \og d'objets invariants\fg : courbes, feuilletages,
fibrations. Vient ensuite une partie plus exp\'erimentale: nous
tra\c{c}ons des orbites de transformations birationnelles de
degr\'e $2$ \`a coefficients r\'eels ainsi que des ensembles
analogues aux ensembles de \textsc{Julia} des polyn\^omes  \`a une
variable. Nous \'etudions aussi les transformations de
\textsc{Cremona} cubiques; en consid\'erant les diff\'erentes
configurations possibles de courbes contract\'ees nous en donnons
\og  la classification \fg. Ceci nous permet de montrer que
l'ensemble des transformations birationnelles exactement de
degr\'e $3$ est irr\'eductible, et en fait rationnellement
connexe.
\end{abstract}

\begin{altabstract}
\selectlanguage{english}
Since the end of the XIXth century, we know that each birational
map of the complex projective plane is the product of a
finite number of quadratic birational maps of the
projective plane; this motivates our work which essentially deals
with these quadratic maps. We establish algebraic
properties such as the classification of
one parameter groups of quadratic birational maps or
the smoothness of the set of quadratic birational maps
in the set of rational maps. We prove that a
finite number of generic quadratic birational maps
generates a free group. We show that if $f$ is a quadratic
birational map or an automorphism of the projective
plane, the normal subgroup generated by $f$ is the full group of
birational maps of the projective plane, which implies
that this group is perfect. We study some dynamical
properties: following an idea of \textsc{Guillot}, we translate some
invariants for foliations in our context, in particular we obtain
that if two generic quadratic birational maps are
birationally conjugated, then they are conjugated by an automorphism of
the projective plane. We are also interested in the presence of
"invariant objects": curves, foliations, fibrations. Then follows
a more experimental part: we draw orbits of quadratic birational
maps with real coefficients and sets analogous to
\textsc{Julia} sets for polynomials of one variable. We study
birational maps of degree $3$ and, by considering the
different possible configurations of the exceptional curves, we
give the "classification" of these maps. We can deduce
from this that the set of the birational maps of degree
$3$ exactly is irreducible, in fact rationally connected.
\end{altabstract}

\subjclass{14E07, 14E05, 37F10, 37F50}

\keywords{
\textsc{Cremona} group, rational map, birational
map, algebraic foliation, birational flow, dynamical
degree, algebraic stability, orbit, \textsc{Julia} set.}

\maketitle
\selectlanguage{french}
\vspace*{6cm}

\begin{flushright}
{\sl \`A la m\'emoire de Felice \textsc{Ronga} qui aurait sans doute aim\'e ce texte,\\

D. C.}
\end{flushright}

\vspace*{3cm}

\begin{flushright}
$\begin{array}{ll}
\text{\sl \`A ceux dont, par pudeur, je parle peu,}\\
\text{\sl avec une pens\'ee particuli\`ere pour celui qui s'est \'eclips\'e,}
\end{array}$

{\sl J. D.}
\end{flushright}

\newpage\hspace{1mm}\newpage

\setlength{\baselineskip}{0.6cm}

\begin{small}

\tableofcontents

\end{small}

\mainmatter

\chapter*{Introduction}

L'\'etude syst\'ematique du groupe des transformations
birationnelles du plan projectif complexe dans lui-m\^eme, encore
appel\'e groupe de \textsc{Cremona},
prend son essor dans les ann\'ees $1860.$ Parmi les protagonistes
on trouve \textsc{Cremona}, \textsc{N\oe ther}, \textsc{de
Jonqui\`eres} etc. Leur approche repose sur une correspondance
entre les transformations birationnelles et certains syst\`emes de
courbes rationnelles que nous allons mentionner. Une
\textbf{\textit{transformation rationnelle}}\label{ind1} $f$ de
$\mathbb{P}^2(\mathbb{C})$ dans lui-m\^eme est une transformation de la forme
\begin{align*}
&f\hspace{0.1cm}\colon\hspace{0.1cm}\mathbb{P}^2(\mathbb{C})
\dashrightarrow\mathbb{P}^2(\mathbb{C}), && (x_0:x_1:x_2) \mapsto
(f_0(x_0,x_1,x_2):f_1(x_0,x_1,x_2):f_2(x_0,x_1,x_2))
\end{align*}
les $f_i$ d\'esignant des polyn\^omes homog\`enes de m\^eme
degr\'e. Le degr\'e de $f$ est par d\'efinition
le degr\'e des $f_i.$ Une \textbf{\textit{transformation
birationnelle}}\label{ind3} de $\mathbb{P}^2(\mathbb{C})$ dans
lui-m\^eme ou \textbf{\textit{transformation de \textsc{Cremona}}}
est une transformation rationnelle admettant un inverse rationnel.
Nous noterons $\mathrm{Bir}(\mathbb{P}^2(\mathbb{C}))$ le groupe
des transformations birationnelles de $\mathbb{P}^2(\mathbb{C})$
dans lui-m\^eme. Si
\begin{align*}
&f\hspace{0.1cm}\colon\hspace{0.1cm}\mathbb{P}^2(\mathbb{C})\dashrightarrow\mathbb{P}^2(\mathbb{C})
\end{align*}
est une transformation birationnelle, le \textbf{\textit{r\'eseau
homalo\"idal}}\label{ind4} associ\'e \`a $f$ est
le syst\`eme de cour\-bes~$\mathscr{H}_f$ d\'efini par
\begin{align*}
&\alpha_0f_0+\alpha_1f_1+\alpha_2f_2=0,&& (\alpha_0:\alpha_1:
\alpha_2)\in\mathbb{P}^2(\mathbb{C});
\end{align*}
c'est l'image r\'eciproque par $f$ du r\'eseau de droites
$\alpha_0x_0+\alpha_1x_1+\alpha_2x_2=0.$ En particulier chaque
courbe du r\'eseau $\mathscr{H}_f$ est rationnelle. Les points
base de $\mathscr{H}_f$ sont les points par lesquels
passent toutes les courbes du r\'eseau; on les appelle aussi
\textbf{\textit{points base}}\label{ind5} de $f.$ Ils peuvent \^etre dans
$\mathbb{P}^2(\mathbb{C})$ ou infiniment proches\footnote{\hspace{0.1cm}
Soient $S$ une surface et $m$ un point de~$S.$ Le diviseur
exceptionnel $E$ obtenu en \'eclatant $m$ est appel\'e
\textbf{\textit{premier voisinage infinit\'esimal}}\label{ind6} de
$m$ et les points de $E$ sont dits \textbf{\textit{infiniment
proches}}\label{ind7} de $m.$ Le \textbf{\textit{$k$-i\`eme
voisinage infinit\'esimal}}\label{ind8} de $m$ est l'ensemble des
points contenus dans le premier voisinage d'un certain point du
$(k-1)$-i\`eme voisinage infinit\'esimal de $m.$ Par opposition
aux points infiniment proches, les points de $S$ sont appel\'es
\textbf{\textit{points propres}}\label{ind9}.} de $\mathbb{P}^2(\mathbb{C});$
d\`es que l'un de ces points n'est pas propre il appartient \`a un
$k$-i\`eme voisinage infinit\'esimal d'un point base. Les points
base propres de~$f$ sont les points d'ind\'etermination, {\it i.e.} ce sont les z\'eros
communs des $f_i;$ ces points sont les points \'eclat\'es par $f.$
La multiplicit\'e au point base $m_i$ de $f$ est la multiplicit\'e
d'une courbe g\'en\'erique de $\mathscr{H}_f$ en $m_i,$ {\it i.e.}
l'ordre en $m_i$ d'un \'el\'ement g\'en\'erique de~$\mathscr{H}_f.$

On a vu qu'\`a toute transformation birationnelle on
peut associer un r\'eseau homalo\"idal, r\'eciproquement on a le

\begin{theoreme}[\cite{RS}]
{\sl Un r\'eseau homalo\"idal d\'efinit une infinit\'e de
transformations birationnelles, chacune pouvant \^etre obtenue \`a
partir d'une autre via composition \`a gauche par un automorphisme
de $\mathbb{P}^2(\mathbb{C}).$}
\end{theoreme}

Ainsi le point de vue homalo\"idal conduit \`a ce que nous appellerons
la conjugaison gauche-droite.

\begin{exemples}
\textbf{\textit{i.}} Consid\'erons la transformation
birationnnelle, appel\'ee \textbf{\textit{involution de
\textsc{Cremona}}}\label{ind11}, d\'efinie par
\begin{align*}
&\sigma\hspace{0.1cm}\colon\hspace{0.1cm}\mathbb{P}^2(\mathbb{C})
\dashrightarrow\mathbb{P}^2(\mathbb{C}), &&
(x_0:x_1:x_2)\mapsto(x_1x_2:x_0x_2:x_0x_1).
\end{align*}\label{not1}

Les points d'ind\'etermination de $\sigma$ sont $(1:0:0),$
$(0:1:0)$ et $(0:0:1).$ Les transformations birationnelles
dont le r\'eseau homalo\"{\i}dal est constitu\'e des coniques
passant par les points 
\begin{align*}
&(1:0:0),&& (0:1:0)&& \text{et} &&(0:0:1) 
\end{align*}
sont les transformations de la forme $A\sigma$ avec $A$ dans
$\mathrm{Aut}(\mathbb{P}^2(\mathbb{C})).$

\bigskip

\textbf{\textit{ii.}} Soit $\mathscr{S}_1$ le r\'eseau
homalo\"{\i}dal form\'e des coniques passant par $(1:0:0),$
$(0:1:0)$ et tangentes \`a la droite d'\'equation $x_2=0.$ Les
transformations birationnelles ayant $\mathscr{S}_1$ pour r\'eseau
homalo\"{\i}dal sont du type $A\rho$ avec 
\begin{align*}
& A\in\mathrm{Aut}(\mathbb{P}^2(\mathbb{C})),&& \rho\hspace{0.1cm}
\colon\hspace{0.1cm}\mathbb{P}^2(\mathbb{C})\dashrightarrow\mathbb{P}^2(
\mathbb{C}),\hspace{2mm}(x_0:x_1:x_2) \mapsto(x_0x_1:x_2^2:x_1x_2)
\end{align*}
dont les points d'ind\'etermination sont $(1:0:0)$ et $(0:1:0).$

 \bigskip

\textbf{\textit{iii.}} Enfin soit $\mathscr{S}_2$ le
r\'eseau homalo\"{\i}dal constitu\'e des coniques passant par
$(0:0:1),$ tangentes et osculatrices \`a $x_0=0$ en ce point;
les transformations birationnelles ayant $\mathscr{S}_2$ pour
r\'eseau homalo\"{\i}dal s'\'ecrivent $A\tau$ o\`u $A$ d\'esigne
un automorphisme de $\mathbb{P}^2(\mathbb{C})$ et
\begin{align*}
&\tau\hspace{0.1cm}\colon \hspace{0.1cm}\mathbb{P}^2(\mathbb{C})
\dashrightarrow\mathbb{P}^2(\mathbb{C}),&&(x_0:x_1:x_2)
\mapsto(x_0^2:x_0x_1:x_1^2-x_0x_2).
\end{align*}
\end{exemples}

Comme nous le verrons les involutions $\sigma,$
$\rho$ et $\tau$ jouent un r\^ole essentiel dans le groupe de
\textsc{Cremona}.

Soient $f$ une transformation birationnelle et
$\mathscr{H}_f$ le r\'eseau homalo\"idal associ\'e \`a $f.$ Les
courbes de $\mathscr{H}_f$ satisfont les \'equations suivantes
(\cite{AC})
\begin{equation}\label{anc}
\sum_{i=1}^q\mu_i^2=n^2-1,\hspace{18mm}\sum_{i=1}^q\mu_i
=3n-3
\end{equation}
o\`u $\mu_i$ d\'esigne la multiplicit\'e aux points base (qui
rappelons le ne sont pas n\'ecessairement propres), $q$ le nombre
de points base et $n$ le degr\'e de $f.$ La premi\`ere \'equation
traduit le fait que deux courbes g\'en\'eriques du r\'eseau
homalo\"{\i}dal se coupent en les points base et un unique autre
point; la seconde exprime que les courbes du r\'eseau sont
rationnelles.

\`A partir, entre autres, de ces \'equations
\textsc{N\oe ther} \'enonce en $1871$ le

\begin{theoreme}[Th\'eor\`eme de N\oe ther]\label{nono}
{\sl Le groupe de \textsc{Cremona} est engendr\'e par le 
groupe des automorphismes de $\mathbb{P}^2(\mathbb{C}),$ {\it i.e.}
$\mathrm{PGL}_3(\mathbb{C}),$ et l'involution
\begin{align*}
&\sigma\hspace{0.1cm}\colon\hspace{0.1cm}\mathbb{P}^2(\mathbb{C})\dashrightarrow
\mathbb{P}^2(\mathbb{C}),&& (x_0:x_1:x_2)\mapsto
(x_1x_2:x_0x_2:x_0x_1).
\end{align*}}
\end{theoreme}

En particulier toute transformation birationnelle du
plan projectif complexe qui n'est pas un automorphisme
est compos\'ee d'un certain nombre de
transformations quadratiques, {\it i.e.} de degr\'e $2.$

La d\'emonstration donn\'ee par \textsc{N\oe ther}
comportait des lacunes (\cite{No, No2, No3}); il faudra attendre
$30$ ans et \textsc{Castelnuovo} pour en avoir une compl\`ete
(\emph{voir} \cite{Cas}). Depuis il y a eu de nombreuses preuves
de ce r\'esultat; pour une chronologie d\'etaill\'ee on peut
consulter~\cite{AC}.

Soit $\varphi\hspace{1mm}\colon\hspace{1mm}\mathrm{PGL}_3(
\mathbb{C})\ast\mathbb{Z}/2\mathbb{Z}\to\mathrm{Bir}(\mathbb{P}^2
(\mathbb{C}))$
le morphisme d\'efini par
\begin{align*}
& \varphi(A)=A \text{ sur } \mathrm{PGL}_3(\mathbb{C}) && \text{et}&& \varphi(1)=\sigma
\text{ sur } \mathbb{Z}/2\mathbb{Z}.
\end{align*}
La description du noyau de $\varphi,$
{\it i.e.} la description des relations dans le groupe
de \textsc{Cremona}, n'aura pas lieu avant les ann\'ees $1980$
(\emph{voir} \cite{Gi}); ce noyau est un groupe libre: c'est une
cons\'equence directe de la simplicit\'e des facteurs
$\mathrm{Aut}(\mathbb{P}^2(\mathbb{C}))$ et
$\mathbb{Z}/2\mathbb{Z}$ et d'un th\'eor\`eme de \textsc{Kuro{\v{s}}} (\cite{Se}).

\begin{remarque}
Il n'y a pas d'analogue au Th\'eor\`eme de \textsc{N\oe ther} en
dimension sup\'erieure comme en t\'emoigne l'\'enonc\'e suivant
d\^u \`a \textsc{Pan}.
\begin{theoreme}[\cite{Pa4}]
{\sl Soit $n$ un entier sup\'erieur ou \'egal \`a $3.$ Tout ensemble de
g\'en\'erateurs du groupe des transformations birationnelles de
$\mathbb{P}^n(\mathbb{C})$ dans lui-m\^eme doit contenir un nombre
infini et non d\'enombrable de transformations de degr\'e
strictement sup\'erieur \`a $1.$}
\end{theoreme}

Notons qu'\textsc{Hudson} avait d\'emontr\'e un
r\'esultat analogue pour $n=3$ via des techniques diff\'erentes
(\cite{Co}, Book II, Chapter VII, \S $3$).
\end{remarque}

\bigskip

Le \textbf{\textit{groupe de \textsc{de
Jonqui\`eres}}}\label{ind12}, not\'e $\mathrm{dJ},\label{not2}$
est le groupe d'invariance de la fibration standard $x_1=$ cte; il
est isomorphe au produit semi-direct
$\mathrm{PGL}_2(\mathbb{C}(x_1))\rtimes\mathrm{PGL}_2(\mathbb{C}).$
Une \textbf{\textit{transformation de \textsc{de
Jonqui\`eres}}}\label{ind13} est un \'el\'ement de $\mathrm{dJ}.$ Une
telle transformation s'\'ecrit en carte affine
\begin{align*}
& \left(\frac{a(x_1)x_0+b(x_1)}{c(x_1)x_0+d(x_1)},\frac{\alpha x_1+\beta}
{\gamma x_1+\delta}\right), && \left[\begin{array}{cc}
\alpha & \beta\\
\gamma & \delta
\end{array}
\right]\in\mathrm{PGL}_2(
\mathbb{C}), && \left[
\begin{array}{cc}
a(x_1) & b(x_1)\\
c(x_1) & d(x_1)
\end{array}
\right]\in\mathrm{PGL}_2(\mathbb{C}(x_1)).
\end{align*}

\bigskip

Dans les ann\'ees $1990$ \textsc{Iskovskikh} donne une pr\'esentation du groupe
de \textsc{Cremona}.

\begin{theoreme}[\cite{I2}]
{\sl Le groupe des transformations birationnelles de
$\mathbb{P}^1(\mathbb{C})\times \mathbb{P}^1(\mathbb{C})$ est engendr\'e par le
groupe des automorphismes de
$\mathbb{P}^1(\mathbb{C})\times\mathbb{P}^1(\mathbb{C})$ et le groupe de
\textsc{de Jonqui\`eres}.

De plus les relations dans
$\mathrm{Bir}(\mathbb{P}^1(\mathbb{C})\times\mathbb{P}^1(\mathbb{C}))$ sont les
relations internes au groupe de \textsc{de Jonqui\`eres}, au
groupe $\mathrm{Aut}(\mathbb{P}^1(\mathbb{C}) \times\mathbb{P}^1(\mathbb{C}))$
auxquelles s'ajoute la relation
\begin{align*}
&(\tau e)^3=\left(\frac{1}{x_0},\frac{1}{x_1}\right)&& \text{o\`u} &&
\tau\hspace{0.1cm}\colon\hspace{0.1cm} (x_0,x_1)\mapsto(x_1,x_0) && \text{et} &&
e\hspace{0.1cm}\colon\hspace{0.1cm}(x_0,x_1)\mapsto\left(x_0,\frac{x_0} {x_1}\right).
\end{align*}}
\end{theoreme}

Alors que \textsc{N\oe ther} s'int\'eresse \`a la
g\'en\'eration de $\mathrm{Bir}(\mathbb{P}^2(\mathbb{C})),$
\textsc{Bertini} se consacre \`a la classification des involutions
birationnelles \`a conjugaison birationnelle pr\`es (\cite{Be});
quelques ann\'ees plus tard \textsc{Kantor} entreprend de classer
les sous-groupes finis de $\mathrm{Bir}(\mathbb{P}^2(\mathbb{C}))$ \`a
conjugaison birationnelle pr\`es (\cite{K}). Malgr\'e les
corrections apport\'ees par \textsc{Wiman} (\emph{voir} \cite{Wi}) cette
classification resta longtemps incompl\`ete et redondante. En
$2000$ \textsc{Bayle} et \textsc{Beauville} pr\'esentent, dans
\cite{BaBe}, une nouvelle d\'emonstration de la classification des
involutions birationnelles; d\`es lors de nouveaux r\'esultats
apparaissent dans le m\^eme esprit (citons par exemple \cite{dF, B}). Pour
une chronologie d\'etaill\'ee de l'\'etude des sous-groupes finis
de $\mathrm{Bir}(\mathbb{P}^2(\mathbb{C}))$ on peut consulter~\cite{DI};
dans cet article \textsc{Dolgachev} et \textsc{Iskovskikh} donnent la
classification des sous-groupes finis du
groupe de \textsc{Cremona}. Leur
d\'emonstration repose sur le point suivant (\cite{M, I}): si $G$
d\'esigne un sous-groupe fini de $\mathrm{Bir}(\mathbb{P}^2(\mathbb{C}))$
alors $G$ est birationnellement conjugu\'e \`a un sous-groupe de
$\mathrm{Aut}(S)$ o\`u $S$ est une surface rationnelle
satisfaisant une des propri\'et\'es suivantes
\begin{itemize}
\item il existe un fibr\'e en coniques\footnote{\hspace{1mm} Les fibres
g\'en\'eriques d'un fibr\'e en coniques sont de genre $0,$ les fibres singuli\`eres sont
union de deux courbes rationnelles, {\it i.e.} les fibres sont de
m\^eme nature que celles d'un vrai pinceau de coniques dans
$\mathbb{P}^2(\mathbb{C}).$} $\pi\hspace{0.1cm}\colon\hspace{0.1cm} S\to\mathbb{P}^1(\mathbb{C})$
invariant par $G;$

\item $S$ est de \textsc{del Pezzo}\footnote{\hspace{1mm} La surface $S$
est une \textbf{\textit{surface de \textsc{del
Pezzo}}}\label{ind14} si $S$ est isomorphe \`a l'une des surfaces
suivantes: $\mathbb{P}^2(
\mathbb{C}),$ la quadrique $\mathbb{P}^1(\mathbb{C})\times\mathbb{P}^1
(\mathbb{C}),$ le plan projectif dont on a \'eclat\'e un point, 
$\mathbb{P}^2(\mathbb{C})$ \'eclat\'e en $n\geq 2$ points
$p_i,$ ces points satisfaisant les conditions suivantes (\cite{Dem})
\begin{itemize}
\item il n'y en a pas trois align\'es;

\item il n'y en a pas six sur une conique;

\item si de plus $n=8$ les $p_i$ ne sont pas sur une cubique plane dont le
point singulier serait un des $p_i.$
\end{itemize}}.
\end{itemize}

\bigskip

Simultan\'ement \textsc{Castelnuovo} \'etudie les
transformations birationnelles qui fixent point par point une
courbe de genre $1.$

\begin{theoreme}[\cite{Cas2}]
{\sl Soit $f$ une transformation birationnelle de $\mathbb{P}^2(\mathbb{C})$
non triviale fixant point par point une courbe irr\'eductible de
genre strictement sup\'erieur \`a $1.$ Alors $f$ est conjugu\'ee
\`a une transformation de \textsc{de Jonqui\`eres} ou bien $f$ est
d'ordre $2,$ $3$ ou $4.$}
\end{theoreme}

En $2006$ \textsc{Blanc}, \textsc{Pan} et \textsc{Vust}
donnent une version un peu plus pr\'ecise de ce r\'esultat en
reprenant la m\^eme d\'emarche que \textsc{Castelnuovo}.

\begin{theoreme}[\cite{BPV}]
{\sl Soit $f\hspace{0.1cm}\colon\hspace{0.1cm}\mathbb{P}^2(\mathbb{C})\dashrightarrow\mathbb{P}^2(\mathbb{C})$ une
transformation de \textsc{Cremona} non triviale qui fixe point par
point une courbe irr\'eductible de genre strictement sup\'erieur
\`a $1.$ Alors~$f$ est conjugu\'ee \`a une transformation de
\textsc{de Jonqui\`eres} ou bien $f$ est d'ordre $2$ ou $3.$ De
plus dans le premier cas si $f$ est d'ordre fini, c'est une
involution.}
\end{theoreme}

Ils d\'ecrivent ensuite les sous-groupes de
$\mathrm{Bir}(\mathbb{P}^2 (\mathbb{C}))$ dont tous les \'el\'ements
fixent point par point une courbe irr\'eductible de genre $g>1.$
Avant d'\'enoncer leur r\'esultat rappelons quelques
d\'efinitions.
Comme nous l'avons dit le groupe de \textsc{de
Jonqui\`eres} est le groupe d'invariance de la fibration standard
$x_1=$ cte. Soit $P$ un \'el\'ement de
$\mathbb{C}[x_1]$ qui n'est pas un carr\'e; on d\'esigne par
$\mathrm{dJ}_P$ le sous-groupe de $\mathrm{dJ}$ d\'efini par
\begin{align*}
&\mathrm{dJ}_P=\left\{
\left(\frac{a_1(x_1)x_0+a_2(x_1)P(x_1)}{a_2(x_1)x_0+a_1(x_1)},
x_1\right)\hspace{0.1cm}|\hspace{0.1cm} a_i\in\mathbb{C}(x_1),\hspace{0.1cm} a_1^2-Pa_2^2\not=0\right\}.
\end{align*}
Remarquons que si $f$ est un \'el\'ement de
$\mathrm{dJ}_P$ la courbe hyperelliptique $x_0^2=P(x_1)$
est fix\'ee point par point par $f.$ On peut ainsi r\'ealiser
n'importe quel genre pour les courbes de points fixes
d'\'el\'ements de $\mathrm{Bir}(\mathbb{P}^2(\mathbb{C})).$

Soient $m_1,$ $\ldots,$ $m_7$ sept points de
$\mathbb{P}^2(\mathbb{C})$ en position g\'en\'erale et $L$ le syst\`eme
lin\'eaire form\'e des cubiques passant par les $m_i;$ il est de
dimension $2.$ Consid\'erons le pinceau de courbes de~$L$ passant
par un point g\'en\'erique $m$ de $\mathbb{P}^2(\mathbb{C}).$ Les $m_j$ et
$m$ en sont des points base et il passe par un
neuvi\`eme point base que nous noterons $\gamma(m);$ la
transformation $\gamma$ est une involution birationnelle de
degr\'e $8$ appel\'ee \textbf{\textit{involution de
\textsc{Geiser}}}\label{ind15}. Les involutions de \textsc{Geiser}
fixent point par point une courbe non hyperelliptique de genre $3$
(\emph{voir} \cite{Dem}).

Soient $m_1,$ $\ldots,$ $m_8$ huit points de
$\mathbb{P}^2(\mathbb{C})$ en position g\'en\'erale; consid\'erons le
pinceau de cubiques $\mathscr{C}$ passant par ces huit points. Il
a un neuvi\`eme point base que nous noterons $m_9.$ Soit~$m$ un
point g\'en\'erique de $\mathbb{P}^2(\mathbb{C});$ le pinceau
$\mathscr{C}$ contient une unique cubique $\mathscr{C}(m)$ passant
par $m.$ Consid\'erons la loi de groupe sur $\mathscr{C} (m)$
ayant $m_9$ comme \'el\'ement neutre. Posons $\beta(m)=-m;$ alors
$\beta$ d\'efinit une involution birationnelle de degr\'e $17$
qu'on appelle \textbf{\textit{involution de
\textsc{Bertini}}}\label{ind16}. Les involutions de
\textsc{Bertini} fixent point par point une courbe non
hyperelliptique de genre $4$ \`a mod\`ele lisse sur un c\^one
quadratique.

Pour finir consid\'erons la surface de \textsc{Del
Pezzo} $S$ de degr\'e $1$ donn\'ee dans l'espace projectif \`a
poids $\mathbb{P}(3,1,1,2)$ par
\begin{align*}
& w^2=x_2^3+P(x_0,x_1), && P \text{ polyn\^ome homog\`ene de degr\'e }6.
\end{align*}
La restriction de la transformation
\begin{align*}
&(w:x_0:x_1:x_2)\mapsto(w:x_0:x_1:\mathrm{j}x_2),&& \mathrm{j}\not=1, && \mathrm{j}^3=1,
\end{align*}
d\'efinit un automorphisme $\varsigma$ d'ordre $3$ de $S$ dont
l'ensemble des points fixes contient une courbe irr\'eductible de
genre $2$ (\emph{voir} \cite{BPV}).

Ces transformations produisent des exemples de groupes
finis fixant une certaine courbe.

\begin{theoreme}[\cite{BPV}]
{\sl Soit $\mathrm{G}$ un sous-groupe non trivial de
$\mathrm{Bir}(\mathbb{P}^2(\mathbb{C})).$ Supposons qu'il existe une
courbe irr\'eductible $\mathcal{C}$ de genre strictement
sup\'erieur \`a $1$ invariante point par point par chaque
\'el\'ement de $\mathrm{G}.$ Alors
\begin{itemize}
\item ou bien $\mathrm{G}$ est cyclique d'ordre $2$ ou $3,$ engendr\'e,
\`a conjugaison pr\`es, par une involution de \textsc{de
Jonqui\`eres}, une involution de \textsc{Bertini}, une involution
de \textsc{Geiser} ou encore par l'automorphisme $\varsigma$
d'ordre~$3;$

\item ou bien $G$ est conjugu\'e \`a un sous-groupe d'un certain
$\mathrm{dJ}_P.$
\end{itemize}

En particulier lorsque $\mathrm{G}$ est infini, $\mathcal{C}$
est hyperelliptique et $\mathrm{G}$ est un groupe ab\'elien conjugu\'e \`a
un sous-groupe de $\mathrm{dJ}.$}
\end{theoreme}

Soit $\mathcal{C}$ une courbe irr\'eductible de $\mathbb{P}^2(
\mathbb{C}).$ On appelle \textbf{\textit{groupe de
d\'ecomposition}}\label{ind17} de $\mathcal{C}$ le groupe des
transformations de \textsc{Cremona} qui fixent $\mathcal{C};$ le
\textbf{\textit{groupe d'inertie}}\label{ind18} de $\mathcal{C}$
est le groupe des transformations birationnelles qui
fixent $\mathcal{C}$ point par point.

Dans \cite{Pa} \textsc{Pan} d\'emontre le

\begin{theoreme}[\cite{Pa}]
{\sl Soient $\mathcal{C}$ une courbe irr\'eductible lisse non
rationnelle contenue dans le plan projectif complexe. Si $f$ appartient au
groupe de d\'ecomposition de $\mathcal{C}$ et n'appartient pas
\`a $\mathrm{Aut}(\mathbb{P}^2(\mathbb{C})),$ alors $\mathcal{C}$ est de degr\'e $3$ et les points
d'ind\'etermination de $f$ appartiennent \`a~$\mathcal{C}.$}
\end{theoreme}

Soient $\mathcal{C}$ une cubique lisse et $m_1,$ $m_2,$
$m_3$ trois points de $\mathcal{C}$ non align\'es. Si $f$ est une
transformation quadratique telle que $\mathrm{Ind}\hspace{0.1cm} f=\{m_1,\hspace{0.1cm}
m_2,\hspace{0.1cm} m_3\}$ on constate que $f(\mathcal{C})$ est une cubique
lisse isomorphe \`a $\mathcal{C}.$ Soit $A$ un automorphisme de
$\mathbb{P}^2(\mathbb{C})$ qui envoie $f(\mathcal{C})$ sur $\mathcal{C}$
alors $Af$ appartient au groupe de d\'ecomposition de
$\mathcal{C};$ une telle transformation est appel\'ee
\textbf{\textit{transformation quadratique
$\mathcal{C}$-g\'en\'erique}}\label{ind19}. Remarquons que les
$m_i$ \'etant choisis arbitrairement le groupe de d\'ecomposition
de $\mathcal{C}$ est infini non d\'enombrable.

\begin{theoreme}[\cite{Pa}]
{\sl Si $\mathcal{C}$ est une cubique lisse, le groupe de
d\'ecomposition de $\mathcal{C} $ est engendr\'e par les
transformations quadratiques $\mathcal{C}$-g\'en\'eriques.}
\end{theoreme}

\textsc{Blanc} d\'ecrit dans \cite{Bl} les \'el\'ements
d'ordre fini du groupe d'inertie d'une courbe de genre~$1;$ il
\'etudie aussi le groupe d'inertie d'une cubique lisse.

\begin{theoreme}[\cite{Bl}]
{\sl Soit $\mathcal{C}$ une cubique lisse. Les \'el\'ements non
triviaux de $\mathrm{Ine}(\mathcal{C})$ sont de degr\'e
sup\'erieur ou \'egal \`a $3.$ Le groupe d'inertie de $\mathcal{C}$
est engendr\'e par ses \'el\'ements de degr\'e~$3.$}
\end{theoreme}

\bigskip

Certaines \'etudes du groupe de \textsc{Cremona} s'inscrivent
dans une th\'ematique classique. Ainsi
il y a eu de nombreux travaux sur les propri\'et\'es
alg\'ebriques des groupes de diff\'eomorphismes des vari\'et\'es;
par exemple \textsc{Dieudonn\'e} a d\'ecrit les automorphismes de
$\mathrm{PGL}_{n+1}(\mathbb{C})=\mathrm{Aut}(\mathbb{P}^n(\mathbb{C}))$
(\emph{voir} \cite{Die}): ils s'obtiennent \`a partir des
automorphismes int\'erieurs, de la contragr\'ediente
($u\mapsto\transp u^{-1}$) et des automorphismes du corps $\mathbb{C}.$
Dans \cite{De} \textsc{D\'eserti} g\'en\'eralise cet \'enonc\'e
\`a $\mathrm{Bir}(\mathbb{P}^2(\mathbb{C})).$

\begin{theoreme}[\cite{De}]\label{autom}
{\sl Soit $\varphi$ un automorphisme du groupe de \textsc{Cremona}; il
existe $\psi$ une transformation birationnelle et $\tau$ un
automorphisme du corps $\mathbb{C}$ tels que
\begin{align*}
&\forall\hspace{0.1cm} f\in\mathrm{Bir}(\mathbb{P}^2(\mathbb{C})),
&&\varphi(f)=\tau(\psi f\psi^{-1}).
\end{align*}
Dit autrement le groupe des automorphismes ext\'erieurs de
$\mathrm{Bir}(\mathbb{P}^2(\mathbb{C}))$ s'identifie au groupe des
automorphismes du corps $\mathbb{C}.$}
\end{theoreme}

La preuve repose sur la description des sous-groupes
ab\'eliens maximaux non d\'enombrables~$\mathrm{G}$ de
$\mathrm{Bir}(\mathbb{P}^2(\mathbb{C}));$ en utilisant entre autres les
travaux de \textsc{Cantat} et \textsc{Favre} sur les sym\'etries
birationnelles des surfaces feuillet\'ees (\cite{CaFa}) on peut
montrer qu'ils satisfont l'une des propri\'et\'es suivantes
\begin{itemize}
\item $\mathrm{G}$ poss\`ede des \'el\'ements de torsion;

\item $\mathrm{G}$ est conjugu\'e \`a un sous-groupe du groupe de
\textsc{de Jonqui\`eres}.
\end{itemize}

\bigskip

Continuons avec une autre propri\'et\'e alg\'ebrique
satisfaite par $\mathrm{Bir}(\mathbb{P}^2(\mathbb{C})).$ Dans les ann\'ees
$1970$ \textsc{Tits} d\'emontre le

\begin{theoreme}[\cite{Ti}]
{\sl Soient $\Bbbk$ un corps de caract\'eristique nulle et
$\Gamma$ un sous-groupe de type fini de $\mathrm{GL}_n(\Bbbk).$
Alors
\begin{itemize}
\item ou bien $\Gamma$ contient un groupe libre non ab\'elien;

\item ou bien $\Gamma$ contient un sous-groupe r\'esoluble d'indice fini.
\end{itemize}}
\end{theoreme}

On dit que $\mathrm{GL}_n(\Bbbk)$ satisfait
l'alternative de \textsc{Tits}. R\'ecemment \textsc{Cantat} a
\'etudi\'e les sous-groupes de type fini du groupe de \textsc{Cremona}
et montr\'e qu'il en est de m\^eme pour le groupe de \textsc{Cremona}
(\cite{Ca, Fa2}). Un r\'esultat analogue avait \'et\'e
d\'emontr\'e par \textsc{Lamy} pour le sous-groupe de
$\mathrm{Bir}(\mathbb{P}^2(\mathbb{C}))$ des
\og automorphismes polynomiaux de $\mathbb{C}^2$\fg\hspace{1mm}
(\emph{voir} \cite{La}).

\bigskip

Mentionnons quelques notions de dynamique.
Soit $f$ la transformation d\'efinie par
\begin{align*}
& f\hspace{0.1cm}\colon\hspace{0.1cm}\mathbb{P}^2(\mathbb{C})\dashrightarrow\mathbb{P}^2(\mathbb{C}), &&(x_0:x_1:x_2)\mapsto(
(x_0+x_1)x_2:x_0x_1:x_0x_2).
\end{align*}
En composant $f$ avec elle-m\^eme, nous obtenons
\begin{align*}
f^2\hspace{1mm}\colon\hspace{1mm}(x_0:x_1:x_2)\mapsto(bx_0x_1+x_0x_2+x_1x_2:x_1
(x_0+x_1):x_2(x_0+x_1)).
\end{align*}
Nous avons donc $\deg(f^2)=\deg(f)^2 -2.$ Lors du calcul de $f^2$
nous avons pu \og factoriser\fg\hspace{1mm} par $x_0x_2$~; ceci traduit,
entre autres, le ph\'enom\`ene g\'eom\'etrique suivant~: la droite
d'\'equation $x_2=0$ est contract\'ee sur $(0:1:0)$ qui est un
point d'ind\'etermination de $f$
\begin{figure}[H]
\begin{center}
\input{intro.pstex_t}
\end{center}
\end{figure}
Ceci a conduit \`a introduire
la notion suivante (\cite{FoSi}): $f$ est
alg\'ebriquement stable si l'une
des propri\'et\'es \'equivalentes suivantes est v\'erifi\'ee

\vspace{1mm}

\begin{itemize}
\item $(\deg f^n)=(\deg f)^n$\hspace{6mm} $\forall\hspace{0.1cm}n\geq 0$~;

\item il n'existe pas de courbe $\mathcal{C}$ de $S$ telle que
$f^n(\mathcal{C})$ soit un point d'ind\'etermination de $f$ pour
un certain $n.$
\end{itemize}

Il y a d'autres fa\c{c}ons de caract\'eriser la stabilit\'e alg\'ebrique 
(\emph{voir} par exemple \cite{DJS}).

\bigskip

\textsc{Diller} et \textsc{Favre} montrent que pour
toute transformation birationnelle $f,$ il existe une modification
propre $\varepsilon\hspace{0.1cm}\colon\hspace{0.1cm} S\to\mathbb{P}^2(\mathbb{C})$ telle que
$\varepsilon^{-1}f\varepsilon$ soit alg\'ebriquement stable
(\cite{DiFa}, th\'eor\`eme~$0.1$).

Introduisons un invariant birationnel. Soient $f,$ $g$
deux transformations de \textsc{Cremona}; en g\'en\'eral on a $\deg
f\not=\deg(gfg^{-1}).$ Par contre il existe deux constantes
positives $\alpha,$ $\beta$ telles que
\begin{align*}
&\forall\hspace{0.1cm} n\in\mathbb{N},&& \alpha\hspace{0.1cm}\deg
f^n\leq\deg(gf^ng^{-1})\leq\beta\hspace{0.1cm} \deg f^n;
\end{align*}
la croissance des degr\'es de $f^n$ est donc un invariant
birationnel. D'o\`u la d\'efinition suivante: soit $f$ une
transformation de \textsc{Cremona}; on appelle
premier degr\'e dynamique de $f$ la
quantit\'e
\begin{align*}
\lambda(f)=\liminf_{n\to+\infty}(\deg f^n)^{1/n}.
\end{align*}
\label{not3}

On dit
qu'une transformation de \textsc{Cremona} $f$ est un automorphisme sur la surface $S$ s'il existe une transformation birationnelle 
$\phi\hspace{1mm}\colon\hspace{1mm}\mathbb{P}^2(\mathbb{C})\dashrightarrow S$
telle que $\phi f\phi^{-1}$ soit un automorphisme de $S;$ on notera
$\mathrm{Aut}(S)$ le groupe des automorphismes de $S.$

\textsc{Diller} et \textsc{Favre} ont donn\'e une
classification des transformations birationnelles \`a conjugaison
pr\`es.

\begin{theoreme}[\cite{DiFa}, th\'eor\`eme 0.2]\label{DiFa}
{\sl Soit $f$ une transformation de \textsc{Cremona}. \`A conjugaison
birationnelle pr\`es elle satisfait une et une seule des
propri\'et\'es suivantes
\begin{itemize}
\item la croissance des degr\'es est born\'ee ({\it i.e.} $\deg f^n=o(1)$), $f$ est un
automorphisme sur une certaine surface rationnelle $S$ et un
it\'er\'e de $f$ appartient \`a la composante neutre de $\mathrm{Aut}(S);$

\item la croissance des degr\'es est lin\'eaire ({\it i.e.} $\deg f^n\sim n$), auquel cas $f$ pr\'eserve
une unique fibration rationnelle et n'est pas birationnellement
conjugu\'e \`a un automorphisme;

\item la croissance des degr\'es est quadratique ({\it i.e.} $\deg f^n\sim n^2$) alors $f$ est un
automorphisme qui pr\'eserve une unique fibration elliptique;

\item la croissance des degr\'es est exponentielle ({\it i.e.} $\deg f^n\sim\alpha^n,$
$\alpha>1$).
\end{itemize}}
\end{theoreme}

Dans les trois premi\`eres \'eventualit\'es $\lambda(f)=
1,$ dans la derni\`ere $\lambda(f)>1.$

Dans ce m\^eme article on trouve une estimation du
nombre de points p\'eriodiques d'un \'el\'ement g\'en\'erique du
groupe de \textsc{Cremona} de degr\'e $d\geq 2.$

\begin{theoreme}[\cite{DiFa}]
{\sl Soit $f$ une transformation de \textsc{Cremona}
alg\'ebriquement stable. Supposons que $f$ ne soit pas un automorphisme
de $\mathbb{P}^2(\mathbb{C})$ et n'ait pas de courbe de
points p\'eriodiques. Notons $\mathrm{Per}_k$ le nombre de points
p\'eriodiques de p\'eriode (divisant) $k.$ Il existe une constante
$C>0$ telle que pour tout $k\geq 0$ on ait
\begin{align*}
|\mathrm{Per}_k-\lambda(f)^k|\leq C.
\end{align*}}
\end{theoreme}

Nous allons voir que sous certaines hypoth\`eses il
y a abondance de points
p\'eriodiques hyperboliques (\emph{voir} Chapitre~\ref{dyndyn}).
Une transformation birationnelle $f$ de $\mathbb{P}^2(\mathbb{C})$
v\'erifie la condition de \textsc{Bedford} et
\textsc{Diller}\label{ind22}
 (\emph{voir} \cite{BD}) si
\begin{align*}
\sum_{n\geq 0}\frac{1}{\lambda(f)^n}|\log(\textrm{dist}(f^n(
\mathrm{Ind}\hspace{0.1cm} f^{-1}),\mathrm{Ind}\hspace{0.1cm} f))|<\infty
\end{align*}
o\`u $\textrm{dist}$ d\'esigne la distance de \textsc{Fubini}-\textsc{Study}.

On peut \'enoncer le

\begin{theoreme}[\cite{BD, Du}]\label{BDDu}
{\sl Soit $f$ un \'el\'ement du groupe de
\textsc{Cremona} n'appartenant pas \`a $\mathrm{Aut}(\mathbb{P}^2(
\mathbb{C})).$ Supposons que $f$ v\'erifie la condition de
\textsc{Bedford} et \textsc{Diller}. Alors $f$ poss\`ede une
infinit\'e de points p\'eriodiques hyperboliques qui
s'\'equidistribuent suivant une mesure de probabilit\'e
$f$-invariante.}
\end{theoreme}

\bigskip

Bien que de nature alg\'ebrique les \'enonc\'es
qui suivent ont \'et\'e \'etablis \`a partir de techniques de dynamique
complexe. Les repr\'esentations lin\'eaires des r\'eseaux de
groupes de \textsc{Lie} r\'eels simples de rang r\'eel strictement
sup\'erieur \`a $1$ ont \'et\'e \'etudi\'ees par \textsc{Margulis}
(\cite{Mar, VGS}); afin de g\'en\'eraliser les travaux de celui-ci
aux repr\'esentations non lin\'eaires des r\'eseaux sur les
vari\'et\'es compactes \textsc{Zimmer} \'etablit un programme
(\cite{Zi1, Zi2, Zi3, Zi4}). \textsc{D\'eserti}
s'int\'eresse \`a une version birationnelle de ce programme.

\begin{theoreme}[\cite{De3}]\label{Zimmer}
{\sl Soient $\Gamma$ un sous-groupe d'indice fini de
$\mathrm{SL}_3(\mathbb{Z})$ et $\rho$ un morphisme injectif de
$\Gamma$ dans le groupe de \textsc{Cremona}. Alors, \`a
conjugaison birationnelle pr\`es, $\rho$ co\"incide avec le
plongement canonique ou la contragr\'ediente ({\it i.e.}
l'involution $u\mapsto\transp u^{-1}$).}
\end{theoreme}

On obtient comme cons\'equence le

\begin{corollaire}[\cite{De3}]
{\sl Soit $\Gamma$ un sous-groupe d'indice fini de $\mathrm{SL}_n(\mathbb{Z}).$
D\`es que $n\geq 4$ le groupe $\Gamma$ ne se plonge pas dans
$\mathrm{Bir}(\mathbb{P}^2(\mathbb{C})).$}
\end{corollaire}

Dans le m\^eme esprit mais avec des
m\'ethodes diff\'erentes \textsc{Cantat} montre le

\begin{theoreme}[\cite{Ca}]
{\sl Soient $\Gamma$ un groupe de \textsc{Kazdhan}\footnote{\hspace{1mm}
Les groupes de \textsc{Kazdhan} sont caract\'eris\'es par
la propri\'et\'e suivante: un
groupe $\mathrm{G}$ est de \textsc{Kazdhan} (ou a la propri\'et\'e (T) de
 \textsc{Kazdhan}) si toute action de $\mathrm{G}$ par isom\'etrie sur un
espace hyperbolique r\'eel ou complexe a un point fixe global
\cite{Delorme, Guich}.} discret et $\rho$
un morphisme de~$\Gamma$ dans le groupe de \textsc{Cremona}. Ou bien
$\rho(\Gamma)$ est fini, ou bien $\rho(\Gamma)$ est conjugu\'e \`a
un sous-groupe de
$\mathrm{Aut}(\mathbb{P}^2(\mathbb{C}))=\mathrm{PGL}_3(\mathbb{C}).$}
\end{theoreme}

\bigskip

Le Th\'eor\`eme \ref{Zimmer} appliqu\'e \`a $\Gamma=
\mathrm{SL}_3(\mathbb{Z})$ permet de d\'ecrire le semi-groupe des
endomorphismes de $\mathrm{Bir}(\mathbb{P}^2(\mathbb{C}));$
en particulier

\begin{theoreme}[\cite{De4}]
{\sl Le groupe de
\textsc{Cremona} est hopfien, {\it i.e.} tout endomorphisme
surjectif de $\mathrm{Bir}(\mathbb{P}^2(\mathbb{C}))$ est
bijectif.

Plus pr\'ecis\'ement soit $\varphi$ un endomorphisme non 
trivial de $\mathrm{Bir}(\mathbb{P}^2(\mathbb{C}));$ il 
existe un plongement de corps $\tau$ de $\mathbb{C}$ 
dans lui-m\^eme et une transformation birationnelle 
$\psi$ telle que
\begin{align*}
&\varphi(f)=\tau(\psi f\psi^{-1}), &&\forall\hspace{1mm} f\in \mathrm{Bir}(\mathbb{P}^2
(\mathbb{C})).
\end{align*}}
\end{theoreme}

Notons qu'on retrouve en particulier le Th\'eor\`eme \ref{autom}.

\bigskip

Certaines des m\'ethodes utilis\'ees pour d\'emontrer le
Th\'eor\`eme \ref{Zimmer} associ\'ees \`a d'autres arguments
permettent de d\'ecrire les sous-groupes nilpotents de
$\mathrm{Bir}(\mathbb{P}^2(\mathbb{C})).$

\begin{theoreme}[\cite{De5}]
{\sl Soient $\mathrm{G}$ un sous-groupe nilpotent, non ab\'elien \`a
indice fini pr\`es et~$\rho$ un morphisme injectif de $\mathrm{G}$
dans $\mathrm{Bir}(\mathbb{P}^2(\mathbb{C})).$ Le groupe
$\mathrm{G}$ satisfait l'une des propri\'et\'es suivantes
\begin{itemize}
\item \`a indice fini pr\`es le premier groupe
d\'eriv\'e\footnote{\hspace{0.1cm} Le premier groupe d\'eriv\'e de $\mathrm{G}$ est
$\langle aba^{-1}b^{-1}\hspace{0.1cm}|\hspace{0.1cm} a,\hspace{0.1cm} b\in\mathrm{G}\rangle.$} de $\mathrm{G}$ est
ab\'elien;

\item tous les \'el\'ements de $\mathrm{G}$ sont de torsion.
\end{itemize}}
\end{theoreme}

Ces \'enonc\'es confortent l'id\'ee qu'en un certain
sens le groupe $\mathrm{Bir}(\mathbb{P}^2(\mathbb{C}))$
se comporte comme un groupe lin\'eaire (de rang $2$):
outre l'alternative de \textsc{Tits} on ne peut pas immerger dans
$\mathrm{Bir}(\mathbb{P}^2(\mathbb{C}))$ les r\'eseaux
$\mathrm{SL}_n(\mathbb{Z})$ et leurs sous-groupes
d'indice fini pour $n\geq 4$ au m\^eme titre que dans
$\mathrm{PGL}_3(\mathbb{C}).$

Ce qui pr\'ec\`ede ne pr\'etend pas bien s\^ur \^etre
une pr\'esentation exhaustive des travaux concernant le groupe de
\textsc{Cremona} mais un point de vue partiel proche des
pr\'eoccupations des auteurs.

Plusieurs m\'emoires de synth\`ese sont disponibles;
citons en particulier \cite{Coble, Coolidge, Godeaux, Hu, RS}. Par
exemple les chapitres $6,$ $7$ et $8$ du tome $4$ de
\cite{Coolidge}  sont consacr\'es aux transformations de
\textsc{Cremona}; on y trouve les notions de points \'eclat\'es,
courbes contract\'ees, le th\'eor\`eme de factorisation de
\textsc{N\oe ther}. Les transformations birationnelles d'ordre $\leq 7,$ celles ayant une
courbe de points fixes de genre $>1,$ les sous-groupes finis de
$\mathrm{Bir}(\mathbb{P}^2(\mathbb{C}))$ sont des
sujets abord\'es. L'ouvrage d'\textsc{Hudson} qui date de $1927$
contient un chapitre sur les transformations birationnelles
quadratiques; toutefois l'\'evolution du vocabulaire et des
techniques en rend sa lecture difficile, de plus certains
arguments \og g\'en\'eriques\fg\hspace{0.1cm} propres \`a l'\'epoque nous
semblent devoir parfois \^etre pr\'ecis\'es. Le livre
d'\textsc{Alberich-Carrami\~{n}ana} est le seul ouvrage r\'ecent
consacr\'e aux transformations de \textsc{Cremona}; l'auteur
reprend avec le langage et la rigueur actuels des r\'esultats
mentionn\'es et/ou d\'emontr\'es par les math\'ematiciens
pr\'ec\'edemment cit\'es et d\'emontre des r\'esultats originaux.
Dans tous ces m\'emoires ce sont essentiellement des propri\'et\'es
alg\'ebriques qui sont \'etudi\'ees, les questions de nature
\og dynamique \fg\hspace{0.1cm} sont bien plus r\'ecentes.

Venons en au contenu de ce texte dont le but
essentiel est finalement de fournir au lecteur suffisamment
d'exemples - parfois innocents - pour permettre d'aborder
des probl\`emes de nature plus difficiles voire conceptuels.
Il est motiv\'e en particulier par l'importance des transformations
quadratiques ne serait-ce qu'au travers du th\'eor\`eme de
\textsc{N\oe ther}; mais elles interviennent profond\'ement
aussi en g\'eom\'etrie ag\'ebrique classique (\og d\'esingularisation
des courbes \`a la \textsc{N\oe ther}\fg\hspace{1mm} par exemple) et du point de
vue de la dynamique nous offrent le premier degr\'e de complexit\'e
\`a \'etudier dans le m\^eme ordre d'id\'ee que les applications
de \textsc{H\'enon} que l'on peut \'evidemment voir dans
$\mathrm{Bir}(\mathbb{P}^2(\mathbb{C})).$

\vspace{12mm}

Dans le premier Chapitre on cherche \`a
comprendre dans l'espace des applications rationnelles
quadratiques de $\mathbb{P}^2(\mathbb{C})$ dans lui-m\^eme le sous-ensemble
constitu\'e des transformations birationnelles. Soyons plus
pr\'ecis; consid\'erons l'espace $\mathrm{Rat}_2$ qui est le projectivis\'e
de l'espace vectoriel des triplets de formes quadratiques en $3$
variables complexes
\begin{align*}
\mathrm{Rat}_2:=\mathbb{P}\{(f_0,f_1,f_2)\hspace{0.1cm}|\hspace{0.1cm}
f_i\in\mathbb{C}[x_0,x_1,x_2]_2\}.
\end{align*}
On constate que $\mathrm{Rat}_2$ s'identifie \`a
$\mathbb{P}^{17}(\mathbb{C}).$ \`A un \'el\'ement $f$ de
$\mathrm{Rat}_2$ on peut faire correspondre la transformation
rationnelle $f^{\bf
\bullet}\hspace{0.1cm}\colon\hspace{0.1cm}\mathbb{P}^2(\mathbb{C})
\dashrightarrow\mathbb{P}^2(\mathbb{C})$\label{not3a} d\'efinie
avec les notations pr\'ec\'edentes par
\begin{align*}
&x=(x_0:x_1:x_2)\dashrightarrow f^{\bullet}(x)=\delta
(f_0:f_1:f_2)(x), && \delta=\frac{1}{\text{pgcd}(f_0,f_1,f_2)}
\end{align*}
On note $\mathrm{Rat}_2^{\bullet}$ l'ensemble des transformations
$f^{\bullet}$ o\`u $f$ appartient \`a $\mathrm{Rat}_2;$ visiblement
$\mathrm{Rat}_2^{\bullet}$ s'identifie au quotient de $\mathrm{Rat}_2$ par la
relation d'\'equivalence $f=(f_0:f_1:f_2)\sim g=(g_0:g_1:g_2)$ si et seulement si la
matrice
\begin{align*}
\left(
\begin{array}{ccc}
f_0 & f_1 & f_2 \\
g_0 & g_1 & g_2
\end{array}
\right)
\end{align*}
est partout de rang \'egal \`a $1.$

Notons que sur l'ouvert de \textsc{Zariski}
$\mathring{\mathrm{R}}\mathrm{at}_2,$ constitu\'e des
$f=(f_0:f_1:f_2)$ tels que $\text{pgcd}(f_0,f_1,f_2) =~1,$
l'application $f\mapsto f^{\bullet}$ est injective. Remarquons
aussi que l'on peut mettre une structure alg\'ebrique sur
$\mathrm{Rat}_2^{\bullet}$ et parler en particulier d'adh\'erence
de \textsc{Zariski} d'un sous-ensemble de
$\mathrm{Rat}_2^{\bullet}.$

Dans $\mathrm{Rat}_2$ nous consid\'erons 
\begin{align*}
\mathrm{Bir}_2:=\{f\in\mathrm{Rat}_2\hspace{0.1cm}|\hspace{0.1cm}f^{\bullet} \text{ soit birationnelle }\}
\end{align*}
et nous notons $\mathrm{Bir}_2^{\bullet}\subset\mathrm{Rat}_2^{\bullet}$
le sous-ensemble correspondant.

Il y a sur $\mathrm{Rat}_2$ (resp.
$\mathrm{Rat}_2^{\bullet}$) deux actions naturelles, toutes les
deux laissant $\mathrm{Bir}_2$ (resp. $\mathrm{Bir}_2^{\bullet}$)
invariant. La premi\`ere est l'action de
$\mathrm{Aut}(\mathbb{P}^2(\mathbb{C}))\simeq\mathrm{PGL}_3(\mathbb{C})$
par conjugaison dynamique sur~$\mathrm{Rat}_2$
\begin{align*}
& \mathrm{PGL}_3(\mathbb{C})\times\mathrm{Rat}_2\to\mathrm{Rat}_2, && (A,f)\mapsto AfA^{-1}.
\end{align*}
Cette action passe bien s\^ur au quotient $\mathrm{Rat}_2^{\bullet}:$
$(A,f^{\bullet})\mapsto Af^{\bullet}A^{-1}.$
L'orbite d'un \'el\'ement $f$ de~$\mathrm{Rat}_2$ (resp. $\mathrm{Rat}_2^{\bullet}$)
pour cette action sera not\'ee $\mathcal{O}_{dyn}(f)$ (resp.
$\mathcal{O}_{dyn}(f^{\bullet})$).

La seconde de nature alg\'ebrico-g\'eom\'etrique est
attach\'ee au point de vue \og r\'eseaux homalo\"{\i}daux\fg.
C'est une action du produit $\mathrm{PGL}_3(\mathbb{C})\times
\mathrm{PGL}_3(\mathbb{C}),$ appel\'ee action gauche-droite
(g.d.), qui se d\'efinit comme suit
\begin{align*}
&\mathrm{PGL}_3(\mathbb{C})\times\mathrm{PGL}_3(\mathbb{C})\times\mathrm{Rat}_2\to\mathrm{Rat}_2,&&
(A,B,f)\mapsto AfB^{-1}.
\end{align*}
Comme pr\'ec\'edemment cette action passe au quotient
$\mathrm{Rat}_2^{\bullet}.$ L'orbite d'un \'el\'ement $f$ est not\'ee
ici $\mathcal{O}_{g.d.}(f).$

Consid\'erons le projectivis\'e de l'espace
$\mathbb{C}[x_0, x_1,x_2]_3$ des polyn\^omes homog\`enes de
degr\'e $3.$ C'est un $\mathbb{P}^9(\mathbb{C})$ sur lequel
$\mathrm{PGL}_3(\mathbb{C})$ agit de fa\c{c}on naturelle: c'est la
projectivisation de l'action de $\mathrm{GL}_3(\mathbb{C})$ sur
$\mathbb{C}[x_0,x_1,x_2]_3$
\begin{align*}
(B,\varphi)\mapsto \varphi\circ B^{-1}.
\end{align*}
Les orbites g\'en\'eriques sont de dimension $8,$ ainsi
cette action induit un feuilletage singulier sur $\mathbb{P}^9
(\mathbb{C}),$ feuilletage qui poss\`ede une int\'egrale
premi\`ere rationnelle $J$ qui n'est autre que l'invariant
des courbes elliptiques (\cite{Ha}). Introduisons maintenant
l'application rationnelle
\begin{align*}
\det\mathrm{jac}\hspace{0.1cm}\colon\hspace{0.1cm}\mathrm{Rat}_2\simeq\mathbb{P}^{17}
(\mathbb{C})\to\mathbb{P}^9(\mathbb{C})
\end{align*}
provenant de l'application homog\`ene
\begin{align*}
(f_0,f_1,f_2)\mapsto\det\left(\frac{\partial f_i}{\partial x_j}\right).
\end{align*}
L'action gauche-droite induit sur $\mathrm{Rat}_2$ un feuilletage
singulier de codimension $1$ not\'e $\mathscr{R},$ les
feuilles g\'en\'eriques de $\mathscr{R}$ \'etant les g.d.
orbites g\'en\'eriques. Il se trouve que l'application rationnelle
$J\circ\det\mathrm{jac}$ en est une int\'egrale premi\`ere rationnelle.
Nous \'etudierons relativement pr\'ecis\'ement ce feuilletage
$\mathscr{R}.$ En particulier nous verrons que l'adh\'erence
ordinaire $\overline{\mathrm{Bir}_2}$ de $\mathrm{Bir}_2$ dans $\mathrm{Rat}_2$ est
une certaine composante irr\'eductible de l'ensemble
singulier du lieu singulier de~$\mathscr{R}.$ L'espace
$\overline{\mathrm{Bir}_2},$ qui est comme on vient de le dire
irr\'eductible, n'est pas lisse. Il s'agit en fait d'une \og
vari\'et\'e rationnelle d\'eterminantielle\fg\hspace{0.1cm} de dimension
$14.$ Par contre l'espace $\mathrm{Bir}_2$ est lui lisse: donnons-en
une br\`eve description. On v\'erifie qu'un \'el\'ement
$f$ de $\mathrm{Rat}_2$ induit une transformation birationnelle
si et seulement si (ceci sera pr\'ecis\'e dans le texte) l'application
correspondante $f^{\bullet}$ contracte les z\'eros de
$\det\mathrm{jac} f.$ La configuration de ces z\'eros est alors l'une
des suivantes
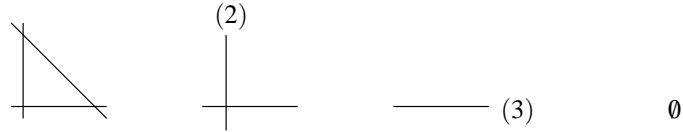
\begin{figure}[H]
\begin{center}
\input{quad.pstex_t}
\end{center}
\end{figure}
\noindent{\it i.e.} correspond aux configurations de triplets de droites
(except\'e la configuration de trois droites concourantes)
auxquelles on ajoute la configuration vide (cette configuration
est associ\'ee aux \'el\'ements $f$ de $\overline{\mathrm{Bir}_2
\setminus\mathring{\mathrm{B}}\mathrm{ir}_2}$ tels que $f^{\bullet}$ soit un isomorphisme
lin\'eaire). De fa\c{c}on duale on constate que les configurations
de points d'ind\'etermination sont

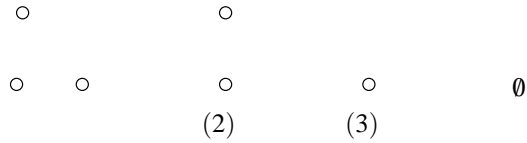
\begin{figure}[H]
\begin{center}
\input{quad2.pstex_t}
\end{center}
\end{figure}

Pour $i=1,$ $2,$ $3$ on note $\Sigma^i$ l'ensemble
des \'el\'ements $f$ de $\mathrm{Bir}_2$ tels que $f^{\bullet}$ ait $i$
points d'ind\'etermination; $\Sigma^0$ d\'esigne l'ensemble
des \'el\'ements de $\mathrm{Rat}_2$ du type $f=\ell(\ell_0,\ell_1,\ell_2)$
o\`u $\ell$ et les $\ell_i$ sont des formes lin\'eaires, les
$\ell_i$ \'etant ind\'ependantes. En fait
\begin{align*}
\Sigma^0=\{f\in\mathrm{Rat}_2\hspace{0.1cm}|\hspace{0.1cm}
f^{\bullet}\in\mathrm{Aut}(\mathbb{P}^2(\mathbb{C}))\}
\end{align*}
correspond exactement aux configurations vides ci-dessus.
Avec les notations pr\'ec\'edentes on~a
\begin{align*}
&(\Sigma^0)^{\bullet}=\mathrm{PGL}_3(\mathbb{C}), &&
(\Sigma^i)^{\bullet}\simeq \Sigma^i\text{ pour }i=1,2,3.
\end{align*}
Les isomorphismes $(\Sigma^i)^\bullet\simeq\Sigma^i$ justifieront
 la confusion parfois abusive des notations $f$ et $f^\bullet.$
Les $\Sigma^i$ sont en fait des orbites de l'action gauche-droite,
plus pr\'ecis\'ement
\begin{align*}
& \Sigma^0=\mathcal{O}_{g.d.}(x_0(x_0:x_1:x_2)), && \Sigma^1=
\mathcal{O}_{g.d.}(\tau), && \Sigma^2=\mathcal{O}_{g.d.}(\rho), &&
\Sigma^3=\mathcal{O}_{g.d.}(\sigma).
\end{align*}
Les dimensions de ces orbites sont les suivantes
\begin{align*}
& \dim\Sigma^0=10, &&\dim\Sigma^1=12, &&\dim\Sigma^2=13,
&&\dim\Sigma^3=14
\end{align*}
et l'on a
\begin{align*}
& \mathrm{Bir}_2=\Sigma^0\cup\Sigma^1\cup\Sigma^2\cup\Sigma^3.
\end{align*}
On peut de plus mentionner les relations d'incidence suivantes
(les adh\'erences sont prises dans~$\mathrm{Bir}_2$)
\begin{align*}
& \overline{\Sigma^0}=\Sigma^0, &&\overline{\Sigma^1}=
\Sigma^0\cup\Sigma^1,&&\overline{\Sigma^2}=\Sigma^0
\cup\Sigma^1\cup\Sigma^2,&&\mathrm{Bir}_2=\overline{\Sigma^3}
=\Sigma^0\cup\Sigma^1\cup\Sigma^2\cup\Sigma^3.
\end{align*}
Ces relations d'incidence se prolongent \'evidemment aux
$(\Sigma^i)^{\bullet}$ dans $\mathrm{Bir}_2^{\bullet}.$

\vspace{6mm}

Dans une session de probl\`emes (\cite{Mu}), lors du congr\`es
international de Moscou, \textsc{Mumford} propose ce qui
suit\footnote{\hspace{0.1cm}
Nous remercions Daniel \textsc{Panazzolo} pour cette r\'ef\'erence.}.

\selectlanguage{english}
\begin{quotation}
\hspace{0.1cm}\og Let $\mathrm{G}= \mathrm{Aut}_\mathbb{C}\mathbb{C}
(x_0,x_1)$ be the \textsc{Cremona} group (...) The problem is to
topo\-logize $\mathrm{G}$ and associate to it a \textsc{Lie} algebra
consisting, roughly, of those meromorphic vector fields $D$ on
$\mathbb{P}^2(\mathbb{C})$ wich ``integrate'' into an analytic
family of \textsc{Cremona} transformations\fg.
\end{quotation}

\selectlanguage{french}
Dans cette optique
nous d\'ecrivons au Chapitre \ref{germgerm} les sous-groupes \`a
un param\`etre de transformations birationnelles quadratiques.
Bien que $\mathrm{Bir}_2^{\bullet}$ ne soit pas un groupe, il
contient de nombreux groupes, en particulier le groupe lin\'eaire
$\mathrm{PGL}_3(\mathbb{C}).$ Comme on l'a dit pr\'ec\'edemment
l'adh\'erence de \textsc{Zariski} d'un tel sous-groupe est
naturellement munie d'une structure de groupe alg\'ebrique. Outre
le probl\`eme pos\'e par \textsc{Mumford}, la description
syst\'ematique des groupes de transformations birationnelles
quadratiques conduit naturellement \`a d\'ecrire les germes de
flots dans $\mathrm{Bir}_2^{\bullet}.$ Un tel objet est la
donn\'ee d'un germe d'application holomorphe
\begin{align*}
&\mathbb{C},0\to\mathrm{Bir}(\mathbb{P}^2(\mathbb{C})), &&
t\mapsto\phi_t
\end{align*}
tel que
\begin{align*}
&\phi_0^{\bullet}=\mathrm{id}, &&\phi_{t+s}^{\bullet}=\phi_t^\bullet
\circ\phi_s^\bullet.
\end{align*}
Si $\phi_t^{\bullet}$ est non lin\'eaire, {\it i.e.} si pour $t$
g\'en\'erique $\phi_t^{\bullet}$ n'est pas dans
$\mathrm{PGL}_3(\mathbb{C}),$ alors $\phi_t^{\bullet}$ contracte
au moins une droite $\mathcal{D}_t$ qui est ou bien immobile
(c'est-\`a-dire ind\'ependante de $t$), ou bien mobile. Par
exemple on montre d'abord que si un germe de flot contracte une
droite \og mobile\fg, ce germe laisse une fibration en droites
invariante; dans le m\^eme ordre d'id\'ee on obtient que si un
germe de flot contracte une unique droite qui est \og immobile\fg,
ce germe de flot est polynomial... En fait comme nous le verrons,
dans tous les cas un germe de flot dans $\mathrm{Bir}_2$
pr\'eserve une fibration en droites.

Donnons un r\'esultat qui s'av\`erera essentiel pour la
classification des flots (\S \ref{cla}): soit $\phi_t^{\bullet}$
un germe de flot quadratique. Le champ de vecteurs rationnel
$\chi= \frac{\partial\phi_t^{\bullet}}{\partial t}\Big|_{t=0}$
s'appelle le g\'en\'erateur infinit\'esimal de $\phi_t^{\bullet};$
\`a ce champ est associ\'e un feuilletage alg\'ebrique not\'e
$\mathcal{F}_\chi.$ On introduit $\langle\phi_t^{\bullet}\rangle\subset
\mathrm{Bir}_2^{\bullet}$ le groupe engendr\'e par les
$\phi_t^{\bullet}$ et
$\overline{\langle\phi_t^{\bullet}\rangle}^{\hspace{0.1cm} \mathsf{Z}}$ son
adh\'erence de \textsc{Zariski} dans $\mathrm{Bir}_2^{\bullet};$
ce sont des groupes ab\'eliens. Enfin $\mathrm{G}(\chi)$ d\'esigne
le groupe ab\'elien alg\'ebrique maximal inclus dans
$\mathrm{Bir}_2^{\bullet}$ et contenant
$\overline{\langle\phi_t^{\bullet}\rangle}^{\hspace{0.1cm}\mathsf{Z}}.$
Alors
\begin{itemize}
\item $\mathcal{F}_\chi$ est d\'efini par une forme ferm\'ee rationnelle;

\item $\phi_t^{\bullet}$ pr\'eserve une fibration rationnelle
(en fait un pinceau de droites);

\item $\dim\mathrm{G}(\chi)\geq 2.$
\end{itemize}

Notons qu'un tel \'enonc\'e est satisfait par n'importe
quel groupe \`a un param\`etre de $\mathrm{PGL}_3(\mathbb{C}).$ Il
nous permet d'obtenir la classification \`a conjugaison pr\`es des
groupes \`a un param\`etre de transformations birationnelles
quadratiques.

En fait ce r\'esultat se g\'en\'eralise en tout degr\'e:
tout germe de flot birationnel laisse une fibration rationnelle
invariante et se rel\`eve en un flot d'automorphismes sur
une surface minimale.

\vspace{6mm}

Dans un troisi\`eme Chapitre on implante, suivant une
id\'ee de \textsc{Guillot}, aux transformations birationnelles
des invariants propres aux feuilletages.
\`A une transformation rationnelle $f=(f_0:f_1:f_2)$ on
associe la $1$-forme diff\'erentielle
\begin{align*}
\omega_f=(x_2f_1-x_1f_2)\mathrm{d}x_0+(x_0f_2-x_2f_0)\mathrm{d}x_1
+(x_1f_0-x_0f_1)\mathrm{d}x_2.
\end{align*}
Lorsque $f$ n'est pas un multiple de l'identit\'e, la forme $\omega_f$
est non identiquement nulle et d\'efinit donc un feuilletage $\mathcal{F}(f)$
sur $\mathbb{P}^2(\mathbb{C}).$ Si $f$ est de degr\'e $n,$ le feuilletage $\mathcal{F}(f)$ est
de degr\'e inf\'erieur ou \'egal \`a $n.$ Dans le cas quadratique
($n=2$), l'application $\mathcal{F}\hspace{0.1cm}\colon\hspace{0.1cm} f\mapsto\omega_f$ induit
une application lin\'eaire
\begin{align*}
\mathcal{F}\hspace{0.1cm}\colon\hspace{0.1cm}\mathrm{Rat}_2\simeq\mathbb{P}^{17}(\mathbb{C})\dashrightarrow\mathcal{F}_2
\simeq\mathbb{P}^{14}(\mathbb{C})
\end{align*}
o\`u $\mathcal{F}_2$ est le projectivis\'e de l'espace
vectoriel des $1$-formes homog\`enes de degr\'e $3$
\begin{align*}
\omega=A_0\mathrm{d}x_0+A_1\mathrm{d}x_1+A_2\mathrm{d}x_2
\end{align*}
satisfaisant la condition d'\textsc{Euler} $\sum x_iA_i=0.$
L'espace des feuilletages proprement dit est
\begin{align*}
&\mathcal{F}_2^{\bullet}=\mathcal{F}_2/\sim \text{ o\`u }\sim 
\text{ est la relation d'\'equivalence }\omega_1\sim\omega_2\Leftrightarrow
\omega_1\wedge\omega_2=0.
\end{align*}

Visiblement l'application $\mathcal{F}\hspace{0.1cm}\colon\hspace{0.1cm}\mathrm{Rat}_2\dashrightarrow
\mathcal{F}_2$ est surjective. Ce qui est plus curieux c'est que sa
restriction $\mathcal{F}_{|\mathrm{Bir}_2}\hspace{0.1cm}\colon\hspace{0.1cm}\mathrm{Bir}_2\dashrightarrow\mathcal{F}_2$
est dominante et g\'en\'eriquement finie: un \'el\'ement g\'en\'erique
de $\mathcal{F}_2$ poss\`ede exactement $35$ ant\'ec\'edents dans
$\mathrm{Bir}_2.$

Le feuilletage $\mathcal{F}(f)$ n'a pas de sens dynamique particulier
vis \`a vis de $f:$ en g\'en\'eral $\mathcal{F}(f)$ n'est pas invariant par
$f.$ Les transformations quadratiques $f$ pour lesquelles
cette propri\'et\'e est satisfaite sont exactement celles qui
pr\'eservent une fibration en droites fibre \`a fibre: dans ce
cas $\mathcal{F}(f)$ est pr\'ecis\'ement cette fibration.

Soit $f$ une transformation birationnelle quadratique, {\it i.e.}
un \'el\'ement g\'en\'erique de $\Sigma^3\simeq
(\Sigma^3)^{\bullet}.$ Alors $f$ poss\`ede $7$ points sp\'eciaux:
$3$ points d'ind\'etermination et $4$ points fixes. Ces $7$ points
sont exactement les points singuliers du feuilletage $\mathcal{F}(f);$ la
position des $3$ points d'ind\'etermination et des $4$ points
fixes d\'etermine compl\`etement l'application $f.$

Un th\'eor\`eme de \textsc{Guillot} (\cite{Gu}) affirme que pour
un \'el\'ement g\'en\'erique $f$ de $\Sigma^3$ on a les relations
suivantes
\begin{align*}
&\displaystyle\sum_{m\in\mathrm{Fix} f} \frac{\mathrm{tr}(Df_{(m)}
-\mathrm{id})}{\det(Df_{(m)}-\mathrm{id})}=-4,
&& \displaystyle\sum_{m\in\mathrm{Fix} f}\frac{1}{\det(Df_{(m)}
-\mathrm{id})}=1,
\end{align*}
$Df_{(m)}$\label{not3d} \'etant la diff\'erentielle de $f$ en $m,$ $\mathrm{tr}(
Df_{(m)})$\label{not3e} sa trace et $\det Df_{(m)}$ son d\'eterminant.
Ce sont des relations qui lient les valeurs propres
de la diff\'erentielle de $f$ aux points fixes; de plus ce sont \og les
seules\fg\hspace{0.1cm} (\cite{Gu2}). On peut tirer des cons\'equences int\'eressantes
de ces formules. Par exemple pour une transformation quadratique
g\'en\'erique, {\it i.e.} avec $7$ points sp\'eciaux, les valeurs
propres de la partie lin\'eaire de $f$ aux points fixes ne peuvent
toutes \^etre de module inf\'erieur strictement \`a $1:$ les points
fixes ne peuvent tous \^etre des attracteurs.

Les relations aux valeurs propres interviennent de fa\c{c}on forte
dans les probl\`emes de classification; mentionnons par
exemple le fait suivant: si deux transformations birationnelles
quadratiques g\'en\'eriques sont birationnellement conjugu\'ees
alors elles le sont lin\'eairement. On peut se demander si deux
transformations birationnelles quadratiques g\'en\'eriques
topologiquement conjugu\'ees le sont par un automorphisme
holomorphe (ou anti-holomorphe); c'est une question ouverte.

\vspace{6mm}

Dans le Chapitre \ref{dyndyn} on aborde quelques exemples de
propri\'et\'es dynamiques des transformations birationnelles, en
particulier quadratiques. Les notions fondamentales de stabilit\'e
alg\'ebrique et de premier degr\'e dynamique d\'ej\`a
mentionn\'ees ci-dessus sont introduites. Nous donnons des
exemples explicites d'\'el\'ements de $\Sigma^3$ qui sont
alg\'ebriquement stables, exemples qui permettent de confirmer la
g\'en\'ericit\'e intuitive de cette propri\'et\'e. Nous nous
int\'eresserons ensuite \`a la pr\'esence ou non \og d'objets
invariants\fg: courbes, feuilletages. Si $\phi$ est un \'el\'ement
du groupe de \textsc{Cremona} et $\mathcal{F}$ un feuilletage holomorphe
(singulier) sur $\mathbb{P}^2(\mathbb{C}),$ on dit que $\mathcal{F}$ est invariant par
$\phi$ si pour toute feuille g\'en\'erique $\mathcal{L}$ de
$\mathcal{F}$ alors $\phi(\mathcal{L})$ est encore une feuille de~$\mathcal{F}.$ 
Si $\mathcal{F}$ est d\'efini par la $1$-forme rationnelle
$\omega,$ alors $\mathcal{F}$ est invariant par $\phi$ si et seulement si
$\phi^*\omega \wedge\omega=0.$ Consid\'erons un \'el\'ement $A$ de
$\mathrm{Aut}(\mathbb{P}^2(\mathbb{C}));$ \'evidemment $A$ poss\`ede un point fixe~$q,$ 
{\it i.e.} $Aq=q,$ et l'image par $A$ d'une droite passant par $q$ est
encore une droite passant par~$q.$ De sorte que $A$ poss\`ede au
moins un feuilletage invariant: la fibration rationnelle form\'ee
par les droites passant par $q.$ En fait $A$ pr\'eserve une
infinit\'e de feuilletages distincts. Par exemple si~$A$ est
diagonalisable, on peut supposer, \`a conjugaison lin\'eaire
pr\`es, que $A$ s'\'ecrit dans la carte affine $x_2=1$
\begin{align*}
(x_0,x_1)\mapsto(\alpha_0x_0,\alpha_1x_1).
\end{align*}
On constate alors que pour chaque nombre complexe
$\beta$ la forme diff\'erentielle
\begin{align*}
\omega_\beta:=\frac{\mathrm{d}x_0}{x_0}+\beta\frac{\mathrm{d}x_1}{x_1}
\end{align*}
est invariante par $A,$ ce qui, comme annonc\'e, produit une
infinit\'e de feuilletages invariants. Qu'en est-il pour les
transformations birationnelles quadratiques ? En invoquant
l'irr\'eductibilit\'e de $\mathrm{Bir}_2$ dans $\mathrm{Rat}_2$ et
le fait qu'un automorphisme de \textsc{H\'enon} quadratique ne
poss\`ede pas de feuilletage invariant (\cite{Br2}) on constate
que g\'en\'eriquement un \'el\'ement de $\mathrm{Bir}_2$ n'en
poss\`ede pas non plus.

De la m\^eme fa\c{c}on, \`a l'inverse des automorphismes, une
transformation birationnelle quadratique g\'en\'erique ne poss\`ede
pas de courbe invariante. Toutefois, comme nous l'avons d\'ej\`a
mentionn\'e, l'\'etude des transformations ayant une courbe
invariante poss\`ede une longue histoire et pr\'esente un
int\'er\^et certain. Nous y apportons quelques enrichissements
en utilisant les concepts r\'ecents de stabilit\'e alg\'ebrique et
premier degr\'e dynamique.

Le paragraphe \S\ref{dynpt} de ce chapitre est consacr\'e \`a la
r\'epartition des points p\'eriodiques et des points
d'ind\'etermination des it\'er\'es d'une transformation
birationnelle quadratique. A priori ces ensembles sont tr\`es
compliqu\'es; nous montrons la densit\'e (au sens de
\textsc{Zariski}) g\'en\'erique de ces deux ensembles. En
utilisant la condition introduite par \textsc{Bedford} et
\textsc{Diller} pr\'ec\'edemment mentionn\'ee on constate que g\'en\'eriquement une
transformation quadratique poss\`ede une infinit\'e de points
p\'eriodiques hyperboliques.

La non stabilit\'e alg\'ebrique d'une transformation birationnelle $f$
est reli\'ee \`a un comportement \og anormal\fg\hspace{0.1cm} du degr\'e de
ses it\'er\'es $f^n;$ dans le cas quadratique ce degr\'e
doit \^etre en principe~$2^n.$ Mais ce n'est pas toujours le cas,
par exemple pour les transformations qui appartiennent \`a un groupe
\`a un param\`etre.
Nous consacrons un paragraphe \`a l'\'etude du cas extr\^eme,
celui o\`u le degr\'e de $f^2$ reste \'egal \`a $2;$ nous en pr\'esentons
la classification \`a conjugaison lin\'eaire pr\`es.

\vspace{6mm}

Dans le Chapitre \ref{algalg} on s'int\'eresse \`a quelques
propri\'et\'es alg\'ebriques du groupe $\mathrm{Bir}(\mathbb{P}^2
(\mathbb{C}))$ mettant le plus souvent en jeu les transformations
quadratiques.
Le fait que $\mathrm{Bir}(\mathbb{P}^2(\mathbb{C}))$ pr\'esente beaucoup de
propri\'et\'es de groupes lin\'eaires am\`ene \`a se poser la
question suivante: existe-t-il une repr\'esentation fid\`ele de
$\mathrm{Bir}(\mathbb{P}^2(\mathbb{C}))$ dans un groupe lin\'eaire
$\mathrm{GL}_n(\Bbbk)$ pour un certain corps~$\Bbbk$ de
caract\'eristique nulle ? Nous constatons en utilisant un argument
de \textsc{Birkhoff} que la r\'eponse est n\'egative.

On montre ensuite que le centralisateur dans $\mathrm{Bir}(
\mathbb{P}^2(\mathbb{C}))$
d'un \'el\'ement g\'en\'erique $f$ de $\Sigma^3$ est isomorphe
\`a $\mathbb{Z},$ {\it i.e.} se r\'eduit aux it\'er\'es de $f.$ Puis on \'etablit
qu'un nombre fini de transformations quadratiques choisies
g\'en\'eriquement engendrent un groupe libre.

Le probl\`eme de la simplicit\'e de $\mathrm{Bir}(\mathbb{P}^2(
\mathbb{C}))$ est
un probl\`eme ouvert. Nous montrons en particulier que si $f$ est
une transformation birationnelle quadratique ou un automorphisme
lin\'eaire non trivial, le sous-groupe normal engendr\'e par $f$
est le groupe de \textsc{Cremona} tout entier; nous en d\'eduisons
que $\mathrm{Bir}(\mathbb{P}^2(\mathbb{C}))$ est parfait.
Nous montrons par exemple que le sous-groupe normal
engendr\'e par une transformation de \textsc{de Jonqui\`eres}
dans $\mathrm{Bir}(\mathbb{P}^2(\mathbb{C}))$ est $\mathrm{Bir}(\mathbb{P}^2(
\mathbb{C}))$ tout entier.

Nous pr\'ecisons le th\'eor\`eme de \textsc{N\oe ther} de la
fa\c{c}on suivante: tout \'el\'ement $f$ de $\mathrm{Bir}(\mathbb{P}^2(
\mathbb{C}))$ s'\'ecrit comme un produit d'involutions
quadratiques standards
\begin{align*}
&f=A_1\sigma A_1^{-1}\ldots A_n\sigma A_n^{-1}, && A_i\in
\mathrm{PGL}_3(\mathbb{C}).
\end{align*}

\vspace{6mm}

Le Chapitre \ref{expexp} est exp\'erimental. On se propose
dans des cas tr\`es sp\'eciaux de donner une id\'ee de ce
que peut \^etre la dynamique d'une application birationnelle.
\`A cet effet nous essayons de repr\'esenter quelques orbites
comme par exemple

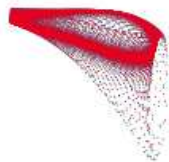
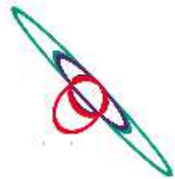
\vspace{6mm}
\begin{tabular}{cc}
\hspace{20mm}\input{des9scan.pstex_t}\hspace{18mm} & \hspace{18mm}
\input{iter8scan.pstex_t}
\hspace{5mm}
\end{tabular}

\bigskip

\noindent et nous reprenons un protocole bien connu pour le
trac\'e des ensembles de \textsc{Julia} des transformations
rationnelles en une variable. Nous appliquons ce proc\'ed\'e \`a
des transformations birationnelles r\'eelles ({\it i.e.} commutant
\`a l'involution anti-holomorphe standard)
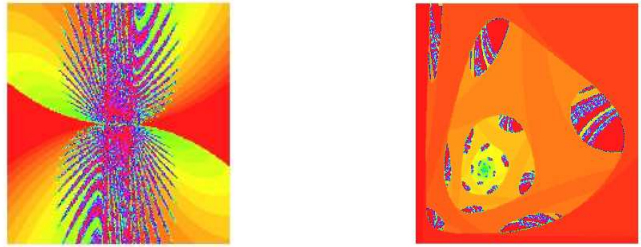
\begin{center}
 \begin{tabular}{cc}
\input{julia49a.pstex_t} &\input{bd48.pstex_t}
\end{tabular}
\end{center}

\vspace{6mm}

Comme on l'a dit l'espace des transformations birationnelles
quadratiques pures est lisse et connexe. Qu'en est-il en degr\'e
quelconque ? Le premier degr\'e int\'eressant est le degr\'e $3$
pour la raison suivante. D'apr\`es le th\'eor\`eme de \textsc{N\oe
ther} toute transformation birationnelle s'\'ecrit comme un
produit 
\begin{align*}
&A_0\sigma A_1\sigma\ldots\sigma A_n, && A_i\in\mathrm{Aut}(
\mathbb{P}^2(\mathbb{C})); 
\end{align*}
lorsque les $A_i$ sont choisis
g\'en\'eriquement le degr\'e d'un tel \'el\'ement est $2^n.$ Pour
obtenir des applications birationnelles de degr\'e $d$ qui n'est
pas une puissance de $2,$ il est donc n\'ecessaire que les
transformations $A_i$ soient assujetties \`a des contraintes. Le
premier degr\'e pour lequel ce ph\'enom\`ene appara\^it est $d=3.$
Dans les trait\'es anciens on trouve une \og description\fg\hspace{0.1cm} des
transformations birationnelles cubiques, description qui repose
sur des consid\'erations de g\'eom\'etrie \'enum\'erative (nombre
de points d'ind\'etermination, configuration de courbes
contract\'ees). Au Chapitre \ref{troistrois} nous proposons une liste de
formes normales \`a conjugaison gauche-droite pr\`es, la
connexit\'e apparaissant en dernier lieu comme sous-produit. Les
techniques utilis\'ees sont pour la plupart classiques: topologie
du compl\'ement de certaines courbes planes, contraction des
z\'eros du d\'eterminant jacobien etc. Malheureusement d\`es le
degr\'e $3$ on ne dispose pas de crit\`ere de birationnalit\'e
analogue \`a celui utilis\'e en degr\'e $2$ comme l'atteste
l'exemple suivant: si $f$ est la transformation
$(x_0^2x_1:x_0x_2^2:x_1^2x_2),$ le lieu d'annulation de
$\mathrm{det}\hspace{0.1cm}\mathrm{jac}\hspace{0.1cm} f$ est bien contract\'e alors que
$f$ n'est pas inversible. Toutefois si $f$ est birationnelle, la
courbe $\mathrm{det}\hspace{0.1cm}\mathrm{jac}\hspace{0.1cm} f=0$ est contract\'ee par
$f$ et ceci s'av\`ere dans bon nombre de cas int\'eressant \`a
exploiter. Nous montrons qu'en degr\'e $3$ les configurations possibles de courbes
contract\'ees sont certaines unions de droites et coniques
sch\'ematis\'ees par les figures suivantes.
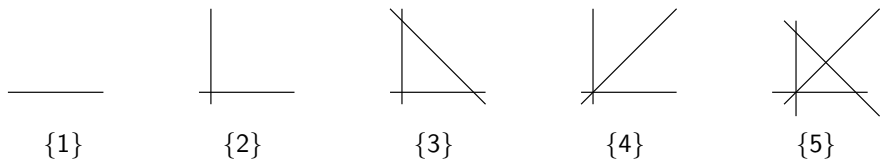
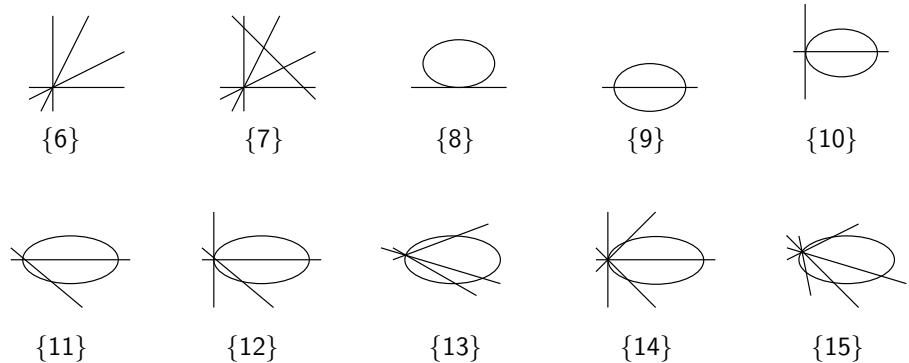
\begin{figure}[H]
\begin{center}
\input{gal.pstex_t}
\end{center}
\end{figure}

\begin{figure}[H]
\begin{center}
\input{gal2.pstex_t}
\end{center}
\end{figure}

\begin{figure}[H]
\begin{center}
\input{gal3.pstex_t}
\end{center}
\end{figure}

On donnera la liste explicite, \`a conjugaison
gauche-droite pr\`es, des transformations correspondant \`a ces
configurations.
En particulier toute transformation birationnelle cubique est \`a
conjugaison dynamique pr\`es du type $Af$ o\`u $A$ d\'esigne un
\'el\'ement de $\mathrm{PGL}_3(\mathbb{C})$ et $f$ un \'el\'ement de la
liste \'evoqu\'ee ci-dessus. Cette classification permettra
d'affirmer que l'\'el\'ement \og g\'en\'erique\fg\hspace{0.1cm} pr\'esente la
derni\`ere configuration et d'\'etablir que l'espace
$\mathring{\mathrm{B}}\mathrm{ir}_3$ des transformations
birationnelles purement de degr\'e $3$ est de dimension $18.$ \`A
conjugaison gauche-droite pr\`es les \'el\'ements qui pr\'esentent la
configuration g\'en\'erique $\mathsf{\{15\}}$ forment une famille
\`a deux param\`etres: alors qu'en degr\'e $2$ il y a trois
orbites gauche-droite, ici il y en a une infinit\'e.

Nous d\'eterminons aussi les dimensions des orbites sous
l'action gauche-droite des mod\`eles;
toujours \`a l'aide de ceux-ci on obtient que la
d\'ecomposition de \textsc{N\oe ther} de tout \'el\'ement de
l'espace $\mathrm{Bir}_3$ des transformations birationnelles de
degr\'e inf\'erieur ou \'egal \`a $3$ fait intervenir au plus huit
fois l'involution de \textsc{Cremona}.

On note que toutes les configurations de
d\'eg\'en\'erescence de la figure $\mathsf{\{15\}}$ n'apparaissent
pas. En degr\'e $2$ il y a d\'ej\`a un ph\'enom\`ene analogue: la
configuration de trois droites concourantes n'est pas r\'ealis\'ee
comme ensemble exceptionnel d'une transformation birationnelle
quadratique.

Notons $\mathscr{X}$ l'ensemble des transformations birationnelles
purement cubiques de configuration $\mathsf{\{15\}}.$ Nous
\'etablissons que l'adh\'erence ordinaire de $\mathscr{X}$ dans $\mathring{\mathrm{B}}\mathrm{ir}_3$
est $\mathring{\mathrm{B}}\mathrm{ir}_3.$ Nous montrons que $\mathring{\mathrm{B}}\mathrm{ir}_3$ est irr\'eductible (en fait
rationnellement connexe); par contre alors que $\mathrm{Bir}_2$ est lisse et
irr\'eductible nous verrons qu'il n'en est
pas de m\^eme pour $\mathrm{Bir}_3$ vu dans $\mathbb{P}^{29}(\mathbb{C})\simeq
\mathrm{Rat}_3.$

En utilisant la \og correspondance
feuilletages-transformations birationnelles\fg\hspace{0.1cm} on obtient qu'un
\'el\'ement g\'en\'erique de $\mathrm{Bir}_3$ poss\`ede $5$ points
fixes, tous de multiplicit\'e $1,$ et $5$ points
d'ind\'etermination, $1$ de multiplicit\'e $4,$ les autres de
multiplicit\'e $1.$ De plus toute transformation birationnelle
cubique g\'en\'erique est d\'etermin\'ee par la position de ses
$5$ points fixes et de ses $5$ points d'ind\'etermination
affect\'es de leur multiplicit\'e.\vspace{0.9cm}

\subsection*{Remerciements.}\hspace{1mm}

Il a \'et\'e tr\`es plaisant de travailler
\`a l'IMPA et au CIRM; nous tenons  \`a remercier ces
deux institutions. Merci \`a l'IUF et l'ANR Symplexe (ANR BLAN$06$-$3$\_$137237$) qui ont contribu\'e au bon d\'eroulement de cette collaboration.

Nous avons profit\'e de discussions agr\'eables et fructueuses
 avec Ivan \textsc{Pan} lors de son s\'ejour \`a Rennes.

Merci \`a Mark \textsc{Baker} et Charles \textsc{Favre} pour leur
disponibilit\'e.

\clearemptydoublepage
\chapter{Transformations rationnelles et
birationnelles quadratiques}\label{gen}

\section{Quelques d\'efinitions et notations}\hspace{0.1cm}

Nous g\'en\'eralisons et pr\'ecisons les notations
pr\'esent\'ees dans l'Introduction. On d\'esigne par~$\mathrm{Rat}_k\label{not4}$ le projectivis\'e de l'espace des
triplets de polyn\^omes homog\`enes de degr\'e $k$ en $3$
variables
\begin{align*}
\mathrm{Rat}_k=\mathbb{P}\{(f_0:f_1:f_2)\hspace{0.1cm}|
\hspace{0.1cm} f_i\in\mathbb{C}[x_0,x_1,x_2]_k\}.
\end{align*}

Le \textbf{\textit{degr\'e}}\label{ind2} d'un \'el\'ement $f=(f_0:f_1:f_2)$ de $\mathrm{Rat}_k$ est
le degr\'e des $f_i.$

\`A un \'el\'ement $f=(f_0:f_1:f_2)$ de $\mathrm{Rat}_k$ on
associe la transformation birationnelle
\begin{align*}
&f^\bullet=\delta(f_0:f_1:f_2), &&
\delta=\frac{1}{\text{pgcd}(f_0,f_1,f_2)}
\end{align*}

Soit $f$ un \'el\'ement de $\mathrm{Rat}_k;$ on dit que $f=
(f_0:f_1:f_2)$ est \textbf{\textit{purement de degr\'e
$k$}}\label{ind23} si les $f_i$ n'ont pas de facteur commun.
L'ensemble des \'el\'ements purement de degr\'e $k$ est not\'e
$\mathring{\mathrm{R}} \mathrm{at}_k\label{not5}.$ Alors 
que~$\mathrm{Rat}_k$ s'identifie \`a un espace projectif, l'espace
$\mathring{\mathrm{R}}\mathrm{at}_k$ en est un ouvert de
\textsc{Zariski}. Un \'el\'ement de $\mathrm{Rat}_k \setminus
\mathring{\mathrm{R}}\mathrm{at}_k$ s'\'ecrit $\psi g=(\psi
g_0:\psi g_1:\psi g_2)$ o\`u $g$ appartient \`a un certain
$\mathrm{Rat}_\ell,$ avec $\ell<~k,$ et $\psi$ est un polyn\^ome
homog\`ene de degr\'e $k-\ell.$ On note $\mathrm{Rat}\label{not6}$
l'ensemble des transformations rationnelles de $\mathbb{P}^2(\mathbb{C})$ dans
lui-m\^eme: c'est l'union $\displaystyle\bigcup_{k\geq
1}\mathring{\mathrm{R}}\mathrm{at}_k.$ C'est aussi la limite
projective des $\mathrm{Rat}^\bullet_k$ avec
\begin{align*}
\mathrm{Rat}_k^\bullet=\{f^\bullet\hspace{0.1cm}|\hspace{0.1cm} f\in\mathrm{Rat}_k\}.
\end{align*}\label{not3b}
Remarquons que si l'\'el\'ement $f$ de $\mathrm{Rat}_k$ est purement
de degr\'e $k$ alors $f$ s'identifie \`a $f^\bullet.$ Ceci
signifie que l'application
\begin{align*}
\mathring{\mathrm{R}}\mathrm{at}_k\to\mathrm{Rat}_k^\bullet
\end{align*}
est injective. Dans la suite lorsqu'il n'y a aucune
ambigu\"{\i}t\'e on utilisera donc la notation $f$ aussi bien pour
les \'el\'ements de $\mathrm{Rat}_k$ que pour ceux de
$\mathrm{Rat}_k^\bullet.$ De m\^eme nous dirons abusivement qu'un
\'el\'ement de $\mathrm{Rat}_k$ \og est\fg\hspace{0.1cm} une transformation
rationnelle.

L'espace $\mathrm{Rat}$ contient le groupe des
transformations birationnelles de $\mathbb{P}^2(\mathbb{C}),$ encore appel\'e
\textbf{\textit{groupe de \textsc{Cremona}}}\label{ind0}, que l'on note
$\mathrm{Bir}(\mathbb{P}^2(\mathbb{C}))\label{not7}.$ Il est parfois plus commode
de travailler dans les espaces $\mathrm{Rat}_k$ que dans les
$\mathrm{Rat}_k^\bullet.$ Aussi nous d\'esignons par
$\mathrm{Bir}_k\subset \mathrm{Rat}_k\label{not8}$ l'ensemble des
transformations $f$ de $\mathrm{Rat}_k$ telles que $f^\bullet$
soit inversible et par $\mathring{\mathrm{B}}\mathrm{ir}_k\subset
\mathrm{Bir}_k\label{not9}$ l'ensemble des transformations
birationnelles qui sont purement de degr\'e $k.$ On pose
\begin{align*}
\mathrm{Bir}_k^\bullet=\{f^\bullet\hspace{1mm}|\hspace{1mm}
f\in\mathrm{Bir}_k\}
\end{align*}\label{not3c}
Le groupe de
\textsc{Cremona} s'identifie alors \`a
$\displaystyle\bigcup_{k\geq
1}\mathring{\mathrm{B}}\mathrm{ir}_k.$ On constate que
$\mathring{\mathrm{B}}\mathrm{ir}_1\simeq \mathrm{PGL}_3(\mathbb{C})$ est
le groupe des automorphismes de $\mathbb{P}^2(\mathbb{C})~;$ visiblement
$\mathring{\mathrm{B}}\mathrm{ir}_1\simeq\mathrm{Bir}_1^\bullet
=\mathrm{Bir}_1.$ L'espace
$\mathrm{Rat}_1$ s'identifie \`a $\mathbb{P}^8(\mathbb{C})$ et
$\mathring{\mathrm{R}} \mathrm{at}_1$ est le projectivis\'e des
matrices de rang sup\'erieur ou \'egal \`a $2.$

Dans ce m\'emoire nous allons d\'ecrire les espaces
$\mathrm{Bir}_k$ et $\mathring{\mathrm{B}}\mathrm{ir}_k$ pour des petits
entiers, pr\'ecis\'ement pour $k=2$ et $k=3.$ Par exemple pour
$k=2$ l'inclusion $\mathring{\mathrm{B}}\mathrm{ir}_2\subset\mathrm{Bir}_2$ est stricte. En effet si
$A$ est un \'el\'ement de $\mathrm{PGL}_3(\mathbb{C})$ et $\ell$ une
forme lin\'eaire, $\ell A$ est dans $\mathrm{Bir}_2$ mais pas dans
$\mathring{\mathrm{B}}\mathrm{ir}_2.$

Sur l'espace $\mathrm{Rat}_k$ il y a deux actions naturelles. La
premi\`ere est l'action de $\mathrm{PGL}_3(\mathbb{C})$ par \textbf{\textit{conjugaison
dynamique}}\label{ind24}
\begin{align*}
&\mathrm{PGL}_3(\mathbb{C})\times\mathrm{Rat}_k\to\mathrm{Rat}_k, &&(A,Q) \mapsto AQA^{-1}
\end{align*}
et la seconde l'action statique que nous appellerons plut\^ot
\textbf{\textit{conjugaison gauche-droite}}\label{ind25} (g.d.) 
de~$(\mathrm{PGL}_3(\mathbb{C}))^2$
\begin{align*}
&\mathrm{PGL}_3(\mathbb{C})\times\mathrm{Rat}_k\times\mathrm{PGL}_3(\mathbb{C})\to\mathrm{Rat}_k,
&&(A,Q,B)\mapsto AQB^{-1}.
\end{align*}
Notons que les ensembles $\mathring{\mathrm{R}}\mathrm{at}_k,$
$\mathrm{Bir}_k$ et $\mathring{\mathrm{B}}\mathrm{ir}_k$ sont
invariants sous ces deux actions. Nous allons plus loin nous
int\'eresser \`a quelques orbites de ces actions. L'orbite d'un
\'el\'ement $Q$ de $\mathrm{Rat}_k$ sous l'action de
$\mathrm{PGL}_3(\mathbb{C})$ par conjugaison dynamique (resp. sous
l'action g.d.) est not\'ee $\mathcal{O}_{dyn} (Q)\label{not10}$
(resp. $\mathcal{O}_{g.d.} (Q)\label{not11}$).

Consid\'erons $f=(f_0:f_1:f_2)$ dans
$\mathring{\mathrm{B}}\mathrm{ir}_k.$ Nous notons $\mathrm{Ind}\hspace{0.1cm}
f\label{not12}$ (ou $\mathrm{Ind}\hspace{0.1cm} f^\bullet$) le
\textbf{\textit{lieu d'ind\'etermination}}\label{ind26} de $f,$
c'est le lieu d'annulation des polyn\^omes $f_i.$ Les
points de  $\mathrm{Ind}\hspace{0.1cm}f$ sont les
\textbf{\textit{points d'ind\'etermination}}\label{ind10} de $f.$
L'ensemble
$\mathrm{Exc}\hspace{0.1cm} f\label{not13}$ (ou $\mathrm{Exc}\hspace{0.1cm} f^\bullet$)
d\'esignera \textbf{\textit{l'ensemble exceptionnel}}\label{ind27}
de~$f,$ {\it i.e.} l'ensemble des \textbf{\textit{courbes contract\'ees}}\label{ind27aa}  par $f;$
c'est \og l'image\fg\hspace{0.1cm} par $f^{-1}$ de $\mathrm{Ind}\hspace{0.1cm} f^{-1}.$
On dit que $f$ \textbf{\textit{\'eclate}}\label{ind27a} les points
de $\mathrm{Ind}\hspace{0.1cm} f$ et
\textbf{\textit{contracte}}\label{ind27b} les courbes de
$\mathrm{Exc}\hspace{0.1cm} f.$

\begin{egs}
Soit $\sigma$ la transformation de \textsc{Cremona}
d\'efinie par
\begin{align*}
& \mathbb{P}^2(\mathbb{C})\dashrightarrow\mathbb{P}^2(\mathbb{C}), && (x_0:x_1:x_2)\mapsto
(x_1x_2:x_0x_2:x_0x_1).
\end{align*}
Notons que $\sigma$ est une involution dont le lieu
d'ind\'etermination et l'ensemble exceptionnel sont
\begin{align*}
&\mathrm{Ind}\hspace{0.1cm} \sigma=\{(1:0:0),\hspace{0.1cm} (0:1:0),\hspace{0.1cm} (0:0:1)\},
&& \mathrm{Exc}\hspace{0.1cm} \sigma=\{x_0=0,\hspace{0.1cm} x_1=0,\hspace{0.1cm} x_2=0\}.
\end{align*}

L'\'el\'ement $\rho=(x_0x_1:x_2^2:x_1x_2)\label{not14}$ de
$\mathrm{Bir}(\mathbb{P}^2(\mathbb{C}))$ a deux points d'ind\'etermination qui sont
$(0:1:0)$ et $(1:0:0);$ les droites contract\'ees par $\rho$ ont
pour \'equation $x_2=0,$ resp. $x_1=0.$

Soit $\tau\label{not15}$ la transformation d\'efinie
par $(x_0^2:x_0x_1:x_1^2-x_0x_2);$ on a
\begin{align*}
&\mathrm{Ind}\hspace{0.1cm}\tau=\{(0:0:1)\},&& \mathrm{Exc}\hspace{0.1cm}\tau=\{x_0=0\}.
\end{align*}

Remarquons que $\rho$ et $\tau$ sont aussi des involutions.

Les trois \'el\'ements $\sigma,$ $\rho$ et
$\tau$ sont birationnellement conjugu\'es \`a des involutions de
$\mathrm{PGL}_3(\mathbb{C})$ (\emph{voir} \cite{DI} par exemple).
\end{egs}

\section{Transformations rationnelles quadratiques}\hspace{0.1cm}

Soit $\mathbb{C}[x_0,x_1,x_2]_\nu\label{not16}$ l'ensemble
des polyn\^omes homog\`enes de degr\'e $\nu$ dans $\mathbb{C}^3.$
Consid\'erons l'application rationnelle $\mathrm{det}\hspace{0.1cm}\mathrm{jac}$ d\'efinie par
\begin{align*}
&\mathrm{det}\hspace{0.1cm}\mathrm{jac}\hspace{0.1cm}\colon\hspace{0.1cm}\mathrm{Rat}_2\to\mathbb{P}(\mathbb{C}[x_0,x_1,x_2]_3)
\simeq\{\text{ courbes de degr\'e $3$ }\}, &&
\hspace{0.1cm}[Q]\mapsto[\mathrm{det}\hspace{0.1cm}\mathrm{jac}\hspace{0.1cm} Q=0].
\end{align*}\label{not16b}

\begin{rem}\label{padef}
L'application $\mathrm{det}\hspace{0.1cm}\mathrm{jac}$ n'est pas d\'efinie pour les applications
$[Q]$ telles que $\mathrm{det}\hspace{0.1cm}\mathrm{jac}\hspace{0.1cm} Q\equiv~0;$ une telle application est
\`a conjugaison g.d. pr\`es du type 
\begin{align*}
&(Q_0:Q_1:0) &&\text{ou bien}&& (x_0^2:x_1^2:x_0x_1).
\end{align*}
En effet soit $Q$ une transformation
rationnelle quadratique telle que $\mathrm{det}\hspace{0.1cm}\mathrm{jac}\hspace{0.1cm} Q\equiv~0;$ notons $Q_i$ les composantes
de $Q.$

Supposons que $Q_0$ soit de rang maximal, par exemple
$Q_0=\frac{x_0^2+x_1^2+ x_2^2}{2};$ on peut alors se ramener \`a
\begin{align*}
& Q_1= ax_0x_2+bx_0x_1+q_1(x_1,x_2), && Q_2=cx_0x_1+q_2(x_1,x_2),
&&q_i\in\mathbb{C}[x_1,x_2]_2.
\end{align*}
Un calcul montre que le
coefficient de $x_0^3$ dans $\mathrm{det}\hspace{0.1cm}\mathrm{jac}\hspace{0.1cm} Q$ est $-ac;$ nous sommes
donc dans l'une des situations suivantes
\begin{align*}
& a\not=0,\hspace{0.1cm} c=0; && a=0,\hspace{0.1cm} c\not=0; && a=c=0.
\end{align*}
Si $a\not=0$ et $c=0$ l'\'egalit\'e $\mathrm{det}\hspace{0.1cm}\mathrm{jac}\hspace{0.1cm} Q\equiv 0$
entra\^ine $Q_2=0.$ Lorsque $a=0$ et $c\not=0$ on a,
\`a partir de $\mathrm{det}\hspace{0.1cm}\mathrm{jac}\hspace{0.1cm} Q\equiv 0,$ l'\'egalit\'e $bQ_2=Q_1.$
Enfin lorsque $a$ et $c$ sont nuls $\mathrm{det}\hspace{0.1cm}\mathrm{jac}\hspace{0.1cm} Q\equiv 0$ implique
que
\begin{itemize}
\item ou bien $q_2=0$ d'o\`u $Q_2=0;$

\item ou bien $b=0$ et $q_1$ et $q_2$ sont $\mathbb{C}$-colin\'eaires, par
suite $Q_1$ et $Q_2$ sont $\mathbb{C}$-colin\'eaires.
\end{itemize}

Le cas o\`u $Q_0$ n'est pas de rang maximum se traite de
la m\^eme fa\c{c}on.
\end{rem}

\begin{pro}\label{imdelta}
{\sl L'application $\mathrm{det}\hspace{0.1cm}\mathrm{jac}$ est surjective.}
\end{pro}

\begin{proof}[\sl D\'emonstration]
Le d\'eterminant jacobien de l'involution de \textsc{Cremona}
s'annule sur l'union de trois droites en position g\'en\'erale.

Soit $Q$ la transformation rationnelle d\'efinie par
$(x_0^2:x_1^2:(x_0-x_1)x_2).$ La courbe donn\'ee par $\mathrm{det}\hspace{0.1cm}\mathrm{jac}\hspace{0.1cm} Q=[x_0x_1(x_0
-x_1)=0]$ est l'union de trois droites concourantes.

Le d\'eterminant jacobien de $\rho$ (resp. $\tau$) est
l'union d'une droite double et d'une droite simple (resp. une
droite triple).

On remarque que
\begin{align*}
\mathrm{det}\hspace{0.1cm}\mathrm{jac}\hspace{0.1cm}\left(-\frac{1}{\alpha}x_0^2+x_2^2:-\frac{\alpha}{2}x_0x_2+
\frac{1+\alpha}{4}x_0^2-\frac{1}{4}x_1^2:x_0x_1\right)=[x_1^2x_2=x_0
(x_0-x_2)(x_0-\alpha x_2)]
\end{align*}
de sorte que l'on atteint toutes les cubiques ayant une forme normale de
\textsc{Weierstrass}.

Posons $Q:=(x_0x_1:x_0x_2:x_0^2+x_1x_2);$ alors $\mathrm{det}\hspace{0.1cm}\mathrm{jac}\hspace{0.1cm} Q
=[x_0(x_0^2-x_1x_2)=0]$ est l'union d'une conique et d'une droite
en position g\'en\'erale.

Ensuite pour la transformation $Q$ d\'efinie par
\begin{align*}
(x_1^2:x_0^2+2x_0x_2:x_0^2+x_0x_1+x_1x_2)
\end{align*}
on a $\mathrm{det}\hspace{0.1cm}\mathrm{jac}\hspace{0.1cm} Q=[x_1(2x_0^2-x_1x_2)=0]$ qui est l'union d'une conique et
d'une droite tangente \`a cette conique.

Consid\'erons $E$ l'espace vectoriel des matrices $3
\times 3$ dont les coefficients sont des formes lin\'eaires; $E$
est de dimension $27.$ Soit $F$ le sous-espace lin\'eaire de $E$
form\'e des matrices jacobiennes d'applications rationnelles quadratiques; on
remarque que $\dim F=18.$ Soit $\mathrm{det}$ le morphisme alg\'ebrique
d\'efini par
\begin{align*}
&\mathrm{det}\hspace{0.1cm}\colon\hspace{0.1cm} E\to\mathbb{C}[x_0,x_1,x_2]_3, && A\mapsto\mathrm{det}\hspace{0.1cm} A.
\end{align*} D'apr\`es ce qui pr\'ec\`ede l'image de $\mathrm{det}$
contient toutes les cubiques planes sauf \'eventuellement la
cubique cuspidale. Consid\'erons la matrice
\begin{align*}
A=\left[\begin{array}{ccc}
x_0 & 0 & x_1\\
x_2 & x_0 & 0\\
0 & x_2 & x_0
\end{array}
\right];
\end{align*}
ce n'est pas la matrice jacobienne d'une transformation
rationnelle quadratique mais 
\begin{align*}
\mathrm{det}\hspace{0.1cm} A=x_0^3+x_1x_2^2.
\end{align*}
L'application
$\mathrm{det}$ est donc surjective; en particulier ses fibres sont non
vides et donc de dimension au moins $27-10=17.$ Puisque $F$ est de
dimension $18$ et $E$ de dimension $27$ la fibre de $\mathrm{det}_{|F}$ au dessus de
toute cubique plane est de dimension strictement positive; en
particulier la cubique cuspidale est atteinte par $\mathrm{det}\hspace{0.1cm}\mathrm{jac}.$

Ainsi on obtient \`a conjugaison pr\`es toutes les
cubiques planes. Pour conclure il suffit de faire agir l'action gauche-droite.
\end{proof}

\section{Crit\`ere de birationnalit\'e}\label{critbir}\hspace{0.1cm}

Nous allons donner une pr\'esentation de la
classification des transformations birationnelles quadratiques qui
repose sur des arguments s'appuyant sur des remarques
\'el\'ementaires. Notons que la Proposition \ref{chaus} qui suit
n'est pas un \'enonc\'e original (\emph{voir} par exemple
\cite{AC}, Proposition~3.5.3.); le point de vue abord\'e ici est
toutefois diff\'erent du point de vue classique.

Soient $g$ une transformation rationnelle et $\psi$ un
polyn\^ome homog\`ene en trois variables. Rappelons que $g$
contracte $\psi$ si l'image par $g$ de la
courbe $[\psi=0]\setminus\mathrm{Ind}\hspace{0.1cm} g$ est un ensemble fini.

Soit $\widetilde{f}\hspace{0.1cm}\colon\hspace{0.1cm}\mathbb{C}^3\to\mathbb{C}^3$ une application
polynomiale homog\`ene de degr\'e $\nu;$ on note $\widetilde{f}_i$ les
composantes de $\widetilde{f}$ et on suppose que $\text{pgcd}
(\widetilde{f}_0, \widetilde{f}_1,\widetilde{f}_2)=1.$ Les
\textbf{\textit{points critiques}}\label{ind56} de
$\widetilde{f}$ sont les z\'eros de $\mathrm{det}\hspace{0.1cm}\mathrm{jac}\hspace{0.1cm}\widetilde{f};$ ils forment
l'ensemble $\mathcal{C}(\widetilde{f}).$ Soit $f=f^\bullet\hspace{0.1cm}\colon\hspace{0.1cm}\mathbb{P}^2(\mathbb{C})
\dashrightarrow\mathbb{P}^2(
\mathbb{C})$ l'application rationnelle induite par
$\widetilde{f};$ par hypoth\`ese $f$ est dominante. On d\'esigne
encore par
$\mathcal{C}(f)$\label{not16c} l'adh\'erence du lieu critique de la restriction
de $f$ \`a $\mathbb{P}^2(\mathbb{C}) \setminus\mathrm{Ind}\hspace{0.1cm} f;$ si $\nu\geq 2$ cet
ensemble est une courbe alg\'ebrique propre. Un calcul local
\'el\'ementaire montre que $\mathcal{C}(f)$ est exactement le
projectivis\'e de $\mathcal{C}(\widetilde{f}).$

\begin{rem}
En g\'en\'eral une transformation rationnelle ne contracte pas
$\mathrm{det}\hspace{0.1cm}\mathrm{jac}\hspace{0.1cm} f$ (c'est par exemple le cas pour $f=(x_0^2:x_1^2:x_2^2)$).
\end{rem}

\begin{pro}\label{chaus}
{\sl Soit $f$ une transformation birationnelle de $\mathbb{P}^2(\mathbb{C})$ dans
lui-m\^eme; $\mathrm{det}\hspace{0.1cm}\mathrm{jac}\hspace{0.1cm} f$ est contract\'e par $f.$}
\end{pro}

\begin{proof}[\sl D\'emonstration]
Nous allons raisonner par l'absurde. Soit $\Gamma$ une composante
irr\'eductible de~$\mathcal{C}(f).$ Supposons que $\Gamma$ ne soit
pas contract\'ee par $f;$ alors $\overline{f(\Gamma\setminus
\mathrm{Ind}\hspace{0.1cm} f)}$ est une courbe $\Gamma'.$ On se donne un point
g\'en\'erique $m$ de $\Gamma;$ en $m,$ la transformation $f$ est
holomorphe, $\Gamma$ est lisse et~$f_{|\Gamma}$ est une submersion
de $\Gamma_{,m}$ sur $\Gamma'_{,m'}$ o\`u $m'=f(m).$ \`A
conjugaison par des diff\'eomorphismes locaux \`a la source et au
but pr\`es, on peut supposer que $m=m'=0$ et que $\Gamma_{,m}
=\Gamma'_{,m'}=(x_1=0).$ Le th\'eor\`eme de la fonction inverse
assure qu'on peut se ramener au voisinage de $0$ \`a 
\begin{align*}
f(x_0,0)=(x_0,0);
\end{align*}
dit autrement il existe $g_1,$ $g_2$ deux fonctions
holomorphes telles que
\begin{align*}
f(x_0,x_1)=(x_0+x_1g_1(x_0,x_1),x_1g_2(x_0,x_1)).
\end{align*}
Toujours par inversion locale on peut supposer que $f(x_0,x_1)=(x_0,x_1
\psi(x_0,x_1)).$ Puisque
\begin{align*}
\mathrm{det}(\mathrm{jac}\hspace{0.1cm} f)=\psi(x_0,x_1)+x_1 \frac{\partial\psi}
{\partial x_1}(x_0,x_1)
\end{align*}
doit s'annuler sur $x_1=0$ on a $\psi(x_0,0)=0.$ Il
en r\'esulte que $f$ s'\'ecrit
\begin{align*}
&(x_0,x_1^k\zeta(x_0,x_1)), && \text{ avec }k\geq 2,
\hspace{3mm}\zeta(x_0,0) \not\equiv 0.
\end{align*}
Quitte \`a se placer en un point
voisin de $0$ sur $x_1=0,$ on peut supposer que $\zeta(0,0)$ est non
nul de sorte qu'\`a conjugaison locale \`a la source et au but
pr\`es $f$ est de la forme $(x_0, x_1^k).$ Or une application
birationnelle est g\'en\'eriquement localement inversible ce qui
n'est pas le cas pour~$(x_0,x_1^k)$ lorsque $k\geq 2.$
\end{proof}

Si $A,$ $B$ sont deux \'el\'ements de $\mathrm{PGL}_3(\mathbb{C})$ et $Q=A\sigma
B$ (resp. $Q=A\rho B,$ resp. $Q=A\tau B$), alors $\mathrm{det}\hspace{0.1cm}\mathrm{jac}\hspace{0.1cm} Q$ est
l'union de trois droites en position g\'en\'erale (resp. l'union
d'une droite double et d'une droite simple, resp. une droite
triple).
Nous allons donner un crit\`ere qui permet de d\'eterminer si une transformation
quadratique rationnelle est birationnelle.

\begin{thm}\label{cri}
{\sl Soit $Q$ une transformation rationnelle purement quadratique
non d\'eg\'en\'er\'ee\label{ind290} ({\it i.e.} $\mathrm{det}\hspace{0.1cm}\mathrm{jac}\hspace{0.1cm} Q\not\equiv 0$).
Supposons que $Q$ contracte $\mathrm{det}\hspace{0.1cm}\mathrm{jac}\hspace{0.1cm} Q;$ alors $\mathrm{det}\hspace{0.1cm}\mathrm{jac}\hspace{0.1cm} Q$ est
l'union de trois droites non concourantes et $Q$ est
birationnelle.

De plus
\begin{itemize}
\item si $\mathrm{det}\hspace{0.1cm}\mathrm{jac}\hspace{0.1cm} Q$ est l'union de trois
droites en position g\'en\'erale, $Q$ est, \`a \'equivalence g.d.
pr\`es, l'involution de \textsc{Cremona} $\sigma;$

\item si $\mathrm{det}\hspace{0.1cm}\mathrm{jac}\hspace{0.1cm} Q$ est l'union d'une droite double et d'une
droite simple, $Q$ co\"incide, \`a conjugaison g.d. pr\`es, avec
$\rho;$

\item enfin si $\mathrm{det}\hspace{0.1cm}\mathrm{jac}\hspace{0.1cm} Q$ est une droite triple, $Q$ appartient
\`a $\mathcal{O}_{g.d.} (\tau).$
\end{itemize}}
\end{thm}

\begin{cor}
{\sl Une transformation rationnelle quadratique de $\mathbb{P}^2(\mathbb{C})$ dans lui-m\^eme
appartient \`a $\Sigma^3$ si et seulement si elle admet trois
points d'ind\'etermination.}
\end{cor}

\begin{rem}
Une transformation birationnelle $Q$ de $\mathbb{P}^2(
\mathbb{C})$ dans lui-m\^eme
contracte $\mathrm{det}\hspace{0.1cm}\mathrm{jac}\hspace{0.1cm} Q$ et ne contracte pas d'autre courbe. On
peut se demander si le crit\`ere de birationnalit\'e du
Th\'eor\`eme \ref{cri} est encore valable en degr\'e strictement
sup\'erieur \`a $2;$ la r\'eponse est non. On a d\`es le degr\'e
$3$ des exemples d'\'el\'ements $Q$ qui contractent $\mathrm{det}\hspace{0.1cm}\mathrm{jac}\hspace{0.1cm} Q$
mais ne sont pas birationnels comme l'atteste la transformation
$(x_0^2x_1:x_0x_2^2:x_1^2x_2).$ Remarquons que si ce crit\`ere de
birationnalit\'e \'etait valable en tout degr\'e, il impliquerait
la conjecture du jacobien de fa\c{c}on plus ou moins triviale.
\end{rem}

\begin{proof}[\sl D\'emonstration du Th\'eor\`eme \ref{cri}]
Montrons que $\mathrm{det}\hspace{0.1cm}\mathrm{jac}\hspace{0.1cm} Q$ est l'union de trois droites.

Supposons $\mathrm{det}\hspace{0.1cm}\mathrm{jac}\hspace{0.1cm} Q$
irr\'eductible. Posons $Q:=(Q_0:Q_1:Q_2).$ \`A composition \`a
gauche pr\`es on peut supposer que $\mathrm{det}\hspace{0.1cm}\mathrm{jac}\hspace{0.1cm} Q$ est contract\'e
sur $(1:0:0);$ alors $\mathrm{det}(\mathrm{jac}\hspace{0.1cm} Q)$ divise $Q_2$ et $Q_3:$
impossible.

De la m\^eme fa\c{c}on si $\det\hspace{0.1cm}\mathrm{jac}\hspace{0.1cm}
Q=Lq$ avec $L$ est lin\'eaire et $q$ quadratique non d\'eg\'en\'er\'ee, on peut se ramener \`a: $q=0$ est contract\'e sur
$(1:0:0);$ alors $Q$ s'\'ecrit $(q_1:q:\alpha q)$ et donc est
d\'eg\'en\'er\'ee.

Ainsi $\det(\mathrm{jac}\hspace{1mm}Q)$ est le produit de trois formes lin\'eaires.

\bigskip

\'Etudions le cas o\`u, \`a conjugaison pr\`es, $\mathrm{det}\hspace{0.1cm}\mathrm{jac}\hspace{0.1cm} Q=x_0x_1x_2.$ Si les droites $x_0=0$ et $x_1=0$
sont contract\'ees sur un m\^eme point, disons $(1:0:0),$ alors $Q$ est
de la forme $(q:x_0x_1:\alpha x_0x_1)$ qui est d\'eg\'en\'er\'ee.
Les droites $x_0=0,$ $x_1=0$ et $x_2=0$ sont donc contract\'ees
sur trois points distincts. Un calcul direct montre qu'ils ne peuvent \^etre
align\'es. On se ram\`ene alors au cas o\`u~$x_0=0$ (resp. $x_1=0,$
resp. $x_2=0$) est contract\'ee sur
$(1:0:0)$ (resp. $(0:1:0),$ resp. $(0:0:1)$); on constate alors que $Q$ est
l'involution de \textsc{Cremona} \`a conjugaison pr\`es.

\bigskip

Consid\'erons l'\'eventualit\'e
o\`u $\mathrm{det}\hspace{0.1cm}\mathrm{jac}\hspace{0.1cm} Q$ a deux branches $x_0=0$ et $x_2=0.$
D'apr\`es ce que l'on vient de voir, les droites $x_0=0$
et $x_2=0$ sont contract\'ees sur deux points distincts, par exemple~$(1:0:0)$ et $(0:1:0).$ La transformation $Q$ s'\'ecrit donc, \`a
\'equivalence g.d. pr\`es,
\begin{align*}
& Q=(x_2(\alpha x_1+\beta x_2):x_0(\gamma x_0+
\delta x_1):x_0x_2).
\end{align*}
Pour que $\mathrm{det}\hspace{0.1cm}\mathrm{jac}\hspace{0.1cm} Q$ ait deux branches il faut que $\alpha\delta=0.$
Un calcul direct montre que $Q$ est birationnelle d\`es que $\beta\delta
-\alpha\gamma\not=0$ et en fait g.d. \'equivalente \`a $\rho.$

\bigskip

On s'int\'eresse au cas o\`u $\mathrm{det}(\mathrm{jac}\hspace{0.1cm} Q)=x_2^3.$
On peut supposer que $x_2=0$ est contract\'ee sur $(1:0:0);$ alors $Q=(q: x_2\ell_1:x_2\ell_2)$ o\`u $q$ d\'esigne une
forme quadratique et les $\ell_i$ des formes lin\'eaires.

\begin{itemize}
\item Si $(x_2,\ell_1,\ell_2)$ est un
syst\`eme de coordonn\'ees, on peut \'ecrire \`a conjugaison pr\`es
\begin{align*}
&Q=(q:x_0x_2:x_1x_2), &&q=ax_0^2+bx_1^2+cx_2^2+dx_0x_1.
\end{align*}
Le calcul explicite de
$\mathrm{det}(\mathrm{jac}\hspace{0.1cm} Q)s$ conduit \`a $a=b=d=0,$ {\it i.e.} ou
bien $Q$ est d\'eg\'en\'er\'ee, ou bien $Q$ repr\'esente une transformation
lin\'eaire, toutes choses exclues.

\item Supposons que $(x_2,\ell_1,\ell_2)$
ne soit pas un syst\`eme de coordonn\'ees, {\it i.e.}
\begin{align*}
& \ell_1=ax_2+\ell(x_0,x_1), &&\ell_2=bx_2+\varepsilon\ell(x_0,x_1).
\end{align*}
Notons que $\ell$ n'\'etant pas identiquement nulle (sinon $Q$ serait
d\'eg\'en\'er\'ee), on peut se
ramener \`a $\ell=x_0$. \`A \'equivalence g.d. pr\`es $Q$ est de la
forme $(q:x_0x_2:x_2^2).$ Un calcul explicite conduit \`a $\mathrm{det}
(\mathrm{jac}\hspace{0.1cm}Q)=-2x_2^2\frac{\partial q} {\partial x_1};$ par suite $\frac{\partial
q}{\partial x_1}$ est divisible par $x_2$ autrement dit $q$ est du
type $\alpha x_2^2+\beta x_0x_2+\gamma x_0^2+\delta x_1x_2.$ Toujours \`a
\'equivalence g.d. pr\`es, on obtient $Q=\tau.$\bigskip
\end{itemize}

Pour finir abordons la possibilit\'e: $\mathrm{det}(\mathrm{jac}\hspace{0.1cm} Q)=x_0x_1(x_0-x_1).$
Comme on l'a vu les droites $x_0=0$ et $x_1=0$ sont
contract\'ees sur deux points distincts que l'on va supposer
\^etre $(1:0:0)$ et $(0:1:0).$ Il s'en suit que $Q$ est de la
forme
\begin{align*}
&(x_1(ax_0+bx_1+cx_2): x_0(\alpha x_0+\beta x_1+\gamma x_2):x_0x_1),
&& a,\hspace{0.1cm} b,\hspace{0.1cm} c,\hspace{0.1cm} \alpha,\hspace{0.1cm}\beta,\hspace{0.1cm}\gamma\in\mathbb{C}.
\end{align*}
On constate que l'image de la droite $x_0=x_1$ par $Q$ est
\begin{align*}
((a+b)x_0+cx_2:(\alpha+\beta)x_0+\gamma x_2:x_0);
\end{align*}
c'est un point si et
seulement si $c$ et $\gamma$ sont nuls auquel cas $Q$ ne d\'epend
pas de $x_2.$
\end{proof}

Posons
\begin{align*}
& \Sigma^3:=\mathcal{O}_{g.d.}(\sigma),&&
\Sigma^2:=\mathcal{O}_{g.d.}(\rho),&&
\Sigma^1:=\mathcal{O}_{g.d.}(\tau).
\end{align*}\label{not17}\label{not18}\label{not19}

Consid\'erons une transformation birationnelle repr\'esent\'ee par
$Q=\ell(\ell_0:\ell_1:\ell_2)$ o\`u $\ell$ et les~$\ell_i$
d\'esignent des formes lin\'eaires, les $\ell_i$ \'etant
ind\'ependantes. La droite $\ell=0$ est une \og droite
contract\'ee apparente\fg\hspace{0.1cm}; en effet l'action de $Q$ sur
$\mathbb{P}^2(\mathbb{C})$ est \'evidemment celle de l'automorphisme $(\ell_0:
\ell_1:\ell_2)$ de $\mathbb{P}^2(\mathbb{C}).$ On notera
$\Sigma^0\label{not20}$ l'ensemble de ces transformations
\begin{align*}
\Sigma^0=\{\ell(\ell_0:\ell_1:\ell_2)\hspace{0.1cm}|\hspace{0.1cm}\ell,\hspace{0.1cm}\ell_i \text{
formes lin\'eaires, les $\ell_i$ \'etant ind\'ependantes}\}.
\end{align*}

Les \'el\'ements de $\Sigma^0$ seront abusivement
appel\'es lin\'eaires; en fait c'est l'ensemble
\begin{align*}
(\Sigma^0)^\bullet=\{f^\bullet\hspace{0.1cm}|\hspace{0.1cm} f\in\Sigma^0\}
\end{align*}
qui s'identifie \`a $\mathrm{PGL}_3(\mathbb{C}).$ On a
$\Sigma^0=\mathcal{O}_{g.d.}(x_0 (x_0:x_1:x_2)):$ \`a conjugaison
gauche-droite pr\`es un \'el\'ement de la forme $\ell A$ s'\'ecrit
$x_0A'$ puis $x_0\mathrm{id}.$ Cette approche permet de voir les
d\'eg\'en\'erescences de transformations quadratiques sur des
applications lin\'eaires.

Un automorphisme polynomial de $\mathbb{C}^2$ est
appel\'e \textbf{\textit{automorphisme de
\textsc{H\'enon}}}\label{ind29} s'il s'\'ecrit (\`a conjugaison
affine pr\`es)
\begin{align*}
&(x_1,P(x_1)-\delta x_0), &&P\in\mathbb{C}[x_1],\hspace{0.1cm} \deg P\geq
2,\hspace{0.1cm}\delta\in\mathbb{C}^*.
\end{align*}
D'apr\`es le Th\'eor\`eme \ref{cri} le prolongement \`a $\mathbb{P}^2(\mathbb{C})$
de tout automorphisme de \textsc{H\'enon} quadratique appartient
\`a $\Sigma^1.$

Remarquons qu'un \'el\'ement de $\Sigma^i$ compte $i$
points d'ind\'etermination et $i$ courbes contract\'ees.

Un \'el\'ement de $\Sigma^i$ ne peut \^etre
lin\'eairement conjugu\'e \`a un \'el\'ement de $\Sigma^j$ o\`u
$j\not=i;$ par contre ils peuvent l'\^etre birationnellement: on a
par exemple signal\'e pr\'ec\'edemment que les involutions
$\sigma,$ $\rho$ et $\tau$ sont dynamiquement conjugu\'ees \`a des
involutions de $\mathrm{PGL}_3(\mathbb{C}).$ Toutefois des arguments que nous
d\'evelopperons au Chapitre \ref{feuilfeuil} montrent qu'un
\'el\'ement g\'en\'erique de $\Sigma^i,$ pour $i\geq 1,$ n'est pas
birationnellement conjugu\'e \`a une transformation lin\'eaire.

\begin{cor}
{\sl On a
\begin{align*}
& \mathring{\mathrm{B}}\mathrm{ir}_2=\Sigma^1\cup\Sigma^2\cup\Sigma^3, &&
\mathrm{Bir}_2=\Sigma^0\cup\Sigma^1\cup\Sigma^2\cup \Sigma^3.
\end{align*}}
\end{cor}

\begin{rems}
\begin{itemize}
\item Une d\'ecomposition de \textsc{N\oe ther} de
$\rho$ est
\begin{align*}
(x_2-x_1:x_1-x_0:x_1)\sigma(x_1+x_2:x_2:x_0)\sigma(x_0+x_2:x_1-x_2:x_2)
\end{align*}
qu'il faudrait en toute l\'egitimit\'e \'ecrire
\begin{align*}
\rho^\bullet=((x_2-x_1:x_1-x_0:x_1)\sigma(x_1+x_2:x_2:x_0)
\sigma(x_0+x_2:x_1-x_2:x_2))^\bullet.
\end{align*}

On retrouve donc le fait classique suivant d\'ej\`a
\'enonc\'e dans \cite{Hu, AC}: pour toute transformation
birationnelle quadratique $Q$ ayant deux points
d'ind\'etermination il existe $\ell_1,$ $\ell_2$ et $\ell_3$ dans
$\mathrm{PGL}_3(\mathbb{C})$ tels que $Q=\ell_1\sigma\ell_2 \sigma\ell_3.$ En
particulier la d\'etermination des transformations lin\'eaires
$\ell$ telle que $\sigma\ell\sigma$ soit quadratique m\'erite
attention (\emph{voir} Chapitre \ref{dyndyn}).

\item La transformation
$\tau=(x_0^2:x_0x_1:x_1^2-x_0x_2)$, qui s'\'ecrit aussi
$\ell_1\sigma\ell_2\sigma \ell_3\sigma\ell_4\sigma\ell_5$
avec
\begin{align*}
& \ell_1=(x_1-x_0:2x_1-x_0:x_2-x_1+x_0), && \ell_2=(x_0+x_2:x_0:x_1), \\
& \ell_3=(-x_1:x_0+x_2-3x_1:x_0), && \ell_4=(x_0+x_2:x_0:x_1),\\
& \ell_5= (x_1-x_0:-2x_0+x_2:2x_0-x_1), &&
\end{align*}
est dans $\Sigma^1.$ Il
s'en suit que tout \'el\'ement de $\Sigma^1$ est de la forme
$\ell_1\sigma\ell_2\sigma\ell_3 \sigma \ell_4 \sigma\ell_5$ avec
$\ell_i$ dans $\mathrm{PGL}_3(\mathbb{C})$ (\emph{voir} \cite{Hu, AC}). La
r\'eciproque est fausse au sens o\`u tout \'el\'ement du type
\begin{align*}
\ell_1\sigma\ell_2 \sigma\ell_3 \sigma \ell_4 \sigma\ell_5
\end{align*}
n'appartient pas \`a $\Sigma^1$; g\'en\'eriquement (sur le choix
des $\ell_i$) un tel \'el\'ement est de degr\'e $16.$
\end{itemize}
\end{rems}

\section[Relations, second crit\`ere de birationnalit\'e]{Relations
et transformations birationnelles, second crit\`ere de
birationnalit\'e}\label{seccritbir}\hspace{0.1cm}

Soit $Q=(Q_0:Q_1:Q_2)$ une transformation rationnelle de
$\mathbb{P}^2(\mathbb{C});$ on appelle \textbf{\textit{relation
lin\'eaire}}\label{ind30} de $Q$ tout triplet
$L=(L_0,L_1,L_2)$ de formes lin\'eaires satisfaisant
\begin{align*}
L_0Q_0+L_1Q_1+L_2Q_2=0.
\end{align*}\label{not21}

On note $\mathrm{RL}(Q)\label{not22}$ le
$\mathbb{C}$-espace vectoriel des relations lin\'eaires de $Q.$ On observe
que la dimension $\mathrm{e}(Q)\label{not23}$ de
$\mathrm{RL}(Q)$ est invariante sous l'action gauche-droite.
Rappelons le r\'esultat classique suivant.

\begin{pro}[\cite{To}]
{\sl Soient $f_1,\ldots,f_n$ des germes de fonctions holomorphes
\`a l'origine de $\mathbb{C}^n$ tels que $\{f_1=\ldots=f_n=0\}=\{0\}.$
Alors le $\mathcal{O}(\mathbb{C}^n,0)$-module des relations
\begin{align*}
R(f):=\{a=(a_1,\ldots,a_n)\in\mathcal{O}(\mathbb{C}^n,0)\hspace{0.1cm}|\hspace{0.1cm}a_1f_1+\ldots+a_nf_n=0\}
\end{align*}
est engendr\'e par les relations triviales $f_{ij}=(0,\ldots,0,f_j,0,
\ldots,0,-f_i,0,\ldots, 0).$}
\end{pro}

Une cons\'equence de cette proposition est la suivante.

\begin{pro}\label{zero}
{\sl Soit $f=(f_0:f_1:f_2)$ un endomorphisme de $\mathbb{P}^2(\mathbb{C}),$ 
c'est-\`a-dire $f$ satisfait $\{f_0=f_1=f_2=0\}=\{0\}.$ Si $\deg f\geq 2$ alors
$\mathrm{RL}(f)=\{0\}.$}
\end{pro}

D\'esormais nous consid\'erons plus particuli\`erement
les transformations quadratiques pures.

\subsection{Calculs explicites d'espaces de relations lin\'eaires}\hspace{0.1cm}

Nous allons calculer $\mathrm{RL}(Q)$ pour quelques $Q$
typiques.\bigskip

Nous commen\c{c}ons par l'involution de \textsc{Cremona}
$\sigma.$ Visiblement $\mathrm{RL}(\sigma)$ est engendr\'e par les
relations $(-x_0,x_1,0)$ et $(x_0,0,-x_2).$ On a donc
$\mathrm{e}(\sigma) =2;$ on remarque que $\mathrm{RL}(\sigma)$
s'identifie naturellement (on repr\'esente les $L$ par des
matrices $3 \times 3$ \`a coefficient dans $\mathbb{C}$) \`a l'alg\`ebre
de \textsc{Lie} commutative des matrices diagonales de trace
nulle. Notons que si l'on tord $\sigma$ par conjugaison diagonale,
{\it i.e.} on consid\`ere
$\sigma_\alpha=(\alpha_0x_1x_2:\alpha_1x_0x_2:\alpha_2 x_0x_1),$
alors $\mathrm{RL}(\sigma_\alpha)$ s'identifie \`a l'alg\`ebre des
matrices diagonales
\begin{align*}
\left[\begin{array}{ccc}
\beta_0 & 0 & 0 \\
0 & \beta_1 & 0 \\
0 & 0 & \beta_2
\end{array}\right]
\end{align*}
telles que $\alpha_0\beta_0+\alpha_1\beta_1+\alpha_2\beta_2=0.$
L'invariance g.d. de l'entier $\mathrm{e}$ nous permet d'affirmer que pour
tout \'el\'ement $Q$ dans $\Sigma^3$ on a $\mathrm{e}(Q)=2.$ \bigskip

Nous consid\'erons maintenant la transformation $\rho'$
de $\Sigma^2$ d\'efinie par 
\begin{align*}
\rho'=(x_1x_2:x_1x_2-x_0x_2:~-x_1^2).
\end{align*}
Un calcul \'el\'ementaire montre que
\begin{align*}
\mathrm{RL}(\rho')=\{(\alpha x_0+\beta
x_1,\alpha x_1,(\alpha+\beta)x_2)\hspace{0.1cm}|\hspace{0.1cm}(\alpha,\beta)\in
\mathbb{C}^2\}.
\end{align*}
Ainsi $\mathrm{e}(\rho')=2$ et $\mathrm{RL}(\rho')$ s'identifie
\`a l'alg\`ebre de \textsc{Lie} commutative de dimension $2$ des
matrices de \textsc{Jordan} du type
\begin{align*}
\left[\begin{array}{ccc}
\alpha & \beta & 0\\
0 & \alpha & 0 \\
0 & 0 & \alpha+\beta
\end{array}
\right].
\end{align*}
Comme pr\'ec\'edemment pour tout \'el\'ement $Q$ de
$\Sigma^2$ on a $\mathrm{e}(Q)=2.$ \bigskip

Nous choisissons dans $\Sigma^1$ l'\'el\'ement $\tau'$
d\'efini par $\tau'=(x_2^2:-x_1x_2:x_1^2-x_0x_2).$ On v\'erifie que
\begin{align*}
\mathrm{RL}(\tau')=\{(\alpha x_0+\beta x_1,\alpha x_1+\beta x_2,\alpha x_2)\hspace{0.1cm}
|\hspace{0.1cm}(\alpha,\beta)\in\mathbb{C}^2\}.
\end{align*}
Ici encore $\mathrm{e}(\tau')=2$ et $\mathrm{RL}(\tau')$
s'identifie \`a l'alg\`ebre de \textsc{Lie} commutative
constitu\'ee des matrices de \textsc{Jordan} de la forme
\begin{align*}
\left[\begin{array}{ccc}
\alpha & \beta & 0\\
0 & \alpha & \beta \\
0 & 0 & \alpha
\end{array}
\right].
\end{align*}

Cette correspondance entre les orbites de $\Sigma^3,$
$\Sigma^2$ et $\Sigma^1$ et les diff\'erents types d'alg\`ebres de
\textsc{Lie} ab\'eliennes de dimension $2$ des matrices $3\times
3$ \`a coefficients dans $\mathbb{C}$ s'av\`ere pour le moins curieuse,
mais peut-\^etre seulement anecdotique.\bigskip

Passons maintenant au calcul de l'invariant $\mathrm{e}$
pour une transformation de $\Sigma^0.$ Consid\'erons par exemple
l'application identit\'e repr\'esent\'ee dans $\mathrm{Bir}_2$ par
$\mathrm{id}_\ell=(\ell x_0:\ell x_1:\ell x_2)$ o\`u $\ell$ est
une forme lin\'eaire non triviale. On constate que $\mathrm{RL}(
\mathrm{id}_\ell)$ est ind\'ependant du choix de $\ell$ et est
engendr\'e par $(x_1,-x_0,0),$ $(x_2,0,-x_0)$ et $(0,x_2,-x_1).$
En r\'esulte que $\mathrm{e}(\mathrm{id}_\ell)=3$ et $\mathrm{RL}
(\mathrm{id}_\ell)$ s'identifie \`a l'alg\`ebre de \textsc{Lie}
simple de dimension $3$ des matrices $3\times 3$
antisym\'etriques. \'Evidemment $\mathrm{e}(Q)$ vaut $3$ pour tout
\'el\'ement $Q$ de $\Sigma^0.$
\bigskip

Nous avons \'etabli la:

\begin{pro}
{\sl Soit $Q$ une transformation birationnelle quadratique. Alors
\begin{itemize}
\item $Q$ appartient \`a $\Sigma^0$ si et seulement si
$\mathrm{e}(Q)=3;$

\item $Q$ est purement quadratique si et seulement
$\mathrm{e}(Q)=2.$
\end{itemize}}
\end{pro}

Comme on l'a vu il y a dans l'adh\'erence de $\mathrm{Bir}_2$
dans $\mathbb{P}^{17}(\mathbb{C})$ des \'el\'ements d\'eg\'en\'er\'es $Q$ pour
lesquels les espaces $\mathrm{RL}(Q)$ se calculent ais\'ement. En
voici quelques~uns
\begin{align*}
\mathrm{RL}(x_0^2:x_1^2:0)=\{(0,0,\ell)\hspace{0.1cm}|\hspace{0.1cm}\ell \text{ forme
lin\'eaire}\}
\end{align*}
est de dimension $3$ et c'est encore une alg\`ebre de
\textsc{Lie};
\begin{align*}
\mathrm{RL}(x_0^2:x_1^2:x_0x_1)=\{( \alpha x_1,\beta
x_0,-\alpha x_0-\beta x_1)\hspace{0.1cm}|\hspace{0.1cm}(\alpha,\beta)\in\mathbb{C}^2\}
\end{align*}
est de dimension $2.$
\bigskip

\subsection{Crit\`ere de birationnalit\'e}\hspace{0.1cm}

Remarquons que si $Q$ est birationnelle purement
quadratique, $\mathrm{RL}(Q)$ contient un \'el\'ement $L$ qui, vu
comme endomorphisme de $\mathbb{C}^3,$ est inversible. Ce n'est plus le
cas pour les transformations $Q$ de $\Sigma^0;$ les \'el\'ements
de $\mathrm{RL}(Q)$ sont alors tous non inversibles.

\begin{lem}\label{fatigue}
{\sl Soit $Q$ un \'el\'ement de $\mathrm{Rat}_2$ tel que
$\mathrm{e}(Q)=2.$ Si $\mathrm{RL}(Q)$ contient un \'el\'ement
inversible, $Q$ appartient \`a
$\Sigma^1\cup\Sigma^2\cup\Sigma^3.$}
\end{lem}

\begin{proof}[\sl D\'emonstration]
Soit $\{L=(L_0,L_1,L_2),\hspace{0.1cm} L'=(L'_0,L'_1,L'_2)\}$ une base de
$\mathrm{RL}(Q);$ on peut supposer que $L$ et $L'$ sont
inversibles. Consid\'erons la transformation quadratique $f=(f_0:
f_1:f_2),$ dite d\'eterminantielle, d\'efinie par $f=L\wedge L'$
\begin{align*}
\left[\begin{array}{ccc} f_0\\ f_1\\
f_2\end{array}\right]=\left[\begin{array}{ccc}
L_1L'_2-L_2 L'_1\\ L_2L'_0-L_0L'_2\\
L_0L'_1-L_1L'_0\end{array}\right].
\end{align*}

Commen\c{c}ons par quelques remarques
\begin{itemize}
\item[\textbf{\textit{i.}}] Supposons que $f$ soit identiquement
nulle. Comme $L$ est inversible, $L'$ est un multiple de~$L$
ce qui contredit le fait que
$\mathrm{RL}(Q)$, de dimension $2,$ soit engendr\'e par $L$ et
$L'.$

\item[\textbf{\textit{ii.}}] Si $Q$ n'est pas purement
quadratique, $\mathrm{e} (Q)$ est minor\'e par $3.$ Ceci se
constate au cas par cas en utilisant le fait que $Q$ est \`a
\'equivalence g.d. pr\`es de l'un des types suivants
\begin{itemize}
\item $(\ell x_0,\ell x_1,\ell x_2),$ $\ell$ d\'esignant une forme
lin\'eaire;

\item $(\ell x_0,\ell x_1,0),$ $\ell$ \'etant une forme lin\'eaire;

\item $(q,0,0)$ avec $q$ forme quadratique.
\end{itemize}
\end{itemize}

D'apr\`es \textbf{\textit{i.}} en un point g\'en\'erique
$m$ les vecteurs $L(m)$ et $L'(m)$ sont ind\'ependants, {\it
i.e.} $f$ est non identiquement nulle et satisfait
$f \wedge Q=0.$ D'apr\`es \textbf{\textit{ii.}} on a $\mathrm{codim}
\hspace{0.1cm} Q^{-1}(0)\geq 2;$ la Proposition
\ref{zero} assure que $Q^{-1}(0)$ ne peut se r\'eduire \`a $\{0\}$
de sorte que $Q^{-1}(0)$ est constitu\'e d'un nombre fini de
droites. En particulier $Q$ satisfait la propri\'et\'e de division
de \textsc{de Rham}-\textsc{Saito}~(\cite{Sa}), {\it i.e.} il existe $h$
dans $\mathcal{O}(\mathbb{C}^3)$ tel que $f$ s'\'ecrive $hQ.$ Pour des
raisons d'homog\'en\'eit\'e~$h$ est une constante non nulle et on
peut supposer que $Q=L\wedge L'.$ Les z\'eros de $Q$ (qui
correspondent en coordonn\'ees homog\`enes aux points
d'ind\'etermination) sont exactement les points $m$ o\`u $L$ et
$L'$ sont colin\'eaires; dit autrement ce sont les directions
propres de $L^{-1}L'.$ Il y a \'evidemment trois possibilit\'es
\begin{itemize}
\item $L^{-1}L'$ a trois vecteurs propres non colin\'eaires;

\item $L^{-1}L'$ compte deux vecteurs propres non colin\'eaires;

\item $L^{-1}L'$ poss\`ede un unique vecteur propre.
\end{itemize}

Dans la premi\`ere \'eventualit\'e $Q$
a trois points d'ind\'etermination et est donc birationnelle, en
fait dans $\Sigma^3.$ Dans les autres cas pour \'etablir la birationnalit\'e on peut supposer \`a
\'equivalence g.d. pr\`es que $L=\mathrm{id};$ il suffit alors d'examiner les
diff\'erents types de \textsc{Jordan} pour~$L'.$ Ainsi dans le second cas
il y a deux points d'ind\'etermination, {\it i.e.} $Q$ appartient \`a $\Sigma^2.$
Enfin dans la derni\`ere possibilit\'e, $Q$ compte un unique point
d'ind\'etermination: $Q$ est un \'el\'ement de $\Sigma^1.$
\end{proof}

Inversement on a le:

\begin{lem}
{\sl Soit $Q$ un \'el\'ement de $\mathbb{P}^{17}(\mathbb{C})$ repr\'esentant une transformation
rationnelle tel que $\mathrm{e}(Q)=~2.$ Supposons que les
\'el\'ements de $\mathrm{RL}(Q)$ soient non inversibles. Alors $Q$ co\"incide,
\`a conjugaison gauche-droite pr\`es, avec $(x_0^2:x_1^2:x_0x_1);$ en
particulier $Q$ appartient \`a $\overline{\mathrm{Bir}_2}.$}
\end{lem}

\begin{proof}[\sl D\'emonstration]
On proc\`ede comme dans le lemme pr\'ec\'edent. La transformation
$Q$ est n\'ecessairement purement quadratique en particulier
$\dim Q^{-1}(0)\leq 1.$ Soit $\{L,L'\}$ une
base de~$\mathrm{RL}(Q).$

Montrons que $f:=L\wedge L'$ ne peut \^etre
identiquement nul. Si $f$ \'etait nul, l'argument de division de
\textsc{de Rham}-\textsc{Saito} (\cite{Sa}) imposerait que $L$ et $L'$ s'annulent
sur un hyperplan de sorte que
\begin{align*}
&L=\ell(a,b,c),&&L'= \ell'(a',b',c')
\end{align*}
o\`u $\ell,$ $\ell'$ d\'esignent des formes lin\'eaires et
$(a,b,c),$ $(a',b',c')$ des \'el\'ements de $\mathbb{C}^3;$ ceci
impliquerait que $Q$ n'est pas purement quadratique. Ainsi comme
dans le Lemme \ref{fatigue} on a $L\wedge L'\not\equiv 0$ et on peut
supposer que $Q$ s'\'ecrit $L\wedge L'.$ Notons que chaque
\'el\'ement de $\mathrm{RL}(Q)$ est, par l'argument pr\'ec\'edent,
pr\'ecis\'ement de rang $2.$ En particulier $\ker(L-tL')$ est,
pour tout $t$ dans $\mathbb{C},$ une droite~$\mathcal{D}_t$ d\'ependant
holomorphiquement de $t.$ On observe que chaque $\mathcal{D}_t$
appartient \`a $Q^{-1}( 0);$ il s'en suit que $\mathcal{D}_t$ est
constant et $L,$ $L'$ ont leur noyau commun, disons $\mathbb{C}(0:0:1).$
Par cons\'equent $Q$ est du type $(Q_0(x_0,x_1):Q_1(x_0,x_1):Q_2(x_0,x_1)).$
Si les $Q_i$ forment une base des formes quadratiques en deux
variables, $Q$ est g.d. conjugu\'e \`a $(x_0^2:x_1^2:x_0x_1)$ comme
annonc\'e. Sinon $Q$ est, toujours \`a conjugaison g.d. pr\`es, du
type $(Q_0:Q_1:0);$ mais pour une telle transformation la
dimension de $\mathrm{RL}(Q)$ est sup\'erieure ou \'egal \`a $3$
ce qui, par hypoth\`ese, est impossible.
\end{proof}

\begin{lem}
{\sl Soit $Q$ une transformation rationnelle quadratique. Si
$\mathrm{e}(Q)\geq 3,$ alors $Q$ appartient \`a
$\overline{\mathrm{Bir}_2}.$}
\end{lem}

\begin{proof}[\sl D\'emonstration]
Supposons qu'il existe $L$ et $L'$ dans $\mathrm{RL}(Q)$ tels que
$L\wedge L'\not\equiv 0.$ Si $L\wedge~L'$ s'annulait sur un
ensemble de codimension $2,$ alors, d'apr\`es \cite{Sa}, la
dimension de $\mathrm{RL}(Q)$ serait~$2$ ce qui n'est pas le cas.
Par suite $L\wedge L'$ est divisible par une forme lin\'eaire ou
une forme quadratique. Dans la premi\`ere \'eventualit\'e $Q$ est
du type $\ell f$ o\`u $\ell$ d\'esigne une forme lin\'eaire et $f$
un endomorphisme de $\mathbb{C}^3;$ visiblement $Q$ est dans
$\overline{\mathrm{Bir}_2}.$ Dans le second cas $Q$ est \`a conjugaison
pr\`es de l'un des types suivants
\begin{align*}
&(x_2^2+x_0x_1:0:0), &&(x_0x_1:0:0),&&(x_0^2:0:0)
\end{align*}
qui sont tous dans $\overline{\mathrm{Bir}_2}:$ ils sont de la forme
$A\sigma,$ $A\rho$ ou $A\tau$ lorsque $A$
d\'eg\'en\`ere sur un \'el\'ement non inversible. \bigskip

Supposons que pour tout choix de $L$ et $L'$ dans
$\mathrm{RL}(Q)$ on ait $L\wedge L'\equiv 0.$ Tout \'el\'ement 
de~$\mathrm{RL}(Q)$ est de rang $1$ sinon on aurait
$\mathrm{e}(Q)=1.$ En particulier il existe $(a,b,c)$ dans $\mathbb{C}^3$
tel que
\begin{align*}
\mathrm{RL}(Q)=\{(\alpha x_0+\beta x_1+\gamma x_2)(a,b,c)\hspace{0.1cm}|\hspace{0.1cm}(\alpha,
\beta,\gamma)\in\mathbb{C}^3\}
\end{align*}
et $Q$ est d\'eg\'en\'er\'ee; on peut d'ailleurs se ramener, \`a \'equivalence g.d. pr\`es, \`a
\begin{align*}
Q=(Q_0(x_0,x_1,x_2):Q_1(x_0,x_1,x_2):0).
\end{align*}
Montrons qu'un tel \'el\'ement est
dans $\overline{\mathrm{Bir}_2}.$ On peut supposer que $Q$ est purement
quadratique. Les formes normales des couples $(Q_0,Q_1)$ sont
donn\'ees par la classification des pinceaux de coniques que l'on
trouve par exemple dans \cite{HP}. \`A \'equivalence g.d. pr\`es
on se ram\`ene aux mod\`eles suivants
\begin{align*}
& (x_1x_2+x_0x_2+x_0x_1:\alpha x_1x_2+\beta x_0x_2+\gamma x_0x_1:0), &&
(x_2^2:x_0x_1:0), \\
& (x_0^2:x_0x_1+x_2^2:0), && (x_0^2+x_1^2:x_0x_1+x_2^2:0).
\end{align*}

Visiblement $(x_1x_2+x_0x_2+x_0x_1:\alpha x_1x_2+\beta x_0x_2
+\gamma x_0x_1:0)$ est une d\'eg\'en\'erescence de transformations de type
$A\sigma;$ elle est donc dans $\overline{\mathrm{Bir}_2}.$ De m\^eme les
transformations $(x_2^2:x_0x_1:0)$ et $(x_0^2:x_0x_1+x_2^2:0)$ le sont
aussi en consid\'erant des transformations du type $A\rho,$
$A\tau$ avec $\mathrm{det}\hspace{0.1cm} A=0.$ Pour montrer que $(x_0^2+x_1^2:x_0x_1+x_2^2:0)$ est dans
l'adh\'erence de $\mathrm{Bir}_2$ on proc\`ede comme suit. Les $2$ formes
quadratiques ci-dessus s'annulent simultan\'ement sur quatre droites en
position g\'en\'erale plus pr\'ecis\'ement sur
\begin{align*}
&\mathrm{Vect}(\mathrm{i},1,e^{3\mathrm{i}\pi/4}),&&
\mathrm{Vect}(\mathrm{i},1,-e^{3\mathrm{i}\pi/4}),&&
\mathrm{Vect}(-\mathrm{i},1,e^{\mathrm{i}\pi/4}),&&
\mathrm{Vect}(-\mathrm{i},1,-e^{\mathrm{i}\pi/4}).
\end{align*}
Par conjugaison \`a droite pr\`es on peut ramener trois de ces
droites en position standard
\begin{align*}
&\mathrm{Vect}(1,0,0),&&\mathrm{Vect}(0,1,0),&&
\mathrm{Vect}(0,0,1).
\end{align*}
Apr\`es cette modification $(x_0^2+x_1^2:x_0x_1+x_2^2:0)$ est du type $A\sigma,$ la
derni\`ere ligne de $A$ \'etant nulle. La transformation $(x_0^2+x_1^2:x_0x_1+x_2^2:0)$ est
donc une d\'eg\'en\'erescence d'un \'el\'ement de~$\Sigma^3.$
\end{proof}

On obtient finalement le:

\begin{thm}
{\sl Soit $Q$ une transformation rationnelle quadratique.

Si $\mathrm{e}(Q)\geq 2$ alors $Q$ appartient \`a
$\overline{\mathrm{Bir}_2}.$

Si $Q$ est non d\'eg\'en\'er\'ee et $\mathrm{e}(Q)=2,$ alors $Q$ appartient \`a
$\Sigma^1\cup\Sigma^2\cup\Sigma^3.$}
\end{thm}

\subsection{Param\'etrisation de l'adh\'erence de $\mathrm{Bir}_2$}\hspace{0.1cm}

En fait en suivant les raisonnements pr\'ec\'edents on a
le:

\begin{thm}\label{bb}
{\sl Soit $Q$ un \'el\'ement de $\mathring{\mathrm{B}}\mathrm{ir}_2;$ il existe $L$ et $L'$
deux applications lin\'eaires telles que~$Q=L\wedge L'.$ Inversement une
transformation non d\'eg\'en\'er\'ee du type $L\wedge L'$ est dans
$\mathrm{Bir}_2.$ Les \'el\'ements de $\mathrm{Bir}_2$ sont donc exactement les
applications quadratiques d\'eterminantielles.\label{ind30a}}
\end{thm}

Ceci permet d'obtenir de nouveau $\overline{\mathrm{Bir}_2}$
comme l'adh\'erence de l'image d'un morphisme (l'ancienne
\'etant $(A,B)\mapsto A\sigma B$). Plus pr\'ecis\'ement introduisons
les notations suivantes
\begin{small}
\begin{align*}
& L=(a_0x_0+b_0x_1+c_0x_2:a_1x_0+b_1x_1+c_1x_2:a_2x_0+b_2x_1+c_2x_2),\\
& L'=(\alpha_0x_0+\beta_0x_1+\gamma_0x_2:\alpha_1x_0+\beta_1x_1+\gamma_1x_2:
\alpha_2x_0+\beta_2x_1+\gamma_2x_2),
\end{align*}
\begin{align*}
& a=\left[\begin{array}{ccc}a_0\\ a_1\\ a_2\end{array}\right], &&
b=\left[\begin{array}{ccc}b_0\\ b_1\\ b_2\end{array}\right], &&
c=\left[\begin{array}{ccc}c_0\\ c_1\\ c_2\end{array}\right],&& \alpha=\left[\begin{array}{ccc}\alpha_0\\ \alpha_1\\ \alpha_2
\end{array}\right], && \beta=\left[\begin{array}{ccc}
\beta_0\\ \beta_1\\ \beta_2\end{array}\right], &&
\gamma=\left[\begin{array}{ccc}\gamma_0\\ \gamma_1\\ \gamma_2
\end{array}\right],
\end{align*}
\begin{align*}
& Q=(A_0x_0^2+B_0x_1^2+C_0x_2^2+D_0x_1x_2+E_0x_0x_2+F_0x_0x_1:A_1x_0^2+B_1x_1^2+C_1x_2^2\\
&
\hspace{9mm}+D_1x_1x_2+E_1x_0x_2+F_1x_0x_1:A_2x_0^2+B_2x_1^2+C_2x_2^2+D_2x_1x_2+E_2x_0x_2
+F_2x_0x_1),
\end{align*}
\begin{align*}
& A=\left[\begin{array}{ccc}A_0\\ A_1\\ A_2\end{array}\right], &&
B=\left[\begin{array}{ccc}B_0\\ B_1\\ B_2\end{array}\right], &&
C=\left[\begin{array}{ccc}C_0\\ C_1\\ C_2\end{array}\right],
&& D=\left[\begin{array}{ccc}D_0\\ D_1\\ D_2\end{array}\right], &&
E=\left[\begin{array}{ccc}E_0\\ E_1\\ E_2\end{array}\right], &&
F=\left[\begin{array}{ccc}F_0\\ F_1\\ F_2
\end{array}\right].
\end{align*}
\end{small}

L'\'egalit\'e $Q=L\wedge L'$ s'exprime alors comme suit
\begin{align*}
& A=a\wedge\alpha,&& B=b\wedge\beta,&& C=c\wedge\gamma,\\
& D=b\wedge\gamma+c\wedge\beta,&& E=a\wedge\gamma+c\wedge
\alpha,&& F=a\wedge\beta+b\wedge\alpha
\end{align*}

\noindent ce qui, comme annonc\'e, permet de param\'etrer
$\overline{\mathrm{Bir}_2}$ comme l'adh\'erence de l'image d'un
morphisme rationnel.

Consid\'erons la matrice $M$ \`a $10$ lignes et $9$
colonnes d\'efinie par
\begin{align*}
M(Q)=\left[\begin{array}{ccccccccc}
A_0 & 0 & 0 & A_1 & 0 & 0 & A_2 & 0 & 0 \\
0 & B_0 & 0 & 0 & B_1 & 0 & 0 & B_2 & 0 \\
0 & 0 & C_0 & 0 & 0 & C_1 & 0 & 0 & C_2 \\
B_0 & F_0 & 0 & B_1 & F_1 & 0 & B_2 & F_2 & 0\\
F_0 & A_0 & 0 & F_1 & A_1 & 0 & F_2 & A_2 & 0\\
C_0 & 0 & E_0 & C_1 & 0 & E_1 & C_2 & 0 & E_2 \\
E_0 & 0 & A_0 & E_1 & 0 & A_1 & E_2 & 0 & A_2 \\
0 & C_0 & D_0 & 0 & C_1 & D_1 & 0 & C_2 & D_2 \\
0 & D_0 & B_0 & 0 & D_1 & B_1 & 0 & D_2 & B_2 \\
D_0 & E_0 & F_0 & D_1 & E_1 & F_1 & D_2 & E_2 & F_2
\end{array}\right].
\end{align*}

Remarquons que si $Q=L\wedge L',$ les vecteurs
\begin{align*}
&\transp(a_0,a_1,a_2,b_0,b_1,b_2,c_0,c_1,c_2),&& \transp(\alpha_0,
\alpha_1,\alpha_2,\beta_0,\beta_1,\beta_2,\gamma_0,\gamma_1,\gamma_2)
\end{align*}
sont dans le noyau de $M(Q).$ Le Th\'eor\`eme \ref{bb} s'\'enonce
en partie comme suit.

\begin{pro}\label{herge}
{\sl Soit $Q$ une transformation rationnelle quadratique non
d\'eg\'en\'er\'ee. Alors~$Q$
appartient \`a $\mathrm{Bir}_2$ si et seulement si $\textrm{rg } M(Q)\leq
7.$}
\end{pro}

Cette proposition est tr\`es pratique pour tester des
exemples effectifs, en particulier avec \texttt{Maple}, puisqu'il
s'agit de d\'ecider de l'annulation de certains polyn\^omes \`a
l'inverse du Th\'eor\`eme \ref{cri} qui lui demande un calcul de
racines de polyn\^omes.

\section{Quelques orbites sous l'action gauche-droite}\hspace{0.1cm}

Au m\^eme titre que $\sigma$
les transformations rationnelles
\begin{align*}
&(x_0^2:x_1^2:x_2^2) && \text{et} &&(x_0^2:x_1^2+x_0x_2:x_2^2)
\end{align*}
ont leur d\'eterminant jacobien qui s'annule sur trois droites
en position g\'en\'erale.

\begin{pro}\label{mod}
{\sl Les orbites de $\sigma=(x_1x_2:x_0x_2:x_0x_1)$ et $(x_0^2:x_1^2:x_2^2)$ sont de dimension~$14,$
l'orbite de $(x_0^2:x_1^2+x_0x_2:x_2^2)$ est de dimension~$15.$}
\end{pro}

\begin{proof}[\sl D\'emonstration]
Calculons le groupe d'isotropie de $(x_0^2:x_1^2:x_2^2):$ soient $A$,
$B$ dans $\mathrm{SL}_3(\mathbb{C})$ et $\eta$ dans $\mathbb{C}^*$ tels que
\begin{align*}
A(x_0^2:x_1^2:x_2^2)B=\eta(x_0^2:x_1^2:x_2^2).
\end{align*}
Notons $B_i$ les composantes de $B;$ \`a partir de
\begin{align*}
\mathrm{det}\hspace{0.1cm}\mathrm{jac}\hspace{0.1cm} A(x_0^2:x_1^2:x_2^2)B=\mathrm{det}\hspace{0.1cm}\mathrm{jac}\hspace{0.1cm} \eta(x_0^2:x_1^2:x_2^2)
\end{align*}
on obtient $B_1B_2B_3=\eta^3 x_0x_1x_2.$
Comme $(x_0^2:x_1^2: x_2^2)$ commute aux permutations, la composante
neutre du groupe d'isotropie est form\'ee d'\'el\'ements
$(A,B)$ diagonaux. Posons
\begin{align*}
&A:=\left(\mu_1x_0: \mu_2x_1:\frac{x_2}{\mu_1 \mu_2}\right),&&
B:=\left(\alpha_1x_0: \alpha_2x_1:\frac{x_2}{\alpha_1\alpha_2}\right).
\end{align*}
L'\'egalit\'e $A(x_0^2:x_1^2:x_2^2)B=\eta(x_0^2:x_1^2:x_2^2)$ entra\^ine
\begin{align*}
&\mu_1 \alpha_1^2=\eta,&&\mu_2 \alpha_2^2=\eta, &&
1=\mu_1\mu_2\alpha_1^2\alpha_2^2 \eta.
\end{align*}
Il en r\'esulte que $\eta^3=1,$ {\it i.e.} \`a indice fini pr\`es
$\eta$ vaut $1$ et $\mu_1\alpha_1^2=\mu_2\alpha_2^2=1.$ Ainsi le
groupe d'isotropie de $(x_0^2:x_1^2:x_2^2)$ est de dimension $2$
donc $\dim\mathcal{O}_{g.d.}(x_0^2:x_1^2:x_2^2)=14.$

Par ailleurs soient $A,$ $C$ dans $\mathrm{SL}_3(\mathbb{C})$ et
$\eta$ dans $\mathbb{C}^*$ tels que
\begin{align*}
A(x_1^2+x_0x_2:x_0^2:x_2^2)=\eta(x_1^2+x_0x_2:x_0^2:
x_2^2)C.
\end{align*}
Un calcul montre qu'on est n\'ecessairement dans l'une
des situations suivantes
\begin{align*}
& A=\left(\eta\beta^2x_0:\eta\alpha^2x_1:\eta\frac{\beta^4}{\alpha^2}x_2\right), &&
C=\left(\alpha x_0:\beta x_1:\frac{\beta^2}{\alpha}x_2\right)&& \text{et} && \eta^3=
\beta^3=1;\\
& A=\left(\eta\beta^2x_0:\eta\alpha^2x_2:\eta\frac{\beta^4}{\alpha^2}x_1\right), &&
C=\left(\alpha x_2:\beta x_1:\frac{\beta^2}{\alpha}x_0\right)&& \text{et} &&
\eta^3=\beta^3=1.
\end{align*} Le groupe d'isotropie de $(x_1^2+x_0x_2:x_0^2:x_2^2)$ est donc de
dimension $1;$ on en d\'eduit que $\dim\mathcal{O}_{g.d.}(x_1^2+
x_0x_2:x_0^2:x_2^2)=15.$

Bien s\^ur on peut calculer de la m\^eme fa\c{c}on le
groupe d'isotropie de $\sigma,$ mais nous aurons besoin plus loin
de l'espace tangent \`a $\mathcal{O}_{g.d.}(\sigma)$ en $\sigma$

\begin{small}
\begin{eqnarray}
{\rm T}_{\sigma}\mathcal{O}_{g.d.}(\sigma)&=&\{(a_2x_1^2+a_3x_2^2+a_4x_0x_1+a_5x_0x_2
+a_6x_1x_2:b_1x_0^2+b_3x_2^2+b_4x_0x_1+b_5x_0x_2+b_6x_1x_2:\nonumber\\
&
&\hspace{3mm}c_1x_0^2+c_2x_1^2+c_4x_0x_1+c_5x_0x_2+c_6x_1x_2)\hspace{0.1cm}|\hspace{0.1cm}
a_i,\hspace{0.1cm} b_i,\hspace{0.1cm} c_i\in\mathbb{C}\}.\nonumber
\end{eqnarray}
\end{small}

On constate ainsi que ${\rm T}_{\sigma}\mathcal{O}_{g.d.}(\sigma)$ est
de codimension $3$ et l'orbite de $\sigma$ est de dimension~$14.$
\end{proof}

\begin{rem}
Posons $Q_0:=(x_0^2:x_1^2:x_2^2)$ et $Q_1:=(x_0^2:x_1^2+x_0x_2:x_2^2).$ Pour
$\varepsilon\not=0$ la transformation
\begin{align*}
Q_\varepsilon:=\left(x_0:x_1:\frac{x_2}{\varepsilon^2}\right)Q_1(x_0:x_1:\varepsilon
x_2)=(x_0^2:x_1^2+\varepsilon x_0x_2:x_2^2)
\end{align*}
est g.d. conjugu\'ee \`a $Q_1;$ lorsque $\varepsilon$ tend
vers $0,$ cette transformation $Q_\varepsilon$ tend vers $Q_0.$
Puisque la dimension du groupe d'isotropie de $Q_0$ est
strictement inf\'erieure \`a celle du groupe d'isotropie de $Q_1,$ on a
$\mathcal{O}_{g.d.}\hspace{0.1cm}(Q_0)\subsetneq
\overline{\mathcal{O}_{g.d.}\hspace{0.1cm}(Q_1)}.$
\end{rem}

\subsection{Orbites g\'en\'eriques et feuilletage
par les orbites}\hspace{0.1cm}

Soit $\mathcal{C}$ une courbe de degr\'e $3$ poss\'edant
une forme normale de \textsc{Weierstrass}
\begin{align*}
x_1^2=x_0(x_0-1)(x_0-~\eta).
\end{align*}
On sait que \vspace{1mm}
\begin{itemize}
\item si $\eta\not\in\{0,\hspace{0.1cm} 1,\hspace{0.1cm}\infty\}$ alors $\mathcal{C}$ est une
courbe elliptique lisse;

\item si $\eta=0$ ou $\eta=1$, alors $\mathcal{C}$ est une
cubique \`a point double;

\item enfin si \og $\eta=\infty$\fg, alors $\mathcal{C}$ est
l'union de trois droites concourantes.
\end{itemize}
\vspace{1mm}

Notons $j$ l'invariant des courbes elliptiques
\begin{align*}
j(\eta)=256\hspace{0.1cm}\frac{(\eta^2-\eta+1)^3}{\eta^2(1-\eta)^2}.
\end{align*}
Consi\-d\'erons l'action de $\mathrm{PGL}_3(\mathbb{C})$ sur $\mathbb{P}^9
(\mathbb{C})=\mathbb{P}(\mathbb{C}[x_0,x_1,x_2]_3)$ donn\'ee par
\begin{align*}
&\mathrm{PGL}_3(\mathbb{C})\times\mathbb{P}^9(\mathbb{C})\to\mathbb{P}^9(\mathbb{C}),
&&(g,P)\mapsto P\circ g^{-1}.
\end{align*}

Elle d\'efinit un feuilletage $\mathscr{F}$ sur
$\mathbb{P}^9(\mathbb{C})$ dont une int\'egrale premi\`ere rationnelle est induite
par~$j$ (extension de $j$ le long des orbites) et not\'ee $J.$
L'\'egalit\'e $\mathrm{det}(\mathrm{jac}\hspace{0.1cm} AQB^{-1})=\mathrm{det}(\mathrm{jac}\hspace{0.1cm} Q) \circ B^{-1}$ assure
la compatibilit\'e de cette action avec l'action g.d. sur $\mathrm{Rat}_2.$

Remarquons que $j(\infty)=j(0)=j(1)=~\infty.$ Un point
g\'en\'erique de $J^{-1}(\infty)$ est une cubique \`a point
double; cette feuille est de dimension maximale $8$ car le groupe
d'isotropie d'une cubique \`a point double est fini. Par
cons\'equent la restriction de $j$ aux courbes elliptiques et
cubiques \`a point double prend toutes les valeurs; il en
r\'esulte que $\overline{J^{-1}(c)}\setminus\mathrm{Sing}\hspace{0.1cm}
\mathscr{F},$ o\`u $\mathrm{Sing}\hspace{0.1cm}\mathscr{F}$\label{not23b}
d\'esigne le lieu singulier de $\mathscr{F},$ est, pour tout $c$, l'orbite d'une courbe elliptique
lisse ou d'une cubique \`a point double. Les autres configurations
de courbes sont contenues dans $\mathrm{Sing}\hspace{0.1cm}\mathscr{F}$.\bigskip

\begin{pro}
{\sl La dimension d'une orbite g\'en\'erique pour l'action
gauche-droite sur $\mathrm{Rat}_2\simeq\mathbb{P}^{17}(\mathbb{C})$ est $16$. En
particulier l'action gauche-droite induit un feuilletage
$\mathscr{R}$ de codimension $1$ sur $\mathbb{P}^{17}(\mathbb{C})$ d'int\'egrale
premi\`ere $J\circ\mathrm{det}\hspace{0.1cm}\mathrm{jac}.$}
\end{pro}

\begin{proof}[\sl D\'emonstration]
Soit $Q$ dans $\mathrm{Rat}_2$ tel que $\mathcal{C}=\mathrm{det}\hspace{0.1cm}\mathrm{jac}\hspace{0.1cm} Q$ soit
g\'en\'erique (l'existence d'une telle transformation est
assur\'ee par la Proposition \ref{imdelta}). Calculons le groupe
d'isotropie de~$Q;$ soient $A$, $B$ dans $\mathrm{SL}_3(\mathbb{C})$ et
$\eta$ dans $\mathbb{C}^*$ tels que $AQB^{-1}=\eta Q$. Cette \'egalit\'e
conduit \`a 
\begin{align*}
\mathrm{det}(\mathrm{jac}\hspace{0.1cm} Q)\circ B^{-1}=\eta\mathrm{det}\hspace{0.1cm}\mathrm{jac}\hspace{0.1cm} Q.
\end{align*}
La courbe
$\mathcal{C}$ est donc invariante par $B$; puisqu'elle est
g\'en\'erique, $B$ appartient \`a un groupe fini. \`A indice fini
pr\`es on a $A=B=\mathrm{id};$ le groupe d'isotropie de $Q$ est
donc fini.
\end{proof}

\subsection{Lieu singulier}\hspace{0.1cm}

Nous allons pr\'eciser la nature du lieu singulier de
$\mathscr{R}.$ D'apr\`es la d\'emonstration de la Proposition \ref{mod} l'ensemble param\'etr\'e par l'application
\begin{align*}
&\psi=\psi_{s,t,u}\hspace{0.1cm}\colon\hspace{0.1cm}(s,t,u)\mapsto(x_1x_2+sx_0^2:x_0x_2+tx_1^2:x_0x_1
+ux_2^2), && (s,t,u)\in\mathbb{C}^3
\end{align*}
est transverse \`a l'orbite de $\sigma$ en $\sigma.$ On constate
que
\begin{align*}
\mathrm{det}(\mathrm{jac}\hspace{0.1cm}\psi)=(2+8stu)x_0x_1x_2-2sx_0^3-2tx_1^3-2ux_2^3.
\end{align*}
En particulier $\mathrm{det}(\mathrm{jac}\hspace{0.1cm}\psi)$ est une submersion en $(0,0,0);$ le
point $(0,0,0)$ correspond \`a $\sigma$ et l'image de $\mathrm{det}(\mathrm{jac}\hspace{0.1cm}\psi)$ est
une transversale \`a l'orbite de $x_0x_1x_2$ dans $\mathbb{P}^9(\mathbb{C}).$ L'union
de trois droites en position g\'en\'erale est
limite d'une des configurations suivantes: courbe elliptique
lisse, cubique \`a point double ou union d'une droite et d'une
conique en position g\'en\'erale. Rappelons que les orbites des
courbes elliptiques ou des cubiques nodales sont lisses de
dimension~$8$ donc ne sont pas contenues dans le lieu singulier de
$\mathscr{F};$ celle d'une conique et d'une droite en position
g\'en\'erale est de dimension $7.$ En fait l'image de $\mathrm{det}(\mathrm{jac}\hspace{0.1cm}\psi)$
s'identifie \`a la sous-vari\'et\'e param\'etr\'ee par
\begin{align*}
&\varphi=\varphi_{s,t,u}\hspace{0.1cm}\colon\hspace{0.1cm}(s,t,u)\mapsto x_0x_1x_2+sx_0^3+tx_1^3
+ux_2^3, && (s,t,u)\in\mathbb{C}^3
\end{align*}
qui est une transversale \`a l'orbite de $x_0x_1x_2$ en $\varphi_{0,0,
0}=x_0x_1x_2.$ Pour les param\`etres $(s,t,u)$ dans $\mathbb{C}^3\setminus\{(0,0,
0)\}$ satisfaisant $s=t$ et $u=27s^2$ la cubique $\varphi_{s,t,u}$
a deux composantes; il en est de m\^eme lorsqu'on fait une
permutation circulaire sur $(s,t,u).$ Notons qu'une cubique
form\'ee d'une conique et d'une droite en position g\'en\'erale
est dans $\mathrm{Sing}\hspace{0.1cm}\mathscr{F}$ si bien que le lieu
singulier de $\mathscr{F}$ dans $\mathbb{P}^9(\mathbb{C})$ au point $x_0x_1x_2$ a trois
branches locales de dimension $7.$ La courbe~$x_0x_1x_2$ \'etant \`a
l'intersection de ces trois branches on en d\'eduit que l'orbite
de $x_0x_1x_2,$ et en fait l'adh\'erence de cette orbite, est contenue
dans le lieu singulier de $\mathrm{Sing}\hspace{0.1cm} \mathscr{F}.$

Comme $\psi$ est submersive en $(0,0,0)$ on obtient que
$\sigma,$ et donc $\Sigma^3,$ est
dans $\mathrm{Sing}\hspace{0.1cm}\mathrm{Sing}\hspace{0.1cm}\mathscr{R}.$ Nous verrons
plus loin que l'adh\'erence de $\Sigma^3$ dans $\mathbb{P}^{17}(\mathbb{C})$ est
aussi $\overline{\mathrm{Bir}_2}.$ On obtient le:

\begin{thm}
{\sl L'adh\'erence de $\mathrm{Bir}_2$ est une composante irr\'eductible de
$\mathrm{Sing}\hspace{0.1cm}\mathrm{Sing}\hspace{0.1cm}\mathscr{R}.$}
\end{thm}

On peut mener des calculs identiques en les \'el\'ements
\begin{align*}
& Q_0=(x_0^2:x_1^2:x_2^2) &&\text{et} &&\widetilde{Q_1}=(x_0^2+x_1x_2:x_1^2:x_2^2)
\end{align*}
pour lesquels on a
$\mathrm{det}(\mathrm{jac}\hspace{0.1cm} Q_0)=\mathrm{det}(\mathrm{jac}\hspace{0.1cm}\widetilde{Q_1})=8x_0x_1x_2.$ Comme on l'a dit
l'orbite de $\widetilde{Q_1}$ est de dimension $15;$ de plus
$\mathrm{det}\hspace{0.1cm}\mathrm{jac}$ est de rang deux dans une transversale \`a
$\widetilde{Q_1}.$ On en d\'eduit, comme pr\'ec\'edemment, que
l'orbite de $\widetilde{ Q_1}$ est contenue dans
$\mathrm{Sing}\hspace{0.1cm}\mathscr{R}$ et pour des raisons de dimension son
adh\'erence est exactement une branche de
$\mathrm{Sing}\hspace{0.1cm}\mathscr{R}.$ Ainsi

\begin{pro}
{\sl L'adh\'erence de $\mathcal{O}_{g.d.}(x_0^2+x_1x_2:x_1^2:x_2^2)$ est une
branche $\mathcal{S}$ du lieu singulier $\mathrm{Sing}\hspace{0.1cm}
\mathscr{R}.$}
\end{pro}

\begin{rem}
Un calcul explicite montre que dans une transversale \`a
$\mathcal{O}_{g.d.}(Q_0)$ en $Q_0$ l'application $\mathrm{det}(\mathrm{jac})$ est
\'equivalente \`a l'application
\begin{align*}
&(s,t,u)\mapsto(tu,su,st),
\end{align*}
\noindent{\it i.e.} \`a l'application induisant $\sigma!$

Posons
\begin{align*}
&\widetilde{Q_1}:=(x_0^2+x_1x_2:x_1^2:x_2^2),&&\widetilde{Q'_1}:=(x_0^2:x_1^2:x_2^2
+x_0x_1),&&\widetilde{Q''_1}:=(x_0^2:x_1^2+x_0x_2:x_2^2);
\end{align*}
voici le mod\`ele local de $\mathcal{S}$ en $Q_0$ dans une section
$3$-plane:

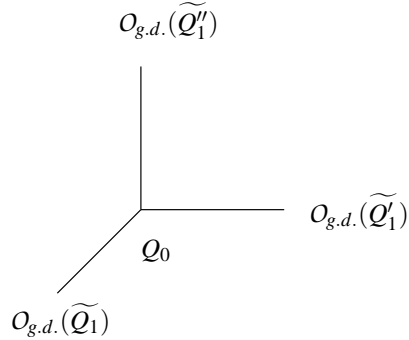
\begin{figure}[H]
\begin{center}
\input{orb1.pstex_t}
\end{center}
\end{figure}

Si $T$ d\'esigne l'orbite de $x_0x_1x_2$ sous l'action de
$\mathrm{PGL}_3(\mathbb{C}),$
on note au passage que $\mathcal{S}$ est une composante
de $\mathrm{det}\hspace{0.1cm}\mathrm{jac}^{-1}(\overline{T})$
qui est donc de dimension $15.$ Remarquons
aussi que $\mathcal{O}_{g.d.}(\widetilde{Q_1}),$ $\mathcal{O}_{
g.d.}(\widetilde{Q'_1})$ et $\mathcal{O}_{g.d.}(\widetilde{
Q''_1})$ sont en fait \'egales, les $\widetilde{Q_1},$  $\widetilde{Q'_1}$
et $\widetilde{Q''_1}$ \'etant \'echang\'ees par permutation de coordonn\'ees.
Nous les avons distingu\'ees car elles apparaissent de fa\c{c}on naturelle
localement en $Q_0$ au travers des d\'eformations de $Q_0$
\begin{align*}
& (x_0^2+tx_1x_2:x_1^2:x_2^2), && (x_0^2:x_1^2+tx_0x_2:x_2^2), &&
(x_0^2:x_1^2:x_2^2+tx_0x_1).
\end{align*}

En $\sigma$ on a la description suivante dans une
section $3$-plane:

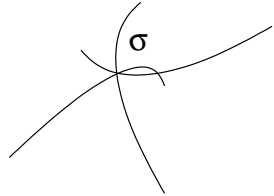
\begin{figure}[H]
\begin{center}
\input{orb2.pstex_t}
\end{center}
\end{figure}

\noindent les trois branches locales faisant partie de l'image
r\'eciproque par $\mathrm{det}(\mathrm{jac})$ de l'adh\'erence de l'orbite d'une conique et d'une
droite en position g\'en\'erale.
\end{rem}

\subsection{Orbites sp\'eciales}\label{spec}\hspace{0.1cm}

On a vu pr\'ec\'edemment que $\mathcal{O}_{g.d.}
(\sigma),$ c'est-\`a-dire $\Sigma^3,$ est de dimension $14.$ Nous
l'avons obtenu par un calcul d'espace tangent. On peut aussi le
faire en explicitant le groupe d'isotropie de $\sigma$ not\'e
\textbf{\textit{${\rm Isot}\hspace{1mm}\sigma.\label{not24}$}} Ce calcul
peut s'av\'erer int\'eressant dans la mesure o\`u il donne des
exemples de relations dans le groupe de \textsc{Cremona}. En effet
soient $(A,B)$ dans $\mathrm{SL}_3(\mathbb{C})\times\mathrm{SL}_3(\mathbb{C})$
tels que $A\sigma=\sigma B;$ alors $(A,B)$ appartient \`a
\begin{align*}
\langle\left(\left(\frac{x_0}{\alpha}:\frac{x_1}{\beta}:\alpha\beta
x_2\right),\hspace{0.1cm}\left(\alpha x_0:\beta
x_1:\frac{x_2}{\alpha\beta}\right)\right), \hspace{0.1cm} \mathscr{S}_6\times
\mathscr{S}_6\hspace{0.1cm}|\hspace{0.1cm}\alpha,\hspace{0.1cm}\beta\in\mathbb{C}^*\rangle
\end{align*}\label{not250}

\noindent o\`u
$\mathscr{S}_6=\{\mathrm{id},\hspace{0.1cm}(x_0:x_2:x_1),\hspace{0.1cm}(x_2:x_1:x_0),\hspace{0.1cm}(x_1:x_0:x_2),
\hspace{0.1cm}(x_1:x_2:x_0),\hspace{0.1cm}(x_2:x_0:x_1)\}.\label{not25}$ Ceci implique
que $\dim{\rm Isot}\hspace{1mm}\sigma=2.$\bigskip

\subsubsection{Orbite de $\rho$}\hspace{0.1cm}

\begin{pro}\label{iso}
{\sl La dimension de $\Sigma^2=\mathcal{O}_{g.d.}(\rho)$ est
$13.$}
\end{pro}

\begin{proof}[\sl D\'emonstration]
Calculons ${\rm Isot}\hspace{1mm}\rho$ sous l'action g.d., {\it
i.e.} cherchons $A$ et $C$ dans $\mathrm{SL}_3(\mathbb{C})$ tels que
$A\rho=\eta \rho C$ o\`u $\eta$ d\'esigne un complexe non nul.
Rappelons que
\begin{align*}
\mathrm{Ind}\hspace{0.1cm}\rho=\{(0:1:0), \hspace{0.1cm} (1:0:0)\};
\end{align*}
la relation $A\rho=\eta\rho C$ assure que $C$ laisse
$\mathrm{Ind}\hspace{0.1cm}\rho$ invariant. Or les points d'ind\'etermination
de $\rho$ \og ne sont pas de m\^eme nature\fg\hspace{0.1cm} donc $C$ fixe
$(0:1:0)$ et $(1:0:0)$; par suite $C$ est de la forme
\begin{align*}
&(ax_0+bx_2:cx_1+dx_2:ex_2), && ace\not=0.
\end{align*}
Un calcul conduit \`a
\begin{align*}
&A=(\eta\gamma\delta x_0+\eta\beta\delta x_2:\eta\alpha^2 x_1:
\eta\alpha\delta x_2), && C=(\gamma x_0+\beta x_2:\delta x_1:\alpha x_2),&&
\eta^3\alpha^2\delta= \alpha \gamma\delta=1.
\end{align*}
La dimension du groupe d'isotropie est donc $3$.
\end{proof}

Le calcul de ${\rm Isot}\hspace{1mm}\rho$ assure qu'on a les
relations suivantes
\begin{align*}
&(\gamma\delta x_0+\beta\delta x_2:\alpha^2x_1:\alpha\delta x_2)\rho=\rho(
\gamma x_0+\beta x_2:\delta x_1:\alpha x_2), &&\alpha,\hspace{0.1cm}\gamma,\hspace{0.1cm}
\delta\in\mathbb{C}^*,\hspace{0.1cm}\beta\in\mathbb{C}.
\end{align*}

\bigskip

\subsubsection{Orbite de $\tau$}\hspace{0.1cm}

Calculons le groupe d'isotropie de $\tau=(x_0^2:x_0x_1:x_1^2
-x_0x_2).$ Rappelons que le lieu d'ind\'etermination de $\tau$ est r\'eduit \`a $\{(0:0:1)\};$ la
relation $A\tau=\eta\tau C,$ o\`u $\eta$ d\'esigne un \'el\'ement de $\mathbb{C}^*$
et $A,$ $C$ deux \'el\'ements de $\mathrm{SL}_3(\mathbb{C}),$ assure que $C$ laisse $(0:0:1)$
invariant. Par suite $C$ est de la forme
\begin{align*}
&(a_0x_0+b_0x_1:a_1x_0+b_1x_1:a_2x_0+b_2x_1 +c_2x_2), &&
c_2(a_0b_1-a_1b_0)\not=0.
\end{align*}
Un calcul \'el\'ementaire conduit \`a

\begin{small}
\begin{align*}
&A=(\eta\varepsilon^2 x_0:\eta\varepsilon (\beta-\alpha+
\varepsilon) x_0+\eta\alpha\varepsilon x_1:\eta(\varepsilon(\beta
+\delta-\alpha-\gamma-\varepsilon)-2\alpha\beta+2\alpha^2+\beta^2)
x_0\\
& \hspace{6mm}+\eta(\alpha\varepsilon+\varepsilon\gamma+2\alpha
\beta-\alpha^2)x_1+\eta\alpha^2x_2),
\\
&C=\left(\varepsilon x_0:(\beta+\varepsilon-\alpha)x_0+\alpha x_1:
\left(\beta-\alpha+2\varepsilon+\gamma-\delta-\frac{\alpha}{\varepsilon}
\right)x_0 +\left(\alpha-\gamma-\frac{\alpha}{\varepsilon}\right)x_1+
\frac{\alpha}{\varepsilon} x_2\right)
\end{align*}
\end{small}
o\`u $\alpha^2=\eta^3\varepsilon^3\alpha=1.$ Par suite $\dim{\rm
Isot}\hspace{1mm}\tau=4.$ On en d\'eduit la:

\begin{pro}\label{iso2}
{\sl La dimension de $\Sigma^1$ est $12.$}
\end{pro}

D'apr\`es ce qui pr\'ec\`ede on a les
relations suivantes: $A\tau=\tau B$ lorsque
\begin{align*}
&A=\left[
\begin{array}{ccc}
\alpha\varepsilon & 0 & \beta\varepsilon \\
\varepsilon\gamma+2\alpha\beta & \alpha^2 & (\varepsilon
\delta+ \beta^2)\\
0 & 0 & \varepsilon^2
\end{array}
\right],&& B=\left[
\begin{array}{ccc}
\alpha & \beta & 0\\
0 & \varepsilon & 0\\
\gamma & \delta & \alpha/\varepsilon
\end{array}
\right], &&\beta,\hspace{0.1cm}\gamma,\hspace{0.1cm}\delta\in\mathbb{C},\hspace{0.1cm}\alpha,\hspace{0.1cm}\varepsilon
\in\mathbb{C}^*.
\end{align*}
\bigskip

\subsubsection{Orbite de $x_0(x_0:x_1:x_2)$}\hspace{0.1cm}

\begin{pro}
{\sl La dimension de $\Sigma^0=\mathcal{O}_{g.d.}(x_0(x_0:x_1:x_2))$ est~$10.$}
\end{pro}

\begin{proof}[\sl D\'emonstration]
On peut calculer le groupe d'isotropie de $x_0(x_0:x_1:x_2)$
mais on peut aussi calculer ${\rm T}_{x_0(x_0:x_1:x_2)}\mathcal{O}_{g.d.}
(x_0(x_0:x_1:x_2));$ cet espace tangent est pr\'ecis\'ement
\begin{eqnarray}
&
&\{(\alpha_1x_0^2+\alpha_4x_0x_1+\alpha_5x_0x_2:\beta_1x_0^2+\beta_2x_1^2+\beta_4
x_0x_1+\beta_5x_0x_2+\beta_6x_1x_2:\nonumber\\
& &\hspace{6mm}\gamma_1x_0^2+\beta_6x_2^2+\gamma_4x_0x_1+\gamma_5x_0x_2+
\beta_2x_1x_2)\hspace{0.1cm}|\hspace{0.1cm}\alpha_i,\hspace{0.1cm}\beta_i,\hspace{0.1cm}\gamma_i
\in\mathbb{C}\}\nonumber
\end{eqnarray}
qui est de codimension $7$ dans l'espace des triplets de formes
quadratiques.
\end{proof}

\subsubsection{Orbites dans $\overline{\mathrm{Bir}_2}\setminus \mathrm{Bir}_2$}
(\emph{voir} \cite{HP}, tome $2,$ page $304$)\hspace{0.1cm}

Soit $f$ une transformation de $\overline{\mathrm{Bir}_2} \setminus
\mathrm{Bir}_2;$ il existe une suite de transformations birationnelles
quadratiques $(f_i)_i$ dont $f$ est la limite. On peut d'ailleurs, comme
nous le verrons dans la Proposition \ref{lili}, choisir les $f_i$ dans
$\Sigma^3.$ Montrons par l'absurde que $\mathrm{det}(\mathrm{jac}\hspace{0.1cm} f)\equiv 0;$
supposons donc que $\mathrm{det}(\mathrm{jac}\hspace{0.1cm} f)\not\equiv 0.$ Les $f_i$ contractent
$\mathrm{det}( \mathrm{jac}\hspace{0.1cm} f_i)\equiv 0$ qui est l'union de trois droites
(Th\'eor\`eme \ref{cri}); par continuit\'e $f$ contracte donc
trois droites qui sont n\'ecessairement distinctes et concourantes sinon $f$
serait birationnelle. En reprenant la d\'emonstration du
Th\'eor\`eme \ref{cri} on en d\'eduit que $f$ s'\'ecrit \`a
\'equivalence g.d. pr\`es $(x_0^2:x_1^2:x_0x_1)$ d'o\`u $\mathrm{det}(\mathrm{jac}\hspace{0.1cm} f)\equiv
0:$ contradiction.

Puisque $\mathrm{det}(\mathrm{jac}\hspace{0.1cm} f)\equiv 0,$ on sait d'apr\`es la
Remarque \ref{padef} que $f$ est, \`a \'equivalence g.d. pr\`es,
de la forme $(Q_0:Q_1:0)$ ou du type $(x_0^2:x_1^2:x_0x_1).$ En
d\'eformant les transformations du type $A\sigma$ (resp. $A\rho$,
resp. $A\tau$), o\`u $A$ d\'esigne un automorphisme de $\mathbb{P}^2(\mathbb{C}),$
on obtient des transformations du type $\widetilde{A}\sigma$
(resp. $\widetilde{A}\rho,$ resp.  $\widetilde{A}\tau$) avec
$\mathrm{det}\hspace{0.1cm}\widetilde{A} =0;$ les orbites de tels \'el\'ements sont dans
$\overline{\mathrm{Bir}_2} \setminus\mathrm{Bir}_2;$ par exemple
\begin{align*}
&\mathcal{O}_{g.d.}(x_2(x_1:x_0:0)), && \mathcal{O}_{g.d.}(x_0x_1(0:0:1)),
&&
\mathcal{O}_{g.d.}(x_2(x_2:x_0:0)), \\
& \mathcal{O}_{g.d.}(x_2^2:0:x_0x_1), &&\mathcal{O}_{g.d.}(x_0^2:0:0), && \mathcal{O}_{g.d.}(x_0x_1:x_0^2+x_1x_2:0),\\
& \mathcal{O}_{g.d.}(0:x_0^2+x_1x_2:x_2^2), && \mathcal{O}_{g.d.}(x_2^2+
x_0x_1:0:0), &&\mathcal{O}_{g.d.}(x_0x_1:x_0^2:x_1^2),\\
&\mathcal{O}_{g.d.}(0: x_0^2:x_1^2). && &&
\end{align*}

\section{Conditions d'incidence; lissit\'e de $\mathrm{Bir}_2$}\hspace{0.1cm}

On se propose dans ce paragraphe d'\'etudier
les conditions d'incidence entre les $\Sigma^i$ et la lissit\'e
de $\mathrm{Bir}_2.$

\begin{pro}\label{lili}
{\sl On a les inclusions
\begin{align*}
&\Sigma^0\subset\overline{\Sigma^1},&&\Sigma^1\subset \overline{
\Sigma^2},&&\Sigma^2\subset \overline{\Sigma^3};
\end{align*}
en particulier $\Sigma^3$ est dense dans $\mathrm{Bir}_2.$}
\end{pro}

\begin{proof}[\sl D\'emonstration]
En composant $\sigma$ \`a droite par $(x_2:x_1:\varepsilon x_0+x_2)$ on a
\begin{align*}
\sigma_1^\varepsilon=(x_1(\varepsilon x_0+x_2):x_2(\varepsilon
x_0+x_2):x_1x_2)
\end{align*}
qui est donc pour $\varepsilon\not=0$ dans
$\mathcal{O}_{g.d.}(\sigma).$ Mais $\sigma_1^\varepsilon$ est g.d.
conjugu\'e \`a
\begin{align*}
\sigma_2^\varepsilon= (x_0x_1:(\varepsilon x_0+x_2)x_2:x_1x_2).
\end{align*}
On remarque que $\displaystyle\lim_{\varepsilon\to 0}\sigma_2^\varepsilon=
(x_0x_1:x_2^2:x_1x_2)=\rho;$ d'o\`u l'inclusion
$\Sigma^2\subset\overline{\Sigma^3}.$\bigskip

Composons $\rho$ \`a droite par $(x_2:x_0+x_1:x_0);$ cela donne
\`a \'equivalence g.d. pr\`es 
\begin{align*}
(x_1x_2+x_0x_2:x_0^2:x_0x_1).
\end{align*}
Composant de
nouveau \`a droite par $(x_0:x_1:x_1+x_2),$ on obtient \`a conjugaison
g.d. pr\`es la transformation
\begin{align*}
f=(x_1x_2+x_1^2+x_0x_2:x_0^2:x_0x_1).
\end{align*}
Soit $g_\varepsilon:=f(x_0/\varepsilon:x_1:-\varepsilon x_2);$ \`a
conjugaison g.d. pr\`es $g_\varepsilon$ s'\'ecrit $(-\varepsilon
x_1x_2+x_1^2-x_0x_2:x_0^2:x_0x_1).$ Pour $\varepsilon=0$ on obtient la
transformation $\tau.$ Par suite $\Sigma^1$ est inclus dans
$\overline{\Sigma^2}.$

Si $\varepsilon$ est non nul, \`a \'equivalence g.d. pr\`es, $\tau$ s'\'ecrit
\begin{align*}
(x_0^2:x_0x_1:\varepsilon^2x_1^2+x_0x_2);
\end{align*}
pour $\varepsilon=0$ on obtient $x_0(x_0:x_1:x_2)$ qui est dans $\Sigma^0.$ Il en
r\'esulte que $\Sigma^0\subset\overline{ \Sigma^1}.$
\end{proof}

R\'ecapitulons

\begin{thm}
{\sl Les adh\'erences \'etant prises dans $\mathrm{Bir}_2,$ on a
\begin{align*}
& \overline{\Sigma^0}=\Sigma^0,&&\overline{\Sigma^1}=\Sigma^0
\cup\Sigma^1,&&\overline{\Sigma^2}=\Sigma^0\cup\Sigma^1\cup
\Sigma^2,
\end{align*}
\begin{align*}
&\mathring{\mathrm{B}}\mathrm{ir}_2=\Sigma^1\cup\Sigma^2\cup\Sigma^3, && \mathrm{Bir}_2=\overline{
\Sigma^3}=\Sigma^0\cup\Sigma^1\cup \Sigma^2\cup\Sigma^3
\end{align*}
avec
\begin{align*}
&\dim\Sigma^0=10,&&\dim\Sigma^1=12,&&\dim\Sigma^2=13 &&\text{et}
&&\dim\Sigma^3=14.
\end{align*}}
\end{thm}
\bigskip

\begin{thm}\label{lis}
{\sl L'ensemble des transformations birationnelles quadratiques
est lisse dans l'ensemble des transformations rationnelles.}
\end{thm}

\begin{proof}[\sl D\'emonstration]
Le fait que chaque $\Sigma^i$ soit une orbite et les
conditions d'incidence impliquent qu'il suffit de montrer
que l'adh\'erence de $\Sigma^3$ est lisse le long de
$\Sigma^0.$

L'espace tangent \`a $\Sigma^0$ en $x_0(x_0:x_1:x_2)$
est donn\'e par
\begin{eqnarray}
{\rm T}_{x_0(x_0:x_1:x_2)}\Sigma^0&=&\{(\alpha_1x_0^2+\alpha_4x_0x_1+\alpha_5
x_0x_2:\beta_1x_0^2+\beta_2x_1^2+\beta_4x_0x_1+\beta_5x_0x_2+\beta_6x_1x_2:\nonumber\\
&
&\hspace{3mm}\gamma_1x_0^2+\beta_6x_2^2+\gamma_4x_0x_1+\gamma_5x_0x_2+\beta_2x_1x_2
)\hspace{0.1cm}|\hspace{0.1cm}\alpha_i,\hspace{0.1cm}\beta_i,\hspace{0.1cm}\gamma_i\in\mathbb{C}\}. \nonumber
\end{eqnarray}

L'espace vectoriel $S$ engendr\'e par

\begin{small}
\begin{align*}
& (x_1^2:0:0), && (x_2^2:0:0), && (x_1x_2:0:0), && (0:x_2^2:0), && (0:0:x_1^2),
&& (0:0:x_2^2), &&(0:0:x_1x_2)
\end{align*}
\end{small}

\noindent est un suppl\'ementaire de ${\rm T}_{x_0(x_0:x_1:x_2)}\Sigma^0$ dans $\mathrm{Rat}_2.$
Soit $f$ dans $\Sigma^3\cap\{x_0(x_0:x_1:x_2)+S\},$ il s'\'ecrit
\begin{align*}
(x_0^2+Ax_1^2+Bx_2^2+Cx_1x_2:x_0x_1+ax_2^2:x_0x_2+\alpha x_1^2+\beta
x_2^2+\gamma x_1x_2).
\end{align*}
N\'ecessairement $f$ a trois points d'ind\'etermination.

Supposons $a\not=0;$ on remarque que la seconde
composante d'un point d'ind\'etermination de~$f$ est
n\'ecessairement non nulle. Si $(x_0:x_1:x_2)$ appartient \`a
$\mathrm{Ind}\hspace{0.1cm} f$ alors $x_0=-ax_2^2/x_1.$ Calculons
$f(-ax_2^2/x_1:x_1:x_2)$
\begin{eqnarray}
f(-ax_2^2/x_1:x_1:x_2)&=&(a^2x_2^4+Ax_1^4+Bx_1^2x_2^2+Cx_1^3x_2:
0:-ax_2^3+\alpha x_1^3+ \beta x_1x_2^2+\gamma x_1^2x_2)\nonumber\\
&=&(P:0:Q).\nonumber
\end{eqnarray}
Comme $f$ doit avoir trois points d'ind\'etermination,
les polyn\^omes $P$ et $Q$ doivent s'annuler sur trois droites distinctes.
En particulier $Q$ divise $P$
\begin{align*}
a^2x_2^4+Ax_1^4+Bx_1^2x_2^2+Cx_1^3x_2=(Mx_1+Nx_2)(-ax_2^3+\alpha x_1^3
+\beta x_1x_2^2+\gamma x_1^2x_2).
\end{align*}
De sorte que
\begin{equation}\label{trois}
B=-\beta^2-a\gamma,\hspace{6mm}C=-\beta\gamma-a\alpha,\hspace{6mm}A=-
\alpha\beta.
\end{equation}
Ces trois \'equations d\'efinissent un graphe lisse, passant par $f$ et
$x_0(x_0:x_1:x_2),$ de codimension~$3$ comme $\Sigma^3.$

Si maintenant $a$ est nul un point d'ind\'etermination $(x_0:x_1:x_2)$
de $f$ satisfait $x_0x_1=0.$ Si~$x_0=0$ on aura
\begin{align*}
 &f(0:x_1:x_2)=(Ax_1^2+Bx_2^2+Cx_1x_2:0:\alpha x_1^2+\beta x_0^2+\gamma
x_1x_2)
\end{align*}
et si $x_1=0$ on a
\begin{align*}
&f(x_0:0:x_2)=(x_0^2+Bx_2^2:0:x_0x_2+\beta x_2^2).
\end{align*}
Pour que $f$ ait un point d'ind\'etermination de type $(x_0:0:x_2)$
il faut que $B=-\beta^2$ et c'est en fait une \'equivalence. Si tel
est le cas $f$ ne poss\`ede qu'un point d'ind\'etermination de ce
type. Comme
$f$ doit avoir trois points d'ind\'etermination, deux sont de
type $(0:x_1:x_2)$ et les polyn\^omes $Ax_1^2+Bx_2^2+Cx_1x_2$ et $\alpha
x_1^2+\beta x_2^2+\gamma x_1x_2$ sont $\mathbb{C}$-colin\'eaires. On obtient
les conditions
\begin{itemize}
\item $a=0,$ $B=-\beta^2,$ $A=-\alpha\beta$ et $C=-\beta\gamma $ si $\beta$
est non nul;

\item $a=B=\beta=A\gamma-\alpha C=0$ sinon.
\end{itemize}
On remarque que dans ce second cas $f$ ne peut avoir trois points
d'ind\'etermination.
Finalement on constate que $\Sigma^3\cap\{x_0(x_0:x_1:x_2)+S\}$ est
contenu dans le graphe d\'etermin\'e par les \'equations~(\ref{trois}). Il
en est de m\^eme de l'adh\'erence $\overline{\Sigma^3}\cap\{x_0(x_0:x_1:
x_2)+S\}$ qui, pour des raisons de dimension, co\"{\i}ncide donc avec
ce graphe. Par suite $\overline{\Sigma^3}$ est lisse le long de
$\Sigma^0.$
\end{proof}

\begin{rem}
Comme nous l'avons dit le fait que $\overline{\Sigma^3}$ soit lisse
le long de $\Sigma^0$ implique \`a cause, entre autres, des
conditions d'incidence que $\overline{\Sigma^3}$ est lisse le
long de $\Sigma^2$ et $\Sigma^1.$ Nous allons toutefois
d\'emontrer ces deux affirmations en construisant des familles
lin\'eaires de transformations birationnelles.

La d\'eformation
$\rho_\epsilon:=\rho+ \varepsilon(0:x_1^2:0)$ de $\rho$ compte, pour
$\varepsilon$ non nul, trois points d'ind\'etermination donc est
dans $\Sigma^3.$ Par ailleurs $\frac{\partial\rho}{\partial
\varepsilon}=(0:x_1^2:0)$ n'appartient pas \`a ${\rm T}_\rho\Sigma^2.$
On retrouve ainsi la lissit\'e de $\overline{\Sigma^3}$ en $\rho$ et par
suite le long de $\Sigma^2.$\bigskip

Montrons que $\overline{\Sigma^3}$ est lisse le long de
$\Sigma^1.$

Posons
\begin{small}
\begin{align*}
&\widetilde{\tau}:=(x_0x_1:x_1x_2+x_0^2:x_1^2)=(-x_0+x_1:2x_0-x_1+x_2:x_0)\tau
(-x_1-x_0-x_1:-x_0-2x_1+x_2).
\end{align*}
\end{small}
L'espace tangent ${\rm T}_{\widetilde{\tau}}\Sigma^1$ \`a $\Sigma^1$ en $\widetilde{\tau}$ est donn\'e par
\begin{eqnarray}
&&\{(\alpha_1x_0^2+\alpha_2x_1^2+\alpha_4x_0x_1+\alpha_5x_0x_2
+\alpha_6x_1x_2:\beta_1x_0^2+\beta_2x_1^2+\alpha_5x_2^2+\beta_4x_0x_1\nonumber\\
&&\hspace{3mm}+\beta_5x_0x_2 +\beta_6x_1x_2:\gamma_1x_0^2+\gamma_2x_1^2+
\gamma_4x_0x_1+(\gamma_1+2 \alpha_5)x_1x_2)\hspace{0.1cm}|\hspace{0.1cm}\alpha_i,\hspace{0.1cm}\beta_i,
\hspace{0.1cm}\gamma_i\in\mathbb{C}\}\nonumber
\end{eqnarray}
Consid\'erons le sous-espace vectoriel $S$ de $\mathrm{Rat}_2$ engendr\'e
par
\begin{align*}
&(x_2^2:0:0), &&(x_0x_2:0:0),&&(0:0:x_2^2),&&(0:0:x_0x_2),&& (0:0:x_0x_1).
\end{align*}
On constate que $\dim S=5$ et ${\rm T}_{\widetilde{\tau}}\Sigma^1\cap
S=\{0\},$ {\it i.e.} $S$ est un suppl\'ementaire de
${\rm T}_{\widetilde{\tau}}\Sigma^1$ dans $\mathrm{Rat}_2.$ D\'eterminons
$\overline{\Sigma^3}\cap\{\widetilde{\tau}+S\};$ soit
\begin{align*}
f=\widetilde{\tau}+(\alpha x_2^2+\beta x_0x_2:0:\gamma x_2^2+\delta x_0x_2+\varepsilon
x_0x_1)
\end{align*}
un \'el\'ement de $\Sigma^3\cap \{\widetilde{\tau}+S\}.$ On note
que tout point de $\mathbb{P}^2(\mathbb{C})$ de la forme $(x_0:x_1:0)$ ne peut \^etre
d'ind\'etermination pour $f;$ ainsi un point de $\mathrm{Ind}\hspace{0.1cm} f$
est du type $(x_0:-x_0^2/x_2:x_2).$ Calculons $f(x_0:-x_0^2/x_2:x_2)$
\begin{align*}
f(x_0:-x_0^2/x_2:x_2)=(-x_0^3x_2+\alpha x_2^4+\beta x_0x_2^3:0:\gamma x_2^4+\delta x_0x_2^3
-\varepsilon x_0^2x_2^2+x_0^4).
\end{align*}
Comme $f$ a trois points d'ind\'etermination; les polyn\^omes
\begin{align*}
& -x_0^3x_2+\alpha x_2^4+ \beta x_0x_2^3 && \text{et} &&
\gamma x_2^4+\delta x_0x_2^3 -\varepsilon x_0^2x_2^2+x_0^4
\end{align*}
s'annulent sur trois droites communes. Pour que ceci ait lieu
il faut que $\gamma=0,$ $\delta=-\alpha$ et $\varepsilon=\beta.$
Il en r\'esulte que $\overline{\Sigma^3}\cap \{\widetilde{\tau}+S\}$
est contenu dans le $2$-plan $\mathcal{P}$ d\'efini par
\begin{align*}
\mathcal{P}=\{(x_0x_1+\alpha x_2^2+\beta x_0x_2:x_0^2+x_1x_2:-\alpha x_0x_2+\beta
x_1x_2+x_1^2)\hspace{0.1cm}|\hspace{0.1cm} \alpha,\hspace{0.1cm}\beta\in\mathbb{C}\}.
\end{align*}
Puisque $\dim\Sigma^3\cap\{\widetilde{\tau}+S\}=2,$ on a $\overline{\Sigma^3}\cap\{\widetilde{
\tau}+S\}= \mathcal{P}.$ Comme $\Sigma^1=\mathcal{O}_{g.d.}
(\widetilde{\tau})$ on v\'erifie de nouveau que $\overline{\Sigma^3}$ est
lisse le long de $\Sigma^1.$
\end{rem}

\begin{pro}
{\sl L'adh\'erence de $\mathrm{Bir}_2$ dans $\mathbb{P}^{17}(\mathbb{C})\simeq\mathrm{Rat}_2$
n'est pas lisse.}
\end{pro}

\begin{proof}[\sl D\'emonstration]
Soit $f$ la transformation rationnelle d\'eg\'en\'er\'ee d\'efinie par $x_2(x_0:x_1:0).$
L'espace tangent \`a $\mathcal{O}_{g.d.}(f)$ en $f$ est donn\'e
par

\begin{small}
\begin{eqnarray}
{\rm T}_f\mathcal{O}_{g.d.}(f)&=&\{(\alpha_1x_0^2+\alpha_3x_2^2+
\alpha_4x_0x_1+\alpha_5x_0x_2+\alpha_6x_1x_2:\alpha_4x_1^2+\beta_3x_2^2\nonumber\\
& &\hspace{3mm}+
\alpha_1x_0x_1+\beta_5x_0x_2+\beta_6x_1x_2:\gamma_5x_0x_2+\gamma_6x_1x_2)\hspace{0.1cm}|\hspace{0.1cm}\alpha_i,
\hspace{0.1cm}\beta_i,\hspace{0.1cm}\gamma_i\in\mathbb{C}\}.\nonumber
\end{eqnarray}
\end{small}

Un suppl\'ementaire $S$ de ${\rm T}_f\mathcal{O}_{g.d.}(f)$ est
l'espace de dimension $8$ engendr\'e par
\begin{align*}
&(x_1^2:0:0),&&(0:x_0^2:0),&&(0:x_1^2:0),&&(0:x_0x_1:0),
\end{align*}
\begin{align*}
& (0:0:x_0^2), && (0:0:x_1^2), && (0:0:x_2^2), && (0:0:x_0x_1).
\end{align*}
Nous allons montrer que
l'intersection de $\{f+S\}$ avec $\overline{\Sigma^3}$ contient un
sous-ensemble analytique non lisse de codimension $3.$ Puisque
$\overline{\Sigma^3}$ est aussi de codimension $3$ on constatera,
en utilisant l'action g.d., la non lissit\'e de
$\overline{\Sigma^3}$ le long de l'orbite de $f.$ Pour cela nous donnons une condition
suffisante pour qu'un \'el\'ement de $\{f+S\}$ soit dans
$\Sigma^3,$ {\it i.e.} ait trois points d'ind\'etermination. Un
\'el\'ement $Q$ de $\{f+S\}$ s'\'ecrit
\begin{align*}
(x_0x_2+ax_1^2:x_1x_2+bx_0^2+cx_1^2+dx_0x_1:ex_0^2+fx_1^2+
gx_2^2+hx_0x_1).
\end{align*}
Les points d'ind\'etermination sont donn\'es par les trois \'equations
suivantes
\begin{align*}
\left\{\begin{array}{lll}
x_0x_2+ax_1^2=0\\
x_1x_2+bx_0^2+cx_1^2+dx_0x_1=0\\
ex_0^2+fx_1^2+hx_0x_1=0
\end{array}
\right.
\end{align*}
ce qui conduit apr\`es \'elimination de $x_2$ \`a $P_1=
P_2=0$ o\`u
\begin{align*}
& P_1=-ax_1^3+bx_0^3+cx_0x_1^2+dx_0^2x_1, &&
P_2=ex_0^4+fx_0^2x_1^2+a^2gx_1^4+hx_0^3x_1.
\end{align*}
Notons que si, pour certaines valeurs des
param\`etres, $P_1$ s'annule sur trois droites distinctes et
divise $P_2,$ la transformation $Q$ correspondante aura trois
points d'ind\'etermination et sera donc birationnelle, en fait
dans $\Sigma^3.$ La divisibilit\'e de $P_2$ par $P_1$ s'exprime
par
\begin{equation}\label{sys}
P_2=(Ax_0+Bx_1) P_1\hspace{5mm}\Leftrightarrow\hspace{5mm}\left\{\begin{array}{lllll}
e=bA \\
f=cA+dB\\
a^2g=-aB\\
h=dA+bB\\
0=-aA+cB
\end{array}\right.
\end{equation}
On remarque que l'ensemble $\Lambda$ des param\`etres tels que
\begin{align*}
& a=0,&& bf-ce=0,&& bh-de=0
\end{align*}
satisfait le syst\`eme (\ref{sys}) (avec $A=e/b$ et $B=0$).
L'ensemble $\Lambda$ est de codimension $3$ et n'est pas lisse,
l'intersection $\Lambda'$ des quadriques $bf-ce=0$ et $bh-de=0$
ne l'\'etant pas. On constate en effet que $\Lambda'$ contient l'espace
lin\'eaire $E$ donn\'e par $b=e=0$ mais ne se r\'eduit
pas \`a $E:$ par exemple l'espace d\'efini par $b=c=d=e=f=h$ est contenu dans
$\Lambda'$ et pas dans $E.$ Comme $\text{ codim }E=\text{codim }\Lambda'$
l'ensemble $\Lambda'$ n'est pas irr\'eductible et par suite n'est
pas lisse; il en est donc de m\^eme pour $\Lambda.$ Maintenant si
$a=b=e=0$ (resp. $b=c=d=e =f=h=1,$ $a=0$) le polyn\^ome $P_1$
vaut $cx_0x_1^2+dx_0^2x_1$ (resp. $x_0^3+x_0x_1^2 +x_0^2x_1$) et s'annule
g\'en\'eriquement sur trois droites distinctes. Ainsi on a
construit dans $\overline{\Sigma^3}\cap\{f+S\}$ un ensemble
analytique non lisse de codimension $3.$
\end{proof}

\clearemptydoublepage
\chapter{Germes de flots birationnels quadratiques}\label{germgerm}

\section{G\'en\'eralit\'es sur les germes de flots birationnels
quadratiques}\label{floflo}\hspace{0.1cm}

\subsection{Quelques rappels utiles}\hspace{0.1cm}

\subsubsection{N\oe ud col}\hspace{0.1cm}

Soit $\mathcal{F}$ un germe de feuilletage holomorphe singulier \`a
l'origine de $\mathbb{C}^2$ donn\'e par un germe de champ holomorphe $X=A
\frac{\partial}{\partial x_0}+B\frac{\partial}{\partial x_1}$ \`a singularit\'e
isol\'ee en l'origine, $A(0)=B(0)=0.$ On note $JX(0)$ la matrice jacobienne
de $X$ en $0$
\begin{align*}
JX(0)=\left[\begin{array}{cc}
\frac{\partial A}{\partial x_0}(0) & \frac{\partial A}{\partial x_1}(0)\\
\frac{\partial B}{\partial x_0}(0) & \frac{\partial B}{\partial x_1}(0)
\end{array}
\right].
\end{align*}
On dit que le point singulier est de type \textbf{\textit{n\oe ud
col}}\label{ind31} si l'une des valeurs propres de $JX(0)$ est nulle
et l'autre non. La th\'eorie des formes normales (\cite{Ma}) dit que
le champ $X$ est alors formellement conjugu\'e \`a un champ rationnel du
type
\begin{align*}
& Y=\eta x_0\frac{\partial}{\partial x_0}+\frac{x_1^{k+1}}{1-\varepsilon x_1^k}
\frac{\partial}{\partial x_1},&& k\geq 1,\hspace{0.1cm}\varepsilon\in\mathbb{C},\hspace{0.1cm}\eta
\in\mathbb{C}^*.
\end{align*}

Un calcul \'el\'ementaire montre que $Y$ ne peut poss\'eder
d'int\'egrale premi\`ere (formelle) m\'eromorphe; en particulier il en
est de m\^eme pour $X.$\bigskip

Nous utiliserons ce fait sous la forme suivante: {\it
un feuilletage alg\'ebrique $\mathcal{F}$ de $\mathbb{P}^2(\mathbb{C})$ poss\'edant un
point singulier de type n\oe ud col n'a pas d'int\'egrale
premi\`ere rationnelle non constante.}

\subsubsection{Travaux de \textsc{Cantat} et
\textsc{Favre}}\label{cantatfavre}\hspace{0.1cm}

Dans \cite{CaFa} \textsc{Cantat} et \textsc{Favre}
\'etudient les paires $(f,\mathcal{F})$ o\`u $f$ est un \'el\'ement de
$\mathrm{Bir}(\mathbb{P}^2(\mathbb{C}))$ et $\mathcal{F}$ un feuilletage alg\'ebrique sur
$\mathbb{P}^2(\mathbb{C})$ laiss\'e invariant par $f,$ {\it i.e.} $f^*\mathcal{F}=\mathcal{F}.$ Ils
mentionnent l'exemple suivant

\begin{eg}\label{monomial}
Consid\'erons dans la carte affine $x_2=1$ la transformation
$f$ d\'efinie par
\begin{align*}
& f(x_0,x_1)=(x_0^ax_1^b,x_0^cx_1^d), && \left[\begin{array}{cc}
a & b\\
c & d\\
\end{array}
\right]\in\mathrm{GL}_2(\mathbb{Z})
\end{align*}
et $\mathcal{F}$ est associ\'e \`a l'une des $1$-formes $\alpha x_1\mathrm{d}x_0+
\beta x_0\mathrm{d}x_1$ o\`u
$(\alpha,\beta)$ est un vecteur propre de la matrice $\left[\begin{array}{cc}
a & b\\
c & d\\
\end{array}
\right].$
\end{eg}

Soient $\mathcal{F}$ un feuilletage sur $\mathbb{P}^2(\mathbb{C})$ 
et $\pi\hspace{0.1cm}\colon\hspace{0.1cm}S\dashrightarrow\mathbb{P}^2(
\mathbb{C})$ un \textbf{\textit{mod\`ele birationnel}}\label{ind32}
de $\mathbb{P}^2(\mathbb{C})$ ({\it i.e.}~$S$ est une surface et $\pi$ un 
morphisme birationnel). On dit
alors que $\pi^*\mathcal{F},$ not\'e encore $\mathcal{F}$ par abus de langage, est un
mod\`ele birationnel de $\mathcal{F};$ dans cette situation on introduit les
groupes

\begin{small}
\begin{align*}
&\mathrm{Aut}(S,\mathcal{F})=\{f\in\mathrm{Aut}(S)\hspace{0.1cm}|\hspace{0.1cm} f^*\pi^*\mathcal{F}=\pi^*\mathcal{F}\}
&&\text{et}
&&\mathrm{Bir}(S,\mathcal{F})=\{f\in\mathrm{Bir}(S)\hspace{0.1cm}|\hspace{0.1cm} f^*\pi^*\mathcal{F}=\pi^*\mathcal{F}\}.
\end{align*}
\end{small}

\noindent \textsc{Cantat} et \textsc{Favre} d\'emontrent alors le:

\begin{thm}[\cite{CaFa}, Th\'eor\`eme 1.2]
{\sl Soit $\mathcal{F}$ un feuilletage sur $\mathbb{P}^2(\mathbb{C})$ tel que l'inclusion $\mathrm{Aut}(S,\mathcal{F})
\subset\mathrm{Bir}(S,\mathcal{F})$ soit stricte pour tout mod\`ele birationnel de $\mathcal{F}.$
Alors $\mathrm{Bir}(S,\mathcal{F})$ poss\`ede un \'el\'ement birationnel d'ordre
infini et
\begin{itemize}
\item soit $\mathcal{F}$ est une fibration rationnelle;

\item soit $\mathcal{F}$ est birationnellement conjugu\'e \`a
l'exemple \ref{monomial} ou bien \`a l'un de ses quotients.
\end{itemize}}
\end{thm}

Dit autrement si $\mathcal{F}$ d\'esigne un feuilletage sur $\mathbb{P}^2(\mathbb{C})$
non birationnellement conjugu\'e \`a l'un des exemples pr\'ec\'edents
alors ou bien $\mathcal{F}$ est une fibration rationnelle, ou bien il y a un mod\`ele
birationnel de $\mathcal{F}$ pour lequel $\mathrm{Aut}(S,\mathcal{F})=\mathrm{Bir}(S,\mathcal{F}).$

\subsection{Propri\'et\'es des germes de flots dans
$\mathrm{Bir}_2$}\hspace{0.1cm}

On appelle \textbf{\textit{germe de flot}}\label{ind33} dans $\mathrm{Bir}_2$ un germe d'application holomorphe
$t\mapsto\phi_t\in \mathrm{Bir}_2$ satisfaisant
\begin{align*}
\phi_{t+s}^\bullet=\phi_t^\bullet\phi_s^\bullet, && \phi_0^\bullet=\mathrm{id}.
\end{align*}

Pour abr\'eger nous parlerons de flot, ceci sera
justifi\'e par la classification qui montrera qu'un germe de flot
se globalise en un certain sens. Comme toujours lorsqu'il n'y
aura pas d'ambigu\"{\i}t\'e nous identifierons les notations
$\phi_t$ et $\phi_t^\bullet.$

L'ensemble des droites contract\'ees par le germe de
flot $\phi_t$ forme un germe d'ensemble analytique dans la
Grassmanienne des droites de $\mathbb{P}^2(\mathbb{C}),$ {\it i.e.} dans l'espace
dual $\check{\mathbb{P}}^2(\mathbb{C}).$ De m\^eme l'ensemble des points
d'ind\'etermination des $\phi_t^\bullet$ constitue un germe d'ensemble
analytique cette fois dans $.$ Appelons
\textbf{\textit{famille de droites contract\'ees}}\label{ind34} une application continue (en fait
analytique) d\'efinie sur un germe de secteur ferm\'e $\Delta$ de
sommet $0$ dans $\mathbb{C},$
\begin{align*}
\mathcal{D}\hspace{0.1cm}\colon\hspace{0.1cm}\Delta\to\check{\mathbb{P}}^2(\mathbb{C}),
\end{align*}
telle que pour chaque $t$ dans $\Delta$ la droite $\mathcal{D}_t$
co\"{\i}ncide avec une droite $\mathcal{D}(t)$ contract\'ee par
$\phi_t^\bullet.$ Cette d\'efinition est rendue pertinente par le
th\'eor\`eme de \textsc{Puiseux}. De m\^eme par \textbf{\textit{famille de
points d'ind\'etermination}}\label{ind35} on entendra une
application $t\mapsto m_t$ continue d\'efinie sur un secteur~$\Delta,$ chaque $m_t$ \'etant d'ind\'etermination pour $\phi_t^\bullet.$

\begin{defi}
Soient $\phi_t$ un flot et $\mathcal{D}_t$ une famille de
droites contract\'ees par $\phi_t$. Si $\mathcal{D}_t$ est
ind\'ependante de $t$, la famille est dite
\textbf{\textit{immobile}}\label{ind36}, sinon elle est dite
\textbf{\textit{mobile}}\label{ind37}.

On a une notion analogue pour les familles de points
d'ind\'etermination.
\end{defi}

\begin{defi}
Soit $\chi$ un champ de vecteurs rationnel sur $\mathbb{P}^2(\mathbb{C});$ on dira
que $\chi$ est \textbf{\textit{rationnellement
int\'egrable}}\label{ind38}
si son flot $\phi_t$ est un flot de transformations birationnelles.
\end{defi}

Un germe de flot $\phi_t$ dans $\mathrm{Bir}_2$ est le flot
d'un champ rationnellement int\'egrable $\chi=\frac{\partial
\phi_t^\bullet}{\partial t}\Big|_{t=0}$ appel\'e
\textbf{\textit{g\'en\'erateur infinit\'esimal}}\label{ind39} de
$\phi_t.$

\begin{pro}\label{fib}
{\sl Soit $\phi_t$ un germe de flot dans $\mathrm{Bir}_2$. Supposons
que $\phi_t$ contracte une droite mobile; alors $\phi_t$
laisse une fibration en droites invariante.}
\end{pro}

\begin{proof}[\sl D\'emonstration]
Notons $\mathcal{D}_t$ la droite contract\'ee par $\phi_t$ et
$m_t$ son image, $\phi_t(\mathcal{D}_t)=~m_t.$

\begin{itemize}
\item Commen\c{c}ons par montrer que ou
bien $m_t$ est immobile, ou bien $\phi_t$ laisse une fibration
invariante fibre \`a fibre. \'Ecrivons que $\phi_t$ et
$\phi_s$ commutent
\begin{align*}
\phi_t^\bullet(\phi_s^\bullet(\mathcal{D}_t))=\phi_s^\bullet(\phi_t^\bullet
(\mathcal{D}_t))=\phi_s^\bullet(m_t).
\end{align*}
Comme $\phi_s^\bullet$ contracte un
nombre fini de courbes, g\'en\'eriquement sur $t,$ \`a $s$ fix\'e,
$\phi_s^\bullet(\mathcal{D}_t)$ est une courbe. Si cette courbe est
contract\'ee par $\phi_t^\bullet$ on aura $\phi_s^\bullet(\mathcal{D}_t)=
\mathcal{D}_t$ (toujours parce que~$\mathrm{Exc}\hspace{0.1cm}
\phi_t^\bullet$ est fini et $\phi_0^\bullet=\mathrm{id}$). Ainsi tous les $\phi_s^\bullet$
laissent $\mathcal{D}_t$ invariante; $\mathcal{D}_t$ \'etant
mobile,~$\phi_s^\bullet$ laisse une fibration en droites (dont font
partie les $\mathcal{D}_t$) invariante fibre \`a fibre comme nous
le pr\'eciserons dans le point qui suit.

Si $\mathcal{D}_t$ n'est pas invariante par $\phi_s^\bullet,$
alors $\phi_s^\bullet(\mathcal{D}_t)$ n'est pas contract\'ee par
$\phi_t^\bullet$ (toujours parce que $\mathrm{Exc}\hspace{0.1cm}
\phi_t^\bullet$ est
fini). Il en r\'esulte que $m_t$ appartient \`a $\mathrm{Ind}\hspace{0.1cm}
\phi_s^\bullet;$ la transformation $\phi_s^\bullet$ ayant un nombre fini
de points d'ind\'etermination $m_t$ est immobile.

\bigskip

\item Montrons que si $m_t$ est immobile
$\phi_t^\bullet$ pr\'eserve une fibration en droites; $m_t$ \'etant
immobile, nous le noterons $m.$ Puisque $\phi_s^\bullet$ et $\phi_t^\bullet$
commutent on a
\begin{align*}
\phi_s^\bullet(\mathcal{D}_t)=\phi_s^\bullet(\phi_{-t}^\bullet(m))=\phi_{-t}^\bullet(
\phi_s^\bullet(m))=\phi_{-t}^\bullet(\mathcal{D}_{s})
\end{align*}
que l'on peut r\'e\'ecrire sous la forme $\phi_{t+s}^\bullet(\mathcal{D}_t)=
\mathcal{D}_s.$ Par suite chaque $\phi_t^\bullet$ pr\'eserve la famille
infinie de droites $\mathcal{D}=(\mathcal{D}_s).$

\bigskip

\begin{itemize}
\item Supposons que $\phi_t$ appartienne \`a $\Sigma^1$ pour $t$
g\'en\'erique; alors $\phi_t$ s'\'ecrit $A\tau B$ o\`u $A$,~$B$
d\'esignent des \'el\'ements de $\mathrm{PGL}_3(\mathbb{C})$ (d\'ependants
de $t$). Un calcul direct montre qu'un ensemble de droites dont
l'image par $\tau$ est encore un ensemble de droites fait partie
du pinceau de droites passant par $(0:0:1).$ Il en r\'esulte que
$\mathcal{D}$ est contenue dans ce pinceau de droites qui est donc
invariant par le flot.

\bigskip

\item Si $\phi_t$ est dans $\Sigma^2$, alors, pour $t$
g\'en\'erique, $\phi_t$ est de la forme $A\rho B$ avec $A$, $B$
deux automorphismes de $\mathbb{P}^2(\mathbb{C}).$ Comme pr\'ec\'edemment on constate
qu'un ensemble de droites dont l'ima\-ge par $\rho$ est encore un
ensemble de droites est inclus dans les deux pinceaux aux points
d'ind\'etermination $(1:0:0)$ et $(0:1:0)$ de $\rho.$ Ici encore
$\mathcal{D}$ est contenue dans un pinceau de droites
invariant par le flot. \bigskip

\item Lorsque $\phi_t$ admet trois points d'ind\'etermination, il
est du type $A\sigma B$. Les droites dont l'image par l'involution
de \textsc{Cremona} sont encore des droites font partie des trois
pinceaux de droites passant par les points $(1:0:0)$, $(0:1:0)$ et
$(0:0:1).$ Il s'en suit que $\mathcal{D}$ est contenue dans
l'union de ces trois pinceaux de droites. Le flot $\phi_t$ laisse donc
ici encore un pinceau de droites invariant.
\end{itemize}
\end{itemize}
\end{proof}

\begin{eg}
Soit $f$ un \'el\'ement du groupe ab\'elien
\begin{align*}
\mathrm{G}=\{(x_0(x_2+\beta x_1):x_1(x_2+\alpha x_0): (x_2+\beta x_1)(x_2+\alpha
x_0))\hspace{0.1cm}|\hspace{0.1cm}\alpha,\hspace{0.1cm}\beta\in\mathbb{C}\};
\end{align*}
on a $\mathrm{Exc}\hspace{0.1cm} f=\{x_2=0,\hspace{0.1cm} x_2+\beta x_1=0,\hspace{0.1cm} x_2+\alpha x_0=0\}$. La
droite mobile d'\'equation $x_2+\beta x_1=0$
(resp. $x_2+\alpha x_0=0$) est contract\'ee sur le point immobile
$(0:1:0)$ (resp. $(0:0:1)$). Dans la carte affine $x_1=1$, un
\'el\'ement $f$ de $\mathrm{G}$ est du type
\begin{align*}
\left(\frac{x_0(x_2+\beta)}{x_2+\alpha x_0},x_2+\beta\right);
\end{align*}
le groupe $\mathrm{G}$ laisse donc la fibration $x_2=$ cte invariante.

Dans le groupe $\mathrm{G}$ ci-dessus il y a beaucoup de flots par
exemple
\begin{align*}
\phi_t=(x_0(x_2+tx_1):x_1(x_2+tx_0):(x_2+tx_1)(x_2+tx_0)).
\end{align*}
Ce flot $\phi_t$ contracte les deux familles mobiles $x_2+tx_0=0$
et $x_2+tx_1=0.$ Il laisse ainsi les deux fibrations $x_2/x_0$ et
$x_2/x_1$ invariantes avec permutation des fibres.
\end{eg}

\begin{lem}
{\sl Soit $\phi_t$ un germe de flot dans $\mathrm{Bir}_2.$ Un point
mobile \'eclat\'e par $\phi_t$ l'est sur une droite immobile.}
\end{lem}

\begin{proof}[\sl D\'emonstration]
Notons $m_t$ le point mobile \'eclat\'e par $\phi_t^\bullet$ et
$\mathcal{D}_t$ la droite sur laquelle il est \'eclat\'e.
\'Ecrivons la commutation de $\phi_t^\bullet$ et $\phi_s^\bullet$
\begin{align*}
\phi_t^\bullet(\phi_s^\bullet(m_t))=\phi_s^\bullet(\phi_t^\bullet(m_t))=
\phi_s^\bullet(\mathcal{D}_t).
\end{align*}
Puisque $\mathrm{Ind}\hspace{0.1cm}\phi_s^\bullet$ est
fini et $m_t$ mobile, l'image de $m_t$ par $\phi_s^\bullet$ est un
point $p_s$ pour $s$, $t$ g\'en\'eriques. Si~$\mathcal{D}_t$ est
mobile, $\phi_s^\bullet(\mathcal{D}_t)$ est en g\'en\'eral une courbe
(car $\mathrm{Exc}\hspace{0.1cm}\phi_s^\bullet$ est fini); alors l'\'egalit\'e
\begin{align*}
\phi_t^\bullet(p_s)=~\phi_s^\bullet(\mathcal{D}_t)
\end{align*}
assure que $p_s$ est un
point d'ind\'etermination de $\phi_t^\bullet$ ce qui contredit le fait
que~$\phi_t^\bullet$ ne compte qu'un nombre fini de points
d'ind\'etermination.
\end{proof}

\begin{eg}
Soit $\mathrm{G}$ le sous-groupe ab\'elien de $\mathrm{Bir}(\mathbb{P}^2(\mathbb{C}))$ d\'efini par
\begin{align*}
\mathrm{G}=\{((x_0+\alpha x_2)^2:x_0x_1 +\beta x_2^2:x_2(x_0+\alpha x_2))\hspace{0.1cm}|\hspace{0.1cm}\alpha,
\hspace{0.1cm}\beta\in\mathbb{C}\}.
\end{align*}
Chaque transformation $((x_0+\alpha x_2)^2:x_0x_1+\beta x_2^2:
x_2(x_0+\alpha x_2))$ contracte le point
$(-\alpha^2:\beta:\alpha)$ qui est donc mobile sur la droite d'\'equation
$x_0=0$ qui est fixe.
\end{eg}

\begin{lem}\label{mobimob}
{\sl Soit $\phi_t$ un germe de flot dans $\mathrm{Bir}_2.$ Il y a au
plus une droite immobile contract\'ee par $\phi_t.$}
\end{lem}

\begin{proof}[\sl D\'emonstration]
\'Ecrivons $\phi_t$ sous la forme $(Q_{1,t}:Q_{2,t}: Q_{3,t})$
o\`u les $Q_{i,t}$ sont des formes quadratiques d\'ependant
analytiquement de $t.$ Comme $\phi_0^\bullet=\mathrm{id},$ les
$Q_{i,0}$ sont divisibles par une m\^eme forme lin\'eaire que l'on
peut supposer \^etre $x_2,$ {\it i.e.} on peut supposer que
$\phi_t$ est du type
\begin{align*}
(x_0x_2+tA_1(x_0,x_1,x_2,t):x_1x_2+tA_2(x_0,x_1,x_2,t):
x_2^2+tA_3(x_0,x_1,x_2,t)).
\end{align*}
On constate en coordonn\'ees homog\`enes que
lorsque $t=0$, l'unique droite contract\'ee est $x_2=0$. Il suffit,
pour terminer, de remarquer que si une droite $\mathcal{D}$ est
contract\'ee pour une infinit\'e de $t$ elle l'est par tous les
$\phi_t^\bullet.$
\end{proof}

\begin{rems}
\begin{itemize}
\item On peut avoir cohabitation entre droites mobiles et droites
immobiles contract\'ees. Consid\'erons le flot
\begin{align*}
\phi_t=(x_0x_2:x_1(tx_0+x_2):x_2(tx_0+x_2));
\end{align*}
on constate que $\phi_t^\bullet$ contracte les droites
$x_2=0$ (immobile) et $tx_0+x_2=0$ (mobile), ces derni\`eres sur
le point $(1:0:0).$

\item On pourrait penser que le Lemme
\ref{mobimob} se g\'en\'eralise en
degr\'e quelconque, \`a savoir qu'un flot birationnel a au plus
une courbe rationnelle immobile contract\'ee. Il n'en est rien
comme le montre l'exemple qui suit. Soit $e^{tA}$ un flot
lin\'eaire g\'en\'erique; alors $\sigma e^{tA}\sigma$ est un flot
de transformations de degr\'e $4$ ayant trois droites fixes et
trois coniques mobiles contract\'ees.
\end{itemize}
\end{rems}

\begin{defi}
Soit $\phi_t$ un flot de $\mathrm{Bir}_2;$ nous dirons que $\phi_t$
est un \textbf{\textit{flot polynomial}}\label{ind40}
s'il existe une carte affine $\mathbb{C}^2$ invariante telle que
$\phi_{t_{|\mathbb{C}^2}}^\bullet \hspace{0.1cm}\colon\hspace{0.1cm}
\mathbb{C}^2\to\mathbb{C}^2$ soit polynomial pour
chaque $t$.
\end{defi}

\begin{pro}\label{pol}
{\sl Soit $\phi_t$ un germe de flot dans $\mathrm{Bir}_2$. Si
$\phi_t$ contracte une unique droite qui de plus est immobile,
$\phi_t$ est un germe de flot polynomial.}
\end{pro}

\begin{proof}[\sl D\'emonstration]
S'il s'agit d'un flot lin\'eaire il n'y a rien \`a faire.
On peut supposer \`a conjugaison pr\`es que l'unique droite
contract\'ee par $\phi_t$ est la droite d'\'equation $x_0=0.$
Pour chaque~$t$ tel que $\phi_t^\bullet$ soit vraiment quadratique
$\phi_t^\bullet$ est g.d. conjugu\'e \`a $\tau,$ {\it i.e.} pour ces
$t$ le flot $\phi_t$ s'\'ecrit 
\begin{align*}
A_t(x_0^2:x_0x_1:x_1^2-x_0x_2)B_t
\end{align*}
o\`u $A_t$, $B_t$ d\'esignent des \'el\'ements de $\mathrm{PGL}_3(\mathbb{C})$ tels que
$B_t$ pr\'eserve la droite d'\'equation $x_0=0$. Par suite l'unique
point d'ind\'etermination de
\begin{align*}
\phi_t=A_t(x_0^2:x_0x_1:x_1^2-x_0x_2)B_t
\end{align*}
est situ\'e sur la
droite d'\'equation $x_0=0.$ Ainsi les points d'ind\'etermination et
les courbes contract\'ees sont \`a l'infini dans la carte $x_0=1$
pour presque tout $t$ donc pour tout $t,$ ce qui signifie que
$\phi_t^\bullet$ y est polynomial.
\end{proof}

Le groupe $\mathrm{Aut}[\mathbb{C}^2]$ des automorphismes polynomiaux de
$\mathbb{C}^2$ contient le groupe, dit \textbf{\textit{groupe
affine}}\label{ind41a},
\begin{align*}
\mathrm{A}=\{(x_0,x_1)\mapsto(a_1x_0+b_1x_1+
c_1,a_2x_0+b_2x_1+c_2)\hspace{0.1cm}|\hspace{0.1cm} a_i,\hspace{0.1cm} b_i,\hspace{0.1cm} c_i\in\mathbb{C},\hspace{0.1cm}
a_1b_2-a_2b_1 \not=0\}
\end{align*}\label{not25a}
des automorphismes affines et le groupe des automorphismes
pr\'eservant la fibration $x_1=$ cte, encore appel\'e
\textbf{\textit{groupe \'el\'ementaire}}\label{ind41b}
\begin{align*}
\mathrm{E}=\{(x_0,x_1)\mapsto(\alpha x_0+P(x_1),\beta x_1+
\gamma)\hspace{0.1cm}|\hspace{0.1cm}\alpha,\hspace{0.1cm} \beta\in\mathbb{C}^*,\hspace{0.1cm}\gamma\in\mathbb{C},\hspace{0.1cm} P\in
\mathbb{C}[x_1]\}.
\end{align*}\label{not25b}
Plus pr\'ecis\'ement \textsc{Jung} a montr\'e que
$\mathrm{Aut}[\mathbb{C}^2]$ est le produit amalgam\'e de $\mathrm{E}$ et
$\mathrm{A}$ le long de leur intersection (\cite{Ju, La2}). Le
groupe des automorphismes polynomiaux a fait l'objet de nombreuses
\'etudes. En particulier \textsc{Friedland} et \textsc{Milnor}
(\emph{voir} \cite{FrMi}) ont montr\'e qu'un \'el\'ement~$f$ de
$\mathrm{Aut}[\mathbb{C}^2]$ est ou bien conjugu\'e \`a un automorphisme
\'el\'ementaire, ou bien conjugu\'e \`a un automorphisme
s'\'ecrivant comme compos\'ee d'un nombre fini de transformations
de \textsc{H\'enon}; dans ce dernier cas on dit que $f$ est
\textbf{\textit{de type \textsc{H\'enon}}}\label{ind40a}.

\bigskip

La classification des groupes ab\'eliens de
transformations polynomiales a \'et\'e \'etablie par
\textsc{Moldavandski} (\cite{Mo, Wr}). \`A conjugaison pr\`es dans
le groupe des automorphismes polynomiaux, un tel groupe
$\mathrm{H}$ v\'erifie l'une des assertions suivantes
\begin{itemize}
\item $\mathrm{H}$ est contenu dans le groupe \'el\'ementaire
$\mathrm{E}$ ou le groupe affine $\mathrm{A};$

\item $\mathrm{H}=\cup_i\mathrm{H}_i$ o\`u les $\mathrm{H}_i,$
$i\in\mathbb{N},$ satisfont $\mathrm{H}_i\subset\mathrm{H}_{i+1}$
et sont individuellement conjugu\'es \`a des sous-groupes finis de
$\mathrm{A}\cap \mathrm{E};$

\item $\mathrm{G}=\mathrm{F}\times\langle f\rangle$ o\`u $f$
d\'esigne un automorphisme de type \textsc{H\'enon} et
$\mathrm{F}$ un sous-groupe fini de l'intersection $\mathrm{A}\cap\mathrm{E}.$
\end{itemize}

Parmi les trois types ci-dessus
seul le premier est susceptible de contenir un groupe \`a un
param\`etre, les deux derniers \'etant d\'enombrables.

Nous allons utiliser le r\'esultat suivant qui fait
partie du folklore.

\begin{lem}\label{elhe}
{\sl Soit $f$ un automorphisme polynomial de degr\'e $2$ de $\mathbb{C}^2$ qui
ne soit pas affine. Alors, \`a conjugaison affine pr\`es,
\begin{align*}
&f=(\alpha x_0+P(x_1),ax_1+b)&&\text{ou} && f=(x_1,P(x_1)-\delta
x_0);
\end{align*}
dit autrement, $f$ est, \`a conjugaison affine pr\`es, un
automorphisme \'el\'ementaire ou une transformation de
\textsc{H\'enon}.}
\end{lem}

\begin{proof}[\sl D\'emonstration]
Posons $f:=(P,Q)$ o\`u $P$ et $Q$ d\'esignent des polyn\^omes de
degr\'e inf\'erieur ou \'egal \`a $2$. \'Ecrivons $P$ (resp. $Q$)
sous la forme $P_0+P_1 +P_2$ (resp. $Q_0+Q_1+Q_2$) avec
$P_0,$ $Q_0$ constant et $\deg P_i=\deg Q_i=i$ pour $i\not=0.$

Puisque le d\'eterminant jacobien de $f$ est constant,
{\it i.e.} $dP\wedge dQ=$ cte, la quantit\'e $dP_2\wedge dQ_2$ est
nulle. Il en r\'esulte que $P_2=\varepsilon_1 q$ et
$Q_2=\varepsilon_2q$ o\`u les $\varepsilon_i$ d\'esignent des
complexes et $q$ une forme quadratique; notons que les
$\varepsilon_i$ sont non tous nuls sinon $f$ serait affine.
L'\'egalit\'e $dP\wedge dQ=$ cte implique alors que
$(\varepsilon_2dP_1-\varepsilon_1dQ_1) \wedge dq=0;$ la forme
quadratique $q$ et la forme lin\'eaire non triviale
$\varepsilon_2P_1- \varepsilon_1Q_1$ ont donc les m\^emes niveaux.
Par suite
\begin{align*}
&q=\varepsilon(\varepsilon_2P_1-\varepsilon_1Q_1)^2, && \varepsilon\in\mathbb{C}^*.
\end{align*}
On peut supposer que l'automorphisme $f$ est de la forme
\begin{align*}
(P_0+P_1+\varepsilon_1x_1^2,Q_0+Q_1+\varepsilon_2 x_1^2).
\end{align*}
Si $\varepsilon_2$ est nul, le d\'eterminant jacobien de $f$
\'etant constant, $Q_1$ ne d\'epend que de $x_1$ et $f$ est
\'el\'ementaire. Si $\varepsilon_2$ est non nul, alors, \`a
conjugaison lin\'eaire pr\`es par $(x_0+\varepsilon_1x_1/
\varepsilon_2,x_1),$ l'automorphisme $f$ s'\'ecrit
$(\widetilde{P}_0+\widetilde{P}_1,\widetilde{Q}_0+
\widetilde{Q}_1+\varepsilon_2 x_1^2).$ Puisque le d\'eterminant
jacobien de $f$ est constant, on obtient que $\widetilde{P}_1=
\widetilde{P}(x_1);$ \`a conjugaison lin\'eaire pr\`es $f$ est
alors un automorphisme de \textsc{H\'enon}.
\end{proof}

Soit $\phi_t$ un germe de flot polynomial quadratique non affine.
On choisit une valeur du param\`etre $t_0$ telle que
$\phi_{t_0}^\bullet=f$ soit purement quadratique. Puisque le
centralisateur d'un automorphisme de type \textsc{H\'enon} dans
$\mathrm{Aut}[\mathbb{C}^2]$ est d\'enombrable (\cite{La}), le Lemme
\ref{elhe} assure qu'\`a conjugaison affine pr\`es $f$ s'\'ecrit
de la fa\c{c}on suivante
\begin{align*}
&f=(\alpha x_0+\beta +x_1^2,ax_1+b), &&\alpha a\not=0.
\end{align*}

\begin{lem}\label{tous}
{\sl Soit $g$ un automorphisme quadratique commutant \`a
l'automorphisme polynomial
\begin{align*}
&f=(\alpha x_0+\beta +x_1^2,ax_1+b), &&\alpha a\not=0.
\end{align*}
Alors $g$ pr\'eserve la fibration $x_1=$ cte.}
\end{lem}

\begin{proof}[\sl D\'emonstration]
Posons $g=(g_1,g_2).$ La commutation de $g$ \`a $f$
nous assure en particulier que
\begin{equation}\label{commutation}
g_2\circ f=ag_2+b
\end{equation}
\'Ecrivons $g_2$ sous la forme $p_0+p_1x_0+p_2x_1+p_3x_0^2+p_4
x_0x_1+p_5x_1^2.$ En d\'eveloppant (\ref{commutation}) on obtient les
\'egalit\'es suivantes
\begin{align*}
p_3=p_4=(\alpha-a)p_1=bp_5=p_1+a(a-1)p_5=0.
\end{align*}
Comme $p_3=p_4=0$, la fibration $x_1=$ cte est
invariante par $g$ si $p_1=0;$ c'est le cas si $\alpha\not=a.$
Supposons $p_1\not=0$ alors $\alpha= a$ et $a(a-1)p_5\not=0;$ en
particulier $a\not=1$ et $b=0.$ Par suite $f$ s'\'ecrit 
\begin{align*}
(ax_0+\beta+x_1^2,ax_1).
\end{align*}
Comme $a$ est distinct de $1$, la transformation $f$
est affinement conjugu\'ee \`a $\widetilde{f}=(ax_0+x_1^2,ax_1).$ 
Un calcul direct montre que si $g$ commute \`a
$\widetilde{f}$, alors $g$ est du type suivant
\begin{align*}
((\mu^2+a(1-a)a_2)x_0+a_1x_1+a_2x_1^2,\mu x_1)
\end{align*}
et en particulier pr\'eserve la fibration $x_1=$ cte.
\end{proof}

\begin{pro}\label{flotaut}
{\sl Soit $\phi_t$ un germe de flot polynomial quadratique.
Alors $\phi_t$ laisse un pinceau de droites invariant. Le flot
$\phi_t$ est affine ou, \`a conjugaison lin\'eaire pr\`es, de
l'une des formes suivantes
\begin{align*}
&(x_0+t(\alpha+x_1^2),x_1);\\
& \left(\left(x_0+\frac{x_1^2}{\alpha}\right)e^{\alpha t}-\frac{x_1^2}
{\alpha},x_1\right),\hspace{0.1cm}\alpha\not=0; \\
&\left(x_0+t(p_0+p_1x_1+x_1^2)+t^2\left(\frac{p_1}{2}+x_1\right)+\frac{t^3}{3},x_1+t
\right);
\end{align*}

\begin{align*}
& \left(\left(x_0+\frac{2}{\alpha^3}+\frac{2x_1}{\alpha^2}+\frac{
x_1^2}{\alpha}\right)e^{\alpha t}-\left(\frac{2}{\alpha^3}
+\frac{2x_1}{\alpha^2}+\frac{x_1^2}{\alpha}\right)-\left(\frac{2}
{\alpha^2}+\frac{2x_1}{\alpha}\right)t-\frac{t^2}{\alpha},x_1+t\right),
\hspace{0.1cm}\alpha\not=0;\\
&\left(x_0+\varepsilon t+\frac{x_1^2}{2\mu}(e^{2\mu t}-1),x_1e^{\mu
t}\right),\hspace{0.1cm}\varepsilon\in\{0,1\},\hspace{0.1cm}\mu\not=0; \\
& \left(x_0e^{\alpha t}+\frac{x_1^2}{2\mu-\alpha}(e^{2\mu
t}-e^{\alpha t}),x_1e^{\mu t}\right),\hspace{0.1cm} \alpha\not=2\mu;\\
& \left(x_0e^{\alpha t}+\varepsilon x_1te^{\alpha t}+\frac{x_1^2}{
\alpha}(e^{2\alpha t}-e^{\alpha t}),x_1e^{\alpha t}\right),\hspace{0.1cm}
\alpha\not=0; \\
& \left(x_0e^{\alpha t}+x_1^2te^{\alpha t},x_1e^{\alpha t/2} \right).
\end{align*}}
\end{pro}

\begin{proof}[\sl D\'emonstration]
On suppose $\phi_t$ non affine. Puisque le centralisateur dans
$\mathrm{Aut}[\mathbb{C}^2]$ d'un automorphisme de \textsc{H\'enon} est
d\'enombrable (\cite{La}), $\phi_t$ est, d'apr\`es les Lemmes
\ref{elhe} et \ref{tous}, \`a conjugaison affine pr\`es, de la
forme suivante
\begin{align*}
\phi_t=(\alpha(t)x_0+P(x_1,t),\mu(t)x_1+
b(t))
\end{align*}
o\`u $\alpha,$ $\mu$ et $b$ sont holomorphes en $t$ et
$P(x_1,t)$ est un polyn\^ome du second degr\'e en $x_1$ \`a
param\`etre~$t.$ Le g\'en\'erateur infinit\'esimal $\chi=\frac{
\partial\phi_t}{\partial t}\big|_{t=0}$ est donc de la forme
\begin{align*}
\chi=(\alpha x_0+p_0+p_1x_1+p_2x_1^2)\frac{\partial}{\partial x_0}+(\mu
x_1+b)\frac{\partial}{\partial x_1}.
\end{align*}
Comme on peut exclure le cas
o\`u $\phi_t$ est affine pour tout $t,$ on suppose $p_2\not=0$
et par homoth\'etie on se ram\`ene \`a $p_2=1.$ Nous allons \`a
conjugaison affine pr\`es donner des formes normales des champs
$\chi$ ci-dessus suivant les valeurs des param\`etres $\alpha,$
$\mu,$ $\ldots$ \`A transformation affine pr\`es en~$x_1$ on peut
supposer que la seconde composante de $\chi$ est de l'un des trois
types $0,$ $\frac{\partial}{\partial x_1},$ $\mu x_1\frac{\partial}
{\partial x_1}.$\bigskip

\begin{itemize}

\item Supposons que $\chi=(\alpha
x_0+p_0+p_1x_1+x_1^2) \frac{\partial} {\partial x_0};$ si $\alpha$ est
nul $\phi_t$ est donn\'e, \`a conjugaison pr\`es par un
\'el\'ement du type $(x_0,x_1+\varepsilon),$ par
\begin{align*}
\phi_t=(x_0+t(\beta+x_1^2),x_1).
\end{align*}

Si $\alpha$ est non nul, on se ram\`ene, \`a
conjugaison pr\`es par des transformations du 
type $(x_0+a+bx_1,x_1),$ \`a $p_0=p_1=0.$ Dans ce cas le flot s'\'ecrit
\begin{align*}
\phi_t=\left(\left(x_0+\frac{x_1^2}{\alpha}\right)e^{\alpha
t}-\frac{x_1^2}{\alpha},x_1\right).
\end{align*}
\bigskip

\item Supposons que $\chi$ soit du type
$(\alpha x_0+p_0+ p_1x_1+x_1^2)\frac{\partial}{\partial
x}+\frac{\partial} {\partial x_1}~;$ comme pr\'ec\'edemment on
distingue les cas $\alpha$ nul et $\alpha$ non nul. Si
$\alpha=0,$ alors
\begin{align*}
\phi_t=\left(x_0+t(p_0+p_1x_1+x_1^2)+t^2\left(\frac{p_1}{2}
+x_1\right)+\frac{t^3}{3},x_1+t\right).
\end{align*}
Lorsque $\alpha\not=0$ on
peut de fa\c{c}on analogue au cas pr\'ec\'edent supposer
$p_0=p_1=0,$ ce qui conduit \`a int\'egrer $(\alpha
x_0+x_1^2)\frac{\partial}{\partial x_0}+\frac{\partial} {\partial x_1};$
on obtient

\begin{footnotesize}
\begin{align*}
\phi_t=\left(\left(x_0+\frac{2}
{\alpha^3}+\frac{2x_1}{\alpha^2}+\frac{x_1^2}{\alpha}\right)e^{\alpha
t}-\left(\frac{2}{\alpha^3}+\frac{2x_1}{\alpha^2}+\frac{x_1^2}{\alpha}
\right)-\left(\frac{2}{\alpha^2}+\frac{2x_1}{\alpha}\right)t-\frac{t^2}
{\alpha},x_1+t\right).
\end{align*}
\end{footnotesize}
\bigskip

\item On consid\`ere le cas o\`u
\begin{align*}
&\chi=(\alpha x_0+p_0 +p_1x_1+x_1^2)\frac{\partial}{\partial x_0}+\mu
x_1\frac{\partial} {\partial x_1},&&\mu\not=0.
\end{align*}
Si $\alpha=0$ le champ $\chi$ est lin\'eairement conjugu\'e \`a
$(\varepsilon+x_1^2)\frac{
\partial}{\partial x_0}+\mu x_1\frac{\partial}{\partial x_1};$ le flot
d'un tel champ est de la forme
\begin{align*}
\phi_t=\left(x_0+\varepsilon t+\frac{x_1^2}{2\mu}(e^{2\mu t}-1
),x_1e^{\mu t}\right).
\end{align*}
Lorsque $\alpha$ est non nul et distinct
de $\mu$ et $2\mu,$ alors \`a conjugaison lin\'eaire pr\`es $\chi$
s'\'ecrit $(\alpha x_0+x_1^2)\frac{\partial}{\partial x_0}+\mu x_1\frac{\partial}{
\partial x_1}$ et son flot est
\begin{align*}
\phi_t=\left(x_0e^{\alpha t}+\frac{x_1^2} {2\mu-\alpha}
(e^{2\mu t} -e^{\alpha t}),x_1e^{\mu t}\right).
\end{align*}
Pour $\alpha=\mu$ on constate que $\chi$ est
lin\'eairement conjugu\'e \`a $(\alpha x_0+\varepsilon x_1+x_1^2)
\frac{\partial}{\partial x_0}+\alpha x_1\frac{\partial}{\partial x_1}$
avec $\varepsilon$ dans $\{0,1\}$ qui conduit \`a
\begin{align*}
\phi_t=\left(x_0e^{\alpha t}+\varepsilon x_1te^{\alpha
t}+\frac{x_1^2}{\alpha}(e^{2\alpha t}-e^{ \alpha t}),x_1e^{\alpha
t}\right) .
\end{align*}
Enfin pour $\alpha=2\mu$ il vient \`a conjugaison lin\'eaire
pr\`es
\begin{align*}
\phi_t=\left(x_0 e^{\alpha t}+x_1^2te^{\alpha t},x_1
e^{\alpha t/2}\right).
\end{align*}
\end{itemize}
\end{proof}

D'apr\`es les \'enonc\'es \ref{fib}, \ref{mobimob}, \ref{pol} et
\ref{flotaut} on a le:

\begin{thm}\label{flot}
{\sl Un germe de flot dans $\mathrm{Bir}_2$ pr\'eserve une fibration en
droites.}
\end{thm}

\begin{defi}
Soit $\phi_t$ un flot birationnel quadratique de g\'en\'erateur
infinit\'esimal $\chi.$ Une \textbf{\textit{sym\'etrie
forte}}\label{ind41} $Y$ de $\chi$ est un champ de
vecteurs rationnellement int\'egrable de flot $\psi_s$ tel que
\begin{itemize}
\item $\phi_t^\bullet$ et $\psi_s^\bullet$ commutent, {\it i.e.} $[\chi,Y]=0;$

\item $\psi_s$ appartient \`a $\mathrm{Bir}_2$ pour tout $s;$

\item $\chi$ et $Y$ ne sont pas $\mathbb{C}$-colin\'eaires.
\end{itemize}
\end{defi}

Soient $\phi_t$ un flot dans $\mathrm{Bir}_2$ et $\chi$
(resp. $\mathcal{F}_\chi$) le champ (resp. feuilletage) associ\'e. On note
$\langle\phi_t^\bullet\rangle\subset~\mathrm{Bir}_2^\bullet$ le groupe engendr\'e
par les $\phi_t^\bullet$ et $\overline{\langle\phi_t^\bullet\rangle}^{\hspace{0.1cm}
\mathsf{Z}}$ son adh\'erence de \textsc{Zariski} dans
$\mathrm{Bir}_2.$ Enfin $\mathrm{G}(\chi)$ d\'esigne le groupe
ab\'elien alg\'ebrique maximal contenu dans $\mathrm{Bir}_2^\bullet$ et
contenant $\overline{\langle\phi_t^\bullet\rangle}^{\hspace{0.1cm}\mathsf{Z}}.$ Avec
les notations pr\'ec\'edentes on a le

\begin{thm}\label{symforte}
{\sl Soit $\phi_t$ un germe de flot dans $\mathrm{Bir}_2$ de
g\'en\'erateur infinit\'esimal $\chi.$ On a les assertions
suivantes
\begin{itemize}
\item[\textbf{\textit{i.}}] Si $\mathrm{G}(\chi)$ est de dimension $1$,
alors $\mathcal{F}_\chi$ est une fibration rationnelle.

\item[\textbf{\textit{ii.}}] Si $\dim\mathrm{G}(\chi)\geq 2,$ alors $\chi$
poss\`ede une sym\'etrie forte.
\end{itemize}

Dans les deux cas $\mathcal{F}_\chi$ est d\'efini par une forme
ferm\'ee rationnelle.}
\end{thm}

\begin{proof}[\sl D\'emonstration]
\textbf{\textit{i.}} Si $\dim\mathrm{G}(\chi)=1,$ alors
$\overline{ \langle \phi_t^\bullet\rangle}^{\hspace{0.1cm} \mathsf{Z}}$ est la
composante neutre de $\mathrm{G}(\chi).$ Ce groupe en tant que
groupe de \textsc{Lie} est isomorphe \`a $\mathbb{C},$ ou \`a $\mathbb{C}^*$ ou
encore \`a un tore $\mathbb{C}/\Lambda.$ D'apr\`es le Th\'eor\`eme
\ref{flot} le groupe $\overline{\langle \phi_t^\bullet\rangle}^{\hspace{0.1cm}
\mathsf{Z}}$ pr\'eserve une fibration en droites, disons $x_1=$
cte. Ceci induit un morphisme
\begin{align*}
\pi\hspace{0.1cm}\colon\hspace{0.1cm}\overline{\langle\phi_t^\bullet
\rangle}^{\hspace{0.1cm}\mathsf{Z}}\to\mathrm{PGL}_2(\mathbb{C})
\end{align*}
d\'ecrivant en particulier l'action de $\phi_t$ sur les fibres. Si
$\phi_t$ pr\'eserve la fibration $x_1=$ cte fibre \`a fibre ({\it
i.e.} $\pi$ est triviale), alors $\mathcal{F}_\chi=\{x_1=\text{cte}\}$ et
la conclusion est \'evidente. Sinon $\overline{\langle
\phi_t^\bullet\rangle}^{\hspace{0.1cm} \mathsf{Z}}$ ne peut \^etre un tore car il n'y
a pas de tore complexe dans les sous-groupes de
$\mathrm{PGL}_2(\mathbb{C});$ par suite l'adh\'erence (topologique) de
$\overline{ \langle \phi_t^\bullet\rangle}^{\hspace{0.1cm} \mathsf{Z}}$ dans
$\mathbb{P}^{17}(\mathbb{C})\simeq\mathrm{Rat}_2$ est une courbe rationnelle. On
sait depuis \textsc{Darboux} (\cite{Jo}) qu'un feuilletage de
$\mathbb{P}^2(\mathbb{C})$ dont toutes les feuilles sont d'adh\'erence des courbes
alg\'ebriques poss\`ede une int\'egrale premi\`ere rationnelle non
constante. Dans notre cas ces courbes sont rationnelles donc
$\mathcal{F}_\chi$ est une fibration rationnelle.

\bigskip

\textbf{\textit{ii.}} Supposons que $\dim\mathrm{G}(\chi)\geq 2.$
On peut alors trouver un germe de groupe \`a un param\`etre
$\psi_s$ dans $\mathrm{G}(\chi)$ non inclus dans
$\langle\phi_t^\bullet\rangle.$ Soit $Y$ le g\'en\'erateur
infinit\'esimal de $\psi_s;$ alors $\chi$ et $Y$ commutent et ne
sont pas $\mathbb{C}$-colin\'eaires. Pour terminer soit $\omega$ une
$1$-forme rationnelle d\'efinissant $\mathcal{F}_\chi;$ si $\chi$ et $Y$
sont g\'en\'eriquement ind\'ependants $\Omega= \omega/i_{Y}\omega$
est ferm\'ee et d\'efinit $\mathcal{F}_\chi.$ Si ce n'est pas le cas, $Y$
s'\'ecrit $r\chi$ avec $r$ rationnelle non constante et la
commutation implique $\chi(r)=0;$ alors $dr$ est ferm\'ee et
d\'efinit $\mathcal{F}_\chi.$
\end{proof}

\begin{rem}\label{sf}
Soit $Y$ une sym\'etrie forte de $\chi;$ le flot de $Y$ est
contenu dans $\mathrm{G}(\chi).$ Comme le groupe $\mathrm{G}(\chi)$ est un
groupe de \textsc{Lie}, $\alpha\chi+\beta Y$ est encore une
sym\'etrie forte de $\chi$ d\`es que $\beta$ est non nul.

En reprenant les d\'emonstrations menant au Th\'eor\`eme
\ref{flot}, on constate que $\mathrm{G}(\chi)$ pr\'eserve une fibration en
droites, disons $x_1=$ cte; un \'el\'ement $\phi=(\phi_1,
\phi_2)$ de $\mathrm{G}(\chi)$ est donc du type (en carte affine)
\begin{align*}
\phi(x_0,x_1)=\left(\frac{a(x_1)x_0+b(x_1)}{c(x_1)x_0+d(x_1)},
\nu(x_1)\right)
\end{align*}
o\`u $\nu$ appartient \`a $\mathrm{PGL}_2(\mathbb{C})$ et $a,$ $b,$ $c$ et $d$
sont des polyn\^omes en $x_1.$ Dans la carte affine $x_2=1$ on peut
donc \'ecrire les champs $\chi$ et $Y$ sous l'une des formes
\begin{align*}
&\chi=\chi_1\frac{\partial}{\partial
x_0}+\chi_2(x_1)\frac{\partial} {\partial x_1},&&
Y=Y_1\frac{\partial}{\partial x_0}+Y_2(x_1) \frac{\partial}
{\partial x_1}, && \text{avec } \chi_2(x_1)=0,\hspace{0.1cm} 1 \text{ ou }
\alpha x_1.
\end{align*}
Puisque $\chi$ et $Y$ commutent quitte \`a changer $Y_1$ en
$Y_1-\varepsilon\chi$ on peut supposer lorsque $\chi_2\not\equiv
0$ que~$Y_2\equiv 0.$
\end{rem}

\subsection{Deux r\'esultats sur les germes de flots dans
$\mathrm{Bir}_n$}\hspace{0.1cm}

Pour un germe de flot quelconque, {\it i.e.} dans
$\mathrm{Bir}_n$ (la d\'efinition \'etant la m\^eme que pour un flot dans
$\mathrm{Bir}_2$), le Th\'eor\`eme \ref{symforte} se g\'en\'eralise de
la fa\c{c}on suivante.

\begin{thm}\label{symforte2}
{\sl Soit $\phi_t$ un germe de flot dans $\mathrm{Bir}_n$ de
g\'en\'erateur infinit\'esimal $\chi.$ Soit $\mathrm{G}(\chi)$ le
groupe ab\'elien maximal alg\'ebrique contenu dans
$\mathrm{Bir}_n$ et contenant $\overline{\langle
\phi_t^\bullet\rangle}^{\hspace{0.1cm} \mathsf{Z}}.$ On a
\begin{itemize}
\item Si $\mathrm{G}(\chi)$ est de dimension $1$,
alors $\mathcal{F}_\chi$ est une fibration rationnelle ou elliptique.

\item Si $\dim\mathrm{G}(\chi)\geq 2,$ alors $\chi$
poss\`ede une sym\'etrie forte.
\end{itemize}

Dans les deux cas $\mathcal{F}_\chi$ est d\'efini par une forme
ferm\'ee rationnelle.}
\end{thm}

\begin{proof}[\sl D\'emonstration]
La seconde assertion se d\'emontre comme au Th\'eor\`eme
\ref{symforte}. Dans la premi\`ere situation les feuilles
de $\mathcal{F}_\chi$ sont d'adh\'erence (ordinaire) alg\'ebriques et $\mathcal{F}_\chi$
poss\`ede une int\'egrale premi\`ere rationnelle; ses fibres sont
n\'ecessairement rationnelles ou elliptiques.
\end{proof}

Le Th\'eor\`eme \ref{flot} se g\'en\'eralise comme suit.

\begin{thm}
{\sl Tout germe de flot birationnel laisse une fibration rationnelle
invariante.}
\end{thm}

\begin{proof}[\sl D\'emonstration]
Soient $\phi_t$ un germe de flot birationnel et $\mathcal{F}$ le feuilletage
associ\'e au g\'en\'erateur infinit\'esimal de $\phi_t;$ on a
$\phi_t^*\mathcal{F}=\mathcal{F}.$ Les transformations de type $(x_0^ax_1^b,x_0^cx_1^d)$ ne
se plongent pas dans un flot; en effet par it\'eration leur degr\'e
cro\^it. Par suite, d'apr\`es le Th\'eor\`eme~1.2. de \cite{CaFa}
(rappel\'e au Chapitre \ref{germgerm}, \S\ref{cantatfavre}), ou bien
$\mathcal{F}$ est une fibration
rationnelle, ou bien $\mathrm{Bir}(S,\mathcal{F})=\mathrm{Aut}(S,\mathcal{F})$ pour un certain mod\`ele
birationnel de $\mathcal{F}.$

Lorsque $\mathrm{Bir}(S,\mathcal{F})=\mathrm{Aut}(S,\mathcal{F})$ les $\phi_t$ induisent
un groupe \`a un param\`etre global (le flot du g\'en\'erateur infinit\'esimal
vu sur $S$) d'automorphismes de $S$ qui bien s\^ur sont isotopes \`a
l'identit\'e. Par suite ils fixent tous les diviseurs exceptionnels de
$S$ et on peut donc supposer la surface $S$ minimale, {\it i.e.} $\mathbb{P}^1(\mathbb{C})\times
\mathbb{P}^1(\mathbb{C})$ ou $\mathbb{P}^2(\mathbb{C})$ ou une surface de
\textsc{Hirzebruch}. Un calcul \'el\'ementaire - au cas par cas -
montre qu'un flot d'automorphismes sur une surface minimale laisse une
fibration rationnelle invariante.
\end{proof}

Dans la suite les $^\bullet$ seront la plupart du temps omis.

\section{Classification des germes de flots birationnels
quadratiques}\label{cla}\hspace{0.1cm}

Soit $\phi_t$ un germe de flot dans $\mathrm{Bir}_2;$ alors
$\phi_t$ s'\'ecrit \`a conjugaison lin\'eaire pr\`es
\begin{align*}
\phi_t=\left(\frac{A(x_1,t)x_0+B(x_1,t)}{C(x_1,t)x_0+D(x_1,t)},
\nu(x_1,t)\right)
\end{align*}
avec
\begin{itemize}
\item $\nu(x_1,t)=x_1$ ou $x_1+t$ ou $e^{\beta t}x_1;$

\item $A,$ $B,$ $C$ et $D$ sont polynomiaux en $x_1$ et
\begin{align*}
&\deg_{x_1} A\leq 1, &&\deg_{x_1} B\leq 2, &&\deg_{x_1} C=0,
\textit{i.e. } C(x_1,t)=C(t), && \deg_{x_1} D\leq 1.
\end{align*}

\item $B(x_1,0)=C(x_1,0)=0$ et $A(x_1,0)=D(x_1,0).$
\end{itemize}

\bigskip

Le g\'en\'erateur infinit\'esimal $\chi:=\frac{\partial
\phi_t} {\partial t}\Big|_{t=0}$ de $\phi_t$ est de la
forme
\begin{align*}
\frac{\alpha x_0^2+\ell(x_1)x_0+P(x_1)}{ax_1+b}\frac{\partial}
{\partial x_0}+\varepsilon(x_1) \frac{\partial}{\partial x_1}
\end{align*}
o\`u
\begin{itemize}
\item $\alpha$, $a$, $b$ d\'esignent des \'el\'ements de $\mathbb{C};$

\item $\ell$ et $P$ des polyn\^omes de degr\'e $1$ et $2$
respectivement;

\item \`a conjugaison lin\'eaire pr\`es et multiplication scalaire
pr\`es (renormalisation du temps), $\varepsilon$ vaut $0,$ $1$ ou
$x_1$.
\end{itemize}

Nous allons classifier, \`a transformations lin\'eaires
et renormalisation pr\`es, les flots de $\mathrm{Bir}_2,$ ce qui revient
peu ou prou \`a d\'etecter parmi les champs $\chi$ ci-dessus ceux
qui sont rationnellement int\'egrables.

L'approche se fait en
trois \'etapes suivant les valeurs de $\varepsilon.$

\subsection{La seconde composante du flot est triviale}\hspace{0.1cm}

Dans ce paragraphe on s'int\'eresse aux champs du type
\begin{align*}
\frac{\alpha x_0^2+\ell(x_1)x_0+P(x_1)}{ax_1+b}\frac{\partial} {\partial
x_0}.
\end{align*}
Nous allons nous ramener \`a \'etudier $5$ mod\`eles de champs. \`A homoth\'etie
en $x_0$ pr\`es, on a l'alternative: $\alpha=0$ ou $\alpha=1.$

\bigskip

\begin{itemize}
\item Lorsque $\alpha$ vaut $1,$ on
constate qu'\`a transformation lin\'eaire pr\`es on peut supposer
$\ell\equiv 0.$ De plus quitte \`a faire une translation et/ou une
homoth\'etie en $x_1$ le couple $(a,b)$ prend les valeurs suivantes:
$(0,1)$ ou $(1,0).$ Autrement dit on se ram\`ene aux deux familles
de champs \begin{align*} \chi=(x_0^2+P(x_1))\frac{\partial}{\partial
x_0}&&\text{ou} &&\chi=\frac{x_0^2+P(x_1)}{x_1} \frac{\partial}{\partial
x_0}.
\end{align*}

\bigskip

\item Lorsque $\alpha$ est nul, il apparait trois types suivant en
particulier les valeurs de $a$

\begin{itemize}
\item  Si $a$ est non nul, alors, \`a
translation et homoth\'etie en $x_1$ pr\`es,
\begin{align*}
\chi=\frac{(\beta x_1+\mu)x_0+P(x_1)} {x_1}\frac{\partial}{\partial x_0}.
\end{align*}

\bigskip

\item Lorsque $a$ est nul, on se ram\`ene aux deux cas suivants
\begin{align*}
&\chi=(\mu x_0+P(x_1))\frac{\partial}{\partial x_0} &&\text{ou}&&
\chi=(\beta x_0x_1+P(x_1))\frac{\partial} {\partial x_0}.
\end{align*}
\end{itemize}
\end{itemize}

Nous allons dans les Remarques et Lemmes qui suivent
\'etudier ces diff\'erentes \'eventualit\'es.

\begin{rem}\label{frou}
Soit $\chi$ le champ d\'efini par $(\mu x_0+P(x_1))\frac{\partial}
{\partial x_0}$ o\`u $\mu$ d\'esigne un \'el\'ement de~$\mathbb{C}$ et $P$
un polyn\^ome quadratique. Si $\mu$ est non nul, $\chi$ est
lin\'eairement conjugu\'e \`a $(\mu x_0+x_1^2)\frac{\partial}{\partial
x_0}$ dont le flot s'\'ecrit
\begin{align*}
\left(x_0e^{\mu t}+\frac{x_1^2}{\mu}\left( e^{\mu t}-1\right),x_1\right).
\end{align*}
Si $\mu$ est nul, le flot de
$\chi$ est lin\'eairement conjugu\'e \`a $(x_0+t(\alpha+x_1^2),x_1)$
(\emph{voir} Proposition \ref{flotaut}).

Tous les deux sont des flots \og polynomiaux \fg.
\end{rem}

\begin{lem}\label{froufrou}
{\sl Soit $P$ un \'el\'ement de $\mathbb{C}[x_1]$ de degr\'e inf\'erieur ou
\'egal \`a $2.$ Les champs
\begin{align*}
&\chi_1=\left(x_0^2+P(x_1)\right)\frac{\partial}{\partial x_0},
&&\frac{\partial P}{\partial x_1}\not=0, \\
&\chi_2=\frac{x_0^2+P(x_1)}{x_1}\frac{\partial} {\partial x_0},
&& P\not\equiv 0,\\
&\chi_3=\left(\alpha x_0x_1+P(x_1)\right)\frac{\partial}{\partial
x_0},&&\alpha\in\mathbb{C}^* ,\\
&\chi_4=\frac{(\alpha x_1+\mu)x_0+P(x_1)}{x_1} \frac{\partial}{\partial
x_0},&&\alpha\in\mathbb{C},\hspace{1mm}\mu\in\mathbb{C}^*
\end{align*}
ne sont pas rationnellement int\'egrables.

Le champ $\frac{x_0^2}{x_1}\frac{\partial}{\partial x_0}$ est
rationnellement int\'egrable et a pour flot $\left(\frac{x_0x_1}{x_1
-tx_0},x_1\right).$}
\end{lem}

\begin{proof}[\sl D\'emonstration]
La preuve se fait dans tous les cas par une int\'egration
\'el\'ementaire directe. \`A titre d'exemple l'int\'egration de
$\chi_1$ se fait comme suit. On \'ecrit $\phi_t=(x_0(t),x_1(t
));$ on a $x_1(t)=x_1$ et l'\'equation du flot conduit \`a $\frac{\mathrm{d}x_0(
t)}{\mathrm{d}t}=x_0^2(t)+P(x_1),$ soit
\begin{align*}
\frac{\mathrm{d}x_0(t)}{x_0(t)^2+P(x_1)}=\mathrm{d}t.
\end{align*}
On \'ecrit formellement l'\'equation diff\'erentielle
ci-dessus sous la forme
\begin{align*}
\frac{-1}{2\mathrm{i}\sqrt{P(x_1)}}\left(
\frac{\mathrm{d}x_0(t)}{x_0(t)+\mathrm{i}\sqrt{P(x_1)}}-\frac{\mathrm{d}x_0(t)}{x_0(t)
-\mathrm{i}\sqrt{P(x_1)}}\right)=\mathrm{d}t.
\end{align*}
On obtient par int\'egration directe
\begin{align*}
\log\left(\frac{x_0(t)+\mathrm{i}\sqrt{P(x_1)}}{x_0(t)
-\mathrm{i}\sqrt{P(x_1)}}\right)=\text{cte}-2\mathrm{i}\sqrt{P(x_1)}t
\end{align*}
soit encore en utilisant la condition initiale $(x_0(0),x_1(0))=(x_0,x_1)$
\begin{align*}
\frac{x_0(t)+\mathrm{i}\sqrt{P(x_1)}}{x_0(t)-\mathrm{i}\sqrt{P(x_1)}}=
\frac{x_0+\mathrm{i}\sqrt{P(x_1)}}{x_0-\mathrm{i}\sqrt{P(x_1)}}\hspace{1mm} e^{-2
\mathrm{i}\sqrt{P(x_1)}}.
\end{align*}
Ce qui donne une expression explicite
pour $x_0(t);$ puisque $P$ est non constant le flot $\phi_t$ est
transcendant multivalu\'e.

Les autres cas se traitent de fa\c{c}on analogue.
\end{proof}

\begin{rems}\label{froub}
\begin{itemize}
\item Lorsque $P$ est constant on \'ecrit le champ $\chi_1$ sous la
forme
\begin{align*}
&\left(x_0^2-\alpha^2\right)\frac{\partial}{\partial
x_0},&&\alpha\in~\mathbb{C}.
\end{align*}
Ce champ est rationnellement int\'egrable;
son flot qui est dans $\mathrm{Bir}_2$ s'\'ecrit
\begin{align*}
& \left(\frac{x_0}{1-tx_0}, x_1\right)\hspace{1mm}\text{ pour } \alpha=0,&&
\left(\frac{\alpha(e^{-2\alpha t}+1)x_0+\alpha^2 (e^{-2\alpha
t}-1)}{(e^{-2\alpha t}-1)x_0+\alpha(e^{-2\alpha t}
+1)},x_1\right)\text{ sinon.}
\end{align*}

\item
Lorsque $\mu$ est nul et $\alpha$ non nul le champ
\begin{align*}
\chi_4=\frac{\alpha x_0x_1+P(x_1)}{x_1}\frac{\partial}{\partial x_0}
\end{align*}
est conjugu\'e \`a $\frac{\alpha x_0x_1+1}{x_1}\frac{
\partial}{\partial x_0}$ dont le flot est
\begin{align*}
\left(-\frac{1}{\alpha
x_1}+e^{\alpha t}\left(x_0+\frac{1}{\alpha x_1}\right),x_1\right)
\end{align*}
Enfin si $\alpha$ est nul, le flot de $\chi_4$ est donn\'e par
$\left(x_0+t\frac{P(x_1)}{x_1},x_1\right).$
\end{itemize}
\end{rems}

\bigskip

On d\'eduit des Remarques \ref{frou}, \ref{froub}
et du Lemme \ref{froufrou} la:

\begin{pro}
{\sl Un germe de flot dans $\mathrm{Bir}_2$ qui pr\'eserve une fibration
en droites fibre \`a fibre est, \`a conjugaison lin\'eaire pr\`es,
de l'une des formes suivantes

\begin{small}
\begin{align*}
&(x_0+t(\alpha+x_1^2),x_1); && \left(\frac{x_0}{1-tx_0},x_1\right);
&& \left(x_0+t\frac{P(x_1)}{x_1},
x_1\right); && \left(\frac{x_0x_1}{x_1-tx_0},x_1\right);
\end{align*}
\end{small}

\begin{small}
\begin{align*}
&\left(x_0e^{\mu t}+\frac{x_1^2}{\mu}\left(e^{\mu t}-1\right),
x_1\right);&&
\left(-\frac{1}{\mu x_1}+e^{\mu t}\left(x_0+\frac{1}{\mu x_1}\right),x_1
\right);
&&\left(\frac{\mu(e^{-2\mu t}+1)x_0+\mu^2 (e^{-2\mu
t}-1)}{(e^{-2\mu t}-1)x_0+\mu(e^{-2\mu t} +1)},x_1\right),&
\end{align*}
\end{small}

o\`u $\mu$ d\'esigne un complexe non nul et $P$ un polyn\^ome de
degr\'e inf\'erieur ou \'egal \`a $2.$}
\end{pro}

\subsection{La seconde composante du flot est une
translation}\hspace{1mm}

Pour ce type de flot le g\'en\'erateur infinit\'esimal
est du type
\begin{align*}
\frac{\alpha x_0^2+\ell(x_1)x_0+P(x_1)}{ax_1+b}\frac{
\partial}{\partial x_0}+\frac{\partial}{\partial x_1}
\end{align*}
avec $\ell$ affine
et $P$ polyn\^ome de degr\'e inf\'erieur ou \'egal \`a $2.$ On
cherche donc parmi les champs de vecteurs de la forme ci-dessus
ceux qui sont rationnellement int\'egrables. On peut supposer que
le couple $(a,b)$ prend les valeurs $(0,1)$ ou $(1,0)$ et quitte
\`a faire une homoth\'etie en $x_0$, ou bien $\alpha$ vaut $1$, ou
bien $\alpha$ est nul.

Si $\alpha=1$, alors en faisant agir la transformation
$\left(x_0-\frac{\ell(x_1)}{2},x_1\right)$ on se ram\`ene au cas o\`u
$\ell\equiv 0.$ On a donc les \'eventualit\'es suivantes \`a
consid\'erer
\begin{align*}
&\chi=(x_0^2+P(x_1)) \frac{\partial} {\partial x_0}+\frac{\partial}
{\partial x_1} &&\text{et} &&\chi=\frac{x_0^2+P(x_1)}{x_1}
\frac{\partial}{\partial x_0}+ \frac{\partial}{\partial x_1}.
\end{align*}

Si $\alpha$ est nul, alors, \`a translation en $x_1$
pr\`es, on a $\ell(x_1)=\beta$ ou $\mu x_1;$ on se ram\`ene aux
trois~cas
\begin{align*}
&(a,b)=(0,1), &&(a,b)=(1,0), &&(a,b)=(1,b).
\end{align*}
On doit \'etudier les champs $\chi$ de l'un des types suivants
\begin{align*}
&(\alpha x_0+P(x_1))\frac{\partial} {\partial
x_0}+\frac{\partial}{\partial x_1},&& (\mu x_0x_1+P(x_1))
\frac{\partial}{\partial x_0}+\frac{\partial}{\partial
x_1},&& \frac{\mu x_0x_1+P(x_1)}{x_1} \frac{\partial}{\partial
x_0}+\frac{\partial}{\partial x_1},
\end{align*}

\begin{align*}
&\frac{\alpha x_0+P(x_1)}{x_1} \frac{\partial}{\partial
x_0}+\frac{\partial}{\partial x_1},&& \frac{\mu x_0x_1+P(x_1)}
{x_1+b}\frac{\partial}{\partial x_0}+\frac{\partial}{\partial x_1},
\end{align*}
o\`u $\alpha,$ $\mu,$ $b$ sont des constantes, $b\not=0$ et $P$
un polyn\^ome de degr\'e inf\'erieur ou \'egal \`a $2.$

Lorsque
\begin{align*}
&\chi=(\alpha x_0+p_0+p_1x_1+p_2x_1^2) \frac{\partial}{\partial
x_0}+\frac{\partial}{\partial x_1},&& \alpha\in\mathbb{C}^*,
\hspace{1mm} p_i\in\mathbb{C},
\end{align*} on a
\begin{itemize}
\item ou bien $p_2$ est nul et $\chi$ est lin\'eaire;

\item ou bien $p_2\not=0$ et $\chi$ est lin\'eairement conjugu\'e
\`a
\begin{align*}
&(\alpha x_0+x_1^2)\frac{\partial}{\partial x_0}+\frac{\partial}
{\partial x_1}, &&\alpha\in\mathbb{C}^*
\end{align*}
dont on constate que le flot est polynomial donc rencontr\'e \`a
la Proposition \ref{flotaut}.
\end{itemize}

Supposons que
\begin{align*}
&\chi=(p_0+p_1x_1+p_2x_1^2)\frac{\partial}{\partial x_0}+\frac{\partial}{\partial x_1},
&& p_i\in\mathbb{C}.
\end{align*}
Si $p_2$ est nul, $\chi$ est lin\'eaire; sinon $\chi$ est lin\'eairement
conjugu\'e \`a
\begin{align*}
(p_0+x_1^2)\frac{\partial}{\partial x_0}+\frac{\partial}{\partial x_1}
\end{align*}
\noindent dont le flot est polynomial.

Notons que
\begin{align*}
&\frac{\mu x_0x_1+p_0+p_1x_1+p_2x_1^2}{x_1}\frac{\partial}{\partial
x_0}+\frac{\partial}{\partial x_1}, &&\mu\in\mathbb{C}^*, \hspace{1mm}
p_i\in\mathbb{C}
\end{align*}
est conjugu\'e via $\left(x_0+\frac{p_2}{\mu}+\frac{p_1+1}{\mu}x_1,x_1
\right)$ \`a
\begin{align*}
\frac{\mu x_0x_1+p_0}{x_1}\frac{\partial}{\partial x_0}+\frac{\partial}{\partial x_1}
\end{align*}

Le champ
\begin{align*}
&\chi=\frac{\alpha x_0+p_0+p_1x_1+p_2x_1^2}{x_1}\frac{\partial}{\partial x_0}
+\frac{\partial}{\partial x_1}, && \alpha\in\mathbb{C}^*,\hspace{1mm}p_i\in\mathbb{C};
\end{align*}
est, \`a conjugaison lin\'eaire pr\`es, de l'un des types suivants

\begin{align*}
&\frac{\alpha x_0}{x_1}\frac{\partial}{\partial x_0}+\frac{\partial}{\partial x_1},&&
\frac{x_0+x_1}{x_1}\frac{\partial}{\partial x_0}+\frac{\partial}{\partial x_1},
&&\frac{x_0+p_1x_1+x_1^2}{x_1}\frac{\partial}{\partial x_0}+\frac{\partial}{\partial x_1},\\
&\frac{\alpha x_0+x_1^2}{x_1}\frac{\partial}{\partial x_0}+\frac{\partial}{\partial x_1},&&
\alpha\in\mathbb{C}^*,\hspace{1mm}p_1\in\mathbb{C}.
\end{align*}

Pour finir remarquons que
\begin{align*}
&\frac{\mu x_0x_1+p_1x_1+p_2x_1^2}{x_1+b}\frac{\partial}{\partial x_0}+\frac{
\partial}{\partial x_1}, && \mu,\hspace{1mm}b\in\mathbb{C}^*,\hspace{1mm}
p_i\in\mathbb{C}
\end{align*}
est lin\'eairement conjugu\'e \`a l'un des deux mod\`eles suivants
\begin{align*}
&\frac{\mu x_0x_1}{x_1+b}\frac{\partial}{\partial x_0}+\frac{\partial}{\partial x_1},
&& \frac{\mu x_0x_1}{x_1+b}\frac{\partial}{\partial x_0}+\frac{\partial}{\partial x_1},
&& b\in\mathbb{C}^*.
\end{align*}
Donc finalment on consid\`erera les \'eventualit\'es suivantes
\begin{align*}
&(\mu x_0x_1+P(x_1))
\frac{\partial}{\partial x_0}+\frac{\partial}{\partial
x_1},&&\frac{\mu x_0x_1+p_0}{x_1}\frac{\partial}{\partial x_0}+\frac{
\partial}{\partial x_1},\\
&\frac{\alpha x_0}{x_1}\frac{\partial}{\partial x_0}+\frac{\partial}{\partial x_1},&&
\frac{x_0+x_1}{x_1}\frac{\partial}{\partial x_0}+\frac{\partial}{\partial x_1},\\
&\frac{x_0+p_1x_1+x_1^2}{x_1}\frac{\partial}{\partial x_0}+\frac{\partial}{\partial x_1},
&&\frac{\alpha x_0+x_1^2}{x_1}\frac{\partial}{\partial x_0}+\frac{\partial}{\partial x_1}, \\
&\frac{\mu x_0x_1}{x_1+b}\frac{\partial}{\partial x_0}+\frac{\partial}{\partial x_1},
&& \frac{\mu x_0x_1}{x_1+b}\frac{\partial}{\partial x_0}+\frac{\partial}{\partial x_1},
\end{align*}
avec $\alpha,$ $b$ dans $\mathbb{C}^*,$ $p_i$ dans $\mathbb{C}$ et $P$ dans $\mathbb{C}[x_1],$
$\deg P\leq 2.$

\bigskip

Cette partie et la suivante contiennent de nombreux Lemmes; les
strat\'egies de preuve utilisent l'une des d\'emarches suivantes
\begin{itemize}
\item lorsque cel\`a est possible on calcule explicitement le flot
par int\'egration;

\item pour \'eliminer certains cas on montre que $\mathcal{F}_\chi$ n'a pas d'int\'egrale premi\`ere
rationnelle en cherchant une singularit\'e de type n\oe ud-col
dans un mod\`ele birationnel de $\mathcal{F}_\chi;$ la pr\'esence d'une
telle singularit\'e est \'evidemment une obstruction \`a
l'existence d'une int\'egrale premi\`ere rationnelle. On montre
ensuite (en g\'en\'eral via \texttt{Maple}) qu'il n'y a pas de
sym\'etrie forte. Le Th\'eor\`eme \ref{symforte} assure alors que
$\chi$ n'est pas rationnellement int\'egrable;

\item une autre proc\'edure consiste \`a faire d\'eg\'en\'erer
$\chi$ (par renormalisation) sur un champ $\chi_0$ que l'on sait
ne pas \^etre rationnellement int\'egrable.
\end{itemize}

\bigskip

Commen\c{c}ons par consid\'erer le cas o\`u $\chi=(x_0^2+P(x_1))
\frac{\partial}{\partial x_0}+\frac{\partial}{\partial x_1}.$ Dans
un premier temps supposons $P$ constant.

\begin{lem}\label{cte}
{\sl Le champ
\begin{align*}
(x_0^2-\alpha^2)\frac{\partial}{\partial x_0}+\frac{\partial}{\partial x_1}
\end{align*}
s'int\`egre en
\begin{align*}
\left(\frac{\alpha(e^{-2\alpha t}+1)x_0+ \alpha^2(e^{-2\alpha t}-1)}{(e^{-2
\alpha t}-1)x_0+\alpha(e^{-2 \alpha t}+1)},x_1+t\right).
\end{align*}}
\end{lem}

\begin{proof}[\sl D\'emonstration]
Le champ $(x_0^2-\alpha^2)\frac{\partial}{\partial x_0}$ est
conjugu\'e via $f=\frac{x_0+\alpha}{x_0-\alpha}$ au champ $-2\alpha
x_0\frac{\partial}{\partial x_0}$ dont le flot est $e^{-2\alpha t}x_0$
d'o\`u le r\'esultat.
\end{proof}

\begin{rem}
Si $\phi_t=\left(\frac{\alpha(e^{-2\alpha t}+1)x_0+ \alpha^2(e^{-2\alpha t}-1)}{(e^{-2
\alpha t}-1)x_0+\alpha(e^{-2 \alpha t}+1)},x_1+t\right),$ on a
\begin{align*}
& \phi_0=\mathrm{id}\in\Sigma^0; && \phi_t\in\Sigma^1,\hspace{1mm}
\forall\hspace{1mm} t\in \frac{2\mathrm{i}\pi\mathbb{Z}}{\alpha}\setminus\{0\};
&& \phi_t\in\Sigma^2,\hspace{1mm}\forall\hspace{1mm} t\not\in\frac{2\mathrm{i}\pi\mathbb{Z}}{\alpha}.
\end{align*}
\end{rem}

Puis traitons l'\'eventualit\'e o\`u $P$ n'est pas constant.

\begin{lem}
{\sl Consid\'erons le champ $\chi$ d\'efini par
\begin{align*}
&\chi=(x_0^2+P(x_1))\frac{\partial}{\partial
x_0}+\frac{\partial}{\partial x_1}, && P\in\mathbb{C}[x_1],\hspace{1mm}\deg P\leq 2,
\hspace{1mm}\frac{\partial P}{\partial x_1}\not=0.
\end{align*}
Le flot de $\chi$ n'est pas dans $\mathrm{Bir}_2.$}
\end{lem}

\begin{proof}[\sl D\'emonstration]
Nous allons utiliser le Th\'eor\`eme \ref{symforte}.

\begin{itemize}
\item Commen\c{c}ons par montrer que $\mathcal{F}_\chi$ n'admet pas
d'int\'egrale premi\`ere rationnelle.

Posons $P(x_1):=p_0+p_1x_1+p_2x_1^2$. Le feuilletage
$\mathcal{F}_\chi$ est d\'ecrit par
\begin{align*}
\omega_\chi=(x_0^2+p_0+p_1x_1+p_2x_1^2) \mathrm{d}x_1-\mathrm{d}x_0.
\end{align*}
Quitte \`a changer $x_0$ en $1/x_0$ et $x_1$ en $1/x_1$, on
obtient
\begin{align*}
\omega_\chi=(x_1^2+x_0^2(p_0 x_1^2+p_1x_1+p_2))\mathrm{d}x_1-x_1^4\mathrm{d}x_0.
\end{align*}
Posons $x_1:=tx_0;$ alors
\begin{align*}
\omega_\chi=t(t^2+p_0t^2x_0^2+p_1tx_0+p_2-t^3x_0^2)\mathrm{d}x_0+
x_0(t^2+p_0t^2x_0^2+p_1tx_0+p_2)\mathrm{d}t.
\end{align*}
Supposons $p_2\not=0$. Le $1$-jet au point singulier
$(0,-\sqrt{-p_2})$ est de la forme $t\mathrm{d}x_0;$ la singularit\'e est
de type n\oe ud-col et $\mathcal{F}_\chi$ n'admet donc pas d'int\'egrale
premi\`ere rationnelle.

Lorsque $p_2=0$ la m\^eme suite de transformations
conduit \`a
\begin{align*}
\omega_\chi=(x_1+p_0x_1x_0^2+p_1x_0^2)\mathrm{d}x_1-x_1^3\mathrm{d}x_0.
\end{align*}
On fait alors deux \'eclatements, plus pr\'ecis\'ement on pose
$x_1=tx_0^2,$ et on obtient la forme
\begin{align*}
t(2t+2p_0tx_0+2p_1-t^2x_0^3)\mathrm{d}x_0+x_0(t+p_0tx_0+p_1)\mathrm{d}t;
\end{align*}
on remarque qu'en $(x_0,t)=(0,-p_1)$ celle-ci pr\'esente une
singularit\'e de type n\oe ud-col car $P(x_1)$ est non constant
donc $p_1$ non nul. Ainsi $\mathcal{F}_\chi$ ne peut poss\'eder
d'int\'egrale premi\`ere rationnelle.\bigskip

\item Pour finir montrons que $\chi$ n'a pas de sym\'etrie forte
$Y.$

D'apr\`es la Remarque \ref{sf}, on peut supposer qu'un
tel $Y$ s'\'ecrit $\frac{Q(x_0,x_1)} {Cx_1+D}
\frac{\partial}{\partial x_0}$ avec $Q$ polyn\^ome de degr\'e
inf\'erieur ou \'egal \`a $2$ et $C$, $D$ dans
$\mathbb{C}^2\setminus\{(0,0)\}.$ On remarque que les p\^oles de $Y$
doivent \^etre invariants par $\chi.$ Mais le feuilletage
$\mathcal{F}_\chi$ ne poss\`ede pas de courbe du type $x_1=$ cte invariante
autre que $x_1=0;$ ainsi $(C,D)$ est \`a valeurs dans $\{(0,1),\hspace{1mm}
(1,0)\}.$ Lorsque $(C,D)=(1,0)$ (resp. $(C,D)=(0,1)$) un calcul
direct montre que $Y$ est nul d\`es que $P$ est non constant.
\end{itemize}
\end{proof}

Ensuite on envisage la possibilit\'e o\`u $\chi=\frac{x_0^2+P(x_1)
}{x_1} \frac{\partial}{\partial x_0}+\frac{\partial}{\partial
x_1}.$ Comme pr\'ec\'edemment nous allons distinguer le cas $P$
constant du cas $P$ non constant.

\begin{lem}\label{cte2}
{\sl Soit $\chi$ le champ donn\'e par
\begin{align*}
\frac{x_0^2-\alpha^2}{x_1}\frac{\partial}{\partial x_0}+\frac{\partial}{\partial x_1}.
\end{align*}
Le flot de $\chi$ est birationnel si et seulement si $2\alpha$ est
dans $\mathbb{Z}\setminus \{0\}$ et appartient \`a
$\mathrm{Bir}_2$ si et seulement si $\alpha^2=1/4.$

Lorsque
$\alpha^2=1/4$ le flot $\phi_t$ s'\'ecrit
\begin{align*}
\left(\frac{(t+2x_1)x_0-t/2}{t+2x_1-2tx_0},x_1+t\right).
\end{align*}}
\end{lem}

\begin{proof}[\sl D\'emonstration]
Si $\alpha$ est nul, on a
\begin{align*}
\phi_t=\left(\frac{x_0}{1+x_0(\ln x_1-\ln(x_1+t))},x_1+t\right)
\end{align*}
qui n'est \'evidemment pas birationnel.\bigskip

Supposons d\'esormais $\alpha$ non nul. On \'ecrit
$\phi_t(x_0,x_1)=(x_0(t),x_1+t);$ l'\'equation du flot conduit \`a
$\frac{\mathrm{d}x_0(t)}{\mathrm{d}t}=\frac{x_0^2(t) -\alpha^2}{x_1+t},$ ou encore
\begin{align*}
\frac{\mathrm{d}x_0(t)}{x_0(t)-\alpha}-\frac{\mathrm{d}x_0(t)}{x_0(t)+\alpha}=2\alpha
\frac{\mathrm{d}t}{x_1+t}.
\end{align*}
On en d\'eduit que
\begin{align*}
\phi_t=\left(\frac{\alpha^2(x_1^{2\alpha}-(x_1+t)^{2
\alpha})+\alpha(x_1^{2\alpha}+(x_1+t)^{2\alpha})x_0}{(x_1^{2\alpha}
+(x_1+t)^{2\alpha})\alpha+(x_1^{2\alpha}-(x_1+t)^{2\alpha})x_0},x_1+t
\right).
\end{align*}
Ce flot est rationnel si et seulement si $2\alpha$ est
un entier relatif non nul; il est dans $\mathrm{Bir}_2$ si et seulement si
$\alpha^2=1/4$ auquel cas
\begin{align*}
\phi_t=\left(\frac{(t+2x_1)x_0-t/2}{t+2x_1-2tx_0},x_1+t\right).
\end{align*}
\end{proof}

\begin{lem}\label{sym}
{\sl Le flot $\phi_t$ du champ
\begin{align*}
&\chi=\frac{x_0^2+p_0+p_1 x_1+p_2
x_1^2}{x_1}\frac{\partial}{\partial x_0}+\frac{\partial}{\partial x_1},
&& p_i\in\mathbb{C},\hspace{1mm}(p_1,p_2)\not=0,
\end{align*}
appartient \`a $\mathrm{Bir}_2$ si et seulement si $(p_0,p_1)=(-1/4, 0).$

Si $(p_0,p_1)=(-1/4, 0)$ alors $\phi_t$ s'\'ecrit
\begin{align*}
\left(\frac{f(t)(2tx_0-t+2x_0x_1)+g(t)(1-2x_0+4
p_2tx_1+4p_2x_1^2)}{2(x_1f(t)+g(t)(1-2x_0))},x_1+t \right)
\end{align*}
avec $f(t)=2\sqrt{p_2}\cos(\sqrt{p_2}t)$ et $g(t)= \sin(\sqrt{p_2}t).$}
\end{lem}

\begin{proof}[\sl D\'emonstration]
Cherchons sous quelles conditions $\chi$ poss\`ede une sym\'etrie
forte. On peut la choisir sous la forme $\frac{Q(x_0,x_1)}
{Cx_1+D}\frac{\partial}{\partial x_0}$ o\`u $Q$ d\'esigne un
polyn\^ome de degr\'e $2$ et $C$, $D$ deux complexes non tous deux
nuls. Les \'eventuels p\^oles de $Y$ doivent \^etre invariants par
$\chi$ et la seule droite $x_1=$ cte invariante par $\mathcal{F}_\chi$ est
$x_1=0;$ par suite
\begin{align*}
(C,D)\in\{(0,1),\hspace{1mm}(1,0) \}.
\end{align*}
On cherche alors $Q$ sous la forme $\alpha^2x_0^2+\ell(x_1)x_0+q(
x_1)$ avec $\deg\ell\leq 1$ et $\deg q\leq 2.$ Un calcul direct
montre que si $(C,D)=(0,1)$ alors $P(x_1)=p_0+p_1x_1+p_2x_1^2$ est
constant, ce que l'on a exclu. Lorsque $(C,D)=(1,0)$ c'est encore
un calcul direct qui montre que l'existence d'une sym\'etrie non
triviale implique $(p_0,p_1)=(-1/4,0),$ cas o\`u l'int\'egration
du flot est explicite via \texttt{Maple}. Supposons que $(p_0,p_1)
\not=(-1/4,0)$. La $1$-forme $\omega_\chi=(x_0^2+P(x_1))\mathrm{d}x_1
-x_1\mathrm{d}x_0$ d\'efinit le feuilletage $\mathcal{F}_\chi;$ quitte \`a changer
$x_0$ en $1/x_0$ et $x_1$ en $1/x_1$, cette forme s'\'ecrit
\begin{align*}
(x_1^2+(p_0x_1^2+p_1x_1+p_2)x_0^2)\mathrm{d}x_1-x_1^3\mathrm{d}x_0.
\end{align*}
Posons $x_1=tx_0;$ alors $\omega_\chi$ est du type
\begin{align*}
x_0(t^2+p_0t^2x_0^2+p_1tx_0+p_2)\mathrm{d}t+t(t^2+p_0t^2x_0^2+p_1tx_0+
p_2-t^2x_0)\mathrm{d}x_0.
\end{align*}
Si $p_2\not=0,$ la
singularit\'e $(x_0,t)=(0,\sqrt{-p_2})$ est de type n\oe ud-col et
$\mathcal{F}_\chi$ ne peut poss\'eder d'int\'egrale premi\`ere rationnelle.
Lorsque $p_2=0$ en \'eclatant deux fois (poser $x_1=tx_0^2$) la forme
\begin{align*}
\omega_\chi=(x_1+p_0x_0^2x_1+p_1x_0^2)\mathrm{d}x_1-x_1^2\mathrm{d}x_0
\end{align*}
on obtient
\begin{align*}
x_0(t+p_0tx_0^2+p_1)\mathrm{d}t+(2t^2+2p_0t^2x_0^2+2p_1t-t^2x_0)\mathrm{d}x_0
\end{align*}
qui poss\`ede une singularit\'e de type n\oe ud-col au point $(x_0,t)=(0
,-p_1)$ puisque $p_1$ est non~nul.
\end{proof}

\begin{rem}
On remarque que
\begin{align*}
\left(\frac{f(t)(2tx_0-t+2x_0x_1)
+g(t)(1-2x_0+4p_2tx_1+4p_2x_1^2)}{2(x_1f(t)+g(t)(1-2x_0))},x_1+t \right)
\end{align*}
avec $f(t)=2\sqrt{p_2}\cos(\sqrt{p_2}t)$ et $g(t)=\sin(\sqrt{p_2}
t)$ a pour ensemble d'ind\'etermination
\begin{align*}
&\{(1:0:0),&&(1:0:2), && (tf(t)-g(t):2g(t)t:-2g(t)))\}
\end{align*}
{\it i.e.} $\phi_t$ est dans $\Sigma^3.$ On constate que $g$ s'annule pour
$t\in\mathbb{Z} \pi/\sqrt{p_2}$ et que pour ces valeurs des
param\`etres $\phi_t$ est dans $\Sigma^2.$
\end{rem}

\begin{lem}\label{faittout}
{\sl Le champ $\chi=(\mu x_0x_1+\alpha)\frac{\partial} {\partial
x_0}+\frac{\partial}{\partial x_1}$ n'est pas rationnellement
int\'egrable lorsque $\mu$ est non nul.}
\end{lem}

\begin{proof}[\sl D\'emonstration]
Raisonnons par l'absurde: supposons que le flot de $\chi$ soit
birationnel. Posons $f(x_0, x_1):=(\varepsilon x_0,x_1)$. Notons
\begin{align*}
\chi_\varepsilon=(\mu x_0x_1+\varepsilon\alpha)\frac{\partial}
{\partial x_0}+\frac{\partial}{\partial x_1}
\end{align*}
le conjugu\'e de $\chi$
par $f$ et $\phi_t^\varepsilon=f\phi_tf^{-1}$ le flot de
$\chi_\varepsilon.$ Par le th\'eor\`eme de \textsc{Cauchy}
appliqu\'e en un point g\'en\'erique la limite de
$\phi_t^\varepsilon$ lorsque $\varepsilon$ tend vers $0$, {\it
i.e.} le flot de $\chi_0,$ est birationnel; or le flot de $\chi_0$
s'\'ecrit $(x_0e^{\mu x_1t+\mu t^2/2},x_1+t)$ et n'est pas birationnel.
\end{proof}

\begin{rem}
Consid\'erons le champ $\chi$ d\'efini par
\begin{align*}
&\chi=(\mu x_0x_1+p_0+p_1x_1+p_2x_1^2)\frac{\partial}{\partial x_0}+
\frac{\partial}{\partial x_1}, && \mu\in\mathbb{C}^*,\hspace{1mm} p_i
\in\mathbb{C},\hspace{1mm} (p_1,p_2)\not=(0,0).
\end{align*}
Le champ
$\chi$ est lin\'eairement conjugu\'e \`a un champ du type $(\mu
x_0x_1+\alpha)\frac{\partial}{\partial x_0}+\frac{\partial}{\partial
x_1}$ dont on sait d'apr\`es le Lemme \ref{faittout} qu'il n'est pas
rationnellement int\'egrable.
\end{rem}

\begin{lem}
{\sl Le champ
\begin{align*}
\chi=\frac{\mu x_0x_1+p_0}{x_1}\frac{\partial}{\partial x_0}+\frac{\partial}
{\partial x_1}
\end{align*}
est rationnellement int\'egrable si et seulement
si $p_0=0,$ {\it i.e.} si et seulement s'il est lin\'eaire.}
\end{lem}

\begin{proof}[\sl D\'emonstration]
Supposons $p_0$ non nul et $\chi$ rationnellement int\'egrable.
En conjuguant $\chi$ par $f=(x_0,\varepsilon x_1)$ et en multipliant
par $\varepsilon$, on obtient le champ rationnellement
int\'egrable (pour~$\varepsilon$ non nul)
\begin{align*}
\chi_\varepsilon=\left(\varepsilon\mu x_0+\frac{p_0}{x_1}\right)
\frac{\partial}{\partial x_0}+\frac{\partial}{\partial x_1};
\end{align*}
on note
$\psi_t^\varepsilon$ le flot de $\chi_\varepsilon.$ La limite de
$\psi_t^\varepsilon$ quand $\varepsilon$ tend vers $0$
est birationnelle. Un calcul montre~que
\begin{align*}
\psi_t^0=(x_0+p_0(\ln(x_1+t)-\ln x_1), x_1+t)
\end{align*}
n'est pas birationnel: contradiction.
\end{proof}

\begin{lem}\label{chi1}
{\sl Le flot de
\begin{align*}
& \chi=\frac{\alpha x_0}{x_1}\frac{\partial}{\partial x_0}+\frac{
\partial}{\partial x_1}
\end{align*}
appartient \`a $\mathrm{Bir}_2$ si et seulement
si $\alpha$ appartient \`a $\{-1,1\}.$

Si $\alpha=1,$ resp. $\alpha=-1$ il est du type
\begin{align*}
&\left(\frac{x_0}{x_1}(x_1+t),x_1+t\right),&&\text{resp. }\left(\frac{x_0x_1
} {x_1+t}, x_1+t\right).
\end{align*}}
\end{lem}

\begin{proof}[\sl D\'emonstration]
Par int\'egration directe on obtient
\begin{align*}
\phi_t(x_0,x_1)=\left(\frac{x_0}{x_1^\alpha}(x_1+t)^\alpha,
x_1+t\right);
\end{align*}
le champ $\chi$ est rationnellement int\'egrable si et seulement si
$\alpha$ est entier et dans $\mathrm{Bir}_2$ si et seulement si $\alpha
=\pm 1.$ Les cas $\alpha=1,$ resp. $\alpha=-1$ conduisent \`a
\begin{align*}
&\left(\frac{x_0}{x_1}(x_1+t),x_1+t\right),&& \text{resp.
}\left(\frac{x_0x_1}{x_1+t}, x_1+t\right).
\end{align*}
\end{proof}

\begin{lem}
{\sl Le champ
\begin{align*}
&\chi=\frac{x_0+x_1}{x_1}\frac{\partial}
{\partial x_0}+\frac{\partial}{\partial x_1}
\end{align*}
n'est pas rationnellement int\'egrable.}
\end{lem}

\begin{proof}[\sl D\'emonstration]
On a
\begin{align*}
\phi_t(x_0,x_1)=\left((x_1+t)\ln\left(\frac{x_1+t}{x_1}\right)+\frac{x_0}
{x_1}(x_1+t),x_1+t\right)
\end{align*}
qui n'est pas birationnel.
\end{proof}

\begin{lem}\label{chi3}
{\sl Le flot de \begin{align*}
& \chi=\frac{x_0+p_1x_1+x_1^2}{x_1}\frac{\partial}{\partial x_0}+
\frac{\partial}{\partial x_1}
\end{align*}
appartient \`a $\mathrm{Bir}_2$ si et
seulement si $p_1=0.$

Si $p_1$ est nul, le flot de $\chi$ est de la forme
\begin{align*}
\left(x_0+\frac{x_0t}{x_1}+t^2+x_1t,x_1+t\right).
\end{align*}}
\end{lem}

\begin{proof}[\sl D\'emonstration]
Le flot du champ $\chi$ est donn\'e par
\begin{align*}
\left(x_0+\frac{tx_0}{x_1}+t^2+p_1(x_1+t)\ln\left(\frac{x_1+t}{x_1}\right)+
tx_1,x_1+t\right)
\end{align*}
qui est donc birationnel si et seulement si $p_1=0;$ auquel cas on
a
\begin{align*}
\phi_t=\left(x_0+\frac{x_0t}{x_1}+t^2+x_1t,x_1+t \right).
\end{align*}
\end{proof}

\begin{lem}
{\sl Le flot de
\begin{align*}
 \chi=\frac{\alpha x_0+x_1^2}{x_1}
\frac{\partial}{\partial x_0}+\frac{\partial}{\partial x_1}.
\end{align*}
est dans $\mathrm{Bir}_2$ si et
seulement si $\alpha$ est dans $\{-1,1\}.$

Si $\alpha=1,$ on
retrouve le cas pr\'ec\'edent; si $\alpha=-1$ le flot de $\chi$
s'\'ecrit
\begin{align*}
\left(\frac{tx_1^2+t^2x_1+x_0x_1+t^3/3}{x_1+t},x_1+t\right).
\end{align*}}
\end{lem}

\begin{proof}[\sl D\'emonstration]
 Le champ $\chi$ est lin\'eairement conjugu\'e \`a
\begin{align*}
X_\varepsilon=\frac{\alpha x_0+\varepsilon x_1^2}{x_1}\frac{
\partial}{\partial x_0}+\frac{\partial}{\partial x_1}
\end{align*}
pour $\varepsilon
\not=0.$ Comme $X_0$ est rationnellement int\'egrable et
dans $\mathrm{Bir}_2$ si et seulement si $\alpha$ appartient \`a
$\{-1,1\}$ (Lemme \ref{chi1}) le champ $\chi$ n'est pas dans $\mathrm{Bir}_2$
pour $\alpha\not\in\{-1,1\}.$ Si $\alpha=1$ on retrouve le cas
du Lemme \ref{chi3} qui est rationnellement int\'egrable; lorsque
$\alpha=-1,$ on a
\begin{align*}
\phi_t=\left(\frac{x_1^2t+x_1t^2+x_0x_1+t^3/3}{x_1+t},x_1+t\right).
\end{align*}
\end{proof}

Lorsque $\alpha=0$ on a le r\'esultat suivant.

\begin{lem}\label{meme3}
{\sl Soit $\chi$ le champ
\begin{align*}
&\frac{p_0+p_1x_1+p_2x_1^2}{x_1} \frac{\partial} {\partial
x_0}+\frac{\partial}{\partial x_1}, && p_i\in\mathbb{C}.
\end{align*}

Le champ $\chi$ est rationnellement int\'egrable si et
seulement si $p_0=0.$

Lorsque $p_0=0,$ le champ $\chi$ est lin\'eaire.}
\end{lem}

\begin{proof}[\sl D\'emonstration]
Le flot de $\chi$ s'\'ecrit
\begin{align*}
\left(x_0+p_0\ln\left(\frac{x_1+t}{x_1}\right)+(p_1+p_2x_1)t+
\frac{p_2t^2}{2},x_1+t \right);
\end{align*}
ainsi $\phi_t$ est birationnel si et seulement si $p_0=0.$
\end{proof}

\begin{lem}\label{chi1b}
{\sl Le flot de
\begin{align*}
&\chi=\frac{\mu x_0x_1}{x_1+b}\frac{\partial} {\partial
x_0}+\frac{\partial}{\partial x_1},&& b\in\mathbb{C}^*
\end{align*}
est dans $\mathrm{Bir}_2$ si et seulement
si $b\mu=\pm 1,$ ce qui conduit aux flots
\begin{align*}
&\left(x_0e^{\mu t}\left(\frac{x_1+1/\mu}{x_1+t+1/\mu}\right),x_1+t
\right),&&\left(x_0e^{\mu t}\left(\frac{x_1+t-1/\mu}{x_1-1/\mu}
\right),x_1+t\right).
\end{align*}}
\end{lem}

\begin{proof}[\sl D\'emonstration]
Le flot du champ $\chi$ est
\begin{align*}
\left(x_0e^{\mu t}\left(\frac{x_1+b}{x_1+t+b}\right)^{b\mu},x_1+t\right).
\end{align*}
On constate qu'il
est dans $\mathrm{Bir}_2$ si $b\mu=\pm 1.$ On en d\'eduit les deux flots
de l'\'enonc\'e.
\end{proof}

\begin{lem}
{\sl Le champ
\begin{align*}
&\chi=\frac{ \mu x_0x_1+1}
{x_1+b}\frac{\partial}{\partial x_0}+\frac{\partial} {\partial x_1},
&& b\in\mathbb{C}^*.
\end{align*}
est rationnellement int\'egrable si et
seulement si $b\mu=1$ et dans ce cas
\begin{align*}
\phi_t=\left(\left(x_0(x_1+b) +\frac{1}{\mu}-\frac{e^{-\mu
t}}{\mu}\right)\frac{e^{\mu t}}{x_1+b+t} ,x_1+t\right).
\end{align*}}
\end{lem}

\begin{proof}[\sl D\'emonstration]
Pour $\varepsilon\not=0$ le champ $\chi$ est
conjugu\'e \`a
\begin{align*}
X_\varepsilon=\frac{\mu x_0x_1+\varepsilon}{x_1+b}\frac{\partial}
{\partial x_0}+\frac{\partial}{\partial x_1};
\end{align*}
d'apr\`es le Lemme \ref{chi1b} le champ $X_0$ est rationnellement
int\'egrable si et seulement si $b\mu=\pm 1$ et
cette condition est n\'ecessaire pour que $\chi$ le soit. Il
faut donc int\'egrer l'\'equation
\begin{align*}
&\frac{\partial x_0(t)} {\partial t}=\frac{\mu x_0(t)(x_1+t)+1}{x_1+b+t},&&
b\mu=\pm 1
\end{align*}
que l'on traite
comme une \'equation lin\'eaire avec second membre. On cherche
$x_0(t)$ sous la forme $C(t)e^{\mu t}(x_1+t+b)^{-b \mu}$ (variation de
la constante); on obtient
\begin{align*}
C'(t)=\frac{e^{-\mu t}}{(x_1+t+b)^{1-b\mu}}.
\end{align*}

Si $b\mu =1$ on a le flot suivant
\begin{align*}
\phi_t=\left(\left(x_0(\mu x_1+1)+1-e^{-\mu t}\right) \frac{
e^{\mu t}}{\mu x_1+1+\mu t},x_1+t\right)
\end{align*}
qui est bien dans $\mathrm{Bir}_2.$

Lorsque $b\mu=-1$ on obtient

\begin{eqnarray}
C'(t)&=& \frac{e^{-\mu t}}{(x_1+t+b)^2}=\frac{e^{-\mu(x_1+t
+b)}}{(x_1+t+b)^2}e^{\mu(x_1+b)}\nonumber\\
&=&
e^{\mu(x_1+b)}\left(\frac{1}{(x_1+b+t)^2}-\frac{\mu}{x_1+b+t}+\sum_{n
\geq 2}\frac{(-\mu)^n(x_1+t+b)^{n-2}}{n!}\right).\nonumber
\end{eqnarray}

Il s'en suit que

\begin{small}
\begin{align*}
x_0(t)=e^{\mu(x_1+t+b)}\left(-1-\mu(x_1+t+b)\ln(x_1+t+b)+\sum_{n\geq 2}
\frac{(-\mu)^n(x_1+t+b)^n}{(n-1)n!}\right)+De^{\mu t}(x_1+t+b).
\end{align*}
\end{small}

La condition initiale $x_0(0)=x_0$ permet
de calculer $D;$ on constate alors que le flot obtenu n'est pas
birationnel (pr\'esence de $\ln$ en particulier).
\end{proof}

\subsection{La seconde composante du flot est une
homoth\'etie}\hspace{1mm}

Dans cette situation le champ $\chi$ a pour seconde
composante $\mu x_1\frac{\partial}{\partial x_1}$ avec $\mu\in\mathbb{C}^*.$
Dans la suite nous supposerons que $\mu$
vaut $1$ pour ne pas alourdir l'expos\'e. Nous obtenons donc une
classification \`a \og renormalisation pr\`es du temps\fg\hspace{1mm} par
une homoth\'etie.

Le champ $\chi$ est de la forme
\begin{align*}
\frac{\delta x_0^2+\ell(x_1)x_0+P(x_1)}{Ax_1+B}\frac{\partial}{\partial
x_0}+x_1\frac{\partial}{\partial x_1}.
\end{align*}

Quitte \`a faire agir une application diagonale $(ax_0,bx_1)$,
on peut supposer que $\delta=1$ ou $\delta=0$ et que
$(A,B)$ prend les valeurs $(0,1),$ $(1,0)$ ou $(1,1).$

\bigskip

Lorsque $\delta$ vaut $1,$ on peut supposer, \`a
transformation affine pr\`es, que $\ell$ est nul.

Si $\delta$ est nul, on pose $P=p_0+p_1x_1+p_2x_1^2$ et
$\ell(x_1)=\beta x_1+\gamma.$ On a les \'eventualit\'es
suivantes \`a \'etudier
\begin{itemize}
\item[\textbf{\textit{1.}}] Si $\beta=\gamma=0$, alors
\begin{align*}
\chi=\frac{P(x_1)}{Ax_1+B}\frac{\partial}{\partial x_0}+x_1\frac{\partial}{\partial x_1};
\end{align*}
comme pr\'ec\'edemment on peut se ramener \`a
\begin{align*}
(A,B)\in\{(0,1),\hspace{1mm}(1,0),\hspace{1mm}(1,1)\}.
\end{align*}
Lorsque $(A,B)=(0,1)$ le champ $\chi$ est de la forme
\begin{align*}
&(p_0+p_1x_1+p_2x_1^2)\frac{\partial} {\partial
x_0}+x_1\frac{\partial}{\partial x_1},&& p_i\in\mathbb{C}.
\end{align*}
Si $p_2=0$ le champ $\chi$ est lin\'eaire sinon il est lin\'eairement
conjugu\'e \`a
\begin{align*}
& (\varepsilon
+x_1^2)\frac{\partial}{\partial x_0}+x_1\frac{\partial}{\partial x_1},
&&\varepsilon\in\{0,1\};
\end{align*}
dans ce cas on constate que le flot de $\chi$ est polynomial et
apparait en particulier dans la Proposition \ref{flotaut}.

Si $(A,B)=(1,0),$ on peut supposer que $\chi$ s'\'ecrit
\begin{align*}
&\chi=\frac{p_0+p_1x_1+p_2x_1^2}{x_1} \frac{\partial} {\partial
x_0}+x_1\frac{\partial}{\partial x_1}, && p_0\not=0;
\end{align*}
ce champ est lin\'eairement conjugu\'e \`a
\begin{align*}
&\frac{1+\varepsilon x_1}{x_1}\frac{\partial}{\partial
x_0}+x_1\frac{\partial}{\partial x_1}, &&\varepsilon\in\{0,1
\}.
\end{align*}

L'autre \'eventualit\'e \`a \'etudier sera la suivante
\begin{align*}
& \chi=\frac{p_0+p_1x_1+p_2x_1^2}{x_1+1}\frac{\partial}{\partial x_0}+x_1\frac{\partial}{\partial x_1},
&& p_i\in\mathbb{C}.
\end{align*}

\item[\textbf{\textit{2.}}] Si $\beta=0$ et $\gamma\not=0,$ nous allons
distinguer deux possibilit\'es
\begin{itemize}
\item Lorsque $B-\gamma\not=0,$ quitte \`a faire agir $\left(
x_0+\frac{p_0}{\gamma}+\frac{p_1}{\gamma-B}x_1,x_1\right)$
le champ $\chi$ est du type
\begin{align*}
\frac{\gamma x_0+\mu x_1^2}{Ax_1+B}\frac{\partial} {\partial x_0}+x_1
\frac{\partial}{\partial x_1}.
\end{align*}
Comme toujours le couple $(A,B)$ est \`a valeurs dans
\begin{align*}
\{(1,0),\hspace{1mm}(0,1),\hspace{1mm}(1,1)\}.
\end{align*}

Quand $(A,B)=(0,1)$ on constate par calcul direct que le flot de
$\chi$ est polynomial.

On d\'etaillera donc les \'eventualit\'es

\begin{small}
\begin{align*}
&\hspace{1cm}\frac{b x_0+\mu x_1^2}{x_1}\frac{\partial} {\partial x_0}+x_1
\frac{\partial}{\partial x_1}, &&\frac{b x_0+\mu x_1^2}{x_1+1}
\frac{\partial} {\partial x_0}+x_1
\frac{\partial}{\partial x_1}, && b,\hspace{1mm}\mu\in\mathbb{C}.
\end{align*}
\end{small}

\item Lorsque $B=\gamma$ on peut comme $\gamma$ est non nul
se ramener \`a $B=\gamma=1,$ {\it i.e.}
\begin{align*}
&\chi=\frac{x_0+p_1x_1+p_2x_1^2}{Ax_1+1}\frac{\partial}{\partial x_0}
+x_1\frac{\partial}{\partial x_1}, && p_i,\hspace{1mm} A\in\mathbb{C}.
\end{align*}
Si $p_1$ est nul, $\chi$ est d'un des types pr\'ec\'edents;
si $p_1$ est non nul alors $p_1=1$ \`a homoth\'etie pr\`es.
Supposons donc que $p_1=1.$ On remarque alors que si
$A$ est nul, le flot est polynomial donc apparait dans la
Proposition \ref{flotaut}. Si $A$ est non nul, on a, quitte \`a
faire une homoth\'etie, $A=1.$ Ainsi le seul cas \`a \'etudier
est le suivant
\begin{align*}
&\chi=\frac{x_0+x_1+\mu x_1^2}{x_1+1}\frac{\partial}{\partial x_0}
+x_1\frac{\partial}{\partial x_1}.
\end{align*}
\end{itemize}

\item[\textbf{\textit{3.}}] Si $\beta\not=0,$ on va distinguer plusieurs
\'eventualit\'es
\begin{itemize}
\item Quand $A\not=\beta$ on a, quitte \`a conjuguer $\chi$ par
$\left(x_0+\frac{p_1+bB-\gamma b}{\beta}+\frac{p_2}{\beta-A}x_1,
x_1\right)$ puis par des homoth\'eties
\begin{align*}
\chi=\frac{\beta x_0x_1+\gamma x_0+\mu}{Ax_1+B}\frac{\partial}{\partial x_0}
+x_1 \frac{\partial}{\partial x_1}, && \beta\in\mathbb{C}^*,\hspace{1mm}
\gamma,\hspace{1mm}\mu\in\mathbb{C}
\end{align*}
le couple $(A,B)$ appartenant \`a $\{(0,1),\hspace{1mm}(1,0),\hspace{1mm}(1,1)\}.$

\item Quand $A=\beta$ et $\gamma\not=B$ alors on peut supposer
\`a conjugaison lin\'eaire pr\`es que
\begin{align*}
&\chi=\frac{x_0x_1+\gamma x_0+p_2x_1^2}{x_1+B}\frac{\partial}
{\partial x_0}+x_1\frac{\partial}{\partial x_1},
\end{align*}
le couple $(\gamma,p_2)$ prenant les valeurs
$\{(0,0),\hspace{1mm}(0,1),\hspace{1mm}(1,0),\hspace{1mm}(1,1)\}.$

\item Quand $A=\beta$ et $\gamma\not=B$ on obtient \`a
conjugaison pr\`es
\begin{align*}
&\chi=\frac{x_0x_1+x_0+p_1x_1+p_2x_1^2}{x_1+1}\frac{\partial}
{\partial x_0}+x_1\frac{\partial}{\partial x_1},&& p_i\in
\mathbb{C}.
\end{align*}
\end{itemize}
\end{itemize}

\bigskip

De mani\`ere analogue au Lemme \ref{cte} on \'etablit
le:

\begin{lem}
{\sl On d\'esigne par $\chi$ le champ $(x_0^2-\alpha^2)\frac{
\partial}{\partial x_0}+x_1\frac{\partial}{\partial x_1}.$ Le flot de
$\chi$ appartient \`a~$\mathrm{Bir}_2$ et s'\'ecrit
\begin{align*}
\left(\alpha\frac{(e^{-2\alpha t}+1)x_0+\alpha(e^{-2\alpha t}-1)}
{(e^{-2\alpha t}-1)x_0+\alpha(e^{-2\alpha t}+1)},e^tx_1\right).
\end{align*}}
\end{lem}

\begin{rem}
Le flot
\begin{small}
\begin{align*}
\phi_t=\left(\alpha\frac{(e^{-2\alpha t}+1)x_0+\alpha(e^{-2\alpha t}-1)}
{(e^{-2\alpha t}-1)x_0+\alpha(e^{-2\alpha t}+1)},e^tx_1\right).
\end{align*}
\end{small}

\noindent compte g\'en\'eriquement deux
points d'ind\'etermination qui sont immobiles; on constate que
\begin{align*}
\mathrm{Exc}\hspace{1mm}\phi_t=\{x_2=0, (e^{-2\alpha t}-1)x_0+\alpha(e^{-2
\alpha t}+1)x_2=0\}
\end{align*}
est constitu\'e d'une droite fixe et d'une droite mobile.

On remarque que
\begin{align*}
&\phi_t\in\Sigma^0, \hspace{1mm}\forall\hspace{1mm} t\in\frac{2\mathrm{i}\pi\mathbb{Z}}{\alpha};
&&\phi_t\in\Sigma^2, \hspace{1mm}\forall\hspace{1mm} t\not\in\frac{2\mathrm{i}\pi\mathbb{Z}}{\alpha}.
\end{align*}
\end{rem}

\begin{lem}
{\sl Le flot du champ
\begin{align*}
& (x_0^2+P(x_1))\frac{\partial }{\partial x_0}+x_1
\frac{\partial}{\partial x_1}, && P\in\mathbb{C}[x_1],\hspace{1mm}\deg P\leq 2,
\hspace{1mm} \frac{\partial P}{\partial x_1}\not=0
\end{align*}
n'appartient pas \`a $\mathrm{Bir}_2$.}
\end{lem}

\begin{proof}[\sl D\'emonstration]
\'Ecrivons $P$ sous la forme $p_0+p_1x_1+p_2x_1^2.$

Montrons que $\mathcal{F}_\chi$ n'a pas d'int\'egrale premi\`ere
rationnelle. Le feuilletage $\mathcal{F}_\chi$ est d\'ecrit par la
$1$-forme $\omega_\chi=(x_0^2+p_0+p_1x_1+p_2x_1^2)\mathrm{d}x_1-x_1\mathrm{d}x_0;$ quitte \`a
changer $x_0$ en $1/x_0$ et $x_1$ en $1/x_1$, cette $1$-forme s'\'ecrit
\begin{align*}
(x_1^2+p_0x_0^2x_1^2+ p_1x_0^2x_1+p_2x_0^2)\mathrm{d}x_1-x_1^3\mathrm{d}x_0.
\end{align*}
Si $p_2$ est nul, on obtient en posant $x_1=tx_0^2$
\begin{align*}
t(2p_1+2t+2p_0tx_0^2-tx_0)\mathrm{d}x_0+x_0(p_1+t+p_0
tx_0^2)\mathrm{d}t;
\end{align*}
on remarque la pr\'esence d'une singularit\'e de type
n\oe ud-col au point $(0,-p_1).$ Le feuilletage~$\mathcal{F}_\chi$ n'a donc pas
d'int\'egrale premi\`ere rationnelle. Si $p_2$ est non nul, posons
$x_1:=tx_0;$ on obtient
\begin{align*}
\omega_\chi=x_0(t^2+p_0t^2x_0^2 +p_1tx_0+p_2)
\mathrm{d}t+t(t^2+p_0t^2x_0^2+p_1tx_0+p_2-t^2x_0)\mathrm{d}x_0.
\end{align*}
Au point singulier $(0,\sqrt{-p_2})$ on a encore une singularit\'e de type n\oe
ud-col; $\mathcal{F}_\chi$ ne peut donc admettre d'int\'egrale premi\`ere
rationnelle.

En utilisant la Remarque \ref{sf}, un calcul montre que,
$P$ n'\'etant pas constant, $\chi$ admet une sym\'etrie forte si
et seulement si $(p_0,p_1)=\left(-\frac{1}{4},0\right).$ Si
$(p_0,p_1)=\left(-\frac{1}{4},0\right)$ on peut supposer, puisque
$P$ n'est pas constant, que $p_2=1;$ on constate alors que le flot
de $\chi$ n'est pas birationnel
\begin{align*}
&\phi_t=\left(\frac{(f-2e^tx_1)\cos(e^tx_1)+(1+2fe^tx_1)\sin(e^tx_1)}
{2f\cos(e^tx_1)+2\sin(e^tx_1)},e^tx_1\right), \end{align*}
avec
\begin{align*}
f=\frac{2x_0\sin x_1+2x_1\cos x_1-\sin x_1}{\cos x_1+
2x_1\sin x_1-2x_0\cos x_1}.
\end{align*}
\end{proof}

\begin{lem}
{\sl Consid\'erons le champ $\chi$ d\'efini par
\begin{align*}
&\frac{x_0^2-\alpha^2}{x_1}\frac{\partial}{\partial x_0}+x_1\frac{\partial}
{\partial x_1}, &&\alpha\in\mathbb{C}.
\end{align*}

Le flot de $\chi$ appartient \`a $\mathrm{Bir}_2$ si et
seulement si $\alpha=0.$

Si $\alpha$ est nul, $\phi_t$ s'\'ecrit
\begin{align*}
\left(\frac{x_0x_1}{x_1-x_0+e^{-t}x_0},e^t x_1 \right).
\end{align*}}
\end{lem}

\begin{proof}[\sl D\'emonstration]
Si $\alpha$ est nul, le flot s'obtient directement par
int\'egration.

\bigskip

Lorsque $\alpha$ est non nul, le flot s'\'ecrit
\begin{align*}
\left(\alpha\frac{\exp\left(2\alpha e^{-t}/x_1-2\alpha
/x_1\right)(x_0+\alpha)+(x_0-\alpha)}{\exp\left(2\alpha e^{-t}/x_1-2
\alpha/x_1 \right)(x_0+\alpha)-(x_0-\alpha) },e^tx_1\right);
\end{align*}
il n'est pas rationnel.
\end{proof}

\begin{rem}
Le flot $\phi_t=\left(\frac{x_0x_1}{x_1-x_0+e^{-t}x_0},e^t x_1 \right)$
v\'erifie
\begin{align*}
& \phi_t=\mathrm{id}\in\Sigma^0,\hspace{1mm}\forall\hspace{1mm} t\in 2\mathrm{i}\pi\mathbb{Z};
&& \text{sinon }\phi_t\in\Sigma^2.
\end{align*}
\end{rem}

\begin{lem}
{\sl Soit $\chi$ le champ d\'efini par
\begin{align*}
&\frac{x_0^2+p_0+p_1x_1+p_2x_1^2} {x_1} \frac{\partial}{\partial
x_0}+x_1 \frac{\partial}{\partial x_1}, && p_i\in\mathbb{C}, && (p_1,p_2)
\not=(0,0).
\end{align*}
Le flot $\phi_t$ de $\chi$ est birationnel si et seulement si
$p_0=p_1=0.$

Si $p_0=p_1=0,$ on a
\begin{align*}
\phi_t=\left(-\frac{f(t)
\left(x_0(1+\alpha)-2p_2x_1\right)+g(t)\left(x_0(\alpha-1)+2p_2x_1\right)
}{f(t)\left(x_1(1-\alpha)-2x_0\right)+g(t)\left(2x_0-x_1(1+\alpha)\right)}
x_1e^t,x_1e^t\right)
\end{align*}
avec
\begin{align*}
&\alpha=\sqrt{1-4p_2},&& f(t)= e^{\left(-(1+\alpha) t/2\right)},&&
g(t)=e^{\left((\alpha -1)t/2\right)}.
\end{align*}}
\end{lem}

\begin{proof}[\sl D\'emonstration]
En utilisant la Remarque \ref{sf} on constate que $\chi$ admet une
sym\'etrie forte si et seulement si $p_0=p_1=0.$ Lorsque
$p_0=p_1=0,$ l'int\'egration du flot est imm\'ediate. Supposons
d\'esormais que $(p_0,p_1)\not=(0,0)$. Le feuilletage $\mathcal{F}_\chi$
est d\'ecrit par la $1$-forme
\begin{align*}
\omega_\chi=(x_0^2+p_0+p_1x_1+p_2x_1^2)\mathrm{d}x_1-x_1^2\mathrm{d}x_0.
\end{align*}
Le $1$-jet au point
singulier $(\sqrt{-p_0},0)$ est du type n\oe ud-col si $p_0$ est
non nul ce qui interdit l'existence d'une int\'egrale premi\`ere
rationnelle. Voici comment on proc\`ede lorsque $p_0=0,$ {\it
i.e.} lorsque la $1$-forme $\omega_\chi$ s'\'ecrit
\begin{align*}
(x_0^2+p_1x_1+p_2x_1^2)\mathrm{d}x_1-x_1^2\mathrm{d}x_0.
\end{align*}
Nous faisons d'abord une ramification $(x_0,x_1)\mapsto(x_0
,x_1^2)$ pour obtenir
\begin{align*}
2(x_0^2+p_1x_1^2+p_2x_1^4)\mathrm{d}x_1-x_1^3\mathrm{d}x_0;
\end{align*}
puis nous \'eclatons l'origine $(x_0,x_1)=(0,0)$ en posant
$x_0=sx_1$
\begin{align*}
2(p_1+s^2+p_2x_1^2-sx_1)\mathrm{d}x_1-x_1^2\mathrm{d}s.
\end{align*}
On remarque alors puisque $p_1$ est non nul qu'au point $(s,x_1)=(\sqrt{-p_1},0)$
cette $1$-forme pr\'esente une singularit\'e de type n\oe ud-col
ce qui une fois de plus exclut l'existence d'une int\'egrale
premi\`ere rationnelle.
\end{proof}

\begin{lem}
{\sl Le flot du champ $\chi$ d\'efini par
\begin{align*}
&\frac{x_0^2-\alpha^2}{x_1+1}\frac{\partial}{\partial x_0}+x_1\frac{\partial}{\partial
x_1},&& \alpha\in\mathbb{C},
\end{align*}
appartient \`a $\mathrm{Bir}_2$ si et seulement si $4\alpha^2=1.$

Lorsque $4\alpha^2=1,$ on a
\begin{align*}
\phi_t=\left(\frac{e^t(4x_0x_1+2x_0-1)+2x_0+1}{2\left(e^t(2x_1+1
-2x_0)+1+2x_0\right)},e^tx_1\right).
\end{align*}}
\end{lem}

\begin{proof}[\sl D\'emonstration]
Supposons $\alpha$ nul. L'\'equation satisfaite par $\phi_t$ a
pour solution
\begin{align*}
\left(\frac{x_0}{1+x_0(\ln(e^tx_1+1)-\ln(x_1+1)-t)},e^tx_1\right);
\end{align*}
le champ $\chi$ n'est donc pas rationnellement int\'egrable.

\bigskip

Si $\alpha$ est non nul,
\begin{align*}
\phi_t=\left(-\alpha\frac{(x_0+\alpha)(e^tx_1+1)^{2\alpha}+
(x_0-\alpha)(x_1+1)^{2\alpha}e^{2\alpha t}}{(x_0-\alpha)(x_1+1)^{2
\alpha}e^{2\alpha t}-(x_0+\alpha)(e^tx_1+1)^{2\alpha}},e^tx_1
\right),
\end{align*}
expression invariante lorsque l'on change $\alpha$ en
$-\alpha.$ On remarque que $\phi_t$ est birationnel si et
seulement si $2\alpha$ est entier; on constate qu'il est
quadratique si et seulement si $2\alpha=1$ auquel cas il
s'\'ecrit
\begin{align*}
\left(\frac{e^t(4x_0x_1+2x_0-1)+2x_0+1}{2\left(e^t(2x_1+1-2x_0
)+1+2x_0\right)},e^tx_1\right).
\end{align*}
\end{proof}

\begin{rem}
Un calcul montre que
\begin{align*}
&\phi_t=\left(\frac{e^t(4x_0x_1+2x_0-1)+2x_0+1}{2\left(e^t(2x_1+1
-2x_0)+1+2x_0\right)},e^tx_1\right)
\end{align*}
appartient \`a
\begin{itemize}
\item $\Sigma^0$ lorsque $t\in 2\mathrm{i}\pi\mathbb{Z},$ plus pr\'ecis\'ement
pour de tels $t$ on a $\phi_t=\mathrm{id};$

\item $\Sigma^3$ sinon.
\end{itemize}
\end{rem}

\begin{lem}\label{bibi}
{\sl Soit $\chi$ le champ d\'efini par
\begin{align*}
&\frac{x_0^2+p_0+p_1x_1+p_2x_1^2}{x_1+1}\frac{\partial}{\partial x_0}+x_1\frac{
\partial}{\partial x_1}, && p_i\in\mathbb{C},\hspace{1mm} (p_1,p_2)\not=(0,0).
\end{align*}
Le flot de $\chi$ appartient \`a $\mathrm{Bir}_2$ si et seulement si $(p_0,p_1)
=(-1/4+p_2,2p_2).$

Si $(p_0,p_1)=(-1/4+p_2,2p_2)$ et $p_2\not=1/4$ le flot
de $\chi$ s'\'ecrit

\begin{small}
\begin{align*}
\left(\frac{\sin(\beta t)\left(e^tx_1(x_1-2x_0)+4\beta^2(x_1+1)(e^tx_1+1)
\right)+2\beta\cos (\beta t)\left((e^t-1)x_1+2(e^tx_1+1)x_0\right)}{2
\left(2\beta\cos( \beta t)(x_1+1)-\sin(\beta t)(2x_0-x_1)\right)},e^t
x_1\right)
\end{align*}
\end{small}

\noindent avec $\beta=\sqrt{4p_2-1}/2.$

Si $(p_0,p_1,p_2)=(0,1/2,1/4),$ on obtient
\begin{align*}
\phi_t =\left(\frac{1}{2}\frac{(2e^tx_1tx_0-e^tx_1^2t-
4e^tx_1x_0-2e^tx_1-4x_0+2x_1)}{(-2x_1-2+2tx_0-tx_1)},e^tx_1\right).
\end{align*}}
\end{lem}

\begin{rem}
A priori $\beta$ est d\'efini au signe pr\`es mais on constate que
dans l'expression du flot $\beta$ n'intervient que sous forme $\beta^2,$
$\sin(\beta t)$ ou $\beta\cos(\beta t)$ qui sont invariantes par $\beta
\to -\beta.$
\end{rem}

\begin{proof}[\sl D\'emonstration]
Si $\chi$ a une sym\'etrie forte $Y,$ on peut supposer qu'elle est
de la forme $\frac{Q(x_0,x_1)} {Cx_1+D} \frac{\partial}{\partial x_0}$
o\`u $Q$ d\'esigne un polyn\^ome de degr\'e inf\'erieur ou \'egal
\`a deux et $C$, $D$ des complexes non tous deux nuls (\emph{voir}
Remarque \ref{sf}). Un calcul montre que l'existence d'un tel $Y$
implique que 
\begin{align*}
& (p_0,p_1)=(-1/4+p_2,2p_2) && \text{ou} && (p_0,p_1)=(-1/4,0)
\end{align*}
(dans ce dernier cas $p_2$ est non nul).
Pour $(p_0,p_1)=(-1/4+p_2,2p_2)$ le flot s'obtient par un
calcul direct. Lorsque $(p_0,p_1)=(-1/4,0)$ on distingue le
cas $p_2=1/4$ du cas $p_2\not=1/4.$ Pour~$p_2=1/4$
on constate que

\begin{small}
\begin{align*}
\phi_t=\left(\frac{1}{2}\frac{(2x_0-x_1-1)
(e^tx_1+1)\ln\left(\frac{e^tx_1+1}{x_1+1}\right)-2(2x_0-x_1-1)x_1e^t-2x_1(e^tx_1
+1)}{(2x_0-x_1-1)\ln\left(\frac{e^tx_1+1}{x_1+1}\right)-2x_1},e^tx_1\right)
\end{align*}
\end{small}

\noindent n'est pas birationnel. Si $p_2\not=1/4,$ on obtient

\begin{small}
\begin{align*}
\phi_t=\left(\frac{1}{2}+\frac{e^tx_1}{2}\frac{(1-\gamma)(e^tx_1
+1)^\gamma(2x_0-1-x_1-\gamma x_1)+(1+\gamma)(x_1+1)^\gamma(1+x_1-2x_0-\gamma x_1)}
{(x_1+1)^\gamma(1+x_1-2x_0-\gamma x_1)+(e^tx_1+1)^\gamma(2x_0-1-x_1-\gamma x_1)},
e^tx_1\right),
\end{align*}
\end{small}

\noindent avec $\gamma=\sqrt{1-4p_2},$ ainsi $\phi_t$ n'est pas
quadratique. \bigskip

Supposons d\'esormais que $(p_0,p_1)\not=(-1/4, 0)$ et
$(p_0,p_1)\not=(-1/4+p_2,2p_2)$ (rappelons que $(p_1,p_2)\not=(0,0)$).
 Quitte \`a faire une translation
en $x_1$ et \`a faire agir $(x_0,\varepsilon x_1)$ sur $\chi$, on obtient le
champ
\begin{align*}
Y_\varepsilon:=\frac{x_0^2+p_0+p_1(\varepsilon x_1 -1)+p_2
(\varepsilon x_1-1)^2}{\varepsilon x_1}\frac{\partial} {\partial x_0}+\frac{\varepsilon
x_1-1}{\varepsilon}\frac{\partial} {\partial x_1}.
\end{align*}
Posons
$\widetilde{Y_\varepsilon}:=\varepsilon Y_\varepsilon.$ On constate que
\begin{align*}
\widetilde{Y_0}=\frac{x_0^2+p_0-p_1+p_2}{x_1}\frac{\partial}{\partial x_0}
-\frac{\partial} {\partial x_1}
\end{align*}
et on remarque que si $\chi$ est
rationnellement int\'egrable alors $\widetilde{Y_0}$ aussi. Posons $\alpha^2=
-(p_0-p_1+p_2);$ d'apr\`es le Lemme \ref{cte2} le flot de $Y_0$
appartient \`a $\mathrm{Bir}_2$ si et seulement si $2\alpha$ vaut $1$ ou~$-1.$ 
Par suite pour que $\phi_t$ appartienne \`a $\mathrm{Bir}_2$ il
faut que $2\alpha=\pm 1,$ {\it i.e.} $p_1-p_0-p_2=1/4.$ Le
feuilletage $\mathcal{F}_\chi$ est d\'ecrit par la $1$-forme
\begin{align*}
(x_0^2+p_0+p_1x_1+p_2x_1^2)\mathrm{d}x_1-x_1(x_1+1)\mathrm{d}x_0
\end{align*}
soit \`a translation en $x_1$ pr\`es
\begin{align*}
(x_0^2-1/4+p_1x_1+p_2x_1^2-2p_2x_1)\mathrm{d}x_1-x_1(x_1 -1)\mathrm{d}x_0.
\end{align*}
Au point singulier $(-1/2,0)$ le $1$-jet de la forme s'\'ecrit
\begin{align*}
(-x_0+(p_1-2p_2)x_1+\text{ termes de degr\'e sup\'erieur})\mathrm{d}x_1+x_1\mathrm{d}x_0;
\end{align*}
d'apr\`es \textsc{Poincar\'e} $\mathcal{F}_\chi$ est localement conjugu\'e
au feuilletage lin\'eaire donn\'e par
\begin{align*}
(-x_0+(p_1-2p_2)x_1)\mathrm{d}x_1+x_1\mathrm{d}x_0.
\end{align*}
Si $p_1-2p_2$ est non nul,
$(-x_0+(p_1-2p_2)x_1)\mathrm{d}x_1+x_1\mathrm{d}x_0$ n'a pas d'int\'egrale
premi\`ere m\'eromorphe. En effet par homoth\'etie on peut se
ramener \`a $p_1-2p_2=1;$ on constate alors que
\begin{align*}
(-x_0+x_1)\mathrm{d}x_1+x_1\mathrm{d}x_0
\end{align*}
a pour int\'egrale premi\`ere
$x_1\exp\left(\frac{x_0}{x_1}\right).$ Reste le cas o\`u
$p_1-2p_2=0,$ l'\'egalit\'e $p_1-p_0-p_2=1/4$ conduit \`a
$(p_0,p_1)=(p_2-1/4,2p_2),$ cas que l'on a justement exclu.
\end{proof}

\begin{rems}
\begin{itemize}
\item Soit $\phi_t$ le flot d\'efini par

\begin{small}
\begin{align*}
\left(\frac{\sin(\beta t)\left(e^tx_1(x_1-2x_0)+4\beta^2(x_1+1)(e^tx_1+1)
\right)+2\beta\cos (\beta t)\left((e^t-1)x_1+2(e^tx_1+1)x_0\right)}{2
\left(2\beta\cos( \beta t)(x_1+1)-\sin(\beta t)(2x_0-x_1)\right)},e^t
x_1\right)
\end{align*}
\end{small}

\noindent o\`u $\beta=\sqrt{4p_2-1}/2.$ On constate que
\begin{align*}
& \phi_0=\mathrm{id}\in\Sigma^0; && \phi_t\in\Sigma^2,\hspace{1mm}
\forall\hspace{1mm} t\in 2\mathrm{i}\pi\mathbb{Z}\setminus\{0\}; && \phi_t\in\Sigma^3,
\hspace{1mm}\forall\hspace{1mm} t\not\in 2\mathrm{i}\pi\mathbb{Z}.
\end{align*}

\bigskip

\item Si $(p_0,p_1,p_2)=\left(0,1/2,1/4\right)$ cas o\`u
\begin{align*}
\phi_t=\left(\frac{2e^tx_1tx_0-e^tx_1^2t-
4e^tx_1x_0-2e^tx_1-4x_0+2x_1}{2(-2x_1-2+2tx_0-tx_1)},e^tx_1\right),
\end{align*}
on a
\begin{align*}
& \phi_0=\mathrm{id}\in\Sigma^0; && \phi_t\in\Sigma^2,\hspace{1mm}\forall\hspace{1mm}
t\in 2\mathrm{i}\pi\mathbb{Z}\setminus\{0\}; && \phi_t\in\Sigma^3,\hspace{1mm}\forall
\hspace{1mm} t\not\in 2\mathrm{i}\pi\mathbb{Z}.
\end{align*}
\end{itemize}
\end{rems}

\begin{lem}
{\sl Le flot de $\frac{1+\varepsilon x_1}{x_1}\frac{\partial}{\partial
x_0}+x_1\frac{\partial}{\partial x_1},$ avec $\varepsilon$ dans $\{0,1
\},$ s'\'ecrit
\begin{align*}
\left(x_0+\frac{1}{x_1}+t\varepsilon-\frac{e^{-t}}{x_1},e^tx_1\right).
\end{align*}}
\end{lem}

Une int\'egration directe conduit au:

\begin{lem}
{\sl Le flot du champ
\begin{align*}
&\chi=\frac{p_0+p_1x_1+p_2x_1^2}{x_1+1}\frac{\partial}{\partial x_0}+x_1
\frac{\partial}{\partial x_1}, && p_i\in\mathbb{C},
\end{align*}
appartient \`a $\mathrm{Bir}_2$ si et seulement si $p_1=p_0+p_2.$

Lorsque $p_1=p_0+p_2$ le champ $\chi$ est affine.}
\end{lem}

\begin{lem}
{\sl Consid\'erons le champ
\begin{align*}
&\chi=\frac{bx_0+\mu x_1^2}{x_1}\frac{\partial}{\partial x_0}+x_1\frac{\partial}
{\partial x_1}, && b,\hspace{1mm}\mu\in\mathbb{C}.
\end{align*}
Le flot de $\chi$ appartient \`a $\mathrm{Bir}_2$ si et seulement
si $b=0.$

Si $b$ est nul, $\phi_t$ est lin\'eaire.}
\end{lem}

\begin{proof}[\sl D\'emonstration]
Si $\chi$ admet une sym\'etrie forte $Y,$ on peut supposer
(Remarque~\ref{sf}) qu'elle est du type $\frac{Q(x_0,x_1)}
{Cx_1+D}\frac{\partial}{\partial x_0}$ avec $Q$ polyn\^ome de degr\'e
$2$ et $C$, $D$ deux complexes non tous deux nuls; un calcul
montre qu'un tel $Y$ existe si et seulement si $\mu b=0.$ Si $b$
est nul, on obtient 
\begin{align*}
\phi_t=(x_0-\mu x_1+\mu e^tx_1,e^tx_1).
\end{align*}
Si $\mu=0$, le flot de $\chi$ est donn\'e par $(x_0e^{b/x_1}e^{-b
e^{-t}/x_1},e^tx_1);$ le champ $\chi$ n'est alors pas rationnellement
int\'egrable. Dans la suite on supposera $\mu b\not=0;$ on cherche
\`a appliquer le Th\'eor\`eme~\ref{symforte}. Le feuilletage
$\mathcal{F}_\chi$ est d\'efini par la $1$-forme
\begin{align*}
(bx_0+\mu x_1^2)\mathrm{d}x_1-x_1^2\mathrm{d}x_0;
\end{align*}
au point $(0,0)$ la singularit\'e en $(0,0)$ est de type n\oe
ud-col, $\mathcal{F}_\chi$ n'admet donc pas d'int\'egrale premi\`ere
rationnelle.
\end{proof}

\begin{lem}
{\sl Si $\mu$ est non nul, le champ
\begin{align*}
\chi=\frac{\mu x_1^2}{x_1+1}\frac{\partial} {\partial
x_0}+x_1\frac{\partial}{\partial x_1},
\end{align*}
n'est pas rationnellement int\'egrable.}
\end{lem}

\begin{proof}[\sl D\'emonstration]
Le flot de $\chi$ est
\begin{align*}
\phi_t=\left(\mu e^tx_1+x_0-\mu x_1+\mu
\ln\left(\frac{x_1+1}{e^tx_1+1}\right),e^tx_1\right)
\end{align*}
qui est transcendant.
\end{proof}

C'est par int\'egration directe que l'on montre le

\begin{lem}\label{dejafai}
{\sl Le champ
\begin{align*}
& \chi=\frac{bx_0}{x_1+1}\frac{
\partial}{\partial x_0}+x_1\frac{\partial}{\partial x_1}, && b\not
=0
\end{align*}
a son flot dans $\mathrm{Bir}_2$ si et seulement si $b=\pm 1.$
\begin{itemize}
\item Si $b=1$, alors
\begin{align*}
\phi_t=\left(\frac{x_0(x_1+1)e^{t}}{(e^tx_1+1)},e^tx_1 \right).
\end{align*}

\item Lorsque $b=-1,$ on a
\begin{align*}
\phi_t=\left(\frac{x_0(e^tx_1+1)}{(x_1+1)e^{t}},e^tx_1\right).
\end{align*}
\end{itemize}}
\end{lem}

\begin{rems}
\begin{itemize}
\item Lorsque $\phi_t=\left(\frac{x_0(x_1+1)e^{t}}{(e^tx_1+1)},e^tx_1 \right)$
on constate que
\begin{align*}
& \phi_t=\mathrm{id}\in\Sigma^0,\hspace{1mm}\forall\hspace{1mm} t\in 2\mathrm{i}\pi\mathbb{Z};
&& \phi_t\in\Sigma^2,\hspace{1mm}\forall\hspace{1mm} t\not\in 2\mathrm{i}\pi\mathbb{Z}.
\end{align*}

\item Les flots
\begin{align*}
\left(\frac{x_0(x_1+1)e^{t}}{(e^tx_1+1)}, e^tx_1
\right)&&\text{et} &&\left(\frac{x_0(e^tx_1
+1)}{(x_1+1)e^{t}},e^tx_1\right)
\end{align*}
sont conjugu\'es via $(1/x_0,x_1).$
\end{itemize}
\end{rems}

\begin{lem}
{\sl Soit $\chi$ le champ d\'efini par
\begin{align*}
&\frac{bx_0+\mu x_1^2}
{x_1+1}\frac{\partial}{\partial x_0}+x_1\frac{\partial}{\partial x_1},
&& b,\hspace{1mm}\mu\in\mathbb{C}^*.
\end{align*}
Le flot de $\chi$ est dans $\mathrm{Bir}_2$ si et seu\-lement si $b$ vaut $1$ auquel
cas on a
\begin{align*}
\phi_t=\left(\frac{x_0+x_0x_1 -\mu x_1^2+\mu x_1^2e^t}
{e^tx_1+1} \hspace{1mm} e^t,e^tx_1\right).
\end{align*}}
\end{lem}

\begin{proof}[\sl D\'emonstration]
Si $\chi$ admet une sym\'etrie forte $Y$, on peut
supposer (Remarque \ref{sf}) qu'elle est de la forme
$\frac{Q(x_0,x_1)}{Cx_1+D}\frac{\partial}{\partial x_0}$ o\`u $Q$
d\'esigne un polyn\^ome de degr\'e inf\'erieur ou \'egal \`a $2$
et $C$, $D$ deux complexes non tous deux nuls. On peut trouver un
tel champ $Y$ si et seulement si $b=\pm 1.$ Si $b=1$, le calcul du
flot est imm\'ediat.

\bigskip

Lorsque $b$ vaut $-1$, le flot $\phi_t$ est du type
\begin{align*}
\left((x_1+e^{-t})\left(\mu e^t-\frac{\mu}{x_1(e^tx_1+1)}-
\mu+\frac{x_0x_1+\mu}{x_1(x_1+1)}+\frac{2\mu}{x_1}\ln\left(\frac{x_1+1}{e^tx_1+1}
\right)\right),e^tx_1\right);
\end{align*}
le champ $\chi$ n'est donc pas rationnellement int\'egrable.

\bigskip

Supposons $b$ distinct de $1$ et de $-1.$ Le feuilletage
$\mathcal{F}_\chi$ est d\'ecrit dans la carte $x_1=1$ par la $1$-forme
$\left((b-1)x_0x_2+\mu-x_0\right)\mathrm{d}x_2+x_2(1+x_2)\mathrm{d}x_0.$ Le $1$-jet au point
singulier $(\mu,0)$ est du type
\begin{align*}
x_2\mathrm{d}x_0+(-x_0+(b-1)\mu x_2)\mathrm{d}x_2.
\end{align*}
Par un argument analogue \`a celui utilis\'e dans le Lemme \ref{bibi} on
obtient que $\mathcal{F}_\chi$ n'a pas d'int\'egrale premi\`ere
rationnelle. Le Th\'eor\`eme \ref{symforte} assure alors que le
flot de $\chi$ n'est pas birationnel quadratique.
\end{proof}

\begin{rem}
Un calcul montre que si $\phi_t$ d\'esigne
\begin{align*}
\left(\frac{x_0+x_0x_1 -\mu x_1^2+\mu x_1^2e^t}
{e^tx_1+1} \hspace{1mm} e^t,e^tx_1\right)
\end{align*}
alors
\begin{align*}
& \phi_t=\mathrm{id}\in\Sigma^0,\hspace{1mm}\forall\hspace{1mm} t\in 2\mathrm{i}\pi\mathbb{Z};
&& \phi_t\in\Sigma^2,\hspace{1mm} \forall\hspace{1mm} t\not\in 2\mathrm{i}\pi\mathbb{Z}.
\end{align*}
\end{rem}

Un calcul direct permet de montrer le

\begin{lem}
{\sl Le flot du champ
\begin{align*}
\chi=\frac{x_0+x_1+\mu x_1^2}{x_1+1}\frac{\partial}{\partial x_0}
+x_1\frac{\partial}{\partial x_1}
\end{align*}
est de la forme
\begin{align*}
\left(\frac{(x_1(t+\mu e^tx_1)+x_0(x_1+1)-\mu x_1^2)e^t}{e^tx_1+1},e^tx_1\right).
\end{align*}}
\end{lem}

\begin{rem}
Soit $\phi_t$ le flot d\'efini par
\begin{align*}
\left(\frac{(x_1(t+\mu e^tx_1)+x_0(x_1+1)-\mu x_1^2)e^t}{e^tx_1+1},e^tx_1\right).
\end{align*}
On a
\begin{align*}
& \phi_0=\mathrm{id}, && \phi_t\in\Sigma^1 \hspace{3mm}\forall
\hspace{1mm} t\in 2\mathrm{i}\pi\mathbb{Z}, && \phi_t\in
\Sigma^2\hspace{3mm}\forall\hspace{1mm}t\not\in 2\mathrm{i}
\pi\mathbb{Z}.
\end{align*}
\end{rem}

\begin{lem}
{\sl Le flot du champ
\begin{align*}
&\chi=(ax_0x_1+bx_0+\mu) \frac{\partial}{\partial
x_0}+x_1\frac{\partial}{\partial x_1}, && a\in \mathbb{C}^*,\hspace{1mm}
b,\hspace{1mm} \mu \in \mathbb{C},
\end{align*}
n'est pas birationnel.}
\end{lem}

\begin{proof}[\sl D\'emonstration]
En utilisant la Remarque \ref{sf} on peut montrer que $\chi$ admet
une sym\'etrie forte si et seulement si $b$ vaut $-1$ ou $\mu=0.$

Si $\mu=0$ ou $b=-1$, l'int\'egration est imm\'ediate;
on obtient respectivement
\begin{small}
\begin{align*}
& \phi_t=(x_0e^{-ax_1} e^{ae^tx_1+bt},e^tx_1)&& \text{et}&&
\phi_t=\left(-\frac{\mu}{ax_1}e^{-t}+\frac{ax_0x_1+\mu}{ax_1}
e^{-ax_1}e^{ae^tx_1-t},e^tx_1\right).
\end{align*}
\end{small}
Dans les deux cas le champ $\chi$ n'est pas rationnellement
int\'egrable.

\bigskip

Supposons $\mu\not=0$ et $b\not=-1.$ Le feuilletage est
d\'efini par la $1$-forme
\begin{align*}
\omega_\chi=(ax_0x_1 +bx_0+\mu)\mathrm{d}x_1-x_1\mathrm{d}x_0.
\end{align*}
Quitte \`a changer $x_0$ en $1/x_0$ et $x_1$ en $1/x_1$, la forme
$\omega_\chi$ s'\'ecrit
\begin{align*}
(ax_0+bx_0x_1+\mu x^2x_1)\mathrm{d}x_1-x_1^2\mathrm{d}x_0;
\end{align*}
au point $(0,0)$ la singularit\'e est de type n\oe ud-col et
$\mathcal{F}_\chi$ n'admet pas d'int\'egrale premi\`ere rationnelle. En
particulier $\chi$ ne peut \^etre rationnellement int\'egrable.
\end{proof}

La preuve du Lemme qui suit se fait par un calcul direct.

\begin{lem}
{\sl On d\'esigne par $\chi$ le champ
\begin{align*}
&\frac{ax_0x_1+\mu}{x_1}\frac{\partial}{\partial x_0}+x_1\frac{\partial}
{\partial x_1}, && a\in\mathbb{C},\hspace{1mm}\mu\in\mathbb{C}^*.
\end{align*}

Si $a\not=-1$, alors \`a conjugaison lin\'eaire pr\`es
\begin{align*}
\phi_t=\left(\left(x_0+\frac{(1-e^{-(a+1)t})}{(a+1)x_1}
\right)e^{at},e^tx_1 \right).
\end{align*}

Si $a$ vaut $-1$, on a \`a conjugaison lin\'eaire pr\`es
\begin{align*}
\phi_t=\left(\left(\frac{t}{x_1}+x_0\right)e^{-t},e^tx_1\right).
\end{align*}

Les deux flots pr\'ec\'edents sont dans $\mathrm{Bir}_2.$}
\end{lem}

\begin{lem}
{\sl Le flot de
\begin{align*}
&\chi=\frac{ax_0x_1+bx_0+\mu}{x_1}\frac{\partial}
{\partial x_0} +x_1\frac{\partial}{\partial x_1},&& b\in\mathbb{C}^*,\hspace{1mm}
a,\hspace{1mm}\mu\in\mathbb{C},
\end{align*}
n'appartient pas \`a $\mathrm{Bir}_2.$}
\end{lem}

\begin{proof}[\sl D\'emonstration]
On constate que $\chi$ admet une sym\'etrie forte si et seulement
si $a$ appartient \`a $\{0,1,2\}$ ou $\mu$ est nul.

\begin{itemize}
\item Si $a=0,$ l'\'equation du flot a pour solution
\begin{align*}
\phi_t=\left(-\frac{\mu}{b}+\left(x_0+\frac{\mu}{b}\right)e^{b
(1-e^{-t})/x_1},e^tx_1\right);
\end{align*}
il en r\'esulte dans ce cas que $\chi$ n'est pas rationnellement int\'egrable.

\bigskip

\item Lorsque $a$ vaut $1,$ le flot de $\chi$ s'\'ecrit
\begin{align*}
\phi_t=\left(-\frac{\mu}{b}+\frac{\mu x_1e^t}{b^2}+\left(x_0+\frac{\mu}{b}
-\frac{\mu x_1}{b^2}\right)\exp\left(t+\frac{b}{x_1}-\frac{be^{-t}}{x_1}\right),e^t x_1
\right);
\end{align*}
on remarque que $\phi_t$ n'est pas birationnel.

\bigskip

\item Supposons $a=2.$ Ici $\phi_t$ est de la forme

\begin{small}
\begin{align*}
\left( -\frac{\mu}{b}+\frac{2x_1\mu e^t}{b^2}
-\frac{2x_1^2\mu e^{2t}}{b^3}+\left(x_0+\frac{\mu}{b}-\frac{2x_1\mu}
{b^2}+\frac{2x_1^2\mu}{b^3}\right)\exp\left(2t+\frac{b}{x_1}-
\frac{be^{-t}}{x_1}\right),e^tx_1\right).
\end{align*}
\end{small}

Dans ce cas le champ $\chi$ n'est donc pas rationnellement int\'egrable.

\bigskip

\item Si $\mu$ est nul, $\phi_t=(x_0e^{b/x_1}e^{at-b e^{-t}/x_1},e^tx_1)$
 est transcendant.

\bigskip

\item Dor\'enavant nous supposerons que
$a\not\in \{0,1,2\}$ et $\mu\not=0.$ Le feuilletage $\mathcal{F}_\chi$ est
d\'ecrit par la $1$-forme $(ax_0x_1+bx_0+\mu)\mathrm{d}x_1-x_1^2\mathrm{d}x_0;$
comme $b$ est par hypoth\`ese non nul la singularit\'e au point
$(-\mu/b,0)$ est de type n\oe ud-col. Le
feuilletage $\mathcal{F}_\chi$ n'admet donc pas d'in\-t\'egrale premi\`ere
rationnelle. D'apr\`es le Th\'eor\`eme \ref{symforte} le champ
$\chi$ n'est pas rationnellement int\'egrable.
\end{itemize}
\end{proof}

On a \`a \'etudier le champ
\begin{align*}
&\chi=\frac{bx_0+\mu}{x_1+1}
\frac{\partial}{\partial x_0}+x_1\frac{\partial}{\partial x_1},&&
b,\hspace{1mm}\mu\in\mathbb{C}.
\end{align*}
Si $b$ est non nul, alors, \`a translation
pr\`es, on peut supposer $\mu$ nul, cas d\'ej\`a trait\'e (Lemme
\ref{dejafai}). Le champ
\begin{align*}
&\chi=\frac{\mu}{x_1+1}\frac{\partial}{\partial x_0} +x_1\frac{\partial}
{\partial x_1}, &&\mu\in\mathbb{C},
\end{align*}
est rationnellement int\'egrable si et seulement si $\mu$ est nul;
dans ce cas, $\phi_t$ est lin\'eaire.

\begin{lem}
{\sl Soit $\chi$ le champ d\'efini par
\begin{align*}
&\frac{ax_0x_1+bx_0+\mu}{x_1+1}\frac{\partial}{\partial x_0}+x_1
\frac{\partial}{\partial x_1}, && a\in\mathbb{C}^*,\hspace{1mm} b,\hspace{1mm}
\mu\in\mathbb{C}.
\end{align*}

Si $\mu$ est nul, $\phi_t$ appartient \`a
$\mathrm{Bir}_2$ si et seulement si $a-b$ prend les valeurs $-1,$
$0$ ou $1.$

Lorsque $a-b$ appartient \`a $\{-1,\hspace{1mm}0,
\hspace{1mm}1\}$ on obtient respectivement
\begin{small}
\begin{align*}
& \phi_t=\left(\frac{e^{bt}x_0(x_1+1)}{e^tx_1+1},e^tx_1\right),&&
\phi_t=\left(x_0e^{bt},e^tx_1\right),
&& \phi_t=\left(\frac{e^{bt}x_0(e^tx_1+1)}{x_1+1},e^tx_1\right).
\end{align*}
\end{small}

Si $\mu$ est non nul, le flot de $\chi$ est birationnel
si et seulement si $a-b=-1$ auquel cas
\begin{align*}
& \phi_t=\left(\frac{bx_0(x_1+1)+\mu-\mu e^{-bt}}{b(e^tx_1+1)}\hspace{1mm}
e^{bt},e^tx_1\right) && \text{si } b\not=0;\\
& \phi_t=\left(\frac{x_0(x_1+1)+\mu t}{e^tx_1+1},e^tx_1\right) &&
\text{sinon.}
\end{align*}}
\end{lem}

\begin{proof}[\sl D\'emonstration]
En conjuguant $\chi$ par $(x_0/\varepsilon,x_1)$, on obtient
\begin{align*}
\chi_\varepsilon=\frac{ax_0x_1+bx_0+\varepsilon\mu}{x_1+1}
\frac{\partial}{\partial x_0}+x_1\frac{\partial}{\partial x_1};
\end{align*}
notons $\psi_t^\varepsilon$ le flot de $\chi_\varepsilon$. La limite de
$\psi_t^\varepsilon$ lorsque $\varepsilon$ tend vers $0$, {\it
i.e.} le flot de
\begin{align*}
\chi_0=\frac{ax_0x_1+bx_0}{x_1+1}
\frac{\partial}{\partial x_0}+x_1\frac{\partial}{\partial x_1},
\end{align*}
est birationnel si $\phi_t$ l'est. Or le flot de $\chi_0$ est
\begin{align*}
\left(\frac{e^{bt}(e^tx_1+1)^{a-b}x_0}{(x_1+1)^{a-b}},x_1e^t\right);
\end{align*}
ce flot est dans $\mathrm{Bir}_2$ si et seulement si $a-b$ prend
les valeurs $-1,$ $0$ ou $1.$ On doit donc \'etudier ces trois
possibilit\'es. Notons que, $\chi_0$ correspondant au cas $\mu=0,$
on peut dans la suite supposer~$\mu$ non nul.

\bigskip

\begin{itemize}
\item[\textbf{\textit{i.}}] Consid\'erons le cas o\`u $a=b;$
nous allons montrer que $\chi$ n'est pas rationnellement
int\'egrable. En utilisant la Remarque \ref{sf}, on montre que
$\chi$ a des sym\'etries fortes si et seulement si $a=b=1,$
$a=b=2$ ou $a=b= -1.$

\begin{itemize}
\item Si $a=b=1,$ le flot est
\begin{align*}
\left(e^t(x_0+\mu)-e^t\mu x_1t-\mu+\mu
x_1e^t\ln\left(\frac{e^tx_1+1}{x_1+1}\right), e^tx_1\right);
\end{align*}
qui n'est pas birationnel.

\bigskip

\item Lorsque $a=b=2,$ on obtient

\begin{small}
\begin{align*}
\phi_t= \left(e^{2t}\mu
x_1^2\left(t+\ln\left(\frac{x_1+1}{e^tx_1+1} \right)\right) +e^t\mu
x_1+e^{2t}\left(x_0-\mu x_1+\frac{\mu}{2}\right)-\frac{\mu}
{2},e^tx_1\right).
\end{align*}
\end{small}

On constate que $\chi$ n'est pas rationnellement int\'egrable.

\bigskip

\item Supposons que $a=b=-1;$ un calcul montre que
\begin{align*}
\phi_t=\left(e^{-t}\left(x_0+\frac{\mu}{x_1}\ln
\left(\frac{ e^tx_1+1}{x_1+1}\right)\right),e^tx_1\right).
\end{align*}
qui est transcendant.

\bigskip

\item Supposons d\'esormais que $a$ soit
distinct de $1,$ $2$ et $-1.$ Le feuilletage $\mathcal{F}_\chi$ n'admet pas
d'int\'egrale premi\`ere rationnelle. En effet $\mathcal{F}_\chi$ est
d\'ecrit par
\begin{align*}
\omega_\chi=(ax_0(x_1+1)+\mu)\mathrm{d}x_1-x_1(x_1+1) \mathrm{d}x_0;
\end{align*}
quitte \`a
changer $x_0$ en $1/x_0,$ on obtient $\omega_\chi=x_0(a (x_1+1)+\mu x_0)\mathrm{d}x_1
+x_1(x_1+1)\mathrm{d}x_0$ soit \`a translation en $x_1$ pr\`es $(ax_1+\mu x_0)x_0\mathrm{d}x_1
+x_1(x_1-1)\mathrm{d}x_0.$ Au point $(0,0)$ la singularit\'e est de type n\oe
ud-col, $\mathcal{F}_\chi$ n'a donc pas d'int\'egrale premi\`ere
rationnelle. Le Th\'eor\`eme \ref{symforte} permet de conclure.
\end{itemize}
\bigskip

\item[\textbf{\textit{ii.}}] Envisageons l'\'eventualit\'e
suivante: $a-b= 1.$ On va encore appliquer le Th\'eor\`eme
\ref{symforte}. En utilisant la Remarque \ref{sf} on constate que
$\chi$ a des sym\'etries fortes si et seulement si~$b=0$ ou $1.$
Dans chacun de ces cas, le flot s'\'ecrit respectivement
\begin{align*}
\left((e^tx_1+1)\left( \mu t+\frac{\mu}{e^tx_1+1}+\mu\ln
\left(\frac{ x_1+1}{e^tx_1+1}\right)+\frac{x_0-\mu}{x_1+1}\right),e^tx_1
\right),
\end{align*}

\begin{small}
\begin{align*}
\left(e^t(e^tx_1+1)\left(
2\mu x_1\ln\left(\frac{e^tx_1+1}{x_1+1}\right)-2\mu x_1t-\frac{\mu x_1}{e^t
x_1+1}-\frac{\mu}{e^t}+\frac{x_0+2\mu x_1+\mu}{x_1+1}\right),e^tx_1
\right),
\end{align*}
\end{small}

\noindent aucun des deux n'\'etant rationnellement int\'egrable.

Supposons que $b$ soit distinct de $0$ et $1.$  Le
feuilletage $\mathcal{F}_\chi$ est d\'ecrit par la $1$-forme
\begin{align*}
\left((b+1)x_0x_1+bx_0+\mu\right)\mathrm{d}x_1-x_1(x_1+1)\mathrm{d}x_0;
\end{align*}
en reprenant l'argument du Lemme
\ref{bibi} au point $(\mu,-1)$ et en utilisant la non nullit\'e 
de~$\mu(b+1)$ on montre que $\mathcal{F}_\chi$ n'a pas d'int\'egrale
premi\`ere rationnelle.\bigskip

\item[\textbf{\textit{iii.}}] Pour finir lorsque $a-b=-1$ on
obtient le flot par un calcul direct.
\end{itemize}
\end{proof}

\begin{rems}
\begin{itemize}
\item Pour le flot
\begin{align*}
& \phi_t=\left(\frac{bx_0(x_1+1)+\mu-\mu e^{-bt}}{b(e^tx_1+1)}\hspace{1mm}
e^{bt},e^tx_1\right), && b\not=0;\\
\end{align*}
on observe que
\begin{align*}
& \phi_0=\mathrm{id}\in\Sigma^0; && \phi_t\in\Sigma^1,\hspace{1mm}\forall\hspace{1mm} t\in
2\mathrm{i}\pi\mathbb{Z}\setminus\{0\}; && \phi_t\in\Sigma^2,\hspace{1mm}\forall\hspace{1mm} t\not\in
2\mathrm{i}\pi\mathbb{Z}.
\end{align*}

\item Lorsque $\phi_t=\left(\frac{e^{bt}x_0(e^tx_1+1)}{x_1+1},e^tx_1\right),$
on a
\begin{align*}
& \phi_t\in\Sigma^0,\hspace{1mm}\forall\hspace{1mm} t\in 2\mathrm{i}\pi\mathbb{Z}; && \phi_t\in\Sigma^2,
\hspace{1mm}\forall\hspace{1mm} t\not\in 2\mathrm{i}\pi\mathbb{Z}.
\end{align*}
\end{itemize}
\end{rems}

Par calcul direct on obtient le

\begin{lem}
{\sl Le flot du champ $\chi=\frac{x_0x_1}{x_1+b}\frac{\partial}{\partial x_0}+x_1
\frac{\partial}{\partial x_1},$ avec $b$ dans $\mathbb{C},$
est donn\'e par
\begin{align*}
\left(\frac{x_0(e^tx_1+b)}{x_1+b},e^tx_1\right).
\end{align*}

Le champ
\begin{align*}
&\frac{x_0x_1+x_1^2}{x_1+b}\frac{\partial}{\partial x_0}+x_1\frac{\partial}{\partial x_1},
&& b\in\mathbb{C}
\end{align*}
n'est pas rationnellement int\'egrable.}
\end{lem}

\begin{rem}
Le flot
\begin{align*}
\left(\frac{x_0(e^tx_1+b)}{x_1+b},e^tx_1\right)
\end{align*} satisfait
\begin{align*}
&\phi_t=\mathrm{id},
\hspace{3mm}\forall\hspace{1mm}t\in 2\mathrm{i}\pi
\mathbb{Z}; && \phi_t\in\Sigma^2,\hspace{3mm}\forall\hspace{1mm}t\not
\in 2\mathrm{i}\pi\mathbb{Z}.
\end{align*}
\end{rem}

\begin{lem}\label{helico}
{\sl Le flot du champ
\begin{align*}
&\chi=\frac{x_0x_1+x_0}{x_1+b}\frac{\partial}{\partial
x_0}+x_1\frac{\partial}{\partial x_1}, && b\in\mathbb{C}
\end{align*}
est dans $\mathrm{Bir}_2$ si et seulement si $b=\frac{1}{2}$
auquel cas
\begin{align*}
\left(x_0e^{2t}\left(\frac{x_1+1/2}{e^tx_1+1/2}\right),e^tx_1\right).
\end{align*}}
\end{lem}

\begin{rem}
Soit $\phi_t$ le flot donn\'e par
\begin{align*}
\phi_t=\left(x_0e^{2t}\left(\frac{x_1+1/2}{e^tx_1+1/2}\right),e^tx_1\right).
\end{align*}
On remarque que
\begin{align*}
& \phi_t=\mathrm{id}\hspace{3mm}\forall
\hspace{1mm}t\in 2\mathrm{i}\pi\mathbb{Z}
&& \phi_t\in\Sigma^2\hspace{3mm} \forall\hspace{1mm} t\not\in
2\mathrm{i}\pi\mathbb{Z} .
\end{align*}
\end{rem}

\begin{proof}[\sl D\'emonstration]
En effet par calcul direct on obtient
\begin{itemize}
\item lorsque $b=0$ le flot s'\'ecrit
\begin{align*}
\left(e^tx_0\exp\left(\frac{tx_1-e^{-t}}{x_1}\right)\exp\left(\frac{1-e^{-t}}{x_1}\right),e^tx_1\right)
\end{align*}
et n'est pas rationnellement int\'egrable;

\item lorsque $b\not=0$ on a
\begin{align*}
\phi_t=\left(x_0e^{t/b}\left(\frac{e^tx_1+b}{x_1+b}\right)^{1-1/b},e^tx_1\right).
\end{align*}
\end{itemize}
\end{proof}

On en d\'eduit le

\begin{lem}
{\sl Le flot du champ
\begin{align*}
&\chi=\frac{x_0x_1+x_0+x_1^2}{x_1+b}\frac{\partial}{\partial
x_0}+x_1\frac{\partial}{\partial x_1}, && b\in\mathbb{C}
\end{align*}
appartient \`a $\mathrm{Bir}_2$ si et seulement si $b=\frac{1}{2}$
auquel cas
\begin{align*}
\phi_t=\left(\frac{(2tx_1^2+x_0(2x_1+1))e^{2t}}{2e^tx_1+1},e^tx_1\right).
\end{align*}}
\end{lem}

\begin{proof}[\sl D\'emonstration]
Le champ
\begin{align*}
&\frac{x_0x_1+x_0+\mu x_1^2}{x_1+b}\frac{\partial}{\partial
x_0}+x_1\frac{\partial}{\partial x_1}
\end{align*}
est lin\'eairement conjugu\'e \`a
\begin{align*}
&X_\varepsilon=\frac{x_0x_1+x_0+\varepsilon x_1^2}{x_1+b}\frac{\partial}{\partial
x_0}+x_1\frac{\partial}{\partial x_1}.
\end{align*}
Or d'apr\`es le Lemme \ref{helico} le flot de $X_0$ appartient \`a
$\mathrm{Bir}_2$ si et seulement si $b=\frac{1}{2};$ ainsi pour
que le flot de $\chi$ soit birationnel quadratique il faut que $b=
\frac{1}{2}.$ Par int\'egration directe on obtient pour $b=\frac{1}{2}$
\begin{align*}
\phi_t=\left(\frac{(2tx_1^2+x_0(2x_1+1))e^{2t}}{2e^tx_1+1},e^tx_1\right)
\end{align*}
d'o\`u le r\'esultat.
\end{proof}

\begin{rem}
On constate que
\begin{align*}
\phi_t=\left(\frac{(2tx_1^2+x_0(2x_1+1))e^{2t}}{2e^tx_1+1},e^tx_1\right)
\end{align*}
v\'erifie
\begin{align*}
&\phi_0=\mathrm{id}, && \phi_t\in\Sigma^1\hspace{3mm}\forall
\hspace{1mm}t\in2\mathrm{i}\pi\mathbb{Z}\setminus\{0\}, && \phi_t\in\Sigma^2
\hspace{3mm}\forall\hspace{1mm}t\in 2\mathrm{i}\mathbb{Z}.
\end{align*}
\end{rem}

Pour finir on a le

\begin{lem}
{\sl Le champ
\begin{align*}
\chi=\frac{x_0x_1+x_0+p_1x_1+p_2x_1^2}{x_1+1}\frac{\partial}{\partial
x_0}+x_1\frac{\partial}{\partial x_1}
\end{align*}
est rationnellement int\'egrable si et seulement si $p_1=p_2$
auquel cas il est lin\'eaire.}
\end{lem}

\begin{proof}[\sl D\'emonstration]
Un calcul montre que le flot de $\chi$ s'\'ecrit
\begin{align*}
\phi_t=\left(x_0+p_1tx_1+(p_2-p_1)x_1\ln\left(\frac{e^tx_1+1}{x_1+1}
\right),e^tx_1\right).
\end{align*}
\end{proof}

\subsection{R\'ecapitulatif}\hspace{1mm}

Rappelons qu'\`a conjugaison lin\'eaire pr\`es les
champs lin\'eaires et leur flot sont en carte affine du type
suivant:

\vspace{5mm}
\begin{center}
\begin{tabular}{|c|c|}
   \hline
    &\\
    champ de vecteurs & flot associ\'e \\
    &\\    
    \hline
        &\\    
     $ \frac{\partial}{\partial x_1}$& $(x_0,x_1+t)$\\
       &\\    
     $x_0\frac{\partial}{\partial x_0}+\frac{\partial}{\partial x_1}$ & $(x_0e^t,x_1+t)$\\
        &\\    
     $(x_1+1)\frac{\partial}{\partial x_0}+\alpha x_1\frac{\partial}{\partial x_1}$ &\hspace{6mm} $(x_0+x_1(e^t-1)+t,e^{\alpha t}x_1)$\hspace{6mm} \\
       &\\    
     $ \alpha_1x_0\frac{\partial}{\partial x_0}+\alpha_2x_1\frac{\partial}
{\partial x_1}$ & $(x_0e^{\alpha_1t},x_1e^{\alpha_2 t})$\\
     &\\
    \hspace{6mm}  $(\alpha x_0+x_1)\frac{\partial}{\partial x_0}+\alpha
x_1\frac{\partial}{\partial x_1}$\hspace{6mm} &$ \left((x_0+tx_1)e^{\alpha t},x_1e^{\alpha t}\right)$\\
     &\\
     $ x_1\frac{\partial}{\partial x_0}+\frac{\partial}{\partial x_1}$ & $\left(x_0+tx_1+\frac{t^2}{2},x_1+t\right)$\\
     &\\
   \hline
\end{tabular}
\end{center}

\vspace{0.6cm}

On dira que le flot $\phi_t$ appartient \`a $\Sigma^i$
si c'est le cas pour $t$ g\'en\'erique. Lorsque la seconde
composante du flot est du type homoth\'etie nous avons, pour
simplifier la pr\'esentation, fait une renormalisation
homoth\'etique dans le temps de sorte que la seconde composante du
g\'en\'erateur infinit\'esimal $\chi$ soit $x_1\frac{\partial}
{\partial x_1}.$ Les tableaux qui suivent donnent la liste des
flots dans $\mathrm{Bir}_2;$ pour chaque flot sont mentionn\'es le
g\'en\'erateur infinit\'esimal, une sym\'etrie forte et une int\'egrale
premi\`ere. Leur lecture doit donc \^etre comprise \`a
conjugaison lin\'eaire et renormalisation pr\`es; au niveau des
champs $\chi$ cel\`a se traduit par une lecture \`a conjugaison
lin\'eaire pr\`es et multiplication scalaire pr\`es.

Dans ce qui suit certaines valeurs des param\`etres produisent
des flots lin\'eaires; nous ne les avons pas mentionn\'ees.

\begin{landscape}
\begin{small}
\begin{center}
\begin{tabular}{|*{5}{c|}l r|}
   \hline
    & & & &\\
    & champ de vecteurs & flot associ\'e & sym\'etrie forte & IP\\
   & & & &\\
   \hline
   & & & &\\
    & $\frac{P(x_1)}{x_1} \frac{\partial}{\partial x_0}$ &
$\left(x_0+t\frac{P(x_1)}{x_1},x_1\right)$ & $\frac{\partial}{\partial x_0}$
& $x_1$\\
    & & & &\\
    & $\frac{\alpha x_0x_1+1}{x_1} \frac{\partial}{\partial x_0}$ &
$\left(-\frac{1}{\alpha x_1}+e^{\alpha t}\left(x_0+\frac{1}{\alpha
x_1}\right),x_1\right)$ & $\frac{1}{x_1}\frac{\partial}{\partial x_0}
+\alpha x_1\frac{\partial}{\partial x_1}$ & $x_1$\\
    & & & &\\
    & $(\alpha x_0+x_1^2)\frac{\partial} {\partial x_0},\hspace{1mm}\alpha\not=0$
& $\left(x_0e^{\alpha t}+\frac{x_1^2}{\alpha}\left( e^{\alpha t}-1\right),x_1
\right)$ & $-2x_1\frac{\partial}{\partial x_0}+\alpha\frac{\partial}
{\partial x_1}$ & $x_1$\\
    & & & &\\
    & $(\alpha+x_1^2)\frac{\partial}{\partial x_0}$ & $(x_0+t(\alpha+x_1^2),x_1)$ &
    $\frac{\partial}{\partial x_0}$ & $x_1$\\
    \hspace*{5mm}$\Sigma^1$\hspace*{5mm}& & & &\\
    & \hspace{2mm} $(\alpha x_0+x_1^2)\frac{\partial}{\partial x_0}+\frac{\partial}{\partial x_1},$
$\alpha\not=0$ \hspace{2mm} & $\psi_1$ & \hspace{2mm} $\left(1+\alpha
x_1+\frac{\alpha^2x_1^2}{2}+\frac{\alpha^3x_0}{2}
\right)\frac{\partial}{\partial x_0}$\hspace{2mm}  &
\hspace{2mm} $\left(x_0+\frac{2}{\alpha^3}
+\frac{2x_1}{\alpha^2}+\frac{x_1^2}{\alpha}\right)e^{-\alpha x_1}$\hspace{2mm} \\
    & & & &\\
    & $(\alpha+x_1^2)\frac{\partial}{\partial x_0}+\frac{\partial}{\partial x_1}$&
   \hspace{2mm}  $(x_0+t(\alpha+x_1^2)+t^2x_1+\frac{t^3}{3},x_1+t)$\hspace{2mm}  & $\frac{\partial}
{\partial x_0}$ & $\alpha x_1+\frac{x_1^3}{3}-x_0$\\
    & & & &\\
    & $(1+x_1^2)\frac{\partial}{\partial x_0}+x_1\frac{
\partial}{\partial x_1},$ & $\left(x_0-\frac{x_1^2}{2}
+t+\frac{x_1^2e^{2t}}{2},e^tx_1\right)$ & $\frac{\partial}
{\partial x_0}$ & $x_1\exp\left(\frac{x_1^2}{2}-x_0\right)$\\
    & & & &\\
    & $\frac{1}{x_1}\frac{\partial}{\partial x_0}+x_1\frac{
\partial}{\partial x_1},$ & $\left(x_0+\frac{1}{x_1}
-\frac{e^{-t}}{x_1},e^tx_1\right)$ & $\frac{\partial}
{\partial x_0}$ & $-\frac{1}{x_1}-x_0$\\
    & & & &\\
   \hline
\end{tabular}
\end{center}
\vspace{0.6cm}
$\psi_1=\left(x_0e^{\alpha t}+\left(\frac{2}{\alpha^3}+
\frac{2x_1} {\alpha^2}+ \frac{x_1^2}{\alpha}\right)(e^{\alpha
t}-1)-\left(\frac{2x_1t}{\alpha}+\frac{t^2} {\alpha}+\frac{2t}
{\alpha^2}\right),x_1+t\right)$
\end{small}
\end{landscape}

\begin{landscape}
\vspace*{0.6cm}
\begin{small}
\begin{center}
\begin{tabular}{|*{5}{c|}l r|}
   \hline
    & & & &\\
    & champ de vecteurs & flot associ\'e & sym\'etrie forte & IP\\
   & & & &\\
   \hline
   & & & &\\
   & $\frac{1+x_1}{x_1}\frac{\partial}{\partial x_0}+x_1\frac{
\partial}{\partial x_1},$ & $\left(x_0+\frac{1}{x_1}
+t-\frac{e^{-t}}{x_1},e^tx_1\right)$ & $\frac{\partial}
{\partial x_0}$ & \hspace{4mm} $x_1\exp\left(-\frac{1}{x_1}-x_0\right)$\hspace{4mm} \\
    & & & & \\
    & $\frac{\alpha x_0x_1+1}{x_1} \frac{\partial}{\partial x_0}+x_1\frac{
\partial}{\partial x_1},\hspace{1mm} \alpha\not=-1$ & $\left(\left(x_0+\frac{(1-e^{
-(\alpha+1)t})}{(\alpha+1)x_1}\right)e^{\alpha t},e^tx_1 \right)$ & $\frac{1
+(1+\alpha)x_0x_1}{x_1}\frac{\partial}{\partial x_0}$ & $\frac{x_1^{1+\alpha}}
{1+(1+\alpha)x_0x_1}$\\
    & & & &\\
    & $\frac{1-x_0x_1}{x_1}\frac{\partial}{\partial x_0}+x_1\frac{
\partial}{\partial x_1}$ & $\left(\left(\frac{t}{x_1}+x_0\right)e^{-t},
e^tx_1\right)$ & $\frac{1}{x_1}\frac{\partial}{\partial x_0}$ & $x_1e^{-x_0x_1}$\\
    & & & &\\
    & $(1+x_1^2)\frac{\partial}{\partial x_0}+x_1\frac{
    \partial}{\partial x_1},$
    & $(x_0+t+\frac{x_1^2}{2}(e^{2t}-1),x_1e^t)$& $\frac{\partial}
{\partial x_0}$ & $x_1\exp\left(\frac{x_1^2}{2}-x_0\right)$\\
    & & & & \\
    \hspace*{5mm}$\Sigma^1$\hspace*{5mm}& $x_1^2\frac{\partial}{\partial x_0}+x_1\frac{\partial}{\partial x_1},$
    & $(x_0+\frac{x_1^2}{2}(e^{2t}-1),x_1e^t)$& $\frac{\partial}
{\partial x_0}$ & $\frac{x_1^2}{2}-x_0$\\
    & & & &\\
    & \hspace{4mm} $(\alpha x_0+x_1^2)\frac{\partial}{\partial x_0}+x_1\frac{\partial}
    {\partial x_1},$ $\alpha\not\in\{1,2\}$\hspace{4mm}  & $\left(x_0e^{\alpha t}+
    \frac{x_1^2}{2-\alpha}(e^{2t}-e^{\alpha t}),x_1e^t\right)$ & \hspace{4mm} $(x_1^2+
    (\alpha-2)x_0)\frac{\partial}{\partial x_0}$\hspace{4mm}  & $\frac{x_1^\alpha}{
x_1^2+(\alpha-2)x_0}$\\
    & & & &\\
    & $(x_0+x_1+x_1^2)\frac{\partial}{\partial x_0}+x_1\frac{\partial}
    {\partial x_1},$ & \hspace{4mm} $\left(x_0e^t+x_1te^t+x_1^2(e^{2t}-e^t),x_1e^t\right)$\hspace{4mm}  &
    $x_1\frac{\partial}{\partial x_0}$ & $x_1\exp\left(x_1-\frac{x_0}{x_1}\right)$\\
    & & & &\\
    & $(x_0+x_1^2)\frac{\partial}{\partial x_0}+x_1\frac{\partial}
    {\partial x_1},$ & $\left(x_0e^t+x_1^2(e^{2t}-e^t),x_1e^t\right)$ &
    $x_1\frac{\partial}{\partial x_0}$ & $x_1-\frac{x_0}{x_1}$\\
    & & & &\\
    & $(2x_0+x_1^2)\frac{\partial}{\partial x_0}+x_1\frac{\partial}{\partial x_1}$
    & $\left(x_0e^{2t}+x_1^2te^{2t},x_1e^t \right)$ & $x_1^2\frac{\partial}
{\partial x_0}$ & $\frac{1}{x_1}\exp\left(\frac{x_0}{x_1^2}\right)$\\
    & & & &\\
   \hline
\end{tabular}
\end{center}
\end{small}
\end{landscape}

\begin{landscape}
\begin{small}
\begin{center}
\begin{tabular}{|*{5}{c|}l r|}
   \hline
    & & & &\\
    & champ de vecteurs & flot associ\'e & sym\'etrie forte & IP\\
   & & & &\\
   \hline
   & & & &\\
   & $x_0^2\frac{\partial}{\partial x_0}$ & $\left(\frac{x_0}{1-tx_0}, x_1\right)$
   & $\frac{x_0^2}{x_1}\frac{\partial}{\partial x_0}$ & $x_1$\\
   & & & & \\
   & $\frac{x_0^2}{x_1}\frac{\partial}{\partial x_0}$ & $\left(\frac{x_0x_1}{x_1-tx_0}, x_1\right)$
   & $x_0^2\frac{\partial}{\partial x_0}$ & $x_1$\\
   & & & & \\
    & \hspace{7mm} $\frac{-x_0}{x_1+1}\frac{\partial}{\partial x_0}+x_1\frac{\partial}
{\partial x_1}$\hspace{7mm}  & $\left(\frac{x_0(e^tx_1+1)}{(x_1+1)e^{t}},e^tx_1\right)$ &
$x_0\frac{\partial}{\partial x_0}$ & $\frac{x_0x_1}{x_1+1}$\\
    & & & &\\
    & $-\frac{x_0}{x_1}\frac{\partial}{\partial x_0}+\frac{\partial}{\partial
    x_1}$& $\left(\frac{x_0x_1}{x_1+t},x_1+t\right)$ & $x_0\frac{\partial}{\partial x_0}$
    & $x_0x_1$\\
    & & & &\\
    & $\frac{x_0}{x_1}\frac{\partial}{\partial x_0}+\frac{
\partial}{\partial x_1}$ & $\left(\frac{x_0}{x_1}(x_1+t),x_1+t\right)$ & \hspace{7mm} $(x_0+x_1)\frac{\partial}
{\partial x_0}$\hspace{7mm}  & $\frac{x_0}{x_1}$\\
   \hspace*{7mm}$\Sigma^2$\hspace*{7mm} & & & &\\
    & $\frac{x_0+x_1^2}{x_1}\frac{\partial}{\partial x_0}+
\frac{\partial}{\partial x_1}$& \hspace{7mm} $\left(x_0+\frac{x_0t}{x_1}+t^2+x_1t,x_1+t
\right)$\hspace{7mm}  & $x_1\frac{\partial}{\partial x_0}$ &$x_1-\frac{x_0}{x_1}$\\
    & & & &\\
    & $\frac{x_1^2-x_0}{x_1}\frac{\partial}{\partial x_0}+\frac{
    \partial}{\partial x_1}$ & $\left( \frac{x_1^2t+x_1t^2+x_0x_1+t^3/3}{x_1+t},x_1+t\right)$&
    $\frac{1}{x_1}\frac{\partial}{\partial x_0}$& $\frac{x_1^3}{3}-x_0x_1$\\
    & & & &\\
    & $x_0^2\frac{\partial}{\partial x_0}+\frac{\partial}{\partial x_1}$ &
     $\left(\frac{x_0}{1-tx_0},x_1+t\right)$ & $\frac{\partial}{\partial x_1}$
     & $x_1+\frac{1}{x_0}$\\
    & & & &\\
    & $\frac{\alpha^2 x_0x_1}{\alpha x_1+1}\frac{\partial}
{\partial x_0}+\frac{\partial}{\partial x_1}$ & $\left(x_0 e^{\alpha t}
\left(\frac{\alpha x_1+1}{\alpha x_1+\alpha t+1}\right),x_1+t\right)$ &
$x_0\frac{\partial}{\partial x_0}$ & $\frac{\alpha e^{\alpha x_1}}{x_0(\alpha x_1+1)}$\\
    & & & &\\
    & $\frac{\alpha^2x_0x_1}{\alpha x_1-1}\frac{\partial}
{\partial x_0}+\frac{\partial}{\partial x_1}$ & $\left( x_0e^{\alpha
t}\left(\frac{\alpha x_1+\alpha t-1}{\alpha x_1-1}\right),x_1+t\right)$ &
$x_0\frac{\partial}{\partial x_0}$ & \hspace{7mm} $\frac{(\alpha x_1-1)e^{\alpha x_1}}{\alpha x_0}$\hspace{7mm} \\
    & & & & \\
   \hline
\end{tabular}
\end{center}
\end{small}
\end{landscape}

\begin{landscape}
\begin{small}
\begin{center}
\begin{tabular}{|*{5}{c|}l r|}
   \hline
    & & & &\\
    & champ de vecteurs & flot associ\'e & sym\'etrie forte & IP\\
   & & & &\\
   \hline
    & & & &\\
    & $(x_0^2-\alpha^2)\frac{\partial}{\partial x_0}+
\frac{\partial}{\partial x_1},$ $\alpha\not=0$ & $\left(\frac{
\alpha(e^{-2\alpha t}+1)x_0+ \alpha^2(e^{-2\alpha t}-1)}{(e^{-2
\alpha t}-1)x_0+\alpha(e^{-2 \alpha t}+1)},x_1+t\right)$ &
$\frac{\partial}{\partial x_1}$ & $e^{2\alpha
x_1}\left(\frac{x_0+\alpha}{x_0-\alpha}\right)$\\
     & & & &\\
  & $\frac{\alpha^2 x_0x_1+\alpha}{\alpha x_1+1}\frac{\partial}{\partial x_0}+\frac{\partial}
{\partial x_1}$ &
\hspace{1.5mm} $\left(\left(x_0\left(x_1+\frac{1}{\alpha}\right)+\frac{1}{\alpha}
-\frac{e^{-\alpha t}}{\alpha}\right)\frac{\alpha e^{\alpha t}}{\alpha x_1+1+\alpha
t},x_1+t\right)$\hspace{1.5mm}  & $\frac{1+x_0+\alpha x_0x_1}{1+\alpha
x_1}\frac{\partial}{\partial x_0}$ & \hspace{1.5mm} $e^{-\alpha x_1}(1+x_0+\alpha x_0x_1)$\hspace{1.5mm} \\
    & & & & \\
    & $\frac{(\alpha+1)x_0x_1+\alpha x_0}{x_1+1}\frac{\partial}{\partial x_0}+x_1\frac{\partial}
{\partial x_1}$ & $\left(\frac{(e^tx_1+1)e^{\alpha t}}{x_1+1}\hspace{1mm} x_0,e^tx_1\right)$ & $x_0\frac{\partial}{\partial x_0}$ & $\frac{x_1^\alpha(x_1+1)}{x_0}$ \\
    & & & & \\
  & $\frac{x_0}{x_1+1}\frac{\partial}{\partial x_0}+x_1\frac{\partial}
{\partial x_1}$ & $\left(\frac{x_0(x_1+1)e^{t}}{(e^tx_1+1)},e^tx_1 \right)$
& $x_0\frac{\partial}{\partial x_0}$ & $\frac{x_1}{x_0(x_1+1)}$\\
    & & & &\\
& $\frac{x_0^2-1/4}{x_1+1}\frac{\partial}{\partial x_0}+x_1\frac{\partial}
{\partial x_1}$ & $\left(\frac{e^t(4x_0x_1+2x_0-1)+2x_0+1}{2\left(e^t(2x_1+1-
2x_0)+1+2x_0\right)},e^tx_1\right)$ & $(1-4x_0^2)\frac{\partial}{\partial
x_0}$& $\frac{x_1(1+2x_0)}{(x_1+1)(1-2x_0)}$\\
    & & & &\\
\hspace*{6mm}$\Sigma^2$\hspace*{6mm}& $\frac{x_0^2}{x_1}\frac{\partial}{\partial x_0}+x_1\frac{\partial}
{\partial x_1}$ & $\left(\frac{x_0x_1}{x_1-x_0+e^{-t}x_0},e^t x_1 \right)$ &
$x_0\frac{\partial}{\partial x_0}+x_1\frac{\partial}{\partial x_1}$ &
$\frac{x_0x_1}{x_1-x_0}$\\
    & & & &\\
& $(x_0^2-\alpha^2)\frac{\partial}{\partial x_0}+x_1\frac{\partial}
    {\partial x_1},$ $\alpha\not=0$ & $\left(\alpha\frac{ (e^{-2\alpha t}+1)x_0+
    \alpha(e^{-2\alpha t}-1)}{(e^{-2\alpha t}-1)x_0 +\alpha(e^{-2\alpha t}+1)},e^tx_1
    \right)$ & $x_1\frac{\partial}{\partial x_1}$ & $x_1^{2\alpha}\left(
\frac{x_0+\alpha}{x_0-\alpha}\right)$\\
    & & & &\\
    & $x_0^2\frac{\partial}{\partial x_0}+x_1\frac{\partial}{\partial x_1}$
    & $\left(\frac{x_0}{1-tx_0},e^tx_1\right)$ & $x_0^2\frac{\partial}
{\partial x_0}$ & $x_1\exp\left(\frac{1}{x_0}\right)$\\
    & & & &\\
    & $\frac{\alpha-x_0x_1}{x_1+1} \frac{\partial}{\partial
x_0}+x_1\frac{\partial}{\partial x_1},$ &
$\left(\frac{\alpha t+x_0(x_1+1)}{e^tx_1+1},e^tx_1\right)$ &
$\frac{1}{x_1+1}\frac{\partial}{\partial x_0}$ & $\frac{\exp(x_0(x_1+1))}{x_1^\alpha}$\\
    & & & &\\
     & \hspace{1.5mm} $\frac{(\beta-1)x_0x_1+\beta x_0+\alpha}{x_1+1} \frac{\partial}{\partial
x_0}+x_1\frac{\partial}{\partial x_1},\hspace{1mm} \alpha\beta\not=0$ \hspace{1.5mm} & $\left(\frac{\beta x_0(x_1+1)+\alpha
-\alpha e^{-\beta t}}{\beta(e^tx_1+1)}\hspace{1mm} e^{\beta t}, e^tx_1\right)$ &
\hspace{1.5mm} $\frac{\alpha/\beta+x_0+x_0x_1}{x_1+1}\frac{\partial}{\partial x_0}$\hspace{1.5mm} &
$\frac{\alpha+\beta x_0+\beta x_0x_1}{x_1^\beta}$\\
    & & & &\\
    & $\frac{x_0} {x_1+1}\frac{\partial}{\partial x_0}+x_1\frac{
\partial}{\partial x_1}$ & $\left(\frac{x_0+x_0x_1} {e^tx_1+1} \hspace{1mm} e^t,e^tx_1\right)$
& $x_0\frac{\partial}{\partial x_0}$ & $\frac{(x_1+1)x_0}{x_1}$\\
    & & & & \\
   \hline
\end{tabular}
\end{center}
\end{small}
\end{landscape}

\begin{landscape}
\begin{small}
\begin{center}
\begin{tabular}{|*{5}{c|}l r|}
   \hline
    & & & &\\
    & champ de vecteurs & flot associ\'e & sym\'etrie forte & IP\\
   & & & &\\
   \hline
       & & & &\\
    & $\frac{x_0+x_1^2} {x_1+1}\frac{\partial}{\partial x_0}+x_1\frac{
\partial}{\partial x_1}$ & $\left(\frac{x_0+x_0x_1-x_1^2+x_1^2e^t} {e^tx_1+1} \hspace{1mm} e^t,
e^tx_1\right)$ & $(x_0-x_1)\frac{\partial}{\partial x_0}$ &
$\frac{(x_1+1)(x_0-x_1)}{x_1}$\\
    & & & &\\
    & \hspace{4mm} $\frac{x_0^2+\beta x_1^2} {x_1} \frac{\partial}{\partial x_0}+x_1\frac{
\partial}{\partial x_1},$ $\beta\not\in\left\{0,\frac{1}{4}\right\}$\hspace{4mm} 
& $\psi_3$ & \hspace{4mm} $x_0\frac{\partial}{\partial x_0}+x_1\frac{\partial}
{\partial x_1}$\hspace{4mm}  & \hspace{2mm} $x_1\left(\frac{\sqrt{4\beta-1}-\mathrm{i}\left(
\frac{2\beta x_1}{x_0}-1\right)}{\sqrt{4\beta-1}+\mathrm{i}\left(\frac{2\beta x_1}{x_0}-1
\right)}\right)^{\frac{\mathrm{i}}{\sqrt{4\beta-1}}}$\hspace{4mm} \\
  & & & &\\
    & $\frac{x_0^2+\frac{x_1^2}{4}}{x_1} \frac{\partial}{\partial x_0}+x_1\frac{
\partial}{\partial x_1}$ & $\left(\frac{e^tx_1(2tx_0-4x_0-tx_1)}{2(2tx_0-tx_1-2x_1)},e^tx_1\right)$ &
$x_0\frac{\partial}{\partial x_0}+x_1\frac{\partial}{\partial x_1}$ &
$x_1\exp\left(
\frac{4x_0}{2x_0-x_1}\right)$\\
    & & & &\\
  \hspace*{6mm} $\Sigma^2$ \hspace*{6mm}&\hspace{4mm} $ \frac{x_0+x_1+\mu x_1^2}{x_1+1}\frac{\partial}{\partial x_0}+x_1\frac{
\partial}{\partial x_1}$\hspace{4mm}  &$\left(\frac{(x_1(t+\mu e^tx_1)+x_0(x_1+1)-\mu x_1^2)e^t}{e^tx_1+1},e^tx_1\right)$ &
$\frac{x_1}{x_1+1}\frac{\partial}{\partial x_0}$& $x_1\exp\left(\frac{\mu x_1^2-x_0-x_0x_1}{x_1}\right)$\\
    & & & &\\
    & $ \frac{x_0x_1}{x_1+b}\frac{\partial}{\partial x_0}+x_1\frac{
\partial}{\partial x_1}$& $\left(\frac{x_0(e^tx_1+b)}{x_1+b},e^tx_1\right)$ &
 $x_0\frac{\partial}{\partial x_0}$&$\frac{x_1+b}{x_0}$\\
    & & & &\\
    & $ \frac{x_0x_1+x_0}{x_1+1/2}\frac{\partial}{\partial x_0}+x_1\frac{
\partial}{\partial x_1}$ & $\left(x_0e^{2t}\left(\frac{x_1+1/2}{e^tx_1+1/2}\right),
e^tx_1\right)$& $x_0\frac{\partial}{\partial x_0}$ & $\frac{x_0\left(x_1+\frac{1}{2}\right)}{x_1^2}$\\
    & & & &\\
    & $ \frac{x_0x_1+x_0+x_1^2}{x_1+1/2}\frac{\partial}{\partial x_0}+x_1\frac{
\partial}{\partial x_1}$& $\left(\frac{(2tx_1^2+x_0(2x_1+1))e^{2t}}{2e^tx_1+1},e^tx_1\right)$& $\frac{x_1^2}{2x_1+1}\frac{\partial}{\partial x_0}$ & $\frac{1}{x_1^2}\exp\left(\frac{(2x_1+1)x_0}{x_1^2}\right)$\\
    & & & &\\
   \hline
\end{tabular}
\end{center}
\vspace{1.2cm}
$\psi_3=\left(-\frac{e^{\left(-(1+\alpha)t/2\right)}
\left(x_0(1+\alpha)-2\beta x_1 \right)+e^{\left((\alpha-1)t/2
\right)}\left(x_0(\alpha-1)+2\beta x_1\right)}{e^{\left(-(1+
\alpha)t/2\right)}\left(x_1(1-\alpha)-2x_0\right)+e^{\left(
(\alpha-1)t/2\right)}\left(2x_0-x_1(1+\alpha)\right)}x_1e^t,
x_1e^t\right)$ et $\alpha=\sqrt{1-4\beta}$
\end{small}
\end{landscape}

\begin{landscape}
\vspace*{1.2cm}
\begin{small}
\begin{center}
\begin{tabular}{|*{5}{c|}l r|}
   \hline
    & & & &\\
    & champ de vecteurs & flot associ\'e & sym\'etrie forte & IP\\
   & & & &\\
   \hline
   & & & &\\
   & $\frac{x_0^2-\frac{1}{4}}{x_1}\frac{\partial}{\partial x_0}+\frac{\partial}
{\partial x_1}$ & $\left(\frac{(t+2x_1)x_0-\frac{t}{2}}{-2tx_0+
t+2x_1},x_1+t\right)$ &
$\left(x_0^2-\frac{1}{4}\right)\frac{\partial}{\partial x_0}$ &
$\frac{2x_0-1}{x_1(2x_0+1)}$\\
    & & & &\\
  & \hspace{1.5mm} $\frac{x_0^2-1/4+\alpha x_1^2}{x_1}\frac{\partial}{\partial x_0}+\frac{
\partial}{\partial x_1},\hspace{1mm} \alpha\not=0$\hspace{1.5mm}  & $\psi_2$ & $\frac{x_0-1/2}{x_1}\frac{\partial}{\partial
x_0}+\frac{\partial}{\partial x_1}$ &
\hspace{1.5mm} $\frac{2x_0-1+2\mathrm{i}\sqrt{\alpha}x_1}{2x_0-1-2\mathrm{i}\sqrt{\alpha}x_1}\exp\left(2
\mathrm{i}\sqrt{\alpha}x_1\right)$\hspace{1.5mm} \\
    & & & & \\
    \hspace*{5mm}$\Sigma^3$\hspace*{5mm} & $\frac{x_0^2+\frac{x_1}{2}+\frac{x_1^2}{4}}{x_1+1}\frac{\partial}
    {\partial x_0}+x_1\frac{\partial}{\partial x_1}$ & \hspace{1.5mm} $\left(\frac{(2e^t
    tx_0x_1-e^ttx_1^2-4e^tx_0x_1-2e^tx_1-4x_0+2x_1)}{2(-2x_1-2+2tx_0-tx_1)},e^tx_1\right)$\hspace{1.5mm}  &
    $\frac{-x_0^2-\frac{x_1^2}{4}+x_0x_1}{x_1+1}\frac{\partial}{\partial x_0}$ &
    $x_1\exp\left(\frac{2(x_1+1)}{2x_0-x_1}\right)$\\
    & & & &\\
    & $\frac{x_0^2-\frac{1}{4}+\alpha(x_1+1)^2}{x_1+1}\frac{\partial}
    {\partial x_0}+x_1\frac{\partial}{\partial x_1},$ & $\psi_4$ &\hspace{1.5mm}  $\frac{(x_0+\frac{1}{2})x_1}{x_1+1}\frac{\partial}{\partial x_0}+x_1\frac{\partial}{\partial x_1}$\hspace{1.5mm}  &
    $\frac{2\beta(x_1+1)-\mathrm{i}(2x_0-x_1)}{2\beta
    (x_1+1)+\mathrm{i}(2x_0-x_1)}x_1^{2\mathrm{i}\beta},$ \\
    & \begin{scriptsize}$\alpha\not\in\left\{0,\frac{1}{4}\right\}$\end{scriptsize} & & &
    \begin{scriptsize} $\beta=\sqrt{4\alpha-1}/2$\end{scriptsize}\\
    & & & & \\
   \hline
\end{tabular}
\end{center}
\vspace{0.6cm}
avec
\begin{align*}
&\psi_2=\left(\frac{2
\sqrt{\alpha}\cos(\sqrt{\alpha}t)(2tx_0-t+2x_1x_0)+\sin(\sqrt{\alpha}t)
(1-2x_0+4\alpha tx_1+4\alpha x_1^2)}{2(2x_1\sqrt{\alpha}\cos(\sqrt{\alpha}t)+\sin(\sqrt{\alpha}t)(1-2x_0))},
x_1+t\right), && \\
& && \\
& && \\
&\psi_4=\left(\frac{\sin(\beta t)\left(e^tx_1(x_1-2x_0)+4\beta^2
(x_1+1)(e^tx_1+1) \right)+2\beta\cos (\beta t)\left((e^t-1)x_1+2(e^t
x_1+1)x_0\right)}{2\left(2\beta\cos (\beta t)(x_1+1)+ \sin(\beta
t)(x_1-2x_0)\right)}, e^t x_1\right), && \beta=\sqrt{4\alpha-1}/2.
\end{align*}
\end{small}
\end{landscape}

\vspace{0.6cm}

L'examen au cas par cas de tous les mod\`eles de flots
birationnels quadratiques (\emph{c.f.} tableaux) montre que la
dimension de l'espace des sym\'etries fortes est sup\'erieure ou
\'egale \`a $2.$ Nous pouvons donc maintenant pr\'eciser
l'\'enonc\'e \ref{symforte}.

\begin{thm}
{\sl Soit $\phi_t$ un flot birationnel quadratique de
g\'en\'erateur infinit\'esimal $\chi.$ Alors la dimension de
$\mathrm{G}(\chi)$ est sup\'erieure ou \'egale \`a $2.$}
\end{thm}

\begin{rems}
\begin{itemize}
\item Ceci g\'en\'eralise le fait suivant: les groupes
ab\'eliens maximaux de $\mathrm{PGL}_3(\mathbb{C})$ sont de dimension $2.$

\item Il y a des exemples o\`u $\dim \mathrm{G}(\chi)\geq 3,$
par exemple pour $\chi=x_1^2\frac{\partial}{\partial x_0}.$

\item L'\'eventualit\'e \textbf{\textit{i.}} du Th\'eor\`eme
\ref{symforte} ne se pr\'esente donc pas. Toutefois l'alternative
a \'et\'e utilis\'ee dans les preuves.
\end{itemize}
\end{rems}

\begin{rem}
On constate en fait a posteriori que si $\phi_t$ est
un germe de flot alors pour~$(x_0,x_1)$ g\'en\'erique fix\'e,
l'application $t\mapsto\phi_t(x_0,x_1)$ est une fonction
m\'eromorphe globale.
\end{rem}

\begin{rem}
On v\'erifie sur chaque exemple que si $m_t$ est un point
d'ind\'etermination mobile alors $m_t$ param\`etre une feuille de
$\mathcal{F}_\chi.$ En fait soit $\mathcal{D}_0$ une droite telle que
$\phi_{-t} (\mathcal{D}_0)=m_t$ et fixons une valeur $t_0$ de
$t;$ on a $m_t=\phi_{-t}(\mathcal{D}_0)=\phi_{t_0-t}(m_0)$
ce qui confirme l'affirmation ci-dessus.
\end{rem}

\begin{rem}
Les configurations de droites contract\'ees et points \'eclat\'es
qu'on a rencontr\'ees sont les suivantes (chaque configuration
est pr\'ecis\'ee avec un exemple)
\begin{itemize}
\item une droite immobile, un point immobile
\begin{align*}
& ((tx_2^2+x_0x_1)e^{-t}:e^tx_1^2:x_1x_2);
\end{align*}

\item une droite immobile, une droite mobile, deux points immobiles
\begin{align*}
& (x_0x_2:x_1(x_2-tx_0):x_2(x_2-tx_0));
\end{align*}

\item une droite immobile, une droite mobile, un point immobile,
un point mobile
\begin{align*}
& (x_0x_1:(x_1+tx_2)^2:x_2(x_1+tx_2));
\end{align*}

\item une droite immobile, deux droites mobiles, deux points immobiles,
un point mobile
\begin{align*}
& (e^t(4x_0x_1+2x_0x_2-x_2^2)+x_2(2x_0+x_2):2e^tx_1P:2x_2P),
\end{align*}
avec $P=e^t(2x_1+x_2-2x_0)+x_2+2x_0.$
\end{itemize}
\end{rem}

\begin{rem}
Soit $f$ une transformation birationnelle quadratique non
p\'eriodique telle que, pour tout $n$ dans $\mathbb{N}$,
l'\'el\'ement $f^n$ soit encore dans $\mathrm{Bir}_2.$ Alors
$\mathrm{G}=\overline{\langle f^{\bullet n}\rangle}^{\hspace{0.1cm}\mathsf{Z}}$ est un
groupe alg\'ebrique ab\'elien non d\'enombrable. Le groupe
$\mathrm{G}$ a un nombre fini de composantes connexes. Comme dans
un groupe de \textsc{Lie} ab\'elien connexe l'application
exponentielle est surjective, il y a un it\'er\'e $f^\ell$ de $f$
qui se plonge dans un flot; par suite $f^\ell$ pr\'eserve une
fibration en droites.
\end{rem}

\noindent\textbf{\textit{Probl\`eme.}} L'adh\'erence dans $\overline{\mathrm{Bir}_2}$
de l'ensemble des flots birationnels quadratiques forme une vari\'et\'e
alg\'ebrique singuli\`ere invariante sous l'action dynamique de
$\mathrm{PGL}_3(\mathbb{C}).$ On peut s'interroger sur la nature de
cette vari\'et\'e.

\clearemptydoublepage
\chapter{Transformations rationnelles, feuilletages, conjugaison
dans les $\Sigma^i$}\label{feuilfeuil}

\section{G\'en\'eralit\'es}

Un feuilletage de degr\'e $\nu$ sur $\mathbb{P}^2(\mathbb{C})$ est
donn\'e par une $1$-forme diff\'erentielle homog\`ene
\begin{align*}
\omega=F_0\mathrm{d}x_0+F_1\mathrm{d}x_1+F_2\mathrm{d}x_2,
\end{align*}
les $F_i$ d\'esignant des \'el\'ements de $\mathbb{C}[x_0,x_1,x_2]_{\nu+1}$ satisfaisant
\begin{align*}
&\text{pgcd }(F_0, F_1,F_2)=1&&\text{et}&& x_0F_0+x_1F_1+x_2F_2=0
\text{ (identit\'e d'\textsc{Euler})}.
\end{align*}
Le lieu singulier $\mathrm{Sing}\hspace{0.1cm}\mathcal{F}_\omega$ du feuilletage
$\mathcal{F}_\omega$ associ\'e \`a $\omega$ est le projectivis\'e de
\begin{align*}
\mathrm{Sing}\hspace{0.1cm}\omega:=\{F_0=F_1=F_2=0\}.
\end{align*}
L'identit\'e d'\textsc{Euler} montre que cet ensemble est toujours non vide;
plus pr\'ecis\'ement
\begin{align*}
\#\hspace{0.1cm}\mathrm{Sing}\hspace{0.1cm}\mathcal{F}_\omega=\nu^2+\nu+1,
\end{align*}
chaque point singulier \'etant compt\'e avec
multiplicit\'e (\cite{Sa}).

Remarquons que le th\'eor\`eme de division de
\textsc{de Rham}-\textsc{Sa\"{\i}to} (\cite{Sa}) sous sa version la plus
\'el\'ementaire implique l'existence de polyn\^omes homog\`enes
$G_0,$ $G_1$ et $G_2$ de degr\'e $\nu$ tels que l'on ait avec des
notations \'evidentes
\begin{align*}
\left[
\begin{array}{c}
F_0\\
F_1\\
F_2
\end{array}
\right] =\left[
\begin{array}{c}
x_0\\
x_1\\
x_2
\end{array}
\right]\wedge \left[
\begin{array}{c}
G_0\\
G_1\\
G_2
\end{array}
\right].
\end{align*}
\'Evidemment les coefficients de la $2$-forme $d\omega$
produisent en fait directement ces $G_i.$

Une courbe $\mathcal{C}$ d'\'equation r\'eduite $h=0$
est dite invariante par le feuilletage $\mathcal{F}_\omega$ si la
$2$-forme~$\omega\wedge dh$ est divisible par $h.$

On note $\mathcal{F}_\nu\label{not26}$ le projectivis\'e de
l'espace vectoriel des formes $\omega$ satisfaisant l'identit\'e
d'\textsc{Euler}
\begin{align*}
\mathcal{F}_\nu=\mathbb{P}\{\omega=F_0\mathrm{d}x_0+F_1\mathrm{d}x_1+F_2\mathrm{d}x_2\hspace{0.1cm}
|\hspace{0.1cm}x_0F_0+x_1F_1+x_2F_2=0,\hspace{0.1cm} \deg F_i=\nu+1\}.
\end{align*}

\`A chaque \'el\'ement de $\mathcal{F}_\nu$ on peut associer un
feuilletage de degr\'e inf\'erieur ou \'egal \`a $\nu$ d\'efini
par $\frac{\omega}{\text{pgcd}(F_0,F_1,F_2)}.$ L'espace des
feuilletages de degr\'e inf\'erieur ou \'egal \`a $\nu$ est not\'e
$\mathcal{F}_\nu^\bullet\label{not26a};$ les feuilletages de degr\'e
pr\'ecis\'ement $\nu$ correspondent aux formes $\omega$ telles que
$\text{pgcd}(F_0,F_1,F_2)=1.$ Il s'agit donc d'un ouvert de
\textsc{Zariski} not\'e $\mathring{\mathcal{F}}_\nu\label{not26aa}.$
L'application naturelle
\begin{align*}
\mathring{\mathcal{F}}_\nu\to\mathcal{F}_\nu^\bullet
\end{align*}
est injective ce qui permettra des abus de langage.

\`A toute transformation rationnelle $f$ on peut associer un feuilletage
$\mathcal{F}(f)$ de la fa\c{c}on suivante
\begin{align*}
&f=(f_0:f_1:f_2)\mapsto\mathcal{F}(f)=(x_1f_2-x_2f_1)\mathrm{d}x_0+(x_2f_0-x_0f_2)\mathrm{d}x_1+
(x_0f_1-x_1f_0)\mathrm{d}x_2
\end{align*}\label{not26b}

\begin{pro}
{\sl Soient $f$ une transformation de \textsc{Cremona} et $\mathcal{F}(f)$ le
feuilletage associ\'e. Si $A$ est un automorphisme de $\mathbb{P}^2(\mathbb{C}),$
alors $\mathcal{F}(A^{-1}fA)=A^*\mathcal{F}(f).$}
\end{pro}

\begin{proof}[\sl D\'emonstration]
Le feuilletage $\mathcal{F}(f)$ est d\'etermin\'e par les champs
de vecteurs associ\'es \`a~$\mathrm{id}$ et $f;$ son conjugu\'e
$A^*\mathcal{F}(f)$ est donc donn\'e par $A^{-1}\mathrm{id} A=
\mathrm{id}$ et $A^{-1}fA.$ Or ces deux champs sont tangents au
feuilletage $\mathcal{F}(A^{-1}fA)$ d'o\`u le r\'esultat.
\end{proof}

Soit $f$ un \'el\'ement du groupe de \textsc{Cremona}; un
\textbf{\textit{point fixe}}\label{ind42} de $f$ est un point $p$
tel que $p$ n'appartienne pas \`a $\mathrm{Ind}\hspace{0.1cm} f\cup\mathrm{Ind}\hspace{0.1cm} f^{-1}$ et
$f(p)=p.$ On note $\mathrm{Fix}\hspace{0.1cm} f\label{not27}$ l'ensemble des
points fixes de $f.$ Si $f$ est un endomorphisme de $\mathbb{P}^2(\mathbb{C}),$
{\it i.e.} une application rationnelle dont l'ensemble
d'ind\'etermination est vide, $\mathrm{Fix}\hspace{0.1cm} f$ est l'ensemble
des points fixes au sens usuel, {\it i.e.}
\begin{align*}
\mathrm{Fix}\hspace{0.1cm} f=\{m\in\mathbb{P}^2(\mathbb{C})\hspace{0.1cm}|\hspace{0.1cm} f(m)=m\}.
\end{align*}

\begin{eg}
La transformation $f=(x_0^2:x_1^2:x_2^2)$ d\'efinit un endomorphisme
quadratique du plan projectif complexe ayant sept points fixes, qui sont exactement
les points singuliers de son feuilletage associ\'e $\mathcal{F}(f).$

Par suite ce ph\'enom\`ene persiste pour une transformation
rationnelle g\'en\'erique.
\end{eg}

\begin{eg}
L'involution de \textsc{Cremona} compte trois points
d'ind\'etermi\-nation 
\begin{align*}
&(1:0:0), && (0:1:0),&& (0:0:1)
\end{align*}
et quatre points fixes
\begin{align*}
&(1:1:1),&&(1:-1:1),&&(-1:1:1),&&(1:1:-1).
\end{align*}
Le feuilletage associ\'e $\mathcal{F}(\sigma)$ est donn\'e par la
forme
\begin{align*}
& x_0(x_2^2-x_1^2)\mathrm{d}x_0+x_1(x_0^2-x_2^2)\mathrm{d}x_1+x_2(x_1^2-x_0^2)\mathrm{d}x_2\\
&\hspace{2mm}
=\frac{1}{2}(x_0^2,x_1^2,x_2^2)^*((x_2-x_1)\mathrm{d}x_0+(x_0-x_2)\mathrm{d}x_1+ (x_1-x_0)\mathrm{d}x_2).
\end{align*}
On en d\'eduit qu'il poss\`ede l'int\'egrale premi\`ere
$\frac{x_0^2-x_2^2}{x_1^2-x_2^2}.$
Le feuilletage associ\'e $\mathcal{F}(\sigma)$ est de degr\'e~$2;$ il poss\`ede
$7$ points singuliers qui sont exactement les points fixes et
les points d'ind\'etermination de $\sigma.$ On en d\'eduit la:
\end{eg}

\begin{pro}
{\sl Si $f$ est une transformation birationnelle quadratique
g\'en\'erique $\mathcal{F}(f)$ poss\`ede $7$ points singuliers distincts;
quatre sont des points fixes de $f,$ les trois autres sont des points
d'ind\'etermination.}
\end{pro}

\begin{eg}
Consid\'erons maintenant l'involution $f$ de $\Sigma^3$ donn\'ee par
$f=(x_0x_2:x_1x_2:x_0x_1).$ On a
\begin{align*}
& \mathrm{Ind}\hspace{0.1cm} f= \mathrm{Ind}\hspace{0.1cm} f^{-1}=\{(1:0:0),\hspace{0.1cm}(0:1:0),\hspace{0.1cm}(0:0:1)\}.
\end{align*}
L'ensemble des points fixes de $f$ est pr\'ecis\'ement la conique
$\{x_2^2-x_0x_1=0\}$ priv\'ee des deux points $(1:0:0)$ et $(0:1:0)$
qui sont d'ind\'etermination. En fait $f$ est bien d\'efini en
restriction \`a la conique $\overline{\mathrm{Fix}\hspace{0.1cm} f}$ sur laquelle
elle co\"{\i}ncide avec l'identit\'e.

Dans une telle situation nous dirons que $f$ poss\`ede une
\textbf{\textit{courbe de points fixes}}\label{ind43}.

Notons que le feuilletage $\mathcal{F}(f)$ est ici le feuilletage en droites
$\frac{x_1}{x_0}=$ cte dont la seule singularit\'e est le point $(0:0:1).$
\end{eg}

\begin{pro}\label{degfeuil}
{\sl Soient $f$ un \'el\'ement du groupe de \textsc{Cremona} et $\mathcal{F}(f)$
le feuilletage associ\'e \`a $f.$ La transformation $f$ a une
courbe de points fixes si et seulement si $\deg\mathcal{F}(f)<\deg f.$

Si $\mathcal{C}$ est une courbe de points fixes de $f$ alors
$\deg\mathcal{C}\leq 2.$}
\end{pro}

\begin{proof}[\sl D\'emonstration]
On \'ecrit $f$ sous la forme $(f_0:f_1:f_2);$ dire que $\deg\mathcal{F}(f)<\deg f$
c'est dire que la $1$-forme
\begin{align*}
\omega=(x_1f_2-x_2f_1)\mathrm{d}x_0+(x_2f_0-x_0f_2)\mathrm{d}x_1
+(x_0f_1-x_1f_0)\mathrm{d}x_2
\end{align*}
est divisible par un polyn\^ome homog\`ene
$\varphi$ de degr\'e strictement sup\'erieur \`a $1.$ La courbe
$\varphi=0$ est visiblement une courbe de points fixes (un calcul
imm\'ediat montre que $\varphi=0$ ne peut pas \^etre contract\'ee).

On peut \'ecrire $\omega=\varphi\omega',$ la forme $\omega'$
s'annulant sur un ensemble de codimension $2$ (une union de droites).
Comme $\omega$ et donc $\omega'$ annulent le champ radial on a
$1\leq\deg\varphi\leq 2.$
\end{proof}

Les points fixes peuvent \^etre de multiplicit\'e
strictement sup\'erieure \`a $1:$ ainsi la transformation
$(x_0^2+x_1x_2:x_1^2+x_0x_2:x_2^2-x_0x_1)$ compte trois points
d'ind\'etermina\-tion
\begin{align*}
&(1:1:-1),&&\left(1-\mathrm{i}\sqrt{3}: 1+\mathrm{i}
\sqrt{3}:2\right),&&\left(1+\mathrm{i}\sqrt{3}:
1-\mathrm{i}\sqrt{3}:2\right),
\end{align*}
et trois points fixes $(1:0: 0),$
$(0:1:0),$ $(0:0:1);$ on peut v\'erifier que $(0:0:1)$ est double.
La transformation
\begin{align*}
f=(x_0(x_2-x_1):x_1(x_2-x_0):x_2^2+\alpha x_1x_2+\beta x_0x_2-(1+
\alpha+\beta)x_0x_1)
\end{align*}
satisfait
\begin{align*}
&\mathrm{Ind}\hspace{1mm} f=\{(1:0:0),\hspace{1mm}(0:1:0),
\hspace{1mm}(1:1:1)\}, && \mathrm{Fix}\hspace{1mm} f=\{(0:0:1)\}
\end{align*}
et $(0:0:1)$ est de
multiplicit\'e $4.$ Les droites d'\'equation $x_0=x_2,$ $x_1=x_2$ et
$x_1=x_0-x_2$ sont contract\'ees et celles d'\'equation $x_0=0,$ $x_1=0$ et
$x_1=x_0$ sont invariantes. On peut v\'erifier qu'un point fixe peut
aussi \^etre de multiplicit\'e $3.$

En g\'en\'eral un \'el\'ement $f$ du groupe de
\textsc{Cremona} ne pr\'eserve pas son feuilletage associ\'e; il
en est ainsi pour l'involution de \textsc{Cremona}. Ceci signifie
en un certain sens qu'il n'y a pas de lien dynamique entre $f$ et
$\mathcal{F}(f).$
Il y a quelque cas sp\'eciaux int\'eressants. Supposons
que~$f$ pr\'eserve une fibration $\ell$ en droites fibre \`a
fibre, {\it i.e.} $\ell\circ f=\ell;$ on peut se ramener \`a
$\ell=\frac{x_1}{x_0}.$ Alors~$f$ s'\'ecrit
\begin{align*}
& f=(x_0\varphi:x_1\varphi:\psi), &&\deg\psi=\deg\varphi+1.
\end{align*}
Un calcul \'el\'ementaire montre que $\mathcal{F}(f)$ est la fibration
$\frac{x_1}{x_0}=$ cte, et dans ce cas $f$ pr\'eserve le
feuilletage $\mathcal{F}(f).$ On peut d\'emontrer (c'est long et
fastidieux) que dans le cas quadratique c'est la seule
possibilit\'e: si $f\in\mathring{\mathrm{B}}\mathrm{ir}_2$ laisse
$\mathcal{F}(f)$ invariant $f$ pr\'eserve une fibration en droites fibre
\`a fibre.

Toutefois certains invariants de $\mathcal{F}(f),$ en particulier
aux points singuliers, vont poss\'eder une traduction naturelle
pour $f.$

\bigskip

\section{Transformations birationnelles quadratiques et feuilletages}

Lorsque $n=2$ on a $\mathrm{Rat}_2\simeq\mathbb{P}^{17}(\mathbb{C})$ et $\mathcal{F}_2\simeq\mathbb{P}^{14}(\mathbb{C});$
ainsi l'application $f\mapsto\mathcal{F}(f)$ induit une application lin\'eaire not\'ee encore $\mathcal{F}(.)$
de $\mathbb{P}^{17}(\mathbb{C})$ dans $\mathbb{P}^{14}(\mathbb{C}).$ Nous allons d\'emontrer que l'application $\mathcal{F}
(.)$ est dominante en restriction \`a $\mathrm{Bir}_2.$

Pour cel\`a on se donne un feuilletage $\mathcal{F}$ g\'en\'erique au sens
o\`u son ensemble singulier $\mathrm{Sing}\hspace{0.1cm}\mathcal{F}$
est constitu\'e de $7$ points en position g\'en\'erale (pas
d'alignement $3$ \`a $3$). Comme on l'a dit on peut trouver une
transformation quadratique $F=(F_0:F_1:F_2)$ telle que $\mathcal{F}$
soit d\'efini par
\begin{align*}
(x_1F_2-x_2F_1)\mathrm{d}x_0+(x_2F_0-x_0F_2)\mathrm{d}x_1+(x_0F_1-x_1F_0)\mathrm{d}x_2;
\end{align*}
l'ensemble $\mathrm{Sing}\hspace{0.1cm}\mathcal{F}$ est alors pr\'ecis\'ement
l'ensemble des points $m$ pour lesquels il existe $\eta$ dans $\mathbb{C}$
tels que $F(m)=\eta m.$ Choisissons trois points $m_1,$ $m_2$ et
$m_3$ dans $\mathrm{Sing}\hspace{0.1cm}\mathcal{F}.$ Comme les $m_i$ ne sont pas
align\'es on peut trouver $\ell$ une forme lin\'eaire telle
qu'avec les notations ci-dessus $\ell(m_i)=-\eta_i;$ l'application
$f=F+\ell\mathrm{id}$ satisfait $\mathcal{F}(f)= \mathcal{F}(F)=\mathcal{F}$ et les points
$m_i$ sont d'ind\'etermination pour $f.$ En particulier $f$ est
birationnelle. On en d\'eduit que la restriction de~$\mathcal{F}(.)$ \`a
$\mathrm{Bir}_2$ est dominante et par suite \`a fibre
g\'en\'erique finie puisque $\dim\mathrm{Bir}_2=\dim\mathcal{F}_2.$

On sait d'apr\`es \cite{GMK} qu'un
feuilletage de degr\'e $2$ est d\'etermin\'e par la position de
ses $7$ points singuliers et que toute configuration g\'en\'erique
est r\'ealis\'ee. On h\'erite d'une propri\'et\'e analogue pour les
\'el\'ements de $\Sigma^3.$

\begin{pro}
{\sl Une transformation birationnelle
quadratique g\'en\'erique est d\'etermin\'ee par la position de ses points
d'ind\'etermination et de ses points fixes.}
\end{pro}

\begin{proof}[\sl D\'emonstration]
Si $Q_1$ et $Q_2$ ont trois points d'ind\'etermination communs on
peut supposer qu'elles sont du type $A_1\sigma$ et $A_2\sigma;$ si
de plus elles ont les m\^emes points fixes, on a, puisque les sept
points sont en position g\'en\'erale, $\mathcal{F}(Q_1)=\mathcal{F}(Q_2).$ Il existe
donc $\eta$ dans $\mathbb{C}$ et $\ell$ une forme lin\'eaire tels que
$(A_1- \eta A_2)\sigma=\ell.\mathrm{id}.$ Un calcul direct montre
que $\ell=0$ et donc $Q_1=Q_2.$
\end{proof}

On en d\'eduit que la fibre g\'en\'erique de l'application
\begin{align*}
\mathcal{F}(.)_{|\mathrm{Bir}_2}\hspace{0.1cm}\colon\hspace{0.1cm}\mathrm{Bir}_2\to\mathcal{F}_2
\end{align*}
a exactement $35=C_7^3$ points; il s'agit en effet de d\'esigner
parmi $7$ points $3$ points d'ind\'etermination.

Remarquons que $\mathcal{F}(.)$ a des fibres non finies.
Par exemple si $\widetilde{\mathcal{F}}$ est le feuilletage donn\'e par
$x_1\mathrm{d}x_0-x_0\mathrm{d}x_1$ qui d\'efinit la fibration $x_1/x_0=$ cte, toutes les
transformations qui laissent invariante fibre \`a fibre cette
fibration sont dans $\mathcal{F}^{-1}(\widetilde{\mathcal{F}});$ en particulier toutes
les transformations $A\sigma$ de la forme
\begin{align*}
(*x_1x_2+*x_0x_2:
*x_1x_2+*x_0x_2:*x_1x_2+*x_0x_2+*x_0x_1)
\end{align*}
sont dans cette fibre.

Il y a aussi des fibres multiples {\it i.e.} ayant moins
de $35$ points. C'est le cas de $\mathcal{F}^{-1}(\mathcal{F}
(\sigma))$ que nous allons d\'eterminer. Remarquons qu'un
\'el\'ement de $\mathcal{F}^{-1}(\mathcal{F}(\sigma))$ est n\'ecessairement dans
$\Sigma^3,$ {\it i.e.} du type $A\sigma B.$ Ainsi les \'el\'ements $Q_i$
de $\mathcal{F}^{-1}(\mathcal{F}(\sigma))$ sont du type
\begin{align*}
& Q_i=A_i\sigma B_i=\sigma+(x_0\ell_i:x_1\ell_i:x_2\ell_i),&& A_i,
\hspace{1mm}B_i\in\mathrm{PGL}_3(\mathbb{C}), \hspace{1mm}
\ell_i\text{ forme lin\'eaire.}
\end{align*}

On obtient donc la:
\begin{pro}
{\sl L'ensemble $\mathcal{F}^{-1}(\mathcal{F}(\sigma))$ est constitu\'e des transformations
birationnelles quadratiques suivantes
\begin{align*}
&Q_0=\sigma,&&
Q_1=\sigma+(x_0+x_1+x_2)\mathrm{id},&&
Q_2 =\sigma+(x_0-x_1-x_2)\mathrm{id},
\end{align*}
\begin{align*}
&Q_3=\sigma+(-x_0-x_1+x_2)\mathrm{id},&&
Q_4=\sigma+(-x_0+x_1-x_2)\mathrm{id}.
\end{align*}}
\end{pro}

Remarquons que les $Q_i$ sont conjugu\'es entre eux par des
\'el\'ements qui commutent avec~$\sigma$

\begin{small}
\begin{align*}
& Q_2=(-x_0:x_1:x_2)Q_1(-x_0:x_1:x_2), && Q_3=(x_0:x_1:-x_2)Q_1(x_0: x_1:-x_2),
\end{align*}
\begin{align*}
&
Q_4=(x_0:-x_1:x_2)Q_1(x_0:-x_1:x_2).
\end{align*}
\end{small}

La configuration des points fixes et d'ind\'etermination de
$\sigma$ est sch\'ematis\'ee par la figure suivante:
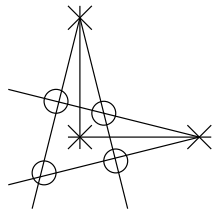
\begin{figure}[H]
\begin{center}
\input{fixind.pstex_t}
\end{center}
\end{figure}

\noindent o\`u les croix repr\'esentent les points d'ind\'etermination
et les cercles les points fixes.

\begin{rem}
Soit $f$ dans $\mathrm{Bir}_2$ telle que $\deg\mathcal{F}(f)=2.$ Supposons que trois
points fixes de $f$ soient align\'es, par exemple situ\'es sur la droite
$x_2=0;$ alors $\mathcal{F}(f)$ a trois
points singuliers sur $x_2=0$ qui est donc invariante par $\mathcal{F}(f).$
Ceci implique que $x_2=0$ est aussi invariante par~$f;$ la restriction
de $f$ \`a $x_2=0$ s'identifie \`a une transformation de \textsc{M\"{o}bius}
avec trois points fixes c'est donc l'identit\'e. Par suite la droite $x_2=0$
est une courbe de points fixes de $f.$
\end{rem}

\bigskip

Revenons \`a la figure ci-dessus. Si $f$ est un \'el\'ement de
$\mathcal{F}^{-1}(\mathcal{F}(\sigma))$ alors
\begin{align*}
\mathrm{Fix}\hspace{0.1cm} f\cup\mathrm{Ind}\hspace{0.1cm} f=\mathrm{Fix}\hspace{0.1cm}
\sigma\cup\mathrm{Ind}\hspace{0.1cm}\sigma
\end{align*}
avec les contraintes suivantes: trois points
de $\mathrm{Ind}\hspace{0.1cm} f$ (resp. $\mathrm{Fix}\hspace{0.1cm} f$) ne peuvent \^etre align\'es.
On retrouve ainsi $\#\mathcal{F}^{-1}(\mathcal{F}(\sigma))=5.$

Des arguments du m\^eme ordre montrent que $\#\mathcal{F}^{-1}(\mathcal{F}(\rho))=15.$
Notons que dans $\mathcal{F}^{-1}(\mathcal{F}(\tau))$ on trouve \`a la fois des
\'el\'ements de $\Sigma^1$ (par exemple $\tau,$ $\tau-x_2\mathrm{id},$
$\tau+x_2\mathrm{id}$) et des \'el\'ements de $\Sigma^2$ comme $\tau+
(-2x_1+x_2)\mathrm{id}.$ Ainsi la donn\'ee du feuilletage $\mathcal{F}(f)$ ne permet pas
la d\'etermination de l'orbite g.d. de $f.$

\section{Relations de type \textsc{Lefschetz}, \textsc{Baum}-\textsc{Bott}}\hspace{0.1cm}

Certains invariants de conjugaison initialement
introduits pour les feuilletages s'adaptent aux transformations
birationnelles. En $2004$ \textsc{Guillot} a montr\'e le:

\begin{thm}[\cite{Gu}]\label{adolfog}
{\sl Soit $B$ un polyn\^ome sur l'espace des transformations lin\'eaires de $\mathbb{C}^n,$
homog\`ene de degr\'e $\nu\leq n$ et  invariant par conjugaison. Soit $f\hspace{0.1cm}
\colon\hspace{0.1cm}\mathbb{P}^n(\mathbb{C})\to\mathbb{P}^n(\mathbb{C})$ une fonction
rationnelle de degr\'e alg\'ebrique $d\geq 2$ d\'efinie en tout point et dont les points
fixes sont isol\'es et simples. Alors
\begin{align*}
\sum_{x\in\mathrm{Fix}f}\frac{B(Df_{(x)}-\mathrm{id})}{\mathrm{det}\hspace{0.1cm}(Df_{(x)}-
\mathrm{id})}
\end{align*}
est un complexe qui ne d\'epend que du polyn\^ome
$B$ et des entiers $n$ et $d.$}
\end{thm}

Lorsque $B=1$ on retrouve un cas particulier du
th\'eor\`eme de \textsc{Lefschetz} suivant.

\begin{thm}[\cite{Le, GrHa}]
{\sl Soient $M$ une vari\'et\'e complexe compacte de dimension $n$
et $f$ une application holomorphe de $M$ dans elle-m\^eme dont les
points fixes sont simples et isol\'es. Alors
\begin{align*}
\sum_{x\in\mathrm{Fix}f}\frac{1}{\mathrm{det}\hspace{0.1cm} (Df_{(x)}-\mathrm{id})}=(-1)^nL(f,0)
\end{align*}
o\`u $L(f,0)$ est le nombre de \textsc{Lefschetz} qui ne d\'epend
que de l'action de $f$ sur la cohomologie de \textsc{Dolbeault} de
$M.$}
\end{thm}

Pour un endomorphisme $f$ de $\mathbb{P}^n(\mathbb{C})$ v\'erifiant les
hypoth\`eses pr\'ec\'edentes, la quantit\'e $(-1)^nL(f,0)$ vaut
$(-1)^n$ (\emph{voir} \cite{GrHa}).

Comme l'indique \textsc{Guillot} le th\'eor\`eme
\label{adolfo} s'adapte aux transformations birationnelles ayant
des points d'ind\'etermination simples pourvu que le degr\'e de
$B$ soit inf\'erieur ou \'egal \`a $1.$ En effet en un tel point
on peut consid\'erer les valeurs propres comme infinies. On
obtient en particulier le:

\begin{thm}\label{adolfo2}
{\sl Soit $f$ un \'el\'ement de $\Sigma^3$ ayant quatre points
fixes en position g\'en\'erale; alors
\begin{equation}\label{rel1}
\displaystyle\sum_{x\in\mathrm{Fix}f} \frac{\mathrm{tr}(Df_{(x)}
-\mathrm{id})}{\mathrm{det}\hspace{0.1cm}(Df_{(x)}-\mathrm{id})}=-4 \hspace{6mm}\text{et}\hspace{6mm}
\displaystyle\sum_{x\in\mathrm{Fix}f}\frac{1}{\mathrm{det}\hspace{0.1cm}(Df_{(x)}
-\mathrm{id})}=1.
\end{equation}}
\end{thm}

\begin{proof}[\sl D\'emonstration]
Il suffit de le v\'erifier pour l'involution $\sigma$ de
\textsc{Cremona}.
\end{proof}

\begin{rem}
Si $f$ satisfait les hypoth\`eses du th\'eor\`eme, en chaque point fixe
$m,$ on a deux valeurs
propres $\delta_1(m),$ $\delta_2(m)$ pour $Df_{(m)}.$ Il y a donc en
tout $8$ valeurs propres qui sont reli\'ees par les relations
(\ref{rel1}). Ces valeurs propres sont invariantes
par conjugaison dynamique et on sait que l'orbite dynamique d'un
\'el\'ement g\'en\'erique de $\Sigma^3$ est de codimension $6$
(\emph{voir} Chapitre \ref{gen}). Nous verrons
que la donn\'ee de trois couples de valeurs propres g\'en\'eriques
peut \^etre r\'ealis\'ee (d\'emonstration du Th\'eor\`eme \ref{conjbir}).
\end{rem}

\bigskip

De plus on a le:

\begin{thm}[\cite{Gu2}]
{\sl Les relations (\ref{rel1}) sont \og les
seules relations\fg\hspace{0.1cm} entre les valeurs propres aux points fixes
d'un \'el\'ement g\'en\'erique de $\Sigma^3.$}
\end{thm}

La multiplicit\'e des points fixes d'une transformation
birationnelle se lit sur celle du feuilletage associ\'e. Par
exemple on a la proposition suivante dont la preuve est
\'el\'ementaire.

\begin{pro}\label{deg}
{\sl Soient $f$ une transformation birationnelle de $\mathbb{P}^2(\mathbb{C})$ et
$m$ un point fixe de~$f.$ On note $\delta_1$ et $\delta_2$ les
valeurs propres de la partie lin\'eaire de $f$ en $m.$ Si les
$\delta_i$ sont distinctes de $1,$ le point $m$ est une
singularit\'e non d\'eg\'en\'er\'ee ({\it i.e.} de multiplicit\'e
$1$) du feuilletage $\mathcal{F}(f)$ associ\'e \`a $f.$}
\end{pro}

En appliquant le Th\'eor\`eme \ref{adolfog} pour
$B=\mathrm{tr},$ on obtient le:

\begin{cor}
{\sl Soit $f$ une transformation de $\Sigma^3$ ayant quatre points
fixes distincts. Les valeurs propres de la partie lin\'eaire
de $f$ aux points fixes ne peuvent pas toutes \^etre de module
strictement inf\'erieur \`a $1;$ autrement dit les points fixes ne
peuvent pas tous \^etre des attracteurs.}
\end{cor}

\begin{proof}[\sl D\'emonstration]
D'apr\`es le Th\'eor\`eme \ref{adolfo2}, on a
\begin{equation}\label{adolpho}
\sum_{m\in\mathrm{Fix} f}\frac{\mathrm{tr}(Df_{(m)}-\mathrm{id})}
{\mathrm{det}\hspace{0.1cm}(Df_{(m)}-\mathrm{id})}=-4
\end{equation}
Notons $\delta_1(m),$ $\delta_2(m)$ les valeurs propres de
$Df_{(m)};$ on remarque que (\ref{adolpho}) se r\'e\'ecrit
\begin{equation}\label{2mimi}
\sum_{m\in\mathrm{Fix} f}\frac{1}{1-\delta_1(m)}+\frac{1}{1-\delta_2(m)}
=4.
\end{equation}
Pour $j=1,2$ posons $\eta_j(m):=(1-\delta_j(m))^{-1}$ puis
$\eta_j(m)=y_j(m)+\mathrm{i}z_j(m).$ L'\'egalit\'e (\ref{2mimi})
est \'equivalente \`a
\begin{align*}
&\displaystyle\sum_{m\in\mathrm{Fix}f}\hspace{0.1cm}\sum_{j=1}^2 y_j(m)=4,
&&\displaystyle\sum_{m\in\mathrm{Fix}f}\hspace{0.1cm}\sum_{j=1}^2 z_j(m)=0
\end{align*}
Supposons que pour tout point fixe $m$ de
$f$ on ait $|\delta_1(m)|,|\delta_2(m)|\leq 1;$ les $1-\delta_j(m)$
sont alors dans le disque de rayon $1$ centr\'e en $1.$ Par
l'application $w\mapsto 1/w$ ce disque est transform\'e en le
demi-plan $\mathrm{Re}\hspace{0.1cm} w\geq 1/2.$\label{not27ab} De sorte que
\begin{align*}
& y_1(m)\geq\frac{1}{2}, && y_2(m)\geq\frac{1}{2}, && \forall\hspace{0.1cm} m\in
\mathrm{Fix}\hspace{0.1cm} f;
\end{align*}
il en r\'esulte que $|\delta_1(m)|=|\delta_2(m)|=1$ pour tout point fixe
$m$ de $f.$
\end{proof}

Par contre on peut produire des exemples o\`u les points fixes
sont soit contractants, soit dilatants ({\it i.e.} aucun des
points fixes n'est hyperbolique). De m\^eme il existe des exemples
o\`u toutes les valeurs propres sont de module $1.$\bigskip

Soit $\mathscr{F}$ un feuilletage sur $\mathbb{P}^2(\mathbb{C}).$
Pla\c{c}ons-nous dans $\mathbb{C}^2,$ carte affine de $\mathbb{P}^2(
\mathbb{C});$ on
d\'esigne par~$\chi$ un champ polynomial d\'efinissant
$\mathscr{F}$ et par $m$ un point singulier de $\mathscr{F}.$ La
matrice jacobienne de $\chi$ au point $m,$ not\'ee $\mathrm{jac}\hspace{0.1cm}\chi(m),$
appartient \`a $M_2(\mathbb{C})$ et admet $2$ valeurs propres
$\eta_1(m)$ et $\eta_2(m).$ La quantit\'e
\begin{align*}
\mathrm{BB}(\mathscr{F}(m))=
\frac{\mathrm{tr}^2(\mathrm{jac}\hspace{0.1cm}\chi(m))}{\mathrm{det}(\mathrm{jac}
\hspace{0.1cm}\chi(m))}=\frac{\eta_1(m)}{\eta_2(m)}+\frac{\eta_2(m)}{\eta_1(m)}+2
\end{align*}\label{not27b}
 ne d\'epend pas du choix du champ (\cite{GMLu});
 c'est l'\textbf{\textit{indice de \textsc{Baum}-\textsc{Bott}}}\label{ind43b} de $\mathscr{F}$
 au point $m$ (\emph{voir} \cite{BB}).

\begin{lem}
{\sl Soient $f,$ $g$ deux transformations rationnelles
quadratiques g\'en\'eriques et $\mathcal{F}(f),$ $\mathcal{F}(g)$ les
feuilletages induits par $f$ et $g.$ Soit $m$ (resp. $p$) un point
fixe pour $f$ (resp.~$g$). Supposons que $f_{,m}$ et $g_{,p}$
soient holomorphiquement conjugu\'es; alors
\begin{align*}
\mathrm{BB}(\mathcal{F}(f)(m))=\mathrm{BB}(\mathcal{F}(g)(p)).
\end{align*}}
\end{lem}

La condition de g\'en\'ericit\'e est en particulier n\'ecessaire
pour que les indices de \textsc{Baum}-\textsc{Bott} soient bien
d\'efinis.

\begin{proof}[\sl D\'emonstration]
C'est un calcul imm\'ediat. Quitte \`a faire une translation on
peut supposer que $m=p=(0:0:1).$ Dans la carte affine $x_2=1,$ la
transformation $f$ s'\'ecrit
\begin{align*}
\left(a_0x_0+a_1x_1,b_0x_0+b_1x_1
\right)+\text{termes de plus haut degr\'e}
\end{align*}
ce qui en d'autres termes dit que $f$ est du type suivant

\begin{footnotesize}
\begin{align*}
& (x_2^{n-1}(a_0x_0+a_1x_1)+x_2^{n-2}q_2+\ldots+q_n:x_2^{n-1}
(b_0x_0+b_1x_1)+x_2^{n-2}r_2+\ldots+
r_n:x_2^n+ x_2^{n-1}s_1+\ldots+s_n)
\end{align*}
\end{footnotesize}
avec $q_i,\hspace{0.1cm}r_i,\hspace{0.1cm} s_i$ dans $\mathbb{C}[x_0,x_1]_i.$

Ainsi
\begin{align*}
&\mathrm{tr}(\mathrm{jac}\hspace{0.1cm} f_{(0)})=a_0+a_1,&& \mathrm{det}(\mathrm{jac}\hspace{0.1cm} f_{(0)})=a_0
b_1-a_1b_0.
\end{align*}
Le feuilletage $\mathcal{F}(f)$ est d\'ecrit dans la carte affine
$x_2= 1$ par
\begin{align*}
\omega=(x_1-b_0x_0-b_1x_1)\mathrm{d}x_0+(a_0x_0+a_1x_1-x_0)\mathrm{d}x_1+ \text{termes de plus
haut degr\'e}.
\end{align*}
Ici encore intervient la condition de g\'en\'ericit\'e:
si $\omega=a\mathrm{d}x_0+b\mathrm{d}x_1$ alors $\text{pgcd}(a,b)=1.$ Il en r\'esulte que

\begin{eqnarray}
\mathrm{BB}(\mathcal{F}(f)(0))&=&\frac{(a_0+b_1-2)^2}{(a_0-1)
(b_1-1)-a_1b_0}\nonumber\\
&=& \left(\mathrm{tr}(\mathrm{jac}\hspace{0.1cm} f_{(0)})-2\right)^2\left(\mathrm{det}(\mathrm{jac}\hspace{0.1cm}
f_{(0)})-\mathrm{tr}(\mathrm{jac}\hspace{0.1cm} f_{(0)})+1\right)^{-1}.\nonumber
\end{eqnarray}
Puisque $f_{,0}$ et $g_{,0}$ sont holomorphiquement conjugu\'es,
on a les \'egalit\'es
\begin{align*}
&\mathrm{tr}(\mathrm{jac}\hspace{0.1cm} f_{(0)})=\mathrm{tr}(\mathrm{jac}\hspace{0.1cm} g_{(0)}),
&&\mathrm{det}(\mathrm{jac}\hspace{0.1cm}f_{(0)})=\mathrm{det}(\mathrm{jac}\hspace{0.1cm} g_{(0)});
\end{align*}
par suite $\mathrm{BB}(\mathcal{F}(f)(0))=\mathrm{BB}(\mathcal{F}
(g)(0)).$
\end{proof}

Pour tout feuilletage g\'en\'erique $\mathscr{F}$ de degr\'e $d$ sur $\mathbb{P}^2(\mathbb{C}),$
on a (\cite{BB})
\begin{align*}
\sum_{m\in\mathrm{Sing}\mathscr{F}}\mathrm{BB}((\mathscr{F})(m))=
(d+2)^2.
\end{align*}
On en d\'eduit une formule pour tout $Q$ g\'en\'erique dans $\mathrm{Rat}_2$
tel que $\deg\mathcal{F}(Q)=2$

\begin{eqnarray}\label{bbb}
\sum_{m\in\Omega}
\frac{\mathrm{tr}^2(\mathrm{jac}\hspace{0.1cm}\chi_Q(m))} {\mathrm{det}(\mathrm{jac}\hspace{0.1cm}\chi_Q(m))}=16
\end{eqnarray}
o\`u $\Omega=\mathrm{Ind}\hspace{0.1cm} Q\cup\mathrm{Fix}\hspace{0.1cm} Q$
et $\chi_Q$ d\'esigne un champ local d\'efinissant le feuilletage
$\mathcal{F}(Q)$ en $m.$

De la m\^eme mani\`ere que pr\'ec\'edemment cette formule est
valable pour tout \'el\'ement g\'en\'erique de $\Sigma^3.$

\begin{rem}
La formule (\ref{bbb}) n'est pas valable pour les \'el\'ements $Q$
de $\Sigma^2$ et $\Sigma^1:$ au point base $m$ de multiplicit\'e
strictement sup\'erieure \`a $1,$ la quantit\'e $\frac{\mathrm{
tr}^2(\mathrm{jac}\hspace{0.1cm}\chi_Q(m))}{\mathrm{det}(\mathrm{jac}\hspace{0.1cm}\chi_Q(m))}$ n'est pas bien
d\'efinie.
\end{rem}

\section{Conjugaison birationnelle entre \'el\'ements de
$\mathrm{PGL}_3(\mathbb{C})$}\hspace{0.1cm}

Le probl\`eme de la classification des automorphismes
de $\mathbb{P}^2(\mathbb{C})$ modulo conjugaison birationnelle est trait\'e dans
\cite{Bl2}; l'auteur y d\'emontre que tout automorphisme de $\mathbb{P}^2(
\mathbb{C})$
est birationnellement conjugu\'e \`a une transformation diagonale
$(\alpha_0x_0:\alpha_1x_1:\alpha_2x_2)$ ou presque diagonale
$(x_0:x_0+x_1:\alpha x_2).$ Il montre par exemple que deux
transformations diagonales birationnellement conjugu\'ees le
sont par une \textbf{\textit{transformation monomiale}}\label{ind19a}, {\it i.e.} une
transformation de la forme
\begin{align*}
&(\gamma x_0^px_1^q,\delta x_0^rx_1^s),&&\gamma,\hspace{0.1cm}\delta\in\mathbb{C}^*,\hspace{0.1cm}
\left[\begin{array}{cc}p & q\\ r & s\end{array}\right]\in\mathrm{GL}_2(\mathbb{Z}).
\end{align*}
Voici un cas sp\'ecial (g\'en\'erique).

\begin{pro}
{\sl Soit $A$ un automorphisme diagonal de $\mathbb{P}^2(\mathbb{C});$ on suppose que
\begin{align*}
\overline{\langle A\rangle}^{\hspace{0.1cm} \mathsf{Z}}=\{(\alpha x_0,\beta
x_1)\hspace{0.1cm}|\hspace{0.1cm} \alpha,\hspace{0.1cm} \beta \in\mathbb{C}^*\}.
\end{align*}
Soit $B$ un automorphisme de $\mathbb{P}^2(\mathbb{C})$ birationnellement
conjugu\'e \`a $A.$ Alors
\begin{itemize}
\item $B$ est lin\'eairement conjugu\'e \`a un \'el\'ement diagonal
$(cx_0,dx_1);$

\item $A$ et $(cx_0,dx_1)$ sont conjugu\'es via une transformation
monomiale.
\end{itemize}}
\end{pro}

\begin{proof}[\sl D\'emonstration]
\'Ecrivons $A$ sous la forme $(ax_0,bx_1).$
Il suffit d'\'etudier les cas o\`u
\begin{align*}
& B=(cx_0,dx_1) && \text{ et } && B=(cx_0,x_0+cx_1).
\end{align*}
Supposons que $B=(cx_0,dx_1);$ soit $f=(f_1,f_2)$ une transformation
birationnelle conjuguant $A$ \`a $B.$ On~a
\begin{align*}
& (f_1,f_2)(a^nx_0,b^nx_1)=(c^nf_1,d^nf_2), && \forall\hspace{0.1cm} n\in\mathbb{Z}.
\end{align*}
Puisque $\overline{
\langle A\rangle}^{\hspace{0.1cm} \mathsf{Z}}=\{(\alpha x_0,\beta x_1)\hspace{0.1cm}|\hspace{0.1cm}
\alpha,\hspace{0.1cm}\beta\in\mathbb{C}^*\},$ on a
\begin{align*}
&(f_1,f_2)(\alpha x_0,\beta x_1)=(\eta(\alpha,\beta)f_1,
\zeta(\alpha,\beta)f_2), && \forall\hspace{0.1cm}\alpha,\hspace{0.1cm}\beta\in\mathbb{C}^*.
\end{align*}
Un calcul \'el\'ementaire (du m\^eme type que celui qui suit) montre
que
\begin{align*}
&(f_1,f_2)=(\gamma x_0^px_1^q,\delta x_0^rx_1^s), &&\gamma,\hspace{0.1cm}
\delta\in\mathbb{C}^*,\hspace{0.1cm}\left[\begin{array}{cc}p & q\\ r & s
\end{array}\right]\in\mathrm{GL}_2(\mathbb{Z}).
\end{align*}

\bigskip

Reste \`a consid\'erer le cas o\`u $B=(cx_0,cx_1+ x_0).$
Comme l'adh\'erence de \textsc{Zariski} $\overline{\langle
A\rangle}^{\hspace{0.1cm} \mathsf{Z}}$ est de dimension $2$ l'\'egalit\'e
$f\circ B^n=A^n\circ f$ implique que l'adh\'erence de
\textsc{Zariski} $\overline{\langle B\rangle}^{\hspace{0.1cm} \mathsf{Z}}$
est aussi de dimension $2.$ Il en r\'esulte que
\begin{align*}
\overline{\langle B\rangle}^{\hspace{0.1cm} \mathsf{Z}}=\{(sx_0,sx_0+tx_1)\hspace{0.1cm}|\hspace{0.1cm} s\in\mathbb{C}^*,
\hspace{0.1cm} t\in\mathbb{C}\};
\end{align*}
par suite on a
\begin{align*}
& f(sx_0,sx_1+tx_0)=(\alpha(s,t)x_0,\beta(s,t)x_1) \phi, &&\forall\hspace{0.1cm} s\in\mathbb{C}^*,
\hspace{0.1cm} t\in\mathbb{C}.
\end{align*}
Consid\'erons la premi\`ere composante
\begin{equation}\label{com}
f_1(sx_0,sx_1+tx_0)=\alpha(s,t)f_1(x_0,x_1);
\end{equation}
comme le membre de gauche est d\'erivable par rapport \`a $s$
celui de droite aussi. On obtient en d\'erivant par rapport \`a $s$ et
en \'evaluant en $(s,t)=(1,0)$
\begin{align*}
x_0\frac{\partial f_1}
{\partial x_0}+ x_1\frac{\partial f_1}{\partial x_1}=\frac{\partial
\alpha}{\partial s}f_1;
\end{align*}
il s'en suit que $f_1$
est homog\`ene. D\'erivons (\ref{com}) par rapport \`a $t$
\begin{align*}
x_0 \frac{\partial f_1}{\partial x_1}(sx_0,sx_1+tx_0)=\frac{\partial
\alpha}{\partial t}(s,t)f_1(x_0,x_1).
\end{align*}
Puisque $f_1$ est homog\`ene on a pour $t=0$
\begin{align*}
x_0\frac{\partial f_1}{\partial x_1}
(sx_0,sx_1)=s^\nu x_0\frac{\partial f_1}
{\partial x_1}(x_0,x_1)=\frac{\partial\alpha}{\partial t}(s,0)f_1(x_0,x_1)
\end{align*}
autrement dit
\begin{align*}
\frac{\partial f_1/\partial x_1}{f_1}=\frac{
\text{cte}}{x_0}.
\end{align*}
On constate alors que n\'ecessairement
$\partial f_1/\partial x_1\equiv 0,$ {\it i.e.} $f_1=f_1(x_0).$ Un
raisonnement analogue conduit \`a $f_2=f_2
(x_0)$ autrement dit $f=f(x_0)$ ce qui est impossible.
\end{proof}

Dans \cite{Se2} Serre \'etablit un r\'esultat analogue
en toute dimension.

\begin{thm}[\cite{Se2}]
{\sl Soient $A$ et $B$ deux parties du tore $T$ de
$\mathrm{PGL}_{n+1}(\mathbb{C})$ et soit $f$ un \'el\'ement
de $\mathrm{Bir}(\mathbb{P}^n(\mathbb{C}))$ tel que $fAf^{-1}=B.$
Alors $A$ et $B$ sont conjugu\'ees via une transformation
monomiale.}
\end{thm}

\section{Conjugaison birationnelle entre \'el\'ements de
$\Sigma^3$}\hspace{0.1cm}

\'Evidemment deux transformations birationnellement conjugu\'ees ne sont
pas n\'ecessairement lin\'eairement conjugu\'ees. Par exemple

\begin{small}
\begin{align*}
&Q_1=(2(x_0+x_1)x_2:x_1(x_2-x_0):(x_0+x_1)(x_2-x_0))&& \text{et} && Q_2
=(x_0-x_2:2 (x_1+x_0):2x_0)
\end{align*}
\end{small}

\noindent sont birationnellement conjugu\'ees par
$\sigma;$ par contre elles ne peuvent pas l'\^etre via un
automorphisme de $\mathbb{P}^2(\mathbb{C})$ (elles n'ont pas m\^eme degr\'e).

N\'eanmoins les transformations $A\sigma$ et $\sigma A$
sont birationnellement conjugu\'ees via $\sigma$ et lin\'eairement
conjugu\'ees via $A^{-1}.$ Ainsi pour tout \'el\'ement $Q$ de
$\Sigma^3$ il existe une transformation de $\Sigma^3,$ en g\'en\'eral distincte
de $Q$, birationnellement conjugu\'ee \`a $Q$ par une
transformation quadratique mais aussi par un automorphisme de
$\mathbb{P}^2(\mathbb{C}).$\bigskip

\begin{pro}\label{tracedet}
{\sl Soient $A$ et $B$ deux automorphismes de $\mathbb{P}^2(\mathbb{C})$ tels que les
transformations $A\sigma$ et $B\sigma$ aient quatre points fixes
en position g\'en\'erale. Supposons que $A\sigma$ (resp.
$B\sigma$) ne laisse pas de courbe rationnelle invariante. Si
$A\sigma$ et $B\sigma$ sont birationnellement conjugu\'ees par~$\varphi,$ 
alors $\varphi$ est un diff\'eomorphisme local
holomorphe en chaque point fixe de $A\sigma.$ De plus, \`a
conjugaison lin\'eaire pr\`es, on a
\begin{itemize}
\item $\mathrm{Fix}(A\sigma)=\mathrm{Fix}(B\sigma);$

\item $\varphi_{|\mathrm{Fix}(A\sigma)=\mathrm{Fix}(B\sigma)}=
\mathrm{id};$

\item et pour tout point fixe $m$ de $A\sigma$
\begin{align*}
& \mathrm{tr}(\mathrm{jac}(A\sigma)_{(m)})=\mathrm{tr}(\mathrm{jac}(B\sigma)_{(m)}),
&&\mathrm{det}(\mathrm{jac}(A\sigma)_{(m)})=\mathrm{det}(\mathrm{jac}(B\sigma)_{(m)}).
\end{align*}
\end{itemize}}
\end{pro}

\begin{rem}
En particulier si chaque matrice $A$ et $B$ a ses coefficients
$\mathbb{Q}$-alg\'ebriquement ind\'ependants, alors, comme nous le verrons au
Chapitre \ref{dyndyn}, les transformations $A\sigma$ et $B\sigma$
satisfont les conclusions de la Proposition \ref{tracedet}.
\end{rem}

\begin{proof}[\sl D\'emonstration]
Soit $m$ un point fixe pour $A\sigma;$ alors $\varphi(m)=
B\sigma\varphi(m).$ Par hypoth\`ese $B\sigma$ ne fixe pas de
courbe rationnelle; il en r\'esulte que $m$ n'est pas
d'ind\'etermination pour $\varphi.$ En particulier $\varphi$ est
holomorphe en chaque point fixe de $A\sigma$ et l'image par
$\varphi$ d'un point fixe de $A\sigma$ est un point fixe de
$B\sigma.$ De m\^eme $\varphi^{-1}$ est holomorphe aux points
fixes de $B\sigma.$ Par suite
\begin{align*}
\varphi(\mathrm{Fix}(A\sigma))=\mathrm{Fix}(B\sigma)
\end{align*}
et $\varphi$ est un diff\'eomorphisme local en chaque point de
$\mathrm{Fix}(A\sigma).$ Comme deux quadruplets de points en
position g\'en\'erale s'\'echangent par un automorphisme, on
obtient le r\'esultat annonc\'e.
\end{proof}

Pour une transformation de type $A\sigma$ ayant quatre points
fixes en position g\'en\'erale on obtient, en normalisant la
position des points fixes, des formes normales qui s'av\`erent
utiles.

\begin{pro}\label{posgal}
{\sl Soit $A$ un automorphisme de $\mathbb{P}^2(\mathbb{C})$ tel que les points fixes et les
points d'ind\'etermination de $A\sigma$ soient en position g\'en\'erale.
Alors $A\sigma$ est lin\'eairement conjugu\'e \`a une
transformation de la forme

\begin{small}
\begin{align*}
(x_0^2+A_0 x_1x_2+B_0x_0x_2+C_0x_0x_1:x_1^2+A_1x_1x_2+B_1x_0x_2+
C_1x_0x_1:x_2^2+A_2x_1x_2+B_2x_0x_2+C_2x_0x_1)
\end{align*}
\end{small}

\noindent les $A_i,$ $B_i,$ $C_i$ \'etant des \'el\'ements de $\mathbb{C}$
satisfaisant

\begin{small}
\begin{align*}
&A_1B_0-A_0B_1\not=0, && A_2=\frac{B_0-A_0C_1}{A_1B_0-A_0B_1}, && B_2=
\frac{A_1-B_1C_0}{A_1B_0-A_0B_1}, &&
C_2=\frac{1-C_0C_1}{A_1B_0-A_0B_1}.
\end{align*}
\end{small}}
\end{pro}

\begin{proof}[\sl D\'emonstration]
Soit $f$ une transformation lin\'eairement conjugu\'ee \`a
$A\sigma$ de sorte que
\begin{align*}
(1:0:0),\hspace{0.1cm}(0:1:0),\hspace{0.1cm}(0:0:1)\in\mathrm{Fix}\hspace{0.1cm} f;
\end{align*}
alors $f$ est du type $(Q_0:Q_1:Q_2)$ o\`u
\begin{align*}
&Q_0=\eta_0x_0^2+A_0x_1x_2+B_0x_0x_2+C_0x_0x_1,&&Q_1=\eta_1
x_1^2+A_1x_1x_2 + B_1x_0x_2+C_1x_0x_1,
\end{align*}
\begin{align*}
Q_2=\eta_2x_2^2+A_2x_1x_2+B_2x_0x_2+C_2x_0x_1.
\end{align*}
Comme les points fixes
sont distincts, $\eta_0\eta_1\eta_2$ est non nul et, \`a conjugaison par homoth\'etie
pr\`es, on peut supposer que $\eta_0=\eta_1=\eta_2=1.$ Consid\'erons
l'intersection $\mathcal{C}$ des deux coniques $Q_0=0$ et $Q_1=0;$ on
remarque que $\mathcal{C}$ contient le point $(0:0:1)$ et les trois
points d'ind\'etermination $p_1$, $p_2$ et $p_3$ de $f.$ Comme ces
quatre points sont en position g\'en\'erale les deux coniques sont
transversales. Par suite si $P$ est un polyn\^ome homog\`ene de
degr\'e $\nu$ s'annulant sur $\mathcal{C},$ alors $P$ s'\'ecrit
\begin{align*}
&D_0Q_0+D_1Q_1, && D_i\in\mathbb{C}[x_0,x_1,x_2]_{\nu-2}.
\end{align*}
On remarque que $x_0Q_2$ s'annule en $(0:0:1)$ et en les
$p_i;$ il en r\'esulte que $x_0Q_2$ s'\'ecrit
\begin{align*}
(\alpha_0x_0+\beta_0x_1+\gamma_0x_2)Q_0+(\alpha_1x_0+\beta_1
x_1+\gamma_1x_2)Q_1.
\end{align*}

Remarquons que $(A_0,B_0)$ est non nul: si $A_0$ et $B_0$
\'etaient nuls, l'intersection de $Q_0=0$ et $Q_1=0$ en $(0:0:1)$ serait
multiple ce qui contredirait l'hypoth\`ese sur les points fixes;
de m\^eme $(A_1,B_1)\not=(0,0).$ Ainsi $B_0A_1-B_1A_0=0$ si
et seulement s'il existe $t$ dans $\mathbb{C}^*$ tel que $(A_1,B_1)=t(A_0,B_0);$
si tel \'etait le cas l'intersection de $tQ_0-Q_1=0$ et $Q_1=0$
en $(0:0:1)$ serait
multiple. Il s'en suit que $B_0A_1-B_1A_0\not=0.$

\`A partir de
\begin{align*}
& (\alpha_0x_0+\beta_0x_1+\gamma_0
x_2)Q_0+(\alpha_1x_0+\beta_1x_1+\gamma_1x_2)Q_1=x_2^2+A_2
x_1x_2+B_2x_0x_2+C_2x_0x_1
\end{align*}
on obtient
\begin{align*}
& A_2=\delta(B_0-A_0C_1),&& B_2=\delta(A_1-B_1C_0), &&
C_2=\delta(1-C_0C_1)
\end{align*}
avec $\delta=(B_0A_1-B_1A_0)^{-1}$ d'o\`u la forme normale
annonc\'ee.
\end{proof}

L'\'enonc\'e suivant est un r\'esultat de rigidit\'e; il dit que
deux transformations quadratiques g\'en\'eriques birationnellement
conjugu\'ees sont lin\'eairement conjugu\'ees.

\begin{thm}\label{conjbir}
{\sl Il existe une hypersurface $\mathcal{H}$ dans $\Sigma^3$
ayant la propri\'et\'e suivante: soient~$f$ et $g$ deux
\'el\'ements de $\Sigma^3\setminus\mathcal{H}$ sans courbe
invariante et ayant leurs sept points sp\'eciaux en position
g\'en\'erale; si $f$ et $g$ sont birationnellement conjugu\'es ils
le sont lin\'eairement.}
\end{thm}

\begin{proof}[\sl D\'emonstration]
D'apr\`es la Proposition \ref{tracedet} les traces et
d\'eterminants des matrices jacobiennes de $f$ et $g$ aux quatre
points fixes correspondants co\"incident. La Proposition
\ref{posgal} nous permet d'\'ecrire $f$ sous la forme normale

\begin{small}
\begin{align*}
&(x_0^2+A_0x_1x_2+B_0x_0x_2+C_0x_0x_1:x_1^2+A_1x_1x_2+B_1x_0x_2+C_1x_0x_1:\\
&x_2^2+\eta(B_0-A_0C_1)x_1x_2+ \eta(A_1-B_1C_0)x_0x_2+\eta (1-C_0C_1)x_0x_1)
\end{align*}
\end{small}

\noindent avec $\eta=(B_0A_1-B_1A_0)^{-1}.$ Les points $(1:0:0),$
$(0:1:0)$ et $(0:0:1)$ appartiennent \`a $\mathrm{Fix}\hspace{0.1cm}f;$ notons $m$ le
quatri\`eme point fixe de $f.$ Pour d\'emontrer le th\'eor\`eme il
suffit d'\'etablir que $f$ est uniquement d\'etermin\'e par la
connaissance des traces et d\'eterminants aux quatre points fixes.
Posons
\begin{align*}
 & t_1:=\mathrm{tr}(\mathrm{jac}\hspace{0.1cm} f_{(0:0:1)})=B_0+A_1,&& \delta_1:=
\mathrm{det}(\mathrm{jac}\hspace{0.1cm} f_{(0:0:1)})=1/\eta,\\
&t_2:=\mathrm{tr}(\mathrm{jac}\hspace{0.1cm} f_{(0:1:0)})=C_0+\eta(B_0-A_0C_1),&&
\delta_2:=\mathrm{det}(\mathrm{jac}\hspace{0.1cm} f_{(0:1:0)})=\eta(C_0B_0-A_0),\\
& t_3:=\mathrm{tr}(\mathrm{jac}\hspace{0.1cm} f_{(1:0:0)})=C_1+\eta(A_1-C_0
B_1),&& \delta_3:=\mathrm{det}(\mathrm{jac}\hspace{0.1cm} f_{(1:0:0)})=\eta(C_1A_1-B_1),\\
&t_4:=\mathrm{tr}(\mathrm{jac}\hspace{0.1cm} f_{(m)}),&&\delta_4:=\mathrm{det}(\mathrm{jac}\hspace{0.1cm} f_{(
m)}).
\end{align*}
Consid\'erons l'application $\Theta$ d\'efinie par
\begin{align*}
&\Theta\hspace{0.1cm}\colon\hspace{0.1cm}\mathbb{C}^6\to\mathbb{C}^6,&&(A_0,A_1,B_0,B_1,C_0,C_1)
\mapsto(t_1,\delta_1,t_2,\delta_2,t_3,\delta_3).
\end{align*}

Le calcul de
$\Theta^{-1}(\alpha,\beta,\gamma,\varepsilon,\mu, \zeta)$ se
ram\`ene \`a r\'esoudre le syst\`eme constitu\'e des \'equations suivantes
\begin{align*}
& B_0+A_1=\alpha, && B_0A_1-B_1A_0=\beta, && C_0\beta+B_0-A_0C_1=\gamma, \\
&  C_0B_0-A_0=\varepsilon, && C_0\beta+A_1-B_1C_0=\mu, && C_1A_1-B_1=\zeta.
\end{align*}
\`A l'aide de \texttt{Maple}, on obtient

\begin{footnotesize}
\begin{align*}
& A_0=\frac{-\alpha\beta^2+\beta(\alpha^2+\varepsilon+\beta)A_1
-\alpha(2\beta+\varepsilon)A_1^2+(\varepsilon+\beta)A_1^3}
{-\alpha\beta\zeta+(\beta\zeta-\beta+\alpha\mu)A_1-\mu A_1^2},
&& \\
&
B_0=-A_1+\alpha,
&& B_1= \frac{-\alpha\beta\zeta+(\beta\zeta-\beta+\alpha\mu)A_1-\mu A_1^2}{
\alpha\beta-(\beta+\varepsilon)A_1},\\
& C_0=\frac{-\beta(\zeta\varepsilon+\beta)+(\alpha\beta+\mu\varepsilon)A_1-(\varepsilon+\beta) A_1^2
}{-\alpha\beta\zeta+(\beta\zeta-\beta+\alpha\mu)
A_1-\mu A_1^2},
&&C_1=\frac{\mu\alpha-\beta-\varepsilon\zeta-\mu A_1}{\alpha\beta-(\beta+\varepsilon)A_1}.
\end{align*}
\end{footnotesize}

Quant \`a $A_1$ il satisfait l'\'equation quadratique

\begin{small}
\begin{align*}
& \beta(-\beta\zeta\varepsilon+\mu\alpha-\beta+\alpha\gamma\zeta-\beta^2-\zeta
\varepsilon)+(\alpha\beta^2+\beta\mu\varepsilon+2\beta\alpha\zeta-\beta\gamma
\zeta-\alpha\gamma\mu-\beta\mu-\alpha^2\beta\zeta+\beta\gamma+\alpha\zeta
\varepsilon)A_1\\
&\hspace{6mm}+(\gamma\mu-\beta\zeta-\varepsilon\zeta-\beta\varepsilon
-\beta^2)A_1^2=0.
\end{align*}
\end{small}

Ceci montre en particulier que pour $\chi=(\alpha,\beta,
\gamma,\varepsilon,\mu,\zeta)$ g\'en\'erique dans l'image de
$\Theta,$ la fibre~$\Theta^{-1}(\chi)$ contient exactement deux
points $\chi_1$, $\chi_2:$ c'est cette condition qui d\'etermine
l'hypersurface $\mathcal{H}=\{\chi\hspace{0.1cm}|\hspace{0.1cm} \#\hspace{0.1cm}\Theta^{-1}(\chi)\not=2\}.$
Pour conclure il suffit de trouver
une valeur explicite de~$\chi$ pour laquelle on aura $t_4(\chi_1)
\not=t_4(\chi_2).$

Un calcul conduit \`a
\begin{align*}
\Theta^{-1}(0,1,1,1,-1,2)=\left(
-\frac{\xi}{3},\xi,-\xi,-\frac{1}{2}-\frac{\xi}{2},
\frac{7\xi}{6}+\frac{1}{2},-\frac{7\xi}{4}-\frac{3}{4}\right)
\end{align*}
o\`u
\begin{align*}
\xi\in\left\{\xi_1=\frac{-1+\mathrm{i}\sqrt{167}}{14},
\xi_2=-\frac{1+\mathrm{i}\sqrt{167}}{14}\right\}.
\end{align*}
On constate
que $t_4(\chi_1)$ (correspondant au choix $\xi=~\xi_1$) et $t_4(\chi_2)$ sont distincts:
\begin{align*}
t_4(\chi_1)\sim 3.67+0.16.10^{-8}\mathrm{i}\not=t_4(\chi_2)\sim 3.67-0.16.
10^{-8}\mathrm{i}.
\end{align*}
\end{proof}

\begin{rems}
L'hypersurface $\mathcal{H}$ n'est pas explicite dans notre
\'enonc\'e; un examen plus pr\'ecis de l'application $\Theta$
permettrait sans doute d'affiner notre r\'esultat.

\bigskip

La d\'emonstration pr\'ec\'edente
montre deux choses
\begin{itemize}
\item si l'on se donne trois couples $(\xi_j,\mu_j)$
g\'en\'eriques, on peut les r\'ealiser comme valeurs propres en
$3$ points fixes $p_j$ d'une certaine transformation quadratique.
Par exemple les~$\xi_j,\hspace{1mm}\mu_j$ peuvent \^etre des complexes
de module $1$ g\'en\'eriques.

\item soit $g=A\sigma$ g\'en\'erique, $m_i$ ses points fixes,
$i=1,\ldots,4.$ On note $t(g,m_i)$ et $\delta(g,m_i)$ les traces
et d\'eterminants de la matrice jacobienne de $A\sigma$ en $m_i.$
Si un \'el\'ement $f$ de $\Sigma^3$ a quatre points fixes $p_i$ en
lesquels
\begin{align*}
&\delta(f,p_i)=\delta(g,m_i),&& t (f,p_i)=t(g,m_i)
\end{align*}
alors $f$ et $g$ sont dynamiquement conjugu\'es via un
automorphisme de $\mathbb{P}^2(\mathbb{C}).$ En effet comme $f$ est dans
$\Sigma^3$ il ne peut avoir trois points fixes align\'es donc les
$p_i$ peuvent \^etre envoy\'es sur les $m_i$ par un automorphisme
et on conclut ais\'ement.
\end{itemize}
\end{rems}

\section{Conjugaison birationnelle entre automorphismes de
\textsc{H\'enon}}\hspace{0.1cm}

Ce qui suit \'etend le Th\'eor\`eme \ref{conjbir} aux
transformations de \textsc{H\'enon}.

\begin{pro}
{\sl Soient $f$ et $g$ deux automorphismes de \textsc{H\'enon}
quadratiques. Si $f$ et $g$ sont birationnellement conjugu\'es,
ils le sont lin\'eairement.}
\end{pro}

\begin{proof}[\sl D\'emonstration]
\`A conjugaison lin\'eaire pr\`es $f$ est du type
$(1+x_1-ax_0^2,bx_0)$ et $g$ de la forme $(1+x_1-\alpha
x_0^2,\beta x_0).$ Notons $\varphi$ la conjuguante birationnelle.
Soit $m$ un point fixe de $f;$ on a $\varphi(m)=g \varphi(m).$
Puisqu'un automorphisme de type \textsc{H\'enon} ne fixe pas de
courbe, $m$ n'est pas d'ind\'etermination pour $\varphi.$ Par
cons\'equent $\varphi$ est holomorphe au voisinage de $m$ et
\begin{align*}
\varphi(\mathrm{Fix}\hspace{0.1cm} f)=\mathrm{Fix}\hspace{0.1cm} g.
\end{align*}
Le m\^eme argument
montre que $\varphi^{-1}$ est holomorphe en $\varphi(m).$ Ainsi
puisque $D\varphi^{-1}_{(\varphi (m))}D\varphi_{(m)}=~\mathrm{id},$
la transformation $\varphi$ est un diff\'eomorphisme local aux
points fixes $p_1$ et $p_2$ de $f.$ Si les $q_j=~\varphi(p_j)$ sont les
points fixes de $g$ on a
\begin{equation}\label{eg}
(\mathrm{det}(\mathrm{jac}\hspace{0.1cm} f_{(p_j)}),\mathrm{tr}(\mathrm{jac}\hspace{0.1cm} f_{(p_j)}))=(\mathrm{det}(
\mathrm{jac}\hspace{0.1cm} g_{(q_j)}),\mathrm{tr}(\mathrm{jac}\hspace{0.1cm} g_{(q_j)})).
\end{equation}
Or un calcul imm\'ediat montre que
\begin{small}
\begin{align*}
& \mathrm{jac}\hspace{0.1cm} f_{(p_1)}=\left[\begin{array}{cc}
1-b-\delta & 1 \\
b & 0
\end{array}
\right], && \mathrm{jac}\hspace{0.1cm} f_{(p_2)}=\left[\begin{array}{cc}
1-b+\delta & 1 \\
b & 0
\end{array}
\right], &&\delta=\sqrt{(1-b)^2+4a},\\
& \mathrm{jac}\hspace{0.1cm} g_{(q_1)}=\left[\begin{array}{cc}
1-\beta-\widetilde{\delta} & 1 \\
\beta & 0
\end{array}
\right],&& \mathrm{jac}\hspace{0.1cm} g_{(q_2)}=\left[\begin{array}{cc}
1-\beta+\widetilde{\delta} & 1 \\
\beta & 0
\end{array}
\right], &&\widetilde{\delta}= \sqrt{(1-\beta)^2+4\alpha}.
\end{align*}
\end{small}

\`A partir de (\ref{eg}) on obtient $b=\beta$ et
\begin{equation}\label{interrogation}
(1-\beta-\widetilde{
\delta})(1-\beta+\widetilde{\delta})=(1-b-\delta)(1-b+\delta)
\end{equation}
soit $(a,b)=(\alpha,\beta).$ Remarquons que $\delta$ et $\widetilde{
\delta}$ sont d\'efinis \`a $\pm 1$ pr\`es mais les produits
(\ref{interrogation}) sont eux bien d\'efinis.
\end{proof}

\clearemptydoublepage
\chapter[Propri\'et\'es dynamiques des transformations
birationnelles quadratiques]{Quelques propri\'et\'es dynamiques
des transformations birationnelles quadratiques}\label{dyndyn}

\section{Stabilit\'e alg\'ebrique}\hspace{0.1cm}

Rappelons une notion introduite dans \cite{FoSi} et
reprise dans \cite{DiFa}. La transformation $f$ est dite
\textbf{\textit{alg\'ebriquement stable}}\label{ind20} s'il
n'existe pas de courbe $\mathcal{C}$ dans $S$ telle
que $f^k(\mathcal{C})$ appartienne \`a $\mathrm{Ind}\hspace{0.1cm} f$ pour un
certain entier $k\geq 0.$ Dit autrement une transformation est
alg\'ebriquement stable si la situation suivante n'arrive pas

\begin{figure}[H]
\begin{center}
\input{as.pstex_t}
\end{center}
\end{figure}

\begin{thm}[\cite{DiFa}]
{\sl Soit $f$ un \'el\'ement du groupe de \textsc{Cremona}; $f$ n'est
pas alg\'ebriquement stable si et seulement s'il existe un entier
$k$ pour lequel $\deg f^k<(\deg f)^k.$}
\end{thm}

\begin{proof}[\sl D\'emonstration]
Nous ne d\'emontrerons que l'implication facile. Soit $k$ le plus
petit entier pour lequel $\deg f^k<(\deg f)^k;$ il existe alors
$h$ une transformation birationnelle de degr\'e $\deg h$ strictement
inf\'erieur \`a $(\deg
f)^{k+1}$ et $\varphi$ un polyn\^ome homog\`ene tels que en tant
que triplet de formes homog\`enes $f^{k+1}=\varphi h.$ On constate
que l'image par $f^k$ de tout point de la courbe d'\'equation
$\varphi=0$ est d'ind\'etermination pour $f.$

L'autre implication figure dans \cite{DiFa}.
\end{proof}

Il y a d'autres fa\c{c}ons de caract\'eriser la stabilit\'e alg\'ebrique
(\cite{DJS}).

Introduisons un invariant birationnel pour les
transformations de \textsc{Cremona}. Si $f$ et $g$ d\'esignent
deux \'el\'ements de $\mathrm{Bir}(\mathbb{P}^2(\mathbb{C}))$ alors en
g\'en\'eral $\deg(gfg^{-1})\not =\deg f;$ par contre il existe
deux constantes positives $\alpha$ et $\beta$ telles que pour tout
entier $n$ on ait (\cite{DiFa}, Proposition 1.15)
\begin{align*}
\alpha\hspace{0.1cm}\deg f^n\leq\deg(gf^ng^{-1})\leq\beta\hspace{0.1cm}\deg f^n.
\end{align*}
Le \og taux de croissance des degr\'es\fg\hspace{0.1cm} est donc un invariant
birationnel: on appelle \textbf{\textit{premier degr\'e
dynamique}}\label{ind21} de $f$ la quantit\'e
\begin{align*}
\lambda(f):=\liminf(\deg f^n)^{1/n}.
\end{align*}
On peut aussi d\'efinir cette
notion pour une transformation birationnelle d'une surface
complexe compacte k\"ahl\'erienne (\cite{DiFa}).

Lorsque $f$ est alg\'ebriquement stable, on a
$\lambda(f)= \deg f.$ En particulier pour $f=A\sigma$ avec $A$
g\'en\'erique on a $\lambda(A\sigma)=2;$ ceci est pr\'ecis\'e
ci-apr\`es.

Dans \cite{DiFa} \textsc{Diller} et \textsc{Favre}
montrent que toute transformation birationnelle d'une surface
complexe compacte est, \`a conjugaison birationnelle pr\`es,
alg\'ebriquement stable. Nous n'utiliserons pas vraiment ce
r\'esultat important parce que l'on souhaite travailler \`a
degr\'e fix\'e ce qui interdit la plupart des conjugaisons
birationnelles.

Soient
$A=(a_0x_0+b_0x_1+c_0x_2:a_1x_0+b_1x_1+c_1x_2:a_2x_0+b_2x_1+c_2x_2)$ un
automorphisme de $\mathbb{P}^2(
\mathbb{C})$ et $\mathcal{C}$ une courbe
contract\'ee par $(A\sigma)^{p+1}$ sur un point $q,$
 l'entier $p$ \'etant
minimal. La courbe $\widetilde{\mathcal{C}}=(A\sigma)^p(
\mathcal{C})$ est contract\'ee par $A\sigma;$ par suite
$\widetilde{\mathcal{C}}$ est l'une des droites de
$\mathrm{Exc}\hspace{0.1cm}\sigma,$ disons $x_2=0,$ et donc $q=(c_0:c_1:c_2).$ Supposons qu'il existe un
entier $k$ tel que $q$ soit \'eclat\'e par $(A\sigma)^{k+1}$ (avec
$k$ minimal) dit autrement $(A\sigma)^k(q)$ est un point
d'ind\'etermination de $\sigma,$ {\it i.e.} deux des composantes
de $(A\sigma)^k(q)$ sont nulles. \'Ecrivons $(A\sigma)^k$ sous la
forme $(P_0:P_1:P_2);$ visiblement les $P_i$ sont des polyn\^omes
en $(x_0,x_1,x_2)$ dont les coefficients sont eux-m\^emes des
polyn\^omes universels \`a coefficients entiers positifs en les variables
$a_0,$ $b_0,$ $c_0,$ $a_1,$ $b_1,$ $c_1,$ $a_2,$ $b_2,$ $c_2.$ Ceci permet de produire
des exemples explicites de transformations $A\sigma$ alg\'ebriquement
stables; en effet si $A$ est une matrice \`a coefficients r\'eels
tous strictement positifs les $P_i(c_0,c_1,c_2)$ sont strictement
positifs. On obtient donc la:

\begin{pro}\label{ass}
{\sl Soit $A$ un automorphisme de $\mathbb{P}^2(
\mathbb{C})$ dont les coefficients
sont des nombres r\'eels strictement positifs; la transformation
$A\sigma$ (resp. $A\rho$) est alg\'ebriquement stable.}
\end{pro}

\begin{rem}
C'est la \og positivit\'e des coefficients\fg\hspace{0.1cm} de $\sigma$ et
$\rho$ qui permet d'obtenir l'\'enon\-c\'e~\ref{ass}. Dans la
transformation $\tau$ on trouve des coefficients de signes
contraires ce qui ne permet pas d'obtenir de fa\c{c}on rapide un
\'enonc\'e de type \ref{ass}. Voici comment on peut l'obtenir. Si
l'\'el\'ement $h$ de $\Sigma^1$ est un automorphisme de
\textsc{H\'enon}, $\deg h^n=2^n=(\deg h)^n$ de sorte que $h$ est
alg\'ebriquement stable; c'est par exemple le cas si $h$ est de la
forme $(y,y^2+x).$ Alors pour $A$ lin\'eaire g\'en\'erique $Ah$
est alg\'ebriquement stable de sorte qu'un \'el\'ement
g\'en\'erique de $\Sigma^1$ est aussi alg\'ebriquement stable.
\end{rem}

\begin{defi}
Nous dirons qu'un automorphisme de $\mathbb{P}^2(
\mathbb{C})$ a ses
\textbf{\textit{coefficients alg\'ebriquement ind\'ependants sur
$\mathbb{Q}$}}\label{ind44} s'il poss\`ede un repr\'esentant $A$ dans
$\mathrm{GL}_3(\mathbb{C})$ dont
les coefficients sont alg\'ebriquement ind\'ependants sur $\mathbb{Q}.$
\end{defi}

\begin{rem}
La condition ci-dessus est invariante sous l'action des
automorphismes du corps $\mathbb{C}.$ En fait deux automorphismes de
$\mathbb{P}^2(\mathbb{C})$ \`a coefficients alg\'ebriquement ind\'ependants sont
\'echang\'es par un automorphisme du corps $\mathbb{C};$ en particulier
toutes leurs propri\'et\'es de nature \og\hspace{0.1cm} alg\'ebriques\fg\hspace{0.1cm} sont
conserv\'ees.
Par contre la condition d'ind\'ependance des coefficients ne
r\'esiste pas \`a la conjugaison dynamique. En effet si $A\sigma$
a ses coefficients rationnellement ind\'ependants sur $\mathbb{Q}$ elle
poss\`ede quatre points fixes distincts tous situ\'es en dehors de
$x_0x_1x_2=0.$ En conjuguant par une homoth\'etie ad-hoc
$(\alpha_0x_0:\alpha_1x_1:\alpha_2x_2)$ on peut supposer
par exemple que $(1:1:1)$ est point fixe: on obtient ainsi une nouvelle
transformation $\widetilde{A}\sigma$ conjugu\'ee \`a la pr\'ec\'edente
mais dont les coefficients ne sont plus rationnellement
ind\'ependants.
\end{rem}

Notons $\mathrm{Aut}(\mathbb{C},+,.)\label{not28}$ l'ensemble
des automorphismes du corps $\mathbb{C}.$ Si $A$ appartient \`a
$\mathrm{PGL}_3(\mathbb{C})$ et $\kappa$ \`a $\mathrm{Aut}(
\mathbb{C},+,.)$ on note
$A^\kappa\label{not29}$ l'\'el\'ement obtenu en faisant agir
$\kappa$ sur les coefficients de $A.$

\begin{cor}\label{coras}
{\sl Soit $A$ un automorphisme de $\mathbb{P}^2(
\mathbb{C})$ dont les coefficients
sont alg\'ebriquement ind\'ependants sur $\mathbb{Q};$ alors $A\sigma$ et
$(A\sigma)^{-1}$ sont alg\'ebriquement stables.}
\end{cor}

\begin{proof}[\sl D\'emonstration]
Soit $B$ un \'el\'ement de $\mathrm{GL}_3(\mathbb{C})$ dont tous
les coefficients sont r\'eels positifs et alg\'ebriquement
ind\'ependants sur $\mathbb{Q}.$ Les coefficients de $B$ font partie d'une
base de transcendance de $\mathbb{C}$ sur $\mathbb{Q}.$ Donnons-nous un
\'el\'ement $A$ de $\mathrm{GL}_3(\mathbb{C})$ dont les coefficients sont
alg\'ebriquement ind\'ependants sur $\mathbb{Q};$ il existe un
automorphisme $\kappa$ du corps $\mathbb{C}$ tel que $A=B^\kappa.$ Alors
$A\sigma=(B\sigma)^\kappa$ est alg\'ebriquement stable puisque
$B\sigma$ l'est.

Comme une transformation birationnelle $f$ d'une surface
$S$ dans elle-m\^eme est alg\'ebriquement stable si et seulement
si $f^{-1}$ l'est (\cite{Ja}) l'inverse de $A\sigma$ est
alg\'ebriquement stable.
\end{proof}

\begin{rem}
On a un \'enonc\'e analogue lorsqu'on remplace $A\sigma$ par
$A\rho$ ou $A\tau.$

Sous les m\^emes hypoth\`eses $(A\sigma)^n$ (resp.
$(A\rho)^n$, resp. $(A\tau)^n$) est alg\'ebriquement stable.
\end{rem}

\section{Feuilletages et courbes invariants pour
des transformations du type $A\sigma$}\hspace{0.1cm}

\subsection{Feuilletages invariants}\hspace{0.1cm}

Dans l'introduction nous avons mentionn\'e que tout
automorphisme de $\mathbb{P}^2(
\mathbb{C})$ poss\`ede des feuilletages invariants;
nous avons d'autre part rencontr\'e des transformations
birationnelles qui pr\'eservent une fibration rationnelle. Comme
l'indique l'exemple suivant il existe des transformations
birationnelles qui laissent invariant un feuilletage ne
d\'efinissant pas une fibration. Consid\'erons dans la carte
affine $x_2=1$ la transformation $f$ de degr\'e $3$ d\'efinie par
\begin{align*}
f(x_0,x_1)=(x_0^2x_1,x_0x_1).
\end{align*}
L'\'el\'ement $f$ pr\'eserve les
feuilletages associ\'es aux formes diff\'erentielles
\begin{align*}
&\omega_\eta=x_1\mathrm{d}x_0+\eta x_0\mathrm{d}x_1, && \eta= \frac{-1\pm \sqrt{5}}{2}
\end{align*}
dont les feuilles sont
transcendantes; par ailleurs un calcul facile assure que les seuls
feuilletages invariants par $f$ sont les deux pr\'ec\'edents.

La Proposition qui suit montre que la pr\'esence de
feuilletages invariants est rare.

\begin{pro}\label{pasfib}
{\sl Un \'el\'ement g\'en\'erique de $\mathrm{Bir}_2$ ne poss\`ede pas de
fibration invariante. Plus pr\'ecis\'ement un \'el\'ement
g\'en\'erique de $\mathrm{Bir}_2$ ne poss\`ede pas de feuilletage
alg\'ebrique invariant.}
\end{pro}

\begin{proof}[\sl D\'emonstration]
Les notations \'etant celles du Chapitre \ref{feuilfeuil}, dans le
produit $\mathbb{P}\mathcal{F}_n\times\mathrm{Bir}_2$ on consid\`ere le
sous-ensemble ferm\'e
\begin{align*}
\Lambda_n=\{([\omega],\phi)\in
\mathbb{P}\mathcal{F}_n\times\mathrm{Bir}_2\hspace{0.1cm}|\hspace{0.1cm}\phi^*\omega\wedge\omega
=0\}.
\end{align*}
Si $([\omega],\phi)$ appartient \`a $\Lambda_n,$ le feuilletage
associ\'e \`a $\omega$ est invariant par $\phi.$ Soit $\pi_2$ la
projection de~$\mathbb{P}\mathcal{F}_n\times \mathrm{Bir}_2$ sur
$\mathrm{Bir}_2.$ La projection de $\Lambda_n$ par $\pi_2$ est un
ferm\'e car les fibres de $\pi_2$ sont compactes. Par suite
$\pi_2(\Lambda_n)$ est ferm\'e alg\'ebrique dans $\mathrm{Bir}_2.$
Dans \cite{Br2} \textsc{Brunella} montre qu'un automorphisme
quadratique de \textsc{H\'enon} ne poss\`ede pas de feuilletage
invariant; le compl\'ement~$\mathcal{B}_n$ de~$\pi_2(\Lambda_n)$
est donc un ouvert de \textsc{Zariski} non vide dans
$\mathrm{Bir}_2.$  En r\'einvoquant qu'un automorphisme de
\textsc{H\'enon} n'a pas de feuilletage invariant et que
$\mathrm{Bir}_2$ est irr\'eductible, on obtient que l'intersection
$\cap\mathcal{B}_n$ est non vide, d'o\`u le r\'esultat.
\end{proof}

\begin{rem}\label{raisonnable}
Il semble raisonnable de penser que si les coefficients de la
matrice $A$ sont ind\'ependants sur $\mathbb{Q}$ alors $f=A\sigma$
ne poss\`ede pas de feuilletage invariant; mais l'\'elimination
directe de $\omega$ dans l'identit\'e $f^*\omega\wedge\omega=0$
s'av\`ere p\'enible. Le Corollaire \ref{mivi} qui \'etablit ce
fait est cons\'equence de la remarque qui suit.
\end{rem}

\begin{rem}
Dans la Proposition \ref{pasfib} le terme \og g\'en\'erique\fg\hspace{0.1cm}
signifie dans le compl\'ement d'une union d\'enombrable
d'ensembles alg\'ebriques propres de $\mathrm{Bir}_2,$ en fait de chaque
$\Sigma^i$ (pour~$i=1,\hspace{1mm}2,\hspace{1mm}3$). En particulier pour $A$ dans le
compl\'ement d'une union d\'enombrable d'ensembles alg\'ebriques
propres de $\mathrm{GL}_3(\mathbb{C})$ la transformation $A\sigma$ ne
poss\`ede pas de feuilletage invariant. Puisqu'on ne peut \'ecrire
$\mathrm{PGL}_3(\mathbb{C})$ comme une union d\'enombrable d'ensembles
alg\'ebriques propres il existe des matrices $A_0$ \`a
coefficients ind\'ependants sur $\mathbb{Q}$ telles que $A_0\sigma$ n'ait
pas de feuilletage invariant. Mais alors (toujours en utilisant un
automorphisme de corps ad hoc) pour toute matrice $A$
\`a coefficients ind\'ependants sur $\mathbb{Q}$ on a la m\^eme propri\'et\'e
d'o\`u l'\'enonc\'e sugg\'er\'e dans la Remarque
\ref{raisonnable}.
\end{rem}

\begin{cor}\label{mivi}
{\sl Soit $A$ un automorphisme de $\mathbb{P}^2(
\mathbb{C})$ dont les coefficients
sont $\mathbb{Q}$-alg\'ebri\-quement ind\'ependants; $A\sigma$ ne poss\`ede pas de feuilletage
invariant.}
\end{cor}

\noindent \textbf{\textit{Probl\`eme.}} Classifier les transformations birationnelles
quadratiques poss\'edant un feuilletage invariant.

\subsection{Courbes invariantes}\hspace{0.1cm}

Nous allons prouver qu'une transformation
quadratique g\'en\'erique de type $A\sigma$ ne poss\`ede pas de
courbe invariante. Cet \'enonc\'e fait partie du folklore, mais il
n'y a pas de r\'ef\'erence connue des auteurs. Nous proposons deux
d\'emarches. La premi\`ere ne concerne que les courbes invariantes
irr\'eductibles; elle donne une classification des transformations
ayant certains types de courbes invariantes. On pourrait penser
proc\'eder comme nous l'avons fait dans le cas des feuilletages,
par d\'eg\'en\'erescence sur un automorphisme de \textsc{H\'enon},
applications pour lesquelles cet \'enonc\'e est connu. Plus
pr\'ecis\'ement soit $(h_\varepsilon)_\varepsilon$ une famille
analytique de transformations birationnelles telle que
$h_{\varepsilon=0}=h_0$ soit un automorphisme de
\textsc{H\'enon}; si $\mathcal{C}_\varepsilon$ est la courbe
invariante par $h_\varepsilon$ la famille
$(\mathcal{C}_\varepsilon)_\varepsilon$ produit une courbe limite
$\mathcal{C}_0$ a priori invariante par $h_0$ mais rien ne dit que
$\mathcal{C}_0$ n'est pas la courbe exceptionnelle de $h_0.$

Pour la premi\`ere approche nous allons invoquer des
arguments de nature dynamique. La seconde, qui s'appuie sur
la premi\`ere, traite
des courbes non irr\'eductibles ({\it i.e.} des courbes
invariantes par les it\'er\'es) et invoque des arguments
d'ind\'ependance alg\'ebrique et de points fixes. Avant toute chose il faut rappeler
qu'une courbe contract\'ee par une transformation birationnelle
est une courbe rationnelle, autrement dit de genre arithm\'etique
$0.$ C'est une cons\'equence directe du th\'eor\`eme de
\textsc{N\oe ther}.

De nombreux auteurs se sont int\'eress\'es au probl\`eme
des courbes invariantes. Dans sa th\`ese
\textsc{Jackson} d\'ecrit les transformations birationnelles
quadratiques $f$ qui pr\'eservent une courbe $\mathcal{C}$ et qui
satisfont $\mathrm{Ind}\hspace{0.1cm} f\cup\mathrm{Ind}\hspace{0.1cm}
f^{-1}\subset\mathcal{C}$ (\emph{voir} \cite{Ja}).

Comme une surface de \textsc{Riemann} compacte de genre $g \geq
2$ ne poss\`ede qu'un nombre fini d'automorphismes, on devine que
pour une transformation birationnelle le fait de poss\'eder une
courbe invariante de genre $g\geq 2$ est tr\`es restrictif; ceci
est pr\'ecis\'e par la:

\begin{rem}
Soit $\mathcal{C}$ une courbe de genre arithm\'etique sup\'erieur
ou \'egal \`a $2.$ Si $f$ est une transformation birationnelle
laissant $\mathcal{C}$ invariante alors $f_{|\mathcal{C}}$ est un
automorphisme de $\mathcal{C},$ en particulier il existe un
it\'er\'e $f^n$ de $f$ tel que $\mathcal{C}$ soit une courbe de
points fixes de $f^n.$
\end{rem}

Commen\c{c}ons par mentionner un cas particulier d'un
r\'esultat de \textsc{Diller}, \textsc{Jackson} et
\textsc{Sommese}.

\begin{thm}[\cite{DJS}]
{\sl Soit $f$ une transformation de \textsc{Cremona} alg\'ebriquement
stable telle que $\lambda(f)>1.$ Si $\mathcal{C}$ est une courbe
invariante par $f,$ alors le genre (arithm\'etique) de
$\mathcal{C}$ est $0$ ou~$1.$}
\end{thm}

Ce r\'esultat est pr\'ecis\'e par \textsc{Pan} dans \cite{Pa3}.

\begin{thm}[\cite{Pa3}]
{\sl Soit $\mathcal{C}$ une courbe irr\'eductible de $\mathbb{P}^2(
\mathbb{C})$ de
genre (arithm\'etique) strictement plus grand que $1.$ Soit $f$
une transformation de \textsc{Cremona} qui pr\'eserve
$\mathcal{C}.$ On a
\begin{itemize}
\item $\lambda(f)=1$ et la suite $(\deg f^n)_n$ cro\^it au plus
lin\'eairement;

\item si $f$ est conjugu\'e \`a un automorphisme d'une surface
rationnelle (lisse), $f$ est d'ordre fini;

\item si la normalisation de $\mathcal{C}$ n'est pas
hyperelliptique, $f$ est d'ordre fini.
\end{itemize}}
\end{thm}

\bigskip

\begin{defi}
Soit $M\in\mathrm{GL}_2(\mathbb{C})$ une matrice $2\times 2.$ On dit que
les valeurs propres $\eta_1,$ $\eta_2$ de~$M$ sont
r\'esonnantes s'il existe des entiers $m_1$ et $m_2$ non tous deux
nuls tels que $\eta_1^{m_1}\eta_2^{m_2}=1.$ Un point fixe
$p$ d'une transformation birationnelle $f$ sera dit
\textbf{\textit{r\'esonnant}}\label{ind45} si $Df_{(p)}$ a
ses valeurs propres r\'esonnantes.
\end{defi}

\begin{pro}\label{res}
{\sl Soit $A$ un automorphisme de $\mathbb{P}^2(
\mathbb{C})$ dont les coefficients
sont $\mathbb{Q}$-alg\'ebrique\-ment ind\'ependants. Les points fixes de $A
\sigma$ sont non r\'esonnants.}
\end{pro}

\begin{proof}[\sl D\'emonstration]
On constate que
\begin{align*}
\{B\sigma\hspace{0.1cm}|\hspace{0.1cm} B\in\mathrm{PGL}_3(\mathbb{C})\}=
\overline{\{A^\kappa\sigma\hspace{0.1cm}|\hspace{0.1cm}\kappa\in\mathrm{Aut}(\mathbb{C},+,.)\}}
\end{align*}
o\`u l'adh\'erence est prise au sens ordinaire.
Remarquons que les points fixes de $A\sigma$ sont en position
g\'en\'erale; en effet si trois d'entre eux \'etaient align\'es,
alors cette propri\'et\'e, invariante par action d'un
automorphisme de corps sur $A,$ serait v\'erifi\'ee par tous les
$B\sigma$ ce qui est absurde. Soit $C$ un automorphisme de
$\mathbb{P}^2(\mathbb{C})$ tel que les points fixes de $g:=CA\sigma C^{-1}$ soient
\begin{align*}
&(0:0:1), && (0:1:0), && (1:0:0) && \text{et} &&(1:1:1).
\end{align*}
On note $\Omega$
l'ensemble des \'el\'ements de $\mathrm{Bir}_2$ fixant ces $4$ points;
l'orbite de $g$ sous l'action de $\mathrm{Aut}(\mathbb{C},+,.)$ est dense
dans $\Omega.$ Supposons que $A\sigma$ poss\`ede un point fixe
r\'esonnant; alors $g$ et~$g^\kappa,$ avec $\kappa$ dans
$\mathrm{Aut}(\mathbb{C},+,.)$, aussi. Il s'en suit que tout \'el\'ement de
$\Omega$ est r\'esonnant en ce point ce qui n'est pas correct.
\end{proof}

\begin{pro}\label{fix}
{\sl Soit $A$ un automorphisme de $\mathbb{P}^2(
\mathbb{C})$ dont les coefficients
sont alg\'ebriquement ind\'ependants sur $\mathbb{Q}.$ La transformation
$(A\sigma)^n$ n'a pas de courbe de points fixes pour chaque
$n$ dans $\mathbb{Z}\setminus\{0\}.$}
\end{pro}

\begin{proof}[\sl D\'emonstration]
Supposons que $\mathcal{C}$ soit une courbe de points fixes pour
$(A\sigma)^n;$ si $d$ est le degr\'e de $\mathcal{C},$ pour tout
$\kappa$ dans $\mathrm{Aut}(\mathbb{C},+,.),$ l'it\'er\'e $n$-i\`eme de
$A^\kappa\sigma$ a une courbe de points fixes de degr\'e $d.$ En
passant \`a l'adh\'erence de \textsc{Zariski} on obtient que pour
tout $f$ dans $\Sigma^3$ l'it\'er\'e $n$-i\`eme de $f$ admet une
courbe de points fixes de degr\'e $d;$ puisque $\Sigma^3$ est
dense dans $\mathrm{Bir}_2$ on en d\'eduit que toute
transformation birationnelle quadratique, et en particulier tout
automorphisme de \textsc{H\'enon}, v\'erifie une telle
propri\'et\'e ce qui est faux (\cite{Br2}).
\end{proof}

\begin{rem}
Tous ces \'enonc\'es se g\'en\'eralisent aux
transformations du type $(A\rho)^n$ et~$(A\tau)^n.$
\end{rem}

\begin{rem}
L'hypoth\`ese d'ind\'ependance sur les coefficients de $A$ est
indispensable: on remarque par exemple que si $A$ est du type
\begin{align*}
& (\alpha x_0-\alpha x_1+x_2:\beta x_0-\beta x_1+x_2:\delta x_0+(1-\delta)x_1),
&& \alpha\not=\beta
\end{align*}
l'\'el\'ement $A \sigma$ laisse la droite $x_0=x_1$ invariante.
\end{rem}

\begin{pro}\label{cc}
{\sl Soit $f=A\sigma$ une transformation birationnelle quadratique
laissant une courbe rationnelle $\mathcal{C}$ invariante.
Supposons que $f$ soit alg\'ebriquement stable et que les
it\'er\'es de $f$ n'aient ni courbe de points fixes, ni point
fixe r\'esonnant; alors $\mathcal{C}$ est une droite, une
conique lisse ou une cubique \`a point double.}
\end{pro}

\begin{proof}[\sl D\'emonstration]
Si $\mathcal{C}$ est lisse, la formule d'\textsc{Hurwitz} assure que
$\mathcal{C}$ est une droite ou une conique lisse.

Supposons d\'esormais que $\mathcal{C}$ soit non lisse
et que $p$ soit un point singulier de $\mathcal{C}.$

Montrons que $A\sigma$ est un diff\'eomorphisme local en
$p,$ {\it i.e.} que $p$ n'est pas dans $\mathrm{Exc}\hspace{0.1cm}
\sigma\cup\mathrm{Ind}\hspace{0.1cm}\sigma.$ D'apr\`es \cite{DJS} le point
$p$ n'est pas d'ind\'etermination pour $\sigma.$ Si $p$
appartient \`a $\mathrm{Exc}\hspace{0.1cm}\sigma\setminus\mathrm{Ind}\hspace{0.1cm}\sigma$
alors~$A\sigma(p)$ est un point d'ind\'etermination de $(A
\sigma)^{-1}.$ Le germe de $A\sigma_{,p}$ est du type 
\begin{align*}
& (u,v)\mapsto(u,uv)
\end{align*}
\`a conjugaison locale g.d. pr\`es. En particulier
il envoie courbe singuli\`ere irr\'eductible sur courbe
singuli\`ere; $A\sigma(p)$ est donc un point d'ind\'etermination
de $(A\sigma)^{-1}$ mais aussi un point singulier ce qui est exclu
(toujours d'apr\`es \cite{DJS}). Ainsi $A\sigma$ est un
diff\'eomorphisme local en $p$ et $A\sigma(\mathcal{C}_{,p})$ est
singuli\`ere.

Si $\pi\hspace{0.1cm}\colon\hspace{0.1cm}\widetilde{\mathcal{C}}\to\mathcal{C}$
est la normalisation de
$\mathcal{C}$ et $\widetilde{A\sigma}$ le relev\'e \`a
$\widetilde{\mathcal{C}}$ de la restriction $A\sigma_{|\mathcal{C}},$
on note
\begin{align*}
\Lambda=\pi^{-1}(\mathrm{Sing}\hspace{0.1cm}(\mathcal{C}))=\{p_1,\ldots, p_s\}.
\end{align*}
Comme $A\sigma$ est un diff\'eomorphisme local en chaque point
singulier de $\mathcal{C},$ l'ensemble $\Lambda$ est invariant par
$\widetilde{A\sigma}.$ Puisque $\widetilde{\mathcal{C}}$ est
isomorphe \`a $\mathbb{P}^1(
\mathbb{C}),$ $\widetilde{A\sigma}$ est conjugu\'e
\`a une translation ou une homoth\'etie; on en d\'eduit que
$\Lambda$ compte au
maximum deux points. En effet si $\# \Lambda\geq 3,$ la
transformation $A\sigma$ est p\'eriodique en restriction \`a
$\mathcal{C};$ il existe alors un entier $n$ tel que
$\mathcal{C}\subset \mathrm{Fix}(A\sigma)^n$ ce qui est exclu
d'apr\`es l'hypoth\`ese. Notons qu'\`a chaque branche locale de
$\mathcal{C}_{,p}$ correspond un~$p_i;$ par suite on a
l'alternative suivante

\begin{itemize}
\item $\mathcal{C}$ a un seul point singulier $p$ (qui
est donc fixe) et $\mathcal{C}$ a au plus deux branches en~$p;$

\item $\mathcal{C}$ a deux points singuliers
distincts $q_1,$ $q_2$ et $\mathcal{C}_{,q_1},$
$\mathcal{C}_{,q_2}$ sont irr\'eductibles.
\end{itemize}

\'Etudions ces diff\'erentes possibilit\'es.

Remarquons que si un germe de diff\'eomorphismes
holomorphe de $\mathbb{C}^2,0$ fixe une courbe singuli\`ere irr\'eductible
ou bien deux courbes lisses tangentes ce point fixe est
r\'esonnant.

\begin{itemize}
\item Puisque les valeurs propres de
la partie lin\'eaire de $A\sigma$ en $p$ ne peuvent \^etre
r\'esonnantes, $\mathcal{C}$ a un unique point
singulier et $\mathcal{ C}_{,p}$ est un croisement ordinaire;
la formule du genre assure que $\mathcal{C}$ est une cubique \`a
point double.

\item On remarque que $(A\sigma)^2$ fixe
les points $q_1$ et $q_2;$ comme $\mathcal{C}_{,q_1}$ est
irr\'eductible non lisse, les valeurs propres de la partie
lin\'eaire de $(A\sigma)^2$ en $q_1$ (resp. $q_2$) sont alors
r\'esonnantes ce qui est impossible.
\end{itemize}
\end{proof}

Comme toujours il y a un \'enonc\'e analogue pour les
transformations g\'en\'eriques de $\Sigma^2$ et $\Sigma^1.$

\bigskip

Comme on le voit les cubiques cuspidales n'apparaissent pas dans
la Proposition \ref{cc}. Certaines transformations quadratiques
poss\`edent de telles cubiques invariantes;
dans \cite{BK} \textsc{Bedford} et \textsc{Kim} \'etudient la
famille de transformations de \textsc{Cremona} d\'efinie par
\begin{align*}
g_{a,b}\hspace{0.1cm}\colon\hspace{0.1cm}\mathbb{P}^2(\mathbb{C})
\dashrightarrow\mathbb{P}^2(\mathbb{C}), && (x_0:x_1:x_2)
\mapsto(x_0(x_1+bx_0):x_2(x_1+bx_0):x_0(x_2+ax_0)),
\end{align*}
famille dont nous reparlerons au Chapitre \ref{expexp}. Notons que les $g_{a,b}$
ne sont pas alg\'ebriquement stables: la droite $x_0=0$ est
contract\'ee par $g_{a,b}$ sur $(0:1:0),$ point qui est \'eclat\'e par $g_{a,b}$
sur la droite $x_2=0.$

Cette famille poss\`ede, pour certains param\`etres, une cubique
cuspidale invariante. Posons
\begin{align*}
\varphi(t)=\left(\frac{t-t^3-t^4}{(1+t)^2},\frac{1-t^5}{t^2(1+t)}\right).
\end{align*}
\textsc{Bedford} et \textsc{Kim} \'etablissent en particulier le:
\begin{thm}[\cite{BK}]
{\sl Soit $t$ dans $\mathbb{C}\setminus\{0,\hspace{0.1cm} 1,\hspace{0.1cm} -1,\hspace{0.1cm}\mathrm{j},\hspace{0.1cm}\mathrm{j}^2\}.$
Il existe un polyn\^ome $P$ dans $\mathbb{C}[x_0,x_1,x_2]_3$ satisfaisant
\begin{align*}
P\circ g_{a,b}= t \det(\mathrm{jac}\hspace{0.1cm}
g_{a,b})\hspace{0.1cm} P
\end{align*}
si et seulement si $(a,b)=\varphi(t).$

Le polyn\^ome $P$ est donn\'e par
\begin{align*}
& P_{t,a,b}(x_0,x_1,x_2)=a(t-1)t^4x_0^3+t(t-1)x_1x_2(tx_1+x_2)\\
& \hspace{28mm} +x_0\left(2bt^3x_1x_2+t^3(t-1)x_1^2+(t-1)(1+bt)x_2^2\right)\\
& \hspace{28mm} + t^3(t-1)x_0 ^2\left(a(x_1+tx_2)+t(x_1+(t-2b)x_2)\right)
\end{align*}}
\end{thm}

Supposons que $(a,b)$ co\"incide avec $\varphi(t)$ pour un
certain $t$ dans $\mathbb{C}\setminus\{0,\hspace{0.1cm} 1,\hspace{0.1cm} -1,\hspace{0.1cm} \mathrm{j},\hspace{0.1cm}\mathrm{j}^2\}.$
La cubique irr\'eductible d'\'equation $P_{t,a,b}=0$ invariante par $g_{a,b}$ contient les deux points fixes de~$g_{a,b}$
\begin{align*}
&p_1=\left(\frac{t^3}{1+t},\frac{t^3}{1+t}\right) && \text{et} && p_2=\left(\frac{-1+t^2+t^3}{t^2(1+t)},\frac{-1+t^2+t^3}{t^2(1+t)}\right)
\end{align*}
et a son cusp en $p_1.$

\begin{pro}
{\sl Soit $f=A\sigma$ une transformation birationnelle quadratique
laissant invariante une courbe irr\'eductible $\mathcal{C}$ de
genre arithm\'etique $1.$ Supposons que $f$ soit alg\'ebriquement
stable et que les it\'er\'es de $f$ n'aient pas courbe de points
fixes; alors $\mathcal{C}$ est une cubique lisse.}
\end{pro}

\begin{proof}[\sl D\'emonstration]
En utilisant des arguments de la d\'emonstration de la
Proposition \ref{cc} on montre que $\mathcal{C}$ est lisse. La restriction de
$f$ \`a $\mathcal{C}$ induit un automorphisme de la normalis\'ee
qui laisse invariant le relev\'e de l'ensemble singulier de
$\mathcal{C}$ (puisque $A\sigma$ est un diff\'eomorphisme local
en chaque point singulier de $\mathcal{C}$). Mais sur une courbe
elliptique un automorphisme qui
poss\`ede une orbite finie est d'ordre fini, ce qui est exclu par
hypoth\`ese; ainsi $\mathcal{C}$ est une courbe elliptique lisse
de $\mathbb{P}^2(\mathbb{C})$ donc une cubique lisse.
\end{proof}

Nous allons dans un premier temps nous int\'eresser aux
transformations birationnelles quadratiques pr\'eservant une
cubique lisse.

\begin{lem}\label{cubcub}
{\sl Soit $\mathcal{C}$ une cubique lisse. On suppose que
$\mathcal{C}'=\sigma(\mathcal{C})$ est aussi une cubique lisse;
alors $\mathcal{C}$ passe par les trois points d'ind\'etermination
de $\sigma.$}
\end{lem}

\begin{proof}[\sl D\'emonstration]
Soient $h$ et $h'$ dans $\mathbb{C}[x_0,x_1,x_2]_3$ d\'efinissant respectivement
$\mathcal{C}$ et $\mathcal{C}';$ l'hypoth\`ese implique que
\begin{align*}
&h\circ\sigma=x_0^\alpha x_1^\beta x_2^\gamma h', &&\alpha, \hspace{1mm}
\beta, \hspace{1mm}\gamma\in\mathbb{N},&&\alpha+\beta+\gamma=3.
\end{align*}
\`A permutation des coordonn\'ees pr\`es
le triplet $(\alpha,\beta,\gamma)$ prend les valeurs $(3,0,0),$
$(2,1,0)$ ou~$(1,1,1).$ Dans le premier cas un calcul direct
montre que $\mathcal{C}$ est l'union de trois droites et dans le
second que $\mathcal{C}$ est singuli\`ere. On v\'erifie, dans la
derni\`ere \'eventualit\'e, que $h$ est du type
\begin{align*}
ax_0^2x_1+bx_0^2x_2+cx_0x_1^2+dx_0x_2^2+ex_1^2x_2+fx_1x_2^2+
gx_0x_1x_2.
\end{align*}
\end{proof}

Notons que si $\mathcal{C}$ et $\mathcal{C}'$ sont comme dans
l'\'enonc\'e $\sigma$ r\'ealise un automorphisme entre
$\mathcal{C}$ et $\mathcal{C}'$ de sorte que ces deux cubiques ont
m\^eme invariant $j;$ on en d\'eduit
l'existence d'un automorphisme de $\mathbb{P}^2(
\mathbb{C})$ \'echangeant
$\mathcal{C}$ et $\mathcal{C}'$ (voir \cite{Ha}).

Soit $\mathcal{T}$ l'espace des polyn\^omes homog\`enes
de degr\'e $3$ s'annulant aux trois points d'ind\'etermination de
$\sigma.$ L'application $\sigma$ induit un isomorphisme
lin\'eaire $\widetilde{\sigma}\hspace{0.1cm}\colon\hspace{0.1cm} \mathcal{T}\to\mathcal{T}$
d\'efini dans la base~$\{x_0^2x_1,\hspace{0.1cm} x_0^2 x_2,\hspace{0.1cm} x_0x_1^2,\hspace{0.1cm} x_0x_2^2,\hspace{0.1cm}
x_1^2x_2,\hspace{0.1cm} x_1x_2^2,\hspace{0.1cm}x_0x_1x_2\}$  par
\begin{align*}
\widetilde{\sigma}\hspace{0.1cm}\colon\hspace{0.1cm}(a,b,c,d,e,f,g)\mapsto(f,e,d,c,b,a,g).
\end{align*}
Les cubiques invariantes par $\sigma$ correspondent \`a l'espace
lin\'eaire $\mathrm{Fix}\hspace{0.1cm}\widetilde{\sigma}$ qui est de dimension projective $3.$

Soient $\mathcal{C}$ une cubique lisse associ\'ee \`a un
\'el\'ement de $\mathcal{T}$ et
$\mathcal{C}'=\sigma(\mathcal{C}).$ Soit $A$ un isomorphisme
lin\'eaire tel que $A(\mathcal{C}')=\mathcal{C};$ la cubique
$\mathcal{C}$ est alors invariante par $A\sigma.$ Notons qu'il
existe un nombre fini de tels isomorphismes.
L'ensemble des \'el\'ements $A$ de $\mathrm{PGL}_3(\mathbb{C})$ tels que $A
\sigma$ laisse une cubique lisse invariante est d'adh\'erence une
vari\'et\'e alg\'ebrique d\'efinie sur $\mathbb{Q}$ de dimension
$\dim\mathbb{P}\mathcal{T}=6.$

Plus g\'en\'eralement soit $\mathcal{C}$ une cubique
lisse; donnons-nous trois points $p_1,$ $p_2$ et $p_3$ non
align\'es sur $\mathcal{C}.$ \`A ce triplet de points
$p=(p_1,p_2,p_3)$ on peut associer $\sigma_p$ une transformation
quadratique de \textsc{Cremona} telle que
$\mathrm{Ind}\hspace{0.1cm}\sigma_p=\{p_1,\hspace{0.1cm} p_2,
\hspace{0.1cm} p_3\}.$ L'image de $\mathcal{C}$ par $\sigma_p$ est
encore une cubique isomorphe \`a $\mathcal{C}$ et il existe une
application lin\'eaire $A_p$ (en fait il en existe un nombre fini)
telle que $A_p \sigma_p(\mathcal{C})=\mathcal{C}.$ Nous noterons
$\varphi_p$ les transformations $A_p\sigma_p.$ L'espace des
cubiques lisses est de dimension $9,$ celui des triplets de points
sur une cubique de dimension~$3;$ ainsi l'espace des \'el\'ements
de $\Sigma^3$ qui laissent invariante une cubique lisse est de
dimen\-sion~$9+~3=~12.$

En particulier on obtient le:

\begin{cor}
{\sl G\'en\'eriquement une transformation de $\Sigma^3$ ne
poss\`ede pas de courbe invariante de genre arithm\'etique $1.$

Soit $A$ un automorphisme de $\mathbb{P}^2(\mathbb{C});$ si
les coefficients de $A$ sont ind\'ependants sur $\mathbb{Q}$ alors~$A\sigma$ ne
pr\'eserve pas de telle courbe.}
\end{cor}

\begin{proof}[\sl D\'emonstration]
Le premier point r\'esulte des consid\'erations de dimension
qui pr\'ec\`edent.

Pour le second, si $A\sigma$ poss\'edait une courbe
elliptique invariante il en serait de m\^eme pour~$A^\kappa
\sigma$ et ce pour chaque $\kappa$ dans $\mathrm{Aut}(\mathbb{C},+,.);$
comme les orbites dynamiques des $A^\kappa\sigma$
sont denses dans~$\Sigma^3$ ceci contredit les calculs de
dimensions ci-dessus.
\end{proof}

Dans \cite{Pa} \textsc{Pan} a d\'emontr\'e que le groupe
$\mathrm{Bir}(\mathcal{C})$ des transformations birationnelles qui
laissent invariante la cubique lisse $\mathcal{C}$ est engendr\'e
par les transformations $\varphi_p=A_p\sigma_p.$

Montrons que tout automorphisme abstrait de
$\mathcal{C}$ est la restriction ou bien d'un $\varphi_p,$ ou bien
d'un isomorphisme lin\'eaire. On note $\mathrm{Aut}(\mathcal{C})$
le groupe des automorphismes de $\mathcal{C}$ vue comme courbe
elliptique abstraite $\mathbb{C}/\Lambda.$ En suivant une id\'ee de
\textsc{Pan} (\cite{Pa2}) nous allons montrer que le morphisme de
restriction
\begin{align*}
\mathrm{Bir}(\mathcal{C})\to\mathrm{Aut}(\mathcal{C})
\end{align*}
est surjectif. Pour
cel\`a fixons deux points g\'en\'eriques $m_1,$ $m_2$ de
$\mathcal{C}$ et $T$ un \'el\'ement de $\mathrm{Aut}(\mathcal{C}).$ Soit
$M$ un point g\'en\'erique du plan projectif complexe; on note $m_1,$
$M'_1$ et $M''_1$ (respective\-ment~$m_2,$ $M'_2$ et $M''_2$) les points
d'intersection de la droite $(Mm_1)$ (resp. $(Mm_2)$) avec la cubique~$\mathcal{C}$

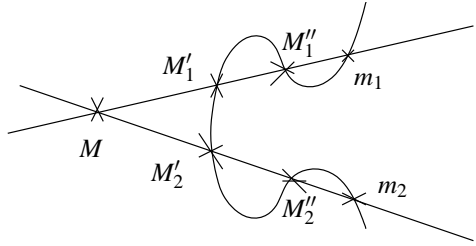
\begin{figure}[H]
\begin{center}
\input{courbell.pstex_t}
\end{center}
\end{figure}

Les couples de points $T(M'_i)$ et $T(M''_i)$
d\'eterminent deux droites dont l'intersection est not\'ee~$\widetilde{T}(M).$ L'application $\widetilde{T}$ co\"{\i}ncide avec $T$
sur $\mathcal{C}$ et est visiblement birationnelle puisque
$\widetilde{T}^{-1}$ inverse $\widetilde{T},$ d'o\`u la surjectivit\'e de
l'application de restriction: $\mathrm{Bir}(\mathcal{C})\to\mathrm{Aut}(\mathcal{C}).$

Revenons aux transformations quadratiques $\varphi_p.$ En
fait chaque triplet $p$ produit un nombre fini de telles applications,
autant qu'il y a d'automorphismes de $\mathbb{P}^2(\mathbb{C})$ qui fixent $\mathcal{C}.$
L'application $R\hspace{0.1cm}\colon\hspace{0.1cm} p\mapsto \varphi_{p_{|\mathcal{C}}}$
peut donc \^etre vue comme une application alg\'ebrique multivalu\'ee de
l'espace des triplets de points distincts de $\mathcal{C}$ dans
$\mathrm{Aut}(\mathcal{C}).$ En g\'en\'eral $\mathrm{Aut}(\mathcal{C})$
s'identifie topologiquement \`a l'union de deux exemplaires
$\mathcal{C}_0$ et $\mathcal{C}_-$ de $\mathcal{C}$ (sauf pour les
cubiques sp\'eciales~$\mathbb{C}/\mathbb{Z}[\mathrm{i}]$ et $\mathbb{C}/\mathbb{Z}[\mathrm{j}])$
\begin{align*}
&\mathcal{C}_0\simeq \{z+a\hspace{0.1cm}|\hspace{0.1cm} a\in\mathbb{C}/\Lambda\}&&\text{ et }
&&\mathcal{C}_-\simeq\{-z+a\hspace{0.1cm}|\hspace{0.1cm} a\in\mathbb{C}/\Lambda\}.
\end{align*}
Si on \'ecrit $\mathcal{C}$ sous forme normale de \textsc{Weierstrass} $x_1^2=x_0
(x_0-1)(x_0-\eta)$ on constate que l'automorphisme $\ell\hspace{0.1cm}\colon\hspace{0.1cm}
(x_0,x_1)\mapsto(x_0,-x_1)$ laisse $\mathcal{C}$ invariante et sa
restriction $\ell_{|\mathcal{C}}$ est dans la composante
$\mathcal{C}_-;$ l'image de $R$ intersecte donc $\mathcal{C}_0$ et
$\mathcal{C}_-.$

Remarquons que $R$ est non constante en tant
qu'application multivalu\'ee; si elle l'\'etait
$\{\varphi_{p_{|\mathcal{C}}}\}$ engendrerait un groupe
d\'enombrable dans $\mathrm{Aut}(\mathcal{C})$ ce qui contredirait la
surjectivit\'e de 
\begin{align*}
\mathrm{Bir}(\mathcal{C})\to\mathrm{Aut}(\mathcal{C})
\end{align*}
puisque
$\mathrm{Bir}(\mathcal{C})$ est engendr\'e par les $\varphi_p.$

Il en r\'esulte, dans le cas des cubiques g\'en\'eriques, que
l'image de $R$ \'evite au plus un nombre fini d'\'el\'ements de
$\mathrm{Aut}(\mathcal{C})$ (en particulier l'identit\'e). Ainsi tout
automorphisme de $\mathcal{C},$ sauf peut-\^etre un nombre fini,
est restriction d'un certain $\varphi_p.$

\begin{rem}
Soient $p=(p_1,p_2,p_3)$ un triplet de points distincts et $a$
un \'el\'ement de la cubique $\mathcal{C}$ que l'on consid\`ere
avec sa structure de groupe. Alors, pour chaque $a,$ l'application
$\varphi_{p+a},$ avec $p+a=(p_1+a,p_2+a,p_3+a),$ est bien d\'efinie.
Comme $\mathcal{C}$ est compacte l'application $a\mapsto\varphi_{{p
+a}_{|\mathcal{C}}}$ est constante (sinon $\varphi_p$ serait surjective).
\end{rem}

Nous passons maintenant \`a l'\'etude des transformations
birationnelles de type $A\sigma$ laissant une courbe de genre nul
invariante.

\begin{pro}\label{cat}
{\sl Soit $A$ un automorphisme de $\mathbb{P}^2(
\mathbb{C})$ d\'efini par
\begin{align*}
(a_0x_0+b_0x_1+c_0x_2:a_1x_0+b_1x_1+c_1x_2:a_2x_0+b_2x_1+c_2x_2).
\end{align*}
La transformation $A\sigma$ laisse une droite invariante
$\mathcal{D}$ si, \`a conjugaison dynamique par une permutation de
coordonn\'ees pr\`es, on a
\begin{align*}
a_0c_1^2+c_0c_1(b_0-a_1)-b_1c_0^2=0.
\end{align*}}
\end{pro}

\begin{proof}[\sl D\'emonstration]
Remarquons que $\mathcal{D}$ passe n\'ecessairement par un des
points de $\mathrm{Ind}\hspace{0.1cm}\sigma,$ par exemple $(0:0:1).$
Pla\c{c}ons-nous dans la carte $x_2=1$ et supposons que
$\mathcal{D}$ soit param\'etr\'ee par~$x_0\mapsto(x_0,tx_0);$ en \'ecrivant
l'invariance de $\mathcal{D}$ on obtient
\begin{align*}
x_0(a_0t^2+t(b_0-a_1)-b_1)+tx_0^2(tc_0-c_1)=0.
\end{align*}
Comme $t$ ne peut
\^etre nul ($x_1=0$ est une droite contract\'ee par $A\sigma$) on
constate que $c_1=tc_0$ d'o\`u l'alternative
\begin{itemize}
\item $c_0=c_1=0$ et $a_0t^2+(b_0-a_1)t-b_1=0;$

\item $t=c_1/c_0;$ en substituant cette valeur de $t$ dans
$(a_0t^2+t(b_0 -a_1)-b_1)$ on obtient
\begin{align*}
a_0c_1^2+c_0c_1(b_0-a_1) -b_1c_0^2=0.
\end{align*}
\end{itemize} Dans les deux
cas la condition annonc\'ee est satisfaite.
\end{proof}

\begin{rem}
Le cas $c_0=c_1=0,$ qui appara\^it en cours de preuve, correspond
aux transformations de type $A\sigma$ qui laissent la fibration
$x_1/x_0$ invariante.

Remarquons aussi que si une transformation $A\sigma$
laisse invariante une fibration en droites, alors, \`a permutation
pr\`es des coordonn\'ees, il s'agit de la fibration $x_1/x_0=$ cte.
Les $A$ convenables sont du type
\begin{align*}
(a_0x_0+b_0x_1:a_1x_0+b_1x_1:a_2x_0+b_2
x_1+c_2x_2);
\end{align*}
ces transformations forment un groupe.
\end{rem}

\begin{rem}
Les \'el\'ements de $\mathrm{Bir}_2$ poss\'edant une droite invariante
forment un ensemble alg\'ebrique dont la trace sur $\Sigma^3$
est une hypersurface.
\end{rem}

Passons maintenant \`a l'\'etude des transformations pr\'eservant
une conique. Soit $\mathcal{C}$ une conique lisse; $\sigma(
\mathcal{C})$ est une conique lisse si et seulement si, \`a
permutation pr\`es des coordonn\'ees,~$\mathcal{C}$ est du type
suivant
\begin{align*}
&\mathcal{C}_{ab}\hspace{0.1cm}\colon\hspace{2mm} x_2^2+ax_1x_2+bx_0x_2
+x_0x_1=0, && a,\hspace{1mm}b\in\mathbb{C},\hspace{1mm}2ab-1
\not=0.
\end{align*}
L'image par $\sigma$ de $\mathcal{C}_{ab}$ est
la conique $\mathcal{C}_{ba}$ de sorte que la permutation
\begin{align*}
P\hspace{0.1cm}\colon\hspace{0.1cm}(x_0:x_1:x_2)\mapsto(x_1:x_0:x_2)
\end{align*}
\'echange ces deux coniques. Soit $f=A\sigma$ une transformation
laissant invariante une conique $\mathcal{C}.$ \`A permutation
pr\`es des coordonn\'ees $\mathcal{C}$ est de type
$\mathcal{C}_{ab};$ de $\mathcal{C}_{ab}=A\sigma
(\mathcal{C}_{ab})=A(\mathcal{C}_{ba})=AP(\mathcal{C}_{ab})$ on d\'eduit
la:

\begin{pro}\label{coni}
{\sl Soient $f=A\sigma$ laissant une conique lisse invariante et $P$
l'involution d\'efinie par $(x_1:x_0:x_2).$ \`A permutation pr\`es
des coordonn\'ees $\mathcal{C}$ est donn\'ee par
\begin{align*}
x_2^2+ax_1x_2+ bx_0x_2+x_0x_1=0
\end{align*}
et $AP$ appartient au groupe orthogonal ${\bf
O}(x_2^2+ax_1x_2+bx_0x_2+x_0x_1).$}
\end{pro}

\begin{rem}
Lorsque $a=b$ l'involution $\sigma$ laisse $\mathcal{C}_{aa}$
invariante et dans ce cas $A$ est dans le groupe orthogonal ${\bf
O} (\mathcal{C}_{aa}).$
\end{rem}

\bigskip

Pour terminer le cas des courbes rationnelles apparaissant dans la
Proposition \ref{cc}, \'etudions les cubiques \`a point double.
L'image par $\sigma$ d'une cubique \`a point double $\mathscr{C}$
d'\'equation
\begin{equation}\label{cub1}
 ax_0^2x_1+bx_0^2x_2+cx_0x_1^2+dx_0x_2^2+ex_1^2x_2+gx_1x_2^2+
 hx_0x_1x_2=0
\end{equation}
est encore une cubique \`a point double (il s'agit ici de cubiques passant
par les trois points d'ind\'etermination de $\sigma$). Mais il y a d'autres
formes normales ayant la m\^eme
propri\'et\'e (correspondant \`a d'autres choix que $(1,1,1)$ pour
le triplet $(\alpha,\beta,\gamma)$ apparaissant dans la
d\'emonstration du Lemme \ref{cubcub}). Ainsi la cubique
d'\'equation
\begin{equation}\label{cub2}
ax_2^3+bx_0x_1^2+cx_0x_2^2+dx_1^2x_2+ex_1x_2^2+gx_0x_1x_2=0
\end{equation}
est envoy\'ee par $\sigma$ sur
\begin{align*}
ax_0x_1^2+bx_2^3+cx_1^2x_2+dx_0x_2^2+ex_0x_1x_2+gx_1x_2^2=0
\end{align*}
(induisant d'ailleurs l'application $(a,b,c,d,e,g)\mapsto
(b,a,d,c,g,e)$). Pour des valeurs g\'en\'eriques de $a,$ $b,$ $c,$
$d,$ $e,$ $g$ la cubique correspondante poss\`ede un point double
en $(1:0:0)$ qui est pr\'ecis\'ement d'ind\'etermination pour $\sigma.$

On v\'erifie facilement que les cubiques donn\'ees par
(\ref{cub1}) et (\ref{cub2}) sont, \`a permutation de
coordonn\'ees pr\`es, les seules cubiques \`a point double ayant
leur image par $\sigma$ encore cubique.

Soit $f=A\sigma$ une transformation birationnelle
quadratique laissant invariante une cubique \`a point double
$\mathscr{C};$ alors $\sigma(\mathscr{C})$ est encore une cubique
\`a point double. Comme \`a conjugaison lin\'eaire pr\`es toutes
les cubiques \`a point double sont les m\^emes, il existe $A_0$
dans $\mathrm{PGL}_3(\mathbb{C})$ tel que $\sigma(\mathscr{C})=A_0( \mathscr{C})$
de sorte que $\mathscr{C}=AA_0(\mathscr{C}).$ Le groupe des
transformations lin\'eaires qui laissent invariante une cubique
\`a point double est fini. Par exemple si $\mathscr{C}_0$ est la
cubique de \textsc{Weierstrass}: $x_1^2x_2=x_0^2(x_0-x_2);$ on v\'erifie que
ce groupe est form\'e des transformations
\begin{align*}
& \mathrm{id},&& \left(4\alpha x_0+4x_1:4x_0-4\alpha x_1:9\alpha
x_0+3x_1-8\alpha x_2\right),\\
& (x_0:-x_1:x_2),&& \left(4\alpha x_0-4x_1:4x_0+4\alpha x_1:9\alpha
x_0-3x_1-8\alpha x_2\right)
\end{align*}
o\`u $\alpha$ satisfait $3\alpha^2=1.$

Pour r\'esumer nous obtenons la:

\begin{pro}
{\sl Soit $\mathscr{C}$ une cubique \`a point double. L'ensemble
des transformations $A\sigma$ qui laissent invariante
$\mathscr{C}$ est fini.}
\end{pro}

\begin{rem}
Puisque deux cubiques \`a point double dans $\mathbb{P}^2(\mathbb{C})$ sont
lin\'eairement conjugu\'ees une telle cubique est
invariante par une infinit\'e de transformations quadratiques de
$\Sigma^3;$ on peut pr\'eciser cela par un argument
essentiellement identique \`a celui des cubiques lisses (choix des
triplets $(p_1,p_2,p_3)$ l'un des $p_i$ pouvant \^etre au point
double).

De la m\^eme fa\c{c}on on \'etablit que g\'en\'eriquement un
\'el\'ement de $\Sigma^3$ ne poss\`ede pas de cubique \`a point
double invariante.
\end{rem}

Finalement on en d\'eduit la:

\begin{pro}\label{courbinv}
{\sl Soit $A$ un automorphisme de $\mathbb{P}^2(
\mathbb{C})$ dont les coefficients sont
alg\'ebriquement ind\'ependants sur $\mathbb{Q}.$ La
transformation quadratique $A\sigma$ ne poss\`ede pas de courbe
irr\'eductible invariante.}
\end{pro}

En utilisant ce fait nous allons montrer un r\'esultat plus
g\'en\'eral.

\begin{thm}\label{courbinv2}
{\sl Soit $A$ un automorphisme de $\mathbb{P}^2(
\mathbb{C})$ dont les coefficients sont
$\mathbb{Q}$-alg\'ebrique\-ment ind\'ependants. La transformation
quadratique $A\sigma$ ne poss\`ede pas de courbe (\'eventuellement
non irr\'eductible) invariante.}
\end{thm}

\begin{proof}[\sl D\'emonstration]
Raisonnons par l'absurde; nous r\'eutilisons ici les notations $f$ pour
l'application quadratique de $\mathbb{C}^3$ dans $\mathbb{C}^3$ et
$f^\bullet$ pour l'application induite sur l'espace projectif $\mathbb{P}^2
(\mathbb{C}).$ Soient $\Gamma$ une courbe invariante
par $f^\bullet=A^\bullet\sigma^\bullet$ et
\begin{align*}
\Gamma=\Gamma_1\cup\ldots\cup\Gamma_s
\end{align*}
la d\'ecomposition de $\Gamma$ en facteurs irr\'eductibles, chaque
$\Gamma_i$ \'etant d\'ecrit par une \'equation irr\'eductible $h_i=0.$
Bien s\^ur la surface $h_i=0$ est invariante par l'application polynomiale
$f.$
Les coefficients de $A$ \'etant alg\'ebriquement ind\'ependants
sur $\mathbb{Q},$ la transformation $f^\bullet$ est alg\'ebriquement
stable; il s'en suit que
\begin{equation}\label{dem}
\forall\hspace{1mm} i\in\{1,\ldots,s\}, \hspace{3mm}
\exists\hspace{1mm} \ell \in\{1,\ldots,s\}\hspace{3mm}
\text{tel que}\hspace{3mm}f^\bullet(\Gamma_i)=\Gamma_\ell.
\end{equation}
\'Ecrivons que $\Gamma$ est invariante par $f^\bullet$
\begin{align*}
&h\circ f=h\eta && \text{avec} && h=h_1\ldots h_s;
\end{align*}
notons que le th\'eor\`eme de la fonction implicite, appliqu\'e en un
point g\'en\'erique de $\Gamma,$ assure que
$\eta$ n'est pas divisible par les $h_i.$ Remarquons que $\eta=0$
est form\'e de branches de
\begin{align*}
\mathrm{Exc}\hspace{1mm}f^\bullet=\{x_0=0,
\hspace{1mm}x_1=0,\hspace{1mm}x_2=0\}.
\end{align*}
En effet supposons que
$\gamma$ soit une branche des z\'eros de $\eta,$ alors $h(f(
\gamma))=0.$ Si $\gamma$ n'est pas contract\'ee par $f^\bullet$ alors
$f^\bullet(\gamma)$ co\"{\i}ncide avec l'un des $\Gamma_k;$ d'apr\`es
(\ref{dem}) il existe $1\leq i\leq s$ tel que $f^\bullet(\gamma)=
f^\bullet(\Gamma_i)$ ce qui n'est pas possible. L'invariance de $\Gamma$
se r\'e\'ecrit donc
\begin{align*}
&h(A\sigma)=\mu h x_0^px_1^qx_2^r, && (p,
q,r)\in\mathbb{N}^3\setminus\{0\},&&\mu\in\mathbb{C}^*.
\end{align*}
Puisque $A$ est g\'en\'erique $f^\bullet$ a quatre points fixes en position
g\'en\'erale. Apr\`es une conjugaison dynamique du type $(\alpha x_0:\beta
x_1:\gamma x_2)$ le point $(1:1:1)$ est fixe par $f=A\sigma\hspace{1mm}
\colon\hspace{1mm}\mathbb{C}^3\to\mathbb{C}^3.$ Ceci revient \`a supposer que la
matrice $A=\left[\begin{array}{ccc}
a_0 & b_0 & c_0 \\
a_1 & b_1 & c_1 \\
a_2 & b_2 & c_2
\end{array}
\right]$ est stochastique
\begin{equation}\label{stochastique}
a_i+b_i+c_i=1, \hspace{40mm} i=0,\hspace{1mm}1,\hspace{1mm}2.
\end{equation}
\'Evidemment cette conjugaison d\'etruit l'hypoth\`ese de
$\mathbb{Q}$-ind\'ependance des coefficients de $A$ mais puisque
cette modification est faite par conjugaison on garde par exemple la
propri\'et\'e qu'il n'y a pas de courbe irr\'eductible invariante. \`A partir de
\begin{align*}
h\circ f=\mu(x_0^px_1^qx_2^r)h
\end{align*}
on obtient
\begin{align*}
h(1,1,1)=\mu h(1,1,1).
\end{align*}
Si $h(1,1,1)$ est nul, $f$ laisse invariante la composante irr\'eductible
de $h=0$ passant par $(1,1,1)$ ce qui implique que $f^\bullet$ pr\'eserve
une courbe irr\'eductible, cas exclu par g\'en\'ericit\'e. Il en r\'esulte 
que $h(1,1,1)$ 
est non nul et que $\mu$ vaut $1.$ Comme l'un des entiers $p,$ $q,$ $r$
est non nul, nous allons supposer que $r\not=0.$ Consid\'erons
maintenant la transformation $\widetilde{f}=\widetilde{A}
\sigma$ d\'efinie par
\begin{align*}
(a_0(x_1-x_0)x_2+x_0x_1,a_1(x_1-x_0)x_2+x_0x_1,a_2(x_1-x_0)x_2
+x_0x_2).
\end{align*}
Pour $a_0,$ $a_1$ et $a_2$ g\'en\'eriques $\widetilde{f}$ induit une
transformation birationnelle et on constate que
\begin{align*}
&\widetilde{f}(1,1,1)=(1,1,1), && \widetilde{f}(1,1,2)=(1,1,2).
\end{align*}
On remarque que le point $(1,1,2)$ n'est pas sur la surface $x_0^p
x_1^qx_2^r=1.$ En particulier parmi les transformations $g=B\sigma$
satisfaisant (\ref{stochastique}) la propri\'et\'e suivante
est g\'en\'erique
\begin{equation}\label{pointfixe}
\text{$g$ poss\`ede un point fixe $m$ situ\'e en dehors de $x_0^px_1^q
x_2^r=1.$}
\end{equation}
Il en est donc ainsi pour $f;$ en effet si ce n'\'etait pas le cas, alors,
pour tout automorphisme de corps $\kappa,$ les \'el\'ements $f^\kappa=
A^\kappa\sigma$ auraient tous leurs points fixes sur $x_0^px_1^qx_2^r=1$
ce qui contredirait la propri\'et\'e (\ref{pointfixe}). Pour un tel point fixe
$m$ on a
\begin{align*}
h(m)=h(f(m))=(x_0^px_1^qx_2^r)(m)h(m),
\end{align*}
\'egalit\'e qui conduit \`a $h(m)=0.$ Le m\^eme raisonnement que
pr\'ec\'edemment assure que $f$ pr\'eserve la composante de
$h=0$ passant par $m$ et donc que $f^\bullet$ a une courbe
irr\'eductible invariante: contradiction avec la Proposition \ref{courbinv}.
\end{proof}

\begin{thm}
{\sl Soit $A\sigma$ une transformation de $\mathring{\mathrm{B}}
\mathrm{ir}_2$ alg\'ebriquement stable. Supposons que~$\Gamma_1,
\hspace{1mm}\ldots,\hspace{1mm}\Gamma_s$ soient des courbes (irr\'eductibles ou non)
invariantes par $A\sigma.$ D\`es que $s\geq 5,$ l'\'element
$A\sigma$ pr\'eserve une fibration.}
\end{thm}

\begin{proof}[\sl D\'emonstration]
Chaque $\Gamma_i$ est donn\'e par une \'equation r\'eduite
$h_i=0;$ l'invariance de $\Gamma_i$ par la transformation $f=A\sigma$ se traduit par
\begin{align*}
&h_i\circ f=\mu_i x_0^{p_i}x_1^{q_i}x_2^{r_i}h_i, && \mu_i\in\mathbb{C}^*,
\hspace{1mm}(p_i,q_i,r_i)\in\mathbb{N}^3\setminus\{0\}.
\end{align*}
On en d\'eduit pour tous les entiers $n_1,$ $\ldots,$ $n_s$
l'\'egalit\'e
\begin{align*}
h_1^{n_1}\ldots h_s^{n_s}\circ f=\mu_1^{n_1}\ldots\mu_s^{n_s}
h_1^{n_1}\ldots h_s^{n_s}x_0^{\sum_{i=1}^sn_ip_i}
x_1^{\sum_{i=1}^sn_iq_i}x_2^{\sum_{i=1}^sn_ir_i}.
\end{align*}
D\`es que $s\geq 5$ on peut trouver des entiers $n_i$ dans
$\mathbb{Z}$ non triviaux tels que
\begin{align*}
&\displaystyle\sum_{i=1}^s n_i\left[
\begin{array}{c}
\deg h_i \\
p_i\\
q_i\\
r_i
\end{array}
\right]=0.
\end{align*}
L'\'el\'ement $f$ pr\'eserve alors la fibration d\'efinie par
$h_1^{n_1} \ldots h_s^{n_s}.$
\end{proof}

En faisant un raisonnement analogue sur les branches de $\det
(\mathrm{jac} f)$ pour une transformation birationnelle $f$ de
degr\'e quelconque on obtient le:

\begin{thm}
{\sl Soit $f$ une transformation de \textsc{Cremona} alg\'ebriquement
stable de degr\'e~$\nu.$ Il existe un entier $N(\nu)$ tel que: d\`es que $f$ poss\`ede
plus de $N(\nu)$ courbes invariantes $f$ pr\'eserve une fibration.}
\end{thm}

\section[Transformations quadratiques pr\'eservant une fibration
rationnelle]{Exemples de transformations quadratiques
pr\'eservant une fibration rationnelle}\hspace{0.1cm}

Donnons quelques exemples de transformations
birationnelles pr\'eservant des fibrations rationnelles qui
peuvent \^etre des fibrations en coniques, en cubiques \`a point
double, en cubiques cuspidales, etc. Ces exemples proviennent de
la classification des flots (\emph{voir} Chapitre
\ref{germgerm}).\bigskip

\begin{itemize}
\item Soit $\beta$ un rationnel; la famille de transformations
\begin{align*}
g_\beta=\left(\left(x_0+ \frac{1-e^{-(\beta+1)t}}{(\beta+1)x_1}
\right)e^{\beta t},e^tx_1 \right)
\end{align*}
pr\'eserve la fibration rationnelle donn\'ee par
$\frac{x_1^{1+\beta}}{1+(1+\beta)x_0x_1}=\text{ cte}$ fibre \`a
fibre. Elle pr\'eserve aussi $x_1=$ cte mais en d\'epla\c{c}ant
les fibres. On constate que $g_2$ laisse une fibration en cubiques
\`a point double invariante: dans la carte $x_0=1$ la famille des
courbes invariantes est donn\'ee par $\alpha
x_1^3-x_2^3-3x_1x_2=0.$ Pour $1+\beta=n$ entier cela donne la
fibration 
\begin{align*}
\alpha x_1^n-x_2^n -nx_1x_2^{n-2}=0.
\end{align*}
\bigskip

\item Pour tout $\beta$ rationnel la famille $\left(
\beta\frac{(e^{-2\beta t} +1)x_0+\beta(e^{-2\beta t}-1)}{(e^{-2\beta t}-1)x_0
+\beta(e^{-2\beta t}+1)},e^tx_1\right)$ pr\'eserve la fibration
rationnelle $x_1^{2\beta}\left(\frac{x_0+\beta}{x_0-\beta}
\right)=\text{ cte}$ fibre \`a fibre; lorsque $\beta$ vaut $1$ c'est une fibration
en cubiques \`a point double.\bigskip

\item On constate que les transformations $\left( \frac{x_1^2t
 +x_1t^2+x_0x_1+\frac{t^3}{3}}{x_1+t},x_1+t\right)$ laissent invariante la fibration
en cubiques \`a point double $x_1^3-3x_0x_1=\text{ cte}$ fibre \`a
fibre.\bigskip

\item La famille de transformations $\left(\frac{x_0(e^tx_1+1)}
{(x_1+1)e^{t}}, e^tx_1\right)$ laisse la fibration en coniques
$\frac{x_0x_1}{x_1+1}=\text{ cte}$ invariante fibre \`a fibre; dans
la carte $x_1=1,$ elle
est donn\'ee par $\frac{x_0}{x_2(1+x_2)}=\text{ cte}$ qui est une
famille de paraboles tangentes. \bigskip

\item L'\'el\'ement $\left(x_0+t(\beta+x_1^2)+t^2x_1+\frac{t^3}
{3},x_1+t \right)$ de $\mathrm{Bir}_2$ pr\'eserve la fibration
en cubiques cuspidales $3\beta x_1+x_1^3-3x_0=\text{ cte}$
(se placer dans la carte affine $x_0=1$).\bigskip

\item Soit $\beta$ dans $\mathbb{Q}\setminus\{2\};$ on note que $\left(x_0e^{\beta t}+
\frac{x_1^2}{2-\beta}(e^{2t}-e^{\beta t}),x_1e^t\right)$ pr\'eserve la
fibration rationnelle $\frac{x_1^\beta}{x_1^2+(\beta-2)x_0}=\text{
cte}$ fibre \`a fibre.\bigskip

\item Les transformations

\begin{footnotesize}
\begin{align*}
&\left(\frac{x_0(x_1+1)e^{t}}{(e^tx_1+1)},e^tx_1 \right), && \text{resp.
} \left(\frac{e^t(4x_0x_1+2x_0-1)+2x_0+1}{2\left(e^t(2x_1+1-2x_0)+1+2x_0
\right)},e^t x_1\right), &&\text{resp.
}\left(\frac{x_0x_1}{x_1-x_0+e^{-t}x_0}, e^t x_1 \right)
\end{align*}
\end{footnotesize}

\noindent pr\'eservent fibre \`a fibre les fibrations en coniques
\begin{align*}
&\frac{x_1}{x_0(x_1+1)}=\text{ cte},&& \text{resp. } \frac{x_1(1+2x_0)}
{(x_1+1)(1-2x_0)}=\text{ cte},&&\text{resp. }\frac{x_0x_1}{x_1-x_0}= \text{
cte}.
\end{align*}
\end{itemize}

\section{Transformations ayant une courbe de points fixes}\hspace{0.1cm}

Dans ce paragraphe nous examinons les diff\'erentes courbes
de points fixes possibles pour une transformation quadratique.

Comme on l'a vu au Lemme \ref{fix} d\`es que $A$ est un
automorphisme de $\mathbb{P}^2(
\mathbb{C})$ dont les coefficients sont
alg\'ebriquement ind\'ependants sur $\mathbb{Q},$ la transformation $(A
\sigma)^n$ (resp. $(A\rho)^n$) n'a pas de courbe de points fixes,
ph\'enom\`ene que l'on observe dej\`a  pour les automorphismes.
Commen\c{c}ons par \'enoncer le:

\begin{thm}[\cite{Pa}, Th\'eor\`eme 1.3., \cite{BPV}, Th\'eor\`eme 1.1.]
{\sl Soient $\mathcal{C}$ une courbe alg\'ebrique irr\'eductible dans
$\mathbb{P}^2(
\mathbb{C})$ et $f$ une transformation de \textsc{Cremona} non
lin\'eaire laissant $\mathcal{C}$ invariante.
\begin{itemize}
\item Si $\mathcal{C}$ est lisse de genre plus grand ou \'egal \`a
$1,$ alors $\mathcal{C}$ est une cubique lisse.

\item Si le genre de $\mathcal{C}$ est sup\'erieur ou \'egal \`a
$2$ et si $f$ est l'identit\'e sur $\mathcal{C},$ alors ou bien
$f$ est conjugu\'ee \`a une transformation de \textsc{de
Jonqui\`eres}, ou bien $f$ est p\'eriodique de p\'eriode
inf\'erieure ou \'egale \`a $3.$
\end{itemize}}
\end{thm}

\begin{eg}
Soit $P$ un polyn\^ome de degr\'e $2d+2$ sans racine multiple; la
transformation d\'efinie par $\left(\frac{x_0+P(x_1)}{x_0+1},x_1\right)$ laisse la courbe
hyperelliptique $x_0^2=P(x_1)$ invariante.
\end{eg}

Rappelons qu'une courbe invariante point par point par
une transformation birationnelle quadratique est n\'ecessairement
de degr\'e inf\'erieur ou \'egal \`a $2$ (Proposition
\ref{degfeuil}). En fait dans~\cite{Bl} \textsc{Blanc} montre le:

\begin{thm}[\cite{Bl}]
{\sl Soient $\mathcal{C}$ une cubique plane lisse et $\mathrm{G}$
le groupe des transformations birationnelles qui fixent
$\mathcal{C}$ point par point. Un \'el\'ement non trivial de
$\mathrm{G}$ est de degr\'e au moins~$3.$

De plus les \'el\'ements de degr\'e $3$ de $\mathrm{G}$
engendrent $\mathrm{G}.$}
\end{thm}

On d\'ecrit maintenant les transformations quadratiques
qui ont une courbe de points fixes; les trois \'enonc\'es qui
suivent s'obtiennent par des calculs \'el\'ementaires.

Comme on l'a d\'ej\`a dit pour que l'image d'une droite par $\sigma$ soit encore une
droite il faut qu'elle passe par l'un des trois points
d'ind\'etermination de $\sigma$ et soit distincte de $x_0=0,$
$x_1=0$ et $x_2=0;$ de m\^eme l'image par $\sigma$ d'une
conique $\mathcal{C}$ est encore une conique si $\mathcal{C}$
passe par deux des points d'ind\'etermination de $\sigma.$ D'o\`u
la:

\begin{pro}
{\sl \`A conjugaison dynamique par permutation et homoth\'eties pr\`es
on a
\begin{itemize}
\item les transformations $A\sigma$ qui poss\`edent une droite de
points fixes s'\'ecrivent

\begin{small}
\begin{align*}
&(\alpha(x_1-x_0)x_2+x_0x_1:(x_1-x_0)x_2+x_0x_1: \gamma(x_1-x_0)
x_2+x_0x_2), &&\alpha,\hspace{0.1cm}\gamma\in\mathbb{C}\setminus\{1\},\hspace{0.1cm}\alpha\not=
\gamma
\end{align*}
\end{small}

(ici la droite de points fixes est la droite d'\'equation $x_0=x_1$);

\item la transformation $A\sigma$ a une conique lisse de points
fixes si elle est du type
\begin{align*}
&(x_0x_1:\alpha x_1x_2-x_0x_2+\beta x_0x_1:x_1x_2), &&\alpha,\hspace{0.1cm}
\beta\in\mathbb{C},\hspace{0.1cm}\alpha\beta\not=1.
\end{align*}
La conique de points fixes est alors donn\'ee par: $x_1^2-\alpha x_1x_2+
x_0x_2-\beta x_0x_1=0.$ Cette conique d\'eg\'en\`ere en deux droites
lorsque $\alpha\beta=1.$
\end{itemize}}
\end{pro}

\begin{rem}
On constate que l'\'el\'ement $(x_0x_1:\alpha x_1x_2-x_0x_2+\beta x_0x_1:x_1x_2)$ pr\'eserve la
fibration $x_0/x_2=\text{ cte}$ fibre \`a fibre et ceci montre que tout \'el\'ement de
$\Sigma^3$ ayant une conique de points fixes pr\'eserve fibre \`a fibre une
fibration en droites. En fait en utilisant les
deux propositions qui suivent nous allons voir que
toute transformation birationnelle
quadratique ayant une conique lisse de points fixes pr\'eserve une
fibration en droites fibre \`a fibre.
\end{rem}

\bigskip

La caract\'erisation des droites (resp. coniques) qui sont transform\'ees
par $\rho$ en des droites (resp. coniques) permet d'\'etablir la:

\begin{pro}\label{courbptfix2}
{\sl La transformation $A\rho$ a une droite de points fixes si elle
s'\'ecrit
\begin{align*}
(a_0x_0x_1+b_0x_2^2+c_0x_1x_2:a_1x_0x_1+b_1x_2^2+c_1x_1x_2
:a_2x_0x_1+b_2x_2^2+ c_2x_1x_2),
\end{align*}
avec
\begin{align*}
& b_1a_2^2+a_2a_1c_1-a_1a_2b_2-a_1^2c_2=0
&& \text{ou}  && a_0b_0b_2+c_0b_2^2-a_2b_0^2-b_0b_2c_2=0.
\end{align*}

La transformation $A\rho$ a deux droites distinctes de
points fixes si elle est de l'une des formes suivantes
\begin{align*}
& (\alpha x_0x_1+x_1x_2:\beta x_0x_1+\gamma x_2^2+\delta x_1x_2:\varepsilon x_0x_1+\eta
x_1x_2), && (\eta-\alpha)^2+4\varepsilon\not=0;\\
& (\alpha x_0x_1+\beta x_2^2+\gamma x_1x_2:x_2^2+\delta x_1x_2:\varepsilon
x_2^2+\eta x_1x_2),&& (\varepsilon-\delta)^2+4\eta \not=0.
\end{align*}

Enfin la transformation $A\rho$ poss\`ede une conique lisse de
points fixes si elle est de l'un des types suivants
\begin{align*}
& (x_0x_1:-x_0x_1-x_2^2-\beta x_1x_2:x_1x_2), 
&& (x_0x_1+\alpha x_2^2+\alpha x_1x_2:x_1x_2:x_2^2), && \beta\in
\mathbb{C},&&\alpha\in\mathbb{C}^*.
\end{align*}
Dans le premier cas la conique est donn\'ee par $x_0x_1+x_2^2+x_1^2+\beta
x_1x_2=0,$ dans le second par~$x_0x_1+\alpha x_2^2+\alpha x_1x_2-x_0x_2=0.$}
\end{pro}

En ce qui concerne les transformations de $\Sigma^1$ on a la:

\begin{pro}\label{courbptfix3}
{\sl Soit $A=(a_0x_0+b_0x_1+c_0x_2:a_1x_0+b_1x_1+c_1x_2:a_2x_0+b_2x_1
+c_2x_2)$ un automorphisme de $\mathbb{P}^2(\mathbb{C}).$

La transformation $A\tau$ poss\`ede une droite de points fixes si
$a_0c_0c_1+b_0c_1^2-a_1c_0^2-b_1c_0c_1=0.$

La transformation $A\tau$ a une conique lisse de points
fixes si $A\tau$ est de la forme
\begin{align*}
&(2x_0^2:2x_0x_1:\alpha x_0^2+\beta x_0x_1-x_1^2+x_0x_2), &&
\alpha,\hspace{0.1cm}\beta\in\mathbb{C};
\end{align*}
la conique en question a pour \'equation $x_0x_2-\alpha x_0^2-\beta x_0x_1
+x_1^2=0.$}
\end{pro}

\begin{rem}
Les it\'er\'es de $(2x_0^2:2x_0x_1:\alpha x_0^2+\beta x_0x_1-x_1^2+x_0x_2)$ sont
dans $\mathrm{Bir}_2.$ En effet dans la carte $x_0=1$ cette application induit la
transformation polynomiale
\begin{align*}
\left(x_1,x_2+\frac{\alpha}{2}+\frac{\beta}{2}x_1-x_1^2\right).
\end{align*}
\end{rem}

\section{Points d'ind\'etermination des it\'er\'es, points p\'eriodiques}\label{dynpt}\hspace{0.1cm}

\subsection{Points d'ind\'etermination, ensembles exceptionnels}\hspace{0.1cm}

Soit $f$ un \'el\'ement du groupe de \textsc{Cremona}.
On note $\mathrm{Ind}^+ \hspace{0.1cm} f=\displaystyle\bigcup_{n\geq
1}\mathrm{Ind}\hspace{0.1cm} f^n\label{not39}$ l'union des ensembles
d'ind\'etermination des it\'er\'es positifs de $f.$ De m\^eme on
introduit $\mathrm{Ind}^- \hspace{0.1cm} f=\displaystyle\bigcup_{n\geq
1}\mathrm{Ind}\hspace{0.1cm} f^{-n}.\label{not40}$ De fa\c{c}on similaire on
d\'efinit $\mathrm{Exc}^+\hspace{0.1cm} f=\displaystyle\bigcup_{n\geq
1}\mathrm{Exc}\hspace{0.1cm} f^n\label{not41}$ l'ensemble des courbes
contract\'ees par les it\'er\'es positifs de $f$ ainsi que
$\mathrm{Exc}^-\hspace{0.1cm} f=\displaystyle \bigcup_{n\geq
1}\mathrm{Exc}\hspace{0.1cm} f^{-n}.\label{not42}$

Soit $A$ un automorphisme de $\mathbb{P}^2(\mathbb{C})$ dont les
coefficients sont alg\'ebriquement ind\'ependants sur
$\mathbb{Q}.$ La transformation $A\sigma$ est alors
alg\'ebriquement stable; soit $p_1,$ $p_2$ et $p_3$ les points
d'ind\'etermination de $A\sigma.$ On remarque que pour $n$ dans
$\mathbb{N}$ les $(A\sigma)^{-n}(p_i)$ sont des points
d'ind\'etermination de $(A\sigma)^{n+1}.$ En particulier la
$\mathbb{Q}$-ind\'ependance implique que l'ensemble
$\mathrm{Ind}^+(A \sigma)$ est infini; il en est de m\^eme
pour $\mathrm{Ind}^-(A\sigma).$

\begin{thm}\label{orbneg}
{\sl Soit $A$ un automorphisme de $\mathbb{P}^2(
\mathbb{C})$ dont les coefficients sont
alg\'ebriquement ind\'ependants sur $\mathbb{Q}.$ Alors
$\mathrm{Ind}^-(A\sigma)$ est \textsc{Zariski} dense.}
\end{thm}

Nous proposons deux fa\c{c}ons d'\'etablir ce r\'esultat.

\begin{proof}[\sl D\'emonstration]
Supposons que $\mathrm{Ind}^-(A\sigma)$ soit contenu dans un
sous-ensemble alg\'ebrique propre $\Gamma$ de
$\mathbb{P}^2(\mathbb{C}).$ Alors il en est de m\^eme pour
$\mathrm{Ind}^-(A^\kappa\sigma)$ pour tout automorphisme $\kappa$
du corps $\mathbb{C}.$ Par densit\'e il en est encore de m\^eme pour tout
\'el\'ement de $\Sigma^3.$ Or nous donnons par la suite des
exemples du type
\begin{align*}
\left(\frac{ax_1+b}{cx_1+d},\frac{\alpha x_2+\beta}{\gamma
x_2+\delta}\right)
\end{align*}
pour lesquels ce n'est pas vrai donc $\mathrm{Ind}^-(A\sigma)$
est \textsc{Zariski} dense.

La deuxi\`eme approche est la suivante.
Soit $\overline{\mathrm{Ind}^-(A\sigma)}^{\hspace{0.1cm}\mathsf{Z}}$
l'adh\'erence de \textsc{Zariski} de $\mathrm{Ind}^-(A\sigma).$
Notons que $\overline{\mathrm{Ind}^-(A\sigma)}^{\hspace{0.1cm}\mathsf{Z}}$
est invariant par $A\sigma.$ Supposons que
$\overline{\mathrm{Ind}^-(A\sigma)}^{\hspace{0.1cm}\mathsf{Z}}$ soit strictement
contenu dans $\mathbb{P}^2(\mathbb{C});$ comme il est infini on a
\begin{align*}
\overline{\mathrm{Ind}^-(A\sigma)}^{\hspace{0.1cm}\mathsf{Z}}=\Gamma\cup\{p_1,
\ldots,p_s\}
\end{align*}
o\`u $\Gamma$ d\'esigne une courbe alg\'ebrique et les $p_i$ des
points de $\mathbb{P}^2(\mathbb{C})$ non situ\'es sur $\Gamma.$ La
$\mathbb{Q}$-d\'ependance fait que $\Gamma$ n'est pas contract\'ee
par $A\sigma.$ On constate alors que $\Gamma$ est invariante par
$A\sigma;$ or g\'en\'eriquement une transformation du type
$A\sigma$ ne laisse pas de courbe invariante (Th\'eor\`eme
\ref{courbinv2}): contradiction.
\end{proof}

Comme le montre la remarque qui suit on ne peut esp\'erer sous
les seules hypoth\`eses de la $\mathbb{Q}$-ind\'ependance obtenir la
densit\'e (au sens ordinaire) des points d'ind\'etermination.

\begin{rem}
Soit $A$ un automorphisme de $\mathbb{P}^2(\mathbb{C})$ dont les
coefficients sont alg\'ebriquement ind\'ependants sur $\mathbb{Q}.$
Supposons que $A\sigma$ ait un point fixe $m$ tel que les valeurs
propres de la partie lin\'eaire de $f$ en $m$ soient de module $1$
mais g\'en\'eriques; $A\sigma$ est alors lin\'earisable au
voisinage de ce point. En particulier il existe un ouvert
(polydisque) invariant par $A \sigma;$ pour tout~$n$ les
\'el\'ements de $\mathrm{Ind}(A \sigma)^n$ et $\mathrm{Exc}
(A\sigma)^n$ ne peuvent rencontrer cet ouvert.

Si $A$ est un \'el\'ement de $\mathrm{SL}_3(\mathbb{Z})$
(resp. $\mathrm{SL}_3(\mathbb{R})$) les points d'ind\'etermination
des $(A\sigma)^n$ sont rationnels (resp. r\'eels). En particulier
les $(A\sigma)^n$ sont holomorphes dans le compl\'ement de
$\mathbb{P}^2(\mathbb{R})$ dans $\mathbb{P}^2(\mathbb{C})$; mais
cet ensemble n'est pas invariant puisque les courbes
contract\'ees, qui l'intersectent \'evidemment, le sont sur des
points r\'eels. Par contre le compl\'ement des courbes
contract\'ees priv\'e de $\mathbb{P}^2(\mathbb{R})$ est invariant.
\end{rem}

On se convainc facilement qu'il n'est pas ais\'e de
d\'ecrire les ensembles $\mathrm{Exc}^\pm$ et $\mathrm{Ind}^{\pm};$
c'est possible dans des situations tr\`es sp\'eciales non g\'en\'eriques
en g\'en\'eral. Un des exemples les plus simples est le suivant
\begin{align*}
h=(h_1,h_2)=\left(\frac{ax_1+b}{cx_1+d},\frac{\alpha x_2+\beta}{\gamma x_2+
\delta}\right).
\end{align*}
On note que $h$ fait partie d'un sous-groupe de $\mathrm{Bir}_2$ ici 
$\mathrm{PGL}_2(\mathbb{C})\times
\mathrm{PGL}_2(\mathbb{C})$ et peut \^etre vu comme automorphisme mais 
sur $\mathbb{P}^1(\mathbb{C})\times\mathbb{P}^1(
\mathbb{C}).$ On constate d\'ej\`a dans cette situation des ph\'enom\`enes
\og int\'eressants \fg.

Si les points fixes de chaque $h_i$ sont distincts de
$0$ et $\infty$ dans $\mathbb{P}^1(
\mathbb{C})$ alors $h$ est dans $\Sigma^3,$ les
trois droites contract\'ees par $h$ \'etant $x_1=-\frac{d}{c},$ $x_2=-
\frac{\delta}{\gamma}$ et la droite \`a l'infini $x_0=0.$ Supposons que chaque
$h_i$ soit une rotation irrationnelle; alors les adh\'erences (ordinaires)
\begin{align*}
&\overline{\left\{h_1^{-n}\left(-\frac{d}{c}\right)\hspace{0.1cm}|\hspace{0.1cm}
n\in\mathbb{N}\right\}}
&&\text{et}&&\overline{\left\{h_2^{-n}\left(
-\frac{\delta}{\gamma}\right)\hspace{0.1cm}|\hspace{0.1cm} n\in\mathbb{N}\right\}}
\end{align*}
sont des cercles. Il en r\'esulte que
\begin{align*}
\overline{\mathrm{Exc}^+h}=\overline{\bigcup_{n\in\mathbb{N}}\mathrm{Exc}\hspace{0.1cm}
h^n},
\end{align*}
qui est aussi dans ce cas
\begin{align*}
\overline{\mathrm{Exc}^-h}=\overline{\bigcup_{n\in\mathbb{N}}\mathrm{Exc}\hspace{0.1cm}
h^{-n}},
\end{align*}
est l'union de deux c\^ones r\'eels quadratiques $H_1$ et $H_2$
dans $\mathbb{P}^2(\mathbb{C})$ et de la droite \`a l'infini~$x_0=0$
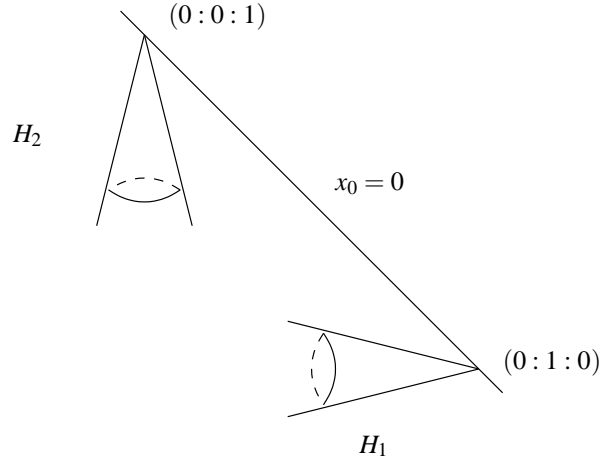
\begin{figure}[H]
\begin{center}
\input{exsimple.pstex_t}
\end{center}
\end{figure}

L'intersection des deux droites $x_1=h_1^{-n}\left(-
\frac{d}{c}\right)$ et $x_2=h_2^{-n}\left(-\frac{\delta}{\gamma}\right)$
contract\'ees par $h^n$ est un point d'ind\'etermination de
$h^n$ (les deux autres \'etant $(0:1:0)$ et $(0:0:1)).$ Par suite l'adh\'erence
\begin{align*}
\overline{\mathrm{Ind}^+h}=\overline{\bigcup_{n\in\mathbb{N}}\mathrm{Ind}\hspace{0.1cm}
h^n}=
\overline{\mathrm{Ind}^-h}=\overline{\bigcup_{n\in\mathbb{N}}\mathrm{Ind}\hspace{0.1cm}
h^{-n}}
\end{align*}
est constitu\'ee des points $(0:1:0),$ $(0:0:1)$ et de
l'intersection des quadriques r\'eelles $H_1$ et~$H_2;$ c'est une
surface quartique r\'eelle qui est topologiquement un tore
$\mathbb{S}^1\times\mathbb{S}^1.$\bigskip

On peut bien s\^ur d\'ecrire sans difficult\'e les $\overline{
\mathrm{Exc}^\pm\hspace{0.1cm} h}$ et $\overline{\mathrm{Ind}^\pm\hspace{0.1cm} h}$
pour n'importe quelles valeurs des param\`etres $a,$ $b,$ $c,$
$d,$ $\alpha,$ $\ldots$ \bigskip

Voici quelques exemples tir\'es de la classification des flots
quadratiques qui permettent de visualiser quelques configurations
d'ensembles $\mathrm{Ind}$ et $\mathrm{Exc}.$

\begin{itemize}
\item On consid\`ere
\begin{align*}
\phi_t=((2x_1+tx_2)x_0-\frac{t}{2}x_2^2:(x_1+tx_2)(-2tx_0+2x_1+tx_2):x_2(-2tx_0+2x_1+tx_2))
\end{align*}
Un calcul montre que
\begin{align*}
\mathrm{Exc}\hspace{0.1cm}\phi_t^n=\mathrm{Exc}\hspace{0.1cm}\phi_{nt}=\{x_1=0,\hspace{0.1cm}
x_1+ntx_2=0,\hspace{0.1cm} -2ntx_0+2x_1+ntx_2=0\}
\end{align*}
et
\begin{align*}
\mathrm{Ind}\hspace{0.1cm}\phi_t^n=\mathrm{Ind}\hspace{0.1cm}\phi_{nt}=\{(1:0:0),\hspace{0.1cm}(1:0:2),\hspace{0.1cm}(1:2nt:-2)\}
\end{align*}
ce qui permet une description imm\'ediate des ensembles
$\mathrm{Ind}^\pm$ et $\mathrm{Exc}^\pm.$  \vspace{0.6cm}

\item Pour le flot $\phi_t$ donn\'e par
\begin{footnotesize}
\begin{align*}
&(2e^t(t-2)x_0x_1-e^ttx_1^2+2(1-e^t)x_1x_2-4x_0x_2:
2e^tx_1(2tx_0-(2+t)x_1-2x_2):2x_2(2tx_0-(2+t)x_1-2x_2))
\end{align*}
\end{footnotesize}
on a
\begin{align*}
\mathrm{Exc}\hspace{1mm}\phi_t^n=\mathrm{Exc}\hspace{1mm}\phi_{nt}=
\{x_1+x_2=0,\hspace{1mm}e^{nt}x_1+x_2=0,\hspace{1mm}2ntx_0-(2+nt)x_1-2x_2=0\}
\end{align*}
et
\begin{align*}
\mathrm{Ind}\hspace{1mm}\phi_t^n=\{(1:0:0),\hspace{1mm}(1:2:-2),\hspace{1mm}
(2+nt-2e^{nt}:2nt:-2nte^{nt})\}
\end{align*}
Rappelons que pour cet exemple $\phi_t$ appartient \`a $\Sigma^3$
sauf pour $t\in 2\mathrm{i}\pi\mathbb{Z}$ o\`u $\phi_t$ est soit
l'identit\'e, soit dans $\Sigma^2.$ Remarquer que suivant les
valeurs de $t$ choisies la description des ensembles
$\mathrm{Ind}^\pm$ n'est pas si simple.
\end{itemize}

\subsection{Points p\'eriodiques}\hspace{0.1cm}

Soit $f$ un \'el\'ement du groupe de \textsc{Cremona}.
Un point \textbf{\textit{p\'eriodique}}\label{ind46} de $f$ est un
point $p$ tel que les germes $f^n_{,p}$ soient holomorphes et
$\mathcal{O}^+(f)=\{p,f(p),\ldots,f^n(p),\ldots\}$ soit un
ensemble fini qui ne rencontre pas $\mathrm{Ind}\hspace{1mm} f\cup
\mathrm{Ind}\hspace{1mm} f^{-1}.$ La p\'eriode de $p$ est par d\'efinition
le cardinal de $\mathcal{O}^+(f).$

\begin{rem}
Un point fixe $p$ de $f^n$ n'est pas n\'ecessairement
un point p\'eriodique de $f;$ il se pourrait en effet que pour
un certain $k<n,$ l'it\'er\'e $k$-i\`eme de $p$ par $f$ soit
d'ind\'etermination pour $f^{-1}$ par exemple.
\end{rem}

Un point p\'eriodique $p$ de p\'eriode $k$ est
\textbf{\textit{hyperbolique}}\label{ind47} si les valeurs propres $\delta_1(p)$ et
$\delta_2(p)$ de~$Df^k_{(p)}$ satisfont
\begin{align*}
|\delta_1(p)|<1<|\delta_2(p)|.
\end{align*}

\bigskip

L'ensemble des points p\'eriodiques hyperboliques d'un
automorphisme quadratique $f$ de \textsc{H\'enon} est
\textsc{Zariski} dense. En effet suivant \cite{FrMi} c'est un
ensemble infini dans $\mathbb{C}^2$ c'est-\`a-dire dans
$\mathbb{P}^2(\mathbb{C});$ s'il n'\'etait pas \textsc{Zariski}
dense $f$ laisserait une courbe invariante ce qui est impossible.
Consid\'erons $t\mapsto f_t$ une famille analytique de
transformations birationnelles telle que
\begin{itemize}
\item $f_0$ soit une transformation de \textsc{H\'enon},

\item pour tout $t$ non nul $f_t$ appartienne \`a $\Sigma^3.$
\end{itemize}
Soit $p_0$ un point p\'eriodique hyperbolique de p\'eriode $N$
pour $f_0.$ Fixons $\mathsf{D}$ un polydisque centr\'e en $p$ tel
que $\mathsf{D}$ ne rencontre pas la droite \`a l'infini. En
particulier
\begin{align*}
&\mathsf{D}\cap\mathrm{Exc}\hspace{1mm}f_0^n=\mathsf{D}\cap\mathrm{Ind}
\hspace{1mm}f_0^n=\emptyset,&&\forall\hspace{1mm}1\leq n\leq N.
\end{align*}
Si $N$ est fix\'e et $|t|$ suffisamment petit, $\mathrm{Exc}\hspace{1mm}f_t^N$ et
$\mathrm{Ind}\hspace{1mm}f_t^N$ ne rencontrent par $\mathsf{D}.$ En
appliquant le th\'eor\`eme de stabilit\'e des points fixes
hyperboliques \`a la famille de diff\'eomorphismes holomorphes
$f^N_{t|\mathsf{D}}$ on constate qu'il existe une application analytique
$p\mapsto p(t)$ telle que pour $|t|$ suffisamment petit $p(t)$
soit un point fixe hyperbolique de $f_t$ de p\'eriode~$N.$

Soit $A$ un automorphisme de $\mathbb{P}^2(
\mathbb{C})$ \`a coefficients
alg\'ebriquement ind\'ependants sur $\mathbb{Q}.$ Puisque $\{A^\kappa
\hspace{1mm}|\hspace{1mm}\kappa\in\mathrm{Aut}(\mathbb{C},+,.)\}$ est topologiquement dense, si $n$
est fix\'e, il existe un automorphisme $\kappa$ du corps $\mathbb{C}$ tel
qu'un conjugu\'e lin\'eaire $g$ de $A^\kappa\sigma$ soit tr\`es
proche de $f_0;$ par suite $g^n$ est proche de $f_0^n$ et
poss\`ede au moins autant de points fixes que $f_0^n$ a de points
fixes hyperboliques (ou plus g\'en\'eralement de points fixes simples).
Il en est de m\^eme pour $(A^\kappa\sigma)^n.$ Cette construction
produit donc pour $A\sigma$ au moins autant de points p\'eriodiques
de p\'eriode $n$ (pas n\'ecessairement tous hyperboliques) que $f_0$
a de points p\'eriodiques simples. On peut \'evidemment faire cela pour
tout~$n.$

\begin{thm}\label{perper}
{\sl Soit $A$ un automorphisme de $\mathbb{P}^2(
\mathbb{C})$ dont les coefficients sont
alg\'ebriquement ind\'ependants sur $\mathbb{Q}.$ L'ensemble des
points p\'eriodiques de $A\sigma$ est \textsc{Zariski} dense.}
\end{thm}

\begin{proof}[\sl D\'emonstration]
Posons $f=A\sigma.$ D'apr\`es ce qui pr\'ec\`ede $f$ a une
infinit\'e de points p\'eriodiques. Notons $\mathrm{Per}(f)$
l'ensemble des points p\'eriodiques de $f.$ L'id\'ee est toujours
la m\^eme. Si l'adh\'erence de \textsc{Zariski}
$\overline{\mathrm{Per}(f)}^{\hspace{1mm}\mathsf{Z}}$ de $\mathrm{Per}(f)$
est strictement contenue dans $\mathbb{P}^2(\mathbb{C}),$ 
alors~$\overline{\mathrm{Per}(f)}^{\hspace{1mm}\mathsf{Z}}$ s'\'ecrit comme
l'union d'une courbe alg\'ebrique $\Gamma$ et de points $p_i$ de
$\mathbb{P}^2(\mathbb{C})$ non situ\'es sur $\Gamma$
\begin{align*}
\overline{\mathrm{Per}(f)}^{\hspace{1mm}\mathsf{Z}}=\Gamma\cup\{p_1,\ldots,
p_s\}.
\end{align*}
Puisque $f(\mathrm{Per}(f))$ est contenu dans $\mathrm{Per}(f),$
la courbe $\Gamma$ est invariante par $f$ ce qui n'est pas
possible en vertu du Th\'eor\`eme \ref{courbinv2}. Il en
r\'esulte que $\mathrm{Per}(f)$ est \textsc{Zariski} dense.
\end{proof}

\begin{rem}
On peut aussi proc\'eder comme suit.

Posons $f=A\sigma.$ D'apr\`es ce qui pr\'ec\`ede $f$ a une
infinit\'e de points p\'eriodiques et une infinit\'e de p\'eriodes
distinctes apparait (ceci est vrai pour les applications de
\textsc{H\'enon}). Notons~$\mathrm{Per}(f)$ l'ensemble des points
p\'eriodiques de $f.$ Si $\mathrm{Per}(f)$ est contenu dans un
sous-ensemble alg\'ebrique propre $\Gamma,$ alors $\Gamma$ est du
type
\begin{align*}
\Gamma=\Gamma_1\cup\ldots\cup\Gamma_n\cup\{p_1, \ldots,p_s\}
\end{align*}
o\`u les $p_i$ sont des points isol\'es et les $\Gamma_i$ des
courbes irr\'eductibles. On en d\'eduit l'existence d'entiers $j$
et $k$ tels que $f^j(\Gamma_k)=\Gamma_k$ et l'it\'er\'e $f^j$ a
une infinit\'e de points p\'eriodiques avec une infinit\'e de
p\'eriodes distinctes; ceci implique que la normalis\'ee de
$\Gamma_k$ est $\mathbb{P}^1(\mathbb{C}).$ Comme un \'el\'ement
de $\mathrm{PGL}_2(\mathbb{C})$ ne peut poss\'eder une
infinit\'e de points p\'eriodiques avec p\'eriodes distinctes,
on obtient une contradiction.
\end{rem}

Dans \cite{Fa} \textsc{Favre} donne une estimation du nombre de
points p\'eriodiques d'un \'el\'ement g\'en\'erique du groupe de
\textsc{Cremona} de degr\'e $d\geq 2;$ puis il montre, avec
\textsc{Diller}, un r\'esultat sur les transformations
bim\'eromorphes d'une surface de \textsc{K\"ahler} compacte que
nous \'enon\c{c}ons dans le cas particulier de $\mathbb{P}^2(\mathbb{C}).$

\begin{thm}[\cite{DiFa}]
{\sl Soit $f$ une transformation de \textsc{Cremona} non lin\'eaire
alg\'ebriquement stable telle que $\lambda(f)>1.$ Supposons que
$f$ n'ait pas de courbe de points p\'eriodiques. Notons
$\mathrm{Per}_k$ le nombre de points p\'eriodiques de p\'eriode
(divisant) $k.$ Il existe une constante~$C>0$ telle que pour tout
$k\geq 0$ on ait
\begin{align*}
\vert\mathrm{Per}_k-\lambda(f)^k\vert\leq C.
\end{align*}}
\end{thm}

Remarquons que les Corollaire \ref{coras} et Lemme
\ref{fix} permettent d'appliquer ce Th\'eor\`eme \`a toute
transformation $A\sigma$ avec $A$ automorphisme g\'en\'erique de
$\mathbb{P}^2(\mathbb{C}).$

\begin{rem}
Il y a des transformations birationnelles de degr\'e $2$ dans $\Sigma^3$
qui sont alg\'ebriquement stable avec un nombre fini de points p\'eriodiques;
par exemple pour $\alpha,$ $\beta$ g\'en\'eriques les applications
$f_{\alpha,\beta}^2,$ o\`u
\begin{align*}
f_{\alpha,\beta}=\left(\frac{\alpha x_0+x_1}{x_0+1},\beta x_1\right)
\end{align*}\label{not35aa}
avec $\alpha,$ $\beta$ deux complexes g\'en\'eriques de module $1,$
sont de ce type. Leur degr\'e dynamique vaut $1.$
\end{rem}

\bigskip

La transformation $f$ v\'erifie la
\textbf{\textit{condition de \textsc{Bedford} et
\textsc{Diller}}}\label{ind22} (\emph{voir} \cite{BD}) si
\begin{align*}
\sum_{n\geq 0}\frac{1}{\lambda(f)^n}|\log(\textrm{dist}(f^n(
\mathrm{Ind}\hspace{1mm} f^{-1}),\mathrm{Ind}\hspace{1mm} f))|<\infty.
\end{align*}

Cette condition implique la stabilit\'e alg\'ebrique.

Le th\'eor\`eme qui suit est tr\`es important; il montre
l'abondance de points p\'eriodiques hyperboliques. Nous
l'\'enon\c{c}ons dans le cadre restreint de $\mathbb{P}^2(\mathbb{C}).$

\begin{thm}[\cite{BD, Du}]\label{BDDu}
{\sl Soit $f$ un \'el\'ement non lin\'eaire du groupe de
\textsc{Cremona}. Supposons que $f$ v\'erifie la condition de
\textsc{Bedford} et \textsc{Diller}. Alors $f$ poss\`ede une
infinit\'e de points p\'eriodiques hyperboliques qui
s'\'equidistribuent suivant une mesure de probabilit\'e
$f$-invariante.}
\end{thm}

Nous allons en d\'eduire la:

\begin{pro}\label{domino}
{\sl Soit $A$ un automorphisme de $\mathbb{P}^2(\mathbb{C})$ \`a coefficients r\'eels
strictement positifs. La transformation birationnelle quadratique
$A\sigma$ v\'erifie la condition de \textsc{Bedford} et
\textsc{Diller}; en particulier $A\sigma$ poss\`ede une infinit\'e
de points p\'eriodiques hyperboliques.}
\end{pro}

\begin{proof}[\sl D\'emonstration]
On note $f=A\sigma$ et $A=\left[
\begin{array}{ccc}
a_0 & b_0 & c_0 \\
a_1 & b_1 & c_1 \\
a_2 & b_2 & c_2
\end{array}
\right].$ Les points d'ind\'etermination~$a,\hspace{1mm}b$ et $c$ de
$f^{-1}$ sont donn\'es par les colonnes de $A$
\begin{align*}
& a=(a_0:a_1:a_2), && b=(b_0:b_1:b_2), && c=(c_0:c_1:c_2).
\end{align*}

Dans la carte affine $x_2=1,$ on note $\Theta$
l'ensemble r\'eel d\'efini par
\begin{align*}
\Theta:=\{x_0\geq 0,\hspace{1mm}x_1\geq 0\}.
\end{align*}
On remarque que $\sigma$ laisse invariant $\Theta.$
Comme les coefficients de $A$ sont positifs on a l'inclu\-sion~$A(\Theta)\subset\Theta;$ par suite
$f(\Theta)\subset\Theta.$ En fait $f(\Theta)$ est
le triangle de sommets $a,$ $b$ et $c.$ Il en r\'esulte que
$f^n(\mathrm{Ind}\hspace{1mm} f^{-1})\subset\Theta$ pour tout $n\geq 0.$
Soit $\textrm{dist}$\label{not35bb} la m\'etrique de
\textsc{Fubini}-\textsc{Study} de diam\`etre $1$ et $0<\varepsilon
<1$ la distance de $f(\Theta)$ \`a
\begin{align*}
\mathrm{Ind}\hspace{1mm}f=\{(1:0:0),\hspace{1mm}(0:1:0),
\hspace{1mm}(0:0:1)\};
\end{align*}
on a
\begin{align*}
\textrm{dist}(f^n(\mathrm{Ind}\hspace{1mm}f^{-1}),\mathrm{Ind}
\hspace{1mm}f)\leq\textrm{dist}(f(\Theta),\mathrm{Ind}\hspace{1mm}
f)=\varepsilon.
\end{align*}
Comme la s\'erie
\begin{align*}
\sum_{n\geq 0}\frac{1}{\lambda(f)^n}|\log(\textrm{dist}(f^n(\mathrm{Ind}
\hspace{1mm} f^{-1}),\mathrm{Ind}\hspace{1mm} f
))|\leq|\log\varepsilon|\sum_{n\geq 0}\frac{1}{2^n}
\end{align*}
est convergente on peut appliquer le Th\'eor\`eme \ref{BDDu}.
\end{proof}

Suite \`a une discussion avec Romain \textsc{Dujardin}, que nous
remercions, nous avons \'etabli le:

\begin{thm}
{\sl L'ensemble des \'el\'ements de $\mathrm{Bir}_2$ poss\'edant une
infinit\'e de points p\'eriodiques hyperboliques est dense dans
$\mathrm{Bir}_2.$}
\end{thm}

\begin{proof}[\sl D\'emonstration]
L'ensemble des transformations de $\mathrm{Bir}_2$ qui v\'erifient
la condition de \textsc{Bedford} et \textsc{Diller} est le
compl\'ement d'un ensemble \og pluri-polaire\fg\hspace{1mm} (\cite{BD}).
Lorsque $A$ est \`a coefficients r\'eels strictement positifs
$A\sigma$ satisfait cette condition, on en d\'eduit que cet
ensemble est non vide et donc dense: c'est une propri\'et\'e des
ensembles pluri-polaires.
\end{proof}

\begin{rem}
On peut d\'emontrer la densit\'e des transformations de
$\mathrm{Bir}_2$ poss\'edant une infinit\'e de points
p\'eriodiques en partant d'une transformation $A\sigma$ avec $A$
dans $\mathrm{PGL}_3(\mathbb{R})$ \`a coefficients
alg\'ebriquement ind\'ependants et en faisant agir
$\mathrm{Aut}(\mathbb{C},+,.).$ Mais cette d\'emarche ne
contr\^ole pas l'hyperbolicit\'e.
\end{rem}

\section{Transformations birationnelles quadratiques de carr\'e
quadratiques}\label{quad}\hspace{1mm}

On note que pour toute transformation de \textsc{Cremona} $f$
telle que
\begin{align*}
& \deg f^n=2&& \forall\hspace{1mm}n\geq 1
\end{align*}
il existe un entier $k$ tel que $f^k$ se plonge dans un flot. En
effet $\overline{\langle f^n\rangle}^{\hspace{1mm}\mathsf{Z}}$ est
un sous-ensemble alg\'ebrique de $\mathbb{P}^{17}(
\mathbb{C})$ et
$\overline{\langle f^n\rangle}^{\hspace{1mm}
\mathsf{Z}}\cap\mathrm{Bir}_2$ est un groupe alg\'ebrique. Par
suite $\overline{ \langle f^n\rangle}^{\hspace{1mm}\mathsf{Z}}
\cap\mathrm{Bir}_2$ a un nombre fini de composantes connexes; une
puissance de $f$ est dans la composante connexe de l'identit\'e et
cette puissance se plonge dans un flot.

Nous avons vu que la non stabilit\'e alg\'ebrique \'etait reli\'ee
\`a une baisse du degr\'e des it\'er\'es. Dans cet ordre d'id\'ee
nous allons d\'ecrire les transformations quadratiques de
$\Sigma^3$ (resp. $\Sigma^2,$ resp. $\Sigma^1$) dont le carr\'e
est encore dans $\mathrm{Bir}_2.$ Ceci se ram\`ene \`a chercher
les automorphismes $A$ de $\mathbb{P}^2(
\mathbb{C})$ pour lesquels $(A\sigma)^2$
(resp. $(A\rho)^2,$ resp. $(A\tau)^2$) est encore dans
$\mathrm{Bir}_2;$ cela produit des exemples de familles de
transformations birationnelles non alg\'ebriquement stables. Les
flots rencontr\'es au Chapitre \ref{germgerm} sont comme cela.

Commen\c{c}ons par un exemple. La famille des transformations
birationnelles $f_{\alpha, \beta}$ d\'efinie par
\begin{align*}
f_{\alpha,\beta}=\left(\frac{\alpha x_0+x_1}{x_0+1},\beta x_1\right)
\end{align*}
o\`u $\alpha,$ $\beta$ d\'esignent des nombres complexes
a \'et\'e \'etudi\'ee dans \cite{De2}
\begin{align*}
\deg f_{\alpha,\beta}^2=2,\hspace{1mm}\deg
f_{\alpha,\beta}^3=3,\hspace{1mm}\deg
f_{\alpha,\beta}^4=3,\hspace{1mm} \deg
f_{\alpha,\beta}^5=4,\hspace{1mm}\deg
f_{\alpha,\beta}^6=4,\hspace{1mm}\deg
f_{\alpha,\beta}^7=5,\hspace{1mm}\ldots
\end{align*}
Lorsque $\alpha$ et $\beta$ sont de la forme $\exp(2\mathrm{i}
\pi\eta)$ et $\exp(2\mathrm{i}\pi \widetilde{\eta})$ avec $\eta,$
$\widetilde{\eta}$ dans $\mathbb{R}\setminus\mathbb{Z},$ cette
famille a une dynamique curieuse (\cite{De2}) que nous
rappelerons bri\`evement au Chapitre \ref{expexp}. Cette famille
montre que la condition $\deg f^2=2$ n'est pas suffisante pour que
tous les it\'er\'es soient de degr\'e~$2.$

\bigskip

Examinons de plus pr\`es les \'el\'ements de $\Sigma^3$
dont le carr\'e est quadratique. La transforma\-tion~$(B\sigma C)^2$
est quadratique si et seulement si $\sigma CB\sigma$ l'est; on se
ram\`ene donc \`a d\'eterminer les \'el\'ements $A$ de
$\mathrm{PGL}_3(\mathbb{C})$ tels que $\sigma A\sigma$ soit
quadratique.

\begin{pro}\label{sigqua}
{\sl Les \'el\'ements $A$ de $\mathrm{PGL}_3(\mathbb{C})$ tels que $\deg\sigma A
\sigma\leq 2$ s'\'ecrivent
$f\ell g$ avec
\begin{align*}
& \ell=(x_0:\alpha x_0+\beta x_1:\gamma x_0+\delta
x_2),&&\hspace{1mm}
\alpha,\hspace{1mm}\gamma\in\mathbb{C},\hspace{1mm}\beta,\hspace{1mm}
\delta\in\mathbb{C}^*&&f,\hspace{1mm} g\in \mathscr{S}_6.
\end{align*}}
\end{pro}

\begin{proof}[\sl D\'emonstration]
Posons $A\sigma:=(F_0:F_1:F_2).$ La transformation
$\sigma A\sigma$ est de degr\'e inf\'erieur ou \'egal \`a $2$ si
et seulement s'il existe $\psi$ et $q_i$ dans $\mathbb{C}[x_0,x_1,x_2]_2$ tels
que
\begin{align*}
(F_1F_2:F_0F_2:F_0F_1)=\psi(q_0:q_1:q_2).
\end{align*}

Si $\psi$ \'etait irr\'eductible, $\psi$ diviserait deux
$F_i$ distincts; $A\sigma$ ne serait alors pas inversible ce
qui est absurde. Ainsi $\psi$ s'\'ecrit $\psi_0\psi_1,$ les $\psi_i$
d\'esignant des formes lin\'eaires.

Si $\psi_0$ et $\psi_1$ co\"incident \`a multiplication
pr\`es par un scalaire, alors ou bien deux des composants
de $A\sigma$ sont multiples l'une de l'autre, ou bien $A\sigma$
est lin\'eaire; ces deux
\'eventualit\'es \'etant exclues, $\psi_0$ et $\psi_1$ ne sont pas
multiples l'une de l'autre.

Quitte \`a r\'eindicer les $\psi_i$ on peut supposer que
$\psi_0$ divise l'un des $F_i$ et $\psi_1$ les deux autres. \`A
permutation pr\`es on a
\begin{align*}
&F_0=\psi_0f_0,&& F_1=\psi_1f_1,&&F_2=\psi_1f_2,
\end{align*}
les $f_i$ \'etant des formes lin\'eaires. On obtient alors
\begin{align*}
\sigma A\sigma= (\psi_1f_1f_2:\psi_0f_0f_2:\psi_1f_0f_1);
\end{align*}
cette transformation est de degr\'e inf\'erieur ou \'egal \`a $2$
si l'un des facteurs de la $i$-\`eme composante divise les autres.
\'Etudions par exemple le cas o\`u $f_1$ divise $\psi_0f_0f_2$ et
$\psi_1f_0f_1.$ Remarquons que si $f_1$ divise $\psi_0$ alors
$A\sigma$ est lin\'eaire donc $f_1$ divise $f_0f_2.$ Consid\'erons
par exemple l'\'eventualit\'e: $f_1$ divise $f_0;$ alors
\begin{align*}
&\sigma A\sigma=(\psi_1f_2:\psi_0f_2:\psi_1f_1),&&
A\sigma=(\psi_0f_0: \psi_1f_0:\psi_0f_2).
\end{align*}
Par d\'efinition les $F_i$ sont du type $*x_1x_2+*x_0x_2+
*x_0x_1:$ ceci impose des conditions sur les $\psi_i$ et $f_i.$
Si $\psi_0$ est de la forme $B_0x_0+C_0x_1$ avec $B_0C_0\not =0,$
n\'ecessairement
\begin{align*}
& f_0=c_0x_2, && f_2=c_2x_2 && \text{et} && \psi_1= B_1x_0+C_1x_1;
\end{align*}
la transformation $A\sigma$ est alors lin\'eaire ce qui est
impossible. Si $\psi_0$ est du type $B_0x_0,$ alors 
\begin{align*}
&f_0=b_0x_1+c_0x_2 && \text{et} &&f_2=b_2x_1+c_2x_2.
\end{align*}
Lorsque $b_0c_0\not=0,$ on constate que
$\psi_1=B_1 x_0$ et $A\sigma$ est lin\'eaire. Si $b_0$ est nul,
$\psi_1=B_1x_0+C_1x_1$ et
\begin{align*}
&\deg A\sigma\leq 2, && A=(x_1: \beta x_0+\gamma x_1:\delta x_1+\varepsilon x_2).
\end{align*}
De m\^eme quand $c_0=0$ on obtient que $A\sigma$ est de degr\'e
inf\'erieur ou \'egal \`a $2$ et $A$ est du type
\begin{align*}
(x_2:\beta x_0+\gamma x_2:\delta x_1+\varepsilon x_2).
\end{align*}

En examinant tous les cas possibles on a le r\'esultat
annonc\'e.
\end{proof}

L'\'enonc\'e pr\'ec\'edent conduit avec les m\^emes notations \`a la:

\begin{pro}
{\sl Soit $Q$ une transformation de $\Sigma^3.$ Si $\deg Q^2\leq 2,$
alors $Q$ est, \`a conjugaison lin\'eaire dynamique pr\`es, de l'un
des types suivants
\begin{align*}
&
\mathcal{Q}_1=\left(\frac{ax_1+1}{x_1},\frac{cx_2+1}{x_2}\right), &&
\mathcal{Q}_2=\left(a+bx_1,\frac{cx_2+x_1}{x_2}\right),
\\
&
\mathcal{Q}_3=\left(\frac{ax_2+1}{x_2},\frac{cx_1+d}{x_1}\right),
&& \mathcal{Q}_4=\left(a+x_2,\frac{cx_1+x_2}{x_1} \right),
\end{align*}
o\`u $a$, $b$, $c$ et $d$ sont des nombres complexes
satisfaisant $bd\not=0.$}
\end{pro}

\begin{rems}
\begin{itemize}
\item Les transformations $\mathcal{Q}_1$ et $\mathcal{Q}_3$ sont
des automorphismes de $\mathbb{P}^1(\mathbb{C})\times~\mathbb{P}^1
(\mathbb{C});$ tous leurs it\'er\'es sont de degr\'e $2.$

\item  On obtient des transformations $f$
de $\Sigma^3$ dont le carr\'e est dans $\Sigma^2;$ par exemple
$\mathcal{Q}_1$ lorsque~$a$ ou $c$ est nul.

\item Les transformations \'etudi\'ees par \textsc{Bedford} et
\textsc{Kim} (\cite{BK}) sont \`a conjugaison dynamique pr\`es de
la forme $\mathcal{Q}_4.$
\end{itemize}
\end{rems}

En g\'en\'eral $\mathcal{Q}_2$ et $\mathcal{Q}_4$ ne sont pas de cube quadratique.

\begin{pro}
{\sl La transformation
$\mathcal{Q}_2=\left(a+bx_1,\frac{cx_2+x_1}{x_2}\right)$ est de cube
quadratique si et seulement si nous sommes dans l'un des deux
cas suivants
\begin{itemize}
\item $c=0;$

\item $c\not=0,$ $b=-1$ et $c^2+a=0.$
\end{itemize}
Dans les deux cas les it\'er\'es de $\mathcal{Q}_2$ sont tous
quadratiques.}
\end{pro}

\bigskip

\begin{proof}[\sl D\'emonstration]
Un calcul montre que $\mathcal{Q}_2^3$ s'\'ecrit
\begin{align*}
\left(a+ab+ab^2+b^3x_1,\frac{c(c^2+2a+ab)x_2+(c^2+ab+a)x_1+bc(1+
b)x_1x_2+b^2x_1^2}{(c^2+a)x_2+cx_1+bx_1x_2}\right).
\end{align*}
Le d\'enominateur de la
seconde composante ne peut \^etre lin\'eaire ($b\not=0$) donc
$\mathcal{Q}_2^3$ est quadratique si et seulement si le
d\'enominateur de la seconde composante de $\mathcal{Q}_2^3$
divise le num\'erateur de celle-ci; cette simplification a lieu si
$(c^2+a)x_2+cx_1+bx_1x_2$ est un produit. C'est le cas si et
seulement si $bc(c^2+a)=0$ et comme $b$ est n\'ecessairement
non nul on doit \'etudier les cas~$c=0$ et $c^2+a=0.$

Si $c=0$ on a
\begin{align*}
\mathcal{Q}_2^3=\left(a+ab+ab^2+b^3x_1,\frac{ab+a+bx_1^2}{bx_2}
\right);
\end{align*}
on remarque que $\mathcal{Q}_2^3$ est \`a conjugaison dynamique
pr\`es du m\^eme type que $\mathcal{Q}_2.$ Ainsi tous les it\'er\'es
de $\mathcal{Q}_2$ sont de degr\'e $2.$

Si $c\not=0$ et $c^2+a=0$ on a
\begin{align*}
\mathcal{Q}_2^3=\left(a+ab+ab^2+b^3x_1,\frac{ac(1+b)x_2+abx_1+bc(1+b)
x_1x_2+b^2x_1^2}{x_1(c+bx_2)}\right)
\end{align*}
Comme $b$ est non nul on constate que $\mathcal{Q}_2^3$ est
quadratique si et seulement si $a(b+1)=0,$ ce qui conduit \`a
$b=-1$ puisque $a\not=0.$ On remarque
qu'alors $\mathcal{Q}_2^4=~\mathrm{id}.$
\end{proof}

Par un calcul \'el\'ementaire on d\'emontre la:

\begin{pro}
{\sl La transformation
$\mathcal{Q}_4=\left(a+x_2,\frac{cx_1+dx_2}{x_1}\right)$ est de
cube quadratique si et seulement si $a$ et $c$ sont nuls; dans ce
cas $\mathcal{Q}_4$ est dynamiquement conjugu\'e \`a $\left(x_2,
\frac{x_2}{x_1}\right)$ qui est p\'eriodique de p\'eriode $6.$
Tous ses it\'er\'es sont quadratiques.}
\end{pro}

Les calculs qui suivent vont nous permettre de distinguer les diff\'erentes
classes de conjugaison des transformations quadratiques de carr\'e
quadratique. Soient $Q$ dans $\mathrm{Bir}_2$ et $m$ un point du plan
projectif complexe; on notera
$\mu_Q(m)$ la multiplicit\'e en $m$ du feuilletage d\'efini par
$Q.$ Nous allons montrer comment on calcule les diff\'erentes
multiplicit\'es par exemple pour les \'el\'ements
\begin{align*}
f_{\alpha,\beta}=\left(\frac{\alpha x_0+x_1}{x_0+1},\beta x_1\right);
\end{align*}
ces transformations comptent trois points d'ind\'etermination
\begin{align*}
&(1:0:0),&&(0:1:0),&&(-1:\alpha:1)
\end{align*}
et deux points fixes $(0:0:1)$ et $
(\alpha-1:0:1).$ Le feuilletage associ\'e \`a $f_{\alpha,\beta}$
est donn\'e par la $1$-forme

\begin{small}
\begin{align*}
(1-\beta)x_1x_2(x_0+x_2)\mathrm{d}x_0+x_2\left(x_2(\alpha
x_0+x_1)-x_0(x_0+x_2) \right)\mathrm{d}x_1+x_1\left(\beta x_0(x_0+x_2)-x_2(\alpha x_0+x_1)\right)
\mathrm{d}x_2.
\end{align*}
\end{small}

Pla\c{c}ons-nous dans la carte $x_1=1$ pour calculer
$\mu_{f_{\alpha,\beta}}(0:1:0);$ le feuilletage associ\'e \`a
$f_{\alpha,\beta}$ est donn\'e par
\begin{align*}
(1-\beta)x_2(x_0+x_2)\mathrm{d}x_0+\left(\beta x_0(x_0+x_2)-x_2(\alpha x_0+1)\right)\mathrm{d}x_2.
\end{align*}
La quantit\'e $(1-\beta)
x_2(x_0+x_2)$ est nulle si et seulement si $x_2=0$ ou $x_0=-x_2.$ Lorsque
$x_2=0$ (resp. $x_0=-x_2$) on a
\begin{align*}
\beta x_0(x_0+x_2)-x_2(\alpha x_0+1)=\beta x_0^2
\end{align*}
(resp. $\beta x_0(x_0+x_2)-x_2(\alpha x_0+1)
=-(\alpha x_0+1)x_0$); il s'en suit que $\mu_{f_{\alpha,\beta}}(0:1:0)=2+1
=3.$ De m\^eme on obtient
\begin{align*}
\mu_{f_{\alpha,\beta}}(1:0:0)=\mu_{f_{\alpha,
\beta}}(0:0:1)=\mu_{f_{\alpha,\beta}}(\alpha-1:0:1)=\mu_{f_{
\alpha,\beta}}(1:-\alpha:1)=1.
\end{align*}

On a la:

\begin{pro}
{\sl Soit $Q$ un \'el\'ement de $\Sigma^3$ dont le carr\'e est de
degr\'e inf\'erieur ou \'egal \`a~$2.$ On suppose que
$\deg\mathcal{F}(Q)=2.$ Notons $p_1,$ $p_2$ et
$p_3$ les points d'ind\'etermination de $Q.$ Le tableau suivant
donne les multiplicit\'es aux points d'ind\'etermination.\bigskip
\begin{small}
\begin{center}
\begin{tabular}{|*{2}{c|}l r|}
 \hline
     & \\
     \begin{small}valeurs de  $(\mu_Q(p_1),\mu_Q(p_2),\mu_Q(p_3))$\end{small} &
      \hspace*{8mm}\begin{small} $Q$ est dynamiquement\end{small}\hspace*{8mm}\\
      \begin{small}\`a permutation pr\`es\end{small} &   \begin{small}conjugu\'ee \`a
       \end{small} \\
     & \\
   \hline
     & \\
     $(1,\hspace{1mm} 1,\hspace{1mm} 1)$ &  $\left(a+\frac{1}{x_1},\hspace{1mm} c+\frac{1}{x_2}\right)$\\
     & \\
     $(1,\hspace{1mm} 1,\hspace{1mm} 3)$ &  $\left(1+bx_1,\hspace{1mm} c+\frac{x_1}{x_2}\right),$
     $b\not\in\{0,\hspace{1mm}1\}$\\
     & \\
     $(1,\hspace{1mm} 5,\hspace{1mm} 1)$ & $\left(1+x_1,\hspace{1mm} c+\frac{x_1}{x_2}\right)$ \\
     & \\
     $(2,\hspace{1mm} 3,\hspace{1mm} 1)$ & $\left(bx_1,\hspace{1mm} 1+\frac{x_1}{x_2}\right),$ $b\not\in\{0,\hspace{1mm} 1\}$\\
     & \\
     $(2,\hspace{1mm} 3,\hspace{1mm} 1)$ & $\left(\frac{b}{x_2},\hspace{1mm} 1+\frac{1}{x_1}\right),$ $b\not\in\{0,\hspace{1mm} 1\}$ \\
     & \\
     $(3,\hspace{1mm} 3,\hspace{1mm} 1)$ &  $\left(bx_1,\hspace{1mm}\frac{x_1}{x_2}\right),$ $b\not\in\{0,\hspace{1mm} 1\}$\\
     & \\
     $(1,\hspace{1mm} 2,\hspace{1mm} 2)$ & $\left(\frac{x_2+b}{x_2},\hspace{1mm} 1+\frac{d}{x_1}\right),$ $bd\not=0$ \\
     & \\
     $(1,\hspace{1mm} 2,\hspace{1mm} 2)$ & $\left(1+x_2,\hspace{1mm} c+\frac{dx_2}{x_1}\right),$ $d\not=0$\\
     & \\
     $(1,\hspace{1mm} 2,\hspace{1mm} 4)$ & $\left(\frac{1}{x_2},\hspace{1mm} 1+\frac{1}{x_1}\right)$ \\
     & \\
     $(2,\hspace{1mm} 2,\hspace{1mm} 2)$ & $\left(x_2,\hspace{1mm} c+\frac{x_2}{x_1}\right),$ $c\not=-1$\\
     & \\
     $(3,\hspace{1mm} 2,\hspace{1mm} 2)$ & $\left(x_2,\hspace{1mm}\frac{x_2}{x_1}-1\right)$ \\
     & \\
  \hline
\end{tabular}
\end{center}
\end{small}}
\end{pro}

\begin{rem}
La condition $\deg\mathcal{F}(Q)=2$ impose que les param\`etres $a,$ $b,$ $c$ et $d$
ci-dessus prennent des valeurs g\'en\'eriques. Lorsque les triplets de
multiplicit\'e de deux transformations sont distincts elles ne peuvent \^etre
dynamiquement conjugu\'ees.
\end{rem}

\begin{rems}
\begin{itemize}
\item On constate que les deux mod\`eles $\left(bx_1,\hspace{1mm}
1+\frac{x_1}{x_2} \right)$ et $\left(\frac{b}{x_2}, \hspace{1mm}
1+\frac{1}{x_1}\right)$ ont, \`a permutation pr\`es, m\^eme
multiplicit\'e aux points d'ind\'etermination mais ne sont pas
dynamiquement conjugu\'es; en effet $\left(bx_1,\hspace{1mm}
1+\frac{x_1} {x_2}\right)$ laisse la fibration $x_1=$ cte
invariante alors qu'un calcul montre que $\left(\frac{b} {x_2},\hspace{1mm}
1+\frac{1}{x_1}\right)$ ne pr\'eserve pas de fibration en droites.

\item On note que $\left(\frac{x_2+b}{x_2}, \hspace{1mm}
1+\frac{d}{x_1} \right)^2$ pr\'eserve deux fibrations en droites
alors que $\left(1+x_2,\hspace{1mm} c+\frac{dx_2}{x_1}\right)^2$
n'en pr\'eserve qu'une; par suite $\left(\frac{x_2+b}{x_2},
\hspace{1mm}1+\frac{d}{x_1}\right)$ et $\left(1+x_2,\hspace{1mm}
c+\frac{dx_2}{x_1}\right)$ ne sont pas dynamiquement conjugu\'es.

\item Chaque fois que la somme des multiplicit\'es
aux points d'ind\'etermination vaut $7$ la transformation correspondante n'a
pas de point fixe.
\end{itemize}
\end{rems}

\subsection{Cas non g\'en\'eriques}\hspace{1mm}

Pour classifier les \'el\'ements de $\Sigma^2$ de
carr\'e quadratique on s'int\'eresse aux \'el\'ements $A$ de~$\mathrm{PGL}_3(
\mathbb{C})$ tels que $\rho A\rho$ appartienne \`a
$\mathrm{Bir}_2.$

\begin{pro}\label{deslee}
{\sl Tout \'el\'ement $A$ de $\mathrm{PGL}_3(\mathbb{C})$ tel que $\rho A\rho$
soit de degr\'e inf\'erieur ou \'egal \`a $2$ est de l'un des
trois types suivants
\begin{align*}
&(\alpha x_0+\beta x_1+\delta x_2:x_1:\varepsilon x_1+\gamma x_2), && (\alpha
x_1+\beta x_2:\gamma x_0+\delta x_2:x_2), \\
& (\alpha x_0+\beta x_2:\gamma x_1+\delta x_2:x_2), &&
\beta,\hspace{1mm}\delta,\hspace{1mm}
\varepsilon\in\mathbb{C},\hspace{1mm}\alpha,\hspace{1mm}\gamma\in\mathbb{C}^*.
\end{align*}}
\end{pro}

\begin{proof}[\sl D\'emonstration]
On reprend la d\'emarche qui nous a permis d'\'etablir la Proposition
\ref{sigqua}.

\'Ecrivons $A\rho$ sous la forme $(F_0:F_1:F_2).$ La
transformation $\rho A\rho$ est dans $\mathrm{Bir}_2$ si et seulement s'il
existe $\psi$ et $q_i$ dans $\mathbb{C}[x_0,x_1,x_2]_2$ satisfaisant
\begin{equation}\label{bruit}
(F_0F_1:F_2^2:F_1F_2)=\psi(q_0:q_1:q_2).
\end{equation}

Si $\psi$ \'etait irr\'eductible, $\psi$ co\"{\i}nciderait, \`a
multiplication par un scalaire pr\`es, avec deux $F_i$ distincts;
$A\rho$ ne serait alors pas inversible ce qui est absurde. Il en
r\'esulte que $\psi$ s'\'ecrit comme le produit de deux formes
lin\'eaires que nous noterons $\psi_0$ et $\psi_1.$

Si $\psi_0=\alpha\psi_1$ avec $\alpha$ dans
$\mathbb{C}^*$ on a l'alternative suivante
\begin{itemize}
\item $q_1=\ell_1^2$ et $F_2=\ell_1\psi_0$ o\`u $\ell_1$ d\'esigne
une forme lin\'eaire non proportionnelle \`a~$\psi_0;$

\item $q_1=\psi_0^2$ et $F_2=\psi_0^2.$
\end{itemize}

Nous allons consid\'erer ces deux \'eventualit\'es au
cas par cas.
\begin{itemize}
\item L'\'egalit\'e (\ref{bruit}) entra\^ine que
\begin{align*}
& F_1=\psi_0^2 && \text{et} && F_2=\ell_1\psi_0.
\end{align*}

Le fait que chaque $F_i$ soit du type $*x_0x_1+*x_2^2+*x_1x_2$
implique que
\begin{align*}
& F_0=\alpha x_0x_1+\beta x_2^2+\delta x_1x_2, && F_1=x_2^2, && F_2=
\varepsilon x_2^2+\gamma x_1x_2,
\end{align*}
autrement dit
\begin{align*}
& A=\left[
\begin{array}{ccc}
\alpha & \beta & \delta \\
0 & 1 & 0 \\
0 & \varepsilon & \gamma
\end{array}
\right], && \alpha\gamma\not=0.
\end{align*}

\item L'\'egalit\'e (\ref{bruit}) conduit \`a $F_0q_2=\psi_0^2q_0.$ On sait que
$\psi_0^2$ ne divise pas $F_0$ (si c'\'etait le cas $F_0$ et $F_2$ seraient
\'egales \`a multiplication par un scalaire pr\`es et $A\rho$ ne serait pas
inversible !) Si $\psi_0$ ne divisait pas $F_0,$ on aurait $q_0=\psi_0^2$
et $F_1,$ $F_2$ seraient multiples l'une de l'autre ce qui est impossible.
Ainsi $F_0$ est du type $\psi_0\ell$ avec $\ell$ forme lin\'eaire non
proportionnelle \`a~$\psi_0;$ en r\'e\'ecrivant (\ref{bruit}) on obtient $\ell q_2=
\psi_0q_0$ d'o\`u: $\psi_0$ divise $q_2.$ Par suite $\psi_0$ divise chaque
$F_i$ ce qui est impossible.
\end{itemize}

Finalement \'etudions le cas o\`u les $\psi_i$ ne sont pas multiples
l'une de l'autre. Il s'en suit que~$q_1=F_2=\psi_0\psi_1.$ \`A partir de
(\ref{bruit}) on a
\begin{align*}
& F_0q_2=\psi_0\psi_1q_0, && F_1=q_2.
\end{align*}
Puisque les $F_i$ ne sont pas proportionnels on a \`a
r\'eindexation pr\`es
\begin{align*}
& F_0=\psi_0\ell, && F_1=\psi_1\widetilde{\ell}, && F_2=\psi_0\psi_1
\end{align*}
o\`u $\ell,$ $\widetilde{\ell}$ d\'esignent deux formes lin\'eaires telles que
$\psi_1$ et $\ell$ (resp. $\psi_0$ et $\widetilde{\ell}$) ne soient pas
proportionnelles.

Les $F_i$ \'etant de la forme $*x_0x_1+*x_2^2+*x_1x_2$ un calcul
montre que $A$ est de l'un des types suivants
\begin{align*}
& \left[
\begin{array}{ccc}
 0 & \gamma & \beta \\
\alpha & 0 & \delta \\
0 & 0 & 1
\end{array}
\right],&& \left[
\begin{array}{ccc}
 \alpha & 0 & \beta \\
0 &\gamma & \delta \\
0 &0 & 1
\end{array}
\right],&&\alpha\gamma\not=0.
\end{align*}
\end{proof}

Cette proposition technique permet d'obtenir le:

\begin{cor}
{\sl Soit $Q$ une transformation de $\Sigma^2.$ Si $\deg Q^2\leq 2,$ alors
$Q$ est, \`a conjugaison lin\'eaire dynamique pr\`es, de l'un des types
suivants
\begin{align*}
& \mathscr{Q}_1=\left(\frac{\alpha x_0x_1+\gamma x_1+\beta}{1+\delta
x_1},\frac{1}{1+\delta x_1}\right), && \mathscr{Q}_2=
\left(\frac{\alpha x_0x_1+x_1+\beta}{x_1},\frac{1}{x_1}\right),\\
& \mathscr{Q}_3=\left(\frac{\alpha}{x_1}+1,\gamma x_0+1\right),
&&
\mathscr{Q}_4=\left(\frac{\alpha}{x_1},\gamma x_0\right),\\
& \mathscr{Q}_5=\left(\frac{\alpha}{x_1}+1,x_0\right), &&
\mathscr{Q}_6=\left(\frac{\alpha}{x_1},x_0+1\right),\\
&
\mathscr{Q}_7=\left(\alpha x_0+1,\frac{\gamma}{x_1}+1\right),
&& \mathscr{Q}_8=\left(\alpha x_0,\frac{\gamma}{x_1}+1\right), \\
& \mathscr{Q}_9=\left(\alpha x_0+1,\frac{1}{x_1}\right),
&&
\mathscr{Q}_{10}=\left(\alpha x_0,\frac{1}{x_1}\right).
\end{align*}}
\end{cor}

\begin{rem}
On remarque que $\deg\mathscr{Q}_i^n\leq 2$ pour tout $n.$
\end{rem}

\bigskip

Toute transformation appartenant \`a $\Sigma^1$
s'\'ecrit $B\tau C;$ elle est de carr\'e quadratique si et
seulement si $\tau CB\tau$ l'est. On s'int\'eresse donc aux
\'el\'ements $A$ de $\mathrm{PGL}_3(\mathbb{C})$ tels que $\tau A\tau$ appartienne
\`a $\mathrm{Bir}_2.$

En suivant la d\'emarche utilis\'ee pour la Proposition
\ref{sigqua} on obtient la:

\begin{pro}
{\sl Tout automorphisme $A$ de $\mathbb{P}^2(
\mathbb{C})$ tel que $\tau A\tau$
soit de degr\'e inf\'erieur ou \'egal \`a $2$ est du type
\begin{align*}
&(x_0:\alpha x_0+\beta x_1:\gamma x_0+\delta x_1+\varepsilon
x_2),&& \alpha,\hspace{1mm} \gamma,\hspace{1mm} \delta\in
\mathbb{C},\hspace{1mm} \beta,\hspace{1mm}\varepsilon\in\mathbb{C}^*.
\end{align*}}
\end{pro}

\begin{proof}[\sl D\'emonstration]
On reprend l'id\'ee de la d\'emonstration de la Proposition \ref{sigqua}.

Notons $F_0,$ $F_1$ et $F_2$ les composantes de $A\tau.$
La transformation $\tau A\tau$ est dans $\mathrm{Bir}_2$ si et seulement si
\begin{equation}\label{bruitb}
(F_0^2:F_0F_1:F_1^2-F_0F_2)=\psi(q_0:q_1:q_2), \hspace{2cm}
q_i,\hspace{1mm}\psi\in\mathbb{C}[x_0,x_1,x_2]_2.
\end{equation}
Si $\psi$ \'etait inversible, $\psi$ diviserait $F_0$ et $F_1$ autrement dit
$A\tau$ ne serait pas inversible ce qui est impossible; $\psi$ s'\'ecrit
donc $\psi_0\psi_1$ avec $\psi_i$ forme lin\'eaire.

Supposons que $\psi_0$ et $\psi_1$ soient multiples l'une de
l'autre; (\ref{bruitb}) se r\'e\'ecrit alors
\begin{align*}
& F_0^2=\psi_0^2q_0, && F_0F_1=\psi_0^2q_1, && F_1^2-F_0F_2=
\psi_0^2q_2
\end{align*}
d'o\`u l'alternative
\begin{itemize}
\item il existe une forme lin\'eaire $\ell$ non proportionnelle \`a $\psi_0$
telle que $q_0=\ell^2$ et $F_0=\ell\psi_0;$

\item $q_0=F_0=\psi_0^2.$
\end{itemize}

Traitons chacune de ces \'eventualit\'es.

\begin{itemize}
\item Commen\c{c}ons par supposer que $q_0=\ell^2$ et $F_0=\ell\psi_0.$
L'\'egalit\'e (\ref{bruitb}) conduit \`a
\begin{align*}
& \ell F_1=\psi_0q_1, && F_1^2=\psi_0^2q_2+\ell\psi_0F_2.
\end{align*}
La premi\`ere \'egalit\'e entra\^ine, puisque $\ell$ et $\psi_0$ ne sont pas
proportionnelles, que $\psi_0$ divise~$F_1;$ par suite $\psi_0^2$ divise
$F_1^2.$ La seconde \'egalit\'e assure alors que $\psi_0$ divise $F_2$
ce qui est impossible.

\item Finalement supposons que $F_0=q_0=\psi_0^2.$ \`A partir de
(\ref{bruitb}) on a
\begin{align*}
& F_1=q_1, && F_1^2=\psi_0^2(q_2+F_2);
\end{align*}
il en r\'esulte que $\psi_0^2$ divise $F_1^2.$ Comme $F_0$ et $F_1$
ne sont pas proportionnelles il existe une forme lin\'eaire $\ell$ non
proportionnelle \`a $\psi_0$ telle que $F_1=\ell\psi_0.$ Chaque
$F_i$ \'etant du type $*x_0^2+*x_0x_1+*(x_1^2-x_0x_2)$ on obtient
\begin{align*}
& F_0=x_0^2, && F_1=x_0(\beta x_0+\alpha x_1), &&\alpha\not=0;
\end{align*}
autrement dit
\begin{align*}
&A=\left[
\begin{array}{ccc}
 1 & 0 & 0 \\
\beta &\alpha & 0 \\
\delta & \varepsilon & \gamma
\end{array}
\right],&&\alpha\gamma\not=0.
\end{align*}
\end{itemize}

Pour finir supposons que les formes lin\'eaires
$\psi_0$ et $\psi_1$ ne soient pas proportionnelles. On en
d\'eduit que $F_0$ s'\'ecrit $\psi_0\psi_1;$ de (\ref{bruitb})
on tire
\begin{align*}
& F_1=q_1, && q_1^2=\psi_0\psi_1(q_2+F_2).
\end{align*}
Par suite $\psi_0\psi_1$ divise $q_1^2;$ puisque $\psi_0$
et $\psi_1$ ne sont pas multiples l'une de l'autre, $q_1$ est
de la forme $\alpha\psi_0\psi_1.$ En particulier $F_0$ et
$F_1$ sont proportionnelles ce qui est impossible.
\end{proof}

\begin{cor}
{\sl Soit $Q$ une transformation de $\Sigma^1.$ Si $Q^2$ est de
degr\'e inf\'erieur ou \'egal \`a $2,$ alors $Q$ est, \`a
conjugaison lin\'eaire dynamique pr\`es, de la forme
\begin{align*}
&(\alpha+ \beta x_0,\gamma+\delta x_0+x_0^2-\varepsilon x_1),
&&\alpha, \hspace{1mm} \gamma,\hspace{1mm}\delta\in\mathbb{C},\hspace{1mm}
\beta,\hspace{1mm}\varepsilon\in\mathbb{C}^*.
\end{align*}}
\end{cor}

N'importe quel it\'er\'e d'une telle transformation est toujours
quadratique.

On a donc le:

\begin{cor}
{\sl Soit $q$ une transformation de $\Sigma^2$ ou $\Sigma^1;$ d\`es
que le carr\'e de $q$ est quadratique, tous les it\'er\'es le sont.}
\end{cor}

Comme on l'a signal\'e ceci n'est pas vrai pour les \'el\'ements
de $\Sigma^3$ (par exemple pour les $f_{\alpha,\beta}$).

\clearemptydoublepage
\chapter{Des propri\'et\'es alg\'ebriques du
groupe de \textsc{Cremona}}\label{algalg}

\section{Le groupe de \textsc{Cremona} ne se
plonge pas dans un $\mathrm{GL}_n(\Bbbk)$}\hspace{1mm}

Bien que $\mathrm{Bir}(\mathbb{P}^2(\mathbb{C}))$ poss\`ede de nombreuses propri\'et\'es des groupes
lin\'eaires on a la:

\begin{pro}\label{nonlin}
{\sl Le groupe de \textsc{Cremona} ne se plonge pas dans $\mathrm{GL}_n
(\Bbbk)$ o\`u $\Bbbk$ d\'esigne un corps de caract\'eristique nulle.}
\end{pro}

Avant de d\'emontrer ce r\'esultat rappelons l'\'enonc\'e suivant
d\^u \`a \textsc{Birkhoff}.

\begin{lem}[\cite{Bi}]
{\sl Soient $\Bbbk$ un corps de caract\'eristique nulle et $A,$ $B,$ $C$ trois
\'el\'e\-ments de~$\mathrm{GL}_n(\Bbbk)$ satisfaisant $[A,B]=C,$
$[A,C]=[B,C]=\mathrm{id}$ et $C$ d'ordre $p$ premier~; alors $p\leq n.$}
\end{lem}

\begin{proof}[\sl D\'emonstration de la Proposition \ref{nonlin}]
Supposons qu'il existe un morphisme injectif $\varsigma$ du groupe
de \textsc{Cremona} dans
$\mathrm{GL}_n(\Bbbk).$ Pour tout nombre premier $p$
consid\'erons, dans la carte affine $x_2=1,$ le groupe
\begin{align*}
\langle\left(\exp\left(-\frac{2\mathrm{i}\pi}{p}\right)x_0,x_1\right),(x_0,x_0x_1),\left(x_0,\exp\left(\frac{2\mathrm{i}\pi}{p}\right)x_1\right)\rangle.
\end{align*}
Les images par $\varsigma$ des trois g\'en\'erateurs satisfont le
lemme de \textsc{Birkhoff} donc $p\leq n~;$ ceci \'etant valable
pour tout premier $p,$ nous obtenons le r\'esultat annonc\'e.
\end{proof}

Toutefois le groupe de \textsc{Cremona} poss\`ede de \og
nombreux sous-groupes lin\'eaires\fg. Il contient, outre les
groupes $\mathrm{PGL}_3( \mathbb{C})$ et
$\mathrm{PGL}_2(\mathbb{C})\times\mathrm{PGL}_2(\mathbb{C}),$ le sous-groupe du
groupe de \textsc{de Jonqui\`eres} pr\'eservant la fibration
$x_1=$ cte fibre \`a fibre qui s'identifie \`a
$\mathrm{PGL}_2(\mathbb{C}(x_1)).$

\section{Centralisateur d'une transformation $A\sigma$ g\'en\'erique}\hspace{1mm}

Dans \cite{Ca} \textsc{Cantat} montre l'\'enonc\'e
suivant:
\begin{thm}[\cite{Ca}, th\'eor\`eme B]
{\sl Soit $f$ une transformation birationnelle d'une surface complexe
compacte $S$ dont le premier degr\'e dynamique est strictement
plus grand que $1.$ Si $g$ est une transformation birationnelle de
$S$ qui commute avec $f$ il existe $m$ dans $\mathbb{N}^*$ et $n$
dans $\mathbb{Z}$ tels que $g^m=f^n.$}
\end{thm}

Nous allons d\'emontrer, seulement pour les
transformations birationnelles du type $A\sigma,$ avec~$A$
g\'en\'erique, un r\'esultat un peu plus pr\'ecis.

\begin{pro}\label{alcides}
{\sl Soit $A$ un \'el\'ement de $\mathrm{PGL}_3(\mathbb{C})$ dont les coefficients
sont rationnellement ind\'ependants; le centralisateur de
$A\sigma$ dans $\mathrm{Bir}(\mathbb{P}^2(\mathbb{C}))$ s'identifie \`a $\mathbb{Z}.$}
\end{pro}

Pour ce faire commen\c{c}ons par \'etablir le:

\begin{lem}\label{abab}
{\sl Soit $A$ un \'el\'ement de $\mathrm{PGL}_3(\mathbb{C})$ dont les coefficients
sont $\mathbb{Q}$-alg\'ebriquement ind\'ependants; le centralisateur de
$A\sigma$ dans $\mathrm{Bir}(\mathbb{P}^2(
\mathbb{C}))$ est ab\'elien.}
\end{lem}

\begin{proof}[\sl D\'emonstration]
Notons $m_1,$ $m_2,$ $m_3$ et $p$ les points fixes de $A\sigma.$
Puisque $A$ est g\'en\'erique, les valeurs propres de la partie
lin\'eaire de $A\sigma$ en les $m_i$ sont \og ind\'ependantes\fg ,
celles au point~$p$ d\'ependant des six autres. D\'esignons par
$\mathrm{Cent}(A\sigma)$ le centralisateur de $A \sigma$ dans
$\mathrm{Bir}(\mathbb{P}^2(\mathbb{C}))$
\begin{align*}
\mathrm{Cent}(A\sigma)=\{f\in\mathrm{Bir}(\mathbb{P}^2(\mathbb{C}))\hspace{1mm}|
\hspace{1mm} fA\sigma=A\sigma f\}.
\end{align*}

Soit $f$ un \'el\'ement de $\mathrm{Cent}(A\sigma).$

Montrons que n\'ecessairement $f$ est holomorphe en tout
point de $\textrm{Fix}(A\sigma).$ Supposons par l'absurde que $f$
ne soit pas holomorphe, par exemple, en $m_1;$ alors $m_1$ est un
point d'ind\'etermination de $f$ dont nous noterons $\Gamma$
l'image. Puisque $A\sigma$ est holomorphe en $m_1$ on constate que
$\Gamma$ est invariante par $A\sigma.$ Or d'apr\`es le Th\'eor\`eme
\ref{courbinv2} la transformation $A\sigma$ n'admet pas de courbe
invariante; $f$ est donc holomorphe en $m_1.$

Puisque $f$ est holomorphe en tout point de $\textrm{Fix}
(A\sigma),$ elle permute les points de $\textrm{Fix}
(A\sigma).$ \og L'ind\'ependance\fg\hspace{1mm} des couples de valeurs
propres de la partie lin\'eaire de $A\sigma$ en les $m_i$ assure
que $f$ fixe au moins deux des $m_i,$ le troisi\`eme pouvant
\^etre a priori permut\'e avec $p$ (on peut penser que cela
n'arrive pas g\'en\'eriquement). Supposons, \`a r\'eindexation
pr\`es, que $m_1$ soit fix\'e par $f.$ Par g\'en\'ericit\'e il
existe un germe de diff\'eomorphisme $\varphi$ d\'efini au voisinage
de $m_1$ tel que
\begin{align*}
\varphi A\sigma_{,m_1}\varphi^{-1}= (\delta x_0,\eta x_1)
\end{align*}
o\`u $\delta,$ $\eta$ sont les valeurs propres de la diff\'erentielle
de $A\sigma$ en $m_1.$ Un calcul \'el\'ementaire montre, puisque
$\delta$ et $\eta$ sont non r\'esonnants, que $\varphi f_{,m_1}
\varphi^{-1}$ est lin\'eaire de la forme $(ax_0,bx_1).$ Le centralisateur
de $A\sigma$ est donc ab\'elien.
\end{proof}

\begin{lem}\label{orbpos}
{\sl Soient $A$ un automorphisme de $\mathbb{P}^2(\mathbb{C})$
g\'en\'erique et $p_1,$ $p_2,$ $p_3$ les points d'ind\'etermination
de $A\sigma.$ Alors pour tout $n\geq 1$ et tout $i$ dans $\{1,2,3\}$
l'image de $p_i$ par $(A\sigma)^n$ est une courbe.}
\end{lem}

\begin{proof}[\sl D\'emonstration]
Supposons que ce ne soit pas le cas, par exemple pour le point
$p=(1:0:~0).$ Il existe alors un entier positif $\ell$ que l'on choisit minimal
tel que $(A\sigma)^\ell(p)$ soit de nouveau un point not\'e $q.$ Alors
pour tout automorphisme de corps $\kappa$ il en est de m\^eme
pour~$A^\kappa\sigma$
\begin{align*}
& (A^\kappa\sigma)^\ell(1:0:0)=q^\kappa, && \forall\hspace{1mm}
\kappa\in\mathrm{Aut}(\mathbb{C},+,.).
\end{align*}
Par densit\'e il en est de m\^eme pour tout $B\sigma$
\begin{align*}
& (B\sigma)^\ell(1:0:0)=m_B, && \forall\hspace{1mm}B\in
\mathrm{PGL}_3(\mathbb{C}),
\end{align*}
$m_B$ d\'esignant un point du plan projectif complexe.
Maintenant consid\'erons la transformation
\begin{align*}
f=(x_1x_2+x_0x_2:x_0x_2:x_0x_1)=(x_0+x_1:x_1:x_2)\sigma.
\end{align*}
On constate que l'image de $(1:0:0)$ par $f$ est la droite
d'\'equation $x_1=x_0$ et que l'image de cette droite par
$f^\ell$ est une droite pour tout entier $\ell:$ contradiction.
\end{proof}

\begin{rem}\label{orbpos2}
Le Lemme \ref{orbpos} implique que si $A$ est un automorphisme
de $\mathbb{P}^2(\mathbb{C})$ \`a coefficients alg\'ebriquement
ind\'ependants sur $\mathbb{Q}$ et si $\mathrm{Ind}\hspace{1mm}
A\sigma=\{p_1,\hspace{1mm}p_2,\hspace{1mm}p_3\}$
alors pour tout $n\geq 1$ et tout $i$ dans $\{1,2,3\}$
l'image de $p_i$ par $(A\sigma)^n$ est une courbe.
\end{rem}

\begin{proof}[\sl D\'emonstration de la Proposition \ref{alcides}]
On reprend une id\'ee de \cite{De2}.

Soit $g$ un \'el\'ement du groupe de \textsc{Cremona}
qui commute \`a $A\sigma.$ Notons $p_i$ les points
d'ind\'etermination de $A\sigma.$ D'apr\`es la Remarque
\ref{orbpos2} l'orbite positive de $p_i$ sous l'action de
$A\sigma$ est constitu\'ee de courbes. Ou bien $p_i$ est
d'ind\'etermination pour $g;$ sinon puisque
$\mathrm{Exc}\hspace{1mm} g$ est fini, il existe un entier positif
$k$ (choisi minimal) tel que $(A\sigma)^k(p_i)$ ne soit pas
contenu dans $\mathrm{Exc}\hspace{1mm} g.$ Quitte \`a remplacer
$g$ par $\widetilde{g}:=g(A \sigma)^k$ on constate que
$\widetilde{g}(p_i)$ est un point d'ind\'etermination de~$A\sigma;$ on a: $\widetilde{g}^{\hspace{1mm}3}(p_i)=p_i$ pour tout~$i.$

Comme l'orbite n\'egative de $p_i$ par
$A\sigma$ est \textsc{Zariski} dense (Th\'eor\`eme \ref{orbneg}), $\widetilde{g}^{\hspace{1mm}3}$
co\"{\i}ncide avec l'identit\'e et $g^3$ est une puissance de
$A\sigma.$ Le centralisateur
de $A\sigma$ dans $\mathrm{Bir}(\mathbb{P}^2(\mathbb{C}))$ s'identifie donc \`a
$\mathbb{Z}\times\mathbb{Z}/p\mathbb{Z}$ avec $p\leq 3.$

Supposons que $p$ soit sup\'erieur ou \'egal \`a $2.$
Il existe alors une transformation birationnelle non triviale $\varphi$ qui
commute \`a $A\sigma$ et telle que $\varphi^p=\mathrm{id};$ c'est
donc le cas pour les \'el\'ements de la forme $B\sigma$ avec
$B=A^\kappa,$ $\kappa$ automorphisme du corps $\mathbb{C},$ puis par
densit\'e pour tout \'el\'ement de type~$B\sigma$ avec $B$ dans
$\mathrm{PGL}_3(\mathbb{C})$ (remarquons qu'une famille \`a un
param\`etre $\varphi_s$ de transformations birationnelles
p\'eriodiques telle que $\varphi_0=\mathrm{id}$ est constamment
l'identit\'e).  Le centralisateur de la transformation
\begin{eqnarray}
f_{\alpha,\beta}&=&(x_0+\alpha x_1:\beta(\alpha x_1+x_2):x_1)\sigma(x_0+x_2:
x_1-\alpha x_2:x_2)\nonumber\\
&=&((\alpha x_0+x_1)x_2:\beta x_1(x_0+x_2):x_2(x_0+x_2))\nonumber
\end{eqnarray}
s'identifie, pour $\alpha,$ $\beta$ g\'en\'eriques, \`a $\mathbb{Z}$
(\emph{voir} \cite{De2}); la d\'emonstration est analogue \`a
celle qu'on vient de pr\'esenter pour une transformation $A\sigma$
avec $A$ \`a coefficients alg\'ebriquement ind\'ependants sur
$\mathbb{Q}.$ L'absence du facteur $\mathbb{Z}/
p\mathbb{Z}$ s'explique par le fait qu'on peut distinguer les
orbites des points d'ind\'etermination de $f_{\alpha,\beta}$
\begin{figure}[H]
\begin{center}
\input{orbind.pstex_t}
\end{center}
\end{figure}

\begin{figure}[H]
\begin{center}
\input{orbind2.pstex_t}
\end{center}
\end{figure}

\begin{figure}[H]
\begin{center}
\input{orbind3.pstex_t}
\end{center}
\end{figure}

Le r\'esultat d\'ecoule du fait qu'\`a conjugaison pr\`es
$f_{\alpha,\beta}$ est du type $B\sigma.$
\end{proof}

\begin{rem}
Il y a toutefois des transformations de la forme $A\sigma$ qui
poss\`edent un \og gros\fg\hspace{1mm} centralisateur; les flots
quadratiques g\'en\'eriques d\'ecrits dans le Chapitre
\ref{germgerm} fournissent de tels exemples.
\end{rem}

\section{Construction de sous-groupes libres}\hspace{1mm}

En g\'en\'eral le groupe engendr\'e par un nombre fini de transformations
birationnelles quadratiques est libre.

\begin{pro}
{\sl Soient $A_1,$ $\ldots,$ $A_n$ des automorphismes g\'en\'eriques de
$\mathbb{P}^2(\mathbb{C}).$ Le groupe $\langle A_1,\hspace{1mm}\ldots,
\hspace{1mm} A_n,\hspace{1mm}\sigma \rangle$ engendr\'e par les
$A_i$ et $\sigma$ est le produit libre
\begin{align*}
\stackrel{n}{\overbrace{\mathbb{Z}\ast\ldots\ast\mathbb{Z}}}\ast(\mathbb{Z}/2\mathbb{Z}).
\end{align*}

En particulier le groupe engendr\'e par les
transformations quadratiques $A_1\sigma,$ $\ldots,$ $A_n\sigma$
est libre.}
\end{pro}

\begin{proof}[\sl D\'emonstration]
On va d\'emontrer la Proposition pour $n=1,$ le m\^eme argument
permettant de conclure pour $n$ quelconque.

Par un argument de \textsc{Baire} il suffit de trouver
un exemple satisfaisant l'\'enonc\'e. Consid\'erons une
transformation $A$ du type suivant
\begin{align*}
A=(\alpha x_0+\beta x_1:\gamma x_0+\delta x_1:x_2).
\end{align*}
Comme le pinceau de droites
$x_0=tx_1$ est invariant par les deux transformations $A$ et
$\sigma$ on h\'erite d'une repr\'esentation
\begin{align*}
&\langle A,\hspace{1mm}\sigma\rangle\to\mathrm{PGL}_2(\mathbb{C}), && A
\hspace{1mm}\colon\hspace{1mm} t\mapsto\frac{\alpha
t+\beta}{\gamma t+\delta},&& \sigma\hspace{1mm}\colon\hspace{1mm}
t\mapsto\frac{1}{t}.
\end{align*}
Dans $\mathrm{PGL}_2(\mathbb{C})$ le sous-groupe engendr\'e par les matrices
\begin{align*}
&\left[\begin{array}{cc}
\alpha & \beta\\
\gamma & \delta
\end{array}
\right] && \text{et} &&\left[\begin{array}{cc}
0 & 1\\
1 & 0
\end{array}
\right]
\end{align*}
est g\'en\'eriquement isomorphe \`a $\mathbb{Z}\ast(\mathbb{Z}/2\mathbb{Z});$
c'est le cas par exemple pour $\left[\begin{array}{cc}
1 & \varepsilon\\
0 & 1
\end{array}
\right] $ et $\left[\begin{array}{cc}
0 & 1\\
1 & 0
\end{array}
\right]$ lorsque $\varepsilon$ est g\'en\'erique. Ceci implique que
$\langle (x_0+\varepsilon x_1:x_1:x_2),\hspace{1mm}\sigma\rangle$ est isomorphe
\`a $\mathbb{Z}\ast\mathbb{Z}/2\mathbb{Z}.$
\end{proof}

\begin{rem}
La preuve ci-dessus montre que la dynamique du groupe engendr\'e
par~$\sigma$ et $A=(\alpha x_0+\beta x_1:\gamma x_0+\delta x_1:x_2)$
est plus ou moins triviale, au sens o\`u elle se r\'eduit \`a la dimension $1.$
En effet $A$ et $\sigma$ se rel\`event par l'application de l'\'eclatement de
l'origine dans la carte affine $x_2=1$ en les applications
\begin{align*}
& \widetilde{A}\hspace{1mm}\colon\hspace{1mm}(t,x_1)\mapsto
\left(\frac{\alpha t+\beta} {\gamma t+\delta}, x_1\right), &&
\widetilde{\sigma}\hspace{1mm}\colon\hspace{1mm}(t,x_1)\mapsto\left(
\frac{1}{t},x_1\right).
\end{align*}

Notons que le groupe engendr\'e par $A$ et $\sigma$ est
l'exemple d'un groupe qui se rel\`eve en un groupe
d'automorphismes d'une surface rationnelle.
\end{rem}

\section{Au sujet de la simplicit\'e}

On ne sait pas si le groupe de \textsc{Cremona} est simple.
Soient $\mathrm{G}$ un groupe et $f$ un \'el\'ement de~$\mathrm{G}.$
On notera $\mathrm{N}(f,\mathrm{G})$\label{not35ba} le \textbf{\textit{sous-groupe normal de
$\mathrm{G}$ engendr\'e par $f$}}\label{ind500}
\begin{align*}
\mathrm{N}(f,\mathrm{G})=\langle hfh^{-1},\hspace{1mm}
hf^{-1}h^{-1}\hspace{1mm}|\hspace{1mm} h\in \mathrm{G}\rangle.
\end{align*}

Nous nous proposons de calculer quelques $\mathrm{N}(f,\mathrm{G})$
pour quelques transformations de \textsc{Cremona} $f$ particuli\`eres.

\subsection{Premiers calculs de $\mathrm{N}(f,\mathrm{Bir}(\mathbb{P}^2(\mathbb{C})))$ et cons\'equences}\hspace{1mm}

Le groupe $\mathrm{PGL}_3(\mathbb{C})$ est simple par suite tout
\'el\'ement non trivial $A$ de $\mathrm{PGL}_3(\mathbb{C})$ satisfait
\begin{align*}
\mathrm{N}(A,\mathrm{PGL}_3(\mathbb{C}))= \mathrm{PGL}_3(\mathbb{C}).
\end{align*}
On en d\'eduit la:

\begin{pro}
{\sl Le sous-groupe normal engendr\'e par $\sigma$ dans le groupe de
\textsc{Cremona} est le groupe de \textsc{Cremona} tout entier:
\begin{align*}
\mathrm{N}(\sigma,\mathrm{Bir}(\mathbb{P}^2(\mathbb{C})))=\mathrm{Bir}(\mathbb{P}^2(\mathbb{C})).
\end{align*}}
\end{pro}

\begin{proof}[\sl D\'emonstration]
Soit $f$ une transformation de \textsc{Cremona}. Le Th\'eor\`eme
de \textsc{N\oe ther} permet d'\'ecrire $f$ sous la forme
\begin{align*}
& f=(A_1)\sigma A_2\sigma A_3\ldots A_n(\sigma),&&A_i\in\mathrm{PGL}_3(\mathbb{C}).
\end{align*}
Puisque $\mathrm{PGL}_3(\mathbb{C})$ est simple on a
\begin{align*}
&\mathrm{N}((-x_0,-x_1),\mathrm{PGL}_3(\mathbb{C}))= \mathrm{PGL}_3(\mathbb{C})
\end{align*}
et tout $A_i$ poss\`ede une \'ecriture de la forme
\begin{align*}
& h_1(-x_0,-x_1)h_1^{-1}h_2(-x_0,-x_1)h_2^{-1}\ldots h_n(-x_0,-x_1)h_n^{-1}, &&
h_i\in\mathrm{PGL}_3(\mathbb{C}).
\end{align*}
Or l'involution $(-x_0,-x_1)$ est conjugu\'ee\footnote{\hspace{1mm} via
l'\'el\'ement $\left(\frac{x_0+1}
{x_0-1},\frac{x_1+1}{x_1-1}\right)$ de
$\mathrm{PGL}_2(\mathbb{C})\times\mathrm{PGL}_2(\mathbb{C})$} \`a
$\left(\frac{1}{x_0}, \frac{1}{x_1}\right);$ on en d\'eduit que
tout \'el\'ement du groupe de \textsc{Cremona} s'\'ecrit comme un
produit de conjugu\'es de $\sigma=\left(\frac{1}{x_0},
\frac{1}{x_1}\right).$
\end{proof}

Ceci implique le:

\begin{cor}\label{vv}
{\sl Toute transformation de \textsc{Cremona} s'\'ecrit comme une
compos\'ee d'involutions toutes conjugu\'ees.}
\end{cor}

Nous donnerons plus loin (Chapitre \ref{algalg}, \S \ref{nonoplus})
une version un peu plus forte du Corollaire~\ref{vv}.

En ce qui concerne le groupe normal engendr\'e par un automorphisme
de $\mathbb{P}^2(\mathbb{C})$ nous avons~la:

\begin{pro}\label{autnorm}
{\sl Soit $A$ un automorphisme non trivial de $\mathbb{P}^2(\mathbb{C});$ alors
\begin{align*}
\mathrm{N}(A,\mathrm{Bir}(\mathbb{P}^2(\mathbb{C})))=\mathrm{Bir}(\mathbb{P}^2(\mathbb{C})).
\end{align*}}
\end{pro}

\begin{proof}[\sl D\'emonstration]
Comme $\mathrm{N}(A,\mathrm{PGL}_3(\mathbb{C}))=\mathrm{PGL}_3(\mathbb{C})$ l'involution $(-x_0,-x_1)$
s'\'ecrit comme un produit de conjugu\'es de $A.$ Il en r\'esulte,
puisque $(-x_0,-x_1)$ et $\sigma$ sont conjugu\'es, que
\begin{align*}
& \sigma=h_1Ah_1^{-1}h_2Ah_2^{-1}\ldots h_nAh_n^{-1},&& h_i\in\mathrm{Bir}
(\mathbb{P}^2(\mathbb{C}));
\end{align*}
d'o\`u l'inclusion $\mathrm{N}(\sigma,\mathrm{Bir}(\mathbb{P}^2(\mathbb{C})))\subset
\mathrm{N}(A,\mathrm{Bir}(\mathbb{P}^2(\mathbb{C}))).$ L'\'egalit\'e
\begin{align*}
\mathrm{N}(\sigma,\mathrm{Bir}(\mathbb{P}^2(\mathbb{C})))=\mathrm{Bir}(\mathbb{P}^2(\mathbb{C}))
\end{align*}
permet de conclure.
\end{proof}

Comme $\rho$ et $\tau$ sont birationnellement
conjugu\'ees \`a des involutions de $\mathrm{PGL}_3(\mathbb{C})$ on a une
propri\'et\'e analogue pour ces transformations.

\begin{cor}
{\sl On a
\begin{align*}
& \mathrm{N}(\rho,\mathrm{Bir}(\mathbb{P}^2(\mathbb{C})))=\mathrm{Bir}(\mathbb{P}^2(\mathbb{C}))&&\text{et}
&&\mathrm{N}(\tau,\mathrm{Bir}(\mathbb{P}^2(\mathbb{C})))=\mathrm{Bir}
(\mathbb{P}^2(\mathbb{C})).
\end{align*}}
\end{cor}

\`A partir de la Proposition \ref{autnorm} et de
\begin{align*}
&\mathrm{N}(A,\mathrm{PGL}_3(\mathbb{C}))=\mathrm{PGL}_3(\mathbb{C}), &&
\forall\hspace{1mm} A\in\mathrm{PGL}_3(\mathbb{C})
\end{align*}
on obtient le:

\begin{cor}\label{envrille}
{\sl Toute transformation de \textsc{Cremona} s'\'ecrit comme un
produit de conjugu\'es de la translation $(x_0,x_1+1).$}
\end{cor}

Puisque la translation $(x_0,x_1+1)=[(x_0,3x_1),\left(x_0,\frac{x_1+1}
{2}\right)]$ est un commutateur, le Corollaire~\ref{envrille} entra\^ine le:

\begin{cor}
{\sl Le groupe de \textsc{Cremona} est parfait, {\it i.e.} le groupe
d\'eriv\'e de $\mathrm{Bir}(\mathbb{P}^2(\mathbb{C}))$ est le groupe
$\mathrm{Bir}(\mathbb{P}^2(\mathbb{C}))$ entier
\begin{align*}
[\mathrm{Bir}(\mathbb{P}^2(\mathbb{C})),\mathrm{Bir}(\mathbb{P}^2(\mathbb{C}))]=\mathrm{Bir}(\mathbb{P}^2(\mathbb{C})).
\end{align*}}
\end{cor}

\subsection{Description du sous-groupe normal de
$\mathrm{Bir}(\mathbb{P}^2(\mathbb{C}))$ engendr\'e par une transformation birationnelle
quadratique}\hspace{1mm}

\begin{pro}\label{norm1}
{\sl Soit $A$ un automorphisme de $\mathbb{P}^2(\mathbb{C});$ on a
\begin{align*}
\mathrm{N}(A\sigma,\mathrm{Bir}(\mathbb{P}^2(\mathbb{C})))=\mathrm{Bir}(\mathbb{P}^2(\mathbb{C})).
\end{align*}}
\end{pro}

\begin{proof}[\sl D\'emonstration]
Supposons que $A^2\not=\mathrm{id}.$ La transformation
$A\sigma$ est conjugu\'ee \`a $\sigma A$ via~$\sigma;$ ainsi
$\sigma A$ appartient \`a $\mathrm{N}(A\sigma,\mathrm{Bir}(\mathbb{P}^2(\mathbb{C}))).$ On
constate alors que $A^2=(A\sigma)(\sigma A)$ est aussi dans~$\mathrm{N}(A\sigma, \mathrm{Bir}(\mathbb{P}^2(\mathbb{C})));$ or la Proposition
\ref{autnorm} assure que
\begin{align*}
\mathrm{N}(A^2, \mathrm{Bir}(\mathbb{P}^2(\mathbb{C})))=\mathrm{Bir}(\mathbb{P}^2(\mathbb{C}))
\end{align*}
donc $\mathrm{N}(A\sigma,\mathrm{Bir}(\mathbb{P}^2(\mathbb{C})))= \mathrm{Bir}(\mathbb{P}^2(\mathbb{C})).$

Supposons que $A$ soit de carr\'e trivial. Soit $B$
l'\'el\'ement de $\mathrm{PGL}_3(\mathbb{C})$ donn\'e par $(\alpha x_0:\beta x_1:~\gamma x_2);$ on a $\mathrm{N}(A\sigma,\mathrm{Bir}(\mathbb{P}^2(\mathbb{C})))=\mathrm{N}
(BA\sigma B^{-1},\mathrm{Bir}(\mathbb{P}^2(\mathbb{C}))).$ Notons que $BA\sigma
B^{-1}=BAB\sigma.$ Un calcul montre qu'il existe $\alpha,$
$\beta$ et $\gamma$ tels que $(BAB)^2$ soit distinct de
l'identit\'e; ce qui pr\'ec\`ede permet de conclure.
\end{proof}

\begin{pro}\label{norm2}
{\sl Si $A$ d\'esigne un automorphisme de $\mathbb{P}^2(\mathbb{C}),$ alors
\begin{align*}
\mathrm{N}(A\rho,\mathrm{Bir}(\mathbb{P}^2(\mathbb{C})))=\mathrm{Bir}(\mathbb{P}^2(\mathbb{C})).
\end{align*}}
\end{pro}

\begin{proof}[\sl D\'emonstration]
Lorsque $A^2\not=\mathrm{id},$ on peut reprendre l'id\'ee
de la Proposition~\ref{norm1}.

Supposons que $A^2= \mathrm{id}.$ Soit $B$
l'automorphisme de $\mathbb{P}^2(\mathbb{C})$ donn\'e par $(x_0:\alpha x_1:x_2);$ par
d\'efinition $\mathrm{N}(A\rho, \mathrm{Bir}(\mathbb{P}^2(\mathbb{C})))=\mathrm{N}(BA\rho
B^{-1},\mathrm{Bir}(\mathbb{P}^2(\mathbb{C}))).$ Remarquons que $BA\rho B^{-1}=BAB\rho$
d'o\`u
\begin{align*}
\mathrm{N}(A\rho,\mathrm{Bir}(\mathbb{P}^2(\mathbb{C})))=\mathrm{N}(BAB\rho,\mathrm{Bir}(\mathbb{P}^2(\mathbb{C}))).
\end{align*}
Un calcul assure l'existence d'un complexe $\alpha$ tel que $(BAB)^2$ ne soit pas trivial.
Ce qui pr\'ec\`ede entra\^ine donc $\mathrm{N}(A\rho,\mathrm{Bir}(\mathbb{P}^2(\mathbb{C})))= \mathrm{Bir}(\mathbb{P}^2(\mathbb{C}))).$
\end{proof}

\begin{pro}\label{norm3}
{\sl Si $A$ d\'esigne un \'el\'ement de $\mathrm{PGL}_3(\mathbb{C}),$ on a
\begin{align*}
\mathrm{N}(A\tau,\mathrm{Bir}(\mathbb{P}^2(\mathbb{C})))=\mathrm{Bir}(\mathbb{P}^2(\mathbb{C}))).
\end{align*}}
\end{pro}

\begin{proof}[\sl D\'emonstration]
Si $A^2$ est non trivial, on reprend la
d\'emonstration de la Proposition \ref{norm1}.

Si $A^2= \mathrm{id},$ on consid\`ere l'\'el\'ement $B$
de $\mathrm{PGL}_3(\mathbb{C})$ du type
\begin{align*}
&(b^2x_0:ax_0+bcx_1:dx_0+ex_1+c^2x_2), && b,\hspace{1mm} c\in\mathbb{C}^*,
\hspace{1mm} a,\hspace{1mm} d,\hspace{1mm} e\in\mathbb{C}.
\end{align*}
On a
\begin{align*}
\mathrm{N}(A\tau, \mathrm{Bir}(\mathbb{P}^2(\mathbb{C})))=\mathrm{N}(BA\tau
B^{-1},\mathrm{Bir}(\mathbb{P}^2(\mathbb{C})));
\end{align*}
de plus $\tau B^{-1}$ s'\'ecrit aussi $C\tau$ o\`u $C$ d\'esigne l'\'el\'ement de 
$\mathrm{PGL}_3(\mathbb{C})$ donn\'e par
\begin{align*}
(b^4x_0:ab^2x_0+b^3cx_1:(a^2-b^2d)x_0+b(2ac-be)x_1+b^2c^2x_2).
\end{align*}
En particulier $BA\tau B^{-1}=
BAC\tau.$ D\`es que $A$ n'est pas de la forme $(x_0:x_1:\alpha x_0+
\beta x_1-x_2)$ on peut trouver $a,$ $b,$ $c,$ $d$ et $e$ tels que
$(BAC)^2\not=\mathrm{id};$ on peut alors conclure comme
pr\'ec\'edemment. Reste \`a traiter l'\'eventualit\'e o\`u $A$ est
du type $(x_0:x_1:\alpha x_0+ \beta x_1-x_2).$ Soit $B$ l'automorphisme
de $\mathbb{P}^2(\mathbb{C})$ d\'efini par $B=(x_0:2x_1:-4x_2).$ Un calcul montre
que la transformation $A\tau BA\tau B^{-1}$ de $\mathrm{N}(A\tau,
\mathrm{Bir}(\mathbb{P}^2(\mathbb{C})))$ s'\'ecrit
\begin{align*}
(x_0:x_1:-3\alpha x_0-\beta x_1+x_2).
\end{align*}
D'apr\`es ce qui pr\'ec\`ede $\mathrm{Bir}(\mathbb{P}^2(\mathbb{C})) =\mathrm{N}(A\tau BA\tau
B^{-1},\mathrm{Bir}(\mathbb{P}^2(\mathbb{C})));$ l'inclusion
\begin{align*}
\mathrm{N}(A\tau BA\tau B^{-1},\mathrm{Bir}(\mathbb{P}^2(\mathbb{C})))\subset\mathrm{N}(A\tau,
\mathrm{Bir}(\mathbb{P}^2(\mathbb{C})))
\end{align*}
conduit \`a $\mathrm{N}(A\tau, \mathrm{Bir}(\mathbb{P}^2(\mathbb{C})))=\mathrm{Bir}(\mathbb{P}^2(\mathbb{C})).$
\end{proof}

Les Propositions \ref{norm1}, \ref{norm2} et \ref{norm3} impliquent le:

\begin{thm}\label{simple}
{\sl Pour tout \'el\'ement $f$ de $\mathrm{Bir}_2$ on a
\begin{align*}
\mathrm{N}(f,\mathrm{Bir}(\mathbb{P}^2(\mathbb{C})))=\mathrm{Bir}(
\mathbb{P}^2(\mathbb{C})).
\end{align*}}
\end{thm}

\subsection{Calculs de $\mathrm{N}(f,\mathrm{Bir}(\mathbb{P}^2(\mathbb{C})))$ dans
un cadre un peu plus g\'en\'eral}\hspace{1mm}

Le groupe $\mathrm{PGL}_2(\mathbb{C}(x_1))$ est simple
\begin{align*}
& \mathrm{N}(f,\mathrm{PGL}_2(\mathbb{C}(x_1)))=\mathrm{PGL}_2(\mathbb{C}(x_1)),&&
\forall\hspace{1mm} f\in\mathrm{PGL}_2(\mathbb{C}(x_1)).
\end{align*}

Soit $f$ une transformation de \textsc{Cremona} pr\'eservant la
fibration rationnelle $x_1=$ cte fibre \`a fibre. Alors $f$
s'\'ecrit
\begin{align*}
\left(\frac{a(x_1)x_0+b(x_1)}{c(x_1)x_0+d(x_1)},x_1\right),
\end{align*}
{\it i.e.} s'identifie \`a un \'el\'ement de $\mathrm{PGL}_2(\mathbb{C}(x_1)).$
On \'ecrit alors abusivement: $f$ est dans $\mathrm{PGL}_2(\mathbb{C}(x_1)).$

\begin{pro}\label{j0}
{\sl Pour tout \'el\'ement $f$ dans $\mathrm{PGL}_2(\mathbb{C}(x_1))$ on a
\begin{align*}
\mathrm{N}(f,\mathrm{Bir}(\mathbb{P}^2(\mathbb{C})))=\mathrm{Bir}(
\mathbb{P}^2(\mathbb{C})).
\end{align*}}
\end{pro}

\begin{proof}[\sl D\'emonstration]
Consid\'erons un \'el\'ement $f$ de $\mathrm{PGL}_2(\mathbb{C}(x_1));$ puisque $\mathrm{PGL}_2(
\mathbb{C}(x_1))$ est simple
l'involution~$(-x_0,x_1)$ s'\'ecrit comme un produit de conjugu\'es de
$f.$ Comme
\begin{align*}
& \mathrm{N}(A,\mathrm{Bir}(\mathbb{P}^2(\mathbb{C})))=\mathrm{Bir}(
\mathbb{P}^2(\mathbb{C})),&&
\forall\hspace{1mm} A \in \mathrm{PGL}_3(\mathbb{C})
\end{align*}
on a le r\'esultat annonc\'e.
\end{proof}

On obtient comme cons\'equence la:

\begin{pro}
{\sl Soit $f$ un \'el\'ement du groupe de \textsc{de Jonqui\`eres};
alors
\begin{align*}
\mathrm{N}(f,\mathrm{Bir}(\mathbb{P}^2(\mathbb{C})))=\mathrm{Bir}(
\mathbb{P}^2(\mathbb{C})).
\end{align*}}
\end{pro}

\begin{proof}[\sl D\'emonstration]
\`A conjugaison pr\`es on peut \'ecrire $f$ sous la forme $(f_1(x_0,x_1),
\gamma(x_1))$ o\`u $\gamma$
est une homoth\'etie ou bien une translation. Soit $h=(h_1(x_0,x_1),x_1)$ un
\'el\'ement de $\mathrm{PGL}_2(\mathbb{C}(x_1))$ alors $g=[f,h]$ appartient \`a
$\mathrm{N}(f,\mathrm{Bir}(\mathbb{P}^2(\mathbb{C})))\cap\mathrm{PGL}_2(
\mathbb{C}(x_1)).$ Pour $h$ bien
choisi $g$ est diff\'erent de l'identit\'e. La Proposition \ref{j0} permet de conclure.
\end{proof}

\begin{cor}\label{fibfib}
{\sl Soit $f$ un
\'el\'ement de $\mathrm{Bir}(\mathbb{P}^2(\mathbb{C})).$ S'il existe une
une transformation de \textsc{Cremona} $h$ telle que $[f,h]$
pr\'eserve une fibration rationnelle, alors
\begin{align*}
\mathrm{N}(f,\mathrm{Bir}(\mathbb{P}^2(\mathbb{C})))=\mathrm{Bir}(\mathbb{P}^2(
\mathbb{C})).
\end{align*}}
\end{cor}

On en d\'eduit la:

\begin{pro}
{\sl Soit $f=(x_1,P(x_1)-\delta x_0),$ $\delta\in\mathbb{C}^*,$ $P\in\mathbb{C}[x_1],$
$\deg P \geq 2,$ un automorphisme de \textsc{H\'enon}; on a
\begin{align*}
\mathrm{N}(f,\mathrm{Bir}(\mathbb{P}^2(\mathbb{C})))=\mathrm{Bir}(\mathbb{P}^2(
\mathbb{C})).
\end{align*}}
\end{pro}

\begin{proof}[\sl D\'emonstration]
Posons $g=(x_0,2x_1);$ le commutateur $[f,g]$ de $f$ et $g$ pr\'eserve
la fibration rationnelle $x_0=$ cte. Le Corollaire \ref{fibfib} implique le r\'esultat.
\end{proof}

\section{Une version un peu plus forte du Th\'eor\`eme de \textsc{N\oe ther}}\label{nonoplus}

\begin{thm}
{\sl Toute transformation birationnelle de $\mathbb{P}^2(
\mathbb{C})$ dans lui-m\^eme s'\'ecrit
comme un produit d'involutions standards
\begin{align*}
& A_0\sigma A_0^{-1}A_1\sigma A_1^{-1}\ldots A_p\sigma A_p^{-1}, &&
A_i\in\mathrm{PGL}_3(\mathbb{C}).
\end{align*}}
\end{thm}

\begin{proof}[\sl D\'emonstration]
Soit $f$ dans le groupe de \textsc{Cremona}. \`A conjugaison
pr\`es par $\sigma$ et/ou un \'el\'ement de $\mathrm{PGL}_3(\mathbb{C})$
on a l'alternative suivante
\begin{itemize}
\item ou bien $f$ est un automorphisme de $\mathbb{P}^2(\mathbb{C});$

\item ou bien $f$ s'\'ecrit
\begin{align*}
&\sigma A_0\ldots A_{p-1}\sigma A_p, &&A_i\in\mathrm{PGL}_3(
\mathbb{C}).
\end{align*}
\end{itemize}
Introduisons l'ensemble $\mathrm{G}$ d\'efini par
\begin{align*}
& \mathrm{G}=\{A_0\sigma A_0^{-1}\ldots A_p\sigma
A_p^{-1}\hspace{1mm} |\hspace{1mm}
A_i\in\mathrm{PGL}_3(\mathbb{C}),\hspace{1mm} p\in\mathbb{N}\}.
\end{align*}
Remarquons que $\mathrm{G}$ est un groupe qui contient $\sigma$ et sur lequel
$\mathrm{PGL}_3(\mathbb{C})$ agit par conjugaison.

Notons que si $D$ est un automorphisme diagonal, {\it i.e.} $D$
est du type $(\alpha x_0:\beta x_1:\gamma x_2)$ avec~$\alpha\beta\gamma$ 
non nul, alors
\begin{align*}
& D\sigma D^{-1}=D^2\sigma=\sigma D^{-2}.
\end{align*}
Consid\'erons un \'el\'ement de $\mathrm{G}$ du type $A\sigma
A^{-1}D_1\sigma D_1^{-1}D_2\sigma D_2^{-1}$ o\`u $A$ (resp. $D_i$)
d\'esigne un \'el\'ement de $\mathrm{PGL}_3(\mathbb{C})$ (resp. un
\'el\'ement diagonal de $\mathrm{PGL}_3(\mathbb{C})$). On constate que
\begin{align*}
& A\sigma A^{-1}D_1\sigma D_1^{-1}D_2\sigma D_2^{-1}=A\sigma A^{-1}
D_1^2D_2^{-2}
\end{align*}
donc $A\sigma A^{-1}D_1^2D_2^{-2}$ est dans $\mathrm{G}.$ Puisque
$\mathrm{PGL}_3(\mathbb{C})$ agit par conjugaison sur $\mathrm{G},$ la
transforma\-tion $\sigma A^{-1}DA$~appartient \`a $\mathrm{G}$ pour
tout $A$ dans $\mathrm{PGL}_3(\mathbb{C})$ et tout $D$ diagonal. Dit
autrement si~$B$ est un automorphisme de $\mathbb{P}^2(\mathbb{C})$
diagonalisable, $\sigma B$ est dans $\mathrm{G};$ ce qui implique,
puisque $\sigma$ est dans~$\mathrm{G},$ que $B$ est dans
$\mathrm{G}.$ Maintenant soit $C$ un \'el\'ement quelconque de
$\mathrm{PGL}_3(\mathbb{C}).$ On peut \'ecrire $C$ sous la forme
$NTN^{-1}$ o\`u $T$ est triangulaire. Pour tout automorphisme $D$
de $\mathbb{P}^2(\mathbb{C})$ on~a
\begin{align*}
C=NTN^{-1}=NTN^{-1}NDN^{-1}ND^{-1}N^{-1}=(N(TD)N^{-1})(ND^{-1}N^{-1}).
\end{align*}
En choisissant $D$ diagonal convenable, on constate que $C$ s'\'ecrit
comme produit de matrices diagonalisables, et donc est dans $\mathrm{G}.$

Finalement $\sigma$ et $\mathrm{PGL}_3(\mathbb{C})$ sont dans $\mathrm{G};$
on conclut en appliquant le th\'eor\`eme de \textsc{N\oe ther}.
\end{proof}

\clearemptydoublepage
\chapter{Exp\'eriences (orbites et ensembles de type \textsc{Julia})}\label{expexp}

Ce Chapitre est disponible sur 
\selectlanguage{english}

http://perso.univ-rennes1.fr/dominique.cerveau/ et 

http://people.math.jussieu.fr/$\sim$deserti/publications.

\clearemptydoublepage
\chapter{Transformations birationnelles de degr\'e $3$}\label{troistrois}

Il s'agit dans ce chapitre d'\'etablir la classification g.d. des
transformations cubiques. Alors que dans le cas quadratique
on distingue trois orbites de degr\'e $2$ pur, ici il y a une infinit\'e
d'orbites. L'espace quotient est en fait de dimension $2.$

\section{G\'en\'eralit\'es}\label{gen33}\hspace{1mm}

Commen\c{c}ons par un \'enonc\'e \'el\'ementaire
valable en tout degr\'e.

\begin{lem}
{\sl Soit $f$ un \'el\'ement du groupe de \textsc{Cremona}; on a
l'inclusion: $\mathrm{Ind}\hspace{1mm} f\subset\mathrm{Exc}\hspace{1mm} f.$}
\end{lem}

\begin{proof}[\sl D\'emonstration]
D'apr\`es la Proposition \ref{chaus} il suffit de v\'erifier que
si $m$ est d'ind\'etermination pour $f=(f_0:f_1:f_2)$ alors
$\mathrm{det}(\mathrm{jac}\hspace{1mm} f_{(m)})=0.$ On travaille dans
$\mathbb{C}^3;$ soit $m$ dans $\mathrm{Ind}\hspace{1mm} f$ alors
\begin{align*}
f_0(m)=f_1(m)=f_2(m)=0.
\end{align*}
Comme les hypersurfaces d'\'equation
$f_i=0$ sont invariantes par le champ radial
\begin{align*}
E=x_0\frac{\partial}{\partial x_0}+x_1\frac{\partial}{\partial
 x_1}+x_2\frac{\partial}{\partial x_2}
\end{align*} 
il existe des coordonn\'ees
locales $(u,v,w)$ en $m$ telles que $E$ soit le champ
$\frac{\partial}{\partial w}$ et les hypersurfaces d'\'equation $f_i=0$ soient donn\'ees par
$g_i(u,v)=0$ les $g_i$ \'etant holomorphes:
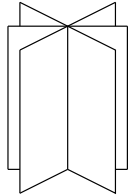
\begin{figure}[H]
\begin{center}
\input{livreouvert.pstex_t}
\end{center}
\end{figure}

Localement $f_i$ s'\'ecrit donc $U_ig_i(u,v),$ avec $U_i$
holomorphes, et ceci implique que $\mathrm{det}(\mathrm{jac}\hspace{1mm} f_{(m)})=~0.$
\end{proof}

Si $f$ est une transformation birationnelle purement cubique,
$\mathrm{det}(\mathrm{jac}\hspace{1mm} f)$ est un polyn\^ome homog\`ene
de degr\'e $6$ dont les z\'eros sont contract\'es par $f.$ Il est
imm\'ediat de constater que les composantes irr\'eductibles de
$\mathrm{det}(\mathrm{jac}\hspace{1mm} f)$ sont toutes de degr\'e
inf\'erieur ou \'egal \`a $2.$ En effet si la courbe
irr\'eductible $\mathcal{C}$ d'\'equation $P=0$ est contract\'ee
par exemple sur $(0:0:1)$ alors les deux premi\`eres composantes
$f_0$ et $f_1$ de $f$ sont divisibles par $P;$ par suite $\deg P
\leq 3.$ Si $\deg P=3$ les composantes $f_0$ et $f_1$ sont
$\mathbb{C}$-colin\'eaires, ceci conduisant \`a la
d\'eg\'en\'erescence de $f;$ on en d\'eduit que les courbes
contract\'ees par $f$ forment un arrangement de droites et de
coniques. Nous verrons que seules les configurations suivantes
peuvent arriver:

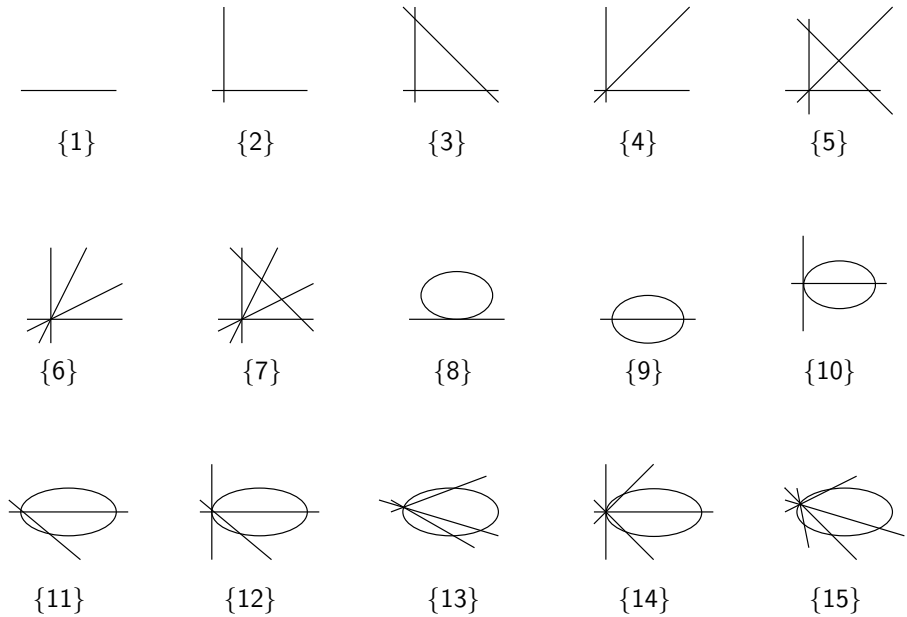
\begin{figure}[H]
\begin{center}
\input{gal.pstex_t}
\end{center}
\end{figure}

\begin{figure}[H]
\begin{center}
\input{gal2.pstex_t}
\end{center}
\end{figure}

\begin{figure}[H]
\begin{center}
\input{gal3.pstex_t}
\end{center}
\end{figure}

\bigskip

Avant de justifier ceci mentionnons le fait suivant qui s'av\`ere
important.\bigskip

\noindent\textbf{\textit{Fait.}} {\sl Soit $f$ une transformation
de \textsc{Cremona}. Puisque $\mathrm{Ind}\hspace{1mm}
f\subset\mathrm{Exc}\hspace{1mm} f$ la restriction de l'application $f$ \`a
$\mathbb{P}^2(\mathbb{C})\setminus\mathrm{Exc}\hspace{1mm} f$ est bien
d\'efinie et holomorphe sur
$\mathbb{P}^2(\mathbb{C})\setminus\mathrm{Exc}\hspace{1mm} f.$ Comme il en
est de m\^eme pour $f^{-1}$ on en d\'eduit que $f$ induit un
biholomorphisme de
$\mathbb{P}^2(\mathbb{C})\setminus\mathrm{Exc}\hspace{1mm} f$ 
sur~$\mathbb{P}^2(\mathbb{C})\setminus\mathrm{Exc}\hspace{1mm} f^{-1}.$}\bigskip

Soit $\Gamma=\Gamma_1\cup\ldots\cup\Gamma_n$ une courbe
plane dont les composantes irr\'eductibles $\Gamma_i$ sont de
degr\'e $d_i.$ On note $\Sigma=\mathbb{P}^2(\mathbb{C})\setminus
\Gamma$ le compl\'ement de $\Gamma.$ Le premier groupe d'homologie
$\mathrm{H}^1(\Sigma,\mathbb{Z})$ s'identifie au groupe
$\mathbb{Z}^n/(d_1, \ldots,d_n);$ lorsque la courbe $\Gamma$ est
\`a singularit\'es ordinaires le groupe fondamental
$\pi_1(\Sigma,*)$ est ab\'elien encore \'egal \`a $\mathbb{Z}^n/(
d_1,\ldots,d_n)$ (th\'eor\`eme de \textsc{Deligne}-\textsc{Fulton},
\emph{voir} \cite{Del}).

Une cons\'equence imm\'ediate est la suivante valable en
tout degr\'e.

\begin{pro}\label{memnbcomp}
{\sl Soit $f$ un \'el\'ement du groupe de \textsc{Cremona}. Les
ensembles $\mathrm{Exc}\hspace{1mm} f$ et $\mathrm{Exc}\hspace{1mm} f^{-1}$ ont le
m\^eme nombre de composantes irr\'eductibles.}
\end{pro}

\begin{rem}\label{deginv}
On note $\textrm{Inv}\label{not35}$ l'application qui \`a une
transformation birationnelle associe son inverse: $\textrm{Inv}(f)
=f^{-1}.$ L'involution $\textrm{Inv}$ respecte la \og
filtration\fg\hspace{1mm} par le degr\'e 
\begin{align*}
\deg(\textrm{Inv}(f))=\deg f.
\end{align*}
En effet, soient $\mathcal{D}$ et $\mathcal{D}'$ deux droites
g\'en\'eriques de $\mathbb{P}^2(\mathbb{C});$ puisque
\begin{align*}
f_*\hspace{1mm},\hspace{1mm} f^*\colon\hspace{1mm}\mathrm{H}^{1,1}(\mathbb{P}^2(\mathbb{C}))\to
\mathrm{H}^{1,1}(\mathbb{P}^2(\mathbb{C}))
\end{align*}
sont auto-adjointes on a
\begin{align*}
\deg f=\deg(f^*\mathcal{D})=(f^*\mathcal{D})
\mathcal{D}'=\mathcal{D}(f_*\mathcal{D}')=\deg(f_*\mathcal{D}')=\deg
f^{-1}.
\end{align*}
\end{rem}

\section{\og Classification\fg\hspace{1mm} des transformations birationnelles
cubiques}\hspace{1mm}

La classification repose en grande partie sur la description de
l'ensemble des courbes contract\'ees.

\subsection{Le lieu exceptionnel contient une conique}\label{conn}\hspace{1mm}

Comme on l'a dit nous allons examiner au cas par cas les
diff\'erentes configurations de ces arrangements de droites et
coniques. Dans ce qui suit on dira qu'une conique est
\textbf{\textit{r\'eduite}}\label{ind49} si elle est
ou bien lisse, ou bien constitu\'ee de deux droites distinctes.
Lorsque la conique $\mathcal{C}$ est constitu\'ee de deux droites
$\mathcal{D}_1$ et $\mathcal{D}_2$ distinctes nous dirons que
$\mathcal{C}$ est contract\'ee si les deux droites $\mathcal{D}_1$
et $\mathcal{D}_2$ le sont sur un m\^eme point. Modulo cette
convention on a le:

\begin{lem}\label{pasdeux}
{\sl Un \'el\'ement $f$ de $\mathring{\mathrm{B}}\mathrm{ir}_3$ ne peut pas contracter deux coniques
r\'eduites distinctes.}
\end{lem}

\begin{proof}[\sl D\'emonstration]
Soit $f$ une transformation birationnelle purement cubique contractant
les deux coniques $Q_1=~0$ et $Q_2=0.$
Remarquons que si les deux coniques ont une branche commune alors
$f$ contracte trois droites sur un m\^eme point ce qui est
impossible. Notons que~$f$ ne peut contracter $Q_1=~0$ et $Q_2=0$
sur un m\^eme point. On se ram\`ene donc au cas o\`u
$Q_1=0$ est contract\'e sur $(1:0:0)$ et $Q_2=0$ sur $(0:1:0);$ on
constate alors que la troisi\`eme composante de $f$ est divisible
par le produit $Q_1Q_2$ ce qui est absurde.
\end{proof}

Il en r\'esulte que $\mathrm{Exc}\hspace{1mm} f$ est une courbe de degr\'e
$6$ constitu\'ee de droites et d'au plus une conique r\'eduite.

Dans le m\^eme ordre d'id\'ee on a le:

\begin{lem}\label{uni}
{\sl Soit $f$ une transformation birationnelle purement cubique; le
lieu exceptionnel de $f$ ne peut \^etre r\'eduit \`a une
conique lisse.}
\end{lem}

\begin{proof}[\sl D\'emonstration]
Raisonnons par l'absurde: supposons que $\mathrm{Exc}\hspace{1mm} f$ soit
r\'eduit \`a une conique~$\mathcal{C}$ d'\'equation $P=0.$ \`A
conjugaison pr\`es on peut supposer que $\mathcal{C}$ est
contract\'ee par $f$ sur~$(1:0:0).$ Il s'en suit que $f$ est de la
forme $(f_0:\ell_1P:\ell_2P),$ les $\ell_i$ d\'esignant des formes
lin\'eaires n\'ecessairement ind\'ependantes. \`A conjugaison g.d.
pr\`es on peut supposer que $\ell_1=x_0$ et $\ell_2=x_1.$ La
transformation $f$ \'etant birationnelle la Proposition
\ref{chaus} assure que
\begin{align*}
\mathrm{det}(\mathrm{jac}\hspace{1mm} f)=P\left( -x_0\frac{\partial f_0} {\partial
x_0}\frac{\partial P}{\partial x_2}-x_1\frac{\partial f_0} {\partial
x_1}\frac{\partial P}{\partial x_2}+\frac{\partial f_0} {\partial
x_2}\left(x_0\frac{\partial P} {\partial x_0}+x_1\frac{\partial
P}{\partial x_1}+P\right)\right)
\end{align*}
est divisible par $P^3.$ Par ailleurs on a d'apr\`es l'identit\'e
d'\textsc{Euler}

\begin{eqnarray}
\mathrm{det}(\mathrm{jac}\hspace{1mm} f)&=&P\left(\frac{\partial P}{\partial x_2}
\left(-3f_0+x_2\frac{\partial f_0} {\partial
x_2}\right)+\frac{\partial f_0}{\partial x_2}\left(3P-
x_2\frac{\partial P}{\partial x_2}\right)\right)\nonumber\\
&=&3P\left(P\frac{\partial f_0}{\partial x_2}-f_0\frac{\partial P}
{\partial x_2}\right).\nonumber
\end{eqnarray}

\noindent On en d\'eduit que $f_0$ est divisible par $P;$ la
transformation $f$ est alors de degr\'e inf\'erieur \`a~$3.$
\end{proof}

\begin{lem}\label{inter}
{\sl Soit $f$ dans $\mathring{\mathrm{B}}\mathrm{ir}_3$ contractant une conique lisse
$\mathcal{C}.$ Supposons que $f$ contracte deux droites
$\mathcal{D}_1$ et $\mathcal{D}_2$ (au moins); alors
$\mathcal{D}_1\cap\mathcal{D}_2$ appartient \`a $\mathcal{C}.$

En particulier si $f$ contracte $\mathcal{D}_1,$
$\ldots,$ $\mathcal{D}_k,$ toutes ces droites sont concourantes
avec point d'intersection sur $\mathcal{C}.$}
\end{lem}

\begin{proof}[\sl D\'emonstration]
On peut se ramener \`a $\mathcal{D}_1:(x_0=0),$ $\mathcal{D}_2:(x_1
=0).$ On peut aussi supposer que $\mathcal{C}:(P=0)$ est contract\'ee
sur $(1:0:0),$ $\mathcal{D}_1$ sur $(0:1:0);$ ces conditions \'etant
satisfaites la transformation
$f$ est du type suivant $(x_0q:\ell P:x_0P)$ avec $q$ forme
quadratique et $\ell$ forme lin\'eaire. La droite
$\mathcal{D}_2$ est alors contract\'ee, \`a conjugaison pr\`es,
sur $(0:0:1)$ ou bien sur $(A:B:0).$  Si $\mathcal{D}_2$ est
contract\'ee sur $(A:B:0),$ alors $P$ est divisible par
$x_1$ ce qui contredit l'irr\'eductibilit\'e. Ainsi $\mathcal{D}_2$ est
contract\'ee sur $(0:0:1)$ et $f=(x_0x_1 \ell:x_1P:x_0P).$ On pose
\begin{align*}
&\ell:=ax_0+bx_1+cx_2, && P:=\alpha x_0^2+\beta x_1^2+\gamma
x_2^2+\delta x_0x_1+\varepsilon x_0x_2+\mu x_1x_2
\end{align*}
et on veut montrer que $\mathcal{D}_1\cap\mathcal{D}_2$ appartient
\`a $\mathcal{C},$ {\it i.e.} que $\gamma=0.$

On va utiliser la birationnalit\'e de $f$ que l'on
\'ecrit dans la carte affine $x_2=1$
\begin{align*}
f(x_0,x_1)=\left(\frac{x_1(ax_0+bx_1+c)}{\alpha x_0^2+\beta x_1^2+
\gamma+\delta x_0x_1+\varepsilon x_0+\mu x_1},\frac{x_1}{x_0}\right).
\end{align*}
On doit exprimer les conditions sur les coefficients pour que $f(x_0,x_1)=(u,v)$
ait g\'en\'eriquement une seule solution. On obtient $x_1=vx_0$ et
\begin{equation}\label{birr}
vx_0(ax_0+bvx_0+c)=u(\alpha x_0^2+\beta v^2x_0^2+\gamma+\delta vx_0^2
+\varepsilon x_0+\mu vx_0).
\end{equation}

On remarque que si le discriminant de l'\'equation (\ref{birr})
est identiquement nul et $\gamma\not=0$ alors~$a=b=c=0,$ cas
exclu. Donc pour que cette \'equation quadratique en $x_0$ ait
g\'en\'eriquement une unique solution il faut que le trin\^ome
ci-dessus poss\`ede une racine $x_0=\eta$ ind\'ependante de $(u,
v).$ C'est le cas par exemple si $\gamma=0$ ($\eta=0$). Une autre
possibilit\'e imm\'ediate~est
\begin{align*}
a=b=\alpha=\beta=\delta=0;
\end{align*}
mais dans ce cas on
v\'erifie que $P$ est divisible par $x_2$ ce qui contredit le
fait que $\mathcal{C}$ soit irr\'eductible.

Supposons que $x_0=\eta\not=0$ soit solution de
(\ref{birr}). On peut voir
\begin{align*}
\eta^2(av+bv^2-u(\alpha+\beta v^2+ \delta v))+\eta(vc-u\varepsilon-
\mu vu)- u\gamma=0
\end{align*}
comme un polyn\^ome en $u,$ $v$ identiquement nul; on obtient
\begin{align*}
& b=\beta=0,&& \mu= -\eta\delta,&& c=-\eta a,&& \gamma=-\eta^2\alpha-
\eta\varepsilon
\end{align*}
de sorte que
\begin{align*}
f=(ax_0x_1:x_1(\alpha(x_0+\eta x_2)+\varepsilon x_2+\delta x_1):x_0(\alpha(x_0
+\eta x_2)+\varepsilon x_2+\delta x_1))
\end{align*}
qui s'av\`ere alors quadratique.

Finalement seule la solution $\gamma=0$ convient.
\end{proof}

\begin{rems}
\begin{itemize}
\item Soit $f$ une transformation rationnelle purement cubique. La
preuve ci-dessus montre que si $f$ contracte une conique lisse et deux
droites qui se coupent sur la conique, alors $f$ est
birationnelle; de plus on dispose de la forme normale
\begin{align*}
(x_0x_1L_3:x_1(x_0L_1+x_1L_2):x_0(x_0L_1+x_1L_2)),
\end{align*}
les $L_i$ d\'esignant des formes lin\'eaires (satisfaisant des conditions
g\'en\'eriques \'evidentes).

\item Soit $f$ une transformation birationnelle purement cubique
donn\'ee par la forme normale ci-dessus. On peut voir que,
g\'en\'eriquement sur $f,$ il y a deux autres droites
contract\'ees n\'ecessairement du type $x_1=\eta x_0.$ Pour cel\`a on
\'evalue $f$ dans la carte affine $x_2=1;$ sur la droite $x_1=\eta x_0$
 on a en reprenant les notations de la d\'emonstration du Lemme
 \ref{inter}
\begin{align*}
f(x_0,\eta x_0)=\left(\frac{\eta
(a+b\eta)x_0+\eta c}{(\alpha+\beta\eta^2+\delta\eta)x_0+\varepsilon+\mu\eta},\eta x_0\right).
\end{align*}
On cherche alors les $\eta$ pour lesquels la premi\`ere composante
est une constante en $x_0,$ {\it i.e.} les $\eta$ pour lesquels
\begin{align*}
\left\vert\begin{array}{cc}
a+b\eta & c\\
\alpha+\beta\eta^2+\delta\eta & \varepsilon+\mu\eta\\
\end{array}
\right\vert=(b\mu-c\beta)\eta^2+(a\mu+b\varepsilon-c\delta)\eta+(a\varepsilon-c\alpha)
\end{align*}
est nul; g\'en\'eriquement il y a deux solutions.
\end{itemize}
\end{rems}

\begin{pro}
{\sl Soit $f$ une transformation birationnelle purement cubique qui
contracte une conique lisse; $f$ a alors au plus $5$ points
d'ind\'etermination distincts.}
\end{pro}

\begin{proof}[\sl D\'emonstration]
Si $f$ contracte une conique $\mathcal{C},$ alors $\mathrm{Exc}\hspace{1mm}
f$ n'est pas r\'eduit \`a $\mathcal{C}$ (Lemme~\ref{uni}) et les
autres courbes contract\'ees par $f$ sont n\'ecessairement des
droites (Lemme \ref{pasdeux}). D\'esignons par $\mathcal{D}_1,$ $\ldots,$
$\mathcal{D}_k$ les droites contract\'ees par $f;$ puisque
$\mathrm{det}(\mathrm{jac}\hspace{1mm} f)$ est de degr\'e $6,$ l'entier $k$ est major\'e
par~$4.$ De sorte que $\mathrm{Exc}\hspace{1mm} f$ a au plus $5$ composantes
et d'apr\`es la Proposition \ref{memnbcomp} l'ensemble $\mathrm{Exc}\hspace{1mm}
f^{-1}$ aussi. Mais les courbes exceptionnelles de $f^{-1}$
sont contract\'ees sur les points d'ind\'etermination de $f.$
\end{proof}

D'apr\`es les Lemmes \ref{pasdeux} et \ref{uni} une
transformation birationnelle cubique qui contracte une conique
lisse contracte au moins une droite. On a l'\'enonc\'e de forme
normale suivant.

\begin{lem}\label{bluesky}
{\sl Soit $f$ une transformation birationnelle purement cubique
qui contracte une conique lisse (et donc au moins une droite). Il existe
$a,$ $b$ dans $\mathbb{C}$ et $L_1,$ $L_2,$ $L_3$ trois formes
lin\'eaires tels qu'on ait \`a conjugaison gauche-droite pr\`es
\begin{align*}
f=(x_1(ax_0+bx_1)L_3:x_0(x_0L_1+x_1L_2):x_1(x_0L_1+x_1L_2)).
\end{align*}}
\end{lem}

\begin{proof}[\sl D\'emonstration]
Notons $f_0,$ $f_1$ et $f_2$ les composantes de $f.$ Soit $P=0$
l'\'equation de la conique lisse $\mathcal{C}$ que l'on peut
supposer \^etre contract\'ee sur $(1:0:0).$ On en d\'eduit que 
l'existence de deux formes lin\'eaires ind\'ependantes $\alpha_1$
et $\alpha_2$ telles que $f=(f_0:\alpha_1P:\alpha_2P).$ Soit $\ell_1=0$
l'\'equation d'une droite $\mathcal{D}_1$ contract\'ee par $f.$ Si
$\mathcal{D}_1$ \'etait contract\'ee sur le point $(1:0:0),$ les
formes lin\'eaires $\alpha_i$ seraient multiples et $f$ serait
d\'eg\'en\'er\'ee; on peut donc supposer, \`a conjugaison g.d.
pr\`es, que $\mathcal{D}_1$ est contract\'ee sur le point
$(0:1:0).$ Il s'en suit que $\alpha_2$ co\"incide, \`a
multiplication par un scalaire pr\`es, avec $\ell_1$ et qu'il
existe une forme quadratique $q$ telle que~$f_0=\alpha_2q$
\begin{align*}
f=(\alpha_2q :\alpha_1P:\alpha_2P).
\end{align*}
Comme $\alpha_1$ et
$\alpha_2$ sont ind\'ependantes $f$ s'\'ecrit \`a conjugaison g.d.
pr\`es $(x_1q:x_0P:x_1P),$ soit~$\left(\frac{q}{P},\frac{x_0}{x_1}\right)$
dans la carte affine $x_2=1.$ La transformation $f$ est
birationnelle si pour tous $\zeta,$ $\mu$ le syst\`eme
\begin{align*}
\left\{\begin{array}{ll}
q-\zeta P=0\\
x_0-\mu x_1=0
\end{array}\right.
\end{align*}
a une unique solution en dehors des points d'ind\'etermination.
Dit autrement tout \'el\'ement du pinceau de coniques d\'efini par
$q-\zeta P=0$ doit couper les \'el\'ements du pinceau de 
droi\-tes~$x_0-~\mu x_1=~0$ en un unique point en dehors de $\mathrm{Ind}\hspace{1mm} f.$
Puisque g\'en\'eriquement une droite coupe une conique en deux
points chaque conique du pinceau $q-\zeta
P=0$ passe par $(0,0),$ {\it i.e.}~$q$ et $P$ s'annulent en
$(0,0);$ ainsi $f$ est du type
\begin{align*}
(x_1(x_0\ell_1+x_1\ell_2):x_0(x_0L_1+x_1L_2):x_1(x_0L_1+x_1L_2)).
\end{align*}
Le pinceau d\'efini par $\eta(x_0\ell_1
+x_1\ell_2)+\nu(x_0L_1+x_1L_2)=0$ contient une conique
d\'eg\'en\'er\'ee; $\mathcal{C}$ \'etant lisse, cette d\'eg\'en\'erescence
n'a pas lieu lorsque $\eta$ est nul et il existe $\nu$ tel que
\begin{align*}
x_0\ell_1+x_1\ell_2+\nu(x_0L_1+x_1L_2)= (ax_0+bx_1)L_3.
\end{align*}
Il en r\'esulte qu'\`a conjugaison g.d. pr\`es $f$ est de la forme
\begin{align*}
(x_1(ax_0+bx_1)L_3:x_0(x_0L_1 +x_1L_2):x_1(x_0L_1+x_1L_2)).
\end{align*}
\end{proof}

Ce Lemme s'av\`ere utile pour \'etablir la:

\begin{pro}\label{condte}
{\sl Soit $f$ un \'el\'ement de $\mathring{\mathrm{B}}\mathrm{ir}_3.$ Supposons que $\mathrm{Exc}\hspace{1mm}
f$ soit constitu\'e d'une seule droite et d'une conique lisse.
\`A conjugaison gauche-droite pr\`es on a
\begin{itemize}
\item lorsque la droite est tangente \`a la conique, $f$ est de la
forme
\begin{align*}
&(x_0(x_0^2+x_1x_2):x_1^3:x_1(x_0^2+x_1x_2))&&\mathsf{\{8\}};
\end{align*}

\item lorsque la droite n'est pas tangente \`a la conique, $f$ est
du type suivant
\begin{align*}
&(x_1^2x_2:x_0(x_0x_2+x_1^2):x_1(x_0x_2+x_1^2))&&\mathsf{\{9\}}.
\end{align*}
\end{itemize}}
\end{pro}

La notation $\mathsf{\{k\}}$ signifie que l'ensemble exceptionnel
de la transformation consid\'er\'ee
pr\'esente la configuration $\mathsf{\{k\}}$ (\emph{voir}
Chapitre \ref{troistrois}, \S\ref{gen33}).

\begin{proof}[\sl D\'emonstration]
D'apr\`es le Lemme \ref{bluesky} et sa d\'emonstration la transformation $f$ est, \`a
conjugaison g.d. pr\`es, du type $(x_1(ax_0+bx_1)L_3:x_0(x_0L_1+x_1L_2):x_1
(x_0L_1+x_1L_2)),$ la conique contract\'ee \'etant d\'efinie par $x_0L_1+x_1L_2=0$
et la droite contract\'ee par $x_1=0.$ Un calcul montre que
\begin{align*}
\mathrm{det}(\mathrm{jac}\hspace{1mm} f)=3(x_0L_1+x_1L_2)x_1(ax_0+bx_1)\left((x_0L_1+x_1L_2)\frac{\partial
L_3}{\partial x_2}-L_3\frac{\partial (x_0L_1+x_1L_2)}{\partial x_2}\right).
\end{align*}
Puisque $\mathrm{det}(\mathrm{jac}\hspace{1mm} f)$ s'annule
uniquement sur les courbes contract\'ees et que $\mathcal{C}$ est
lisse, le coefficient $a$ est nul. De plus
\begin{align*}
(x_0L_1+x_1L_2)\frac{\partial L_3}{\partial x_2}-L_3\frac{\partial
(x_0L_1+x_1L_2)}{\partial x_2}
\end{align*}
est un multiple de $x_1^p(x_0L_1+x_1L_2)^q,$ le couple $(p,q)$ appartenant
\`a $\{(2,0),\hspace{1mm}(1,1)\}.$ Posons
\begin{align*}
&L_i:=a_ix_0+b_ix_1+c_ix_2.
\end{align*}

Si $p=2,$ alors ou bien $(a_1,a_2)=\left(\frac{c_1
a_3}{c_3},\frac{b_3c_1+a_3c_2-c_3b_1}{c_3}
\right),$ ou bien $c_1=a_3=c_3=0.$ Dans le premier cas,
$c_3$ est non nul et on peut supposer que $L_3=x_2;$ on obtient
pour mod\`ele \`a conjugaison g.d. pr\`es
\begin{align*}
(x_1^2 x_2:x_0(x_1^2+x_0x_2):x_1(x_1^2+x_0x_2)).
\end{align*}
La configuration des courbes contract\'ees est donn\'ee par:
\begin{figure}[H]
\begin{center}
\input{config17.pstex_t}
\end{center}
\end{figure}

On constate que les points d'ind\'etermination de $f$ sont les
points d'intersection de la droite~$x_1=~0$ avec la conique $x_1^2+x_0x_2=0.$

La seconde \'eventualit\'e conduit \`a conjugaison g.d.
pr\`es \`a
\begin{align*}
&(x_1^3:x_0P:x_1P), && P=a_1x_0^2+(b_1+a_2)x_0x_1+b_2x_1^2+c_2x_1x_2;
\end{align*}
\`a conjugaison g.d. pr\`es on obtient
\begin{align*}
(x_0(x_0^2+x_1x_2):x_1^3:x_1(x_0^2+x_1x_2));
\end{align*}
la configuration des courbes contract\'ees est alors donn\'ee par:
\begin{figure}[H]
\begin{center}
\input{config18.pstex_t}
\end{center}
\end{figure}

Notons que l'unique point d'ind\'etermination de $f$ est
le point de tangence entre la droi\-te~$x_1=0$ et la conique
$x_1^2+x_0x_2=0.$

Lorsque $q=1,$ on note que $c_1$ et $c_2$ sont nuls
d'o\`u le mod\`ele suivant
\begin{align*}
& (x_1^2(a_3x_0+b_3x_1+c_3x_2):x_0P:x_1P), && P=a_1x_0^2+(b_1+a_2)x_0x_1+b_2x_1^2
\end{align*}
auquel cas $\mathcal{C}$ n'est pas lisse.
\end{proof}

\begin{rem}
\begin{itemize}
\item L'inverse de $f=(x_0(x_0^2+x_1x_2):x_1^3:x_1(x_0^2+x_1x_2))$
est
\begin{align*}
(x_0x_1x_2:x_1x_2^2:x_2^3-x_0^2x_1);
\end{align*}
on constate que $\mathrm{Exc}\hspace{1mm} f^{-1}
= \{x_1=0,\hspace{1mm} x_2=0\}.$ Les lieux exceptionnels de $f$ et $f^{-1}$ \og ne
sont donc pas de m\^eme nature\fg : $f$ contracte une conique
et une droite tangente \`a cette conique et $f^{-1}$
deux droites. Ce ph\'enom\`ene n'arrive pas en degr\'e
deux o\`u dans ce cas les configurations de $\mathrm{Exc}\hspace{1mm}
f$ et $\mathrm{Exc}\hspace{1mm} f^{-1}$ sont lin\'eairement les m\^emes.
Remarquons, comme l'indique l'\'enonc\'e
\ref{memnbcomp}, que les groupes fondamentaux des compl\'ements
des deux ensembles exceptionnels sont isomorphes.

\item L'inverse de $f=(x_1^2x_2:x_0(x_0x_2+x_1^2):x_1(
x_0x_2+x_1^2))$ est
\begin{align*}
(x_1(x_2^2-x_0x_1):x_2(x_2^2-x_0x_1):x_0x_2^2);
\end{align*}
on remarque que
$\mathrm{Exc}\hspace{1mm} f^{-1}=\{x_2=0,\hspace{1mm} x_2^2-x_0x_1=0\},$ dit autrement $f$ et
$f^{-1}$ contractent une conique et une droite non tangente \`a la
conique.
\end{itemize}
\end{rem}

\begin{rem}\label{te2}
Soit $f$ une transformation birationnelle purement cubique dont le
lieu exceptionnel contient une conique r\'eduite $P=0$ et au moins
deux autres droites. Ces courbes ne peuvent pas \^etre
contract\'ees sur $3$ points distincts align\'es. En effet \`a
conjugaison g.d. pr\`es on peut supposer que $\mathrm{Exc}\hspace{1mm} f=\{
x_0=0,\hspace{1mm} x_1=0,\hspace{1mm} P=0\}.$ Raisonnons par l'absurde: supposons que $P=0,$
$x_0=0,$ $x_1=0$ soient respectivement contract\'es sur les points
$(0:0:1),$ $(0:1:0)$ et~$(0:1:1).$ On constate que la premi\`ere
composante de $f$ est divisible par $x_0x_1P:$ ce qui n'est pas possible.
\end{rem}

Le lemme qui suit est l'argument fondamental de la classification
des transformations birationnelles cubiques.

\begin{lem}\label{vent}
{\sl Soit $f$ une transformation birationnelle purement cubique.
On suppose que $f$ contracte une conique r\'eduite $\mathcal{C}$
et au moins deux droites distinctes. Alors $f$ est, \`a
conjugaison gauche-droite pr\`es, de l'un des types suivants

\begin{small}
\begin{align*}
&\mathfrak{(a)}&&(x_0(\alpha x_0^2+\beta x_1^2+\gamma x_0x_1
+\delta x_0x_2+\varepsilon x_1x_2):x_1
(\alpha x_0^2+\beta x_1^2+\gamma x_0x_1+\delta x_0x_2+\varepsilon
x_1x_2):x_0x_1x_2);\\
&\mathfrak{(b)}&&(x_0(\alpha x_0^2+\beta x_1^2+\gamma x_0x_1+
\delta x_0x_2+\varepsilon x_1x_2):x_1
(\alpha x_0^2+\beta x_1^2+\gamma x_0x_1+\delta x_0x_2+\varepsilon
x_1x_2):x_0x_1^2);\\
&\mathfrak{(c)}&&(x_0(x_0^2+x_1x_2+\gamma x_0x_2):x_1(x_0^2+
x_1x_2+\gamma x_0x_2):x_0x_1(x_0-x_1));\\
&\mathfrak{(d)}&&(x_0x_2(x_1+\gamma x_0):x_1x_2(x_1+\gamma x_0):
x_0x_1(x_0-x_1)).
\end{align*}\label{not31}\label{not32}\label{not33}\label{not34}\end{small}
Tous ces mod\`eles sont birationnels d\`es qu'ils sont
purement cubiques.}
\end{lem}

\begin{proof}[\sl D\'emonstration]
Supposons que les deux droites contract\'ees aient pour \'equation
$x_0=~0,$ resp. $x_1=0$ et que $\mathcal{C}$ soit donn\'ee par $P=0.$
On se ram\`ene, gr\^ace au Lemme \ref{pasdeux} et \`a la Remarque~\ref{te2}, au cas o\`u $x_0=0,$ $x_1=0$ et $\mathcal{C}$ sont
contract\'ees sur 
\begin{align*}
&(0:1:0), &&(1:0:0)&& \text{et} &&(0:0:1)
\end{align*}
respectivement. Il existe alors une forme lin\'eaire $\ell$ telle
que $f=(x_0P:x_1P:x_0x_1\ell).$ \`A conjugaison pr\`es on peut supposer
que $\ell$ est l'une des trois formes lin\'eaires $x_2,$ $x_1-x_0,$ $x_1.$
\'Ecrivons~$P$ sous la forme
\begin{align*}
& P=ax_0^2+bx_1^2+cx_0x_1+dx_0x_2+ex_1x_2+\varepsilon x_2^2.
\end{align*}
Puisque $f$ est purement cubique $P$ n'est pas divisible par $x_0,$ resp.
$x_1;$ dans la carte affine $x_0=1$ la transformation $f$ s'\'ecrit
\begin{align*}
& f(x_1,x_2)= \left(x_1,\frac{x_2}{a+bx_1^2+cx_1+dx_2+ex_1x_2+\varepsilon x_2^2}
\right).
\end{align*}
La birationnalit\'e de $f$ implique que les \'equations $f(x_1,x_2)=(u,v)$
poss\`edent g\'en\'eriquement une seule solution; on d\'emontre alors
sans difficult\'e que $\varepsilon=0.$ Lorsque $\ell=x_2$ (resp. $\ell=x_1$)
on obtient la premi\`ere (resp. deuxi\`eme) forme normale. Supposons
que $\ell=x_1-x_0.$ Observons que $d$ et $e$ ne peuvent pas \^etre
simultan\'ement nuls; comme $x_0$ et $x_1$ jouent un r\^ole
sym\'etrique on peut supposer que
\begin{align*}
& e=1 && \text{et} && P=ax_0^2+bx_1^2+cx_0x_1+dx_0x_2+x_1x_2
\end{align*}
\`A composition \`a droite pr\`es par $(x_0:x_1:x_2-(c-bd)x_0-bx_1)$ la
transformation $f$ s'\'ecrit
\begin{align*}
(x_0(\alpha x_0^2+x_1x_2+\gamma x_0x_2):x_1(\alpha x_0^2+x_1x_2+\gamma x_0x_2):
x_0x_1(x_0-x_1)).
\end{align*}
On note que $\alpha$ et $\gamma$ ne sont pas simultan\'ement nuls. Si $\alpha$ est non nul on se
ram\`ene \`a la forme normale
\begin{align*}
& (x_0(x_0^2+x_1x_2+\gamma x_0x_2):x_1(x_0^2 +x_1x_2+\gamma x_0x_2):
x_0x_1(x_0-x_1)), &&\mathfrak{(c)}
\end{align*}
qui contracte la conique lisse $x_0^2+x_1x_2+\gamma x_0x_2=0$ et les
droites
\begin{align*}
&x_0=0,&&x_1=0,&&x_1=x_0,&&x_1=-\gamma x_0.
\end{align*}
Ces droites sont au nombre
de quatre si $\gamma\not\in\{0,\hspace{1mm} 1\};$ de plus la droite $x_1=-\gamma x_0$ est
tangente \`a la conique lisse d'\'equation $x_0^2+x_1x_2+\gamma x_0x_2=0.$
Lorsque $\alpha$ est nul, on obtient la forme normale
\begin{align*}
& (x_0x_2(x_1+\gamma x_0):x_1x_2(x_1+\gamma x_0):x_0x_1(x_0-x_1)), &&
\mathfrak{(d)}
\end{align*}
qui contracte les cinq droites $x_0=0,$ $x_1=0,$ $x_1=x_0,$ $x_1+\gamma x_0=0,$ $x_2=0.$
\end{proof}

Dans l'Appendice \ref{demdem} on r\'eduit le nombre de  param\`etres
apparaissant dans les diff\'erents mod\`eles du Lemme \ref{vent}
\`a conjugaison g.d. pr\`es (\S\ref{reduc}). On compare ensuite
les mod\`eles poss\'edant la m\^eme configuration de courbes
contract\'ees; on constate en utilisant les configurations de
$\mathrm{Ind}\hspace{1mm} f^{-1}$ que certains sont g.d. conjugu\'es,
d'autres non (\S\ref{compa}). Cette \'etude conduit \`a la:

\begin{pro}\label{listee}
{\sl Soit $f$ un \'el\'ement de $\mathring{\mathrm{B}}\mathrm{ir}_3.$ On suppose que $f$
contracte une conique r\'eduite $\mathcal{C}$ et au moins deux
droites distinctes. Alors $f$ est, \`a conjugaison gauche-droite
pr\`es, de l'une des formes suivantes

\begin{align*}
& 1. && (x_0(x_0+x_1)(x_1+x_2):x_1(x_0+x_1)(x_1+x_2):x_0x_1x_2),\hspace{2mm}
\mathsf{\{5\}};\\
& 2. && (x_0(x_0+x_1+x_2)(x_0+x_1):x_1(x_0+x_1+x_2)(x_0+x_1):x_0x_1x_2),\hspace{2mm}
\mathsf{\{5\}};\\
& 3. && (x_0x_2(x_0+x_1):x_1x_2(x_0+x_1):x_0x_1^2),\hspace{2mm} \mathsf{\{5\}};\\
& 4. && (x_0(x_0^2+x_1^2+\gamma x_0x_1):x_1(x_0^2+x_1^2+\gamma x_0x_1):x_0
x_1x_2),\hspace{1mm} \gamma^2\not=4,\hspace{2mm} \mathsf{\{6\}};\\
& 5. && (x_0x_2(x_1+x_0):x_1x_2(x_1+x_0):x_0x_1(x_0-x_1)), \hspace{2mm}
\mathsf{\{7\}}; \\
& 6. && (x_0(x_0^2+x_1^2+\gamma x_0x_1+\gamma_+x_0x_2+x_1x_2):x_1(x_0^2+
x_1^2+\gamma x_0x_1+\gamma_+x_0x_2+x_1x_2):x_0x_1x_2),\\
& && \gamma^2\not=4,\hspace{2mm}\gamma_\pm=\frac{\gamma\pm\sqrt{\gamma^2
-4}}{2},\hspace{2mm} \mathsf{\{7\}};\\
& 7. && (x_0(x_1^2+x_0x_2):x_1(x_1^2+x_0x_2):x_0x_1x_2),\hspace{2mm}
\mathsf{\{10\}};\\
& 8. && (x_0(x_1^2+x_0x_2):x_1(x_1^2+x_0x_2):x_0x_1^2),\hspace{2mm}\mathsf{\{10\}};\\
& 9. &&(x_0(x_0^2+x_1x_2):x_1(x_0^2+x_1x_2):x_0x_1^2),\hspace{2mm}
\mathsf{\{10\}};\\
& 10. &&(x_0(x_0^2+x_1x_2+x_0x_2):x_1(x_0^2+x_1x_2+x_0x_2):x_0x_1x_2),\hspace{2mm}
\mathsf{\{11\}};\\
& 11. &&(x_0(x_0x_1+x_0x_2+x_1x_2):x_1(x_0x_1+x_0x_2+x_1x_2):x_0x_1x_2),\hspace{2mm}
\mathsf{\{11\}};\\
& 12. &&(x_0(x_0^2+x_1x_2):x_1(x_0^2+x_1x_2):x_0x_1(x_0-x_1)),\hspace{2mm}
\mathsf{\{12\}};\\
& 13. &&(x_0(x_0^2+x_0x_1+x_1x_2):x_1(x_0^2+x_0x_1+x_1x_2):x_0x_1x_2),\hspace{2mm}
\mathsf{\{12\}};\\
& 14. &&(x_0(x_1^2+\gamma x_0x_1+x_1x_2+x_0x_2):x_1(x_1^2+\gamma x_0x_1+x_1x_2+x_0
x_2):x_0x_1x_2),\hspace{1mm}\gamma\not\in\{0,1\},\hspace{2mm} \mathsf{\{13\}};\\
& 15. && (x_0(x_0^2+x_1^2+\gamma x_0x_1+x_0x_2):x_1(x_0^2+x_1^2+\gamma x_0x_1+x_0x_2)
:x_0x_1x_2),\hspace{1mm}\gamma^2\not=4, \hspace{2mm} \mathsf{\{14\}}; \\
& 16. && (x_0(x_0^2+x_1x_2+x_0x_2):x_1(x_0^2+x_1x_2+x_0x_2):x_0x_1(x_0-x_1)), \hspace{2mm}
\mathsf{\{14\}};\\
& 17. && (x_0(x_0^2+x_1^2+\gamma x_0x_1+\delta x_0x_2+x_1x_2):x_1(x_0^2+x_1^2+\gamma x_0x_1+\delta x_0x_2+x_1x_2):x_0x_1x_2),\\
& &&\gamma^2\not=4,\hspace{1mm}
\delta\not\in\left\{\frac{\gamma\pm\sqrt{\gamma^2-4}}{2}\right\},
\hspace{2mm}\mathsf{\{15\}}
\end{align*}}
\end{pro}

\begin{rem}
Consid\'erons la transformation $f$ de la Proposition \ref{listee}
d\'efinie par
\begin{align*}
(x_0x_2(x_1+x_0):x_1x_2(x_1+x_0):x_0x_1(x_0-x_1)), \hspace{2mm}
\mathsf{\{7\}};
\end{align*}
On v\'erifie que
\begin{align*}
\mathrm{Exc}\hspace{1mm} f=\{x_2=0,\hspace{1mm} x_1=0,\hspace{1mm} x_0=0,
\hspace{1mm} x_1+x_0=0,\hspace{1mm} x_0=x_1\}; \\
\mathrm{Ind}\hspace{1mm} f=\{(0:1:0),\hspace{1mm}(0:0:1),\hspace{1mm}(1:0:0),\hspace{1mm}(1:1:0)\}.
\end{align*}
On constate en particulier que le nombre de points d'ind\'etermination
diff\`ere du nombre de courbes contract\'ees ce qui n'arrive pas en
degr\'e $2;$ la raison est la suivante: les droites $x_2=0$ et $x_1+
x_0=0$ sont contract\'ees sur un m\^eme point, le point $(0:0:1).$

Puisque $f$ est une involution, son inverse pr\'esente la m\^eme
configuration. En \'eclatant le point $(0:0:1)$ on constate que le diviseur
exceptionnel est envoy\'e sur $x_2=0$ et sur ce diviseur exceptionnel il
y a encore un point d'ind\'etermination qui lui est envoy\'e sur
$x_0+x_1=0.$

Le degr\'e $3$ est le premier degr\'e o\`u l'on constate ce ph\'enom\`ene:
un point d'ind\'etermination est \'eclat\'e sur plusieurs courbes.

Mentionnons un autre \'el\'ement $f$ du groupe de \textsc{Cremona}
pour lequel un point d'ind\'etermination est \'eclat\'e sur plusieurs 
courbes et pour lequel le comportement de $f^{-1}$ diff\`ere de celui
de $f$ (\emph{voir} \cite{BK2})
\begin{small}
\begin{align*}
&((x_0x_1-x_0^2)^nx_2:x_0^{n-1}(x_1-x_0)^{n+1}(x_0+x_2)+x_2
\displaystyle\sum_{j=0}^n a_j(x_0x_1-x_0^2)^{n-j}(x_2^2-x_0x_1
-x_1x_2)^j:\\
&x_2(x_0x_1-x_0^2)^{n-1}(x_2^2-x_0x_1-x_1x_2)).
\end{align*}
\end{small}
Cette transformation a quatre courbes exceptionnelles
\begin{align*}
& x_0=0, && x_0=x_1, && x_2=0, && -x_2^2+x_0x_1+x_1x_2=0;
\end{align*}
les trois premi\`eres sont contract\'ees sur le point $(0:1:0).$
Les composantes de $f^{-1}$ sont donn\'ees par
\begin{align*}
&x_0^nx_2\left(\displaystyle\sum_{j=0}^na_jx_0^{n-j}x_2^j-x_0^{n-1}(x_0+x_1)\right),\\
&(x_0+x_2)\left(\displaystyle\sum_{j=0}^na_jx_0^{n-j}x_2^j-x_0^{n-1}x_2^j- x_0^{n-1}(x_0+x_1)\right)^2,\\
& x_0^{n-1}x_2\left(x_0^{n-1}(x_0^2+x_0x_1+x_1x_2)-(x_0+x_2)\displaystyle\sum_{j=0}^na_jx_0^{n-j}x_2^j\right);
\end{align*}
son ensemble exceptionnel est form\'e de 
\begin{small}
\begin{align*}
&x_0=0, && x_2=0, && x_0^n+x_0^{n-1}x_1-\displaystyle\sum_{j=0}^na_jx_0^{n-j}x_2^j=0, &&
x_0^{n+1}-(x_0+x_2)(x_0^n+x_0^{n-1}x_1-\displaystyle\sum_{j=0}^na_jx_0^{n-j}x_2^j)=0
\end{align*}
\end{small}
Les deux premi\`eres courbes sont contract\'ees sur $(0:1:0),$ la
troisi\`eme sur $(0:0:1)$ et la derni\`ere sur $(1:1:0).$
\end{rem}

\subsection{Le lieu exceptionnel ne contient pas de conique}\label{sconn}\hspace{1mm}\bigskip

Nous continuons la classification en \'etudiant les transformations
cubiques dont l'ensemble exceptionnel est un arrangement
de droites.

\subsubsection{Le lieu exceptionnel est r\'eduit \`a une droite}\hspace{1mm}

Soit $f$ une transformation birationnelle qui contracte
exactement une droite que l'on supposera \^etre $x_2=0.$ On a
\begin{align*}
\mathrm{H}^1(\mathbb{P}^2(
\mathbb{C})\setminus\{x_2=0\},\mathbb{Z})=\mathrm{H}^1(\mathbb{P}^2(
\mathbb{C})
\setminus\mathrm{Exc}\hspace{1mm} f^{-1},\mathbb{Z})=0
\end{align*}
ce qui implique que
$\mathrm{Exc}\hspace{1mm} f^{-1}$ est aussi une droite. \`A conjugaison g.d.
pr\`es on se ram\`ene \`a
\begin{align*}
\mathrm{Exc}\hspace{1mm} f= \mathrm{Exc}\hspace{1mm}
f^{-1}=\{x_2=0\}
\end{align*}
de sorte que $f$ induit un automorphisme
polynomial de degr\'e $3$ dans la carte $\mathbb{C}^2=(x_2~=1).$

\begin{pro}\label{1dte}
Soit $f$ une transformation birationnelle purement cubique
con\-tractant exactement une droite; $f$ est
gauche-droite conjugu\'ee \`a la transformation
\begin{align*}
& (x_0x_2^2+x_1^3:x_1x_2^2:x_2^3), && \mathsf{\{1\}}.
\end{align*}
En particulier ces transformations forment une seule
orbite sous l'action gauche-droite.
\end{pro}

La d\'emonstration de cette proposition utilise le:

\begin{lem}[\cite{MM}]\label{mm}
{\sl Soit $f$ un germe de fonction holomorphe \`a l'origine de
$\mathbb{C}^n;$ \'ecrivons $f$ sous la forme $f_1^{n_1}\ldots f_p^{n_p}.$
On suppose que les $f_i$ sont irr\'eductibles et que les $n_i$
sont premiers entre eux, {\it i.e.} $f$ n'est pas une puissance.
Soit $g$ un germe de fonction holomorphe satisfai\-sant~$df\wedge
dg=0.$ Il existe un germe de fonction holomorphe $\ell$ \`a
l'origine de $\mathbb{C}$ telle que~$g=~\ell\circ~f.$}
\end{lem}

\begin{rem}
Si $f$ et $g$ sont de plus des polyn\^omes homog\`enes, sous les
hypoth\`eses du Lemme \ref{mm} on obtient $g= \varepsilon f^s$
avec $\varepsilon$ dans $\mathbb{C}$ et $s$ dans $\mathbb{N}.$
\end{rem}

\begin{proof}[{\sl D\'emonstration de la Proposition \ref{1dte}}]
D'apr\`es ce qui pr\'ec\`ede on peut supposer que $f$ induit
un automorphisme polynomial dans la carte $x_2=1.$ On
\'ecrit $f$ sous la forme
\begin{align*}
f=(P,Q)=(P_0+P_1+P_2+P_3,Q_0+Q_1+Q_2+Q_3)
\end{align*}
les $P_i$ et $Q_i$
appartenant \`a $\mathbb{C}[x_0,x_1]_i.$ \`A conjugaison
g.d. pr\`es on peut supposer que $P_0=Q_0=0$ et que $P_3$ est non
identiquement nul. Comme $f$ est un automorphisme polynomial, il
existe un complexe non nul $\eta$ tel que $\mathrm{d}P\wedge \mathrm{d}Q=\eta
\mathrm{d}x_0\wedge \mathrm{d}x_1.$ Si $P_3$ n'est pas une puissance le Lemme~\ref{mm}
assure que $Q_3=\eta P_3;$ \`a conjugaison g.d. pr\`es, on se
ram\`ene \`a $Q_3=0$ puis \`a $Q_1=Q_2=0$ ce qui est impossible. De
sorte que $P_3$ est une puissance et la seule possibilit\'e \`a
conjugaison lin\'eaire pr\`es est $P_3=x_1^3;$ on peut ici encore se
ramener \`a $Q_3=0$ et $Q_2=\varepsilon x_1^2.$

Si $\varepsilon$ est nul, $f$ s'\'ecrit
$(P_1+P_2+x_1^3,Q_1)$ et n\'ecessairement $Q_1= \mu x_1$ ce qui fait
qu'\`a conjugaison g.d. pr\`es on se ram\`ene \`a $(x_0+x_1^3,x_1).$

Reste \`a s'assurer que $\varepsilon=0;$ si ce n'est pas
le cas on peut supposer que $\varepsilon$ vaut $1.$ Toujours en
exploitant $\mathrm{d}P\wedge \mathrm{d}Q=\eta \mathrm{d}x_0\wedge \mathrm{d}x_1,$ il vient $\mathrm{d}P_2 \wedge
\mathrm{d}x_1^2+\mathrm{d}x_1^3\wedge \mathrm{d}Q_1=0$ qui conduit \`a 
\begin{align*}
\mathrm{d}(2 P_2-3x_1Q_1)\wedge\mathrm{d}x_1=0
\end{align*}
de sorte que $P_2=\frac{3x_1Q_1}{2}+\beta x_1^2.$ \`A
conjugaison g.d. pr\`es $f$ est du type
\begin{align*}
\left(P_1+\frac{3x_1Q_1}{2}+\beta x_1^3,Q_1+ x_1^2\right)
\end{align*}
avec la condition $(3Q_1^2/4-2P_1x_1)\wedge \mathrm{d}x_1=0;$ celle-ci implique que
$Q_1$ est divisible par $x_1$ et emp\^eche que $f$ soit un
automorphisme.
\end{proof}

\subsubsection{Le lieu exceptionnel est constitu\'e de
deux droites}\hspace{1mm}

Soit $f$ un \'el\'ement de $\mathring{\mathrm{B}}\mathrm{ir}_3$ qui contracte
exactement deux droites. Comme $\mathrm{Exc}\hspace{1mm} f$ a deux
composantes, il en est de m\^eme pour $\mathrm{Exc}\hspace{1mm} f^{-1}.$ Les
arguments topologiques pr\'ec\'edents nous indiquent que
$\mathrm{Exc}\hspace{1mm} f^{-1}$ est constitu\'e soit de deux droites
distinctes, soit d'une conique lisse et d'une droite tangente \`a
la conique. Rappelons que le compl\'ement de ces deux
types de courbes dans $\mathbb{P}^2(\mathbb{C})$ sont biholomorphes.  Quoiqu'il en
soit on peut supposer, \`a conjugaison g.d. pr\`es, que
$\mathrm{Exc}\hspace{1mm} f=\{x_0=0,\hspace{1mm} x_2=0\}$ et que 
$\mathrm{Exc}\hspace{1mm} f^{-1}$
est de l'un des deux types suivants
\begin{itemize}
\item $\mathrm{Exc}\hspace{1mm} f^{-1}=\{x_0=0,\hspace{1mm} x_2=0\};$

\item $\mathrm{Exc}\hspace{1mm} f^{-1}=\{x_2=0,\hspace{1mm} x_1^2-x_0x_2=0\}.$
\end{itemize}

Traitons la premi\`ere possibilit\'e; $f$ induit un
biholomorphisme de $\mathbb{C}^2\setminus\{x_0=0\}$ dans lui-m\^eme. Si on
\'ecrit $f=(f_0,f_1)$ dans cette carte, on note que $f_0$ et $f_1$
n'ont pas de p\^ole en dehors de $x_0= 0$ et que $f_0$ n'y a pas de
z\'ero. Par suite il existe $p,$ $q$ dans $\mathbb{Z}$ et $g_1$ dans
$\mathbb{C}[x_0,x_1]$ tels que
\begin{align*}
&f_0(x_0,x_1)=x_0^p, && f_1(x_0,x_1)=x_0^qg_1(x_0,x_1).
\end{align*}
Mais le d\'eterminant jacobien de $f$ est $px_0^{p-1}\frac{\partial
g_1}{\partial x_1}(x_0, x_1);$ comme $f$ est un automorphisme de $\mathbb{C}^2
\setminus \{x_0=0\},$ on a $\frac{\partial g_1}{\partial x_1}(x_0,x_1)=
ax_0^n.$ Ainsi \`a conjugaison lin\'eaire pr\`es
\begin{align*}
& f=(x_0^p,x_0^{n+q}x_1+x_0^q\psi(x_0)), && \psi\in\mathbb{C}[x_0].
\end{align*}
Visiblement $p$
vaut $\pm 1~;$ il nous reste \`a lister $q,$ $n+q$ et $\psi$ pour
lesquels ce type de transformations est de degr\'e $3$ exactement.
On obtient, \`a conjugaison g.d. pr\`es, les mod\`eles
suivants
\begin{align*}
& (x_0x_2^2:x_0^2x_1:x_2^3), && (x_0x_2^2:x_0^3+x_0x_1x_2:x_2^3), \\
& (x_0^2x_2:x_0^3+x_2^3+x_0x_1x_2:x_0x_2^2), && (x_0^2x_2:x_0^2x_1+x_2^3:x_0x_2^2).
\end{align*}

Consid\'erons pour finir la seconde \'eventualit\'e:
$\mathrm{Exc}\hspace{1mm} f^{-1}=\{x_1=0,\hspace{1mm} x_1^2-x_0x_2=0\}.$ La
Proposition \ref{condte} donne la description de $f^{-1}$
\begin{align*}
(x_1^3:x_0(\alpha x_0^2+\beta x_0x_1+\gamma x_1^2+\delta x_1x_2):x_1(
\alpha x_0^2+\beta x_0x_1+\gamma x_1^2+\delta x_1x_2))
\end{align*}
soit \`a
conjugaison g.d. pr\`es $(x_0(x_0^2+x_1x_2):x_1^3:x_1(x_0^2+x_1x_2));$ l'inverse de
cette transformation est donn\'e par $(x_0x_1x_2:x_1x_2^2:x_2^3-x_0^2x_1).$

Il s'en suit la:

\begin{pro}\label{2dtes}
{\sl Soit $f$ une transformation birationnelle purement cubique
qui contracte exactement deux droites ($\mathsf{\{2\}}$). Alors
$f$ est \`a conjugaison gauche-droite pr\`es de l'une des formes
suivantes
\begin{align*}
& (x_0x_2^2:x_0^2x_1:x_2^3); && (x_0x_2^2:x_0^3+x_0x_1x_2:x_2^3); && (x_0^2x_2:x_0^3
+x_2^3+x_0x_1x_2:x_0x_2^2); \\
& (x_0^2x_2:x_0^2x_1+x_2^3:x_0x_2^2); && (x_0x_1x_2:x_1x_2^2:x_2^3-x_0^2x_1). &&
\end{align*}}
\end{pro}

\subsubsection{Le lieu exceptionnel est constitu\'e de trois
droites}\label{33dtes}\hspace{1mm}

On proc\`ede plus ou moins comme pr\'ec\'edemment. Soit
$f$ un \'el\'ement de $\mathring{\mathrm{B}}\mathrm{ir}_3$ qui contracte exactement trois
droites. Ces trois droites sont
\begin{itemize}
\item ou bien concourantes;

\item ou bien en position g\'en\'erale.
\end{itemize}

Nous allons traiter ces deux \'eventualit\'es l'une apr\`es l'autre.

\begin{itemize}
\item Supposons que $\mathrm{Exc}\hspace{1mm}f$ soit constitu\'e
de trois droites concourantes. Comme $\mathrm{Exc}\hspace{1mm}f$ et
$\mathrm{Exc}\hspace{1mm} f^{-1}$ ont m\^eme nombre de composantes
et ont leurs compl\'ements diff\'eomorphes, on a les deux
possibilit\'es suivantes
\begin{itemize}
\item $\mathrm{Exc}\hspace{1mm}f^{-1}$ est form\'e d'une conique lisse et
de deux droites r\'ealisant la configuration~$\mathsf{\{11\}}.$ En
inversant les mod\`eles r\'ealisant cette configuration ({\sl 10.}
et {\sl 11.} dans la Proposition \ref{listee}) on constate que
cette situation n'arrive pas.

\item $\mathrm{Exc}\hspace{1mm} f^{-1}$ est form\'e de trois droites qui sont
n\'ecessairement concourantes. Dans ce cas on peut supposer, \`a
conjugaison g.d. pr\`es, que $\mathrm{Exc}\hspace{1mm} f=\mathrm{Exc}
\hspace{1mm}f^{-1}.$
\begin{pro}\label{3concourantes}
{\sl Soit $f$ un \'el\'ement de $\mathring{\mathrm{B}}\mathrm{ir}_3$
dont le lieu exceptionnel est constitu\'e de trois droites
concourantes ($\mathsf{\{4\}}$). Alors, \`a conjugaison
gauche-droite pr\`es, $f$ est de l'un des types suivants
\begin{align*}
& \hspace{5mm}&&&(x_0^3:x_0^2x_1:(x_0-x_1)x_1x_2);
&& (x_0^2(x_0-x_1):x_0x_1(x_0-x_1):x_0x_1x_2+x_1^3).
\end{align*}
Dans les deux cas les ensembles $\mathrm{Exc}\hspace{1mm}f^{-1}$
sont aussi form\'es de trois droites concourantes.}
\end{pro}

\begin{proof}[\sl D\'emonstration]
D'apr\`es ce qui pr\'ec\`ede il suffit d'examiner le cas o\`u
\begin{align*}
\mathrm{Exc}\hspace{1mm}f=\mathrm{Exc}\hspace{1mm}f^{-1}=\{x_0x_1(x_0- x_1)=0\}.
\end{align*}
Sous cette hypoth\`ese $f$ respecte la fibration $x_1/x_0=$ cte.
En effet consid\'erons une droite $x_1= \beta x_0$ avec
$\beta\not\in\{0,1,\infty\};$ introduisons l'application
holomorphe
\begin{align*}
\psi\hspace{1mm}\colon\hspace{1mm}\mathbb{C}\to\mathbb{C}^2\setminus\{x_0x_1
(x_0-x_1)=0\}\subset\mathbb{P}^2(\mathbb{C})
\setminus\{x_0x_1(x_0-x_1)=0\}
\end{align*}
d\'efinie par $t\mapsto(e^t,\beta e^t).$
Comme $f$ est par hypoth\`ese un automorphisme
de $\mathbb{P}^2(\mathbb{C})\setminus\{x_0x_1(x_0-x_1)=0\},$
l'application $\varphi=f\circ\psi$ est une application
holomorphe de $\mathbb{C}$ \`a valeurs dans
$\mathbb{P}^2(\mathbb{C})\setminus\{x_0x_1(x_0-x_1)=0\}
\simeq\mathbb{C}\setminus\{0,1\}\times\mathbb{C}.$
Une application directe du th\'eor\`eme de \textsc{Picard} dit que
$\varphi$ est \`a valeurs dans une droite $x_1=\beta'
x_0$ et ceci pour tout $\beta;$ d'o\`u l'affirmation. Par
cons\'equent $f$ s'\'ecrit en coordonn\'ees homog\`enes
\begin{align*}
&f=((cx_0+dx_1)Q:(ax_0+bx_1)Q:C), && ad-bc\not=0,
\end{align*}
avec $Q$ quadratique et $C$ cubique.
Comme $Q$ est contract\'e par $f$ on peut supposer
que~$Q=x_0\ell$ o\`u $\ell$ est l'une des formes
lin\'eaires $x_0$ ou $x_0-x_1.$ Finalement \`a
conjugaison g.d. pr\`es on se ram\`ene \`a
\begin{align*}
&f=(x_0^2\ell:x_0x_1\ell:C), && \deg C=3.
\end{align*}
Mais le calcul explicite de $\det(\mathrm{jac}\hspace{1mm} f)$ conduit
\`a
\begin{align*}
C=x_2C_1(x_0,x_1)+C_2(x_0,x_1).
\end{align*}
En \'ecrivant que $x_1=0$ est contract\'ee on constate que
$C_1(x_0,0)=0$ et
\begin{align*}
f=(x_0^2\ell:x_0x_1\ell:x_1x_2L(x_0,x_2)+C_2(x_0,x_1))
\end{align*}
avec $L$ lin\'eaire et $\ell=x_0$ ou $\ell=x_0-x_1.$ Un
tel $f$ contracte effectivement les droites~$x_0=~0$ et
$x_1=0.$

Maintenant si $\ell=x_0,$ la droite $x_0=x_1$
sera contract\'ee si et seulement si $L(x_0,x_0)=0.$
\`A conjugaison g.d. pr\`es on aura
\begin{align*}
&f=(x_0^3:x_0^2x_1:(x_0-x_1)x_1x_2+C_2(x_0,x_1)),
&& \deg C_2=3
\end{align*}
qui d\'efinit bien une transformation birationnelle cubique
contractant exactement trois droites. On se ram\`ene
alors \`a conjugaison g.d. pr\`es au premier mod\`ele
annonc\'e
\begin{align*}
&(x_0^3:x_0^2x_1:(x_0-x_1)x_1x_2).
\end{align*}

Lorsque $\ell=x_0-x_1$ l'application $f$ est du type
\begin{align*}
f=(x_0^2(x_0-x_1):x_0x_1(x_0-x_1):x_1x_2L(x_0,x_1) +C_2(x_0,x_1)).
\end{align*}
Pour que la droite $x_1=0$ soit contract\'ee il faut que
$L=$ cte $x_0,$ la constante pouvant \^etre
suppos\'ee \'egale \`a $1.$ On se ram\`ene alors
par conjugaison gauche-droite \`a
\begin{align*}
&f=(x_0^2(x_0-x_1):x_0x_1(x_0-x_1):x_0x_1x_2+
x_1^3).
\end{align*}

Les transformations
\begin{align*}
&(x_0^3:x_0^2x_1:(x_0-x_1)x_1x_2)
\end{align*}
et
\begin{align*}
& (x_0^2(x_0-x_1):x_0x_1(x_0-x_1):(x_0-x_1)x_1x_2+x_1^3)
\end{align*}
ne sont pas g.d. conjugu\'ees: contrairement \`a la seconde
la premi\`ere contracte une droite avec multiplicit\'e $4.$
\end{proof}
\end{itemize}

\item Supposons que $\mathrm{Exc}\hspace{1mm}f$ soit form\'e de trois
droites en position g\'en\'erale. Il y a deux possibilit\'es suivant
que $\mathrm{Exc}\hspace{1mm}f^{-1}$ contienne ou non une conique
lisse.

Lorsque $\mathrm{Exc}\hspace{1mm}f^{-1}$ contient une conique lisse
$\mathrm{Exc}\hspace{1mm}f^{-1}$ contient aussi deux droites; on est donc
dans l'une des configurations $\mathsf{\{10\}}$ ou
$\mathsf{\{11\}}$
\begin{figure}[H]
\begin{center}
\input{condte.pstex_t}
\end{center}
\end{figure}

Le compl\'ement de la configuration $\mathsf{\{10\}}$
dans $\mathbb{P}^2(\mathbb{C})$ est diff\'eomorphe au compl\'ement
de trois droites en position g\'en\'erale. Un calcul explicite
montre que les trois mod\`eles $f$ de la Proposition \ref{listee}
pr\'esentant la configuration $\mathsf{\{10\}}$ satisfont:
$\mathrm{Exc}\hspace{1mm}f^{-1}$ est l'union de trois droites en position
g\'en\'erale (configuration $\mathsf{\{3\}}).$ Pour ce qui est de
la possibilit\'e $\mathsf{\{11\}}$ ou bien on remarque que le
compl\'ement de $\mathsf{\{11\}}$ n'est pas diff\'eomorphe au
compl\'ement de~$\mathsf{\{3\}}$ (les groupes fondamentaux sont
diff\'erents), ou bien on inverse les deux mod\`e\-les~$f$ de la
Proposition \ref{listee} pr\'esentant la configuration
$\mathsf{\{11\}}$ pour constater que $\mathrm{Exc}\hspace{1mm}f^{-1}$ est
encore de type $\mathsf{\{11\}}.$ Cette derni\`ere situation est
donc \`a exclure.

\begin{pro}\label{3dtes}
{\sl Soit $f$ un \'el\'ement de $\mathring{\mathrm{B}}\mathrm{ir}_3$ qui contracte trois
droites en position g\'en\'erale ($\mathsf{\{3\}}$). Si
$\mathrm{Exc}\hspace{1mm}f^{-1}$ contient une conique lisse, il
contient deux droites dont l'une est tangente \`a la conique; de
plus $f$ est \`a conjugaison gauche-droite pr\`es de l'une des
formes suivantes
\begin{align*}
&(x_0^2(x_1-x_2):x_0x_1(x_1-x_2):x_1^2x_2); && (x_0^2x_2:x_0x_1x_2:x_0x_1^2-x_1^2x_2); \\
&(x_0x_1x_2:x_1^2x_2:x_0(x_1^2-x_0x_2)).&&
\end{align*}}
\end{pro}

Le cas restant est celui o\`u $\mathrm{Exc}\hspace{1mm}f^{-1}$ est lui aussi form\'e
de trois droites en position g\'en\'erale. \`A conjugaison gauche-droite
pr\`es nous pouvons supposer que
\begin{align*}
\mathrm{Exc}\hspace{1mm}f=\mathrm{Exc}\hspace{1mm}f^{-1}=\{x_0x_1x_2=0\}.
\end{align*}
Dans la carte
affine $\mathbb{C}^2=(x_2 =1)$ la transformation $f$ induit un automorphisme
du compl\'ement de $x_0x_1=0.$ On en d\'eduit que $f$ s'\'ecrit
$(x_0^px_1^q,x_0^rx_1^s)$ avec $\left[\begin{array}{cc} p & q \\ r & s
\end{array}\right]\in\mathrm{GL}_2(\mathbb{Z}).$ Lorsque l'on homog\'en\'eise $f$ on trouve une
expression du type
\begin{align*}
(x_0^{p_0}x_1^{q_0}x_2^{r_0}:x_0^{p_1}x_1^{q_1}x_2^{r_1}:x_0^{p_2}x_1^{q_2}
x_2^{r_2})
\end{align*}
o\`u les sommes $p_i+q_i+ r_i$ valent $3$ et les
produits $p_0p_1p_2,$ $q_0q_1q_2,$ $r_0r_1r_2$ sont nuls. Le
d\'eterminant jacobien de $f$ vaut, \`a multiplication pr\`es par
une constante,
\begin{align*}
x_0^{(\sum p_i)-1}x_1^{(\sum q_i)-1}x_2^{(\sum r_i)-1};
\end{align*}
puisque $\mathrm{Exc}\hspace{1mm}f$ est constitu\'e de trois droites il faut
que
\begin{align*}
&\sum p_i>1, && \sum q_i>1,&& \sum r_i>1.
\end{align*}
Comme la transformation $f$
s'\'ecrit $(x_0^{p_0-p_2}x_1^{q_0-q_2},x_0^{p_1-p_2}x_1^{q_1-q_2})$ dans
la carte affine $x_2=~1,$ on a de plus
\begin{align*}
\det\left[\begin{array}{cc}
p_0-p_2 & q_0-q_2 \\
p_1-p_2 & q_1-q_2
\end{array}
\right]=\pm 1.
\end{align*}
L'examen de toutes les \'eventualit\'es nous conduit \`a la:

\begin{pro}\label{3dtesb}
{\sl Soit $f$ une transformation birationnelle purement cubique
telle que les lieux exceptionnels de $f$ et $f^{-1}$ soient
form\'es de trois droites en position g\'en\'erale~($\mathsf{\{3\}}$).
Alors $f$ est, \`a conjugaison gauche-droite pr\`es, du type suivant
\begin{align*}
(x_0^3:x_1^2x_2:x_0x_1x_2).
\end{align*}}
\end{pro}
\end{itemize}

\subsubsection{Le lieu exceptionnel est form\'e de quatre droites ou
plus}\hspace{1mm}

\begin{lem}\label{ventru}
{\sl Soit $f$ un \'el\'ement de $\mathring{\mathrm{B}}\mathrm{ir}_3.$ On suppose qu'\`a
conjugaison gauche-droite pr\`es~$f$ contracte les droites $x_0=0,$
$x_1=0$ et $x_2=0$ sur trois points distincts $p_1,$ $p_2$ et $p_3;$
on suppose en outre que $f$ contracte une quatri\`eme droite
$\ell=0.$ Alors les $p_i$ ne sont pas align\'es.}
\end{lem}

\begin{proof}[\sl D\'emonstration]
Supposons que $f$ contracte $x_0=0,$ $x_1=0,$ $x_2=0$ sur $(1:0:0),$
resp. $(1:0:1),$ resp. $(0:0:1).$ Il existe $q_0$ et $q_2$ dans $\mathbb{C}[x_0,x_1,
x_2]_2$ tels que
\begin{align*}
& f=(x_2q_0:x_0x_1x_2:x_0q_2) && \text{avec } && x_2q_0(x_0,0,x_2)=x_0q_2(x_0,0,x_2).
\end{align*}
La droite $\ell=0$ est
contract\'ee sur un point $p$ qui ne peut \^etre align\'e avec les
trois pr\'ec\'edents (sinon $x_0x_1x_2$ serait divisible par $\ell$).
Par suite on peut supposer que $p=(0:1:0);$ la transformation $f$
est donc de la forme $(x_2\ell\ell_0:x_0x_1x_2:x_0\ell\ell_2),$ les $\ell_i$
\'etant lin\'eaires. Mais en se souvenant que $x_1=0$ est
contract\'ee sur $(1:0:1),$ on constate que $\ell_0$ et $\ell_2$
sont du type suivant
\begin{align*}
&\ell_0=\alpha x_0+\beta x_1,&& \ell_2= \alpha x_2+\gamma x_1,
&&\alpha\beta\gamma \not=0.
\end{align*}
On se ram\`ene alors par homoth\'etie \`a $\alpha=\beta= \gamma=1$
et
\begin{align*}
f=(x_2(x_0+x_1)(Ax_0+Bx_1+Cx_2): x_0x_1x_2:x_0(x_1+x_2)(Ax_0+Bx_1+Cx_2))
\end{align*}
qui dans la carte affine $x_1=1$ s'\'ecrit
\begin{align*}
\left(\frac{x_0+1}{x_0}(Ax_0+B+Cx_2),
\frac{x_2+1}{x_2}(Ax_0+B+Cx_2)\right).
\end{align*}
En composant par $\psi=\left(\frac{1}{x_0-1},\frac{1}{x_2-1}\right)$ on obtient
\begin{align*}
f\circ\psi=\left(x_0\left(\frac{A}{x_0-1}+B+\frac{C}{x_2-1}\right),x_2\left(
\frac{A}{x_0-1}+B+\frac{C}{x_2-1}\right)\right)
\end{align*}
qui est injective en m\^eme temps que
\begin{align*}
\left(x_0\left(\frac{A}{x_0-1}+B+\frac{C}{x_2-1}\right),\frac{x_2}{x_0}\right)
\end{align*}
Si cette derni\`ere application est injective g\'en\'eriquement
alors $A=C=0$ ce qui implique que~$f$ est quadratique: contradiction.
\end{proof}

Avec les m\^emes techniques on en d\'eduit la:

\begin{pro}
{\sl Une transformation birationnelle purement cubique $f$ ne peut
pas contracter quatre droites en position g\'en\'erale.}
\end{pro}

\begin{proof}[\sl D\'emonstration]
Raisonnons par l'absurde: supposons que $f$ contracte les droites
\begin{align*}
&x_0=0, && x_1=0, && x_2=0 &&  \text{et} &&
x_2=x_0+x_1
\end{align*}
sur $p_0,$ $p_1,$ $p_2$ et $p_3$ respectivement.

Si les $p_i$ sont distincts, on peut se ramener, d'apr\`es le
Lemme \ref{ventru}, au cas o\`u
\begin{align*}
& p_0=(1:0:0), && p_1=(0:1:0),&& p_2=(0:0:1),&& p_3=(1:1:1).
\end{align*}
Un calcul direct montre que $f$ est du type suivant

\begin{small}
\begin{align*}
(x_1x_2(a(x_0+x_1-x_2)+x_2-x_1):x_0x_2(b(x_0+x_1-x_2)+x_2-x_0)
:x_0x_1(c(x_0+x_1-x_2)+x_2)).
\end{align*}
\end{small}

On remarque que si l'une des conditions $a=1,$ $b=1,$ $c=0$ est
satisfaite alors $f$ est quadratique. En carte affine $x_2=1,$ la
transformation $f$ s'\'ecrit
\begin{align*}
\left(\frac{a(x_0+x_1-1)+1-x_1}{
x_0(c(x_0+x_1-1)+1)},\frac{b(x_0+x_1-1)+1-x_0}{x_1(c(x_0+x_1-1)+1)}\right).
\end{align*}
On doit exprimer l'injectivit\'e de $f;$ pour cel\`a on regarde les
\'equations $f=(u,v)$ qui conduisent~\`a
\begin{align*}
&a(x_0+x_1-1)+1-x_1=ux_0(c(x_0+x_1-1)+1), \\
& b(x_0+x_1-1)+1-x_0=vx_1(c(x_0+x_1-1)+1).
\end{align*}

Apr\`es \'elimination on obtient que
\begin{align*}
x_1=-\frac{cux_0^2+(u(1-c)-a)x_0+(a-1)}{ucx_0+1-a}
\end{align*}
et
\begin{small}
\begin{eqnarray}\label{injj}
 & &c^2u(v-u)x_0^3+c(u^2c^2+2u(a-1)-av+uv(2-a-c))x_0^2\\
& &+(a(2-a-b)+b-1+u(2c-b)(1-a)-v(a^2+c+a-ac)+uv(a-1)(c-1))x_0\nonumber\\
& & +(a-1)^2(1-v)=0.\nonumber
\end{eqnarray}
\end{small}
Pour qu'il n'y ait qu'une solution g\'en\'erique, cette \'equation
doit se factoriser de la fa\c{c}on suivante
\begin{align*}
(x_0-x_0(u,v))(x_0-\eta_1)(x_0-\eta_2)(x_0-\eta_3)=0
\end{align*}
les $\eta_i$ \'etant des
constantes. Soit $x_0=\eta$ une solution de (\ref{injj})
ind\'ependante de $u,$ $v.$ Le coefficient de $u^2$ est
$\eta^2c^2(c-1);$ comme $c$ est non nul, $\eta^2c^2(c-1)=0$ si et
seulement si $c=1$ ou~$\eta=0.$ Si $\eta\not=0,$ alors $c=1$ et on
constate que $\eta=a=1;$ dans ce cas $f$ est quadratique.
Si~$\eta=0;$ \`a partir de $(a-1)^2(1-v)=0$ on obtient que
$a$ vaut $1$ et que $f$ est quadratique. \bigskip

Supposons que les $p_i$ ne soient pas distincts;
d'apr\`es le Lemme \ref{pasdeux} il y a au moins parmi ces points
trois points distincts. En particulier on peut invoquer la
Proposition \ref{listee}; on constate que l'arrangement de quatre
droites contract\'ees en position g\'en\'erale n'appara\^it pas
dans celle-~ci.
\end{proof}

Les configurations de quatre droites contract\'ees (au moins)
contiennent donc au moins trois droites concourantes. Comme on l'a
vu si les quatre droites sont contract\'ees par $f$ sur
trois points alors $f$ est \`a \'equivalence pr\`es donn\'ee par
la Proposition \ref{listee}. On peut donc supposer que les quatre
droites sont contract\'ees sur quatre points distincts et que
trois d'entre elles sont les droites $x_0=0,$ $x_1=0$ et $x_1=x_0.$ On
consid\`ere tour \`a tour les cas o\`u les droites $x_0=0,$ $x_1=0$ et
$x_0=x_1$ sont contract\'ees sur trois points align\'es,
respectivement non align\'es et distincts.

\begin{pro}\label{5dtes}
{\sl Soit $f$ un \'el\'ement de $\mathring{\mathrm{B}}\mathrm{ir}_3$ qui, \`a conjugaison
gauche-droite pr\`es, contracte les droites $x_0=0,$ resp. $x_1=0,$
resp. $x_1=x_0$ sur trois points $(1:0:0),$ resp. $(0:1:0),$ resp.
$(1:1:0)$ et (au moins) une quatri\`eme droite. Alors $f$
contracte cinq droites et est, \`a conjugaison gauche-droite
pr\`es, le mod\`ele {\sl 5.} de la Proposition~\ref{listee}.}
\end{pro}

\begin{proof}[\sl D\'emonstration]
La transformation $f$ est de la forme
\begin{align*}
(x_1q+(x_0-x_1)x_1\ell:x_0q:x_0x_1(x_0-x_1))
\end{align*}
o\`u $q$ d\'esigne une forme quadratique et $\ell$ une
forme lin\'eaire. Supposons que $f$ contracte une quatri\`eme
droite d'\'equation $L=0.$ On est ramen\'e \`a conjugaison pr\`es
\`a l'alternative suivante
\begin{align*}
&L=x_1-\beta x_0 \text{ avec } \beta\not\in\{0,1\}, && \text{ou}
&&L=x_2.
\end{align*}

Nous allons voir que n\'ecessairement $L=x_2;$ raisonnons
par l'absurde {\it i.e.} supposons que $L=x_1-\beta x_0$ avec $\beta
\not\in\{0,1\}.$ Sur $x_1-\beta x_0=0$ on a
\begin{align*}
(\beta x_0q(x_0,\beta x_0,x_2)+\beta(1-\beta)x_0^2\ell(x_0,\beta x_0,x_2):
x_0q(x_0,\beta x_0,x_2):\beta(1-\beta)x_0^3);
\end{align*}
ainsi $q(x_0,\beta x_0,x_2)$ est divisible par $x_0^2$ et $\ell(x_0,
\beta x_0,x_1)$ est divisible par $x_0.$ Il s'en suit que $q$ et $\ell$
ne d\'ependent pas de $x_2$ ce qui est impossible.

On suppose donc que $x_2=0$ est contract\'ee; sur $x_2=0$ on a
\begin{align*}
(x_1q(x_0,x_1,0)+x_1(x_0-x_1)\ell(x_0,x_1,0):x_0
q(x_0,x_1,0):x_0x_1(x_0-x_1));
\end{align*}
il en r\'esulte que $q(x_0,x_1,0)=\varepsilon
x_1(x_0-x_1)$ et $\ell(x_0,x_1,0)=ax_0-\varepsilon x_1.$ Ainsi
toujours \`a conjugaison g.d. pr\`es
\begin{align*}
f=(x_1x_2(Ax_0+Bx_1+Cx_2+c(x_0-x_1)):x_0x_2(Ax_0+Bx_1+Cx_2):
x_0x_1(x_0-x_1)).
\end{align*}
Ce mod\`ele est birationnel si $C=Ac=0.$ Or si $C=A=0,$ la transformation est
quadratique donc $C=c=0$ et
\begin{align*}
f=(x_1x_2(Ax_0+Bx_1):x_0x_2(Ax_0+Bx_1):x_0x_1(x_0-x_1))
\end{align*}
qui est le mod\`ele annonc\'e.
\end{proof}

On termine en supposant que les droites $x_0=0,$ $x_1=0$ et $x_1=x_0$ sont
contract\'ees respectivement sur $(1:0:0),$ $(0:1:0)$ et
$(1:0:1).$

\begin{lem}
{\sl Soit $f$ un \'el\'ement de $\mathring{\mathrm{B}}\mathrm{ir}_3.$ Supposons qu'\`a
conjugaison gauche-droite pr\`es $f$ contracte la droite $x_0=0,$
resp. $x_1=0,$ resp. $x_1=x_0$ sur $(1:0:0),$ resp. $(0:1:0),$ resp
$(0:0:1)$ et (au moins) une quatri\`eme droite de la forme $x_1=\eta
x_0$ avec $\eta \not\in\{0,1\}.$ Alors $f$ contracte une conique
r\'eduite; en particulier les hypoth\`eses de la Proposition
\ref{listee} sont v\'erifi\'ees.}
\end{lem}

\begin{proof}[\sl D\'emonstration]
La transformation $f$ est du type
\begin{align*}
& (x_1(x_0-x_1)\ell_0:x_0(x_0-x_1)\ell_1:
x_0x_1\ell_2)
\end{align*}
o\`u les $\ell_i$ sont lin\'eaires. Comme les trois
droites jouent le m\^eme r\^ole et que l'un des $\ell_i$ doit
d\'ependre effectivement de $x_2,$ on peut supposer que $\ell_2=x_2.$
En \'ecrivant que la quatri\`eme droite est contract\'ee on
obtient que $f$ est du type suivant
\begin{align*}
f=(x_1(x_0-x_1)(a(x_1-\eta x_0)+bx_2):x_0(x_0-x_1) (\alpha(x_1-\eta x_0)+\beta
x_2):x_0x_1x_2).
\end{align*}
Un calcul \'el\'ementaire montre que le d\'eterminant
jacobien de $f$ s'\'ecrit
\begin{align*}
\mathrm{det}(\mathrm{jac}\hspace{1mm} f)=3x_0x_1(x_0-x_1)(x_1-\eta x_0)\left(a
\alpha(\eta x_0^2+x_1^2-(\eta+1)x_0x_1)-a\beta x_0x_2+\alpha bx_1x_2\right).
\end{align*}
Notons $\mathcal{C}$ la conique donn\'ee par $a \alpha(\eta
x_0^2+x_1^2-(\eta+1)x_0x_1)-a\beta x_0x_2+\alpha bx_1x_2=0.$ On remarque que
$\mathcal{C}$ est lisse si et seulement si
\begin{align*}
a\alpha(a\beta-b\eta\alpha)(b\alpha-a\beta)\not=0.
\end{align*}
Comme cette conique est contract\'ee le Lemme est d\'emontr\'e dans ce cas.

\bigskip

Si $a$ est nul, $\mathcal{C}$ est l'union des droites
$x_1=0$ et $x_2=0$ qui sont toutes deux contract\'ees sur $(0:1:0):$
la conique $\mathcal{C}$ est r\'eduite et contract\'ee et $f$ satisfait
donc les
hypoth\`eses de la Proposition \ref{listee}. On suppose $a\not=0;$
lorsque $\alpha$ est nul alors $\beta$ ne l'est pas et on constate
de m\^eme que $\mathcal{C}$ est contract\'ee et r\'eduite.

\bigskip

Reste \`a \'etudier l'\'eventualit\'e o\`u $(a
\beta- b\eta\alpha)(b\alpha-a\beta)=0$ avec $a\alpha\not=0.$

Si $b=0$ alors $\beta=0;$ la transformation est alors
g.d. conjugu\'ee \`a
\begin{align*}
(x_1(x_0-x_1)(x_1-\eta x_0): x_0(x_0-x_1)(x_1-\eta x_0):x_0x_1x_2).
\end{align*}
La droite $x_1-\eta x_0$ est contract\'ee sur $(0:0:1),$ comme la
droite $x_1=x_0;$ ici encore le Lemme est d\'emontr\'e et on se
trouve dans la configuration $\mathsf{\{6\}}:$ le lieu
exceptionnel est constitu\'e des quatre droites $x_0=0,$ $x_1=0,$
$x_0=x_1$ et $x_1=\eta x_0.$ On retrouve le mod\`ele {\sl 4.} de
la Proposition~\ref{listee}.

Supposons d\'esormais que $b$ soit non nul;
$\mathcal{C}$ est l'union des deux droites $\mathcal{D}_1,$ resp.
$\mathcal{D}_2$ d'\'equation
\begin{align*}
& x_1=\frac{a\beta}{\alpha b}x_0, && \text{resp. } && \alpha(\alpha
b\eta-a\beta\eta-a\beta)x_1-a\beta^2x_2+a\alpha\beta
\eta x_0=0.
\end{align*}
La droite $\mathcal{D}_1$ est contract\'ee par $f$ si
et seulement si $\eta=\frac{a\beta}{\alpha b}$ ou $a\beta-\alpha
b= 0.$ Commen\c{c}ons par supposer que $\eta=\frac{a\beta}{\alpha
b}.$ On constate que la droite $\mathcal{D}_1$ co\"{\i}ncide avec
la droite $x_1=\eta x_0$ et qu'elle est contract\'ee sur $\left(
a(a\beta-\alpha b):\alpha(a\beta-\alpha b):-a\alpha\right);$ la
droite $\mathcal{D}_2,$ donn\'ee par $ax_1+bx_2-ax_0=0,$ est aussi
contract\'ee sur ce point. La conique $\mathcal{C}$ est donc
r\'eduite et contract\'ee. Si
$a\beta-\alpha b=0,$ on remarque que $\mathcal{D}_1$ co\"{\i}ncide
avec la droite $x_1=x_0$ et que $\mathcal{D}_2,$ d'\'equation
$\alpha\eta x_0-\alpha x_1-ax_2=0,$ est contract\'ee si et seulement si
$a=\beta$ auquel cas elle l'est sur le point $(0:0:1):$ la conique
$\mathcal{C}$ est encore r\'eduite et contract\'ee.
\end{proof}

\begin{pro}\label{aumoins4}
{\sl Soit $f$ une transformation birationnelle purement cubique.
Supposons qu'\`a conjugaison gauche-droite pr\`es $f$ contracte
les droites d'\'equation $x_0=0,$ $x_1=0,$ $x_0=x_1$ sur $(1:0:0),$ resp.
$(0:1:0),$ resp. $(0:0:1)$ et la droite $x_2=0.$ Alors $f$ appartient
\`a l'orbite sous l'action gauche-droite d'une des transformations suivantes
\begin{itemize}
\item mod\`ele {\sl 2.} de la Proposition
\ref{listee},\hspace{10mm}$\mathsf{\{5\}};$

\item
$(x_1(x_0-x_1)(x_0+x_2):x_0(x_0-x_1)(x_2-x_1):x_1x_2(x_0+x_1)),\hspace{10mm}\mathsf{\{7\}}.$
\end{itemize}}
\end{pro}

\begin{rem}
Comme nous le verrons dans la d\'emonstration,
\begin{align*}
&(x_1(x_0-x_1)(x_0+x_2):x_0(x_0-x_1)(x_2-x_1):x_1x_2(x_0+x_1))
\end{align*}
n'est pas g.d. conjugu\'e aux mod\`eles de la Proposition
\ref{listee} pr\'esentant la configuration~$\mathsf{\{7\}}.$
\end{rem}

\begin{proof}[\sl D\'emonstration]
D'apr\`es les hypoth\`eses $f$ est de la forme
\begin{align*}
(x_1(x_0-x_1) (\alpha x_0+\beta x_2):x_0(x_0-x_1)(\gamma x_1+\delta
x_2):x_0x_1(\varepsilon(x_0-x_1) +\mu x_2)).
\end{align*}
Notons que $\beta,$ $\delta$ et $\mu$ ne peuvent \^etre
nuls puisque $f$ est purement cubique; on peut donc les supposer
\'egaux \`a $1.$ On remarque qu'une transformation du type
ci-dessus est birationnelle si et seulement si nous sommes dans
l'une des situations suivantes
\begin{itemize}
\item $\gamma=\varepsilon=0;$

\item $\alpha=\gamma=0;$

\item $\gamma=-1,$ $\varepsilon=\frac{\alpha}{\alpha+1}.$
\end{itemize}
Dans les deux premiers cas on obtient \`a conjugaison g.d. pr\`es
\begin{align*}
& f_1=(x_1(x_0-x_1)(x_0+x_2):x_0x_2(x_0-x_1):x_0x_1x_2), &&
\mathsf{\{5\}};
\end{align*}
la derni\`ere \'eventualit\'e conduit \`a
\begin{align*}
& (x_1(x_0-x_1)(\alpha x_0+x_2):x_0(x_0-x_1)(x_2-x_1):x_1x_2(\alpha
x_0+x_1))
\end{align*}
soit \`a
\begin{align*}
& f_2=(x_1x_2(x_0-x_1):x_0(x_0-x_1)(x_2-x_1):x_1^2x_2), &&
\mathsf{\{5\}}
\end{align*}
ou \`a
\begin{align*}
&
f_3=(x_1(x_0-x_1)(x_0+x_2):x_0(x_0-x_1)(x_2-x_1):x_1x_2(x_0+x_1)),
&& \mathsf{\{7\}}
\end{align*}
suivant que $\alpha$ est nul ou non.

Les \'el\'ements $f_1$ et $f_2$ sont g.d. conjugu\'ees
\begin{align*}
& f_2=(-x_2:-x_0:x_1+x_2)f_1(-x_1:x_0-x_1:x_2).
\end{align*}

La transformation $f_1$ est g.d. conjugu\'ee au
mod\`ele {\sl 2.} de la Proposition \ref{listee} que nous noterons
$g$
\begin{align*}
& f_1=(x_2:x_0:x_1)g(x_1-x_0:-x_1:x_0+x_2).
\end{align*}

L'\'el\'ement
\begin{align*}
f_3=(x_1(x_0-x_1)(x_0+x_2):x_0(x_0-x_1)(x_2-x_1):x_1x_2(x_0+x_1))
\end{align*}
n'est pas g.d. conjugu\'ee au
mod\`ele {\sl 6.} de la Proposition \ref{listee}: sur la composante de
$\mathrm{det}\hspace{1mm}\mathrm{jac}$  de multiplicit\'e $2$
chacun de ces deux \'el\'ements admet deux points d'ind\'etermination
dont la position diff\`ere (dans le cas de $f_3$ ces points ne sont pas
n\'ecessairement des points d'intersection avec d'autres courbes de
$\mathrm{Exc}\hspace{1mm}f_3$ contrairement \`a ce qu'il se passe pour le mod\`ele {\sl 6.}).
Elle n'est pas non plus g.d. conjugu\'ee au mod\`ele {\sl 5.} de cette
m\^eme Proposition; en effet le mod\`ele~{\sl 5.} poss\`ede la
propri\'et\'e suivante que $f_3$ n'a pas:
on peut trouver dans le lieu exceptionnel de cette transformation
trois courbes contract\'ees sur trois points distincts align\'es.
Voici des dessins qui illustrent ce discours;
on adopte les notations suivantes. Les droites correspondent aux courbes
contract\'ees par la transformation consid\'er\'ee; quand la multiplicit\'e
(dans le d\'eterminant jacobien) de l'une d'elle est strictement sup\'erieure
\`a $1$ elle est pr\'ecis\'ee entre parenth\`ese. Les \og cercles\fg\hspace{1mm}
correspondent aux points sur lesquels les \'el\'ements du lieu
exceptionnel sont contract\'es; on pr\'ecise entre crochets le nombre
de droites contract\'ees sur un point lorsque celui-ci est strictement
sup\'erieur \`a $1.$ Les \og croix rouges\fg\hspace{1mm} caract\'erisent les points
d'ind\'etermination de la transformation \'etudi\'ee.

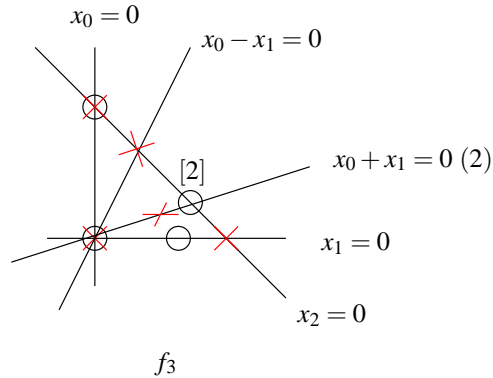
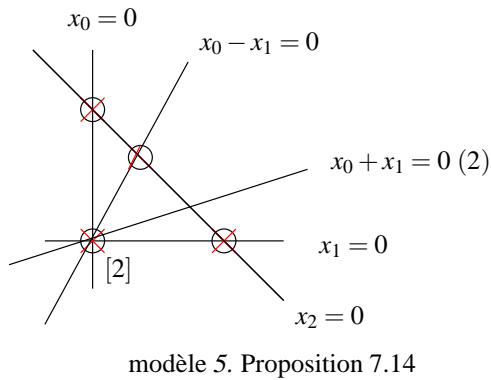
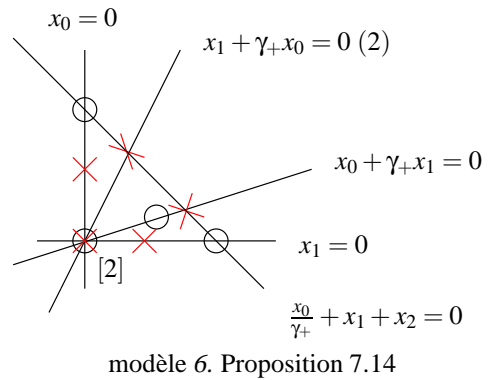
\begin{figure}[H]
\begin{center}
\input{expli4.pstex_t}
\end{center}
\end{figure}

\begin{figure}[H]
\begin{center}
\input{expli5.pstex_t}
\end{center}
\end{figure}

\end{proof}

\begin{rem}
On constate que certaines configurations, par exemple $5$ droites
concourantes ou $6$ droites, n'arrivent pas; on peut d'ailleurs le
prouver de fa\c{c}on directe sans utiliser la classification.
\end{rem}

\subsection{R\'ecapitulatif}\label{dimor}\hspace{1mm}

Le tableau qui suit donne \`a conjugaison g.d. pr\`es
les diff\'erents mod\`eles de transformations birationnelles
purement cubiques; pour chacun d'entre eux on a mentionn\'e la
configuration des courbes contract\'ees (deuxi\`eme colonne) et
la dimension de son orbite sous l'action gauche-droite (troisi\`eme
colonne). Si $\gamma$ d\'esigne un complexe,
on rappelle que
\begin{align*}
\gamma_+:=\frac{\gamma+\sqrt{\gamma^2-4}} {2} && \text{et} &&
\gamma_-:=\frac{\gamma-\sqrt{\gamma^2-4}}{2}
\end{align*}
sont les deux racines du trin\^ome $t^2-\gamma t+1.$

Rappelons qu'on obtient la dimension de l'orbite d'une
transformation birationnelle en d\'eterminant son groupe
d'isotropie. Donnons un exemple.

Consid\'erons la transformation $f=(x_0x_2^2+x_1^3:x_1x_2^2:x_2^3).$
Cherchons $A,$ $B$ dans $\mathrm{SL}_3(\mathbb{C})$ et $\eta$ dans 
$\mathbb{C}^*$ tels que
$Af=\eta fB;$ on obtient que n\'ecessairement
\begin{align*}
&A=(\gamma^3x_0+\alpha\delta^2x_1+\beta\delta^2x_2:\gamma\delta^2x_1:\delta^3x_2),&&
B=\left(\frac{\gamma^3}{\delta^2}x_0+\alpha x_1+\beta x_2:\gamma
x_1:\delta x_2\right)
\end{align*}
avec $\delta^6=1$ et $\gamma^4=\delta.$ La dimension du groupe d'isotropie
est donc $2,$ celle de l'orbite de $f$ sous l'action gauche-droite
$16-2=14.$

\begin{scriptsize}
\begin{landscape}
\begin{center}
\begin{tabular}{|*{3}{c}*{3}{c|}l r|}
   \hline
    & & & & & \\
    & \hspace*{2mm} & $(x_0x_2^2+x_1^3:x_1x_2^2:x_2^3)$  & \hspace*{2mm}&\hspace*{2mm}$\mathsf{\{1\}}$\hspace*{2mm}
    &\hspace*{2mm}14\hspace*{2mm}\\
    & \hspace*{2mm} &  $(x_0x_2^2:x_0^2x_1:x_2^3)$ & \hspace*{2mm}&\hspace*{2mm}$\mathsf{\{2\}}$\hspace*{2mm}
    &\hspace*{2mm}15\hspace*{2mm}\\
    & \hspace*{2mm} & $(x_0x_2^2:x_0^3+x_0x_1x_2:x_2^3)$ & \hspace*{2mm}&\hspace*{2mm}$\mathsf{\{2\}}$\hspace*{2mm}
    &\hspace*{2mm}15\hspace*{2mm}\\
    & \hspace*{2mm} & $(x_0^2x_2:x_0^3+x_2^3+x_0x_1x_2:x_0x_2^2)$ & \hspace*{2mm}&\hspace*{2mm}$\mathsf{\{2\}}$\hspace*{2mm}
    &\hspace*{2mm}14\hspace*{2mm}\\
    & \hspace*{2mm} & $(x_0^2x_2:x_0^2x_1+x_2^3:x_0x_2^2)$ & \hspace*{2mm}&\hspace*{2mm}$\mathsf{\{2\}}$\hspace*{2mm}
    &\hspace*{2mm}15\hspace*{2mm}\\
    & \hspace*{2mm} & $(x_0x_1x_2:x_1x_2^2:x_2^3-x_0^2x_1)$ & \hspace*{2mm}&\hspace*{2mm}$\mathsf{\{2\}}$\hspace*{2mm}
    &\hspace*{2mm}15\hspace*{2mm}\\
    & \hspace*{2mm} & $(x_0^3:x_1^2x_2:x_0x_1x_2)$ & \hspace*{2mm}&\hspace*{2mm}$\mathsf{\{3\}}$\hspace*{2mm}
    &\hspace*{2mm}15\hspace*{2mm}\\
    & \hspace*{2mm} & $(x_0^2(x_1-x_2):x_0x_1(x_1-x_2):x_1^2x_2)$ & \hspace*{2mm}&\hspace*{2mm}$\mathsf{\{3\}}$\hspace*{2mm}
    &\hspace*{2mm}15\hspace*{2mm}\\
    & \hspace*{2mm} & $(x_0^2x_2:x_0x_1x_2:x_1^2(x_0-x_2))$ & \hspace*{2mm}&\hspace*{2mm}$\mathsf{\{3\}}$\hspace*{2mm}
    &\hspace*{2mm}15\hspace*{2mm}\\
    & \hspace*{2mm} & $(x_0x_1x_2:x_1^2x_2:x_0(x_1^2-x_0x_2))$ & \hspace*{2mm}&\hspace*{2mm}$\mathsf{\{3\}}$\hspace*{2mm}
    &\hspace*{2mm}15\hspace*{2mm}\\
    & \hspace*{2mm} & $(x_0^3:x_0^2x_1:(x_0-x_1)x_1x_2)$ & \hspace*{2mm}&\hspace*{2mm}$\mathsf{\{4\}}$\hspace*{2mm}
    &\hspace*{2mm}15\hspace*{2mm}\\
    & \hspace*{2mm} & $(x_0^2(x_0-x_1):x_0x_1(x_0-x_1):x_0x_1x_2+x_1^3)$ & \hspace*{2mm}&\hspace*{2mm}$\mathsf{\{4\}}$\hspace*{2mm}
    &\hspace*{2mm}16\hspace*{2mm}\\
    & \hspace*{2mm} & $(x_0x_2(x_0+x_1):x_1x_2(x_0+x_1):x_0x_1^2)$ & \hspace*{2mm}&\hspace*{2mm}$\mathsf{\{5\}}$\hspace*{2mm}
    &\hspace*{2mm}16\hspace*{2mm}\\
    & \hspace*{2mm} & $(x_0(x_0+x_1)(x_1+x_2):x_1(x_0+x_1)(x_1+x_2):x_0x_1x_2)$ & \hspace*{2mm}&\hspace*{2mm}$\mathsf{\{5\}}$\hspace*{2mm}
    &\hspace*{2mm}16\hspace*{2mm}\\
    & \hspace*{2mm} & $(x_0(x_0+x_1+x_2)(x_0+x_1):x_1(x_0+x_1+x_2)(x_0+x_1):x_0x_1x_2)$ & \hspace*{2mm}&\hspace*{2mm}$\mathsf{\{5\}}$\hspace*{2mm}
    &\hspace*{2mm}16\hspace*{2mm}\\
    & \hspace*{2mm} & $(x_0(x_0^2+x_1^2+\gamma x_0x_1):x_1(x_0^2+x_1^2+\gamma x_0x_1):x_0x_1x_2),$ $\gamma^2\not=4$ & \hspace*{2mm}&\hspace*{2mm}$\mathsf{\{6\}}$\hspace*{2mm}
    &\hspace*{2mm}15\hspace*{2mm}\\
    & \hspace*{2mm} & $(x_0x_2(x_1+x_0):x_1x_2(x_1+x_0):x_0x_1(x_0-x_1))$ & \hspace*{2mm}&\hspace*{2mm}$\mathsf{\{7\}}$\hspace*{2mm}
    &\hspace*{2mm}16\hspace*{2mm}\\
    & \hspace*{2mm} & $(x_0(x_0^2+x_1^2+\gamma x_0x_1+\gamma_+x_0x_2+x_1x_2):x_1(x_0^2+x_1^2+\gamma x_0x_1+\gamma_+x_0x_2+x_1x_2):x_0x_1x_2)$ & \hspace*{2mm}&\hspace*{2mm}$\mathsf{\{7\}}$\hspace*{2mm}
    &\hspace*{2mm}16\hspace*{2mm}\\
    & \hspace*{2mm} & $(x_1(x_0-x_1)(x_0+x_2):x_0(x_0-x_1)(x_2-x_1):x_1x_2(x_0+x_1))$ & &\hspace*{2mm}$\mathsf{\{7\}}$\hspace*{2mm}&\hspace*{2mm}16\hspace*{2mm}\\
    & \hspace*{2mm} & $(x_0(x_0^2+x_1x_2):x_1^3:x_1(x_0^2+x_1x_2))$ & \hspace*{2mm}&\hspace*{2mm}$\mathsf{\{8\}}$\hspace*{2mm}
    &\hspace*{2mm}16\hspace*{2mm}\\
    & \hspace*{2mm} & $(x_1^2x_2:x_0(x_0x_2+x_1^2):x_1(x_0x_2+x_1^2))$ & \hspace*{2mm}&\hspace*{2mm}$\mathsf{\{9\}}$\hspace*{2mm}
    &\hspace*{2mm}15\hspace*{2mm}\\
    & \hspace*{2mm} & $(x_0(x_1^2+x_0x_2):x_1(x_1^2+x_0x_2):x_0x_1x_2)$ & \hspace*{2mm}&\hspace*{2mm}$\mathsf{\{10\}}$\hspace*{2mm}
    &\hspace*{2mm}16\hspace*{2mm}\\
    & \hspace*{2mm} & $(x_0(x_1^2+x_0x_2):x_1(x_1^2+x_0x_2):x_0x_1^2)$ & \hspace*{2mm}&\hspace*{2mm}$\mathsf{\{10\}}$\hspace*{2mm}
    &\hspace*{2mm}16\hspace*{2mm}\\
    & \hspace*{2mm} & $(x_0(x_0^2+x_1x_2):x_1(x_0^2+x_1x_2):x_0x_1^2)$ & \hspace*{2mm}&\hspace*{2mm}$\mathsf{\{10\}}$\hspace*{2mm}
    &\hspace*{2mm}16\hspace*{2mm}\\
    & \hspace*{2mm} & $(x_0(x_0x_1+x_0x_2+x_1x_2):x_1(x_0x_1+x_0x_2+x_1x_2):x_0x_1x_2)$ & \hspace*{2mm}&\hspace*{2mm}$\mathsf{\{11\}}$\hspace*{2mm}
    &\hspace*{2mm}16\hspace*{2mm}\\
    & \hspace*{2mm} & $(x_0(x_0^2+x_1x_2+x_0x_2):x_1(x_0^2+x_1x_2+x_0x_2):x_0x_1x_2)$ & \hspace*{2mm}&\hspace*{2mm}$\mathsf{\{11\}}$\hspace*{2mm}
    &\hspace*{2mm}16\hspace*{2mm}\\
    & \hspace*{2mm} & $(x_0(x_0^2+x_0x_1+x_1x_2):x_1(x_0^2+x_0x_1+x_1x_2):x_0x_1x_2)$ & \hspace*{2mm}&\hspace*{2mm}$\mathsf{\{12\}}$\hspace*{2mm}
    &\hspace*{2mm}16\hspace*{2mm}\\
    & \hspace*{2mm} & $(x_0(x_0^2+x_1x_2):x_1(x_0^2+x_1x_2):x_0x_1(x_0-x_1))$ & \hspace*{2mm}&\hspace*{2mm}$\mathsf{\{12\}}$\hspace*{2mm}
    &\hspace*{2mm}16\hspace*{2mm}\\
    & \hspace*{2mm} & $(x_0(x_1^2+\gamma x_0x_1+x_1x_2+x_0x_2):x_1(x_1^2+\gamma x_0x_1+x_1x_2+x_0x_2):x_0x_1x_2),$ $\gamma\not=0,\hspace{1mm} 1$ & \hspace*{2mm}&\hspace*{2mm}$\mathsf{\{13\}}$\hspace*{2mm}
    &\hspace*{2mm}16\hspace*{2mm}\\
    & \hspace*{2mm} & $(x_0(x_0^2+x_1^2+\gamma x_0x_1+x_0x_2):x_1(x_0^2+x_1^2+\gamma x_0x_1+x_0x_2):x_0x_1x_2),$ $\gamma^2\not=4,$ &\hspace*{2mm} &\hspace*{2mm}$\mathsf{\{14\}}$\hspace*{2mm}
    &\hspace*{2mm}16\hspace*{2mm}\\
    & \hspace*{2mm} & $(x_0(x_0^2+x_1x_2+x_0x_2):x_1(x_0^2+x_1x_2+x_0x_2):x_0x_1(x_0-x_1))$ & \hspace*{2mm}&\hspace*{2mm}$\mathsf{\{14\}}$\hspace*{2mm}
    &\hspace*{2mm}16\hspace*{2mm}\\
    & \hspace*{2mm} & $(x_0(x_0^2+x_1^2+\gamma x_0x_1+\delta x_0x_2+x_1x_2):x_1(x_0^2+x_1^2+\gamma x_0x_1+\delta x_0x_2+x_1x_2):x_0x_1x_2),$ $\gamma^2\not=4,$
    $\delta\not=\gamma_\pm$ & \hspace*{2mm}&\hspace*{2mm}$\mathsf{\{15\}}$\hspace*{2mm}&\hspace*{2mm}16\hspace*{2mm}\\
    & & & & & \\
\hline
\end{tabular}
\end{center}
\end{landscape}
\end{scriptsize}

Parmi ces mod\`eles cinq familles d\'ependent effectivement de
plus de un param\`etre. Ces param\`etres, attach\'es \`a des
configurations de courbes contract\'ees, peuvent \^etre reli\'es
\`a des configurations de points sur la droite
projective, pr\'ecis\'ement $4$ ou $5$ points suivant qu'il y a un
ou deux param\`etres. Le lecteur int\'eress\'e pourra d\'ecrire
l'espace des invariants; il y rencontrera l'invariant $j$ des
courbes elliptiques.

\begin{rem}
Dans le tableau les diff\'erents mod\`eles ne sont pas g.d.
\'equivalents. On le v\'erifie en examinant les invariants: type
et nombre de composantes des z\'eros du d\'eterminant jacobien et
leur multiplicit\'e, configuration des points d'ind\'etermination
(par exemple alignement ou non) etc.
\end{rem}

\begin{rem}
Toute transformation birationnelle cubique est donc \`a
conjugaison dynamique pr\`es du type $Af$ o\`u $A$ d\'esigne un
automorphisme de $\mathbb{P}^2(\mathbb{C})$ et $f$ un \'el\'ement de la liste
donn\'ee ci-dessus.
\end{rem}

Comme nous l'avons fait dans le cas quadratique on peut \'etudier 
l'espace des relations lin\'eaires associ\'e \`a une application 
birationnelle $f$ de $\mathring{\mathrm{B}}\mathrm{ir}_3.$ Si 
$f=(f_0,f_1,f_2)$ on introduit l'espace vectoriel 
\begin{align*}
\mathrm{RL}(f):=\{L=(L_0,L_1,L_2) \text{ application lin\'eaire }
\hspace{1mm}|\hspace{1mm} L_0f_0+L_1f_1+L_2f_2=0\}.
\end{align*}
En utilisant la classification pr\'ec\'edente on v\'erifie sans peine
que pour toute transformation $f$ purement cubique $\dim
\mathrm{RL}(f)=1.$ De plus si $L$ est un \'el\'ement non trivial
de $\mathrm{RL}(f)$ alors $L$ est de rang deux. On obtient alors
un \'enonc\'e semblable au Th\'eor\`eme \ref{bb}.

\begin{thm}\label{bb}
{\sl Soit $f$ un \'el\'ement de $\mathring{\mathrm{B}}\mathrm{ir}_3.$
Il existe une application lin\'eaire $L=(L_0,L_1,L_2)$ et une 
application quadratique $Q=(Q_0,Q_1,Q_2)$ telles que 
\begin{align*}
f=L\wedge Q=(L_1Q_2-L_2Q_1:L_2Q_0-L_0Q_2:L_0Q_1-L_1Q_0).
\end{align*}
En particulier les \'el\'ements de $\mathring{\mathrm{B}}\mathrm{ir}_3$
sont des applications d\'eterminantielles. }
\end{thm}

\begin{eg}
Pour l'\'el\'ement g\'en\'erique pr\'esentant la configuration $\mathsf{\{15\}}$
on constate que 
\begin{align*}
& L=(x_1,-x_0,0) && \text{et} && Q=(0,x_0x_2,-(x_0^2+x_1^2+\gamma x_0x_1
+\delta x_0x_2+x_1x_2))
\end{align*}
conviennent.
\end{eg}

Toutefois l'\'enonc\'e \ref{bbb} est moins complet que le Th\'eor\`eme
\ref{bb}; en effet si $L$ est une application lin\'eaire de rang deux 
et $Q$ une application quadratique, alors $L\wedge Q$ ne d\'efinit
pas n\'ecessairement une application birationnelle.

\section{Irr\'eductibilit\'e de $\mathring{\mathrm{B}}\mathrm{ir}_3$}\hspace{1mm}

Pour $\gamma^2\not=4$ et $\delta\not=\gamma_\pm$ on note $\xi_{\gamma,
\delta}$ la transformation
\begin{align*}
(x_0(x_0^2+x_1^2+\gamma x_0x_1+\delta x_0x_2+x_1x_2):x_1(x_0^2+x_1^2+\gamma
x_0x_1+\delta x_0x_2+x_1x_2):x_0x_1x_2)
\end{align*}
qui pr\'esente la configuration $\mathsf{\{15\}}.$ Consid\'erons
$A$ et $B$ deux automorphismes de $\mathbb{P}^2(\mathbb{C})$ tels que
\begin{equation}\label{lonelx_1}
A\xi_{\gamma,\delta}B=\xi_{\gamma',\delta'}.
\end{equation}
On constate que $\mathrm{Exc}\hspace{1mm}( A\xi_{\gamma,\delta}B)=B^{-1}(\mathrm{Exc}\hspace{1mm}
\xi_{\gamma,\delta})=\mathrm{Exc}\hspace{1mm}\xi_{\gamma',\delta'}.$ Comme
\begin{align*}
\mathrm{Exc}\hspace{1mm}\xi_{\gamma,\delta}=\{x_0=0,\hspace{1mm} x_1=0,\hspace{1mm}x_0^2+x_1^2+\gamma x_0x_1+\delta x_0x_2+x_1x_2=0,\hspace{1mm}x_0^2+x_1^2
+\gamma x_0x_1=0\}
\end{align*}
l'ensemble des automorphismes $A$ et $B$ de $\mathbb{P}^2(\mathbb{C})$ satisfaisant (\ref{lonelx_1}) est fini.

Posons
\begin{align*}
\mathscr{X}:=\{A\xi_{\gamma,\delta}B\hspace{1mm}|\hspace{1mm} A,\hspace{1mm} B\in\mathrm{PGL}_3(\mathbb{C}),\hspace{1mm} \gamma,\hspace{1mm} \delta\in\mathbb{C},
\hspace{1mm}\gamma^2\not=4,\hspace{1mm}\delta\not=\gamma_\pm\}\subset\mathring{\mathrm{B}}\mathrm{ir}_3.
\end{align*}
Il s'agit d'une r\'eunion d'orbites gauche-droite qui
repr\'esentent les transformations dont l'ensemble exceptionnel
est constitu\'e d'une conique lisse $\mathcal{C}$ et quatre
droites non tangentes \`a $\mathcal{C}$ (configuration
$\mathsf{\{15\}}$).

La Proposition qui suit est cons\'equence directe de ce
qui pr\'ec\`ede.

\begin{pro}
{\sl L'adh\'erence de $\mathscr{X}$ dans $\mathrm{Rat}_3\simeq\mathbb{P}^{29}(\mathbb{C})$ est une
vari\'et\'e alg\'ebrique irr\'eductible de dimension $18.$}
\end{pro}

Dans les trait\'es anciens on trouve le discours suivant: les courbes du r\'eseau
homalo\"idal associ\'e \`a une transformation birationnelle cubique satisfont
les \'equations suivantes (\cite{AC})
\begin{align*}
&\sum\mu_i^2=8, &&\sum\mu_i=6, && \mu_i\in\mathbb{N}
\end{align*}
o\`u $\mu_i$ d\'esigne la multiplicit\'e aux points base (qui
rappelons le ne sont pas n\'ecessairement propres) et $q$ le nombre
de points base. Ce syst\`eme d'\'equations a pour solution, \`a
r\'eindexation pr\`es, $(\mu_1,\mu_2,\mu_3,
\mu_4,\mu_5)=(2,1,1,1,1).$ Ces formules sont satisfaites par tout
\'el\'ement \og g\'en\'erique\fg\hspace{1mm} de $\mathscr{X}$ qui compte
cinq points d'ind\'etermination dont $1$ d'ordre $2$ et $4$
d'ordre $1.$

\begin{pro}[D\'ecomposition de \textsc{N\oe ther} des \'el\'ements
de $\mathscr{X}$]\label{decompnoe}

{\sl Soit $f$ un \'el\'ement de $\mathscr{X};$ alors $f$ poss\`ede
une \'ecriture du type
\begin{align*}
&B\sigma A\sigma C, && A,\hspace{1mm} B,\hspace{1mm} C\in\mathrm{PGL}_3(\mathbb{C}).
\end{align*}}
\end{pro}

\begin{proof}[\sl D\'emonstration]
Consid\'erons un automorphisme $A$ de $\mathbb{P}^2(\mathbb{C})$ de la forme
sp\'eciale suivante
\begin{align*}
A=(b_0x_1+c_0x_2:b_1x_1+c_1x_2:a_2x_0+b_2x_1+c_2x_2).
\end{align*}
La transformation $\sigma A\sigma$ s'\'ecrit
\begin{align*}
&((b_1x_1+c_1x_2)Q:(b_0x_1+c_0x_2)Q:x_0(b_0x_1+c_0x_2)(b_1x_1+c_1x_2))
\end{align*}
avec $Q=a_2x_1x_2+b_2x_0x_2+c_2x_0x_1.$ On constate qu'elle est
cubique pure d\`es que $a_2b_2c_2$ est non nul, condition qui
caract\'erise la lissit\'e de la conique $Q=0.$ Visiblement les
courbes 
\begin{align*}
&Q=0, && b_0x_1+c_0x_2=0,&&b_1x_1+c_1x_2=0,&&x_1=0&& 
\text{et}&&x_2=0
\end{align*}
sont contract\'ees. Ainsi pour des valeurs g\'en\'eriques
des param\`etres la configuration des courbes contract\'ees par
$\sigma A\sigma$ est du type $\mathsf{\{15\}}.$ Notons que la
configuration de $\mathrm{Exc}\hspace{1mm}(\sigma A\sigma)$ d\'etermine les
lignes de la matrice $A$ \`a multiplication scalaire pr\`es, ce
qui correspond \`a l'action \`a gauche du groupe diagonal, et
permutation $(x_0:x_2:~x_1)$ pr\`es. On en d\'eduit facilement que
l'ensemble des transformations du type $B \sigma A\sigma C,$ avec
$A$ comme ci-dessus, est de dimension $18.$ Par suite aux
configurations de type $\mathsf{\{15\}}$ est associ\'ee une
d\'ecomposition de la forme
\begin{align*}
& B\sigma A\sigma C, && A=(b_0x_1+c_0x_2:b_1x_1+c_1x_2: a_2x_0+b_2x_1+
c_2x_2), && B,\hspace{1mm} C\in\mathrm{PGL}_3(\mathbb{C}).
\end{align*}
\end{proof}

\begin{rem}
L'\'enonc\'e pr\'ec\'edent donne une nouvelle justification \`a
l'\'etude des transformations $A$ telles que le degr\'e de $\sigma
A\sigma$ soit anormal c'est-\`a-dire plus petit ou \'egal \`a $3,$
\'etude faite au Chapitre \ref{dyndyn}.
\end{rem}

\bigskip

Notons
$\complement\mathscr{X}=\mathring{\mathrm{B}}\mathrm{ir}_3\setminus\mathscr{X};$
chaque \'el\'ement de $\complement\mathscr{X}$ est g.d. conjugu\'e \`a
l'un des mod\`eles repr\'esentant toutes les configurations
except\'ees $\mathsf{ \{15\}}.$ Comme ces mod\`eles d\'ependent au
plus d'un param\`etre, on en d\'eduit que la dimension de
l'adh\'erence de $\complement\mathscr{X}$ est en tout point inf\'erieure
ou \'egale \`a $17.$

D'apr\`es la Remarque \ref{deginv} et pour des raisons de dimension on a
$\textrm{Inv}( \overline{\mathscr{X}})=\overline{\mathscr{X}},$ ce que l'on peut
voir de fa\c{c}on explicite en inversant les $\xi_{\gamma,\delta}.$\bigskip

Pour d\'emontrer l'irr\'eductibilit\'e de
$\mathring{\mathrm{B}}\mathrm{ir}_3$ on proc\`ede comme suit; on
cherche \`a joindre tout \'el\'ement $f_0$ de $\complement\mathscr{X}$ \`a
l'adh\'erence de $\mathscr{X}$ par un chemin rationnel $s\mapsto f_s$ de
sorte que $f_{s=0}=f_0$ et $f_s$ appartiennent \`a
$\overline{\mathscr{X}}$ pour les valeurs g\'en\'eriques de $s.$ Il est
souvent plus commode de travailler avec des formes \og
pr\'enormales\fg\hspace{1mm} qui font intervenir \og beaucoup\fg\hspace{1mm} de
param\`etres, comme celles donn\'ees dans le Lemme \ref{vent},
plut\^ot qu'avec les formes finales (par exemple celles de la
Proposition \ref{listee}). On proc\`ede au cas par cas suivant la
nature du lieu exceptionnel.

\subsection{Le lieu exceptionnel contient une
conique}\hspace{1mm}\label{exccon}

Nous allons essayer de joindre les formes pr\'enormales
du Lemme \ref{vent} de type $\mathfrak{(b)}$ \`a celles de type
$\mathfrak{(a)}$ qui sont n\'ecessairement dans $\overline{\mathscr{X}}.$ Soit
\begin{align*}
&f_0=(x_0Q: x_1Q:x_0x_1^2),&& Q=\alpha x_0^2+\beta x_1^2+
\gamma x_0x_1+\delta x_0x_2+\varepsilon x_1x_2
\end{align*}
un \'el\'ement de type $\mathfrak{(b)}.$ On consid\`ere la
transformation
\begin{align*}
f_s= (x_0Q:x_1Q:x_0x_1(sx_2+x_1)).
\end{align*}
Pour $s=0,$ on atteint la transformation
initiale $f_0$ et pour $s$ non nul $f_s$ est
dans $\mathscr{X}.$ \bigskip

Maintenant on cherche \`a relier un \'el\'ement $f_0$ de
type $\mathfrak{(c)}$
\begin{align*}
f_0=(x_0(x_0^2+\gamma x_0x_2+x_1x_2):x_1(x_0^2+\gamma
x_0x_2+x_1x_2):x_0x_1(x_0 -x_1))
\end{align*}
\`a $\overline{\mathscr{X}}.$ On introduit
\begin{align*}
f_s:=(x_0(x_0^2+\gamma x_0x_2+x_1x_2):x_1(x_0^2+\gamma x_0x_2+
x_1x_2):x_0x_1(sx_2+x_0-x_1)).
\end{align*}
Pour $s=0,$ on obtient la transformation initiale
et pour $s\not=0$ on constate que $f_s$
est dans $\overline{\mathscr{X}};$ en effet $f_s\left(x_0:x_1:\frac{x_2-x_0+x_1}{s}\right)$
s'\'ecrit
\begin{align*}
& \left(x_0P:x_1P:x_0x_1x_2\right),
&& P=\left(1-\frac{\gamma}{s}\right)x_0^2+\frac{x_1^2}{s}+\frac{\gamma
-1}{s}x_0x_1+\frac{\gamma}{s}x_0x_2+\frac{x_1x_2}{s}
\end{align*}

\bigskip

Pour les \'el\'ements $f_0$ de type $\mathfrak{(d)}$
\begin{align*}
f_0=(x_0x_2(x_1+\gamma x_0):x_1x_2(x_1+\gamma x_0):x_0x_1(x_0-x_1))
\end{align*}
on peut proc\'eder
simplement comme suit. Les transformations
\begin{align*}
f_s=(x_0(sx_0^2+\gamma x_0x_2+x_1x_2):x_1(s
x_0^2+\gamma x_0x_2+x_1x_2):x_0x_1(x_0-x_1))
\end{align*}
sont conjugu\'ees pour $s\not=0$ \`a $f_{s=1}$ que
l'on sait \^etre dans $\overline{\mathscr{X}}$ car de la forme
$\mathfrak{(c)}.$ Pour $s=0,$ on a l'\'el\'ement $f_0$
de type $\mathfrak{(d)}.$\bigskip

\subsection{Le lieu exceptionnel est constitu\'e de plus de
$4$ droites}\hspace{1mm}

Il s'agit de montrer que la transformation
\begin{align*}
& f=(x_1(x_0-x_1)(x_0+x_2):x_0(x_0-x_1)(x_2-x_1):(x_0+x_1)x_1x_2)
\end{align*}
obtenue dans la Proposition \ref{aumoins4} est dans $\overline{\mathscr{X}}.$
On constate que
\begin{align*}
(x_0+x_1:x_1:x_2)f(x_0-x_1:x_1:x_2),
\end{align*}
qui s'\'ecrit aussi
\begin{align*}
& (x_0(x_0^2+2x_1^2-3x_0x_1+x_0
x_2-2x_1x_2):x_1(x_0^2+2x_1^2-3x_0x_1+x_0x_2-2x_1x_2):x_0x_1x_2),
\end{align*}
est du type $\mathfrak{(a)}$ \`a homoth\'etie pr\`es donc $f$ appartient bien \`a $\overline{\mathscr{X}}.$

\bigskip

\subsection{Le lieu exceptionnel est constitu\'e de
trois droites}\hspace{1mm}

Les \'el\'ements de $\mathring{\mathrm{B}}\mathrm{ir}_3$ dont le lieu exceptionnel
est constitu\'e d'une conique lisse et deux droites appartiennent
\`a $\overline{\mathscr{X}};$ on d\'eduit de l'invariance de
$\overline{\mathscr{X}}$ par $\textrm{Inv}$ que les mod\`eles obtenus dans
la Proposition \ref{3dtes} sont dans l'adh\'erence de $\mathscr{X}.$ \bigskip

Consid\'erons une transformation de type
$\mathfrak{(a)}.$ Si $\alpha=\varepsilon=1$ et $\beta=\gamma=
\delta=0$ elle s'\'ecrit
\begin{align*}
f=(x_0(x_0^2+x_1x_2): x_1(x_0^2+x_1x_2):x_0x_1x_2).
\end{align*}
Posons
\begin{align*}
g_s:=\left(
\frac{x_0-x_2}{s}:x_1:x_2\right)f(sx_0:sx_1:x_2)
=(x_0^3:x_1(sx_0^2+x_1x_2):x_0x_1x_2);
\end{align*}
on remarque que $g_{s=0}=(x_0^3:x_1^2x_2 :x_0x_1x_2)$ est le
mod\`ele de la Proposition \ref{3dtesb}.\bigskip

Reste \`a consid\'erer les deux mod\`eles dont le lieu
exceptionnel est constitu\'e de trois droites concourantes. Lorsque
$\alpha=\beta=1,$ $\gamma=2$ et $\delta=\varepsilon=0$ une
transformation du type $\mathfrak{(a)}$ est de la forme
\begin{align*}
(x_0(x_0+x_1)^2):x_1(x_0+x_1)^2):x_0x_1x_2);
\end{align*}
on constate alors que
\begin{align*}
(x_0+x_1:x_1:x_2)(x_0(x_0+x_1)^2):x_1(x_0+x_1)^2):x_0x_1x_2)
(x_0-x_1:x_1:x_2)
\end{align*}
co\"incide avec $(x_0^3:x_0^2x_1:x_1x_2(x_0-x_1)),$ l'un des deux
mod\`eles de la Proposition \ref{3concourantes}.

Consid\'erons la seconde transformation obtenue dans la
Proposition \ref{3concourantes}
\begin{align*}
(x_0^2(x_0-x_1):x_0x_1(x_0-x_1):x_0x_1x_2+x_1^3);
\end{align*}
elle s'\'ecrit dans la carte $x_1=1$ comme suit
\begin{align*}
\left(x_0,\frac{x_0x_2+1}{x_0(x_0-1)}\right).
\end{align*}
On peut la pertuber de la fa\c{c}on suivante
\begin{align*}
&\left(x_0,\frac{x_0x_2+1}{x_0(x_0-1)+\varepsilon x_2}\right),
&&\varepsilon\in\mathbb{C};
\end{align*}
on remarque que c'est une homographie en $x_2.$ En
homog\'en\'eisant on obtient donc l'\'el\'ement de~$\mathring{\mathrm{B}}\mathrm{ir}_3$ suivant
\begin{align*}
(x_0(x_0(x_0-x_1)+\varepsilon x_1x_2):x_1(x_0(x_0-x_1)+\varepsilon
x_1x_2):x_0x_1x_2+x_1^3),
\end{align*}
\'el\'ement qui contracte la conique lisse
$x_0(x_0-x_1)+\varepsilon x_1x_2$ d\`es que $\varepsilon$ est non
nul; or un tel \'el\'ement est dans $\overline{\mathscr{X}}$ (\emph{voir}
\S\ref{exccon}) donc
\begin{align*}
(x_0^2(x_0-x_1):x_0x_1(x_0-x_1):x_0x_1x_2+x_1^3);
\end{align*}
l'est aussi.

\subsection{Le lieu exceptionnel est constitu\'e de
deux droites}\hspace{1mm}

Le mod\`ele $(x_0x_1x_2:x_1x_2^2:x_2^3-x_0^2x_1)$ de la Proposition
\ref{2dtes}, dont l'inverse
\begin{align*}
& (x_0(x_1x_2+x_0^2):x_1^3:x_1(x_1x_2+x_0^2))
\end{align*}
contracte une droite et une conique,
est, de part l'invariance de $\overline{\mathscr{X}}$ par $\textrm{Inv},$
dans $\overline{\mathscr{X}}.$

\bigskip

Voici rapidement comment on proc\`ede pour les cas
restants. Consid\'erons une transformation de la forme
$\mathfrak{(b)}$ que l'on sait \^etre dans $\overline{\mathscr{X}}.$
Apr\`es conjugaison dynamique par $(x_1:x_0:x_2)$ on obtient lorsque $\beta
=\delta=1$ et
$\alpha=\gamma=\varepsilon=0$ l'\'el\'ement
\begin{align*}
f=(x_0^3+x_0x_1x_2:x_0^2x_1+x_1^2x_2:x_0^2x_1).
\end{align*}
Posons
\begin{align*}
g_s=\left(x_2:\frac{x_1-x_2}{s}:x_0\right)f(x_2:x_0:
sx_1)=(x_0x_2^2:x_0^2x_1:x_2^3+s x_0x_1x_2);
\end{align*}
alors $g_{s=0}=( x_0x_2^2:x_0^2x_1:x_2^3)$ est le mod\`ele obtenu dans la
Proposition \ref{2dtes}.

\bigskip

Partons encore d'une transformation du type $\mathfrak{(b)}.$
Lorsque $\alpha=\beta=\varepsilon=1$ et $\gamma=\delta=0$ on a
\begin{align*}
(x_0(x_0^2+x_1^2+x_1x_2):x_1(x_0^2+x_1^2+x_1x_2):x_0x_1^2).
\end{align*}
Posons
\begin{align*}
&g_s:=\left(\frac{x_2}{s}:\frac{x_0-x_2}
{s^3}:x_1\right)(x_0(x_0^2+x_1^2+x_1x_2):x_1(x_0^2+x_1^2+x_1x_2):x_0x_1^2)(
sx_0:x_2:s^2x_1),
\end{align*}
soit
\begin{align*}
&g_s=(x_0x_2^2:x_0^3+x_0x_1x_2:x_2(s^2 x_0^2+x_2^2+s^2x_1 x_2)).
\end{align*}
On remarque $g_{s=0}=(x_0x_2^2:x_0^3+x_0x_1x_2:x_2^3)$ est un des cas
de la Proposition \ref{2dtes}; il est donc dans~$\overline{\mathscr{X}}.$\bigskip

Le mod\`ele $(x_0^2x_2:x_0^2x_1+x_2^3:x_0x_2^2)$ de la Proposition
\ref{2dtes} a pour inverse
\begin{align*}
(x_0^3:x_0x_1x_2-x_2^3:x_0^2x_2)
\end{align*}
et
\begin{align*}
& (-x_2:x_1:x_0)(x_0^3:x_0x_1x_2-x_2^3:x_0^2x_2)(-x_2:x_1:x_0)=
(x_0x_2^2:x_0^3+x_0x_1x_2:x_2^3).
\end{align*}
On vient de voir que $(x_0x_2^2:x_0^3+x_0x_1x_2:x_2^3)$ appartient
\`a $\overline{\mathscr{X}};$ comme
$\overline{\mathscr{X}}$ est invariant par $\textrm{Inv}$ la transformation
$(x_0^2x_2:x_0^2x_1+x_2^3:x_0x_2^2)$ est dans l'adh\'erence de $\mathscr{X}.$
\bigskip

Pour finir prenons la transformation $f=(x_0^2x_2: x_0^3+x_2^3
+x_0x_1x_2:x_0x_2^2)$ de la Proposition~\ref{2dtes}. Soit $f_s$ la
transformation birationnelle d\'efinie par
\begin{align*}
f_s=(x_0^2x_2:x_0^3+x_2^3+x_0x_1(x_2-sx_0):x_0x_2(x_2-
sx_0)).
\end{align*}
On constate que pour $s$ non nul $f_s$ contracte trois
droites donc est dans $\overline{\mathscr{X}}$ et que $f_{s=0}=f;$ ceci
implique que $f$ appartient aussi \`a $\overline{\mathscr{X}}.$\bigskip

\subsection{Le lieu exceptionnel est constitu\'e
d'une droite}\hspace{1mm}

Consid\'erons le mod\`ele $(x_0x_2^2:x_0^3+ x_0x_1x_2:x_2^3)$ qui
contracte exactement deux droites. Posons
\begin{align*}
g_s=(x_1-s^3x_2:x_0-sx_2:x_2)(x_0x_2^2
:x_0^3+x_0x_1x_2:x_2^3)\left(x_1+sx_2:\frac{x_0}{s}-3sx_1:
x_2\right);
\end{align*}
on constate que $g_s=\left(
x_1^3+\frac{x_0x_1x_2}{s}+x_0x_2^2:x_1x_2^2:x_2^3\right).$
Lorsque $s$ tend vers $+\infty,$ la transformation~$g_s$ tend vers $(x_1^3+x_0x_2^2:x_1x_2^2:x_2^3);$ on reconnait
l'\'el\'ement obtenu dans la Proposition \ref{1dte} qui contracte
pr\'ecis\'ement une droite.\bigskip

Finalement on obtient le:

\begin{thm}
{\sl L'adh\'erence (ordinaire) de $\mathscr{X}$ dans $\mathring{\mathrm{B}}\mathrm{ir}_3$ est
$\mathring{\mathrm{B}}\mathrm{ir}_3.$ En particulier $\mathring{\mathrm{B}}\mathrm{ir}_3$ est irr\'eductible.}
\end{thm}

\subsection{Non irr\'eductibilit\'e de $\mathrm{Bir}_3$ en tant que
sous-ensemble de $\mathbb{P}^{29}(\mathbb{C})$}\hspace{1mm}

Alors que $\mathrm{Bir}_2$ est lisse et irr\'eductible, il n'en
est pas de m\^eme pour $\mathrm{Bir}_3$ (vu comme sous-ensemble de
$\mathbb{P}^{29}(\mathbb{C}));$ comme nous l'avons vu $\mathring{\mathrm{B}}\mathrm{ir}_3$ est
irr\'eductible (en fait rationnellement connexe). Consid\'erons un
\'el\'ement de $\mathrm{Bir}_3$ du type suivant $f_0=\ell A\sigma$ o\`u
$\ell$ est une forme lin\'eaire g\'en\'erique, disons $x_0+x_1+x_2.$ Il
est clair que $f_0$ en tant que triplet de formes homog\`enes ne
peut \^etre limite de transformations de $\mathring{\mathrm{B}}\mathrm{ir}_3.$ Ceci r\'esulte
du fait qu'aucun \'el\'ement de $\mathring{\mathrm{B}}\mathrm{ir}_3$ ne contracte quatre
droites en position g\'en\'erale. Ceci implique que $\mathrm{Bir}_3$ n'est
pas irr\'eductible; plus pr\'ecis\'ement les composantes
irr\'eductibles de $\overline{\mathrm{Bir}_3}$ sont $\overline{\mathring{\mathrm{B}}\mathrm{ir}_3}$ et
$\overline{\mathrm{Bir}_3'}$ o\`u $\mathrm{Bir}_3'$ d\'esigne l'ensemble des
transformations de type $\ell Q$ avec $\ell$ forme lin\'eaire et
$Q$ dans $\mathrm{Bir}_2.$ On note que $\overline{\mathrm{Bir}_3'}$ est de
dimension $16.$ Toutefois chaque \'el\'ement de $\mathrm{Bir}_2$ se
repr\'esente dans $\overline{\mathring{\mathrm{B}}\mathrm{ir}_3}.$ Plus pr\'ecis\'ement
avec les notations habituelles $\mathrm{Bir}_3^\bullet=(\overline{
\mathring{\mathrm{B}}\mathrm{ir}_3})^\bullet.$

\section{D\'ecomposition de \textsc{N\oe ther}}\hspace{1mm}

Nous avons vu dans la Proposition \ref{decompnoe} la
d\'ecomposition de \textsc{N\oe ther} des transformations cubiques
g\'en\'eriques, {\it i.e.} les \'el\'ements de $\mathscr{X}.$ Dans les
Lemmes qui suivent on propose une d\'ecomposition de \textsc{N\oe ther}
pour les \'el\'ements du tableau obtenu au Chapitre \ref{troistrois},
\S\ref{dimor}. On prendra garde au fait que ces d\'ecompositions ne
sont peut-\^etre pas optimales.

\begin{lem}
{\sl La d\'ecomposition de \textsc{N\oe ther} d'un \'el\'ement de
$\mathring{\mathrm{B}}\mathrm{ir}_3$ ayant la configuration~$\mathsf{\{1\}}$ 
fait intervenir au plus huit fois l'involution de
\textsc{Cremona}.}
\end{lem}

\begin{proof}[\sl D\'emonstration]
On a vu aux Chapitre \ref{troistrois}, \S\ref{conn} et Chapitre \ref{troistrois},
\S\ref{sconn} qu'une transformation birationnelle purement cubique poss\'edant la configuration
$\mathsf{\{1\}}$ s'\'ecrit $(x_0x_2^2+x_1^3:x_1x_2^2:~x_2^3)$ \`a conjugaison
g.d. pr\`es. Or
\begin{align*}
(x_0x_2^2+x_1^3:x_1x_2^2:x_2^3)=(x_0:x_2:x_1)\rho(x_2:x_1:x_0)\tau
(x_2:x_1:-x_0)\rho(x_0:x_2:x_1)
\end{align*}
d'o\`u le r\'esultat.
\end{proof}

\begin{lem}
{\sl Soit $f$ une transformation birationnelle purement cubique dont la
configuration est $\{\mathsf{2}\};$ la d\'ecomposition de
\textsc{N\oe ther} de $f$ fait intervenir au plus huit fois l'involution de
\textsc{Cremona}.}
\end{lem}

\begin{proof}[\sl D\'emonstration]
\`A conjugaison g.d. pr\`es $f$ est de l'un des types suivants
\begin{align*}
& (x_0x_2^2:x_0^2x_1:x_2^3), && (x_0x_2^2:x_0^3+x_0x_1x_2:x_2^3),
&& (x_0^2x_2:x_0^3+x_2^3+x_0x_1x_2:x_0
x_2^2),\\
& (x_0^2x_2:x_0^2x_1+x_2^3:x_0x_2^2), && (x_0x_1x_2:x_1x_2^2:x_2^3-x_0^2x_1). &&
\end{align*}
L'\'enonc\'e r\'esulte des \'egalit\'es suivantes
\begin{align*}
& (x_0x_2^2:x_0^2x_1:x_2^3)=((x_2:x_1:x_0)\rho)^2,\\
& (x_0x_2^2:x_0^3+x_0x_1x_2:x_2^3)=(x_2:x_0:x_1)\rho(x_2:x_1:x_0)\tau(x_2:x_0:-x_1),\\
& (x_0^2x_2:x_0^3+x_2^3+x_0x_1x_2:x_0x_2^2)=(x_1:x_2:x_0)\tau(x_1:x_0:-x_2)\tau(x_0:x_2:-x_1),\\
& (x_0^2x_2:x_0^2x_1+x_2^3:x_0x_2^2)=(x_1:x_0:x_2)\rho(x_0+x_1:x_0:x_2)\rho(x_2:x_1:x_0)\rho,\\
& (x_0x_1x_2:x_1x_2^2:x_2^3-x_0^2x_1)=(x_1:x_0:-x_2)\tau(x_2:x_0:x_1)\rho.
\end{align*}
\end{proof}

\begin{lem}
{\sl Les \'el\'ements de $\mathring{\mathrm{B}}\mathrm{ir}_3$
pr\'esentant la configuration $\{\mathsf{3}\}$ comptent dans leur
d\'ecomposition de \textsc{N\oe ther} au plus cinq fois l'involution de
\textsc{Cremona}.}
\end{lem}

\begin{proof}[\sl D\'emonstration]
On a vu que toute transformation birationnelle purement cubique
poss\'edant la configuration $\{\mathsf{3}\}$ est \`a conjugaison
g.d. pr\`es de l'une des formes suivantes
\begin{align*}
& (x_0^3:x_1^2x_2:x_0x_1x_2), && (x_0^2(x_1-x_2):x_0x_1(x_1-x_2):x_1^2x_2), \\
& (x_0^2x_2:x_0x_1x_2:x_0x_1^2-x_1^2x_2), && (x_0x_1x_2:x_1^2x_2:x_0(x_1^2-x_0x_2)).
\end{align*}

On constate que $(x_0^3:x_1^2x_2:x_0x_1x_2)=(\rho(x_1:x_2:x_0))^2$ donc
toute transformation appartenant \`a $\mathcal{O}_{g.d.}(x_0^3:x_1^2x_2
:x_0x_1x_2)$ s'\'ecrit $A\rho(x_1:x_2:x_0)\rho B$ o\`u $A,$ $B$ d\'esignent
deux automorphismes de~$\mathbb{P}^2(\mathbb{C}).$

Comme $(x_0^2(x_1-x_2):x_0x_1(x_1-x_2):x_1^2x_2)=(x_0+x_2:x_2:x_1)\rho
(-x_0+x_1:x_2-x_1:x_0)\sigma$ on a l'inclusion

\begin{small}
\begin{align*}
\mathcal{O}_{g.d.}(x_0^2(x_1-x_2):x_0x_1(x_1-x_2
):x_1^2x_2)\subset\{A\rho(-x_0+x_1:x_2-x_1:x_0)\sigma B\hspace{1mm}|\hspace{1mm} A,\hspace{1mm} B\in\mathrm{PGL}_3(\mathbb{C})\}.
\end{align*}
\end{small}

Puisque $(x_0^2x_2:x_0x_1x_2:x_0x_1^2-x_1^2x_2)=(x_1:x_2:x_0+x_2)\rho(-x_0-x_1+x_2:x_0:x_1
)\sigma,$ pour toute transformation $f$ de $\mathcal{O}_{g.d.}(
x_0^2x_2:x_0x_1x_2:x_0x_1^2-x_1^2x_2)$ il existe $A$ et $B$ deux automorphismes de~$\mathbb{P}^2(\mathbb{C})$ tels que $f=A\rho(-x_0-x_1+x_2:x_0:x_1)\sigma B.$

Pour finir on note que $(x_0x_1x_2:x_1^2x_2:x_0(x_1^2-x_0x_2))=(x_1:x_0:-x_2)
\tau\sigma$ donc
\begin{align*}
\mathcal{O}_{g.d.}(x_0x_1x_2:x_1^2x_2:x_0(x_1^2-x_0x_2))
\subset\{ A\tau\sigma B\hspace{1mm}|\hspace{1mm} A,\hspace{1mm} B\in\mathrm{PGL}_3(\mathbb{C})\}.
\end{align*}
\end{proof}

\begin{lem}
{\sl La d\'ecomposition de \textsc{N\oe ther} d'une transformation birationnelle
purement cubique de configuration $\mathsf{\{4\}}$ fait intervenir
au plus six fois l'involution $\sigma.$}
\end{lem}

\begin{proof}[\sl D\'emonstration]
On a vu au \S\ref{33dtes} qu'un \'el\'ement de
$\mathring{\mathrm{B}}\mathrm{ir}_3$ de configuration
$\mathsf{\{4\}}$ est de l'un des types suivants
\begin{align*}
& (x_0^3:x_0^2x_1:x_1x_2(x_0-x_1));&&
(x_0^2(x_0-x_1):x_0x_1(x_0-x_1): x_1(x_0x_2+x_1^2)).
\end{align*}
Or
\begin{align*}
&
(x_0^3:x_0^2x_1:x_1x_2(x_0-x_1))=(x_0+x_2:x_2:x_1)\sigma(x_0:x_1:x_2-x_0)
\rho(x_1:x_2:x_0)
\end{align*}
et $(x_0^2(x_0-x_1):x_0x_1(x_0-x_1):x_1(x_0x_2+x_1^2))$ s'\'ecrit
\begin{align*}
(x_1+x_2:x_2:x_0+
x_2)\rho(-x_0+x_1+x_2:x_1:x_0-x_1)\tau(x_0:x_1:-x_2)
\end{align*}
d'o\`u l'\'enonc\'e.
\end{proof}

\begin{lem}
{\sl La d\'ecomposition de \textsc{N\oe ther} des transformations de
\textsc{Cremona} cubiques dont la configuration est
$\{\mathsf{5}\}$ fait intervenir au plus trois fois l'involution
de \textsc{Cremona}.}
\end{lem}

\begin{proof}[\sl D\'emonstration]
Soit $f$ une telle transformation birationnelle; elle est, \`a
conjugaison g.d. pr\`es, de l'un des types suivants
\begin{align*}
& (x_0x_2(x_0+x_1):x_1x_2(x_0+x_1):x_0x_1^2), \\
& (x_0(x_0+x_1)(x_1+x_2):x_1(x_0+x_1)(x_1+x_2):x_0x_1x_2), \\
& (x_0(x_0+x_1+x_2)(x_0+x_1):x_1(x_0+x_1+x_2)(x_0+x_1):x_0x_1x_2).
\end{align*}
On constate que
\begin{footnotesize}
\begin{align*}
& (x_0x_2(x_0+x_1):x_1x_2(x_0+x_1):x_0x_1^2)=(x_1+x_2:-x_2:-x_0-x_2)\rho(-x_0-x_1+x_2:-x_0:x_0+x_1)\sigma,
\\
&  (x_0(x_0+x_1)(x_1+x_2):x_1(x_0+x_1)(x_1+x_2):x_0x_1x_2)=(-x_0:
\beta x_2:\alpha\beta x_1)\sigma(x_2:x_0+\alpha x_1:x_1+\beta x_2)\sigma A,\\
& (x_0(x_0+x_1+x_2)(x_0+x_1):x_1(x_0+x_1+x_2)(x_0+x_1):x_0x_1x_2)=\sigma(x_0+x_2:-x_0+x_2:4x_1+2x_2)\rho B
\end{align*}
\end{footnotesize}
avec
\begin{align*}
&A=\left(-\frac{x_2}{\alpha\beta}:-\frac{x_1}{\beta}:x_0+x_1\right),&&
B=\left(-\frac{x_0}{2}+\frac{x_1}{2}:x_2:\frac{x_0}{2}+\frac{x_1}{2}\right).
\end{align*}
On en d\'eduit le r\'esultat annonc\'e.
\end{proof}

\begin{lem}
{\sl Soit $f$ une transformation birationnelle quadratique avec la
configuration~$\{\mathsf{6}\};$ il existe $A$ dans $\mathrm{PGL}_3
(\mathbb{C})$ tel qu'\`a conjugaison g.d. pr\`es $f$ s'\'ecrive $\sigma
A\sigma.$}
\end{lem}

\begin{proof}[\sl D\'emonstration]
On a vu que $f$ co\"incide \`a conjugaison g.d. pr\`es avec
\begin{align*}
& (x_0(x_0^2+x_1^2+\gamma x_0x_1):x_1(x_0^2+x_1^2+\gamma x_0x_1):x_0x_1x_2),&&
\gamma^2\not=4.
\end{align*}
Or lorsque $\gamma^2\not=4$ on a
$(x_0(x_0^2+x_1^2+\gamma x_0x_1):x_1(x_0^2+x_1^2+\gamma x_0x_1):x_0x_1x_2)
=\sigma A\sigma B$ o\`u
\begin{align*}
&\hspace{1mm} B=\left(x_0+\gamma_+x_1:x_0+\gamma_-x_1:x_2\right),&&
\hspace{1mm}\gamma_{\pm}=\frac{
\gamma\pm\delta^{-1}}{2},\\
&A=\left(-\delta x_0+\delta x_1:\delta\gamma_+x_0-\delta\gamma_-x_1:x_2\right), &&
\delta=\frac{1}{\sqrt{\gamma^2-4}},
\end{align*}
d'o\`u l'\'enonc\'e.
\end{proof}

\begin{lem}
{\sl Soit $f$ un \'el\'ement de $\mathring{\mathrm{B}}\mathrm{ir}_3$
dont la configuration est $\mathsf{\{7\}};$ on compte dans la
d\'ecomposition de \textsc{N\oe ther} de $f$ au plus trois fois
l'involution de \textsc{Cremona}.}
\end{lem}

\begin{proof}[\sl D\'emonstration]
D'apr\`es les  Chapitre \ref{troistrois}, \S\ref{conn} et
Chapitre \ref{troistrois}, \S\ref{sconn} la transformation $f$ est
\`a conjugaison g.d. pr\`es de l'un des types suivants

\begin{small}
\begin{align*}
&(x_0(x_0^2+x_1^2+\gamma x_0x_1+\gamma_+x_0x_2+x_1x_2):x_1(
x_0^2+x_1^2+\gamma x_0x_1
+\gamma_+x_0x_2+x_1x_2):x_0x_1x_2),&& \gamma^2
\not=4, \\
& (x_0x_2(x_1+x_0):x_1x_2(x_1+x_0):x_0x_1(x_0-x_1)),&& \\
& (x_1(x_0-x_1)(x_0+x_2):x_0(x_0-x_1)(x_2-x_1):x_1x_2(x_0+x_1)). &&
\end{align*}
\end{small}

L'\'el\'ement
$(x_0(x_0^2+x_1^2+\gamma x_0x_1+\gamma_+x_0x_2+x_1x_2):x_1(
x_0^2+x_1^2+\gamma x_0x_1+\gamma_+x_0x_2+x_1x_2):x_0x_1x_2)$ s'\'ecrit
lorsque $\gamma^2\not=4$
\begin{align*}
\sigma(\delta\gamma_+x_0-\delta\gamma_-x_1:-\delta x_0+\delta x_1:x_1+x_2)
\sigma(\gamma_+x_0+x_1:\gamma_-x_0+x_1:x_2)
\end{align*}
avec $\delta=\frac{1}{\sqrt{\gamma^2-4}}$ et
$\gamma_{\pm}=\frac{\gamma\pm\delta^{-1}}{2}.$

On remarque que
$(x_0x_2(x_1+x_0):x_1x_2(x_1+x_0):x_0x_1(x_0-x_1))=A\rho B\sigma$ o\`u $A$ d\'esigne
l'automorphisme de $\mathbb{P}^2(\mathbb{C})$ d\'efini par
$(2(x_1+x_2):2(x_1-x_2):2(x_0+x_2))$ et $B$ celui donn\'e par
\begin{align*}
\left(-\frac{x_0}{2}-\frac{x_1}{2}+x_2:-\frac{x_0}{2}+\frac{x_1}{2}:\frac{x_0}{2}+\frac{x_1}{2}
\right).
\end{align*}

Pour finir notons que $(x_1(x_0-x_1)(x_0+x_2):x_0(x_0-x_1)(x_2-x_1):x_1x_2(x_0+x_1))$
s'\'ecrit aussi
\begin{align*}
(x_0+x_2:x_2-x_0:x_1-x_2)
\rho(x_0-x_1+2x_2:-x_0+x_1:x_0+x_1)\sigma.
\end{align*}
\end{proof}

\begin{lem}
{\sl Soit $f$ un \'el\'ement de $\mathring{\mathrm{B}}\mathrm{ir}_3$
ayant la configuration $\{\mathsf{8}\};$ alors \`a conjugaison
gauche-droite pr\`es $f= \rho(x_1:x_2:x_0)\tau.$ Ainsi $\sigma$
apparait au plus six fois dans la d\'ecomposition de \textsc{N\oe ther} de
$f.$}
\end{lem}

\begin{proof}[\sl D\'emonstration]
On sait que $f$ s'\'ecrit $(x_0(x_0^2+x_1x_2):x_1^3:x_1(x_0^2+x_1x_2))$ \`a
conjugaison g.d. pr\`es. Or
\begin{align*}
(x_0(x_0^2+x_1x_2):x_1^3:x_1(x_0^2+x_1x_2))=\rho(x_1:x_2:x_0)\tau( x_1:x_0:-x_2)
\end{align*}
donc $f$ est de la forme
$A\rho(x_1:x_2:x_0)\tau B$ o\`u $A,$ $B$ d\'esignent deux automorphismes
de $\mathbb{P}^2(\mathbb{C}).$
\end{proof}

\begin{lem}
{\sl Toute transformation birationnelle purement cubique ayant pour
configuration $\{\mathsf{9}\}$ s'\'ecrit \`a conjugaison
gauche-droite pr\`es
\begin{align*}
& \rho(x_0:x_0+x_1:x_2)\rho;
\end{align*}
en particulier $\sigma$ apparait au plus quatre fois dans sa
d\'ecomposition de \textsc{N\oe ther}.}
\end{lem}

\begin{proof}[\sl D\'emonstration]
On a vu que toute transformation birationnelle purement cubique
dont le lieu exceptionnel est r\'eduit \`a une conique et une
droite non tangente \`a cette conique s'\'ecrit \`a conjugaison
g.d. pr\`es
\begin{align*}
(x_1^2x_2:x_0(x_0x_2+x_1^2):x_1(x_0x_2+x_1^2)).
\end{align*}
Or
\begin{align*}
(x_1^2x_2:x_0(x_0x_2+x_1^2):x_1(
x_0x_2+x_1^2))=(x_1:x_0:x_2)\rho(x_0:x_0+x_1:x_2)\rho(x_0: x_2:x_1)
\end{align*}
donc tout \'el\'ement de $\mathring{\mathrm{B}}\mathrm{ir}_3$
poss\'edant la configuration $\{\mathsf{9}\}$ est du type
\begin{align*}
& A\rho(x_0:x_0+x_1:x_2)\rho B, && A,\hspace{1mm} B\in\mathrm{PGL}_3(\mathbb{C}).
\end{align*}
\end{proof}

\begin{lem}
{\sl La d\'ecomposition de \textsc{N\oe ther} de toute transformation de
\textsc{Cremona} cubique poss\'edant la configuration
$\{\mathsf{10}\}$ fait intervenir au plus six fois l'involution de
\textsc{Cremona}.}
\end{lem}

\begin{proof}[\sl D\'emonstration]
Une transformation birationnelle cubique dont la configuration est
$\{\mathsf{10}\}$ est, \`a conjugaison g.d. pr\`es de l'une des
formes suivantes
\begin{itemize}
\item $(x_0(x_1^2+x_0x_2):x_1(x_1^2+x_0x_2):x_0x_1x_2);$

\item $(x_0(x_0^2+x_1x_2):x_1(x_0^2+x_1x_2):x_0x_1x_2);$

\item $(x_0(x_1^2+x_0x_2):x_1(x_1^2+\gamma x_0x_1+x_0x_2)
:x_0x_1^2).$
\end{itemize}
On d\'eduit l'\'enonc\'e des \'egalit\'es
\begin{align*}
& (x_0(x_1^2+x_0x_2):x_1(x_1^2+x_0x_2):x_0x_1x_2)=(x_2:x_0:x_1)\sigma(x_0:x_0+x_1:x_2)\rho(x_0:x_2:x_1),\\
& (x_0(x_0^2+x_1x_2):x_1(x_0^2+x_1x_2):x_0x_1x_2)=(x_0:x_2:x_1)\sigma(x_0:x_0+x_1:x_2)\rho(x_1:x_2:x_0),
\end{align*}
et
\begin{align*}
& (x_0(x_1^2+x_0x_2):x_1(x_1^2+x_0x_2):x_0x_1^2)
 =\sigma (x_1:x_2:x_0+x_1)\rho (x_2:x_0:x_1).
\end{align*}
\end{proof}

\begin{lem}
{\sl La d\'ecomposition de \textsc{N\oe ther} d'un \'el\'ement de
$\mathring{\mathrm{B}}\mathrm{ir}_3$ poss\'edant la configuration
$\{\mathsf{11}\}$ fait intervenir au plus trois fois l'involution
de \textsc{Cremona}.}
\end{lem}

\begin{proof}[\sl D\'emonstration]
Soit $f$ une transformation birationnelle purement cubique qui
contracte exactement deux droites et une conique. On a
l'alternative suivante
\begin{itemize}
\item l'une des droites est de multiplicit\'e $3,$ l'autre de
multiplicit\'e $1$ et $f$ s'\'ecrit \`a conjugaison g.d. pr\`es
\begin{align*}
(x_0(x_0^2+x_1x_2+x_0x_2):x_1(x_0^2+x_1x_2+x_0x_2):x_0x_1x_2);
\end{align*}

\item les deux droites sont toutes deux de multiplicit\'e $2.$
\end{itemize}

Dans les deux cas les \'el\'ements pr\'esentant cette
configuration forment une orbite gauche-droite (\emph{voir} \S \ref{dimor}).

Dans la premi\`ere \'eventualit\'e, on constate que la
d\'ecomposition de \textsc{N\oe ther} de $f$ fait intervenir au plus trois
fois l'involution de \textsc{Cremona} puisque

\begin{align*}
(x_0(x_0^2+x_1x_2
+x_0x_2):x_1(x_0^2+x_1x_2+x_0x_2):x_0x_1x_2)=(x_1:x_2:x_0)
\sigma A\rho(x_1:x_2:x_0).
\end{align*}
o\`u $A=(x_0+x_1+x_2:x_0:x_2).$

Int\'eressons-nous \`a la seconde possibilit\'e. Un
calcul montre que $\sigma(x_1:x_2:x_0+x_1+x_2)\sigma,$
qui s'\'ecrit aussi
\begin{align*}
(x_1(x_1x_2+x_0
x_2+x_0x_1):x_2(x_1x_2+x_0x_2+x_0x_1):x_0x_1x_2),
\end{align*}
contracte la conique $x_1x_2+x_0x_2+x_0x_1=0$ et les droites
$x_1=0,$ $x_2=0;$ on note que chacune de ces droites est de
multiplicit\'e $2$ dans la d\'ecomposition de $\mathrm{det}\hspace{1mm}\mathrm{jac}.$
\end{proof}

\begin{lem}
{\sl La d\'ecomposition de \textsc{N\oe ther} des transformations birationnelles
purement cubiques dont le lieu exceptionnel est r\'eduit \`a une
conique, deux droites et une droite tangente \`a la conique
($\mathsf{\{12\}}$) fait intervenir au plus six fois l'involution
de \textsc{Cremona}.}
\end{lem}

\begin{proof}[\sl D\'emonstration]
\`A conjugaison g.d. pr\`es nous sommes dans l'une des situations
suivantes
\begin{itemize}
\item la droite tangente \`a la conique est de multiplicit\'e $2$
et
\begin{align*}
f=(x_0(x_0^2+x_1x_2):x_1(x_0^2+x_1x_2):x_0x_1(x_0-x_1));
\end{align*}

\item la droite tangente \`a la conique est de multiplicit\'e $1$
et $f$ s'\'ecrit
\begin{align*}
(x_0(x_0^2+x_0x_1+x_1x_2):x_1(x_0^2+x_0x_1+x_1x_2):x_0x_1x_2).
\end{align*}
\end{itemize}

\'Etudions la premi\`ere \'eventualit\'e. La
d\'ecomposition de $f$ compte alors au plus six fois l'involution
de \textsc{Cremona} puisque

\begin{small}
\begin{align*}
(x_0(x_0^2+x_1x_2):x_1(x_0^2+x_1x_2):x_0x_1(x_0-x_1))=(x_0:x_0-x_1:x_2)
\sigma(x_1-x_0:x_1:x_2)\tau(x_1:x_0:-x_2).
\end{align*}
\end{small}

Supposons que la droite tangente \`a la conique soit de
multiplicit\'e $1.$ On constate que

\begin{small}
\begin{align*}
\sigma(x_1+x_2:x_2:x_0+x_1+x_2)
\sigma=(x_1(x_1x_2+x_0x_2+x_0x_1):(x_2+x_1)(x_1x_2+x_0x_2+x_0x_1):x_0x_1(x_2+x_1))
\end{align*}
\end{small}

\noindent contracte la
conique $x_1x_2+x_0x_2+x_0x_1=0,$ les droites $x_1+x_2=0,$ $x_2=0$ et $x_1=0;$ cette
derni\`ere droite \'etant de multiplicit\'e $2.$
\end{proof}

\begin{lem}
{\sl Pour chaque \'el\'ement $f$ de
$\mathring{\mathrm{B}}\mathrm{ir}_3$ poss\'edant la configuration
$\{\mathsf{13}\}$ il existe $A$ dans $\mathrm{PGL}_3(\mathbb{C})$ tel que
$f=\sigma A\sigma.$}
\end{lem}

\begin{proof}[\sl D\'emonstration]
Une transformation birationnelle purement cubique $f$ de
configuration $\{\mathsf{13}\}$ est, \`a conjugaison g.d. pr\`es,
du type suivant
\begin{align*}
& (x_0(x_1^2+\gamma x_0x_1+x_0x_2+x_1x_2):x_1(x_1^2+
\gamma x_0x_1+x_0x_2+x_1x_2):x_0x_1x_2), &&\gamma\not=0,\hspace{1mm} 1.
\end{align*}

On constate que $(x_0(x_1^2+\gamma x_0x_1+x_0x_2+x_1x_2):x_1(x_1^2+
\gamma x_0x_1+x_0x_2+x_1x_2):x_0x_1x_2)$ s'\'ecrit aussi
\begin{align*}
&\sigma(-\gamma x_2:x_1+x_2:-\gamma^2x_0+x_1+(1-\gamma)x_2)\sigma
\left(x_2:-\frac{x_1}{\gamma}:x_0+\frac{x_1}{\gamma}\right).
\end{align*}
\end{proof}

\begin{lem}
{\sl Soit $f$ une transformation birationnelle purement cubique de
configuration~$\mathsf{\{14\}};$ la d\'ecomposition de \textsc{N\oe ther}
de $f$ fait intervenir au plus trois fois l'involution $\sigma.$}
\end{lem}

\begin{proof}[\sl D\'emonstration]
D'apr\`es ce qui pr\'ec\`ede on a l'alternative
\begin{itemize}
\item $f=(x_0(x_0^2+x_1^2+\gamma x_0x_1+x_0x_2):x_1(x_0^2+x_1^2+\gamma x_0x_1
+x_0x_2):x_0x_1x_2),$ $\gamma^2\not=4;$

\item $f=(x_0(x_0^2+x_1x_2+x_0x_2):x_1(x_0^2+x_1x_2+x_0x_2):x_0x_1
(x_0-x_1)).$
\end{itemize}

Dans le premier cas de figure on constate que $f=\sigma A
\sigma B$ avec
\begin{align*}
& A=(-x_0-x_1:\gamma_-x_0+\gamma_+x_1:\gamma_-x_0+\gamma_+x_1+(
\gamma^2-4)x_2), \\
& B=(-\delta(x_0+\gamma_-x_1):\delta(x_0+\gamma_+
x_1):x_2), &&\delta=\frac{1}{\sqrt{\gamma^2-4}}.
\end{align*}

Dans la seconde \'eventualit\'e on a $f=\sigma A\rho B,$ les
automorphismes $A$ et $B$ de $\mathbb{P}^2(\mathbb{C})$ \'etant donn\'es par
\begin{align*}
& A=(-2x_1+4x_2:2x_1+4x_2:-4x_0+x_1+4x_2), \\
& B=\left(\frac{x_0}{4}
+\frac{x_1}{4}+x_2:-\frac{x_0}{4}-\frac{x_1}{4}:-\frac{x_0}{2}+\frac{x_1}{2}\right).
\end{align*}
\end{proof}

On en d\'eduit l'\'enonc\'e suivant.

\begin{thm}
{\sl La d\'ecomposition de \textsc{N\oe ther} d'une transformation birationnelle
cubique fait intervenir au plus huit fois l'involution de
\textsc{Cremona}.}
\end{thm}

\section{Transformations birationnelles cubiques
et feuilletages}\hspace{1mm}

On mime ce qui a \'et\'e fait au Chapitre \ref{feuilfeuil}
pour le cas quadratique avec des notations analogues.
Puisque $\mathrm{Rat}_3\simeq\mathbb{P}^{29}(\mathbb{C})$ et $\mathcal{F}_3
\simeq\mathbb{P}^{23}(\mathbb{C})$ on a
\begin{eqnarray}
\mathcal{F}(\cdot)\hspace{1mm}\colon\hspace{1mm}\mathbb{P}^{29}(\mathbb{C})
&\to&\mathbb{P}^{23}(\mathbb{C})\nonumber \\
f=(f_0:f_1:f_2)&\mapsto&(x_1f_2-x_2f_1)\mathrm{d}x_0+(x_2f_0-x_0f_2)\mathrm{d}x_1
+(x_0f_1-x_1f_0)\mathrm{d}x_2\nonumber
\end{eqnarray}
Soit $f$ dans $\mathring{\mathrm{B}}\mathrm{ir}_3.$ Lorsque le
feuilletage $\mathcal{F}(f)$ est vraiment de degr\'e $3,$ {\it
i.e.} lorsque $\mathcal{F}(f)$ ne s'annule pas sur une
hypersurface, les points singuliers de $\mathcal{F}(f)$ sont la
r\'eunion des points d'ind\'etermination et des points fixes de
$f.$ Alors qu'en degr\'e $2$ cette m\^eme application restreinte
\`a $\mathrm{Bir}_2$ est essentiellement surjective, il n'en est pas de
m\^eme en degr\'e $3.$ Comme $\mathring{\mathrm{B}}\mathrm{ir}_3$ est de dimension $17,$ son
image est a priori de dimension au plus $17.$

Rappelons que $\mathscr{X}$ est l'orbite des transformations
$\xi_{\gamma,\delta}$ de la forme
\begin{align*}
(x_0(x_0^2+x_1^2+\gamma x_0x_1+\delta x_0x_2+x_1x_2): x_1(x_0^2
+x_1^2+\gamma x_0x_1+\delta x_0x_2+x_1x_2):x_0x_1x_2)
\end{align*}
induisant la configuration $\mathsf{\{15\}}.$ Un \'el\'ement de
$\mathscr{X}$ s'\'ecrit $B\xi_{\gamma,\delta}C$ pour certains $B,$ $C$
dans~$\mathrm{GL}_3(\mathbb{C})$ et $\gamma,$ $\delta$ dans $\mathbb{C}.$ Notons
que $B\xi_{\gamma,\delta}C$ est conjugu\'e \`a
$CB\xi_{\gamma,\delta}$ et par suite $\mathcal{F}(B \xi_{\gamma,\delta}C)$
est lin\'eairement conjugu\'e \`a $\mathcal{F}(CB \xi_{\gamma,\delta}).$
Dans la suite on se concentre donc sur les transformations du type~$A\xi_{\gamma,\delta}.$ Elles s'\'ecrivent
\begin{align*}
(a_0x_0Q+b_0x_1Q+c_0x_0x_1x_2: a_1x_0Q+b_1x_1Q+c_1x_0x_1x_2:a_2
x_0Q+b_2x_1Q+c_2x_0x_1x_2)
\end{align*}
avec $A=\left[
\begin{array}{ccc}
a_0 & b_0 & c_0 \\
a_1 & b_1 & c_1 \\
a_2 & b_2 & c_2
\end{array}\right]$ et $Q=x_0^2+x_1^2+\gamma x_0x_1+\delta x_0x_2+x_1x_2.$ Un calcul direct
montre que~$\mathcal{F}(A\xi_{\gamma,\delta})$ est donn\'e par la $1$-forme
\begin{eqnarray}
\omega&=& \left(x_1(a_2x_0Q+b_2x_1Q+c_2x_0x_1x_2)-x_2(a_1x_0Q+b_1x_1Q+c_1x_0x_1x_2)
\right)\mathrm{d}x_0+ \nonumber\\
\hspace{3mm}& & \left(x_2(a_0x_0Q+b_0x_1Q+c_0x_0x_1x_2)-x_0(a_2x_0Q+b_2x_1Q+c_2x_0x_1x_2)
\right)\mathrm{d}x_1+\nonumber\\
\hspace{3mm}& & \left(x_0(a_1x_0Q+b_1x_1Q+c_1x_0x_1x_2)-x_1(a_0x_0Q+b_0x_1Q+c_0x_0x_1x_2)
\right)\mathrm{d}x_2\nonumber
\end{eqnarray}
On remarque que $\omega,$ \'ecrit dans la carte affine $x_2=1,$ a
son jet d'ordre $1$ nul en $(0,0);$ ce type de feuilletages a des
propri\'et\'es dynamiques sp\'eciales (\cite{Croco}). En fait son
$2$-jet en $(0,0)$ est donn\'e par
\begin{align*}
(-\delta a_1x_0^2+(a_1+b_1\delta+c_1)x_0x_1+b_1x_1^2)\mathrm{d}x_0
+(\delta a_0x_0^2+(a_0+b_0\delta+c_0)x_0x_1+b_0x_1^2)\mathrm{d}x_1.
\end{align*}
En particulier, toujours pour des param\`etres $\gamma,$ $\delta$ et $A$
g\'en\'eriques, la multiplicit\'e du point singulier $(0,0)$ pour
$\omega$ est $4.$ Ceci semble indiquer qu'il y a g\'en\'eriquement
$10$ points singuliers. Pour confirmer cel\`a on constate que si
$A=\mathrm{id},$ les points singuliers de $\xi_{\mathrm{id},\gamma,\delta}$
pour $\gamma,$ $\delta$ g\'en\'eriques sont au nombre de $9$ (avec la
notation habituelle $\gamma_\pm=\frac{\gamma\pm\sqrt{\gamma^2-4}}{2}$)
\begin{align*}
& (1:1:0), && (1:-1:0), && (1:-\gamma_+:0), && (1:-\gamma_-:0), && (0:0:1),
\end{align*}
\begin{align*}
& (0:-1:1), && (-\delta:0:1), && ((1+\delta):(1+\delta):-(1+\gamma)), && (\delta-1:\delta-1:
1-\gamma).
\end{align*}

On remarque alors que le $2$-jet en $(0,0)$ de
$\omega_{\mathrm{id},\gamma,\delta}$ n'est pas de type g\'en\'eral, en
particulier sa multiplicit\'e $\mu$ est strictement plus grande
que $4.$ Un feuilletage de degr\'e $3$ de $\mathbb{P}^2(\mathbb{C})$ a $13$ points
singuliers compt\'es avec multiplicit\'e. Ceci implique que les
autres points singuliers sont de multiplicit\'e $1$ et que
$\mu=5.$ Si $A$ est voisin de l'identit\'e et g\'en\'erique, les
huit points simples subsistent, le point $(0,0)$ devient de
multiplicit\'e $4$ donc un dixi\`eme point singulier appara\^it.
D'o\`u la:

\begin{pro}
{\sl Un  \'el\'ement g\'en\'erique de $\mathscr{X}$ poss\`ede dix points
sp\'eciaux. Cinq sont des points d'ind\'etermination dont l'un est
de multiplicit\'e $4,$ les cinq autres sont des points fixes
simples.}
\end{pro}

En g\'en\'eral pour un \'el\'ement $f$ de $\mathscr{X}$ au moins
quatre des cinq points fixes sont en position g\'en\'erale (pas
d'alignement $3$ \`a $3$). Pour cel\`a il suffit de construire un
exemple. On se donne une transformation $\xi_{\gamma,\delta}$ et
on choisit quatre points $p_0,$ $p_1,$ $p_2$ et $p_3$ en position
g\'en\'erale tels que $\xi_{\gamma,\delta}$ y soit un
diff\'eomorphisme local. On suppose que les
$q_i:=\xi_{\gamma,\delta}(p_i)$ sont aussi en position
g\'en\'erale. Alors il existe un unique automorphisme lin\'eaire
$A$ tel que $Aq_i=p_i.$ La transformation birationnelle
$A\xi_{\gamma,\delta}\in\mathscr{X}$ poss\`ede quatre points fixes en
position g\'en\'erale.

Consid\'erons maintenant deux transformations
g\'en\'eriques de $\mathscr{X}$ ayant leurs cinq points
d'ind\'etermination communs compt\'es avec multiplicit\'e ({\it
i.e.} le point de multiplicit\'e $4$ est commun aux deux
transformations) ainsi que $4$ points fixes communs; comme les
transformations sont g\'en\'eriques ces $4$ points fixes sont en
position g\'en\'erale. Remarquons que les points
d'ind\'etermination d\'eterminent les $5$ courbes contract\'ees:
ce sont donc les m\^emes pour les deux transformations. On se
ram\`ene bien s\^ur au cas o\`u les deux transformations
s'\'ecrivent $A\xi_{\gamma,\delta}$ et~$B\xi_{\gamma',\delta'}.$ Mais si $\xi_{\gamma,\delta}$
et $\xi_{\gamma',\delta'}$ ont m\^emes point d'ind\'etermination elles sont
\'egales, c'est un calcul \'el\'ementaire. Les quatre points fixes
communs nous permettent d'affirmer que $A=B.$ Comme en degr\'e $2$
on obtient le:

\begin{thm}
{\sl Une transformation birationnelle de degr\'e $3$ g\'en\'erique
est d\'etermin\'ee par la position de ses $5$ points fixes et de
ses $5$ points d'ind\'etermination affect\'es de leur
multiplicit\'e.}
\end{thm}

Dans le cas quadratique
pour $7$ points en position g\'en\'erale on peut
trouver un \'el\'ement de~$\mathrm{Bir}_2$ ayant ces $7$ points
sp\'eciaux; les configurations de $10$ points \'etant de
dimension $20$ il n'est pas possible de r\'ealiser un $10$-uplet
de points g\'en\'eriques comme points sp\'eciaux d'un \'el\'ement
de~$\mathrm{Bir}_3.$ Pour une transformation de $\mathring{\mathrm{B}}
\mathrm{ir}_3$ g\'en\'erique la position de $9$ points sp\'eciaux
d\'etermine celle du $10$\`eme.

\appendix

\chapter{D\'emonstration de la Proposition
7.14}\label{demdem}

Nous pr\'esentons ici une d\'emonstration de la Proposition
\ref{listee}. Il s'agit de r\'eduire le nombre de param\`etres
apparaissant dans les diff\'erents mod\`eles du Lemme \ref{vent}
\`a conjugaison g.d. pr\`es. Pour cel\`a \`a chaque transformation
on associe deux invariants gauche-droite
\begin{itemize}
\item la configuration des courbes contract\'ees avec
multiplicit\'e, celle qui appara\^it dans la d\'ecomposition en
facteurs irr\'eductibles du d\'eterminant jacobien;

\item la configuration des points images des courbes contract\'ees
qui n'est rien d'autre que celle des points d'ind\'etermination de
$f^{-1}.$
\end{itemize}

 Nous allons voir que pour certaines valeurs
sp\'ecifiques des param\`etres quelques mod\`eles faisant partie
de types distincts (on fait ici allusion aux types
$\mathfrak{(a)},$ $\mathfrak{(b)},$ $\mathfrak{(c)}$ et
$\mathfrak{(d)}$ introduits au Lemme \ref{vent}) pr\'esentent des configurations
d'ensembles $\mathrm{Exc}\hspace{1mm} f$ et $\mathrm{Ind} \hspace{1mm} f^{-1}$
identiques. \bigskip

Nous avons d\'ecoup\'e la preuve en deux parties. La premi\`ere
(\S \ref{reduc}) est une approche grossi\`ere de la
classification qui consiste modulo quelques conjugaisons
imm\'ediates \`a lister les ensembles $\mathrm{Exc}\hspace{1mm} f$ lorsque
$f$ d\'ecrit les diff\'erents mod\`eles du Lemme \ref{vent}. La seconde
(\S \ref{compa}) consiste \`a comparer les transformations qui
pr\'esentent une m\^eme configuration d'ensemble exceptionnel.

\section{Premi\`ere \'etape}\label{reduc}\hspace{1mm}

\subsection{\'El\'ements de la forme $\mathfrak{(a)}$}\hspace{1mm}

Soit $f$ une transformation de type $\mathfrak{(a)}.$
Suivant que le param\`etre $\alpha$ est nul ou non, on se ram\`ene
\`a $\alpha=0$ ou $\alpha=1.$

$\mathfrak{(a)}$ \textbf{\textit{i.}} Commen\c{c}ons par
\'etudier l'\'eventualit\'e $\alpha=0.$ La transformation $f$ \'etant
dans $\mathring{\mathrm{B}}\mathrm{ir}_3$ le param\`etre $\delta$ et le couple $(\beta,\varepsilon)$
sont non nuls; \`a conjugaison g.d.
pr\`es il suffit donc de traiter les possibilit\'es suivantes
\begin{align*}
& (\beta,\delta,\varepsilon)=(0,1,1),&& (\beta,\delta,\varepsilon)=(1,1,0),
&& (\beta,\delta,\varepsilon)=(1,1,1).
\end{align*}

$\mathfrak{(a)}$ \textbf{\textit{i.1.}} Supposons que
$(\beta,\delta,\varepsilon)=(0,1,1).$ Notons que si $\gamma=0,$ la
transformation $f$ est quadratique, par suite $\gamma=1$ \`a
conjugaison g.d. pr\`es et
\begin{align*}
f=(x_0(x_0x_1+x_0x_2+x_1x_2):x_1(x_0x_1+x_0x_2+x_1x_2):
x_0x_1x_2)
\end{align*}
pour laquelle on observe la configuration $\mathsf{\{11\}}$
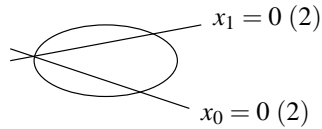
\begin{figure}[H]
\begin{center}
\input{config1.pstex_t}
\end{center}
\end{figure}

Rappelons que les indices (\hspace{1mm}) indiquent la multiplicit\'e qui
appara\^it dans la d\'ecomposition en facteurs irr\'eductibles du
d\'eterminant jacobien. Cette multiplicit\'e, qui est \'evidemment
un invariant g.d., n'est mentionn\'ee que lorsqu'elle est
strictement sup\'erieure \`a $1.$

$\mathfrak{(a)}$ \textbf{\textit{i.2.}} Si
$(\beta,\delta,\varepsilon)=(1,1,0),$ alors
\begin{align*}
f=(x_0(x_1^2+\gamma x_0x_1+x_0x_2):x_1(x_1^2+\gamma x_0x_1+x_0x_2):x_0x_1x_2).
\end{align*}
On a l'alternative suivante
\begin{itemize}
\item ou bien $\gamma=0,$ l'ensemble $\mathrm{Exc}\hspace{1mm} f$ est du
type $\mathsf{\{10\}}$
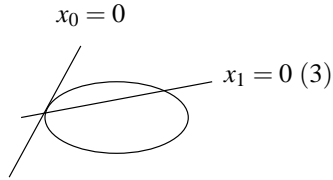
\begin{figure}[H]
\begin{center}
\input{config2.pstex_t}
\end{center}
\end{figure}

\item ou bien $\gamma\not=0,$ on peut supposer que $\gamma$ vaut
$1$ et on constate que $f$ pr\'esente la configuration~$\mathsf{\{12\}}$
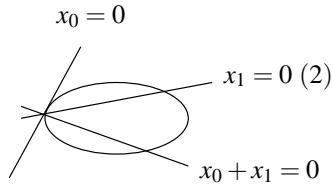
\begin{figure}[H]
\begin{center}
\input{config3.pstex_t}
\end{center}
\end{figure}
\end{itemize}

$\mathfrak{(a)}$ \textbf{\textit{i.3.}} Lorsque
$(\beta,\delta,\varepsilon)=(1,1,1)$ on a
\begin{align*}
f=(x_0(x_1^2+\gamma x_0x_1
+x_0x_2+x_1x_2):x_1(x_1^2+\gamma x_0x_1+x_0x_2+x_1x_2):x_0x_1x_2).
\end{align*}
La conique d'\'equation $x_1^2+\gamma x_0x_1 +x_0x_2+x_1x_2=0$ est
lisse si et seulement si $\gamma\not=1.$ Les courbes contract\'ees
pr\'esentent la configuration $\mathsf{\{11\}},$ resp. $\mathsf{\{13\}}$
lorsque la conique est lisse
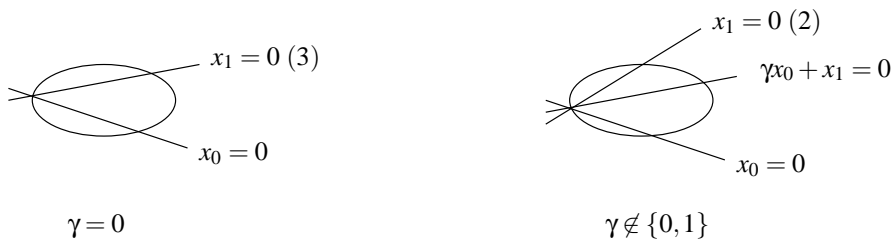
\begin{figure}[H]
\begin{center}
\input{config4.pstex_t}
\end{center}
\end{figure}

\bigskip

\noindent et la configuration $\mathsf{\{5\}}$ lorsqu'elle ne
l'est pas

\bigskip
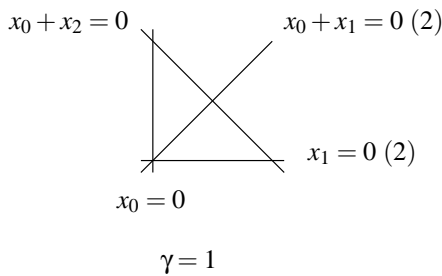
\begin{figure}[H]
\begin{center}
\input{config4b.pstex_t}
\end{center}
\end{figure}

\bigskip

$\mathfrak{(a)}$ \textbf{\textit{ii.}} Supposons que
$\alpha=1.$ On observe que le couple $(\beta,\varepsilon)$
est non nul; comme pr\'ec\'edemment \`a conjugaison g.d. pr\`es il suffit donc de
consid\'erer les possibilit\'es
\begin{align*}
& (\beta,\varepsilon)=(1,0),&& (\beta,\varepsilon)=(0,1), && (\beta,\varepsilon)=(1,1).
\end{align*}

$\mathfrak{(a)}$ \textbf{\textit{ii.1.}} \'Etudions
l'\'eventualit\'e $(\beta,\varepsilon)=(1,0)$
\begin{align*}
f=(x_0(x_0^2+x_1^2+\gamma x_0x_1+\delta x_0x_2):x_1(x_0^2+x_1^2+
\gamma x_0x_1+\delta x_0x_2):x_0x_1x_2).
\end{align*}
Posons
\begin{align*}
& \gamma_+:=\frac{\gamma+\sqrt{\gamma^2-4}}{2},&&\gamma_-:=\frac{\gamma-
\sqrt{\gamma^2-4}}{2}.
\end{align*}
La conique d'\'equation $x_0^2+x_1^2+\gamma x_0x_1+\delta
x_0x_2=0$ est lisse si et seulement si $\delta\not=0;$ nous allons
donc distinguer pour d\'ecrire les ensembles $\mathrm{Exc}\hspace{1mm} f$
le cas $\delta$ nul du cas $\delta$ non nul. Si $\delta$ est non
nul, le lieu exceptionnel de $f$ contient une conique lisse et
on observe les configurations $\mathsf{\{12\}}$ et $\mathsf{\{14\}}$
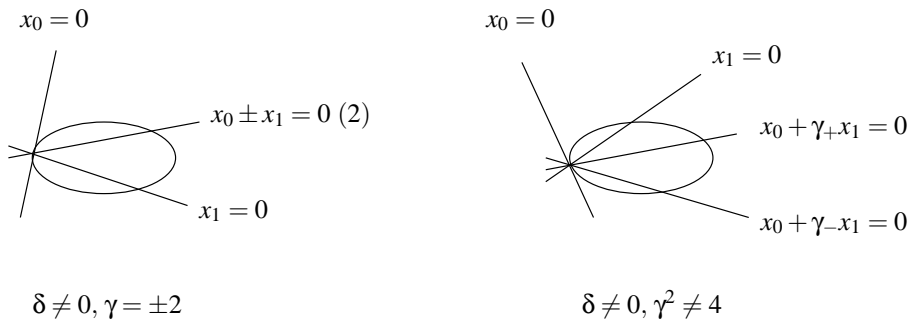
\begin{figure}[H]
\begin{center}
\input{config5.pstex_t}
\end{center}
\end{figure}

Lorsque $\delta$ est nul on constate que $f$ pr\'esente
la configuration $\mathsf{\{6\}}$
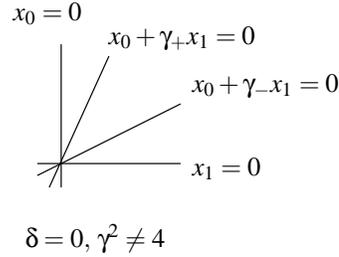
\bigskip
\begin{figure}[H]
\begin{center}
\input{config5b.pstex_t}
\end{center}
\end{figure}

Remarquons que pour $\gamma^2=4$ et $\delta=0,$ la conique est
une droite double, donc non r\'eduite, ce qui n'entre pas dans la situation
\'etudi\'ee.

$\mathfrak{(a)}$ \textbf{\textit{ii.2.}} Si $(\beta,\varepsilon)=(
0,1),$ alors $f$ s'\'ecrit
\begin{align*}
(x_0(x_0^2+\gamma x_0x_1+\delta x_0x_2+x_1x_2):x_1(x_0^2+\gamma x_0x_1
+\delta x_0x_2+x_1x_2):x_0x_1x_2).
\end{align*}
La conique d'\'equation $x_0^2+\gamma x_0x_1+\delta x_0x_2+x_1x_2=0$ est lisse
si et seulement si $\gamma\delta\not=1.$ On observe les configurations
$\mathsf{\{13\}},$ resp. $\mathsf{\{12\}}$
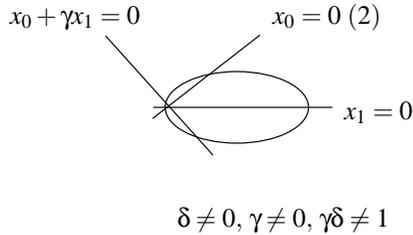
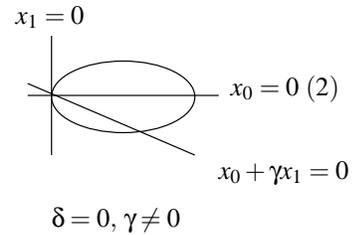
\begin{figure}[H]
\begin{center}
\input{config6.pstex_t}
\end{center}
\end{figure}

\bigskip

Lorsque $\gamma$ est nul, le nombre de droites diminue et on
observe les configurations $\mathsf{\{13\}},$ resp.~$\mathsf{\{12\}}$

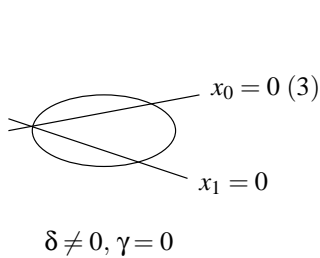
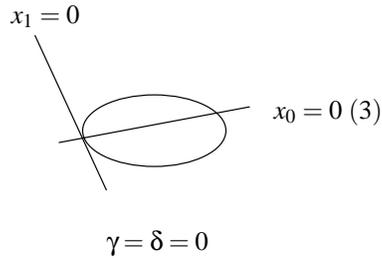
\bigskip

\begin{figure}[H]
\begin{center}
\input{config6a.pstex_t}
\end{center}
\end{figure}

\bigskip

Enfin si $\gamma\delta=1$ on constate que $f$ ne
contracte plus de conique lisse; plus pr\'ecis\'ement $f$ a pour
configuration $\mathsf{\{5\}}$

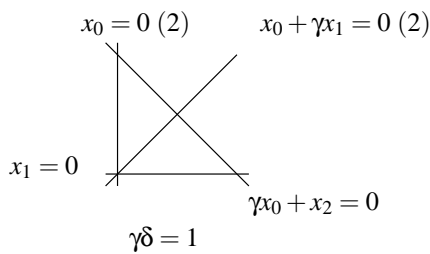
\bigskip

\begin{figure}[H]
\begin{center}
\input{config6b.pstex_t}
\end{center}
\end{figure}

$\mathfrak{(a)}$ \textbf{\textit{iii.}} Enfin on
s'int\'eresse \`a la possibilit\'e $(\beta,\varepsilon)=(1,1);$ on a
\begin{align*}
f=(x_0(x_0^2+x_1^2+\gamma x_0x_1+\delta x_0x_2+x_1x_2):x_1
(x_0^2+x_1^2+\gamma x_0x_1+\delta x_0x_2 +x_1x_2):x_0x_1x_2).
\end{align*}
Notons que la conique $x_0^2+x_1^2+\gamma x_0x_1+\delta
x_0x_2+x_1x_2=0$ est lisse si et seulement si
$\delta\not\in\{\gamma_+,\gamma_-\}.$ Si $\delta$ est distinct de
$\gamma_+$ et $\gamma_-,$ les configurations des courbes
contract\'ees sont les configura\-tions~$\mathsf{\{15\}},\hspace{1mm}
\mathsf{\{13\}},\hspace{1mm}\mathsf{\{14\}}$ et $\mathsf{\{12\}}$
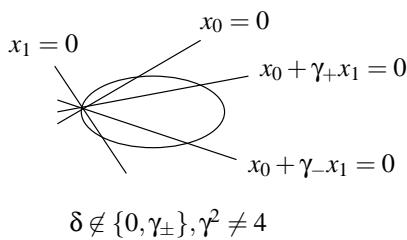
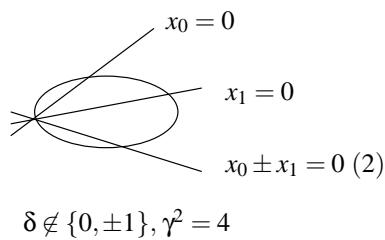
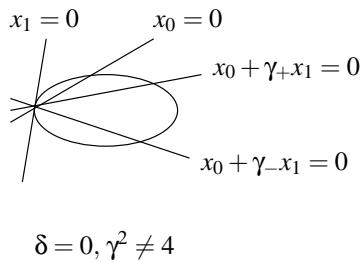
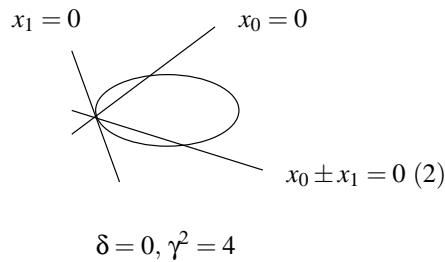
\begin{figure}[H]
\begin{center}
\input{config7.pstex_t}
\end{center}
\end{figure}

\begin{figure}[H]
\begin{center}
\input{config7b.pstex_t}
\end{center}
\end{figure}

Lorsque la conique d\'eg\'en\`ere on a ou bien
$\delta=\gamma_\pm,$ ou bien $(\gamma,\delta)=(\pm 2,\pm 1).$
Si on change~$x_0$ en $-x_0$ on voit que les configurations $(
\gamma,\delta)=(2,1)$ et $(\gamma,\delta)=(-2,-1)$ sont isomorphes,
tout comme les configurations $(\gamma,\delta)=(2,-1)$
et $(\gamma,\delta)=(-2,1);$ de sorte qu'il n'y a que les \'eventualit\'es
$(\gamma,\delta)=(2,\pm 1)$ \`a consid\'erer. On observe les
configurations $\mathsf{\{5\}}$ et $\mathsf{\{7\}}$
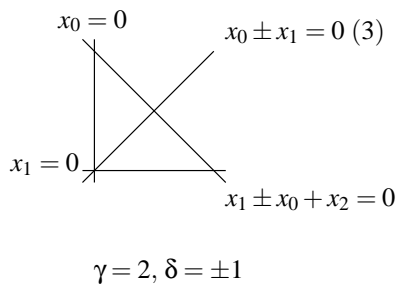
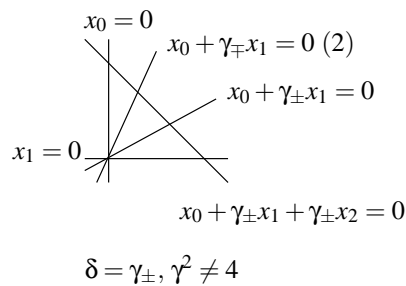
\begin{figure}[H]
\begin{center}
\input{config8.pstex_t}
\end{center}
\end{figure}

\subsection{\'El\'ements du type $\mathfrak{(b)}$}\hspace{1mm}

Soit $f$ un \'el\'ement de la forme $\mathfrak{(b)}.$

$\mathfrak{(b)}$ \textbf{\textit{i.}} Dans un premier
temps supposons que $\alpha$ soit nul: \`a conjugaison pr\`es on a
$\alpha=1.$ La transformation $f$ \'etant purement cubique on peut
supposer que $\delta=1$ et $(\beta, \varepsilon)\not=(0,0);$ on
distingue les possibilit\'es suivantes \`a conjugaison g.d.
pr\`es
\begin{align*}
& (\beta,\delta,\varepsilon)=(1,1,0), && (\beta,\delta,\varepsilon)=(0,1,1),
&& (\beta,\delta,\varepsilon)=(1,1,1).
\end{align*}

$\mathfrak{(b)}$ \textbf{\textit{i.1.}} Si
$(\beta,\delta,\varepsilon)=(1,1,0),$ alors $f=(x_0(x_1^2+\gamma
x_0x_1+ x_0x_2): x_1(x_1^2+\gamma x_0x_1+x_0x_2):x_0x_1^2)$ et la
configuration de $\mathrm{Exc}\hspace{1mm} f$ est la configuration
$\mathsf{\{10\}}$
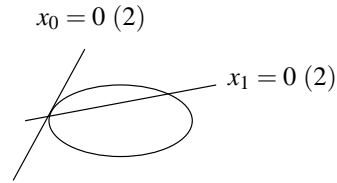
\begin{figure}[H]
\begin{center}
\input{config9.pstex_t}
\end{center}
\end{figure}

$\mathfrak{(b)}$ \textbf{\textit{i.2.}} Lorsque
$(\beta,\delta,\varepsilon)=(0,1,1)$ on a
\begin{align*}
& f=(x_0(\gamma x_0x_1+ x_0x_2+x_1x_2):
x_1(\gamma x_0x_1+x_0x_2+x_1x_2):x_0x_1^2).
\end{align*}
La conique $\gamma x_0x_1+x_0x_2+x_1x_2=0$ est lisse si et seulement
si $\gamma\not=0;$ nous avons donc les configurations $\mathsf{\{12\}}$
et $\mathsf{\{5\}}$
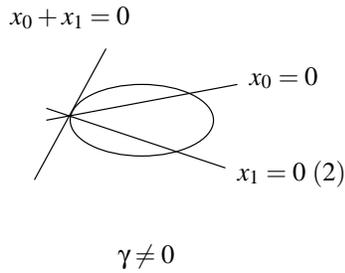
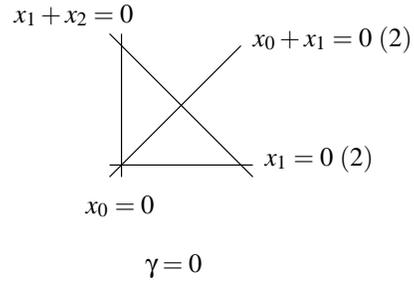
\begin{figure}[H]
\begin{center}
\input{config10.pstex_t}
\end{center}
\end{figure}

$\mathfrak{(b)}$ \textbf{\textit{i.3.}} Pour
$(\beta,\delta,\varepsilon)=(1,1,1)$ la transformation $f$ est du type
\begin{align*}
(x_0(x_1^2+\gamma x_0x_1+x_0x_2+x_1x_2):x_1(x_1^2+\gamma x_0x_1
+x_0x_2+x_1x_2):x_0x_1^2).
\end{align*}
La conique $x_1^2+\gamma x_0x_1+x_0x_2+x_1x_2=0$ \'etant lisse si et
seulement si $\gamma\not=1,$ on observe les configurations $\mathsf{\{12\}}$
et $\mathsf{\{5\}}$
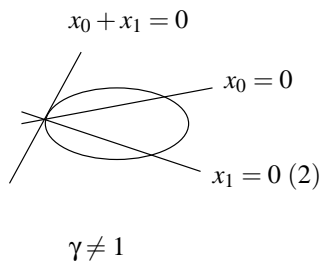
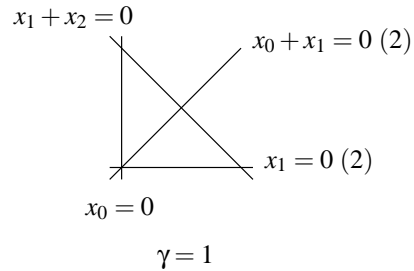
\begin{figure}[H]
\begin{center}
\input{config11.pstex_t}
\end{center}
\end{figure}

$\mathfrak{(b)}$ \textbf{\textit{ii.}} Supposons que
$\alpha=1;$ puisque $f$ appartient \`a $\mathring{\mathrm{B}}\mathrm{ir}_3$ on a $(\beta,\varepsilon)\not=(0,0);$
par suite \`a conjugaison pr\`es
\begin{align*}
(\beta,\varepsilon)\in\{(1,0),\hspace{1mm}(0,1),\hspace{1mm}(1,1)\}.
\end{align*}

$\mathfrak{(b)}$ \textbf{\textit{ii.1.}} Lorsque
$(\beta,\varepsilon)=(1,0)$ la transformation $f$ s'\'ecrit
\begin{align*}
(x_0(x_0^2+x_1^2+\gamma x_0x_1+\delta x_0x_2):x_1(x_0^2+ x_1^2+
\gamma x_0x_1+\delta x_0x_2): x_0x_1^2)
\end{align*}
et pr\'esente la configuration $\mathsf{\{10\}}$
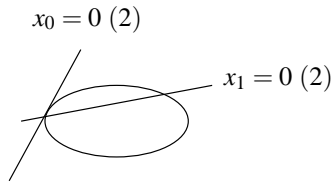
\begin{figure}[H]
\begin{center}
\input{config9.pstex_t}
\end{center}
\end{figure}

$\mathfrak{(b)}$ \textbf{\textit{ii.2.}} Si
$(\beta,\varepsilon)=(0,1)$ on a
\begin{align*}
f=(x_0(x_0^2+\gamma x_0x_1+\delta x_0x_2+x_1x_2):x_1(x_0^2+
\gamma x_0x_1+\delta x_0x_2+x_1x_2):x_0x_1^2).
\end{align*}
La conique d'\'equation $x_0^2+\gamma x_0x_1+\delta x_0x_2+x_1
x_2=0$ \'etant lisse si et seulement si $\gamma\delta\not=1$ on
observe pour $\gamma\delta\not=1$ les configurations $\mathsf{
\{10\}}$ et $\mathsf{\{12\}}$
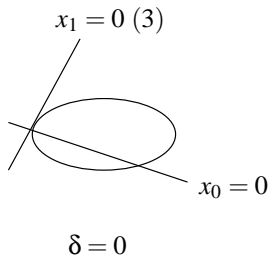
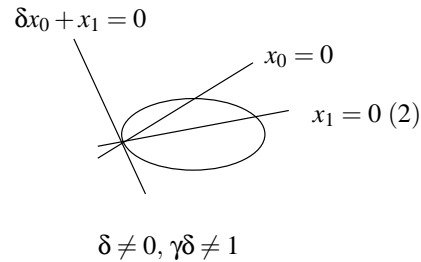
\begin{figure}[H]
\begin{center}
\input{config12.pstex_t}
\end{center}
\end{figure}

\bigskip

Si $\gamma\delta=1$ on obtient la configuration
$\mathsf{\{5\}}$

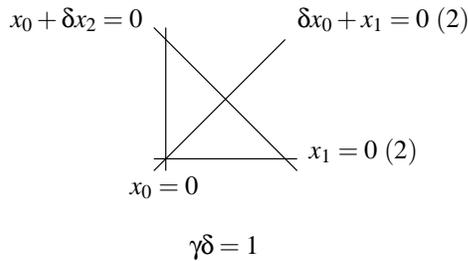
\bigskip
\begin{figure}[H]
\begin{center}
\input{config12b.pstex_t}
\end{center}
\end{figure}

$\mathfrak{(b)}$ \textbf{\textit{ii.3.}} Pour finir
lorsque $(\beta,\varepsilon)=(1,1)$ la transformation $f$ s'\'ecrit
\begin{align*}
(x_0(x_0^2+x_1^2+\gamma x_0x_1+\delta x_0x_2+x_1x_2):x_1(x_0^2
+x_1^2+\gamma x_0x_1+\delta x_0x_2+x_1x_2):x_0x_1^2).
\end{align*}
La conique d'\'equation $x_0^2+x_1^2+\gamma x_0x_1+\delta x_0x_2
+x_1x_2=0$ est lisse si et seulement si $\delta$ n'appartient pas \`a $\{\gamma_+,
\gamma_-\}.$ Dans ce cas $\mathrm{Exc}\hspace{1mm} f$ pr\'esente la
configuration $\mathsf{\{10\}}$ ou $\mathsf{\{12\}}$
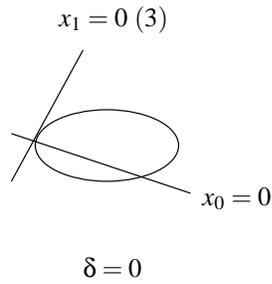
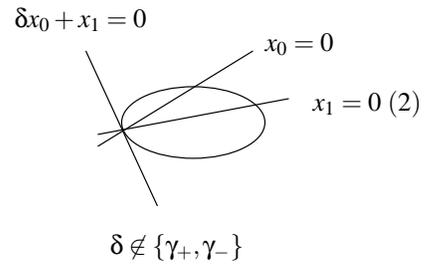
\begin{figure}[H]
\begin{center}
\input{config13.pstex_t}
\end{center}
\end{figure}

\bigskip

Lorsque $\delta=\gamma_\pm$ on
observe la configuration $\mathsf{\{5\}}$

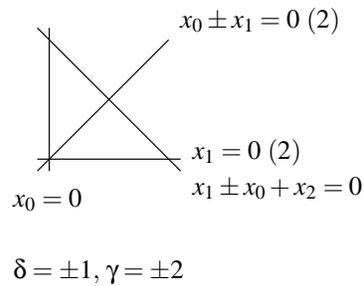
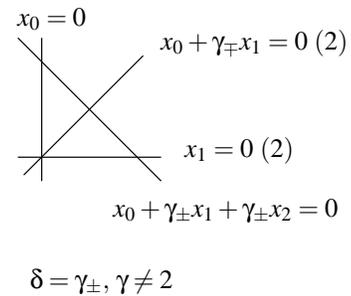
\bigskip
\begin{figure}[H]
\begin{center}
\input{config13b.pstex_t}
\end{center}
\end{figure}

\newpage

\subsection{\'El\'ements de la forme $\mathfrak{(c)}$ et $\mathfrak{(d)}$}\hspace{1mm}

Les transformations de type $\mathfrak{(c)}$ conduisent
aux configurations $\mathsf{\{14\}}$ et $\mathsf{\{12\}}$
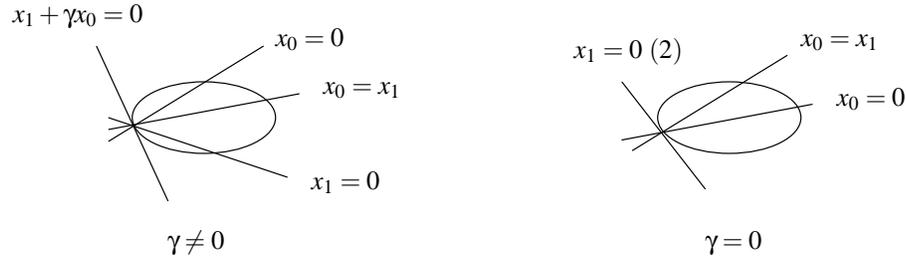
\begin{figure}[H]
\begin{center}
\input{config15.pstex_t}
\end{center}
\end{figure}

\bigskip

Enfin pour un \'el\'ement de la forme $\mathfrak{(d)}$
on a la configuration $\mathsf{\{7\}}$

\begin{figure}[H]
\begin{center}
\input{config14.pstex_t}
\end{center}
\end{figure}

\section{Deuxi\`eme \'etape}\label{compa}\hspace{1mm}

\subsection{Configuration $\mathsf{\{5\}}$}\hspace{1mm}

Les transformations du \S\ref{reduc} ayant la
configuration $\mathsf{\{5\}}$ sont les suivantes
\begin{align*}
& f_1=(x_0(x_1+x_2)(x_0+x_1):x_1(x_1+x_2)(x_0+x_1):x_0x_1x_2); \\
& f_2=(x_0(x_2+\gamma x_0)\left(\frac{x_0}{\gamma}+x_1\right):
x_1(x_2+\gamma x_0)\left(\frac{x_0}{\gamma}+x_1\right):x_0x_1x_2); \\
& f_3=(x_0(x_0+x_1)(x_0+x_1+x_2):x_1(x_0+x_1)(x_0+x_1+x_2):x_0x_1x_2); \\
& f_4=(x_0x_2(x_0+x_1):x_1x_2(x_0+x_1):x_0x_1^2); \\
& f_5=(x_0(x_0+x_1)(x_1+x_2):x_1(x_0+x_1)(x_1+x_2):x_0x_1^2); \\
& f_6=(x_0\left(\frac{x_0}{\gamma}+x_1\right)(x_2+\gamma x_0):
x_1\left(\frac{x_0}{\gamma}+x_1\right)(x_2+\gamma x_0):x_0x_1^2); \\
& f_7=(x_0(x_0+x_1)(x_0+x_1+x_2):x_1(x_0+x_1)(x_0+x_1+x_2):x_0x_1^2); \\
& f_8=(x_0(x_1+\gamma_+x_0)(x_1+\gamma_-x_0+x_2):x_1(x_1+\gamma_+
x_0)(x_1+\gamma_-x_0+x_2):x_0x_1^2).
\end{align*}

On remarque que d'une part $f_4,$ $f_5,$ $f_6,$ $f_7$ et
$f_8$ sont g.d. conjugu\'ees
\begin{align*}
& f_4=f_5(x_0:x_1:x_2-x_1); &&
f_4=\left(\frac{x_0}{\gamma}:x_1:\frac{x_2}{\gamma}\right)f_6(\gamma
x_0:x_1: x_2-\gamma^2x_0);\\
& f_4=f_7(x_0:x_1:x_2-x_0-x_1);&&
f_4=(\gamma_+x_0:x_1:\gamma_+x_2)f_8\left(\frac{x_0}{\gamma_+}:x_1:x_2
-\frac{\gamma_-}{\gamma_+}x_0-x_1\right);
\end{align*}

\noindent et que d'autre part $f_1$ et $f_2$ appartiennent \`a la
m\^eme g.d. orbite
\begin{align*}
f_1=(\gamma x_1:x_0:x_2)f_2\left(x_1:\frac{x_0}{\gamma}:\gamma
x_2\right).
\end{align*}

Les \'el\'ements $f_1$ et $f_4$ ne sont pas g.d.
conjugu\'es: $\#\hspace{1mm}\mathrm{Ind}\hspace{1mm} f_4\not=\#\hspace{1mm}\mathrm{Ind}\hspace{1mm}
f_1.$

Les transformations $f_1$ et $f_3$ ne sont pas g.d. conjugu\'ees: tous
les points d'ind\'etermination de $f_3$ sont des intersections de droites de
$\mathrm{Exc} f_3$ ce qui n'est pas le cas pour $f_1.$

\`A conjugaison g.d. pr\`es il reste donc seulement les
mod\`eles suivants pour la configuration~$\mathsf{\{5\}}$
\begin{align*}
& (x_0(x_1+x_2)(x_0+x_1):x_1(x_1+x_2)(x_0+x_1):x_0x_1x_2); \\
& (x_0(x_0+x_1)(x_0+x_1+x_2):x_1(x_0+x_1)(x_0+x_1+x_2):x_0x_1x_2); \\
& (x_0x_2(x_0+x_1):x_1x_2(x_0+x_1):x_0x_1^2).
\end{align*}

\subsection{Configuration $\mathsf{\{6\}}$}\hspace{1mm}

Un seul mod\`ele pr\'esente la configuration
$\mathsf{\{6\}}$
\begin{align*}
& (x_0(x_0^2+x_1^2+\gamma x_0x_1):x_1(x_0^2+x_1^2+\gamma x_0x_1)
:x_0x_1x_2), && \gamma^2\not=4.
\end{align*}

\subsection{Configuration $\mathsf{\{7\}}$}\hspace{1mm}

Deux mod\`eles du \S\ref{reduc} pr\'esentent la
configuration $\mathsf{\{7\}}$
\begin{footnotesize}
\begin{align*}
& f_1=(x_0(x_0^2+x_1^2+\gamma x_0x_1+\gamma_+ x_0x_2+x_1x_2):x_1(
x_0^2+x_1^2+\gamma x_0x_1+\gamma_+ x_0x_2+x_1x_2):x_0x_1x_2), &&
\gamma^2\not=4;\\
& f_2=(x_0x_2(x_1+\gamma x_0):x_1x_2(x_1+\gamma
x_0):x_0x_1(x_0-x_1)).&&
\end{align*}
\end{footnotesize}

Les transformations $f_1$ et $f_2$ ne sont pas g.d.
conjugu\'ees: trois des points de $f_2(\mathrm{Exc} f_2)$ sont
align\'es ce qui n'est pas le cas pour $f_1(\mathrm{Exc} f_1).$

\subsection{Configuration $\mathsf{\{10\}}$}\hspace{1mm}

Les \'el\'ements obtenus au \S\ref{reduc} ayant la
configuration $\mathsf{\{10\}}$ sont
\begin{align*}
& f_1=(x_0(x_1^2+x_0x_2):x_1(x_1^2+x_0x_2):x_0x_1x_2); \\
& f_2=(x_0(x_0^2+x_1x_2):x_1(x_0^2+x_1x_2):x_0x_1x_2); \\
& f_3=(x_0(x_1^2+\gamma x_0x_1+x_0x_2):x_1(x_1^2+\gamma
x_0x_1+x_0x_2):x_0x_1^2); \\
& f_4=(x_0(x_0^2+x_1^2+\gamma x_0x_1+\delta x_0x_2):x_1(x_0^2
+x_1^2+\gamma x_0x_1+\delta x_0x_2):x_0x_1^2); \\
& f_5=(x_0(x_0^2+\gamma x_0x_1+x_1x_2):x_1(x_0^2+\gamma x_0
x_1+x_1x_2):x_0x_1^2); \\
& f_6=(x_0(x_0^2+x_1^2+\gamma x_0x_1+x_1x_2):x_1(x_0^2+x_1^2
+\gamma x_0x_1+x_1x_2):x_0x_1^2).
\end{align*}

On constate que $f_1$ et $f_2$ (resp. $f_3$ et $f_4,$
resp. $f_5$ et $f_6$) sont dans la m\^eme g.d. orbite
\begin{align*}
& f_1=(x_1:x_0:x_2)f_2(x_1:x_0:x_2);&& f_5=f_6(x_0:x_1:x_2-x_1);&&
 f_3=f_4\left(x_0:x_1:\frac{x_2-x_0}{\delta}\right).
\end{align*}

L'\'el\'ement $f_5$ n'est conjugu\'e ni \`a $f_1,$ ni
\`a $f_2,$ ni \`a $f_4;$ la configuration de $\mathrm{Exc}
 f_5$ se distingue de celle des $f_i$ pour $1\leq i\leq 4$
 par la propri\'et\'e suivante: la droite de $\mathrm{Exc} f_5$
tangente \`a la conique de $\mathrm{Exc} f_5$ est de
multiplicit\'e $3.$

L'\'el\'ement $f_4$ n'est conjugu\'e ni \`a $f_1,$ ni \`a $f_2,$
ni \`a $f_5:$ les deux droites de $\mathrm{Exc} f_4$ apparaissent avec
multiplicit\'e $2$ ce qui n'est pas le cas pour les autres mod\`eles.

La transformation $f_1$ n'est pas conjugu\'ee \`a $f_3:$
la droite de  $\mathrm{Exc} f_1$ tangente \`a la conique de
$\mathrm{Exc} f_1$ est de multiplicit\'e $1$ ce qui n'est pas le
cas pour $f_3.$

Par ailleurs on remarque que $f_3$ est g.d. conjugu\'e \`a
\begin{align*}
(x_0(x_1^2+x_0x_2):x_1(x_1^2+ x_0x_2):x_0x_1^2)
\end{align*}
et que $f_4$ appartient \`a
\begin{align*}
\mathcal{O}_{g.d.}(x_0(x_0^2+x_1x_2):x_1(x_0^2+x_1
x_2):x_0x_1^2).
\end{align*}

Ainsi \`a conjugaison g.d. pr\`es on a les mod\`eles
\begin{align*}
& (x_0(x_1^2+x_0x_2):x_1(x_1^2+x_0x_2):x_0x_1x_2);
&& (x_0(x_1^2+x_0x_2):x_1(x_1^2+x_0x_2):x_0x_1^2);\\
&(x_0(x_0^2+x_1x_2):x_1(x_0^2+x_1 x_2):x_0x_1^2).&&
\end{align*}

\subsection{Configuration $\mathsf{\{11\}}$}\hspace{1mm}

Consid\'erons les \'el\'ements du \S\ref{reduc} ayant
pour configuration $\mathsf{\{11\}}$
\begin{align*}
& f_1=(x_0(x_0x_1+x_1x_2+x_0x_2):x_1(x_0x_1+x_1x_2+x_0x_2):x_0x_1x_2);\\
& f_2=(x_0(x_1^2+x_0x_2+x_1x_2):x_1(x_1^2+x_0x_2+x_1x_2):x_0x_1x_2);\\
& f_3=(x_0(x_0^2+x_0x_2+x_1x_2):x_1(x_0^2+x_0x_2+x_1x_2):x_0x_1x_2).
\end{align*}

Un calcul \'el\'ementaire conduit \`a
\begin{align*}
f_3=(x_1:x_0:x_2)f_2(x_1:x_0:x_2).
\end{align*}

Les transformations $f_1$ et $f_3$ ne sont pas g.d. conjugu\'ees:
les multiplicit\'es des droites de~$\mathrm{Exc} f_3$ sont $3$ et $1$ alors que celles
de $\mathrm{Exc} f_1$ apparaissent avec multiplicit\'e $2.$

\subsection{Configuration $\mathsf{\{12\}}$}\hspace{1mm}

Les transformations de configuration $\mathsf{\{12\}}$
sont
\begin{align*}
& f_1=(x_0(x_1^2+x_0x_1+x_0x_2):x_1(x_1^2+x_0x_1+x_0x_2)
:x_0x_1x_2); &&\\
& f_2=(x_0(x_0^2+x_1^2+2x_0x_1+x_0x_2):x_1(x_0^2+x_1^2+2
x_0x_1+x_0x_2):x_0x_1x_2); &&\\
& f_3=(x_0(x_0^2+\gamma x_0x_1+x_1x_2):x_1(x_0^2+\gamma x_0
x_1+x_1x_2):x_0x_1x_2); &&\\
& f_4=(x_0(x_0^2+x_1^2+2x_0x_1+x_1x_2):x_1(x_0^2+x_1^2+2x_0
x_1+x_1x_2):x_0x_1x_2); &&\\
& f_5=(x_0(x_0x_1+x_1x_2+x_0x_2):x_1(x_0x_1+x_1x_2+x_0x_2):
x_0x_1^2); &&\\
& f_6=(x_0(x_1^2+\gamma x_0x_1+x_0x_2+x_1x_2):x_1(x_1^2+
\gamma x_0x_1+x_0x_2+x_1x_2):x_0x_1^2); &&\\
& f_7=(x_0(x_0^2+\gamma x_0x_1+\delta x_0x_2+x_1x_2):x_1(
x_0^2+\gamma x_0x_1+\delta x_0x_2+x_1x_2):x_0x_1^2); &&\\
& f_8=(x_0(x_0^2+x_1^2+\gamma x_0x_1+\delta x_0x_2+x_1x_2):x_1
(x_0^2+x_1^2+\gamma x_0x_1+\delta x_0x_2+x_1x_2):x_0x_1^2); &&\\
& f_9=(x_0(x_0^2+x_1x_2):x_1(x_0^2+x_1x_2):x_0x_1(x_0-x_1)).&&
\end{align*}

L'\'el\'ement $f_9$ n'est g.d. conjugu\'e \`a aucun des
autres $f_i:$ la droite de $\mathrm{Exc}\hspace{1mm} f_9$ tangente \`a la
conique de $\mathrm{Exc}\hspace{1mm} f_9$ apparait avec multiplicit\'e $2$
ce qui n'est pas le cas pour les autres mod\`eles.

Les huit premi\`eres transformations sont g.d. conjugu\'ees

\begin{small}
\begin{align*}
& f_2=(x_1:x_0:x_2)f_4(x_1:x_0:x_2); && f_7=f_8(x_0:x_1:x_2-x_1);
\end{align*}
\end{small}

\begin{small}
\begin{align*}
& f_1=(\gamma x_1:x_0:x_2)f_3\left(x_1:\frac{x_0}{\gamma}:\gamma
x_2\right);&& f_5=f_6(x_0:x_1:(\gamma-1)x_2-x_1);\\
& f_1=(x_0:-x_0-x_1:-x_0-x_2)f_2(x_0:-x_0-x_1:x_1+x_2); &&f_7=f_5(x_0:x_1:x_0+x_2); \\
& f_5=(x_0+x_1:-x_1:x_2-x_1)f_1(-x_0-x_1:x_1:x_2).&&
\end{align*}
\end{small}

Finalement \`a conjugaison g.d. pr\`es on a deux mod\`eles
\begin{align*}
& (x_0(x_0^2+x_1x_2):x_1(x_0^2+x_1x_2):x_0x_1(x_0-x_1));\\
& (x_0(x_0^2+x_0x_1+x_1x_2):x_1(x_0^2+x_0
x_1+x_1x_2):x_0x_1x_2). &&
\end{align*}

\subsection{Configuration $\mathsf{\{13\}}$}\hspace{1mm}

Les transformations pr\'esentant la configuration
$\mathsf{\{13\}}$ sont
\begin{small}
\begin{align*}
&f_1=(x_0(x_1^2+\gamma x_0x_1+x_1x_2+x_0x_2):x_1(x_1^2+\gamma x_0
x_1+x_1x_2+x_0x_2):x_0x_1x_2), && \gamma\not=0,\hspace{1mm} 1;
\end{align*}
\begin{align*}
&f_2=(x_0(x_0^2+\gamma x_0x_1+x_1x_2+\delta
x_0x_2):x_1(x_0^2+\gamma x_0x_1+x_1x_2+\delta x_0x_2):x_0x_1x_2),
&& \gamma\delta\not=0,\hspace{1mm}\gamma \delta\not=1;
\end{align*}
\begin{align*}
& f_3=(x_0(x_0^2+x_1^2+2x_0x_1+x_1x_2+\delta
x_0x_2):x_1(x_0^2+x_1^2 +2x_0x_1+x_1x_2+\delta x_0x_2):x_0x_1x_2),
&& \delta\not=0,\hspace{1mm} 1,\hspace{1mm} -1.
\end{align*}
\end{small}

Notons que $f_1$ et $f_2$ appartiennent \`a la m\^eme
g.d. orbite
\begin{small}
\begin{align*}
f_1=\left(x_0:x_1:x_2\right)(x_0(x_0^2+\gamma
x_0x_1+x_1x_2+x_0x_2):x_1(x_0^2+\gamma x_0x_1+x_1x_2+
x_0x_2):x_0x_1x_2)\left(x_0:\delta x_1:\frac{x_2}{\delta}\right).
\end{align*}
\end{small}

Les \'el\'ements $f_1$ et $f_3$ sont g.d. conjugu\'es
\begin{small}
\begin{align*}
f_1=A\left(x_0(x_0^2+x_1^2+2x_0x_1+x_1x_2+\frac{\gamma-1}{\gamma}
x_0x_2):x_1(x_0^2+x_1^2 +2x_0x_1+x_1x_2+\frac{\gamma-1}{\gamma}
x_0x_2):x_0x_1x_2\right)B
\end{align*}
\end{small}
avec
\begin{small}
\begin{align*}
& A=\left(\frac{\gamma-1}{\gamma}x_0:(1-\gamma)x_0+(1-\gamma)x_1:
x_0-\frac{x_2}{\gamma}\right) && \text{et} &&B=\left(-\gamma
x_0:\gamma x_0+x_1:\frac{\gamma x_1}{1-\gamma}+\frac{x_2}
{1-\gamma}\right).
\end{align*}
\end{small}

\subsection{Configuration $\mathsf{\{14\}}$}\hspace{1mm}

Les mod\`eles pr\'esentant la configuration
$\mathsf{\{14\}}$ sont les suivants
\begin{align*}
& f_1=(x_0(x_0^2+x_1^2+\gamma x_0x_1+\delta
x_0x_2):x_1(x_0^2+x_1^2+ \gamma x_0x_1+\delta x_0x_2):x_0x_1x_2),
&& \gamma^2\not=4,\hspace{1mm}\delta\not=0;
\end{align*}
\begin{align*}
& f_2=(x_0(x_0^2+x_1^2+\gamma
x_0x_1+x_1x_2):x_1(x_0^2+x_1^2+\gamma x_0x_1+x_1x_2):x_0x_1x_2),
&& \gamma^2\not=4;
\end{align*}
\begin{align*}
& f_3=(x_0(x_0^2+\gamma x_0x_2+x_1x_2):x_1(x_0^2+\gamma x_0x_2+x_1
x_2):x_0x_1(x_0-x_1)), && \gamma\not=0.
\end{align*}

On remarque que $f_1$ et $f_2$ sont g.d. conjugu\'es
\begin{align*}
& f_2=(x_1:x_0:\delta
x_2)f_1\left(x_1:x_0:\frac{x_2}{\delta}\right);
\end{align*}
par contre $f_1$ et $f_3$ ne le sont pas: trois des
points de $f_3(\mathrm{Exc} f_3)$ sont align\'es ce qui n'est pas
le cas pour $f_1(\mathrm{Exc} f_1).$

\subsection{Configuration $\mathsf{\{15\}}$}\hspace{1mm}

Un seul mod\`ele du \S \ref{reduc} poss\`ede la
configuration $\mathsf{\{15\}},$ il s'agit de
\begin{footnotesize}
\begin{align*}
& (x_0(x_0^2+x_1^2+\gamma x_0x_1+\delta x_0x_2+x_1x_2):x_1(
x_0^2+x_1^2+\gamma x_0x_1+\delta x_0x_2+x_1x_2):x_0x_1x_2), &&
\gamma^2\not=4,\hspace{1mm}\delta\not\in\{0,\hspace{1mm} 1\}.
\end{align*}
\end{footnotesize}

\backmatter

\chapter*{Index}

\newlength{\largeur}
\setlength{\largeur}{\textwidth} \addtolength{\largeur}{1cm}

\hspace*{-1.2cm}\begin{tabular}{p{\largeur}p{4.2cm}}

\noteB{alg\'ebriquement stable}{ind20}

\noteB{application quadratique d\'eterminantielle}{ind30a}

\noteB{automorphisme de \textsc{H\'enon}}{ind29}

\noteB{automorphisme de type \textsc{H\'enon}}{ind40a}

\noteB{coefficients alg\'ebriquement ind\'ependants sur
$\mathbb{Q}$}{ind44}

\noteB{condition de \textsc{Bedford} et \textsc{Diller}}{ind22}

\noteB{conjugaison dynamique}{ind24}

\noteB{conjugaison gauche-droite}{ind25}

\noteB{conique r\'eduite}{ind49}

\noteB{contracter une courbe}{ind27b}

\noteB{courbe contract\'ee par une transformation birationnelle}{ind27aa}

\noteB{courbe de points fixes}{ind43}

\noteB{degr\'e d'une transformation rationnelle}{ind2}

\noteB{droite mobile}{ind37}

\noteB{droite immobile}{ind36}

\noteB{\'eclatement d'un point}{ind27a}

\noteB{ensemble exceptionnel d'une transformation birationnelle}{ind27}

\noteB{famille de droites contract\'ees}{ind34}

\noteB{famille de points d'ind\'etermination}{ind35}

\noteB{flot polynomial}{ind40}

\noteB{g\'en\'erateur infinit\'esimal}{ind39}

\noteB{germe de flot}{ind33}

\noteB{groupe affine}{ind41a}

\noteB{groupe de \textsc{Cremona}}{ind0}

\noteB{groupe de d\'ecomposition d'une courbe irr\'eductible de $\mathbb{P}^2(\mathbb{C})$}{ind17}

\noteB{groupe de \textsc{de Jonqui\`eres}}{ind12}

\noteB{groupe d'inertie d'une courbe irr\'eductible de $\mathbb{P}^2(\mathbb{C})$}{ind18}

\noteB{indice de \textsc{Baum}-\textsc{Bott} du feuilletage
$\mathscr{F}$ au point $m$}{ind43b}

\noteB{groupe \'el\'ementaire}{ind41b}

\end{tabular}

\setlength{\largeur}{\textwidth} \addtolength{\largeur}{1cm}

\hspace*{-1.2cm}\begin{tabular}{p{\largeur}p{4cm}}

\noteB{involution de \textsc{Bertini}}{ind16}

\noteB{involution de \textsc{Cremona}}{ind11}

\noteB{involution de \textsc{Geiser}}{ind15}

\noteB{$k$-i\`eme voisinage infinit\'esimal}{ind8}

\noteB{lieu d'ind\'etermination d'une transformation
birationnelle}{ind26}

\noteB{mod\`ele birationnel}{ind32}

\noteB{point base d'une transformation birationnelle}{ind5}

\noteB{point critique d'une application polynomiale
homog\`ene}{ind56}

\noteB{point d'intermination d'une transformation birationnelle}{ind10}

\noteB{point fixe d'une transformation birationnelle}{ind42}

\noteB{point fixe r\'esonnant}{ind45}

\noteB{point p\'eriodique hyperbolique}{ind47}

\noteB{point infiniment proche}{ind7}

\noteB{point immobile}{ind36}

\noteB{point mobile}{ind37}

\noteB{point p\'eriodique d'une transformation birationnelle}{ind46}

\noteB{point propre d'une transformation birationnelle}{ind9}

\noteB{point singulier de type n\oe ud col}{ind31}

\noteB{premier degr\'e dynamique}{ind21}

\noteB{premier voisinage infinit\'esimal}{ind6}

\noteB{purement de degr\'e $k$}{ind23}

\noteB{rationnellement int\'egrable}{ind38}

\noteB{relation lin\'eaire}{ind30}

\noteB{r\'eseau homalo\"idal}{ind4}

\noteB{sous-groupe normal de
$\mathrm{G}$ engendr\'e par l'\'el\'ement $f$ de $\mathrm{G}$}{ind500}

\noteB{surface de \textsc{del Pezzo}}{ind14}

\noteB{sym\'etrie forte}{ind41}

\noteB{transformation quadratique $\mathcal{C}$-g\'en\'erique}{ind19}

\noteB{transformation rationnelle}{ind1}

\noteB{transformation rationnelle non d\'eg\'en\'er\'ee}{ind290}

\noteB{transformation birationnelle de $\mathbb{P}^2(\mathbb{C})$ dans lui-m\^eme,
transformation de \textsc{Cremona}}{ind3}

\noteB{transformation de \textsc{de Jonqui\`eres}}{ind13}

\noteB{transformation monomiale}{ind19a}

\end{tabular}

\chapter*{Index des notations}

\setlength{\largeur}{\textwidth} \addtolength{\largeur}{-1cm}

\hspace*{-1.4cm}\begin{tabular}{p{2cm}p{\largeur}}

\noteAB{\sigma}{involution de \textsc{Cremona}, transformation
birationnelle quadratique d\'efinie par}{}

\noteAC{}{\hspace{4cm}$(x_0:x_1:x_2)\mapsto(x_1x_2:x_0x_2:x_0x_1)$}{not1}

\noteA{\mathrm{dJ}}{groupe de \textsc{de Jonqui\`eres}}{not2}

\noteA{\lambda(f)}{premier degr\'e dynamique de $f$}{not3}

\noteAB{f^{\bullet}}{\`a un \'el\'ement $f=(f_0:f_1:f_2)$ de $\mathrm{Rat}_k$ on
associe la transformation birationnelle}{}

\noteAC{}{\hspace{3cm}$f^\bullet=\delta(f_0:f_1:f_2),\hspace{6mm}
\delta=(\text{pgcd}(f_0,f_1,f_2))^{-1}$}{not3a}

\noteA{Df_{(m)}}{diff\'erentielle de $f$ en $m$}{not3d}

\noteA{\mathrm{tr}(Df_{(m)})}{trace de $Df_{(m)}$}{not3e}

\noteA{\mathrm{Rat}_k}{projectivis\'e de l'espace des triplets de polyn\^omes homog\`enes
de degr\'e $k$ en $3$ variables}{not4}

\noteA{\mathring{\mathrm{R}} \mathrm{at}_k}{ensemble des \'el\'ements de $\mathrm{Rat}_k$
purement de degr\'e $k$}{not5}

\noteA{\mathrm{Rat}}{ensemble des transformations rationnelles de $\mathbb{P}^2(\mathbb{C})$ dans lui-m\^eme}{not6}

\noteA{\mathrm{Bir}(\mathbb{P}^2(\mathbb{C}))}{groupe des transformations
birationnelles de $\mathbb{P}^2(\mathbb{C}),$ groupe de
\textsc{Cremona}}{not7}

\noteA{\mathrm{Bir}_k}{ensemble des transformations de $\mathrm{Rat}_k$ inversibles en tant que
transformations rationnelles}{not8}

\noteAB{\mathrm{Rat}_k^{\bullet}}{ensemble d\'efini par}{}

\noteAC{}{\hspace{5cm}$\{f^\bullet\hspace{1mm}|\hspace{1mm} f\in\mathrm{Rat}_k\}$}{not3b}

\noteAB{\mathrm{Bir}_k^{\bullet}}{ensemble d\'efini par}{}

\noteAC{}{\hspace{5cm}$\{f^\bullet\hspace{1mm}|\hspace{1mm} f\in\mathrm{Bir}_k\}$}{not3c}

\noteA{\mathring{\mathrm{B}}\mathrm{ir}_k}{ensemble des
transformations birationnelles purement de degr\'e $k$}{not9}

\noteA{\mathcal{O}_{dyn} (Q)}{orbite d'un \'el\'ement
$Q$ de $\mathrm{Rat}$ sous l'action de $\mathrm{PGL}_3(\mathbb{C})$ par conjugaison
dynamique}{not10}

\noteA{\mathcal{O}_{g.d.} (Q)}{orbite d'une transformation
rationnelle $Q$ sous l'action gauche-droite}{not11}

\noteA{\mathrm{Ind}\hspace{1mm} f}{lieu d'ind\'etermination de $f$}{not12}

\noteA{\mathrm{Exc}\hspace{1mm} f}{ensemble exceptionnel de $f$}{not13}

\noteAB{\rho}{transformation birationnelle quadratique donn\'ee
par}{}

\noteAC{}{\hspace{4cm}$(x_0:x_1:x_2)\mapsto(x_0x_1:x_2^2:x_1x_2)$}{not14}

\noteAB{\tau}{transformation birationnelle quadratique d\'efinie
par}{}

\noteAC{}{\hspace{4cm}$(x_0:x_1:x_2)\mapsto(x_0^2:x_0x_1:x_1^2-x_0x_2)$}{not15}

\noteA{\mathbb{C}[x_0,x_1,x_2]_\nu}{ensemble
des polyn\^omes homog\`enes de degr\'e $\nu$ dans $\mathbb{C}^3$}{not16}

\noteAB{\mathrm{det}\hspace{1mm}\mathrm{jac}}{application d\'efinie
par}{}

\noteAC{}{\hspace{1.8cm}$\mathrm{Rat}_2\to\mathbb{P}(\mathbb{C}[x_0,x_1,x_2]_3)
\simeq\{\text{ courbes de degr\'e $3$ }\},\hspace{1mm}
\hspace{1mm}[Q]\mapsto[\mathrm{det}\hspace{1mm}\mathrm{jac}\hspace{1mm} Q=0]$}{not16b}

\noteA{\mathcal{C}(f)}{lieu critique de $f$}{not16c}

\noteA{\Sigma^3}{orbite de $\sigma$ sous l'action gauche-droite, {\it i.e.}
$\Sigma^3=\mathcal{O}_{g.d.}(\sigma)$}{not17}

\noteA{\Sigma^2}{orbite de $\rho$ sous l'action gauche-droite, {\it i.e.}
$\Sigma^2=\mathcal{O}_{g.d.}(\rho)$}{not18}

\noteA{\Sigma^1}{orbite de $\tau$ sous l'action gauche-droite, {\it i.e.}
$\Sigma^1=\mathcal{O}_{g.d.}(\tau)$}{not19}
\end{tabular}

\setlength{\largeur}{\textwidth} \addtolength{\largeur}{-1cm}

\hspace*{-1.4cm}\begin{tabular}{p{2cm}p{\largeur}}

\noteAB{\Sigma^0}{ensemble des transformations du type}{}

\noteAC{}{\hspace{2cm}
$\ell(\ell_0:\ell_1:\ell_2),$ $\ell,\hspace{1mm}\ell_i$
formes lin\'eaires, les $\ell_i$ \'etant ind\'ependantes}{not20}

\noteA{\mathrm{RL}(Q)}{$\mathbb{C}$-espace vectoriel des relations lin\'eaires de $Q$}{not22}

\noteA{\mathrm{e}(Q)}{dimension de $\mathrm{RL}(Q)$}{not23}

\noteA{\mathrm{Sing}\hspace{1mm}\mathscr{F}}{lieu singulier
du feuilletage $\mathscr{F}$}{not23b}

\noteA{{\rm Isot}\hspace{1mm}f}{groupe d'isotropie de $f$}{not24}

\noteA{\langle f_1,\hspace{1mm} \ldots,\hspace{1mm} f_n\rangle}{groupe engendr\'e par $f_1,$ $\ldots,$ $f_n$}{not250}

\noteAB{\mathscr{S}_6}{groupe donn\'e par}{}

\noteAC{}{\hspace{1.5cm}$\mathscr{S}_6=\{\mathrm{id},\hspace{1mm}(x_0:x_2:x_1),\hspace{1mm}(x_2:x_1:x_0),\hspace{1mm}(x_1:x_0:x_2),\hspace{1mm}(x_1:x_2:x_0),\hspace{1mm}(x_2:x_0:x_1)\}$}{not25}

\noteAB{\mathrm{A}}{groupe affine}{}

\noteAC{}{\hspace{1cm}$\mathrm{A}=\{(x_0,x_1)\mapsto(a_1x_0+b_1x_1+
c_1,a_2x_0+b_2x_1+c_2)\hspace{1mm}|\hspace{1mm}a_i,\hspace{1mm} b_i,\hspace{1mm}c_i\in\mathbb{C},\hspace{1mm}
a_1b_2-a_2b_1 \not=0\}$}{not25a}

\noteAB{\mathrm{E}}{groupe \'el\'ementaire}{}

\noteAC{}{\hspace{2cm} $\mathrm{E}=\{(x_0,x_1)\mapsto(\alpha x_0+P(x_1),\beta
x_1+ \gamma)\hspace{1mm}|\hspace{1mm}\alpha,\hspace{1mm}\beta\in\mathbb{C}^*,
\hspace{1mm}\gamma\in\mathbb{C},\hspace{1mm} P\in
\mathbb{C}[x_1]\}$}{not25b}

\noteAB{\mathcal{F}_\nu}{projectivis\'e de l'espace vectoriel des formes
$\omega$ satisfaisant l'identit\'e d'\textsc{Euler}:}{}

\noteAC{}{\hspace{1.5cm}$\mathcal{F}_\nu=\mathbb{P}\{\omega=F_0dx_0+F_1dx_1+F_2dx_2\hspace{1mm}|\hspace{1mm}
x_0F_0+x_1F_1+x_2F_2=0,\hspace{1mm}\deg F_i=\nu+1\}$}{not26}

\noteA{\mathcal{F}_\nu^\bullet}{espace des feuilletages de degr\'e
inf\'erieur ou \'egal \`a $\nu$}{not26a}

\noteA{\mathring{\mathcal{F}}_\nu}{ensemble des feuilletages de degr\'e
$\nu$}{not26aa}

\noteAB{\mathcal{F}(.)}{application donn\'ee par}{}

\noteAC{}{
\hspace{7mm}$\mathrm{Rat}_n\to\mathcal{F}_n,\hspace{3mm}
f=(f_0:f_1:f_2)\mapsto(x_1f_2-x_2f_1)dx_0+(x_2f_0-x_0f_2)dx_1+
(x_0f_1-x_1f_0)dx_2$}{not26b}

\noteA{\mathrm{Fix}\hspace{1mm}f}{ensemble des points fixes de $f$}{not27}

\noteA{\mathrm{Re}\hspace{1mm}z}{partie r\'eelle de $z$}{not27ab}

\noteA{\mathrm{BB}(\mathscr{F}(m))}{indice de \textsc{Baum}-\textsc{Bott} du
feuilletage $\mathscr{F}$ au point $m$}{not27b}

\noteA{\mathrm{Aut}(\mathbb{C},+,.)}{groupe des automorphismes du corps $\mathbb{C}$}{not28}

\noteAB{A^\kappa}{\'el\'ement obtenu en faisant agir
l'\'el\'ement $\kappa$ de $\mathrm{Aut}(\mathbb{C},+,.)$ sur les coefficients de}{}

\noteAC{}{l'automorphisme $A$ de $\mathbb{P}^2(\mathbb{C})$}{not29}

\noteAB{\mathrm{Ind}^+ \hspace{1mm}f}{union des ensembles
d'ind\'etermination des it\'er\'es positifs de $f:$}{}

\noteAC{}{
\hspace{5cm}$\bigcup_{n\geq 1}\mathrm{Ind}\hspace{1mm}f^n$}{not39}

\noteAB{\mathrm{Ind}^- \hspace{1mm}f}{union des ensembles
d'ind\'etermination des it\'er\'es n\'egatifs de $f:$}{}

\noteAC{}{
\hspace{5cm}$\bigcup_{n\geq 1}\mathrm{Ind}\hspace{1mm} f^{-n}$}{not40}

\noteAB{\mathrm{Exc}^+ \hspace{1mm}f}{union des ensembles
exceptionnels des it\'er\'es positifs de $f:$}{}

\noteAC{}{\hspace{5cm}
$\bigcup_{n\geq 1}\mathrm{Exc}\hspace{1mm} f^n$}{not41}

\noteAB{\mathrm{Exc}^- \hspace{1mm} f}{union des ensembles
exceptionnels des it\'er\'es n\'egatifs de $f:$}{}

\noteAC{}{\hspace{5cm}
$\bigcup_{n\geq 1}\mathrm{Exc}\hspace{1mm} f^{-n}$}{not42}

\noteAB{f_{\alpha,\beta}}{famille de transformations birationnelles d\'efinie
par}{}

\noteAC{}{
\hspace{4.7cm}$\left(\frac{\alpha x_0+x_1}{x_0+1},\beta x_1\right),$ \hspace{3mm}
$\alpha,\hspace{1mm}\beta\in\mathbb{C}^*$}{not35aa}

\noteA{\textrm{dist}}{m\'etrique de
\textsc{Fubini}-\textsc{Study}}{not35bb}

\noteA{\mathrm{N}(f,\mathrm{G})}{sous-groupe normal de
$\mathrm{G}$ engendr\'e par l'\'el\'ement $f$ de $\mathrm{G}$}{not35ba}

\noteA{\textrm{Inv}}{application qui \`a une transformation
birationnelle associe son inverse}{not35}

\end{tabular}

\setlength{\largeur}{\textwidth} \addtolength{\largeur}{-1cm}

\hspace*{-1.4cm}\begin{tabular}{p{2cm}p{\largeur}}

\noteAB{\mathfrak{(a)}}{transformation du type }{}

\noteAC{}{\hspace{7mm}
$(x_0(\alpha x_0^2+\beta x_1^2+\gamma x_0x_1
+\delta x_0x_2+\varepsilon x_1x_2):x_1
(\alpha x_0^2+\beta x_1^2+\gamma x_0x_1+\delta x_0x_2+\varepsilon
x_1x_2):x_0x_1x_2)$}{not31}

\noteAB{\mathfrak{(b)}}{transformation de la forme}{}

\noteAC{}{\hspace{0.8cm}
$(x_0(\alpha x_0^2+\beta x_1^2+\gamma x_0x_1+
\delta x_0x_2+\varepsilon x_1x_2):x_1
(\alpha x_0^2+\beta x_1^2+\gamma x_0x_1+\delta x_0x_2+\varepsilon
x_1x_2):x_0x_1^2)$}{not32}

\noteAB{\mathfrak{(c)}}{transformation du type }{}

\noteAC{}{\hspace{2.7cm}$(x_0(x_0^2+x_1x_2+\gamma x_0x_2):x_1(x_0^2+
x_1x_2+\gamma x_0x_2):x_0x_1(x_0-x_1))$}{not33}

\noteAB{\mathfrak{(d)}}{transformation de la forme}{}

\noteAC{}{\hspace{3.7cm}$(x_0x_2(x_1+\gamma x_0):x_1x_2(x_1+\gamma x_0):
x_0x_1(x_0-x_1))$}{not34}
\end{tabular}
\vspace{8mm}

\bibliographystyle{alpha}
\bibliography{texte}
\nocite{*}

\end{document}

%% file: des9scan.pstex_t
\begin{picture}(0,0)%
\includegraphics{des9scan.pstex}%
\end{picture}%
\setlength{\unitlength}{3947sp}%
\begingroup\makeatletter\ifx\SetFigFont\undefined%
\gdef\SetFigFont#1#2#3#4#5{%
  \reset@font\fontsize{#1}{#2pt}%
  \fontfamily{#3}\fontseries{#4}\fontshape{#5}%
  \selectfont}%
\fi\endgroup%
\begin{picture}(1200,1200)(1201,-1561)
\end{picture}%

%% file: iter8scan.pstex_t
\begin{picture}(0,0)%
\includegraphics{iter8scan.pstex}%
\end{picture}%
\setlength{\unitlength}{3947sp}%
\begingroup\makeatletter\ifx\SetFigFont\undefined%
\gdef\SetFigFont#1#2#3#4#5{%
  \reset@font\fontsize{#1}{#2pt}%
  \fontfamily{#3}\fontseries{#4}\fontshape{#5}%
  \selectfont}%
\fi\endgroup%
\begin{picture}(1200,1200)(1201,-1561)
\end{picture}%

%% file: julia49a.pstex_t
\begin{picture}(0,0)%
\includegraphics{julia49a.pstex}%
\end{picture}%
\setlength{\unitlength}{3947sp}%
\begingroup\makeatletter\ifx\SetFigFont\undefined%
\gdef\SetFigFont#1#2#3#4#5{%
  \reset@font\fontsize{#1}{#2pt}%
  \fontfamily{#3}\fontseries{#4}\fontshape{#5}%
  \selectfont}%
\fi\endgroup%
\begin{picture}(2400,1950)(1201,-2311)
\end{picture}%

%% file: bd48.pstex_t
\begin{picture}(0,0)%
\includegraphics{bd48.pstex}%
\end{picture}%
\setlength{\unitlength}{3947sp}%
\begingroup\makeatletter\ifx\SetFigFont\undefined%
\gdef\SetFigFont#1#2#3#4#5{%
  \reset@font\fontsize{#1}{#2pt}%
  \fontfamily{#3}\fontseries{#4}\fontshape{#5}%
  \selectfont}%
\fi\endgroup%
\begin{picture}(2400,1950)(1201,-2311)
\end{picture}%